\documentclass[11pt,letterpaper]{article}
\usepackage[margin=1.05in]{geometry}
\usepackage[T1]{fontenc}
\usepackage{lmodern,microtype,amsmath,amssymb,amsthm,mathtools}
\usepackage{booktabs,longtable,array,enumitem,mdframed,xurl}
\usepackage[hidelinks,hypertexnames=false]{hyperref}
\makeatletter
\renewcommand*\l@subsection{\@dottedtocline{2}{1.5em}{3em}}
\makeatother
\newtheorem{theorem}{Theorem}[subsection]
\newtheorem{proposition}[theorem]{Proposition}
\newtheorem{lemma}[theorem]{Lemma}
\newtheorem{corollary}[theorem]{Corollary}
\theoremstyle{definition}

\theoremstyle{remark}
\newtheorem{remark}[theorem]{Remark}

\numberwithin{equation}{subsection}

\hypersetup{pdftitle={Critical characteristic-polynomial chaos for complex iid matrices},pdfauthor={Yin Cai, Xiang Fang, Feng Guo}}
\begin{document}
\title{Critical characteristic-polynomial chaos\\for complex iid matrices}
\author{Yin Cai\thanks{School of Mathematics, Hangzhou Normal University,
Hangzhou 311121, P. R. China.\\
Email: \mbox{\href{mailto:cymath@hznu.edu.cn}{\nolinkurl{cymath@hznu.edu.cn}}}.}
\and Xiang Fang\thanks{Department of Applied Mathematics,
National Yang Ming Chiao Tung University, Hsinchu 30010, Taiwan.\\
Email: \mbox{\href{mailto:xfang@nycu.edu.tw}{\nolinkurl{xfang@nycu.edu.tw}}}.}
\and Feng Guo\thanks{School of Mathematics, South China University of Technology,
Guangzhou 510640, P. R. China.\\
Email: \mbox{\href{mailto:70207994@nuaa.edu.cn}{\nolinkurl{70207994@nuaa.edu.cn}}}.}}
\date{}
\maketitle
 
\begin{abstract}
We prove convergence to critical Gaussian multiplicative chaos for
characteristic-polynomial measures of complex iid matrices.
Let $A_N=N^{-1/2}(\xi_{ij})$, with a fixed law satisfying
$\mathbb E\xi=\mathbb E\xi^2=0$, $\mathbb E|\xi|^2=1$,
and finite moments of every order; no density is assumed.
On each fixed disk $\{|z|\le r\}$, $0<r<1$,
expectation-normalized determinant measures at power $2\sqrt2$,
multiplied by $\sqrt{\tfrac12\log N}$, converge to critical chaos
with a fourth-cumulant covariance correction.
The normalized height-deficit marked measure converges to the product
of this spatial limit and the standard Rayleigh law.
Convergence holds jointly with centered averages of the
$\sqrt2$-scaled original logarithmic field against finitely many
smooth compactly supported bulk tests, and with the two-sided
near-critical curve at powers $q=2\sqrt2-2\lambda/\sqrt{\log N}$.
On fixed compact parameter intervals, the curve converges in law
in the uniform weak topology; all limiting coordinates use
the same Gaussian field.
The limiting curve is an explicit deterministic profile times the
same critical measure.
\end{abstract}


\clearpage
\tableofcontents
\clearpage

\section{Introduction}
\label{chap:introduction}

In this paper we show that, after suitable normalization, critical and
near-critical characteristic-polynomial measures of complex iid random
matrices converge jointly to deterministic multiples of the same
critical Gaussian multiplicative chaos measure.
The paper is organized around a single main result,
Theorem~\ref{main:IMW}, which also includes the height marks and
averages of the original logarithmic field, with all limiting
objects constructed from the same Gaussian field.

\subsection*{The characteristic-polynomial measure}

 {We first explain how to turn the characteristic polynomial into a random spatial measure and explain why the critical exponent requires a different normalization.}
The characteristic polynomial associates a random function of a complex
variable to a random matrix. If the eigenvalues are
$\lambda_1,\ldots,\lambda_N$, then
\[
 |\det(A_N-zI)|=\prod_{j=1}^N|z-\lambda_j|.
\]
Its logarithm is therefore a sum of logarithmic potentials. Taking a
positive power of the absolute determinant gives a spatial weight that
favors points where this function is large. We study the random measure
formed by this weight at its critical exponent. The question is to
identify the limiting measure, including its dependence on the field
from which it is constructed.

The matrices are $A_N=N^{-1/2}(\xi_{ij})$, with independent copies of a
fixed complex random variable satisfying
$\mathbb E\xi=\mathbb E\xi^2=0$, $\mathbb E|\xi|^2=1$, and finite
moments of every order. No density is assumed, and the real and imaginary
parts need not be independent. We work on a fixed closed disk
$K=\{|z|\le r_0\}$, $0<r_0<1$, inside the unit disk. For a positive
field exponent $\gamma$, we write
\[
 \mu_{N,\gamma}(\mathrm dz)
 =\frac{|\det(A_N-zI)|^{\sqrt2\gamma}}
        {\mathbb E|\det(A_N-zI)|^{\sqrt2\gamma}}\,\mathrm dA(z),
 \qquad z\in K,
\]
where $\mathrm dA$ is planar Lebesgue measure. The denominator makes
the expected density equal to one; it does not make the total mass a
deterministic constant. We use $\gamma$ for the exponent of the
$\sqrt2$-scaled logarithmic field. Thus the corresponding power of the
absolute determinant is $q=\sqrt2\gamma$.

Logarithmic correlations explain why such measures are connected with
Gaussian multiplicative chaos (GMC). In the bulk, the limiting covariance of
the logarithmic field grows like the logarithm of inverse distance as
two points approach one another; the precise iid covariance is recorded
by Cipolloni--Landon \cite[(1.7)]{CL}. Gaussian multiplicative chaos
gives a meaning to the exponential of this rough Gaussian field by
first smoothing it, normalizing the exponential and then removing the
smoothing. The subcritical and critical constructions have different
normalizations \cite[Theorems 5 and 10]{DRSVSH},
\cite[Theorems 5.3--5.5]{JSW}. In our exponent convention the critical
value is $\gamma=2$. For every fixed $0<\gamma<2$, convergence of the
original expectation-normalized iid determinant measure follows from
\cite[Definition 1.1 and Theorem 1.2]{CL}.

At the endpoint, the additional normalization becomes part of the
question. We put $t=\frac12\log N$. The critical measure studied here is
\[
 Z_N=\sqrt t\,\mu_{N,2}.
\]
The factor $\sqrt t$ is the square-root normalization associated with
critical Gaussian chaos; its conversion to the derivative-normalized
measure is specified below. A fixed subcritical convergence theorem
does not by itself justify passing to this endpoint or letting the
exponent vary with $N$. Both operations require control on the scale
at which critical mass is retained. Even the mean mass illustrates
the change: $\mathbb E Z_N(K)=\sqrt t\,|K|$, so the ordinary first
moment gives no bound uniform in $N$. The relevant convergence concerns
the random finite measure, with its stated topology and normalization.

The spatial measure describes where the critical weight is located.
We also ask which field values carry its mass, and record their
normalized deficits below the critical height scale in a marked measure.
These marks describe how the weights respond to a change of exponent
of order $1/\sqrt{\log N}$. Joint convergence with averages of the
original logarithmic field identifies the dependence among these
observations, which separate marginal limits do not determine. The resulting
statement connects the critical spatial measure, its height marks and the
nearby exponent curve on each fixed compact interval of the rescaled
parameter through one limiting Gaussian field. The complete
formal statement is Theorem~\ref{main:IMW} in
Section~\ref{sp:common}.

\subsection*{The joint critical limit}

We now describe the conclusions of Theorem~\ref{main:IMW}; that is, {we state the joint limit of the critical measure, its height marks, the field averages, and the curve of nearby measures.}
Its complete
definitions and formal joint statement are in Section~\ref{sp:common}.
We write $\kappa=\mathbb E|\xi|^4-2\ge-1$. The limiting generalized
Gaussian field $G_\kappa$ has covariance kernel
\[
 K_\kappa(z,w)=-\log|z-w|
       +\frac\kappa2(1-|z|^2)(1-|w|^2).
\]
Let $D_\kappa$ be its derivative-normalized critical chaos, with the
normalization specified in Section~\ref{sp:common}, and put
$M_\kappa=\sqrt{2/\pi}\,D_\kappa$. The critical spatial measure $Z_N$
converges in distribution to $M_\kappa|_K$. Thus the critical measure constructed from the Gaussian field is also
the limit of the original iid determinant weights, with the conversion
to derivative normalization specified above. The smooth fourth-cumulant
term remains in the covariance: different entry laws need not give the
same field as complex Ginibre. Section~\ref{chap:gaussian} constructs
the field and chaos together for the full range $\kappa\ge-1$.

To describe the field values carrying this mass, we use the terminal
field before the fourth-cumulant correction and its normalized deficit
\[
 Y_N(z)=\sqrt2\left[\log|\det(A_N-zI)|
                         +\frac N2(1-|z|^2)\right],\qquad
 R_N(z)=\frac{2t-Y_N(z)}{\sqrt t}.
\]
The marked finite measure $\Xi_N(\mathrm dz,\mathrm du)
=Z_N(\mathrm dz)\delta_{R_N(z)}(\mathrm du)$ converges jointly with
$Z_N$ to $M_\kappa|_K\otimes\rho_{\mathrm R}$, where
\[
 \rho_{\mathrm R}(\mathrm du)
       =u e^{-u^2/2}\mathbf1_{\{u>0\}}\,\mathrm du.
\]
Thus the mass has a Rayleigh distribution in the deficit coordinate,
with the same random spatial factor. The statement concerns a
mass-weighted normalized deficit, not the unscaled height at a uniformly
chosen spatial point. At determinant zeros the density vanishes and
the mark is set to zero. The corrected field used inside the proof is
$Y_N^0=Y_N+c_\kappa$, with
$c_\kappa(z)=\sqrt2\kappa(1-|z|^2)^2/4$. This deterministic function
cancels from the expectation-normalized density and changes the
normalized mark by $-c_\kappa/\sqrt t$, which tends uniformly to zero
on $K$. The theorem itself uses the raw mark $R_N$.

The limiting field is retained in this convergence. For every prescribed
finite family of real smooth compactly supported tests $\varphi$ in
the unit disk, consider
\[
 \mathcal Y_N(\varphi)
 =\int\varphi(z)Y_N(z)\,\mathrm dA(z)
    -\mathbb E\!\int\varphi(z)Y_N(z)\,\mathrm dA(z).
\]
These coordinates converge jointly with the spatial and marked measures
to $G_\kappa(\varphi)$. The integral precedes centering, so the definition
also covers discrete entry laws without asking for a finite pointwise
mean of the zero-height logarithm. Test supports may extend beyond
$K$ while remaining inside the unit disk. The chaos and all limiting
field coordinates are constructed from the same $G_\kappa$.
This finite-family statement specifies their dependence; it makes no
additional assertion about convergence in a distribution-space topology.

For a fixed $A<\infty$, we set
\[
 \gamma_N(\lambda)=2-\frac\lambda{\sqrt t},\qquad
 Z_N(\lambda)=\sqrt t\,\mu_{N,\gamma_N(\lambda)},
 \qquad |\lambda|\le A.
\]
The exponents are positive for all sufficiently large $N$. The whole
curve converges in distribution in
$C([-A,A],\mathcal M_+(K))$, equipped with the uniform
bounded-Lipschitz metric, to $(F(\lambda)D_\kappa|_K)_{|\lambda|\le A}$,
where
\[
 F(\lambda)=\sqrt{\frac2\pi}\,e^{-\lambda^2/2}
          \int_0^\infty u e^{-u^2/2+\lambda u}\,\mathrm du.
\]
Here $\mathcal M_+(K)$ is the space of finite positive measures with
its weak topology, including total mass. The marked measure also uses
the weak finite-measure topology on $K\times\mathbb R$. All the
spatial, marked, field-test and curve conclusions hold jointly in the
product space stated in Theorem~\ref{main:IMW}. In particular,
$F(0)=\sqrt{2/\pi}$, consistent with the critical spatial limit.
Positive $\lambda$ approaches criticality from below and negative
$\lambda$ from above. The interval is fixed: neither a growing window
nor a fixed supercritical exponent is included.

The connection between the mark and the curve is an exact change of
exponent. If $\psi_N(\gamma,z)=\log\mathbb E e^{\gamma Y_N(z)}$, then
\[
 Z_N(\lambda)(\mathrm dz)
 =e^{\psi_N(2,z)-\psi_N(\gamma_N(\lambda),z)-2\lambda\sqrt t}
       e^{\lambda R_N(z)}Z_N(\mathrm dz).
\]
The first factor tends to $e^{-\lambda^2/2}$ uniformly on the fixed
spatial and parameter sets, as quantified in Proposition~\ref{cw:inputs}.
Integrating $e^{\lambda u}$ against the limiting Rayleigh kernel then
gives $F(\lambda)D_\kappa$. This step requires the uniform exponential
mark-tail estimate in probability as well as the scalar normalization
asymptotic. Weak convergence of the marked finite measure alone controls
bounded tests and does not justify this unbounded tilt. The required estimate
and the deduction of the continuous curve are in Section~\ref{cw:subsection}.

\subsection*{Relation to previous work}

Cipolloni--Landon \cite[Definition 1.1 and Theorem 1.2]{CL} prove weak
convergence in law of the original expectation-normalized determinant
measures for each fixed $0<q<2\sqrt2$. Their complex iid entry
assumptions agree with ours; they also treat the real symmetry class.
Their sharp one-point moment asymptotics, however, allow every fixed
nonnegative determinant power, including the critical power; the
subcritical restriction concerns the random-measure convergence theorem.
That theorem does not include the endpoint or the joint critical limit
of the measure, its marks and the compact near-critical curve. Their Proposition~7.5 already contains mixed
information involving regularized field observations, and their
Corollary~1.7 gives leading free-energy freezing. The present
contribution is the critical random-measure limit with the coupling,
marks and fixed compact curve specified above. The critical direction
was highlighted in Cipolloni's June 2026 Oberwolfach report
\cite[p.~1497]{OW}. The barrier-truncation approach has antecedents in Berestycki
\cite[Sections~2--4]{Berestycki} and its extension to asymptotically
Gaussian fields by Lambert--Ostrovsky--Simm
\cite[Section~2.2]{LOS}. The general subcritical convergence framework
of Claeys--Fahs--Lambert--Webb \cite[Theorem 2.6]{CFLW} uses a positive
integrability gap that closes at criticality; see their Lemma~2.8.
We reach the endpoint through determinant-weighted estimates
and removal of the barrier in probability.

Bourgade--Dubach--Hartung--Keles \cite[Theorem 1.7]{BDHK} obtain the
strict subcritical determinant-chaos limit for random normal matrix
eigenvalue ensembles under their potential assumptions (A1)--(A4),
including complex Ginibre. Varying the confining potential and varying
an iid entry law are distinct extensions of that common example.
Their Theorems~1.1--1.2 also include smooth linear statistics, and
Corollary~1.13 gives leading free-energy freezing. Earlier, Lambert
\cite[(2.2)--(2.3), Theorem 2.2]{Lam} proved a strict subcritical
Ginibre chaos limit for a mesoscopically smoothed logarithmic field,
with smoothing $\epsilon=N^{-1/2+\alpha}$, fixed $0<\alpha<1/2$, and
a fixed smoothing subtracted. That theorem concerns a regularized
field on a fixed small bulk disk, a distinction from the original
determinant density considered here.

The original logarithmic-field marginal also has established
precedents. For complex Ginibre, Webb
\cite[Propositions 22 and 24]{Webb}, refining Rider--Vir\'ag
\cite[Corollary~2]{RV},
proved convergence of the centered log characteristic polynomial to
a Gaussian field, including its constant mode. For general iid matrices, the marginal
limit of the smooth field averages also follows from
\cite[Theorem 2.2 and Corollary 2.4]{CES}, applied to their
logarithmic-potential tests with an exterior cutoff. Here we retain these averages
jointly with the critical measure, its marks and the compact near-critical
curve for general complex iid laws.

The Gaussian limiting objects and their normalization have established
precedents. The square-root normalization and its factor
$\sqrt{2/\pi}$ relative to derivative chaos are proved in the Gaussian
settings of Duplantier--Rhodes--Sheffield--Vargas
\cite[Theorems 5 and 10]{DRSVSH}. The continuum subcritical-to-critical
relation was proved for the planar zero-boundary Gaussian free field
by Aru--Powell--Sep\'ulveda \cite[Theorem 1.1]{APSCritical}, and for
logarithmic covariance with a sufficiently regular remainder by
Junnila--Saksman--Webb \cite[Theorem 5.5]{JSW}. In derivative
normalization, the continuum subcritical chaos divided by $2-\gamma$
tends to twice the corresponding derivative-normalized critical measure
as $\gamma\uparrow2$.
The canonical construction used here is supplied by
\cite[Theorems 5.3--5.4]{JSW}. The continuum relation concerns the
limiting subcritical measures. Its proof already uses approximation
depths of order $(2-\gamma)^{-2}$ and exponential reweighting of
critical endpoint marks, followed by truncation; see
\cite[proof of Lemma~4.1, (4.3)--(4.4)]{APSCritical} and
\cite[Section~3.2, (3.17)--(3.18)]{PowellReview}.
Here this mechanism is applied to the original iid determinant
measures, using the marked critical limit, uniform control of exact
normalizers, and tail bounds in probability for both signs of the
window parameter.

The Rayleigh endpoint also has a Gaussian precedent. Aru--Powell--Sep\'ulveda
\cite[Theorem 6.7]{APS} obtain a Brownian-meander endpoint mark for
local-set approximations of the planar Gaussian free field. Their
displayed mark has an additional $\sqrt2$ in its scale; after rescaling,
the endpoint is standard Rayleigh. Their theorem uses bounded
continuous mark tests. The exponential tests required for the present
window need the separate tail estimate stated above.

The branching precedents include the square-root normalization of
A\"id\'ekon--Shi \cite[Theorem 1.1]{AS} and the critical Gibbs
trajectory limit of Madaule \cite[Theorem 1.2]{Madaule}.
Pain \cite[Theorems 1.1(ii)--(iii) and 1.2(ii)--(iii)]{Pain} studies a
two-sided window of order $n^{-1/2}$ for branching random walk and
obtains exponential tilts of a Brownian-meander endpoint. In his notation,
$W_{n,\beta}=\sum_{|v|=n}e^{-\beta V(v)}$, where $|v|$ and $V(v)$
are the generation and position of the particle $v$, and
$\Psi(\beta)=\log\mathbb E W_{1,\beta}$.
Under his assumptions (1.1)--(1.3) and stronger moment hypothesis (1.6), set
$\beta_n(\lambda)=1-\lambda/(\sigma\sqrt n)$. Here
$\sigma^2=\Psi''(1)>0$ and $D_\infty$ denotes the branching
derivative-martingale limit. His scalar limits,
including the critical value supplied by \cite[Theorem 1.1]{AS},
imply convergence of the entire scalar curve on every fixed compact
parameter interval. Indeed, the raw functions
$\sigma\sqrt n\,W_{n,\beta_n(\lambda)}$ are convex in $\lambda$;
convergence on a countable dense set therefore gives locally uniform
convergence in probability under the survival-conditioned law.
Only afterward is expectation normalization applied: analyticity of
$\Psi$ near one gives $n\Psi(\beta_n(\lambda))\to\lambda^2/2$
uniformly on compact intervals. Thus the normalized curve has limit
$F(\lambda)D_\infty$, with the same profile $F$ as above.
The expectation in this normalization is under the original branching
law. This compact-curve conclusion is a consequence of the cited
scalar results and convexity. The matrix result here additionally
identifies the spatial measure, raw deficit marks and prescribed
field averages jointly, with all limiting quantities constructed from
the same Gaussian field.

\subsection*{Proof ideas}

 {Next we follow the proof from a fixed barrier to the joint limit, keeping track of the mass removed at each step.}
In general, convergence of a regularized field does not by itself control the large
values carrying the critical mass. We must also identify how this mass
depends on the limiting field. We first restrict the field at successive
scales, compare the restricted matrix measure with a Gaussian measure,
and then show that the restriction loses negligible mass in probability.
The height information obtained in this comparison will also be used
to identify the measures at nearby powers.

\paragraph{The Gaussian target.}
The first task is to specify a Gaussian comparison measure together with
the field from which it is constructed. We construct both on one
probability space for the full range $\kappa\ge-1$; see
Section~\ref{gaussian:gp2:canonical}. The history barrier
keeps only points whose regularized field stays below $2k+B$ at
successive integer logarithmic scales $k$, with $B$ fixed; larger $k$
means finer smoothing. On the matrix side the condition is
$X_{N,k}(z)\le2k+B$, where $X_{N,k}$ denotes the corrected regularized
field. We impose the same condition on the corresponding Gaussian
observations $G_k(z)$; see \eqref{sp:barriers}.
For a fixed starting scale $k_0$, we write
\[
 \begin{aligned}
 E_\infty^B
   &=\{z\in K:G_k(z)\le2k+B
          \text{ for every integer }k\ge k_0\},\\
 Z_B&=\mathbf1_{E_\infty^B}\,M_\kappa|_K.
 \end{aligned}
\]
Thus $Z_B$ is the restriction of the critical Gaussian measure
to points satisfying the barrier at all these scales.

For fixed $B$ and each real continuous spatial test $f$ on $K$, the
integral $Z_B(f)$ is square integrable and is determined by the
integer-scale Gaussian observations
(Proposition~\ref{p1:proposition}). Therefore, for every
$\varepsilon>0$, there are a finite list $D_G$ of these observations
and a bounded smooth function $\psi$, with all derivatives bounded,
such that
\[
 \mathbb E\bigl|Z_B(f)-\psi(D_G)\bigr|^2<\varepsilon;
\]
see \eqref{p1:finite-approximation}.
The list and function are fixed before $N\to\infty$.
The list can be enlarged to retain any prescribed fixed collection
of field values at fixed positive regularization scales in the bulk.
We must now show that $\psi(D_N)$, where $D_N$ is the matching list
of matrix observations, approximates the restricted matrix mass,
as in \eqref{sp:fb:projection}.

\paragraph{Identification at a fixed barrier.}
Convergence of the matrix observations alone does not control
expectations containing the determinant weights. We need comparisons
of these weighted expectations and a second-moment bound with an
additional temporary cutoff, explained below
(Proposition~\ref{p2:matrix}). These estimates
will control the squared approximation error and thereby identify
the dependence on the field.

The estimates build on the Hermitization and dynamical multiscale
approach of Cipolloni--Landon \cite[Section 1.2.1]{CLMax}, the
local/global decomposition in \cite[Section~3.4]{CL}, and the
singular-value relaxation estimates of \cite[Section~4]{BCH}.
This approach separates a regularized initial field from recent
singular-value fluctuations, which are compared with independent
reference systems. Here the comparison is carried out at elapsed time
$t^b/N$, for a fixed sufficiently large $b$, and for the nearby spatial
pairs needed below. It retains the prescribed initial observations
and compares the exact expected normalizing constants, with relative
errors small enough to remain negligible after the critical rescaling.

The second moment involves all spatial pairs, including points
arbitrarily close together. Estimates at fixed positive separation do
not control these pairs. To control the closest pairs, we temporarily
discard points where the corrected final field $Y_N^0$ is too large.
This field is unsmoothed; the history barrier instead concerns regularized
values at successive scales. More precisely, we write
\[
 \widehat Z_N^B(\mathrm dz)
   =\mathbf1_{\{Y_N^0(z)\le2t-t^{1/3}\}}\,Z_N^B(\mathrm dz),
 \qquad t=\tfrac12\log N;
\]
see \eqref{sp:barrier-weights}--\eqref{sp:fb:spatial}.
The Gaussian comparison measure has only the history barrier.
The terminal cutoff bounds the final exponential weights.
Together with estimates at intermediate separations, this bound
controls the contribution of the closest pairs, as quantified in
\eqref{fixed:eq:fb-diagonal}.
The expected discarded mass tends to zero at fixed $B$:
\[
 \mathbb E\bigl[Z_N^B(K)-\widehat Z_N^B(K)\bigr]
      \longrightarrow0
      \qquad(N\to\infty);
\]
see \eqref{p2:e}.
The terminal cutoff stays in place during the second-moment comparison
\eqref{p2:c}; we then remove it using this first-moment loss estimate.

\begin{samepage}
For fixed $B$, a real spatial test $f\in C(K;\mathbb R)$ and $\delta>0$,
we choose matching finite observation lists $D_G,D_N$ and a real
bounded smooth function $\psi$ with bounded derivatives, all fixed
before $N\to\infty$, such that
\[
\begin{aligned}
 \mathbb E|Z_B(f)-\psi(D_G)|^2&<\delta,\\
 \limsup_{N\to\infty}
 \mathbb E|\widehat Z_N^B(f)-\psi(D_N)|^2&\le\delta.
\end{aligned}
\]
\end{samepage}
The expectations are on the respective Gaussian and matrix probability
spaces. The first line is the Gaussian approximation. Expanding the matrix
square, we use convergence of the observations, comparisons for
expectations of the spatial integral times a function of those
observations, and the second-moment upper bound. This gives the second
line; see \eqref{p1:finite-approximation} and \eqref{sp:fb:projection}.
The lists can include any prescribed fixed matching field observations.
Letting $\delta\downarrow0$ after the matrix limit therefore identifies
the mass together with those observations.

Further second-moment estimates identify the height information:
for real-valued tests $f\in C(K)$ and $h\in C_c^\infty(\mathbb R)$,
\[
 \mathbb E\left|
   \widehat\Xi_N^{B,0}(f\otimes h)
      -\rho_{\mathrm R}(h)\widehat Z_N^B(f)
 \right|^2\longrightarrow0;
\]
see \eqref{p2:d}.
Here $\widehat\Xi_N^{B,0}$ attaches the corrected height mark
to $\widehat Z_N^B$, as in \eqref{sp:fb:marked}, and
$\rho_{\mathrm R}(h)=\int_0^\infty h(u)u e^{-u^2/2}\,\mathrm du$
is the Rayleigh average of $h$.
Tightness of the total mass and control of large marks in
\eqref{sp:fb:mark-tightness} complete the finite-measure conclusion.
After removing the terminal cutoff and returning to the mark defined
from the original terminal field, using \eqref{sp:fb:raw-correction},
Corollary~\ref{sp:fb:identification} gives
\[
 (Z_N^B,\Xi_N^B,C_N)
   \Longrightarrow
 (Z_B,Z_B\otimes\rho_{\mathrm R},C_G).
\]
Here $C_N,C_G$ are the prescribed matching observation vectors.
The limiting spatial measure and observations use the same Gaussian
field. The history barrier still excludes mass.

\paragraph{Removing the barrier.}
We must now show that the history barrier discards negligible critical
mass. Two estimates are needed. For $0<\epsilon<1$, write
$J_\epsilon=\lfloor(1-\epsilon)t\rfloor$, and let $E_N(B,\epsilon)$
be the event that the barrier holds at every point of $K$ through
scale $J_\epsilon$, as in \eqref{p34:prefix}.
The first estimate controls the probability of failure under the
original unweighted matrix law: for all sufficiently large $B$,
\[
 \limsup_{N\to\infty}\mathbb P\bigl(E_N(B,\epsilon)^c\bigr)
       \le C(1+B)e^{-2B};
\]
see \eqref{p34:global-bound}. The constant $C$ is independent of
$\epsilon$, although the finite-$N$ comparison rate may depend on it.
Proving the bound requires a growing list of scale observations;
convergence of each fixed list would not suffice.
The second estimate concerns the critical mass lost when the barrier
is extended to the remaining scales. Writing $Z_N^{B,J_\epsilon}$
for the measure restricted only through scale $J_\epsilon$, we have
\[
 \mathbb E\bigl[Z_N^{B,J_\epsilon}(K)-Z_N^B(K)\bigr]
   \longrightarrow
   \bigl((1-\epsilon)^{-1/2}-1\bigr)\mathbb E Z_B(K);
\]
see \eqref{p34:loss-formula}.
This limit is finite at each fixed $B,\epsilon$ and tends to zero
as $\epsilon\downarrow0$ with $B$ fixed.
When $E_N(B,\epsilon)$ holds, all discarded mass comes from the
remaining scales. Splitting off the event's complement and applying
Markov's inequality to the nonnegative mass loss combines the two
estimates; see \eqref{p34:split}.

We take $N\to\infty$ first with $B,\epsilon$ fixed, then
$\epsilon\downarrow0$ with $B$ fixed, and finally $B\to\infty$.
Thus the discarded mass tends to zero in probability, using only
$\mathbb E Z_B(K)<\infty$ at each fixed barrier.
The same loss controls the marked measure: removing a point removes
exactly its spatial weight, whatever its height mark.
The field observations remain unchanged.
Together with removal of the Gaussian barrier in \eqref{p1:removal},
this gives the unrestricted spatial and marked limits jointly with
the prescribed fixed field observations
(Corollary~\ref{p34:joint-corollary}).
These observations are still regularized; the next step recovers
averages of the original field.

\paragraph{Averages of the original field.}
To recover the field averages $\mathcal Y_N(\varphi)$ in
\eqref{sp:of:original}, we first use the regularized integrals
\[
 I_{N,S}(\varphi)
   =\int \varphi(z)
       \bigl[X_{N,S}(z)-\mathbb E X_{N,S}(z)\bigr]\,\mathrm dA(z),
\]
at a fixed integer scale $S\ge1$; see \eqref{p5:centered-scale}.
Here $\varphi$ is a smooth real test function with compact support
in the bulk. We approximate each such integral by a finite sum of
field values (Proposition~\ref{field:prop:of-quadrature}).
These values are among the observations already covered by the
joint limit. We first refine the finite sum and then remove the
regularization, using estimates for spatially integrated logarithms.
In particular, for every $\delta>0$,
\[
 \lim_{S\to\infty}\limsup_{N\to\infty}
 \mathbb P\bigl\{
   |\mathcal Y_N(\varphi)-I_{N,S}(\varphi)|>\delta
 \bigr\}=0;
\]
see \eqref{field:eq:of-clock-approximation}.
Throughout both approximations the spatial and marked measures
are kept unchanged, so their dependence on the same Gaussian field
is preserved (Proposition~\ref{p5:transfer}).
In the original field average, the logarithm is integrated in space
before subtracting its actual expectation, which also permits
discrete entry laws. This gives joint convergence for every
prescribed fixed finite family of such test functions
(Corollary~\ref{p5:field-corollary}); their supports need not lie
in the measure disk $K$.

\paragraph{Nearby determinant powers.}
We return to the height marks to describe how the measure changes when
the determinant power approaches criticality. The exact identity
\eqref{cw:tilt} reads
\[
 Z_N(\lambda)(\mathrm dz)
   =e^{\lambda R_N(z)-\lambda^2/2+e_N(\lambda,z)}
      Z_N(0)(\mathrm dz),
 \qquad Z_N(0)=Z_N.
\]
The deterministic error $e_N$ tends to zero uniformly in the spatial
variable and on each fixed compact parameter interval; see
\eqref{cw:normalizer-control}. Thus the deterministic factor tends
to $e^{-\lambda^2/2}$.
Integrating the resulting exponential weight against the Rayleigh law
would give the stated profile. This passage requires an additional
estimate: weak convergence of the marked measure controls bounded
continuous tests, whereas $e^{\lambda R_N}$ is unbounded.

To control these tests for $|\lambda|\le A$, fix $b>A$.
Separate estimates, obtained before identifying the curve, show that
the total masses at parameters $b$ and $-b$ are tight
(Proposition~\ref{cw:inputs}). The exact identity then gives, for all
sufficiently large $N$,
\[
 \Xi_N\!\left(e^{A|r|}\mathbf1_{\{|r|>L\}}\right)
   \le C_b e^{-(b-A)L}
       \bigl[Z_N(b)(K)+Z_N(-b)(K)\bigr],
\]
where $C_b$ is deterministic; see \eqref{cw:tail-domination}.
The two signs control the two mark tails, giving \eqref{cw:tails}
in probability. After truncating the mark, the kernels are bounded
and vary continuously with the parameter on $[-A,A]$.
The continuous-mapping argument then gives joint convergence of the
truncated curve, uniformly on this interval in the metric for weak
convergence of measures; see \eqref{cw:truncated-joint}.
Removing the truncation after the matrix limit gives the same
convergence for the curve, by \eqref{cw:curve-approximation}.
The Rayleigh integral identifies the limit as
\[
 M_\kappa|_K\,e^{-\lambda^2/2}
       \int_0^\infty e^{\lambda r}\,\rho_{\mathrm R}(\mathrm dr)
     =F(\lambda)D_\kappa|_K;
\]
see \eqref{sp:cw:profile}.
The spatial measure, marked measure and each prescribed finite family
of original-field averages remain unchanged in these approximations,
preserving their joint dependence on the same Gaussian field.

\subsection*{Organization of the paper}

Section~\ref{sp:common} gives the common definitions and states the main
theorem. Sections~\ref{sp:fixed-section}--\ref{cw:subsection}
give the six propositions and deduce the main theorem from them;
Section~\ref{main:limit-order} explains the parameter and limit
orders. Section~\ref{chap:gaussian} constructs the common Gaussian field and
chaos.
Sections~\ref{chap:scalar}--\ref{chap:annuli} develop the scalar,
resolvent, relaxation and weighted-comparison estimates, assembled at a
fixed barrier in Section~\ref{chap:fixed}.
Sections~\ref{chap:aggregate}--\ref{chap:window} treat the barrier through
an initial fraction of the scales, barrier removal, original-field transfer
and the continuous near-critical curve. Readers can follow the proof of
the main theorem first and then use the references following the
propositions to find the required statements in these later sections. The complete
parameter schedule is in Section~\ref{main:parameter-schedule}.

\subsection*{AI Declaration}

\enlargethispage{\baselineskip}
The proof in this work was developed during a long chain of conversations with OpenAI ChatGPT (GPT-6) which also helped with the preparation
of the manuscript and numerical verification.  All results were written and checked by the authors who take full responsibility.

\begingroup

\newcommand{\E}{\mathbb E}
\newcommand{\Prob}{\mathbb P}
\newcommand{\PP}{\mathbb P}
\newcommand{\R}{\mathbb R}
\newcommand{\C}{\mathbb C}
\newcommand{\Disk}{\mathbb D}
\newcommand{\one}{\mathbf1}
\newcommand{\ind}{\mathbf1}
\newcommand{\dd}{\,\mathrm d}
\newcommand{\TV}{\mathrm{TV}}
\newcommand{\BL}{\mathrm{BL}}
\newcommand{\cM}{\mathcal M_+}
\newcommand{\cD}{\mathcal D}
\newcommand{\rhoR}{\rho_{\mathrm R}}
\newcommand{\cSH}{\sqrt{2/\pi}}
\newcommand{\csh}{\sqrt{2/\pi}}
\newcommand{\Ytest}{\mathcal Y}
\newcommand{\norm}[1]{\left\|#1\right\|}
\newcommand{\Cov}{\operatorname{Cov}}
\theoremstyle{plain}
\newtheorem{earlytheorem}{Theorem}
\newtheorem{earlyproposition}{Proposition}
\newtheorem{earlycorollary}{Corollary}
\renewcommand{\theearlycorollary}{\Alph{earlycorollary}}
\newtheorem{earlytail}{Corollary}
\renewcommand{\theearlytail}{T}
\allowdisplaybreaks[1]
\section{Critical chaos, marked heights, and the compact critical window}
\label{chap:main}

The main theorem couples the critical spatial measure, its height
marks, the original field averages and the nearby-exponent curve.
We first identify the measure under a fixed barrier, then remove
the barrier and recover the original field tests. The marked limit,
with exponential tail control, determines the whole compact interval
of nearby exponents.

Six propositions separate the estimates from their joint deduction.
The Gaussian target uses the established critical-chaos construction
with the fixed-barrier integrability and measurability proved below.
The matrix propositions compare the required weighted moments and
control the barrier through a fixed initial fraction of the scales.
They also estimate the mass lost on the remaining scales.
Approximation of the original field and bounds for the two nearby
total masses then preserve the joint limiting law when the exponent
changes. We give the complete deductions
here and a precise reference after each proposition to its later
verification.

\subsection{Common setting and the joint conclusion}
\label{sp:common}

 {We first fix the common definitions; then we state the joint convergence theorem with its normalization, coupling, and topology.}
Let \(A_N=N^{-1/2}(\xi_{ij})_{i,j\le N}\), where the entries are independent
copies of a fixed complex random variable satisfying
\begin{equation}
 \E\xi=\E\xi^2=0,\qquad \E|\xi|^2=1,\qquad
 \E|\xi|^p<\infty\quad\text{for every }p>0.
 \label{main:model}
\end{equation}
No density assumption is imposed, and the real and imaginary parts
of an entry need not be independent.

We say that an event holds with overwhelming probability if, for every
fixed \(D>0\), its complement has probability at most \(C_DN^{-D}\)
for all sufficiently large \(N\). Constants and eventual thresholds may
depend on the parameters held fixed in the relevant statement. The law
and any simultaneous range of parameters are specified locally.
We write \(\kappa=\E|\xi|^4-2\ge-1\) and
\(\Disk=\{z\in\C:|z|<1\}\). Fix a closed disk
\(K=\{|z|\le r_0\}\), \(0<r_0<1\), and an integer \(k_0\ge1\). Put
\[
 t=\tfrac12\log N,\qquad n=\lfloor t\rfloor,
 \qquad \dd A=\dd x\dd y.
\]
The deterministic Hermitization density \(\rho^z\) describes the limiting
distribution of the symmetrized singular values of the shifted matrix.
It is symmetric and has total mass one. Its Stieltjes transform is the branch
\[
 m^z(w)=\int_{\R}\frac{\rho^z(x)}{x-w}\dd x,
 \qquad -\frac1{m^z(w)}=w+m^z(w)-\frac{|z|^2}{w+m^z(w)},
\]
with \(\operatorname{Im}m^z(w)>0\) in the upper half-plane and
\(m^z(w)\sim-1/w\) at infinity
(see the conventions in Section~\ref{chap:scalar}). The following formulas
regularize the logarithm near determinant zeros and subtract its
deterministic leading term. The scale parameter increases as the
regularization height decreases. We also include the fourth-cumulant
correction used in the Gaussian comparison. We set
\begin{align}
 L_{N,\eta}(z)
 &=\frac12\log\det\bigl((A_N-zI)(A_N-zI)^*+\eta^2I\bigr),
       \qquad \eta>0,\nonumber\\
 L_{N,0}(z)&=\log|\det(A_N-zI)|,\label{sp:logdet}\\
 \Phi_N(z,\eta)
 &=L_{N,\eta}(z)
   -\frac N2\int_{\R}\log(x^2+\eta^2)\rho^z(x)\dd x ,
       \label{sp:fb:Phi}\\
 u_s&=\tanh s,\qquad p_z(u)=u(1-u|z|^2),\qquad
 \eta_s(z)=(1-u_s)\sqrt{u_s^{-1}-|z|^2},\nonumber\\
 X_{N,s}(z)
 &=\sqrt2\left[\Phi_N(z,\eta_s(z))
                    +\frac\kappa4p_z(u_s)^2\right].
       \label{sp:fb:scale}
\end{align}
The raw and corrected terminal fields are
\begin{align}
 Y_N(z)&=\sqrt2\left[\log|\det(A_N-zI)|
                         +\frac N2(1-|z|^2)\right],\nonumber\\
 Y_N^0(z)&=Y_N(z)+c_\kappa(z),\qquad
 c_\kappa(z)=\frac{\sqrt2\kappa}{4}(1-|z|^2)^2 .
       \label{sp:fb:endpoints}
\end{align}
We define the critical density and the two endpoint marks by
\begin{equation}
 W_N(z)=\frac{e^{2Y_N(z)}}{\E e^{2Y_N(z)}}
       =\frac{e^{2Y_N^0(z)}}{\E e^{2Y_N^0(z)}},\qquad
 R_N=\frac{2t-Y_N}{\sqrt t},\qquad
 R_N^0=\frac{2t-Y_N^0}{\sqrt t}.
       \label{sp:fb:roots}
\end{equation}
The equality of the two normalized densities follows by cancellation
of the deterministic correction. At a determinant zero, \(W_N=0\)
and both marks are set to zero.

Before imposing the barrier, we define the critical measure and
its marked version with marks from the original terminal field by
\begin{equation}
 Z_N(\dd z)=\sqrt t\,W_N(z)\dd A(z),\qquad
 \Xi_N(\dd z,\dd u)=Z_N(\dd z)\delta_{R_N(z)}(\dd u).
 \label{sp:unbarred}
\end{equation}
All these measures are restricted spatially to \(K\).

\paragraph{One Gaussian space.}
We fix a centered Gaussian family indexed by scale and spatial point,
\(\{G_s(z):s>0,\ z\in\Disk\}\), with its continuous
positive-scale versions and covariance
\begin{align}
 C_\kappa(s,z;r,w)
   &=-\tfrac12\log D(u_s,z;u_r,w)
                       +\tfrac\kappa2p_z(u_s)p_w(u_r),\nonumber\\
 D(u,z;v,w)
   &=1-uv-uv\bigl[(1-u)|z|^2+(1-v)|w|^2\bigr]
                       +uv|z-w|^2 .
       \label{sp:fb:covariance}
\end{align}
Its generalized spatial limit is the field \(G_\kappa\) with covariance
\begin{equation}
 K_\kappa(z,w)=-\log|z-w|
                +\frac\kappa2(1-|z|^2)(1-|w|^2).
 \label{sp:kernel}
\end{equation}
On this same Gaussian probability space, let \(D_\kappa\) be its
derivative-normalized critical chaos and put
\[
 M_\kappa=\sqrt{2/\pi}\,D_\kappa,\qquad
 \rhoR(\dd u)=u e^{-u^2/2}\one_{\{u>0\}}\dd u.
\]
The positive-semidefinite construction, including \(\kappa=-1\),
and the common realization are those of
Section~\ref{chap:gaussian}, in particular
Section~\ref{gaussian:gp2:canonical}. The derivative normalization uses
the approximants
\[
 (2\operatorname{Var}\widehat G_s-\widehat G_s)
 e^{2\widehat G_s-2\operatorname{Var}\widehat G_s}\dd A,
 \qquad \widehat G_s=G_\kappa*\rho_{e^{-s}},
\]
where \(\rho\) is a fixed smooth nonnegative mollifier of unit mass.
These approximants converge to \(D_\kappa\) in probability against
compactly supported continuous spatial tests.

Every auxiliary disk carries a restriction of this same family.
For each prescribed finite observation or test family, choose a closed
bulk disk \(K_*=\{|z|\le r_*\}\) whose interior contains \(K\),
the observation sites, and all test supports. Choose also a fixed
closed bulk disk \(K_+\) whose interior contains \(K\), for the
maximum estimate. These disks are fixed before the matrix limit
and independently of the barrier and terminal fraction.

\paragraph{Field tests and measure topologies.}
Since determinant zeros cause singularities, we use spatial averages to study fluctuations.
For a real \(\varphi\in C_c^\infty(\Disk)\), 
define
\begin{align}
 \Ytest_N(\varphi)
 &=\int\varphi(z)Y_N(z)\dd A(z)
                  -\E\int\varphi(z)Y_N(z)\dd A(z)\nonumber\\
 &=\sqrt2\left[
      \int\varphi(z)\log|\det(A_N-zI)|\dd A(z)
       -\E\int\varphi(z)\log|\det(A_N-zI)|\dd A(z)
         \right].
       \label{sp:of:original}
\end{align}
The spatial integral precedes centering. Spatial logarithmic
integrability, established in Lemma~\ref{field:lem:of-integrability},
makes this
definition well defined under the stated moment assumptions,
including for discrete entry laws. No pointwise expectation of
the zero-height logarithm is used. The deterministic bulk
potential cancels in the second line of \eqref{sp:of:original}.

We write \(\mathcal M_+(S)\) for finite positive Borel measures on a
Polish space \(S\), with weak convergence against bounded continuous
functions, including the constant function. Thus total mass is part
of the topology. For a metric space \(S\), put
\begin{equation}
 d_{\BL,S}(\mu,\nu)
   =\sup_{\|f\|_\infty\le1,\,\operatorname{Lip}(f)\le1}
                            |\mu(f)-\nu(f)|.
 \label{sp:BL}
\end{equation}
We use Euclidean distance on \(K\) and the metric
\(|z-w|+|u-v|\) on \(K\times\R\). The total-variation norm is
\(\|\sigma\|_{\TV}=|\sigma|(S)\).

For each fixed compact parameter interval and all sufficiently
large \(N\), define
\begin{align}
 \gamma_N(\lambda)&=2-\lambda/\sqrt t,\qquad
 \psi_N(\gamma,z)=\log\E e^{\gamma Y_N(z)},\nonumber\\
 Z_N(\lambda)(\dd z)
   &=\sqrt t\,e^{\gamma_N(\lambda)Y_N(z)
                      -\psi_N(\gamma_N(\lambda),z)}\dd A(z).
 \label{sp:window-measure}
\end{align}
More precisely, on \(|\lambda|\le b\), take \(N\) large enough
that \(t\ge4b^2\). Then
\(\gamma_N(\lambda)\in\Gamma=[3/2,5/2]\).
The scalar estimate used in Proposition~\ref{cw:inputs} gives
finite, strictly positive normalizing constants uniformly on this
domain. All matrix-measure definitions are used for sufficiently
large \(N\). At determinant zeros the positive-exponent densities
vanish. In particular, \(Z_N(0)=Z_N\).
Set
\[
 \mathcal C_A=C([-A,A],\mathcal M_+(K)),\qquad
 d_A(\zeta,\zeta')=
   \sup_{|\lambda|\le A}d_{\BL,K}(\zeta_\lambda,\zeta'_\lambda).
\]
The finite-measure space on compact \(K\) is complete and separable
for \(d_{\BL,K}\), so \(\mathcal C_A\) is Polish.

\begin{earlytheorem}[Joint critical limit]
\label{main:IMW}
For every fixed \(0\le A<\infty\) and every fixed finite family
\(\varphi_1,\ldots,\varphi_m\in C_c^\infty(\Disk;\R)\), put
\[
 \boldsymbol{\mathcal Y}_N=(\Ytest_N(\varphi_j))_{j=1}^m,
 \qquad \boldsymbol G=(G_\kappa(\varphi_j))_{j=1}^m,
 \qquad \mathcal Z_N=(Z_N(\lambda))_{|\lambda|\le A}.
\]
For all sufficiently large \(N\), \(\mathcal Z_N\) is a random
element of \(\mathcal C_A\), and
\begin{align}
 &\bigl(Z_N,\Xi_N,\boldsymbol{\mathcal Y}_N,\mathcal Z_N\bigr)
       \Longrightarrow\nonumber\\
 &\qquad\bigl(M_\kappa|_K,
        M_\kappa|_K\otimes\rhoR,\boldsymbol G,
        (F(\lambda)D_\kappa|_K)_{|\lambda|\le A}\bigr)
 \label{sp:main-limit}
\end{align}
in
\(\mathcal M_+(K)\times\mathcal M_+(K\times\R)
  \times\R^m\times\mathcal C_A\), where
\begin{align}
 F(\lambda)
 &=\sqrt{2/\pi}\,e^{-\lambda^2/2}
           \int_0^\infty u e^{-u^2/2+\lambda u}\dd u\nonumber\\
 &=\sqrt{2/\pi}\,e^{-\lambda^2/2}
                   +2\lambda\Phi_{\mathrm{nor}}(\lambda).
 \label{sp:cw:profile}
\end{align}
Here \(\Phi_{\mathrm{nor}}\) is the standard normal distribution
function. Every limiting coordinate uses the Gaussian family fixed
above. The field-test supports need not lie in \(K\).
\end{earlytheorem}

\paragraph{A joint approximation principle.}
We use the following elementary fact in product spaces. Let \(T\)
be Polish with compatible metric \(d\). Suppose, for each fixed
approximation index \(a\), that \(U_N^a\Rightarrow U^a\), with
\(U_N,U_N^a\) defined together and \(U^a,U\) on one limiting
probability space. If \(U^a\to U\) in probability and
\begin{equation}
 \lim_{a\to\infty}\limsup_N
      \Prob\{d(U_N,U_N^a)>\eta\}=0\qquad(\eta>0),
 \label{sp:approximation-hypothesis}
\end{equation}
then \(U_N\Rightarrow U\). Indeed, for any real \(H\) bounded by
one with Lipschitz constant at most one,
\begin{align}
 |\E H(U_N)-\E H(U)|
 \le{}&\E[2\wedge d(U_N,U_N^a)]
       +|\E H(U_N^a)-\E H(U^a)|\nonumber\\
      &+\E[2\wedge d(U^a,U)].
 \label{sp:approximation-proof}
\end{align}
The first term is at most
\(\eta+2\Prob\{d(U_N,U_N^a)>\eta\}\), the middle term vanishes
at fixed \(a\), and the last term tends to zero as \(a\to\infty\).
Taking these limits and then \(\eta\downarrow0\) proves the claim.
An approximation parameter tending to zero is handled by
reparametrization. For product spaces we use sums of compatible
component metrics, truncating them when convenient. Unchanged
coordinates contribute zero to the approximation distance.

\subsection{Fixed-barrier identification}
\label{sp:fixed-section}

 Next, we first identify the critical measure under a fixed integer-scale
history barrier, jointly with its height marks and each prescribed
finite collection of coarse field observations.
The Gaussian target and the matrix estimates
are treated separately before the joint deduction.

For integers \(j\ge k_0\), define
\begin{equation}
 A_{N,j}^B(z)=\bigcap_{k=k_0}^j\{X_{N,k}(z)\le2k+B\},
 \qquad
 A_j^B(z)=\bigcap_{k=k_0}^j\{G_k(z)\le2k+B\}.
 \label{sp:barriers}
\end{equation}
We write \(d=t^{1/3}\), \(Z_N^{B,j}=\one_{A_{N,j}^B}Z_N\), and
\begin{align}
 W_N^B(z)&=W_N(z)\one_{A_{N,n}^B(z)},&
 \widehat W_N^B(z)&=W_N^B(z)\one_{\{Y_N^0(z)\le2t-d\}},
 \label{sp:barrier-weights}\\
 Z_N^B(\dd z)&=\sqrt t\,W_N^B(z)\dd A(z),&
 \widehat Z_N^B(\dd z)&=\sqrt t\,\widehat W_N^B(z)\dd A(z).
 \label{sp:fb:spatial}
\end{align}
Thus \(Z_N^B=Z_N^{B,n}\). The corrected and raw marked measures are
\begin{align}
 \Xi_N^{B,0}(\dd z,\dd u)
   &=Z_N^B(\dd z)\delta_{R_N^0(z)}(\dd u),&
 \widehat\Xi_N^{B,0}(\dd z,\dd u)
   &=\widehat Z_N^B(\dd z)\delta_{R_N^0(z)}(\dd u),\nonumber\\
 \Xi_N^B(\dd z,\dd u)&=Z_N^B(\dd z)\delta_{R_N(z)}(\dd u).
 \label{sp:fb:marked}
\end{align}
The deterministic correction cancels from \(W_N\), but remains
in the endpoint cutoff and the corrected mark.

For real \(s\ge k_0+1\), set
\begin{align}
 v_s(z)&=C_\kappa(s,z;s,z),&
 W_s^G(z)&=e^{2G_s(z)-2v_s(z)},\nonumber\\
 W_s^{G,B}(z)&=W_s^G(z)\one_{A_{\lfloor s\rfloor}^B(z)},&
 R_s^G(z)&=\frac{2s-G_s(z)}{\sqrt s},\nonumber\\
 \mu_s^G(\dd z)&=\sqrt s\,W_s^G(z)\dd A(z),&
 \mu_s^{G,B}(\dd z)&=\sqrt s\,W_s^{G,B}(z)\dd A(z).
 \label{sp:gaussian-measures}
\end{align}
The terminal Gaussian weight and mark use the real scale parameter \(s\),
whereas the barrier is observed only at the integers through
\(\lfloor s\rfloor\). The Gaussian measure has no endpoint cutoff.
We also use the probability measures
\begin{equation}
 \frac{\dd Q_{N,z}^{\mathrm{det}}}{\dd\Prob}=W_N(z),\qquad
 \frac{\dd Q_{s,z}^{G}}{\dd\Prob_G}=W_s^G(z),
 \label{sp:weighted-laws}
\end{equation}
where \(\Prob_G\) is the law of the specified Gaussian family.
We call these one-root laws, with root at $z$. Weighting at two spatial
points and normalizing the product gives a two-root law.

\paragraph{Coarse observations.}
Choose a countable dense set \(\mathcal D\subset K\), and define
the completed sigma-field
\begin{equation}
 \mathcal H=
 \overline{\sigma\{G_k(v):k\ge k_0\text{ integer},\ v\in\mathcal D\}}.
 \label{sp:integer-sigma-field}
\end{equation}
We match matrix and Gaussian observations by using the same points
and scales. We call \(C_N,C_G\) an admissible matching pair if it
consists of a fixed finite sublist of the integer-scale observations
on \(\mathcal D\), adjoined to any prescribed fixed finite list
at sites in \(K_*\) with regularization parameters in \([a_0,1]\)
for some fixed \(a_0>0\). Such a pair can be written as
\begin{equation}
 C_N=(X_{N,s_i}(v_i))_{i=1}^m,\qquad
 C_G=(G_{s_i}(v_i))_{i=1}^m,\qquad
 \eta_{s_i}(v_i)\in[\eta_*,1],\quad s_i>0,\quad v_i\in K_*.
 \label{sp:coarse-vectors}
\end{equation}
Here \(\eta_*>0\) may depend on the chosen finite list. Indeed,
a finite list of integer scales has a positive minimum regularization
parameter, and \(k\ge1\) gives \(\eta_k(v)\le\eta_1(0)<1\).
Repeated observations and height one are allowed. Finite unions
of admissible lists remain admissible. Every list is fixed before
\(N\to\infty\), and all its Gaussian coordinates use the same
family, including those at sites outside \(K\).

\subsubsection{The Gaussian measure at a fixed barrier}
\label{p1:section}

 {Now we construct the Gaussian barrier restriction and approximate its spatial integrals using finitely many field observations.}
The fixed-barrier measure must be constructed from the same Gaussian
family as the coarse observations. We first identify it as a random
spatial restriction of the canonical critical measure. A bound on
three-point expectations then gives the integrability needed for
the matrix identification argument.

In this Gaussian subsection, \(t\ge k_0+1\) may be any real number
and \(n=\lfloor t\rfloor\). Uniformly on \(K\), the exact variance is
\begin{equation}
 v_t(z)=t-\tfrac12\log\bigl(4(1-|z|^2)\bigr)
       +\tfrac\kappa2(1-|z|^2)^2+O(e^{-2t}).
 \label{p1:variance}
\end{equation}
The bounded spatial correction fixes the normalization of the
limiting measure. For integers \(L\ge k_0\), write
\[
 E_L^B=\{z\in K:A_L^B(z)\text{ holds}\},\qquad
 E_\infty^B=\bigcap_{L\ge k_0}E_L^B.
\]
\begin{earlyproposition}[Gaussian fixed-barrier target]
\label{p1:proposition}
For each fixed \(B\in\R\), the measure
\begin{equation}
 Z_B=\one_{E_\infty^B}M_\kappa|_K
 \label{p1:target}
\end{equation}
is finite and atomless. For every \(f\in C(K;\R)\),
\begin{equation}
 \mu_t^{G,B}(f)\longrightarrow Z_B(f)
       \quad\text{in }L^2,\qquad t\to\infty
       \text{ through real values}.
 \label{p1:L2}
\end{equation}
The random measure \(Z_B\) is \(\mathcal H\)-measurable. Moreover,
\begin{equation}
 \norm{Z_B-M_\kappa|_K}_{\TV}\longrightarrow0
       \quad\text{in probability as }B\to\infty.
 \label{p1:removal}
\end{equation}
All Gaussian objects in these assertions use the specified common
probability space.
\end{earlyproposition}

The proof also gives, for each fixed \(B\),
\begin{equation}
 \sup_{t\ge k_0+1}\E[\mu_t^{G,B}(K)^3]<\infty,
 \qquad \E[Z_B(K)^3]<\infty,
 \label{p1:third-moments}
\end{equation}
and \(L^p\) convergence in~\eqref{p1:L2} for every \(1\le p<3\).
The constants in these moment bounds may depend on \(B,K,k_0,\kappa\).
The verification is given by the canonical construction in
Section~\ref{gaussian:gp2:canonical}, the restriction argument in
Proposition~\ref{gaussian:gp2:probability}, and
Theorem~\ref{gaussian:gp2:gaussian-main}. The latter also gives
\(L^2\) convergence when the two real terminal scale parameters differ.
Removal of the barrier on the same Gaussian probability space
follows from Proposition~\ref{gaussian:gp2:max-prop}
and the paragraph following Theorem~\ref{gaussian:gp2:gaussian-main}.

\paragraph{Finite-coordinate approximation.}
Measurability with respect to the integer-scale observations is proved
in Section~\ref{gaussian:gp2:coarse}. Its integer specialization of
Lemma~\ref{gaussian:gp2:cylinder} gives the following consequence.
For each fixed \(B\), \(f\in C(K;\R)\), and \(\delta>0\), there are
a fixed finite matching list \(D_N,D_G\) of integer observations
and a real bounded smooth function \(\psi\), with all derivatives
bounded, such that
\begin{equation}
 \E|Z_B(f)-\psi(D_G)|^2<\delta.
 \label{p1:finite-approximation}
\end{equation}
The list may be enlarged by any prescribed admissible coarse
observations, because \(\psi\) may ignore the added coordinates.
No equality between the integer-history sigma-field and the
sigma-field of the full positive-time family is needed.

\subsubsection{Matrix estimates at a fixed barrier}
\label{p2:section}

 {Next we shall state weighted moment comparisons and an estimate for the expected mass discarded by the terminal cutoff, preparing the matrix identification.}
 
We return to \(t=\tfrac12\log N\), \(n=\lfloor t\rfloor\), and
\(d=t^{1/3}\). The coarse observations keep track of the relation
between the critical measure and the Gaussian field. The endpoint
cutoff makes the contribution of very close spatial pairs small.
We prove that its discarded mass tends to zero, compare first and
second moments, and then identify the centered height marks.
The verification uses the order (a), (e), (b), (c), (d), so the
cutoff estimate is available for the mixed first moment.

Proposition~\ref{p1:proposition} supplies the Gaussian
\(L^2\) limit on the same probability space as every prescribed
coarse vector. We use the weights and probability laws already
defined in \eqref{sp:barrier-weights}--\eqref{sp:weighted-laws}.
\begin{earlyproposition}[Matrix estimates at a fixed barrier]
\label{p2:matrix}
Fix \(B\in\R\), \(f\in C(K;\R)\), and any admissible matching
coarse vectors \(C_N,C_G\in\R^m\). For every
\(\Psi\in C_b^\infty(\R^m;\R)\), with all derivatives bounded,
the following assertions hold:
\begin{enumerate}[label=\textup{(\alph*)},leftmargin=2.1em,itemsep=5pt]
\item The coarse vectors satisfy
\begin{equation}
 C_N\Longrightarrow C_G.\label{p2:a}
\end{equation}
\item The mixed first moments satisfy
\begin{equation}
 \E[\widehat Z_N^B(f)\Psi(C_N)]
       \longrightarrow\E[Z_B(f)\Psi(C_G)].\label{p2:b}
\end{equation}
\item The second moments with the terminal cutoff satisfy
\begin{equation}
 \limsup_{N\to\infty}\E[\widehat Z_N^B(f)^2]
       \le\E[Z_B(f)^2].\label{p2:c}
\end{equation}
\item For every \(h\in C_c^\infty(\R;\R)\),
\begin{equation}
 \E\left|\widehat\Xi_N^{B,0}(f\otimes h)
      -\rhoR(h)\widehat Z_N^B(f)\right|^2
       \longrightarrow0.\label{p2:d}
\end{equation}
\item The expected loss from the endpoint cutoff satisfies
\begin{equation}
 \E\norm{Z_N^B-\widehat Z_N^B}_{\TV}
       \longrightarrow0.\label{p2:e}
\end{equation}
\end{enumerate}
\end{earlyproposition}

The entry law, disks, barrier, test functions, and observation list
are fixed throughout. Each positive separation and each auxiliary
precision parameter is fixed before \(N\to\infty\). Constants
may depend on these fixed choices. Auxiliary comparison parameters
are chosen after the fixed tests and required precision; their
choices may differ between applications. All spatial integrals below
are over \(K\) or \(K^2\).

Clauses (a) and (b) follow from the coarse and first-moment
estimates in Section~\ref{fixed:sec:fb-first}, with the cutoff
estimate \eqref{oneroot:eq:cutoff-tv} in
Section~\ref{chap:oneroot} giving (e).
The pair and diagonal estimates in
Sections~\ref{fixed:sec:fb-cut-pair}--\ref{fixed:sec:fb-quadratic}
give (c) and (d). The full second-moment matching
\eqref{fixed:eq:fb-M3}, together with Proposition~\ref{p1:proposition},
also gives convergence of the second moments in (c).

\subsubsection{Joint identification at a fixed barrier}

 {Now we combine the moment comparisons to identify the restricted matrix measure, its marks, and the retained field observations jointly.}

\begin{earlycorollary}
\label{sp:fb:identification}
For every fixed \(B\) and every fixed admissible coarse vector,
\begin{equation}
 (\widehat Z_N^B,\widehat\Xi_N^{B,0},C_N)
 \Longrightarrow
 (Z_B,Z_B\otimes\rhoR,C_G)
       \label{sp:fb:cut-conclusion}
\end{equation}
in
\(\mathcal M_+(K)\times\mathcal M_+(K\times\R)\times\R^m\).
The same conclusion holds for
\((Z_N^B,\Xi_N^{B,0},C_N)\) and for
\((Z_N^B,\Xi_N^B,C_N)\).
The limiting measure and coarse observations are coupled on the
Gaussian probability space specified above.
\end{earlycorollary}

\begin{proof}
\emph{Gaussian approximation and projection.}
Fix a real \(f\in C(K)\), and write
\[
 T_N=\widehat Z_N^B(f),\qquad T=Z_B(f).
\]
By Proposition~\ref{p1:proposition}, \(T\in L^2(\mathcal H)\).
The finite-coordinate approximation \eqref{p1:finite-approximation}
gives, for every \(\varepsilon>0\), a fixed finite matching pair
\(D_N,D_G\) and a real bounded smooth \(\psi\), with bounded
derivatives, such that
\begin{equation}
 \E|T-\psi(D_G)|^2<\varepsilon.
 \label{sp:fb:cylinder}
\end{equation}
Enlarge the list to include the prescribed vector \(C_N,C_G\).
The function \(\psi\) can ignore the added coordinates, so
\eqref{sp:fb:cylinder} is preserved. The entire list and all
derivative bounds remain fixed before \(N\to\infty\).

Expanding the matrix square and applying
\eqref{p2:a}--\eqref{p2:c} to the enlarged list gives
\begin{align}
 \limsup_{N\to\infty}\E|T_N-\psi(D_N)|^2
 &\le
 \E T^2-2\E[T\psi(D_G)]+\E[\psi(D_G)^2]\nonumber\\
 &=\E|T-\psi(D_G)|^2<\varepsilon .
       \label{sp:fb:projection}
\end{align}
Here convergence of the last term follows from coarse convergence
and boundedness of \(\psi^2\). Only the upper bound
\eqref{p2:c} is used for the first term.

For finitely many real spatial tests \(f_1,\ldots,f_r\), choose one
union of the approximating lists and the prescribed coarse vector.
The same calculation gives a total squared approximation error at
most \(r\varepsilon\) in the limit superior. Coarse convergence
and continuity of the corresponding \(\psi_i\)'s identify the joint
law of their approximants and \(C_N\). For a bounded function of
Lipschitz constant at most one on \(\R^{r+m}\), Cauchy--Schwarz
bounds each of the matrix and Gaussian approximation errors by
\(\sqrt{r\varepsilon}\). Letting \(\varepsilon\downarrow0\) proves
\begin{equation}
 \bigl((\widehat Z_N^B(f_i))_{i=1}^r,C_N\bigr)
 \Longrightarrow
 \bigl((Z_B(f_i))_{i=1}^r,C_G\bigr).
       \label{sp:fb:finite-spatial}
\end{equation}

\emph{Spatial measure convergence.}
Taking \(f=1\) in \eqref{sp:fb:finite-spatial} proves tightness of
\(\widehat Z_N^B(K)\). On compact \(K\), sets of positive measures
with bounded total mass are weakly compact. Thus the laws of
\(\widehat Z_N^B\) are tight, jointly with the tight coarse vectors.
A countable uniformly dense subset of \(C(K)\), together with the
constant function, determines a finite measure. Equation
\eqref{sp:fb:finite-spatial} therefore identifies every
subsequential joint law and gives
\begin{equation}
 (\widehat Z_N^B,C_N)\Longrightarrow(Z_B,C_G).
       \label{sp:fb:spatial-conclusion}
\end{equation}

\emph{Marks and tightness in the mark coordinate.}
For each finite collection of tensors \(f_i\otimes h_i\), where
\(f_i\in C(K;\R)\) and \(h_i\in C_c^\infty(\R;\R)\),
\eqref{p2:d} and \eqref{sp:fb:spatial-conclusion} give
\begin{equation}
 \bigl(\widehat Z_N^B,
        (\widehat\Xi_N^{B,0}(f_i\otimes h_i))_i,C_N\bigr)
 \Longrightarrow
 \bigl(Z_B,(\rhoR(h_i)Z_B(f_i))_i,C_G\bigr).
       \label{sp:fb:finite-marked}
\end{equation}
This application of Slutsky's theorem retains the same coarse
coordinates.

Choose \(0\le h_L\le1\) in \(C_c^\infty(\R)\), increasing
pointwise to one, with \(h_L=1\) on \([-L,L]\) and \(h_L=0\)
outside \([-2L,2L]\). Such a family is obtained by scaling a smooth
even cutoff that is nonincreasing in the absolute value of its
argument. Since the spatial marginal of
\(\widehat\Xi_N^{B,0}\) is exactly \(\widehat Z_N^B\),
\eqref{p2:d} with \(f=1\) yields, for each fixed \(L\),
\begin{equation}
 \widehat\Xi_N^{B,0}(1\otimes(1-h_L))
   =(1-\rhoR(h_L))\widehat Z_N^B(K)+o_{\Prob}(1).
       \label{sp:fb:mark-cutoff}
\end{equation}
The total masses are tight and \(1-\rhoR(h_L)\to0\). Positivity
and the support of \(h_L\) therefore imply, for every \(\delta>0\),
\begin{equation}
 \lim_{L\to\infty}\limsup_{N\to\infty}
 \Prob\!\left\{
   \widehat\Xi_N^{B,0}
      \bigl(K\times(\R\setminus[-2L,2L])\bigr)>\delta
       \right\}=0 .
       \label{sp:fb:mark-tightness}
\end{equation}
Together with tight total masses, this is the tightness criterion
for laws of finite positive measures on \(K\times\R\).
Hence the triples in \eqref{sp:fb:cut-conclusion} are tight in
the stated product space. Choose countable spatial and smooth mark
test classes whose finite tensor span is uniformly dense in continuous
compactly supported tests. The resulting countable tensor family
determines their subsequential marked measures. Equation
\eqref{sp:fb:finite-marked} identifies every such joint limit
as \((Z_B,Z_B\otimes\rhoR,C_G)\).
This proves \eqref{sp:fb:cut-conclusion}, including total mass
and tightness in the mark coordinate.

\emph{Removing the terminal cutoff.}
Removing the terminal cutoff adds positive measures, and the added
marked and spatial measures have equal total mass. Thus
\begin{equation}
 \norm{\Xi_N^{B,0}-\widehat\Xi_N^{B,0}}_{\TV}
 =\norm{Z_N^B-\widehat Z_N^B}_{\TV}.
       \label{sp:fb:cut-restoration}
\end{equation}
By \eqref{p2:e}, both sides tend to zero in probability.
Total variation dominates the bounded-Lipschitz metric, so
Slutsky's theorem transfers \eqref{sp:fb:cut-conclusion} to
\((Z_N^B,\Xi_N^{B,0},C_N)\).

\emph{Returning to the raw mark.}
Off the determinant zeros,
\[
 R_N^0(z)-R_N(z)=-\frac{c_\kappa(z)}{\sqrt t}.
\]
With the convention at zeros in \eqref{sp:fb:roots}, it follows
everywhere that
\(\lvert R_N^0-R_N\rvert\le\norm{c_\kappa}_{\infty,K}/\sqrt t\).
Use the product metric and bounded-Lipschitz distance from
\eqref{sp:BL}. Testing
against functions bounded by one with Lipschitz constant at
most one gives
\begin{equation}
 d_{\BL,K\times\R}(\Xi_N^B,\Xi_N^{B,0})
 \le \frac{\norm{c_\kappa}_{\infty,K}}{\sqrt t}\,Z_N^B(K)
 \longrightarrow0
       \quad\text{in probability}.
       \label{sp:fb:raw-correction}
\end{equation}
The last step uses the tightness of \(Z_N^B(K)\) already obtained
after cutoff restoration. A final application of Slutsky's theorem
gives the claimed limit with the mark from the original terminal field,
and with the same spatial measure and coarse coordinates.
\end{proof}

The \(L^2\) assertion in Proposition~\ref{p1:proposition} is on
the Gaussian probability space. Corollary~\ref{sp:fb:identification}
is joint convergence in distribution for the matrix observables.
The cutoff restoration and mark displacement are approximation
estimates on each original matrix space. In particular, the
second-moment argument retains the terminal cutoff; its removal uses
only the first-moment estimate \eqref{p2:e}.
\subsection{Removal of the barrier}
\label{p34:subsection}

{We next show that the discarded mass vanishes in probability as the barrier rises, giving the unrestricted joint limit.}

The fixed-barrier argument identifies the part of the critical
measure whose integer-scale history stays below a fixed level.
We now show that the discarded mass becomes negligible as that
level increases. We first control all spatial points through a
fixed fraction of the scales. The expected mass lost in the
remaining fraction can then be computed from the Gaussian limit.

We use the fixed-barrier joint conclusion of
Corollary~\ref{sp:fb:identification} and the Gaussian removal in
\eqref{p1:removal}. The auxiliary disks and the prescribed coarse
list remain fixed as in the common setting.

\paragraph{The barrier through an initial fraction of the scales.}
For \(0<\epsilon<1\), define
\begin{equation}
 J_\epsilon=\lfloor(1-\epsilon)t\rfloor,\qquad
 E_N(B,\epsilon)=
 \{A_{N,J_\epsilon}^B(z)\text{ holds for every }z\in K\}.
 \label{p34:prefix}
\end{equation}
For each fixed \(\epsilon\), one has \(k_0\le J_\epsilon\le n\)
for all sufficiently large \(N\).

\begin{earlyproposition}[Uniform barrier bound through an initial fraction of the scales]
\label{p34:global}
There are constants \(C,B_0\), independent of \(\epsilon\), such
that for every fixed \(\epsilon\in(0,1)\) and every \(B\ge B_0\),
\begin{equation}
 \limsup_{N\to\infty}\Prob(E_N(B,\epsilon)^c)
       \le p(B),\qquad p(B)=C(1+B)e^{-2B}.
 \label{p34:global-bound}
\end{equation}
The constants may depend on the fixed ensemble, bulk disks, and
\(k_0\). The finite-\(N\) comparison rate may depend on \(\epsilon\).
\end{earlyproposition}

This is Corollary~\ref{aggregate:cor:au-maximum}. Its proof in
Section~\ref{aggregate:sec:maximum} combines the Gaussian comparison for
a growing list of observations with the deterministic spatial mesh and
the Gaussian maximum estimate.

\begin{earlyproposition}[Expected mass lost on the remaining scales]
\label{p34:loss}
For every fixed \(B\in\R\) and every fixed \(\epsilon\in(0,1)\),
\begin{equation}
 \E[Z_N^{B,J_\epsilon}(K)-Z_N^B(K)]
 \longrightarrow
 \big((1-\epsilon)^{-1/2}-1\big)\E Z_B(K).
 \label{p34:loss-formula}
\end{equation}
The difference is nonnegative for all sufficiently large \(N\).
\end{earlyproposition}

The calculation is \eqref{barrier:eq:br-loss} in
Section~\ref{chap:barrier}. It uses the comparison for one exponentially weighted point with a
growing list of history observations. Lemma~\ref{barrier:lem:br-root}
then changes the Gaussian scale at which this weight is evaluated. The full-barrier expectation follows from
Proposition~\ref{p2:matrix}(b),(e), before any removal of the barrier.

\paragraph{Removal in the full joint law.}

\begin{earlycorollary}
\label{p34:joint-corollary}
For every \(\delta>0\),
\begin{equation}
 \lim_{B\to\infty}\limsup_{N\to\infty}
 \Prob\{\norm{Z_N-Z_N^B}_{\TV}>\delta\}=0,
 \label{p34:spatial-removal}
\end{equation}
and the same assertion holds for \(\Xi_N,\Xi_N^B\).
For every fixed admissible coarse vector, jointly,
\begin{equation}
 (Z_N,\Xi_N,C_N)
 \Longrightarrow
 (M_\kappa|_K,M_\kappa|_K\otimes\rhoR,C_G)
 \label{p34:joint-limit}
\end{equation}
in
\(\mathcal M_+(K)\times\mathcal M_+(K\times\R)\times\R^m\).
The limiting coordinates have the coupling of the specified
exact Gaussian family.
\end{earlycorollary}

\begin{proof}
\emph{The discarded mass.}
Let \(\Delta_N^B=\norm{Z_N-Z_N^B}_{\TV}\).
Fix \(B\ge B_0\) and \(\epsilon\in(0,1)\).
On \(E_N(B,\epsilon)\), \(Z_N=Z_N^{B,J_\epsilon}\), so
\(\Delta_N^B=Z_N^{B,J_\epsilon}(K)-Z_N^B(K)\).
Splitting off the complement of the event and applying Markov's
inequality to the nonnegative loss gives
\begin{align}
 \Prob(\Delta_N^B>\delta)
 &\le\Prob(E_N(B,\epsilon)^c)
       +\delta^{-1}\E[Z_N^{B,J_\epsilon}(K)-Z_N^B(K)].
 \label{p34:split}
\end{align}
Propositions~\ref{p34:global} and~\ref{p34:loss} imply
\begin{equation}
 \limsup_N\Prob(\Delta_N^B>\delta)
 \le p(B)+\frac{(1-\epsilon)^{-1/2}-1}{\delta}\,\E Z_B(K).
 \label{p34:ordered}
\end{equation}
For this fixed \(B\), \(\E Z_B(K)<\infty\). The left side
is independent of \(\epsilon\), so we may send
\(\epsilon\downarrow0\) to obtain
\[
 \limsup_N\Prob(\Delta_N^B>\delta)\le p(B).
\]
Now \(B\to\infty\) proves~\eqref{p34:spatial-removal}.
The order is
\begin{equation}
 N\to\infty\text{ with }B,\epsilon\text{ fixed};\qquad
 \epsilon\downarrow0\text{ with }B\text{ fixed};\qquad
 B\to\infty.
 \label{p34:order}
\end{equation}
Only the fixed-\(B\) moment is used, and no conditioning on the
global event enters~\eqref{p34:split}.

\emph{The marked mass.}
The measures \(\Xi_N\) and \(\Xi_N^B\) use the same raw mark.
Their positive difference has total mass equal to the discarded
spatial mass. Therefore
\begin{equation}
 \norm{\Xi_N-\Xi_N^B}_{\TV}=\Delta_N^B.
 \label{p34:marked-removal}
\end{equation}
This proves marked barrier removal as well.

\emph{Joint convergence.}
For the fixed coarse vector, put
\begin{align*}
 U_N&=(Z_N,\Xi_N,C_N),&
 U_N^B&=(Z_N^B,\Xi_N^B,C_N),\\
 U&=(M_\kappa|_K,M_\kappa|_K\otimes\rhoR,C_G),&
 U^B&=(Z_B,Z_B\otimes\rhoR,C_G).
\end{align*}
Using the bounded-Lipschitz metrics from \eqref{sp:BL},
a bounded metric for the product topology is
\begin{align*}
 d((\mu,\nu,c),(\mu',\nu',c'))
 ={}&(1\wedge d_{\BL,K}(\mu,\mu'))
      +(1\wedge d_{\BL,K\times\R}(\nu,\nu'))\\
    &+(1\wedge|c-c'|).
\end{align*}
Total variation dominates each bounded-Lipschitz distance.
The coarse coordinate is unchanged, so
\begin{equation}
 d(U_N,U_N^B)\le2(1\wedge\Delta_N^B).
 \label{p34:matrix-distance}
\end{equation}
Equation~\eqref{p34:spatial-removal} therefore supplies the
matrix approximation condition \eqref{sp:approximation-hypothesis}.

On the specified Gaussian space, \(M_\kappa|_K-Z_B\) is a
positive measure, and \(\rhoR\) is a probability measure. Hence
\[
 \norm{(M_\kappa|_K-Z_B)\otimes\rhoR}_{\TV}
       =\norm{M_\kappa|_K-Z_B}_{\TV}.
\]
By \eqref{p1:removal}, \(d(U^B,U)\to0\) in probability on
this space. At each fixed \(B\), Corollary~\ref{sp:fb:identification}
gives \(U_N^B\Rightarrow U^B\). The approximation argument
\eqref{sp:approximation-proof}, with \(N\to\infty\) first and
then \(B\to\infty\), proves \eqref{p34:joint-limit}.
The prescribed coarse coordinate is unchanged in both approximations.
\end{proof}

The moment needed in this argument is
\(\E Z_B(K)<\infty\) at each fixed \(B\).
The original normalization gives
\(\E Z_N(K)=\sqrt t\,|K|\), so the final removal is used in
probability. The joint conclusion retains every prescribed fixed
coarse vector, including observations outside the measure disk.
Its extension to original integral field tests is the separate
transfer argument of Proposition~\ref{p5:transfer}.
\subsection{Transfer to the original field}
\label{p5:subsection}

 {Next we recover averages of the unsmoothed field while keeping the spatial and marked measures in the same joint limit.}
Indeed, the preceding joint convergence includes every fixed finite
collection of positive-scale field values. We now recover the
original centered spatial integrals. At a fixed positive scale,
we approximate each integral by finitely many field values.
We then refine this approximation and remove the positive scale.
The spatial and marked measures remain unchanged during these
approximations, which preserves their joint limiting law.

For the fixed enlarged disk \(K_*\), we write
\(a_S(z)=\eta_S(z)\). The original field test is the spatially
integrated, then centered variable \eqref{sp:of:original}.
For every fixed \(S\ge1\), \(a_S\) is continuous with positive
minimum on \(K_*\), and \(\sup_{K_*}a_S\le Ce^{-2S}\).
The positive-scale integral centered by its expectation is
\begin{align}
 I_{N,S}(\varphi)
 &:=\int_{K_*}\varphi(z)
          [X_{N,S}(z)-\E X_{N,S}(z)]\dd A(z)\nonumber\\
 &=\sqrt2\left[
    \int_{K_*}\varphi(z)L_{N,a_S(z)}(z)\dd A(z)
    -\E\int_{K_*}\varphi(z)L_{N,a_S(z)}(z)\dd A(z)\right].
 \label{p5:centered-scale}
\end{align}
All corrections between \(X_{N,S}\) and \(\sqrt2L_{N,a_S}\)
are deterministic, so this identity is exact. The actual
expectation in~\eqref{p5:centered-scale} is distinguished from
the deterministic centering in~\eqref{sp:fb:scale}.

\begin{earlyproposition}[Transfer from coarse observations to the original field]
\label{p5:transfer}
Fix a closed bulk disk \(K_*\Subset\Disk\) and an integer
\(k_*\ge1\). Let \(V_N\) take values in a Polish space
\(\mathcal E\), and let \(V\) be defined on the same
probability space as the exact Gaussian family.
Suppose that, for every fixed finite list of integer scale parameters
\(s_i\ge k_*\) and sites \(z_i\in K_*\),
\begin{equation}
 (V_N,(X_{N,s_i}(z_i))_i)
 \Longrightarrow (V,(G_{s_i}(z_i))_i)
 \label{p5:premise}
\end{equation}
with this same limiting coupling.
Then, for every fixed finite family
\(\varphi_1,\ldots,\varphi_m\in
C_c^\infty(\operatorname{int}K_*;\R)\),
\begin{equation}
 (V_N,(\Ytest_N(\varphi_j))_{j=1}^m)
 \Longrightarrow (V,(G_\kappa(\varphi_j))_{j=1}^m)
 \label{p5:conclusion}
\end{equation}
in \(\mathcal E\times\R^m\).
\end{earlyproposition}

This follows from the proof of Proposition~\ref{field:prop:of-joint},
applied on \(K_*\) using sufficiently large integer scale parameters.
It uses Theorem~\ref{field:thm:of-approximation}, in the form adapted to this
scale profile
\eqref{field:eq:of-clock-approximation}, for removal of the positive
cutoff, Lemma~\ref{field:lem:of-coarse-moments} for moment convergence,
and Proposition~\ref{field:prop:of-quadrature} for finite quadrature.
The retained coordinate \(V_N\) is unchanged in these approximations.

\paragraph{Application after removing the barrier.}
We apply the general transfer result to the spatial and marked measures
already identified in Corollary~\ref{p34:joint-corollary}.

\begin{earlycorollary}
\label{p5:field-corollary}
For every fixed finite family
\(\varphi_1,\ldots,\varphi_m\in C_c^\infty(\Disk;\R)\),
jointly,
\begin{equation}
 (Z_N,\Xi_N,(\Ytest_N(\varphi_j))_{j=1}^m)
 \Longrightarrow
 (M_\kappa|_K,M_\kappa|_K\otimes\rhoR,
                   (G_\kappa(\varphi_j))_{j=1}^m)
 \label{p5:application}
\end{equation}
in
\(\mathcal M_+(K)\times\mathcal M_+(K\times\R)\times\R^m\).
The test supports need not lie in \(K\). All limiting
coordinates have the coupling of the specified Gaussian field.
\end{earlycorollary}

\begin{proof}
Choose one closed bulk disk \(K_*\) whose interior contains
\(K\) and all the test supports. In the Polish space
\(\mathcal E=\mathcal M_+(K)\times\mathcal M_+(K\times\R)\),
take
\[
 V_N=(Z_N,\Xi_N),\qquad
 V=(M_\kappa|_K,M_\kappa|_K\otimes\rhoR).
\]
Corollary~\ref{p34:joint-corollary} gives \eqref{p5:premise} for every permitted
fixed finite coarse list on this enlarged disk, including
sites outside \(K\). For every integer \(s\ge1\),
\[
 0<\eta_s(z)\le
       \frac{1-\tanh1}{\sqrt{\tanh1}}<1.
\]
Each fixed finite list with \(s_i\ge k_0\) therefore has
heights in \([\eta_*,1]\) for some fixed \(\eta_*>0\), and
is a permitted coarse list. Proposition~\ref{p5:transfer}
applies with \(k_*=k_0\). It keeps the entire measure pair
unchanged and uses the same Gaussian field that generates
\(D_\kappa\). This proves~\eqref{p5:application}.
\end{proof}

The conclusion retains total mass in the measure topology and
the already restored raw mark. It is joint convergence for
each prescribed finite family of original field tests.
All regularization and quadrature estimates in the transfer
are taken under the unweighted matrix law.
\subsection{The compact critical window}
\label{cw:subsection}

 {In this subsection we use the height marks and their exponential tails to identify the curve on each fixed compact interval of the near-critical parameter.}

Changing the exponent reweights the normalized height deficits.
The exact ratio of densities therefore expresses the nearby
measures in terms of the marked measure. To pass to its limit
uniformly over a compact parameter interval, we first control
the exponential mark tails. Proposition~\ref{cw:inputs} provides the
normalizing-constant estimate and the bounds for total masses at two
nearby exponents needed for this step.

We define the deterministic error in the normalizing constants by
\begin{equation}
 e_N(\lambda,z)
   =\psi_N(2,z)-\psi_N(\gamma_N(\lambda),z)
                   -2\lambda\sqrt t+\frac{\lambda^2}{2}.
 \label{cw:normalizer-error}
\end{equation}

\begin{earlyproposition}[Normalizing constants and tightness at nearby exponents]
\label{cw:inputs}
For every fixed \(b>0\),
\begin{equation}
 \delta_N(b):=
    \sup_{|\lambda|\le b,\,z\in K}|e_N(\lambda,z)|
       \longrightarrow0.
 \label{cw:normalizer-control}
\end{equation}
Moreover,
\begin{equation}
 S_N(b):=Z_N(b)(K)+Z_N(-b)(K)
 \label{cw:endpoint-sum}
\end{equation}
is tight:
\begin{equation}
 \lim_{H\to\infty}\limsup_{N\to\infty}
          \Prob\{S_N(b)>H\}=0.
 \label{cw:endpoint-tightness}
\end{equation}
\end{earlyproposition}

\begin{proof}[Verification by the later window estimates]
Lemma~\ref{window:lem:normalizer}, in particular
\eqref{window:eq:cw-modulus}, proves the normalizer assertion.
Equation~\eqref{window:eq:cw-two-masses} proves tightness of the
sum using one event under the unweighted matrix law, on which the
barrier holds at every spatial point through the retained fraction $1/2$. Proposition~\ref{window:prop:cw-P} retains the
stronger uniform marginal statement. Its inputs are
Theorems~\ref{window:thm:cw-source} and
\ref{window:thm:cw-survival}, together with
Proposition~\ref{p34:global}. The full joint curve conclusion
is deduced below from these estimates and the marked limit.
\end{proof}

\paragraph{Exact exponent change and exponential mark tails.}
The ratio of the densities in \eqref{sp:window-measure} gives
\begin{equation}
 Z_N(\lambda)(\dd z)
   =e^{\lambda R_N(z)-\lambda^2/2+e_N(\lambda,z)}
                                  Z_N(0)(\dd z).
 \label{cw:tilt}
\end{equation}
Indeed, the logarithmic ratio is
\(-\lambda Y_N/\sqrt t+\psi_N(2,z)
-\psi_N(\gamma_N(\lambda),z)\).
Definition~\eqref{cw:normalizer-error} gives \eqref{cw:tilt}
exactly. At determinant zeros both positive-exponent densities
vanish, so the identity of measures remains valid.

Fix \(0\le A<\infty\) and one \(b>A\), independently of \(N\).
For \(L\ge0\) and \(r\in\R\),
\[
 e^{A|r|}\one_{\{|r|>L\}}
       \le e^{-(b-A)L}(e^{br}+e^{-br}).
\]
Applying \eqref{cw:tilt} at \(b\) and \(-b\) gives
\begin{equation}
 \Xi_N\!\left(e^{A|r|}\one_{\{|r|>L\}}\right)
   \le e^{b^2/2+\delta_N(b)-(b-A)L}S_N(b).
 \label{cw:tail-domination}
\end{equation}
Proposition~\ref{cw:inputs} therefore implies, for every
\(\eta>0\),
\begin{equation}
 \lim_{L\to\infty}\limsup_{N\to\infty}
 \Prob\!\left\{
   \Xi_N\!\left(e^{A|r|}\one_{\{|r|>L\}}\right)>\eta
          \right\}=0.
 \label{cw:tails}
\end{equation}
Without the indicator, the same bound shows that
\begin{equation}
 H_{N,A}:=\Xi_N(e^{A|r|})
       \le e^{b^2/2+\delta_N(b)}S_N(b)
       \quad\text{is tight}.
 \label{cw:weighted-mass}
\end{equation}
It also yields the pathwise estimate
\begin{equation}
 \sup_{|\lambda|\le A}Z_N(\lambda)(K)
       \le e^{\delta_N(b)}H_{N,A}.
 \label{cw:sup-mass}
\end{equation}
The strict inequality \(b>A\) gives the decay in
\eqref{cw:tail-domination}. The two signs control both
prelimit mark tails. These tail conclusions are in probability.

\begin{proof}[Proof of Theorem~\ref{main:IMW}]
Fix \(A\), the test family, and \(b>A\). Write
\(\nu_\kappa=M_\kappa|_K\otimes\rhoR\).
Corollary~\ref{p5:field-corollary} gives
\begin{equation}
 U_N:=(Z_N,\Xi_N,\boldsymbol{\Ytest}_N)
       \Longrightarrow
 U:=(M_\kappa|_K,\nu_\kappa,\boldsymbol G)
 \label{cw:joint-premise}
\end{equation}
in the Polish space
\[
 \mathcal E=\mathcal M_+(K)\times
                 \mathcal M_+(K\times\R)\times\R^m.
\]
We append the curve while retaining all coordinates of \(U_N\).

\emph{Finite-matrix continuity and measurability.}
For sufficiently large fixed \(N\), the exponent range is
positive and compact. Determinant powers are continuous in the
matrix, spatial point, and exponent, including at determinant
zeros. The entry moment assumptions give a finite-\(N\)
integrable dominating power on the compact point/parameter
domain. Their normalizing constants are therefore continuous,
and the uniform scalar estimate \eqref{window:eq:cw-scalar} gives a
positive minimum there. For each matrix realization the
normalized density is continuous and bounded on this compact
domain. Integration gives continuity of
\(\lambda\mapsto Z_N(\lambda)\), even in total variation.
The same argument is locally uniform in the matrix entries,
so the matrix-to-curve map is continuous and hence measurable.
Thus \(\mathcal Z_N\) is a random element of \(\mathcal C_A\).

\emph{The truncated continuous map.}
For \(L>0\), choose a continuous \(0\le\chi_L\le1\), equal
to one on \([-L,L]\) and zero outside \([-L-1,L+1]\).
For a finite positive measure \(\nu\) on \(K\times\R\), define
\begin{equation}
 T_L(\nu)_\lambda(\dd z)
   =\int_\R e^{\lambda r-\lambda^2/2}\chi_L(r)
                                      \nu(\dd z,\dd r),
                 \qquad |\lambda|\le A.
 \label{cw:truncated-map}
\end{equation}
Uniform continuity of the kernel on the compact
parameter/mark rectangle gives \(T_L(\nu)\in\mathcal C_A\).

The map \(T_L:\mathcal M_+(K\times\R)\to\mathcal C_A\)
is continuous. To see this, let \(\nu_n\) converge weakly
to \(\nu\); their total masses are bounded. The functions
\[
 (z,r)\longmapsto f(z)e^{\lambda r-\lambda^2/2}\chi_L(r),
 \quad \norm f_\infty\le1,\quad\operatorname{Lip}(f)\le1,
 \quad |\lambda|\le A,
\]
form a uniformly bounded equicontinuous family with common
compact support. Finite uniform nets for this family, weak
convergence for each member of a net, and bounded total
masses give
\[
 \sup_{|\lambda|\le A}
 d_{\BL,K}\bigl(T_L(\nu_n)_\lambda,T_L(\nu)_\lambda\bigr)
       \longrightarrow0.
\]
Applying this continuous map to \eqref{cw:joint-premise},
and retaining \(U_N\), gives, for each fixed \(L\),
\begin{equation}
 (U_N,T_L(\Xi_N))
       \Longrightarrow (U,T_L(\nu_\kappa)).
 \label{cw:truncated-joint}
\end{equation}

\emph{Uniform approximation of the curve.}
By \eqref{cw:tilt}, total variation bounds the curve distance:
\begin{align}
 \Delta_{N,L}
  &:=d_A(\mathcal Z_N,T_L(\Xi_N))\nonumber\\
  &\le (e^{\delta_N(b)}-1)H_{N,A}
       +\Xi_N\!\left(e^{A|r|}\one_{\{|r|>L\}}\right).
 \label{cw:curve-error}
\end{align}
For the first term, compare \(\mathcal Z_N\) with the
untruncated kernel \(e^{\lambda r-\lambda^2/2}\) applied to
\(\Xi_N\). The uniform bound
\(|e^{e_N}-1|\le e^{\delta_N(b)}-1\) controls this error.
Removing \(\chi_L\) from that kernel gives an error bounded by the
second term.
The first term tends to zero in probability by
\eqref{cw:normalizer-control} and \eqref{cw:weighted-mass};
the second is controlled by \eqref{cw:tails}. Hence
\begin{equation}
 \lim_{L\to\infty}\limsup_{N\to\infty}
       \Prob\{\Delta_{N,L}>\eta\}=0
       \qquad(\eta>0).
 \label{cw:curve-approximation}
\end{equation}

On the specified Gaussian space set
\[
 \mathcal Z_\lambda
   =M_\kappa|_K\,e^{-\lambda^2/2}
                        \int_0^\infty e^{\lambda r}\rhoR(\dd r).
\]
The Rayleigh law has every fixed exponential moment and
\(M_\kappa(K)<\infty\) almost surely. Thus
\(\mathcal Z\in\mathcal C_A\), and
\begin{equation}
 d_A(T_L(\nu_\kappa),\mathcal Z)
   \le M_\kappa(K)
       \int_{|r|>L}e^{A|r|}\rhoR(\dd r)
       \longrightarrow0\quad\text{almost surely}.
 \label{cw:limit-approximation}
\end{equation}
Only almost-sure finiteness of the limiting mass is used here.

\emph{Joint convergence after removing the mark cutoff.}
Choose a compatible metric \(d_{\mathcal E}\le1\) on
\(\mathcal E\), and use
\[
 d((u,\zeta),(u',\zeta'))=
 d_{\mathcal E}(u,u')+(1\wedge d_A(\zeta,\zeta'))
\]
on \(\mathcal E\times\mathcal C_A\). The actual joint
approximants are \((U_N,T_L(\Xi_N))\) and
\((U,T_L(\nu_\kappa))\). Their distances from the respective
targets satisfy
\begin{align}
 d((U_N,\mathcal Z_N),(U_N,T_L(\Xi_N)))
   &=1\wedge\Delta_{N,L},\nonumber\\
 d((U,\mathcal Z),(U,T_L(\nu_\kappa)))
   &=1\wedge d_A(\mathcal Z,T_L(\nu_\kappa)).
 \label{cw:joint-distances}
\end{align}
Thus \eqref{cw:curve-approximation} and
\eqref{cw:limit-approximation} verify both approximation
conditions, while \eqref{cw:truncated-joint} gives convergence
at every fixed \(L\). Applying \eqref{sp:approximation-proof},
first as \(N\to\infty\) and then as \(L\to\infty\), gives
\[
 (U_N,\mathcal Z_N)\Longrightarrow(U,\mathcal Z).
\]
Both approximations keep the spatial measure, marked measure,
and original field-test vector unchanged. The limiting
coordinates retain their specified Gaussian coupling.

\emph{Identification of the profile.}
Integration by parts and completion of the square give
\[
 \int_0^\infty r e^{-r^2/2+\lambda r}\dd r
   =1+\lambda\sqrt{2\pi}\,e^{\lambda^2/2}
        \Phi_{\mathrm{nor}}(\lambda).
\]
Since \(M_\kappa=\cSH D_\kappa\), this identifies
\(\mathcal Z_\lambda=F(\lambda)D_\kappa|_K\) with
\eqref{sp:cw:profile}. The integral formula gives positivity
for every real \(\lambda\), and
\(F(0)=\cSH\) agrees with the critical coordinate of
\eqref{cw:joint-premise}. This proves \eqref{sp:main-limit}.
\end{proof}

The window radius \(A\), the margin \(b>A\), and the fraction
of the history used in Proposition~\ref{cw:inputs} are fixed before
\(N\to\infty\). The mark cutoff is removed after the matrix
limit. The result retains each prescribed finite family of
original field tests and both signs of the window parameter.

\begin{earlytail}[Exponential mark tails in probability]
\label{main:T}
For every fixed $a,\delta>0$,
\[
 \lim_{L\to\infty}\limsup_{N\to\infty}
 \PP\left\{\int e^{a|u|}\ind_{\{|u|>L\}}
                          \Xi_N(\dd z,\dd u)>\delta\right\}=0.
\]
The displayed tail bound holds in probability. It does not assert an
exponential moment bound in expectation for the untruncated mass.
\end{earlytail}

\begin{proof}
Take \(A=a\) in \eqref{cw:tails}.
\end{proof}

\subsection{Parameter choices, limits, and scope}
\label{main:limit-order}

The complete order of parameter choices is given in
Section~\ref{main:parameter-schedule}. The tests, exponential coefficients,
truncation levels, derivative and moment requirements, and required
bounds for weighted errors on exceptional events are fixed before the
logarithmic cutoff and flow exponents.
All exponents are fixed before $N\to\infty$. The retained field history
is fixed before the separate auxiliary coupling for the scalar estimate
is chosen. This later choice does not change the history. The global
event and the nearby masses use a separate fixed initial fraction of
scales. These choices follow the order specified in that section.

Barrier removal takes $N\to\infty$, then $\varepsilon\downarrow0$
at fixed $B$, then $B\to\infty$. Approximation of the original field fixes
the coarse scale parameter before the matrix limit, refines quadrature
next, and then sends that scale parameter to infinity. The compact window is fixed
before the matrix limit and its exponential mark cutoff is removed
afterward. Rates are not asserted uniformly in these later limits.

The later sections give the analytic, Gaussian, and probabilistic
arguments in full and state the hypotheses of the external results
where they are used. We specify the cited version when using a theorem
or formula. The iid GMC theorem of Cipolloni--Landon and the normal-matrix
theorem of Bourgade--Dubach--Hartung--Keles treat strict subcritical
ranges; see \cite{CL,BDHK}. The theorem above concerns the joint critical,
marked, and two-sided compact-window limit.

\endgroup

\begingroup
\newcommand{\E}{\mathbb E}
\newcommand{\Pp}{\mathbb P}
\newcommand{\PP}{\mathbb P}
\newcommand{\R}{\mathbb R}
\newcommand{\C}{\mathbb C}
\newcommand{\D}{\mathbb D}
\newcommand{\dd}{\,\mathrm d}
\newcommand{\ind}{\mathbf1}
\newcommand{\eps}{\varepsilon}
\newcommand{\cF}{\mathcal F}
\newcommand{\cG}{\mathcal G}
\newcommand{\cM}{\mathcal M}
\newcommand{\rhoR}{\rho_{\mathrm R}}
\newcommand{\norm}[1]{\left\|#1\right\|}
\newcommand{\src}[1]{\nolinkurl{#1}}
\newcommand{\Cov}{\operatorname{Cov}}
\newcommand{\Var}{\operatorname{Var}}
\newcommand{\dist}{\operatorname{dist}}
\newcommand{\Tr}{\operatorname{Tr}}
\section{Canonical Gaussian chaos, barriers, and endpoint marks}\label{chap:gaussian}

 {The purpose of this section is to construct the Gaussian comparison measure together with its field and prove the barrier and mark estimates needed for matrix identification.}
In particular, we apply canonical critical GMC to the   Gaussian  
family used in Section~\ref{chap:main} and proves the barrier and mark
estimates needed for its matrix identification. The   family and
its limiting generalized field are realized together. We first identify
the canonical critical measure on this space, then restrict it by the
integer-scale history and establish the fixed-barrier moment and
endpoint estimates.

The covariance and measure normalizations are those fixed in
\eqref{sp:fb:covariance} and Proposition~\ref{p1:proposition}.
In particular $M_\kappa=\sqrt{2/\pi}\,D_\kappa$, with the derivative
normalization specified in Section~\ref{gaussian:gp2:canonical}.
Canonical GMC, covariance uniqueness, the maximum comparisons and
the continuous Brownian bounds are cited with their hypotheses at
the point of use. Matrix identification is the separate deduction
in Corollary~\ref{sp:fb:identification}.

\subsection{The exact Gaussian family and its normalization}
\label{gaussian:gp2:setting}

 {We now specify the exact Gaussian covariance and retain the spatial corrections that determine the normalization of the limiting measure.}
We use the common bulk disk $K\Subset\D=\{z\in\C:|z|<1\}$,
integer $k_0\ge1$ and parameter $\kappa\ge-1$ from Section~\ref{chap:main}.
The barrier height $B$ is real and fixed. All auxiliary bulk disks
carry restrictions of the same Gaussian family; a larger disk is
used when convolution or additional test supports require it.
Area remains $\dd A=\dd x\dd y$.

We recall the covariance because its exact form, including its
bounded spatial correction, will be used throughout this section.

For $s>0$, define
\[
 u_s=\tanh s,\qquad p_z(u)=u(1-u|z|^2).
\]
The centered, jointly Gaussian family $G_s(z)$ has covariance
\begin{align}
 C_\kappa(s,z;t,w)
 &=-\frac12\log D(u_s,z;u_t,w)
       +\frac\kappa2p_z(u_s)p_w(u_t),\label{gaussian:gp2:covariance}\\
 D(u,z;v,w)
 &=1-uv-uv\big((1-u)|z|^2+(1-v)|w|^2\big)
                    +uv|z-w|^2.\label{gaussian:gp2:D}
\end{align}
Write $v_t(z)=C_\kappa(t,z;t,z)$ and
\begin{align}
 W_t(z)&=e^{2G_t(z)-2v_t(z)},\\
 A_n^B(z)&=\bigcap_{k=k_0}^{n}\{G_k(z)\le2k+B\},\\
 F_t^B(z)&=W_t(z)\ind_{A_{\lfloor t\rfloor}^B(z)},\\
 \mu_t^B(\dd z)&=\sqrt t\,F_t^B(z)\dd A(z),
 \qquad \mu_t(\dd z)=\sqrt t\,W_t(z)\dd A(z).
 \label{gaussian:gp2:measures}
\end{align}
Here $\mu_t^B$ is the Gaussian measure denoted by
$\mu_t^{G,B}$ in \eqref{sp:gaussian-measures}; $W_t$ is the variance-normalized
Gaussian exponential denoted by $W_t^G$ there.
The terminal time $t$ is real, whereas the barrier is sampled only
at integers through $\lfloor t\rfloor$. There is no additional
barrier at a fractional terminal time.

\paragraph{Variance normalization.}
The factor $2v_t(z)$ in \eqref{gaussian:gp2:measures} is not replaced by $2t$. In fact,
\[
 v_t(z)=t-\tfrac12\log\big(4(1-|z|^2)\big)
             +\tfrac\kappa2(1-|z|^2)^2+O(e^{-2t}).
\]
Dropping this bounded, spatially varying correction would change the limiting measure by a nontrivial deterministic multiplier.

\subsubsection{Positivity and uniform covariance estimates}

{Now we shall verify that the covariance defines a Gaussian family throughout the stated parameter range and controls its logarithmic correlations.}
Positivity of the limiting spatial kernel
\eqref{gaussian:gp2:limit-kernel} follows from
\cite[Appendix~D]{CL}, after multiplying its complex-case kernel by
two. The expansion below verifies positivity of the full positive-time
family. It also separates a scalar Gaussian, multiplied by a positive
deterministic function, from an independent residual family. We use this
decomposition for the history restrictions.

The algebraic identity
\begin{equation}
 D(u,z;v,w)=|1-uvz\overline w|^2-p_z(u)p_w(v)
 \label{gaussian:gp2:feature-identity}
\end{equation}
expresses the covariance as a sum of positive kernels. This is
particularly useful at $\kappa=-1$, where the cumulant correction
cannot be treated as an independent Gaussian summand. The logarithmic part is
\begin{align}
 -\tfrac12\log D
 &=-\log|1-uvz\overline w|\nonumber\\
 &\quad+\sum_{a\ge1}\frac{p_z(u)^a p_w(v)^a}{2a}
       (1-uvz\overline w)^{-a}(1-uv\overline z w)^{-a}.
 \label{gaussian:gp2:features}
\end{align}
The coefficient of the function
$p_z(u)^a u^{b+c}z^b\overline z^{\,c}$, paired with its conjugate at $(v,w)$, is
\[
 \frac1{2a}\binom{a+b-1}{b}\binom{a+c-1}{c}\ge0.
\]
Pairing conjugate functions yields a real positive semidefinite kernel. The first term in \eqref{gaussian:gp2:features} has the usual expansion in real and imaginary monomials with nonnegative coefficients. The cumulant term changes only the $(a,b,c)=(1,0,0)$ coefficient, from $1/2$ to $(1+\kappa)/2$. Thus the full kernel is positive semidefinite for every $\kappa\ge-1$, including the endpoint. The resulting consistent finite-dimensional laws define the exact Gaussian family.

If $|z|,|w|\le r_*<1$, then
\[
 (1-r_*^2)(1-uv)+uv|z-w|^2
 \le D(u,z;v,w)\le1-uv+uv|z-w|^2.
\]
Here $uv(2-u-v)\le1-uv$ follows from
$1/(uv)+u+v\ge3$. Consequently, for $s,t\ge s_*>0$,
\begin{equation}
 C_\kappa(s,z;t,w)
 =\min\{s,t,\log_+|z-w|^{-1}\}+O(1),
 \label{gaussian:gp2:log-profile}
\end{equation}
uniformly on fixed bulk compacts. At coincident spatial points the logarithm is $+\infty$.

At equal time, put
\[
 d_t(z)=(1-u_t)(1+u_t-2u_t^2|z|^2).
\]
Then
\[
 D(u_t,z;u_t,w)=\frac{d_t(z)+d_t(w)}2+u_t^2|z-w|^2.
\]
Using $d_t\asymp e^{-2t}$ and the bounded spatial derivatives of its normalized coefficient gives
\begin{equation}
 \E|G_t(z)-G_t(w)|^2\le C e^{2t}|z-w|^2.
 \label{gaussian:gp2:metric}
\end{equation}
The smooth covariance also gives continuous versions jointly on compact subsets of $(0,\infty)\times\D$. All countable-coordinate measurability statements below use this version.

\subsection{Canonical identification on the original Gaussian space}
\label{gaussian:gp2:canonical}

 {In this subsection we shall identify the canonical critical measure on the original Gaussian space while preserving its dependence on the field.}
We first construct the canonical critical measure from the limiting
Gaussian field by smooth mollification. We then compare the exact approximants with those
mollifications while retaining every prescribed finite field vector.
A final measurability argument gives convergence in probability on
the specified space. These are separate steps: identifying the law
of the measure alone would not identify its coupling to the field.

We use Junnila--Saksman--Webb
\cite[Theorems 5.3--5.4]{JSW} for canonical critical GMC
and Junnila--Saksman \cite[Theorem 1.1]{JS} for covariance
uniqueness. The smooth remainder below has the required regularity.
The history restriction is treated separately in
Proposition~\ref{gaussian:gp2:probability}.

\subsubsection{The limiting field and its mixed kernel}

Next the limiting field is constructed  and we compute its covariance with the regularized observations that remain in the joint limit.
For a real $\varphi\in C_c^\infty(\D)$, the variables
\[
 G_\kappa(\varphi)=L^2\!\lim_{t\to\infty}
                  \int G_t(z)\varphi(z)\dd A(z)
\]
exist by dominated convergence of mixed covariances. The dominating singularity is $C+\log_+|z-w|^{-1}$, which is locally integrable in the two spatial variables. The resulting generalized field has covariance
\begin{equation}
 K_\kappa(z,w)=-\log|z-w|
              +\frac\kappa2(1-|z|^2)(1-|w|^2).
 \label{gaussian:gp2:limit-kernel}
\end{equation}
It is jointly Gaussian with the whole exact family and realizes the field in \eqref{sp:kernel} on that family's probability space.

Taking the second time to infinity in \eqref{gaussian:gp2:covariance} gives the exact mixed kernel
\begin{align}
 C_{t,\infty}(z,w)
 &=-\tfrac12\log\left(u_t|z-w|^2
              +(1-u_t)(1-u_t|z|^2)\right)\nonumber\\
 &\quad+\tfrac\kappa2p_z(u_t)(1-|w|^2).
 \label{gaussian:gp2:mixed-kernel}
\end{align}
It is interpreted after integration in $w$. With
\[
 a_t(z)^2=\frac{(1-u_t)(1-u_t|z|^2)}{u_t}=O(e^{-2t}),
\]
the difference from \eqref{gaussian:gp2:limit-kernel} consists of $O(e^{-2t})$ terms and
$-\frac12\log(1+a_t(z)^2/|z-w|^2)$. Since
\[
 \int_{|x|\le R}\log(1+a^2/|x|^2)\dd A(x)
       \le C_Ra^2(1+|\log a|),
\]
we obtain
\begin{equation}
 \sup_{z\in K}\left|
 \Cov(G_t(z),G_\kappa(\varphi))
       -\int K_\kappa(z,w)\varphi(w)\dd A(w)\right|
 \le C_\varphi(1+t)e^{-2t}.
 \label{gaussian:gp2:mixed-rate}
\end{equation}

\subsubsection{From covariance uniqueness to the specified coupling}

{Now we retain the field observations in the uniqueness argument so that it identifies the measure's coupling as well as its law.}
Let $\rho$ be a fixed smooth, compactly supported, nonnegative
mollifier of unit mass,
and put $\widehat G_t=G_\kappa*\rho_{e^{-t}}$. On $K$, the
convolution uses only the larger bulk disk for all sufficiently
large $t$. Canonical critical GMC gives
\[
 \sqrt t\,e^{2\widehat G_t-2\Var\widehat G_t}\dd A
       \longrightarrow M_\kappa
       =\sqrt{\frac2\pi}\,D_\kappa
 \quad\hbox{in probability}.
\]
The derivative normalization defining $D_\kappa$ is
$(2\Var\widehat G_t-\widehat G_t)e^{2\widehat G_t-2\Var\widehat G_t}\dd A$.
Here the remainder in \eqref{gaussian:gp2:limit-kernel} is smooth and hence belongs
to $H^s_{\rm loc}$ for every $s>2$, as required in \cite{JSW}.
Positivity is required of the full covariance, already proved above;
the smooth remainder need not itself be positive semidefinite.
The canonical limit is atomless: the stationary reference derivative
measure is atomless by \cite[Theorem 4]{DRSV};
\cite[Theorem 5]{DRSVSH} identifies its Seneta--Heyde limit as
$\sqrt{2/\pi}$ times that measure. The continuous positive multiplier
in the decomposition proof of \cite[Theorem 5.3]{JSW} preserves
atomlessness. This also covers $\kappa=-1$.

The equal-time covariance difference between $G_t$ and $\widehat G_t$ is uniformly bounded and tends uniformly to zero away from the spatial diagonal. This follows by comparing both logarithmic kernels with
$\log(1/(e^{-t}\vee|z-w|))$ and using uniform off-diagonal convergence. Equation \eqref{gaussian:gp2:mixed-rate} also gives
\begin{equation}
 \sup_{z\in K}|\Cov(G_t(z)-\widehat G_t(z),G_\kappa(\varphi))|
       \longrightarrow0.
 \label{gaussian:gp2:cross-condition}
\end{equation}
These bounds compare the normalized exponentials of $2G_t$ and
$2\widehat G_t$ with common deterministic reference measure
$\sqrt t\,\dd A$. The covariance-uniqueness theorem allows this
reference measure to depend on the sequence index. To identify the
coupling as well, we apply that theorem on a larger space.

To retain a finite Gaussian vector $V=(G_\kappa(\varphi_j))_{j\le m}$,
adjoin $m$ disjoint copies of $[0,1]$ to $K$. Rescale the Euclidean
metric on each component to diameter at most one, and set the distance
between distinct components to two. This is a compact doubling metric
space: balls of radius less than two stay in one component, while a
larger ball is covered using finitely many component covers.
On $K$ the two fields are $2G_t$ and $2\widehat G_t$; on the $j$th
interval both are $V_j$. The reference measure is $\sqrt t\,\dd A$ on
$K$ and Lebesgue measure on each interval. For each fixed $t$, smooth
covariances on $K$, constant interval fields, and the positive distances
between components give the H\"older regularity of the covariance
pseudometric required in \cite[Definition 2.1]{JS}; uniform H\"older
constants in $t$ are not required. The canonical limit on $K$ and the
interval densities give an atomless limiting measure. Covariance
differences on $K\times K$ are uniformly bounded and vanish uniformly
off the diagonal; cross-component differences involving $K$ are twice
\eqref{gaussian:gp2:cross-condition}, and interval-pair differences are zero.
Thus all hypotheses of \cite[Theorem 1.1]{JS} hold for every deterministic
sequence $t\to\infty$. The interval masses are
$e^{V_j-\Var(V_j)/2}$, so taking their logarithms recovers the original $V$, not independent copies. Covariance uniqueness therefore proves
\begin{equation}
 (\mu_t,V)\Longrightarrow(M_\kappa,V).
 \label{gaussian:gp2:joint-canonical-law}
\end{equation}

We now pass from the joint law to convergence on the original
space. The measurable-target principle is recorded in
\cite[Lemma~4.1]{JS}; the argument below keeps the weak measure
topology explicit. Choose countably many $\varphi_j$ generating the
limiting generalized field. The canonical target $M_\kappa$ is measurable
with respect to that field. For every bounded continuous function
$F$ of a finite measure, \eqref{gaussian:gp2:joint-canonical-law} gives
\[
 \E[H F(\mu_t)]\longrightarrow\E[H F(M_\kappa)]
\]
first for bounded continuous finite-coordinate $H$, then for all bounded field-measurable $H$ by $L^1$ approximation and increasing finite sigma-fields. Taking $H=F(M_\kappa)$ and using marginal convergence for $F^2$ gives
$\E|F(\mu_t)-F(M_\kappa)|^2\to0$.
To make the topology explicit, take a countable uniformly dense family
$(f_\ell)$ in the unit ball of $C(K)$, including the constant one, and
put $F_\ell(\nu)=\arctan(\nu(f_\ell))$. The bounded metric
\[
 d_*(\nu,\nu')=\sum_{\ell\ge1}2^{-\ell}
       |F_\ell(\nu)-F_\ell(\nu')|
\]
metrizes weak convergence of finite positive measures: the constant
test controls total mass, after which uniform approximation controls
every continuous test. Each summand converges in probability and the
tail is deterministically summable. Hence $d_*(\mu_t,M_\kappa)\to0$
in probability.

This argument applies along every deterministic real sequence $t\to\infty$. It does not need temporal independent increments or an interpolation theorem at the fractional terminal time.

\paragraph{Boundary of the bulk disk.}
The canonical measure gives deterministic area-null sets zero mass. 
One way to check this without a finite first moment is to use the local decomposition in \cite{JSW}: 
the chaos is a positive continuous multiple of a stationary reference critical measure. 
For a locally finite stationary random measure $\nu$ and a fixed area-null set $E$, Tonelli gives
$\int_{[0,1]^2}\nu(E+x)\dd x=0$ pathwise. 
Stationarity then gives $\Pp(\nu(E)>0)=0$. The local positive multiplier preserves nullity, in particular on $\partial K$.

The construction and convergence above precede all history
restrictions. They use neither fixed-barrier moment estimates nor
matrix identification, and require no finite first moment of the
unrestricted critical mass. They therefore supply the canonical
measure for the subsequent Gaussian barrier arguments.

\subsection{A sharp summable maximum estimate}
\label{gaussian:gp2:maximum}

{In this subsection we bound the field maximum sharply enough to sum the probabilities of barrier crossings over all late scales.}
To pass from a finite barrier history to the full history, we must
make crossings after a large integer time unlikely throughout $K$.
The logarithmic correction to the maximum gives a summable bound
for those late crossings. Its dependence on the barrier height
will also be used in Corollary~\ref{aggregate:cor:au-maximum}.

We combine Acosta's continuum comparison and cell-supremum estimate
\cite[Proposition 2.1 and Lemma 2.2]{Acosta} with
Ding--Roy--Zeitouni's subset maximum bound
\cite[Proposition 1.1, inequality (3)]{DRZ}.
The latter retains the sharp exponent $2$ when we sum over the
cell fluctuations; Acosta's general exponential tail alone does
not supply that exponent.

\begin{proposition}[Sharp upper tail and late crossings]
\label{gaussian:gp2:max-prop}
For sufficiently large integer $k$ and $\lambda\ge1$,
\begin{equation}
 \Pp\!\left(\sup_{z\in K}G_k(z)
             >2k-\tfrac34\log k+\lambda\right)
       \le C(1+\lambda)e^{-2\lambda}.
 \label{gaussian:gp2:sharp-max}
\end{equation}
For fixed $B$, if
\[
 \mathcal L_L^B=
 \{\exists k>L,\ z\in K:\ G_k(z)>2k+B\},
\]
then
\begin{equation}
 \Pp(\mathcal L_L^B)
       \le C_B\frac{1+\log(2+L)}{\sqrt{1+L}}\longrightarrow0.
 \label{gaussian:gp2:late-crossing}
\end{equation}
For sufficiently large $B$, the probability of any crossing at any integer $k\ge k_0$ is at most $C(1+B)e^{-2B}$.
\end{proposition}
\begin{proof}
The logarithmic profile \eqref{gaussian:gp2:log-profile} and metric bound \eqref{gaussian:gp2:metric} imply the hypotheses of the continuum comparison. The local metric condition there is weaker: when $|z-w|\le e^{-k}$, the right side of \eqref{gaussian:gp2:metric} is at most $Ce^k|z-w|$.

After a fixed spatial rescaling and a finite cover, write
$\epsilon\asymp e^{-k}$ for the resulting mesh size.
The comparison field on each mesh cell is
$a(x)\xi_v+\psi_v(x)$, where the discrete modified branching field $\xi$ is independent of the cell fields $\psi$, and $1/2\le a(x)\le1$ for large $k$.
To apply the subset estimate, we check \cite[Assumption (A.0)]{DRZ}.
Relabel the mesh points $\epsilon a$ by
$a\in V:=\{0,\ldots,N-1\}^2$, where
$N=\lceil\epsilon^{-1}\rceil$ is the grid side length and
$L_\epsilon=\log(1/\epsilon)$.
Acosta's field has equal variances $L_\epsilon\le\log N$ and,
for distinct indices,
\[
 \Cov(\xi_a,\xi_b)\ge L_\epsilon-\log\|a-b\|_\infty-C,
 \qquad
 \E(\xi_a-\xi_b)^2\le2\log_+|a-b|+2C.
\]
Thus (A.0), including its variance-difference term, holds with a fixed
constant. Since $\log N=k+O(1)$, the centering
$m_N=2\log N-\frac34\log\log N$ differs
from $2k-\frac34\log k$ by $O(1)$.

Classify cells according to the positive part of their cell supremum: $\Gamma_j$ contains those whose supremum lies in $(j-1,j]$, with $j=0$ covering nonpositive suprema. The uniform Gaussian cell tail gives
\[
 \E\frac{|\Gamma_j|}{|V|}\le C e^{-c j^2}
\]
after changing constants. For $M\ge0$, $a x+j\ge M$ and $a\in[1/2,1]$ imply $x\ge M-2j$. Conditionally on $\Gamma_j$, the subset estimate in dimension two therefore bounds the crossing probability by
\[
 C(1+\lambda)e^{-2\lambda+4j+C}
                \sqrt{|\Gamma_j|/|V|}.
\]
Here the subset estimate is used with $z=\lambda$ and
$y=2j+C$ after a fixed lift; its stated range permits negative $z-y$.
The sets $\Gamma_j$ are independent of $\xi$, as required for conditioning. Jensen's inequality and summation over $j$ give \eqref{gaussian:gp2:sharp-max}. A fixed centering shift changes only its constant.

Substitute $\lambda=B+\frac34\log k$ for sufficiently large $k$. The resulting tail is bounded by
$C_B(1+\log k)k^{-3/2}$. Summation gives \eqref{gaussian:gp2:late-crossing}. For sufficiently large $B$, the same substitution gives the uniform
bound $C(1+B+\log k)e^{-2B}k^{-3/2}$ at all sufficiently large
integer times. Sum it and treat the finitely many earlier continuous
Gaussian maxima by Gaussian concentration. This gives the stated
all-time bound. No independence in $k$ is used.
\end{proof}

\subsection{Random coarse boundaries and convergence in probability}
\label{gaussian:gp2:restriction}

Next we prove that random barrier boundaries carry no limiting mass, allowing finite history restrictions to pass through the Gaussian limit.
Let
\[
 E_L^B=\{z\in K:G_k(z)\le2k+B\text{ for }k_0\le k\le L\},
 \qquad E_\infty^B=\bigcap_{L\ge k_0}E_L^B.
\]
The fixed-barrier target is $Z_B=\ind_{E_\infty^B}M_\kappa$.
For a fixed finite history, weak convergence can pass through this
restriction once its boundary has zero limiting mass. The random part of the boundary
is contained in level sets of the same field that generates the measure, so
atomlessness alone does not establish this nullity. To prove this nullity, we represent the field as a residual family plus
an independent scalar Gaussian multiplied by a positive function. We then
condition on the residual family and integrate over that scalar.

\begin{lemma}[Zero mass on barrier levels, including $\kappa=-1$]
\label{gaussian:gp2:positive-feature}
The kernel
$C_\kappa(s,z;t,w)-f_s(z)f_t(w)$ is positive semidefinite, where
\begin{equation}
 f_s(z)=\tfrac12p_z(u_s)^2>0,
 \qquad f_\infty(z)=\tfrac12(1-|z|^2)^2.
 \label{gaussian:gp2:feature}
\end{equation}
For every fixed integer $k\ge k_0$,
\begin{equation}
 M_\kappa\{z\in K:G_k(z)=2k+B\}=0\quad\hbox{almost surely}.
 \label{gaussian:gp2:level-null}
\end{equation}
\end{lemma}
\begin{proof}
In \eqref{gaussian:gp2:features}, the $(a,b,c)=(2,0,0)$ term has coefficient $1/4$ and is exactly $f_s(z)f_t(w)$. It is not the coefficient altered by $\kappa$. Removing it preserves positivity even at $\kappa=-1$.

Represent the family in law as $G_s=R_s+\xi f_s$, where $\xi$ is standard Gaussian and independent of the residual family $R$. This can be done on an enlarged space; the desired assertion concerns the joint law. The residual normalized measures obey the exact identity
\[
 \mu_t^R=e^{-2\xi f_t+2f_t^2}\mu_t.
\]
Uniform convergence $f_t\to f_\infty$ on $K$ and the canonical convergence show that
\[
 M_R=e^{-2\xi f_\infty+2f_\infty^2}M_\kappa
\]
is their limit in probability. Since each residual approximant is measurable with respect to $R$, so is $M_R$; in particular $M_R$ is independent of $\xi$. Equivalently,
\[
 M_\kappa(\dd z)=e^{2\xi f_\infty(z)-2f_\infty(z)^2}M_R(\dd z).
\]
Condition on $R$ and $M_R$. At any fixed $z$, the level equation
$R_k(z)+\xi f_k(z)=2k+B$ fixes one value of $\xi$, because $f_k(z)>0$. Integrating first in the scalar Gaussian, then in $M_R$, gives zero by Tonelli. This proves \eqref{gaussian:gp2:level-null} without assuming
$\E M_R(K)<\infty$. Since it is a property of the joint law, it
holds on the specified Gaussian realization. For each fixed $B$,
a countable intersection gives the assertion simultaneously for
all integers $k\ge k_0$.
\end{proof}

\begin{proposition}
\label{gaussian:gp2:probability}
As real $t\to\infty$,
\begin{equation}
 \mu_t^B\longrightarrow Z_B:=\ind_{E_\infty^B}M_\kappa
 \quad\hbox{in probability in }\cM_+(K).
 \label{gaussian:gp2:prob-conv}
\end{equation}
Here $\cM_+(K)$ has the weak topology of finite positive measures.
\end{proposition}
\begin{proof}
Fix an integer $L\ge k_0$. Continuity of the coarse fields puts
the boundary of $E_L^B$ in the finite union of the level sets in
\eqref{gaussian:gp2:level-null}, together with the deterministic
disk boundary. Both have zero $M_\kappa$ mass by the preceding
lemma and Section~\ref{gaussian:gp2:canonical}.

From any deterministic sequence of real times tending to infinity,
canonical convergence in probability gives a subsequence on which
the measures converge weakly almost surely. On each sample in
this event, the set $E_L^B$ is fixed and its boundary has zero
limiting mass. Portmanteau therefore gives convergence of the
restrictions. By the subsequence principle,
\[
 \ind_{E_L^B}\mu_t\longrightarrow\ind_{E_L^B}M_\kappa
 \quad\hbox{in probability in }\cM_+(K).
\]
In particular this includes total mass.

The late-crossing estimate makes this finite-history approximation
uniform over all later terminal times. On $(\mathcal L_L^B)^c$,
\[
 E_{\lfloor t\rfloor}^B=E_L^B=E_\infty^B
 \qquad\hbox{throughout }K,\quad t\ge L.
\]
For the bounded weak metric $d_*$ used in
Section~\ref{gaussian:gp2:canonical} and every $\delta>0$, it follows that
\begin{align}
 \Pp\{d_*(\mu_t^B,Z_B)>\delta\}
 &\le\Pp(\mathcal L_L^B)\nonumber\\
 &\quad+\Pp\{d_*(\ind_{E_L^B}\mu_t,
                 \ind_{E_L^B}M_\kappa)>\delta\}.
 \label{gaussian:gp2:restriction-error}
\end{align}
First take $t\to\infty$ with $L$ fixed, and then $L\to\infty$
using \eqref{gaussian:gp2:late-crossing}. This proves the stated
weak convergence in probability without any moment assumption
on $M_\kappa(K)$.
\end{proof}

The remaining step to the fixed-barrier $L^2$ target is
integrability. The following three-point and shell estimates provide
a uniform third moment at fixed $B$, before any limiting moment
assertion is used.

\subsection{A multiplicative comparison with a three-leaf Brownian tree}
\label{gaussian:gp2:relative}

In this subsection we compare three weighted Gaussian histories with a Brownian tree to obtain bounds that survive the critical normalization.
Indeed, we obtain a uniform third moment by integrating a three-point
estimate over the spatial separations. We retain the pointwise
bound also when spatial points coincide. The closest pair determines two branching times for a
Brownian comparison tree. We compare the barrier probabilities
multiplicatively, so that the estimate remains useful after the
large normalizing factor for the product of the three exponential
weights is restored.

Let $n=\lfloor t\rfloor$. Relabel three spatial points so that $(z_2,z_3)$ is a closest pair, and set
\begin{align}
 r&=\left\lfloor\min\{n,\log_+|z_2-z_3|^{-1}\}\right\rfloor,\\
 q&=\left\lfloor\min\{n,\log_+|z_1-z_2|^{-1},
                           \log_+|z_1-z_3|^{-1}\}\right\rfloor.
 \label{gaussian:gp2:qr}
\end{align}
Then $0\le q\le r\le n$. The two distances from $z_1$ differ by at most a factor of two. A standard Brownian tree has common stem $[0,q]$, one singleton branch after $q$, a pair stem $[q,r]$, and two separate branches after $r$. We write its three labeled paths as $S^i$.

Equation \eqref{gaussian:gp2:log-profile} bounds the discrepancy
between each exact history covariance and its tree counterpart
by a uniform constant. When spatial points coincide, we retain
a separate coordinate label for each path. The Gaussian covariance
may therefore be singular, and the comparison must also apply in this case.

\begin{lemma}[Relative comparison on the three-leaf tree]
\label{gaussian:gp2:tree-comparison}
Suppose a centered Gaussian array $X_k^i$, $i=1,2,3$, $k_0\le k\le n$, differs in covariance from $S_k^i$ by at most $M\ge1$. Put
\[
 a_k=A\sqrt{\log(2+\dist(k,\{0,q,r,n\}))},
 \qquad A^2\ge24M.
\]
There are uniform $C,c>0$ and a nonnegative random variable $W$, independent of the tree, with $\Pp(W>x)\le Ce^{-cx^2}$, such that for arbitrary deterministic thresholds $b_k^i$,
\begin{equation}
 \Pp(X_k^i\le b_k^i\ \forall i,k)
 \le C\E_W\Pp(S_k^i\le b_k^i+2a_k+W\ \forall i,k).
 \label{gaussian:gp2:relative-bound}
\end{equation}
No covariance inverse or distinctness of the formal coordinates is required.
\end{lemma}
\begin{proof}
Let $v_k^i=\Var X_k^i$. Add independent $N(0,3M)$ noises $\xi_k^i$ to $X$. Compare the result with
\[
 Y_k^i=S_k^i+Z+\eta_k^i,
 \quad Z\sim N(0,2M),
 \quad \eta_k^i\sim N(0,v_k^i-k+M),
\]
where all noises and the tree are independent. The last variances lie in $[0,2M]$. The arrays have equal diagonal variances $v_k^i+3M$, and the off-diagonal covariances of $Y$ dominate those of the enlarged $X$. Gaussian rectangle comparison gives
\[
 \Pp(X_k^i+\xi_k^i\le b_k^i+a_k\ \forall i,k)
 \le\Pp(Y_k^i\le b_k^i+a_k\ \forall i,k).
\]
This version of comparison follows by covariance interpolation with a product of smooth decreasing one-dimensional cutoffs: cross second derivatives are nonnegative and diagonal covariance differences vanish. Degenerate laws follow by regularization.

The four distinguished times in the distance formula are the initial
time, the two branching times and the terminal time. Only a bounded
number of the $3n$ formal coordinates have any given distance from this
set. Hence
\[
 \sum_{i,k}\Pp(\xi_k^i>a_k)
       \le C\sum_{d\ge0}(2+d)^{-A^2/(6M)}<\infty.
\]
Every success probability is at least $1/2$. Independence and
$\log(1-x)\ge-2x$ for $x\le1/2$ therefore give a uniform positive lower bound for the probability that all added noises respect their allowances. This event is independent of $X$.

Finally put
\[
 W=|Z|+\max\{0,\sup_{i,k}(-\eta_k^i-a_k)\}.
\]
Gaussian tails, the same shell count, and $(x+a)^2\ge x^2+a^2$ give its stated uniform tail. The event for $Y$ implies
$S_k^i\le b_k^i+2a_k+W$. Dividing by the uniform positive product probability proves \eqref{gaussian:gp2:relative-bound}.
\end{proof}

\subsubsection{The exact triple tilt}

{Next we compute how three exponential weights shift the Gaussian history while keeping their normalizing factor explicit.}
The normalized triple weight changes only the mean of the Gaussian
history. We keep its scalar normalizer separate from the resulting
barrier probability. It is
\begin{equation}
 \E\prod_{i=1}^3W_t(z_i)
 =\exp\left(4\sum_{i<j}C_\kappa(t,z_i;t,z_j)\right)
 \le C e^{8q+4r}.
 \label{gaussian:gp2:triple-normalizer}
\end{equation}
Under the normalized triple tilt, the mean at history coordinate $(k,z_i)$ is $2\sum_j C_\kappa(k,z_i;t,z_j)$. After removing the own-root part $2k$, its excess, up to a uniformly bounded error, is
\[
 b_1(k)=4\min(k,q),\qquad
 b_2(k)=b_3(k)=2\min(k,q)+2\min(k,r).
\]
Thus the tree drift is $4$ on the common stem, $2$ on the pair stem, and zero on terminal branches. Apply Lemma~\ref{gaussian:gp2:tree-comparison} to the centered tilted history. All bounded mean errors are absorbed by a fixed height lift. The scalar factor in \eqref{gaussian:gp2:triple-normalizer} is retained.

\subsubsection{Integer observations to a continuous tree}

Now we fill the gaps between integer observations with Brownian bridges at a cost bounded independently of the terminal scale. 
The comparison lemma still observes the tree at integers. To use
the continuous Brownian estimates below, fill in the intervening
edges while retaining a uniform positive probability of respecting
the enlarged barrier. Condition on the tree skeleton at integers and insert independent unit Brownian bridges on each distinct edge. Impose bridge constraints only on the observed intervals in $[k_0,n]$; no constraint is added before $k_0$. Count a shared edge once, not once for every leaf label. A bridge's positive maximum has a Gaussian tail. Requiring it to be at most a sufficiently large multiple of
$\sqrt{\log(2+d)}$ on an interval at distance $d$ from these four times has a uniformly positive product probability over all edges.

On this bridge event, the integer constraints imply continuous constraints with a slightly larger pinned envelope. The deterministic drift is affine between integer branch times, so interpolation introduces no additional drift error. Since
$\sqrt{\log(2+d)}\le C_\alpha(1+d^\alpha)$ for a fixed $0<\alpha<1/2$, the resulting continuous barrier has height
\[
 H=C_0+|B|+W\ge1
\]
and envelope $C_1\dist(s,\{0,q,r,n\})^\alpha$. We may enlarge it, on each edge, to the usual endpoint-pinned curve. The additional designated time $r$ along the singleton branch only strengthens the original constraint.

A global $\sqrt{\log n}$ height lift would enlarge the later
exponential height factor, while an additive error would remain
multiplied by the triple normalizer. Instead, the bridge construction
costs only a fixed multiplicative constant. Conditional on the
integer skeleton and on $W$, the
independent bridge event has probability $p_{\rm br}\ge c>0$.
Writing $A_{\rm int}$ and $A_{\rm cont}$ for the two constraint events,
$A_{\rm int}\cap E_{\rm br}\subset A_{\rm cont}$, and therefore
$\Pp(A_{\rm int}\mid W)\le c^{-1}\Pp(A_{\rm cont}\mid W)$.

\subsection{The Brownian three-point estimate}
\label{gaussian:gp2:brownian}

In this subsection we combine survival estimates along the Brownian tree to bound three-point expectations, including configurations with coincident spatial points.

We next estimate survival on each distinct edge of the tree.
Positive-length stems use a transition density, whereas terminal
branches use a one-sided survival bound. A stem of zero length
is omitted; it is never assigned a transition density. This
separation will include branch times that coincide or reach the
terminal time.

For $x\le H$, write $h_H(x)=1+H-x$. The continuous Brownian upper bounds in \cite[Proposition 2.1, (2.1), and Proposition 2.2, (2.3)]{CHL}, after reflection and a height lift of one, give the following estimates. For an edge of length $\ell\ge1$, a transition constrained below
$H+C\min(s,\ell-s)^\alpha$ has density at most
\begin{equation}
 \frac{C h_H(x)h_H(y)}{(1+\ell)^{3/2}}\dd y,
 \qquad x,y\le H.
 \label{gaussian:gp2:killed-density}
\end{equation}
A one-sided continuation of length $\ell\ge0$ survives with probability at most
\begin{equation}
 \frac{C h_H(x)}{\sqrt{1+\ell}}.
 \label{gaussian:gp2:one-sided}
\end{equation}
The latter also bounds a zero-length continuation, because $h_H(x)\ge1$.
The density estimate is never used at zero length. To match the arbitrary
fixed curve coefficient $C$ with the precise Brownian hypotheses, put
$A=\max(1,C)$ and choose
\[
 0<\delta<\min\{\tfrac12-\alpha,A^{-1/(2\alpha)}\},
 \qquad \beta=\tfrac12-\delta>\alpha.
\]
Then $C u^\alpha\le1+\delta^{-1}u^\beta$ for every $u\ge0$:
below $A^{-1/\alpha}$ the first term is at most one, and above it
$\delta^{-1}u^{\beta-\alpha}\ge A$.
Reflecting about $H+1$ therefore gives endpoint distances
$h_H(x),h_H(y)\ge1>\delta$, exactly as required by \cite{CHL}.
Multiplying its bridge probability bound by the free transition density
gives \eqref{gaussian:gp2:killed-density}; its one-sided bound and the trivial
bound for $\ell\le1$ give \eqref{gaussian:gp2:one-sided}.
Only these continuous upper bounds are imported; no decorated-walk
renewal statement is used.

\begin{proposition}[Exact Gaussian three-point envelope]
\label{gaussian:gp2:three-point}
For $q,r$ in \eqref{gaussian:gp2:qr}, uniformly in the three spatial points and real $t\ge k_0+1$,
\begin{equation}
 \E\prod_{i=1}^3F_t^B(z_i)
 \le
 \frac{C_{B,k_0,K,\kappa}\,e^{2q+2r}}
 {(1+q)^{3/2}(1+r-q)^{3/2}
          \sqrt{1+n-q}(1+n-r)}.
 \label{gaussian:gp2:three-envelope}
\end{equation}
\end{proposition}
\begin{proof}
We prove the continuous tree estimate from
Section~\ref{gaussian:gp2:relative}, average over $W$, and restore
\eqref{gaussian:gp2:triple-normalizer}. The three cases distinguish
a positive pair stem, a pair stem of zero length, and the first
branching before a fixed time near the start of the observed history.

\textit{Case 1: $q\ge2k_0+2$ and $r>q$.}
Let $U_q=x$ and $U_r=y$ be the two branch values of the centered tree. Cameron--Martin, on distinct edges only, gives the exact likelihood
\begin{equation}
 \exp\{4U_q-8q+2(U_r-U_q)-2(r-q)\}
       =e^{2x+2y-6q-2r}.
 \label{gaussian:gp2:CM}
\end{equation}
The drifted event has become the centered-tree event under the raised curved barrier.

The delayed start is handled by conditioning at $k_0$, not by imposing constraints earlier. On $[k_0,q]$, its envelope is bounded by a constant depending on $k_0$ plus a curve pinned at $k_0,q$. Absorb the constant into $H$. Integrating the starting Gaussian value $v$ in \eqref{gaussian:gp2:killed-density}, and using $q-k_0\ge q/2$, bounds the common-stem endpoint density by
\[
 \frac{C(1+H)h_H(x)}{(1+q)^{3/2}}\dd x,
 \qquad x\le H.
\]
The singleton branch contributes $C h_H(x)/\sqrt{1+n-q}$. The pair stem has length $r-q\ge1$ and contributes \eqref{gaussian:gp2:killed-density}. The two pair continuations contribute $C h_H(y)^2/(1+n-r)$, including $r=n$.

Their product with \eqref{gaussian:gp2:CM} is bounded by
\begin{align}
 &\frac{C(1+H)e^{-6q-2r}}
 {(1+q)^{3/2}(1+r-q)^{3/2}\sqrt{1+n-q}(1+n-r)}
 \nonumber\\[-1mm]
 &\hspace{6mm}\times
 \int_{-\infty}^H\int_{-\infty}^H
       e^{2x+2y}h_H(x)^3h_H(y)^3\dd x\dd y.
 \label{gaussian:gp2:branch-integral}
\end{align}
The integral equals $\frac{361}{64}e^{4H}$, since
$\int_0^\infty e^{-2u}(1+u)^3\dd u=19/8$. Its height factor is integrable against the uniform Gaussian tail of $W$. Restoring the normalizer changes the exponential factor to
\[
 e^{8q+4r}e^{-6q-2r}=e^{2q+2r},
\]
as required.

\textit{Case 2: $q\ge2k_0+2$ and $r=q$.}
There is no pair-stem transition. The likelihood is $e^{4U_q-8q}$. There are three one-sided continuations, with denominator $(1+n-q)^{3/2}$. The remaining endpoint integral is
\[
 \int_{-\infty}^H e^{4x}h_H(x)^4\dd x
       =\frac{103}{128}e^{4H}.
\]
The normalizer is at most $Ce^{12q}$; after the factor $e^{-8q}$ it gives $e^{4q}$. This is precisely \eqref{gaussian:gp2:three-envelope} with $r=q$, including $q=r=n$. No artificial density on a zero-length edge has been introduced.

\textit{Case 3: $q<2k_0+2$.}
Choose once a fixed integer $a_0=2k_0+3$, so $a_0>q$. If $n\le a_0+1$, all times are bounded and ordinary Gaussian exponential moments give \eqref{gaussian:gp2:three-envelope}, after increasing its constant.

Suppose $n>a_0+1$. Discard constraints before $a_0$ and absorb the bounded envelope value there into $H$. The vector of drifted branch values at $a_0$ is Gaussian with uniformly bounded means and covariances, depending only on $k_0$.

If $r>a_0$, there is a singleton and one pair path at time $a_0$. The singleton continuation has denominator $\sqrt{1+n-a_0}$. On the remaining pair stem use its drift-$2$ likelihood
$e^{2(y-v)-2(r-a_0)}$, where $v$ is the initial pair value. Its killed transition has denominator $(1+r-a_0)^{3/2}$, and the two terminal continuations give $1+n-r$. The endpoint integral in $y$ is $C e^{2H}$ times a fixed polynomial factor. Integrating $e^{-2v}$ and all initial endpoint polynomials against the bounded-time Gaussian vector costs at most $C(1+H)^C$. Restoring the normalizer gives $C e^{2r}$ times a constant depending on the bounded $q,a_0$; this is bounded by the claimed $C e^{2q+2r}$. The denominators with $a_0$ and $q$ are comparable, including $r-a_0=1$.

If $r\le a_0$, all three paths have split by $a_0$. Condition on their values there and use three independent continuations. The resulting probability is at most $C(1+H)^C(1+n)^{-3/2}$; the normalizer is bounded because $q,r$ are bounded. Again this is \eqref{gaussian:gp2:three-envelope}.

All cases have height dependence bounded by $C e^{C H}(1+H)^C$, so averaging over $W$ is legitimate. Finally, for real $t\in[n,n+1)$ the exact triple normalizer and every history tilt differ from the tree expressions at time $n$ by uniformly bounded terms, by \eqref{gaussian:gp2:log-profile}. These were already allowed in the relative comparison. This completes the proof without adding a barrier at $t$.
\end{proof}

\paragraph{Coincident points.}
If two or all three spatial points coincide, keep their formal leaf labels. The covariance comparison remains valid even though the original array is singular. A coincident pair has $r=n$, and three coincident points have $q=r=n$. These cases are included in the proof, not excluded merely because the spatial diagonal has area zero.

\subsection{The spatial sum and the fixed-barrier limit}
\label{gaussian:gp2:moment}

 {Next we shall sum the spatial three-point bounds to obtain a third-moment bound and convergence in mean square at each fixed barrier.}
The exponential factor in \eqref{gaussian:gp2:three-envelope} is
cancelled by the volume of the corresponding spatial shell.
It remains to sum the time denominators uniformly in $n$.
The factor $(1+n)^{3/2}$ in the next lemma accounts for the
three critical mass normalizations.

For a closest-pair labeling, a spatial $(q,r)$ shell has volume at most $C_K e^{-2q-2r}$: fix $z_2$, put $z_3$ in a disk of radius $C e^{-r}$, and $z_1$ in a disk of radius $C e^{-q}$. Shells with $q=0$ or $r=0$ are bounded by the size of $K$; terminal shells collect all smaller distances. There are at most three closest-pair choices.

\begin{lemma}
\label{gaussian:gp2:shell}
For every integer $n\ge0$,
\begin{equation}
 S_n:=(1+n)^{3/2}\sum_{0\le q\le r\le n}
 \frac1{(1+q)^{3/2}(1+r-q)^{3/2}
              \sqrt{1+n-q}(1+n-r)}\le128.
 \label{gaussian:gp2:shell-sum}
\end{equation}
\end{lemma}
\begin{proof}
Write $a=q$, $b=r-q$, $c=n-r$, so $a+b+c=n$. We use
$\sum_{j\ge0}(1+j)^{-3/2}\le3$.

For $a\ge n/2$, the ratio $(1+n)^{3/2}/(1+a)^{3/2}$ is at most $2^{3/2}$, and
\[
 (1+b+c)^{-1/2}(1+c)^{-1}\le(1+c)^{-3/2}.
\]
Summing over $b,c$ bounds this region by $9\,2^{3/2}$.

For $a<n/2$, put $m=n-a$. The convolution estimate
\begin{equation}
 \sum_{b=0}^{m}\frac1{(1+b)^{3/2}(1+m-b)}
       \le\frac{12}{1+m}
 \label{gaussian:gp2:convolution}
\end{equation}
follows by splitting at $b=m/2$. The first half is at most $6/(1+m)$. The other half is at most
$2^{3/2}(1+m)^{-3/2}(1+\log(1+m))\le2^{5/2}/(1+m)$, using
$(1+\log x)/\sqrt x\le2$ for $x\ge1$.
Thus the second region is at most $36\,2^{3/2}$. The total is
$45\,2^{3/2}<128$.
\end{proof}

For continuous-barrier truncations of the martingale and mollified
critical approximations, Lacoin
\cite[Propositions~4.1 and~5.1]{Lacoin} proves $L^2$ convergence.
For the integer-history restriction here, the finite-time bounds in
Proposition~\ref{gaussian:gp2:three-point} and
Lemma~\ref{gaussian:gp2:shell} are independent of any limiting $L^2$
assertion. Combined with the already identified probability limit
in Proposition~\ref{gaussian:gp2:probability}, they give the third-moment
bound and the convergence conclusions of the following theorem.

\begin{theorem}[Fixed-barrier Gaussian limit]
\label{gaussian:gp2:gaussian-main}
For the exact family \eqref{gaussian:gp2:covariance}, any fixed
$B\in\R$, $k_0\ge1$, $\kappa\ge-1$, and bulk disk $K$,
\begin{equation}
 \sup_{t\ge k_0+1}\E[\mu_t^B(K)^3]<\infty.
 \label{gaussian:gp2:third-moment}
\end{equation}
For every $f\in C(K)$ and every $1\le p<3$,
\begin{equation}
 \mu_t^B(f)\longrightarrow
 Z_B(f)=\int_K f(z)\ind_{E_\infty^B}(z)M_\kappa(\dd z)
 \quad\hbox{in }L^p
 \label{gaussian:gp2:Lp}
\end{equation}
as $t\to\infty$ through real values. In particular, $Z_B(K)\in L^3$ and
\begin{equation}
 \lim_{R\to\infty}\sup_{s,t\ge R}
 \E|\mu_t^B(f)-\mu_s^B(f)|^2=0,
 \label{gaussian:gp2:mixed-cauchy}
\end{equation}
where the supremum is over real times. The random measure $Z_B$ is
measurable with respect to the completed sigma-field generated by the
integer histories on $K$.
\end{theorem}
\begin{proof}
We integrate the three-point estimate in
Proposition~\ref{gaussian:gp2:three-point} over $K^3$ and use the shell
volumes and Lemma~\ref{gaussian:gp2:shell}. Since
$t^{3/2}\le C(1+\lfloor t\rfloor)^{3/2}$, this proves
\eqref{gaussian:gp2:third-moment} for all real $t\ge k_0+1$.
Proposition~\ref{gaussian:gp2:probability} gives convergence in
probability to the stated restriction. Applying Fatou's lemma along an
almost surely convergent subsequence gives $\E[Z_B(K)^3]<\infty$.

For $\Delta_t=\mu_t^B(f)-Z_B(f)$, positivity of the measures gives
$\sup_{t\ge k_0+1}\E|\Delta_t|^3\le C_{B,f}$. Thus, for $1\le p<3$,
the variables $|\Delta_t|^p$ are uniformly integrable.
Convergence in probability proves \eqref{gaussian:gp2:Lp} through all
real terminal times. The constants may depend on $B$ and the other
fixed parameters; no $L^3$ convergence is asserted.
The triangle inequality in $L^2$ then gives
\eqref{gaussian:gp2:mixed-cauchy} as both real times tend to infinity.
The random-measure
measurability assertion is proved in Section~\ref{gaussian:gp2:coarse}.
\end{proof}

\paragraph{Dependence on $B$ and removal of the barrier.}
We define all $Z_B$ by the restriction formula in
\eqref{gaussian:gp2:Lp}, on the same realization of the exact Gaussian
family and the canonical measure. Let
\[
 S_K=\sup_{k\ge k_0\text{ integer}}\sup_{z\in K}(G_k(z)-2k).
\]
For a fixed countable dense set $\mathcal D\subset K$, continuity allows
the spatial supremum to be taken over $\mathcal D$.
Proposition~\ref{gaussian:gp2:max-prop} gives
\[
 \Pp(S_K>B)\le C(1+B)e^{-2B}
 \quad\text{for sufficiently large }B.
\]
Here $C$ depends on the fixed $K,k_0,\kappa$. Thus $S_K<\infty$
almost surely. On this event every integer barrier holds throughout
$K$ when $B\ge S_K$, so the common restriction gives
\[
 Z_B=M_\kappa|_K
 \quad\text{for every }B\ge S_K,
 \quad\text{simultaneously, almost surely}.
\]
In particular, for every $\delta>0$ and sufficiently large $B$,
\[
 \Pp\!\left(\norm{Z_B-M_\kappa|_K}_{\mathrm{TV}}>\delta\right)
 \le \Pp(S_K>B)\le C(1+B)e^{-2B}.
\]
The total-variation difference is therefore eventually zero on almost
every sample. The moment bounds above concern each fixed $B$, with
constants allowed to depend on $B$.

\subsection{Measurability and approximation by finitely many observations}
\label{gaussian:gp2:coarse}

The purpose of this subsection is to approximate spatial integrals of the restricted Gaussian measure in mean square by smooth functions of finitely many observations at integer scales.

The matrix identification uses smooth functions of finitely many
Gaussian observations. We explain why they can approximate $Z_B(f)$
while keeping the observation list fixed before the matrix limit.
All sigma-fields in this section are completed by adjoining the null
sets of the fixed ambient Gaussian probability space.

We first use the integer histories. For the countable dense set
$\mathcal D\subset K$ chosen above, we write
\[
 \mathcal H=
 \overline{\sigma\{G_k(v):k\ge k_0\text{ integer},\ v\in\mathcal D\}}.
\]
For each integer $k\ge k_0$, continuity makes $G_k|_K$ an
$\mathcal H$-measurable random element of $C(K)$, whose Borel
sigma-field is generated by evaluations on $\mathcal D$. The
finite-history density defining $\mu_n^B$ is jointly Borel measurable
in these fields and the spatial variable. Thus $\mu_n^B(f)$ is
$\mathcal H$-measurable for every $f\in C(K)$ and integer $n\ge k_0+1$.
Convergence in $L^2$ along integers gives the measurability of $Z_B(f)$.
A countable uniformly dense family in $C(K;\R)$, containing the
constant one, determines finite positive measures and generates their
weak Borel sigma-field. Applying the preceding conclusion to this
family proves that $Z_B$ is an $\mathcal H$-measurable random finite
measure. This proves the measurability assertion in
Theorem~\ref{gaussian:gp2:gaussian-main}.

We also retain the approximation statement for the full positive-time
family. Choose a countable dense set of positive times containing
every positive integer, and enumerate the resulting observations on
$\mathcal D$ as
\[
 C_G=(G_{s_j}(z_j))_{j\ge1},\qquad
 C_G^{(m)}=(G_{s_j}(z_j))_{j\le m},\qquad
 \cG_m=\sigma(C_G^{(m)}).
\]
Any prescribed finite additional coarse observations may be included
in the enumeration. We set
$\cG_\infty=\overline{\sigma(\bigcup_m\cG_m)}$.
Joint continuity gives measurability of the positive-time fields and
their tested finite-history measures with respect to $\cG_\infty$.
In particular, $\mathcal H\subset\cG_\infty$ and $Z_B(f)$ is
$\cG_\infty$-measurable.

The following standard approximation fact uses only a finite Borel law;
Gaussianity is not needed for the approximation argument. We include its
short proof for completeness.

\begin{lemma}[Smooth finite-coordinate approximation]
\label{gaussian:gp2:cylinder}
For every $T\in L^2(\cG_\infty)$ and $\delta>0$, there are a finite
$m$ and $g\in C_b^\infty(\R^m)$, with bounded derivatives of all
orders, such that
\[
 \norm{T-g(C_G^{(m)})}_2<\delta.
\]
No nondegeneracy of the covariance of $C_G^{(m)}$ is required.
For real $T$, the function $g$ may be chosen real.
\end{lemma}
\begin{proof}
The conditional expectations $\E[T\mid\cG_m]$ converge to $T$ in
$L^2$. We choose a sufficiently large fixed $m$ and a Borel
representative of this conditional expectation as a function of
$C_G^{(m)}$. We truncate that function, approximate it in $L^2$ of
the finite Borel law by a compactly supported continuous function,
and then mollify uniformly. These approximations give a smooth
bounded function with bounded derivatives of all orders. For real
$T$ they may all be chosen real. Regularity of a finite Borel measure
on $\R^m$ is sufficient; a density with respect to Lebesgue measure
is not needed.
\end{proof}

For the particular target $T=Z_B(f)$, with $f\in C(K;\R)$, we apply
the same conditional-expectation argument to an enumeration of the
integer observations generating $\mathcal H$. Given $\varepsilon>0$,
we obtain a fixed finite integer observation vector $D_G\in\R^m$ and a real
$\psi\in C_b^\infty(\R^m)$, with all derivatives bounded, such that
\[
 \E|Z_B(f)-\psi(D_G)|^2<\varepsilon.
\]
This proves the consequence stated in \eqref{p1:finite-approximation}.
A prescribed finite coarse vector $V_G$ on the same Gaussian space,
including observations in the enlarged disk $K_*$, may be appended:
the function $\widetilde\psi(x,y)=\psi(x)$ ignores the added
coordinates and retains all the stated bounds. The enlarged list and
the derivative bounds remain fixed before $N\to\infty$. No equality
between $\mathcal H$ and $\cG_\infty$ is required.

The canonical untruncated measure $M_\kappa$ is determined by
$G_\kappa$. The fixed-barrier restriction $Z_B$ also uses the exact
histories $G_k$; the argument above establishes its measurability from
those integer histories. It makes no assertion that $Z_B$ is
determined by the generalized limiting field alone at fixed $B$.

\subsection{A two-point bound beyond the barrier scale}
\label{gaussian:gpm:external}

{Next we shall control pairs of Gaussian weights evaluated beyond the last observed barrier scale, keeping the resulting extra correlation explicit.}

In the matrix annular argument, the exponential weights may be evaluated
at scales later than the last retained barrier scale. We therefore keep
their two evaluation scales separate from the last scale in the barrier
history. This is
the Gaussian estimate used in Lemma~\ref{annuli:lem:gaussian}.
When both exponential weights are evaluated at the last barrier scale,
the bound below has, up to integer rounding, the two-point form proved in dimension two by Lacoin
\cite[Proposition~2.6, (2.21)]{Lacoin}. That estimate uses a continuous
barrier on an independent-increment Gaussian family; here we keep
the integer history of \eqref{gaussian:gp2:covariance} and the two
possibly later evaluation scales.
Retain the covariance \eqref{gaussian:gp2:covariance} and each root's
exact variance normalization. Let $n\ge k_0$ be the final retained
integer, and let the two exponential weights be evaluated at real scales
$\ell_z,\ell_w\ge n$. Put
\[
 q_\infty=\log_+|z-w|^{-1},\qquad
 q=\min(\ell_z,\ell_w,q_\infty),\qquad
 m=\min(n,\lfloor q\rfloor),
\]
with $q_\infty=\infty$ at coincidence. Here $q$ records the
correlation depth visible to the two roots, whereas $m$ is the
integer part retained by the barrier. Their difference must remain
in the general estimate; its effect is absorbed into the constant
under the stated annular hypotheses.

\begin{theorem}[Two-point bound with two evaluation scales]
\label{gaussian:gpm:external-theorem}
Uniformly in these evaluation scales and spatial points,
\begin{equation}
 \E[W_{\ell_z}(z)W_{\ell_w}(w)
          \ind_{A_n^B(z)\cap A_n^B(w)}]
 \le\frac{C_B e^{4q-2m}}{(1+m)^{3/2}(1+n-m)}.
 \label{gaussian:gpm:external-bound}
\end{equation}
If $j$ is an integer with $0\le j\le n$, $|q_\infty-j|\le C_0$,
and $q_\infty\le n+C_0$, then
\begin{equation}
 e^{-2j}\E[W_{\ell_z}(z)W_{\ell_w}(w)
          \ind_{A_n^B(z)\cap A_n^B(w)}]
 \le\frac{C_{B,C_0}}{(1+j)^{3/2}(1+n-j)}.
 \label{gaussian:gpm:external-annular}
\end{equation}
There is no dependence on the unobserved gaps $\ell_v-n$ in this annular
bound. At coincidence \eqref{gaussian:gpm:external-bound}, rather than
\eqref{gaussian:gpm:external-annular}, is the assertion.
\end{theorem}

\begin{proof}
Double-root tilting leaves the centered history covariance unchanged.
For $k\le n$, the mean at either spatial site is
\[
 2C_\kappa(k,z;\ell_z,z)+2C_\kappa(k,z;\ell_w,w)
       =2k+2\min(k,q_\infty)+O(1),
\]
with the analogous formula at $w$. Since
$\min(k,q_\infty)\ge\min(k,m)$, the centered barrier is contained in
$X_k^i\le B+C-2\min(k,m)$. The exact scalar factor is
\begin{equation}
 \E[W_{\ell_z}(z)W_{\ell_w}(w)\ind_A]
 =e^{4C_\kappa(\ell_z,z;\ell_w,w)}Q_{\rm dbl}(A)
 \le Ce^{4q}Q_{\rm dbl}(A).
 \label{gaussian:gpm:double-tilt}
\end{equation}

Compare the centered history to two Brownian paths sharing their stem up
to time $m$ and independent afterward. Their covariance differs by a
fixed constant, including formal duplicate labels. The proof of
Lemma~\ref{gaussian:gp2:tree-comparison} applies verbatim with two leaves and pin
set $\{0,m,n\}$. For clarity, add independent $N(0,3M)$ noises to the
exact array. Compare it with the tree plus a common $N(0,2M)$ and
independent noises of variance $\Var X_k^i-k+M\in[0,2M]$. Diagonal
variances agree and off-diagonal covariances are ordered. Covariance
interpolation for decreasing rectangle factors gives the correct order.
The summable allowances
$A\sqrt{\log(2+\dist(k,\{0,m,n\}))}$ have a uniform positive product
probability, independent of the exact array. Dividing by that probability
and absorbing the remaining noises into an independent height $H\ge1$
with a uniform Gaussian tail yields a multiplicative comparison.
There is no additive error competing with the common-stem probability.

Insert the independent unit Brownian bridges once per distinct edge.
Their pinned allowances again have a uniform positive product
probability. As in the three-leaf proof, the integer event is bounded by
a constant times a continuous event under a curve
$H+C\dist(s,\{0,m,n\})^\alpha$, $0<\alpha<1/2$. The deterministic drift
is $2$ on the common stem and zero afterward.

Suppose first $m\ge2k_0+2$. For a centered tree $U$ the drift likelihood is
$\exp(2U_m-2m)$. Write $x=U_m\le H$ and $h_H(x)=1+H-x$.
Conditioning at the actual barrier start $k_0$ and using
\eqref{gaussian:gp2:killed-density}, the common-stem endpoint density is at most
\[
 \frac{C(1+H)h_H(x)}{(1+m)^{3/2}}\dd x.
\]
This follows by integrating the starting endpoint polynomial against a
fixed-variance Gaussian; it imposes no barrier before $k_0$.
The two independent continuations each contribute at most
$Ch_H(x)/\sqrt{1+n-m}$, also when their length is zero. Thus the drifted
survival is bounded by
\begin{align}
 \frac{C(1+H)e^{-2m}}{(1+m)^{3/2}(1+n-m)}
       \int_{-\infty}^H e^{2x}h_H(x)^3\dd x
 =\frac{19C(1+H)e^{2H-2m}}{8(1+m)^{3/2}(1+n-m)}.
 \label{gaussian:gpm:two-tree}
\end{align}
No transition density is assigned to a zero-length edge.

If $m<2k_0+2$, condition instead at a fixed integer
$a\ge\max(m,k_0)$ and discard earlier constraints. The bounded-time
Gaussian vector $(U_m,U_a^1,U_a^2)$ may be singular; its polynomial and
exponential moments remain uniformly finite. The two future survival
factors give denominator $1+n-a$, comparable to $1+n-m$ when $n$ is
large. If $n$ is bounded, ordinary Gaussian moments suffice. This proves
\eqref{gaussian:gpm:two-tree} after changing constants. It includes $m=0$, delayed
barrier starts and $m=n$.

Average over $H$, whose Gaussian tail controls all height factors, and
restore \eqref{gaussian:gpm:double-tilt}. This proves
\eqref{gaussian:gpm:external-bound}. In the stated annular regime, $q=j+O(1)$ and
$m=j+O(1)$, with truncation at $n$ only on the last bounded number of
annuli. This proves \eqref{gaussian:gpm:external-annular}.
\end{proof}

\begin{remark}[Why the terminal core is separate]
At coincidence with $q>n$, compensation at $j=n$ leaves the factor
$e^{4(q-n)}$ in \eqref{gaussian:gpm:external-bound}. The annular bound cannot be
extended into this region by quoting a same-time full-barrier theorem.
The matrix application controls that core by its separate endpoint cutoff.
\end{remark}

\subsubsection{Exact signed compensation}

 {Now we introduce a negative coefficient at an observed history coordinate and track the exact cancellation that restores the two-point weight.}
The signed Gaussian comparison uses a negative coefficient at an
observed history coordinate in its positive exponential weight. Multiplication by the corresponding
positive exponential restores the two-root weight exactly. This
identity keeps the root self-variances until they cancel.
Throughout this subsection impose the annular hypotheses of
\eqref{gaussian:gpm:external-annular}. In particular, coincidence and a terminal
core with $q_\infty>n+C_0$ are outside the present assertion.
Form a finite Gaussian observation vector $G$ containing the two
external roots, the retained history and the coarse observations.
Let $\Sigma=\operatorname{Cov}(G)$, and let $e_r$ denote its coordinate
unit vectors, with distinct labels retained for repeated observations.
Let $c=(j,z)$ be an observed history coordinate, with $k_0\le j\le n$,
and let $r_z,r_w$ be the external-root labels. Put
\[
 p=2e_{r_z},\quad r=2e_{r_w},\quad a=p+r-2e_c,\quad
 Z_G(a)=e^{a^T\Sigma a/2},\quad
 \mathcal C_{G,j}=e^{2j}Z_G(a)/(Z_G(p)Z_G(r)).
\]
Let $G_a$ denote the Gaussian vector under the normalized
exponential tilt by $a$, with mean $\Sigma a$ and covariance
$\Sigma$. Let $\mathcal B$ be a smooth decreasing barrier product
with $0\le\mathcal B\le1$, supported in the hard barrier of
height $B+1$, and let $\Psi$ be
a measurable coarse factor with $0\le\Psi\le1$. Put
$g(x)=e^{2x_c-4j}\mathcal B(x)\Psi(x)$. Since $c$ is observed by
the barrier, $g\le e^{2(B+1)}$. Exact cancellation gives
\begin{equation}
 \mathcal C_{G,j}\E g(G_a)
 =e^{-2j}\E[W_{\ell_z}(z)W_{\ell_w}(w)\mathcal B(G)\Psi(G)]
 \le\frac{C_B}{(1+j)^{3/2}(1+n-j)}.
 \label{gaussian:gpm:compensation}
\end{equation}
The logarithm of the scalar factor is exactly
\[
 2j+4\Sigma_{r_zr_w}-4\Sigma_{r_zc}-4\Sigma_{r_wc}+2\Sigma_{cc}.
\]
Its four covariances equal $j+O(1)$ in the annular geometry, so
$\mathcal C_{G,j}$ is bounded above and below. No covariance inverse or
independence of the coarse factor is used. The upper bound discards
$\Psi$ using $0\le\Psi\le1$. Smooth matrix comparisons impose
their derivative hypotheses separately; in the annular application
those costs are recorded in \eqref{annuli:eq:g-derivatives}.

\subsection{Endpoint strips and terminal completion}
\label{gaussian:gpm:endpoint}

The purpose of this subsection is to bound endpoints near the barrier so that omitted terminal scales can be restored with errors that vanish after critical normalization.

The one-root comparison needs to complete a short omitted terminal
history with error $o(t^{-1/2})$. The separated two-root comparison
needs error $o(t^{-1})$, and its endpoint-cut removal retains both
barriers. We prove the corresponding one-path and pair estimates
separately. They are used in Theorem~\ref{oneroot:thm:one-root}
and Corollary~\ref{separated:cor:sm-cut}; their different powers
of $t$ arise from normalizing a single critical mass in the one-point
estimate and a product of two critical masses in the two-point estimate.

For the marked statements fix a positive separation $r>0$ and
$|z-w|\ge r$. Constants may depend on $r$, $B$, $k_0$, $K$ and $\kappa$.
The observation horizon $t$ is real and barriers end at
$n=\lfloor t\rfloor$. The two exponential weights may be evaluated
at different scales $\ell_z,\ell_w\ge t$, at or after the terminal scale.
Under the normalized double-root law
\[
 \frac{\dd Q}{\dd\PP}
 =\exp\big(2G_{\ell_z}(z)+2G_{\ell_w}(w)
       -2v_{\ell_z}(z)-2v_{\ell_w}(w)
       -4C_\kappa(\ell_z,z;\ell_w,w)\big),
\]
write $Y_s^a=G_s(a)-2s$, for $a=z,w$. Its mean is $O_r(1)$ for
$s\le t$, uniformly in both external gaps, since the own-root covariance
is $s+O(1)$ and the cross-root covariance is bounded at fixed separation.
Its centered covariance is the original covariance.
Under one exponential weight evaluated at scale $\ell\ge t$, the same conclusions hold
with the cross-root term omitted.
The product of the two exponential weights, each divided by its own
expectation, differs from $Q$ by the uniformly bounded factor
$e^{4C_\kappa(\ell_z,z;\ell_w,w)}$.
All strip, completion and factorization assertions below are uniform in
these external evaluation scales; observations remain at times at most $t$.

\subsubsection{Covariances of later increments and earlier observations}

Now our goal is to quantify how covariances between late increments and earlier observations decay, providing the error bound for separating their laws.
To separate later increments from earlier observations, a covariance
comparison with an $O(1)$ error is insufficient. The exact formula
gives for $s\ge k$
\[
 |\partial_s C_\kappa(s,z;k,z)|\le Ce^{-2(s-k)}.
\]
Indeed $|\partial_s\tanh s|\le Ce^{-2s}$ and the same-point determinant
is bounded below by a constant times $e^{-2k}$. Therefore
\begin{equation}
 |C_\kappa(a+j,z;k,z)-C_\kappa(a,z;k,z)|
       \le Ce^{-2(a-k)},\qquad k\le a,\ j\ge0.
 \label{gaussian:gpm:prefix-decay}
\end{equation}
For a fixed positive-cutoff observation the right side is $Ce^{-2a}$;
for an observation at the other spatial site it is $C_re^{-2a}$ uniformly
in its scale. Taking a second increment shows that two late future blocks
at the separated sites have cross-covariance $O_r(e^{-2a})$.

Every retained early coordinate has uniformly positive variance, even at
$\kappa=-1$: since $D(u,z;u,z)\le1-u^2$ and $p_z(u)^2\le u^2$,
\begin{equation}
 v_s(z)\ge\tfrac12[-\log(1-u_s^2)-u_s^2]\ge u_s^4/4.
 \label{gaussian:gpm:positive-variance}
\end{equation}
For large $s$, $v_s=s+O(1)$. These facts justify all unconditional density
bounds below, without uniform conditional nondegeneracy.

\subsubsection{An additive rectangle comparison for polynomial precision}

Next we shall compare Gaussian barrier probabilities by shifting thresholds, with errors smaller than any prescribed inverse power of the scale.
At fixed separation the two-root normalizer is uniformly bounded.
A prescribed polynomial error therefore suffices after the later
critical normalizations. Fix the desired error power $D>0$ before
choosing the Gaussian noise allowance. Let centered Gaussian vectors $U,V$ of dimension at most $Ct$ have
covariances differing entrywise by at most $M$. Add independent variance
$2M$ noises to $U$; add to $V$ a common variance $M$ Gaussian and
independent noises of variances $M+\Var U_i-\Var V_i\in[0,2M]$.
Diagonal variances agree and the latter off-diagonal covariances dominate.
Gaussian rectangle interpolation, followed by a Gaussian union bound on
the added noises, gives for arbitrary thresholds
\begin{equation}
 \PP(U_i\le b_i\ \forall i)
 \le\PP(V_i\le b_i+C_D\sqrt{\log(2+t)}\ \forall i)+C_Dt^{-D}.
 \label{gaussian:gpm:additive-comparison}
\end{equation}
Interchanging $U,V$ gives the corresponding lower comparison. Singular arrays are handled
by independent regularization and a limit. The bounded covariance
difference itself is not asserted positive semidefinite.

Integer Brownian constraints can likewise be replaced by continuous
constraints with an extra $O(\sqrt{\log t})$ lift and $O(t^{-D})$ error.
Condition on all sampled values; each omitted interval is a Brownian
bridge of length at most $k_0$. Its overshoot above the larger endpoint has
tail at most $e^{-2x^2/k_0}$. Include the deterministic initial endpoint and
specified terminal endpoint in this interpolation. This handles the
unobserved initial interval without imposing an earlier constraint on the
original process. The present additive comparison is used only at
polynomial survival precision; the annular and three-leaf estimates used
the multiplicative comparison instead.

\begin{lemma}[One-root and separated pair endpoint strips]
\label{gaussian:gpm:strips}
Let $A_j^a=\{Y_k^a\le B:k_0\le k\le j\}$ and
$\mathcal A_j=A_j^z\cap A_j^w$. For integers $k_0\le j\le n$
and $1\le R\le j$, under one root,
\begin{equation}
 Q_{\rm one}(A_j,0\le B-Y_j\le R)
 \le\frac{C_DR(1+\sqrt{\log(2+j)})(1+R+\sqrt{\log(2+j)})}{j^{3/2}}
          +C_Dj^{-D}.
 \label{gaussian:gpm:one-strip}
\end{equation}
The same one-path bound holds for either marginal path under the
fixed-separated double-root law, with constants allowed to depend on
the separation. The regression proof uses only the bounded rooted mean
and unchanged one-path covariance, which that marginal also has.
Under the separated double-root law,
\begin{align}
 Q(\mathcal A_j)&\le C_D\frac{1+\log(2+j)}j+C_Dj^{-D},
                                               \label{gaussian:gpm:pair-survival}\\
 Q(\mathcal A_j,0\le B-Y_j^z\le R)
 &\le\frac{C_DR(1+\log(2+j))(R+1+\sqrt{\log(2+j)})}{j^2}
              +C_Dj^{-D}.\label{gaussian:gpm:pair-strip}
\end{align}
The other pair endpoint is symmetric. An arbitrary translated endpoint
strip at a time $J\ge1$ also satisfies
\begin{equation}
 Q(|Y_J-c|\le h)\le Ch/\sqrt J,\qquad c\in\R,
 \label{gaussian:gpm:translate}
\end{equation}
uniformly in its center and the rooted mean.
\end{lemma}
\begin{proof}
Condition in the one-root case on $Y_j=y\in[B-R,B]$. Gaussian regression
gives
\[
 \E(Y_k\mid Y_j=y)=(k/j)y+O_B(1),\qquad
 \Cov(Y_k,Y_l\mid Y_j)=\min(k,l)-kl/j+O(1).
\]
For the mean, $C(k,j)/v_j=k/j+O(1/j)$ and $|y|\le|B|+j$.
Expanding $(k+O(1))(l+O(1))/(j+O(1))$ proves the covariance bound.
Apply \eqref{gaussian:gpm:additive-comparison} to the residual and interpolate the
Brownian bridge. It suffices to bound a Brownian bridge from $0$ to $y$
below a height $H_j$, where
$\max(0,B)+1\le H_j\le C_D(1+\sqrt{\log(2+j)})$.
Reflection gives survival at most $2H_j(H_j-y)/j$.
The endpoint density is at most $Cj^{-1/2}$; integration over the
length-$R$ interval proves \eqref{gaussian:gpm:one-strip}.

The pair survival follows by comparison to two independent Brownian
motions and their one-sided reflection probabilities. For the pair strip,
condition on $Y_j^z=y$. The first path has mean $(k/j)y+O(1)$ and the
second has bounded mean. Its conditional covariance differs by $O_r(1)$
from an independent Brownian bridge and Brownian motion: all cross-site
covariances are bounded, and the rank-one conditioning correction remains
bounded. Reflection bounds the conditional survival by
\[
 \bigl(1-e^{-2H_j(H_j-y)/j}\bigr)
     \bigl(2\Phi_{\mathrm{nor}}(H_j/\sqrt j)-1\bigr)
 \le CH_j^2(H_j+R+|B|)/j^{3/2}.
\]
Integrate against the endpoint density to get \eqref{gaussian:gpm:pair-strip}.
Finally \eqref{gaussian:gpm:translate} is the unconditional Gaussian density bound.
It concerns moving endpoint thresholds; the killed-strip bounds continue
to keep $B$ fixed.
\end{proof}

\paragraph{Completing the terminal history.}
Take an integer $J\le n$ with $t-J=O(\log t)$. On the retained
barrier, a violation at a later integer time requires either a deficit
$O(\log t)$ at $J$ or a later increment exceeding that amount.
Bounded rooted means and increment variance at most $t-J+C$ make the
latter probability $O(t^{-D})$ with a sufficiently large constant.
Lemma~\ref{gaussian:gpm:strips} gives one-root completion
$O((\log t)^{5/2}t^{-3/2})+O(t^{-D})$ and pair completion
$O((\log t)^3t^{-2})+O(t^{-D})$. The real-time final increment is handled
the same way. These are respectively $o(t^{-1/2})$ and $o(t^{-1})$.
Multiplying by any uniformly bounded measurable insertion preserves these
bounds. A bounded Lipschitz endpoint mark changes by
$O(\log t/\sqrt t)$ off the exceptional event. Its contribution
to a critically normalized pair expression also includes the
exceptional event, where the marks are bounded. The total error is at most
\begin{align*}
 &t\,O(\log t/\sqrt t)\,Q(\mathcal A_J)+O(t^{1-D})\\
 &\qquad=O((\log t)^2/\sqrt t)+O(t^{1/2-D}\log t)+O(t^{1-D})=o(1)
\end{align*}
when $D>1$. Thus this mark replacement uses pair survival as well
as the size of the mark displacement. All these estimates concern
the stated rooted Gaussian laws; no Markov property of the exact
scale family is used.

\subsection{Killed future laws and separated Rayleigh marks}
\label{gaussian:gpm:marks}

In this subsection we separate late increments from earlier observations and show that endpoint tests contribute Rayleigh averages at fixed positive spatial separation.

The Rayleigh law is the endpoint law of the standard Brownian meander.
Its occurrence in marked critical GMC is established in
\cite[Theorem~6.7]{APS}, where the endpoint mark has an additional
factor $\sqrt2$, and in \cite[Lemma~2.7]{PowellReview} for star-scale
invariant fields. For the exact family here, we first compare a killed
Gaussian future with killed Brownian motion, then separate the late
futures from the early history, and finally identify the endpoint law.
The comparison is additive: the proof retains the prefix survival
probability before applying the critical normalization. It therefore
does not require the leading constant in the survival probability of
the exact Gaussian path.

We write
\[
 \rhoR(h)=\int_0^\infty h(v)v e^{-v^2/2}\dd v.
\]
This is the standard Rayleigh probability law. Marks below are bounded
Lipschitz functions on $\R$; the limiting law is supported on $[0,\infty)$.
We test the early observations using a bounded function $H$ of class
$C^2$, with bounded derivatives, applied to a fixed finite list $U$ of
observations at positive cutoffs. Smooth bounded test functions are
included. The sites of those observations may
coincide with either root, since their heights stay fixed and positive.
The conventions of the preceding section for fixed separation, a real
terminal time and external evaluation scales remain in force. The endpoint test functions and the test of the early observations may
take either sign; the estimates below control absolute errors.

\begin{lemma}[Uniform killed endpoint comparison]
\label{gaussian:gpm:killed-cdf}
Let $(V_j)_{1\le j\le N}$ be a centered Gaussian vector, with
$N\asymp t$ and covariance $\min(j,l)+O(1)$, and let $W$ be
standard Brownian motion starting at zero. Set
\[
 F_V(d,y)=\PP(\max_{j\le N}V_j\le d,V_N\le y),\qquad
 F_W(d,y)=\PP(\sup_{0\le s\le N}W_s\le d,W_N\le y).
\]
For $d\ge C_D\sqrt{\log t}$,
\begin{equation}
 \sup_{y\in\R}|F_V(d,y)-F_W(d,y)|
 \le C_D\frac{\sqrt{\log t}}{\sqrt t}+C_Dt^{-D}.
 \label{gaussian:gpm:cdf}
\end{equation}
For each fixed bounded Lipschitz function $h:\R\to\R$ and
an auxiliary deterministic scale $m=m(t)>0$ satisfying
$\sqrt m\log t=o(\sqrt t)$, the expectations of
$h((d-V_N)/\sqrt t)$ and $h((d-W_N)/\sqrt t)$, restricted to
their respective barrier events, differ in absolute value by at most
$C_h\log t/\sqrt t+O(t^{-D})$, uniformly for
$C_D\sqrt{\log t}\le d\le\sqrt m\log t$.
\end{lemma}
\begin{proof}
Use \eqref{gaussian:gpm:additive-comparison} in both directions, with final
threshold $\min(d,y)$, and interpolate the Brownian bridges.
For $S=C_D\sqrt{\log t}$, this brackets $F_V(d,y)$ by
$F_W(d\pm CS,y\pm CS)$ up to $O(t^{-D})$. The change in $d$ is at most
\[
 \PP(\sup W\in[d-CS,d+CS])\le CS/\sqrt N,
\]
by $\PP(\sup W\le b)=2\Phi_{\mathrm{nor}}(b/\sqrt N)-1$ for $b\ge0$.
The change in $y$ is at most the unconditional Gaussian endpoint mass
in an interval of length $2CS$, again $CS/\sqrt N$. This proves
\eqref{gaussian:gpm:cdf} uniformly for every real $y$.

For a bounded Lipschitz function of $(d-V_N)/\sqrt t$, truncate that
argument at $C_D\sqrt{\log t}$. Both Gaussian future and Brownian
unconditional endpoint tails outside the truncation are $O(t^{-D})$,
uniformly for the stated $d=o(\sqrt t)$. The truncated function has total
variation $O_h(\sqrt{\log t})$. Integration by parts against
\eqref{gaussian:gpm:cdf} proves the asserted marked error. Since these are
subprobability laws, include the endpoint term involving their total-mass
difference; \eqref{gaussian:gpm:cdf} also bounds it, by taking $y\to\infty$.
A bounded Lipschitz
function need not have finite variation on the full line, which is why
this truncation is needed.
\end{proof}

The next theorem compares a marked expectation directly with its
unmarked counterpart. This formulation also applies to a centered mark
$h-\rhoR(h)$, which remains bounded and Lipschitz even when it loses
compact support. The resulting zero Rayleigh mean is what the spatial
second-moment argument will use.

\begin{theorem}[Factorization of endpoint tests at separated points]
\label{gaussian:gpm:factorization}
At fixed positive separation, uniformly on the fixed spatial compact,
with $H$ and its finite observation list fixed as $t\to\infty$,
\begin{align}
 t\big|\E_Q[H(U)\ind_{\mathcal A_n}
 h_1(-Y_t^z/\sqrt t)h_2(-Y_t^w/\sqrt t)]
 -\rhoR(h_1)\rhoR(h_2)\E_Q[H(U)\ind_{\mathcal A_n}]\big|
 \longrightarrow0.
 \label{gaussian:gpm:pair-marks}
\end{align}
The same assertion holds when the two exponential weights are each
divided by their own expectation, without normalizing their product. There is also the one-root
version with prefactor $\sqrt t$ and one Rayleigh factor.
No ratio with a possibly zero or signed expectation is used.
\end{theorem}

\begin{proof}
We first prove the pair assertion. The argument truncates the prefix,
separates the futures, and identifies their killed endpoint laws,
retaining the probability of the surviving prefix at each step.

\emph{Prefix truncation and future separation.}
For a fixed constant $c_{\mathrm{gap}}>0$ to be chosen below, set
\[
 m=\lfloor t^{3/4}\rfloor,\quad
 g=\lceil c_{\mathrm{gap}}\log t\rceil,\quad
 a=m+g,\quad N=\lfloor t\rfloor-a,
\]
\[
 d_i=B-Y_m^{v_i},\quad(v_1,v_2)=(z,w),\qquad
 d_-=t^{1/16},\qquad d_+=\sqrt m\log t.
\]
Restrict the prefix to
\[
 \mathcal C=\mathcal A_m\cap\{d_-\le d_1,d_2\le d_+\}.
\]
The small-deficit part discarded has, by \eqref{gaussian:gpm:pair-strip}, probability
$O(t^{-11/8}\log t)+O(m^{-D})=o(t^{-1})$. The large-deficit part has
probability $O(e^{-c(\log t)^2})$, by unconditional Gaussian tails.
Keep the prefix-survival estimate
\begin{equation}
 Q(\mathcal A_m)\le C\log(2+m)/m.
 \label{gaussian:gpm:prefix-survival}
\end{equation}

Center $Z_s^v=Y_s^v-\E_QY_s^v$ and define future increments
$V_j^i=Z_{a+j}^{v_i}-Z_a^{v_i}$. Equation
\eqref{gaussian:gpm:prefix-decay} shows that every covariance between a vector of later increments and the
entire earlier history or the coarse observations is $O(e^{-2g})$; the
cross-covariances between the two vectors of later increments have the
same bound. Each future marginal has covariance $\min(j,l)+O(1)$.

We define a product law that retains the joint law of the initial
history and coarse observations, and each future marginal, while
making the two futures independent of each other and of those earlier
observations. Convex covariance interpolation stays
positive semidefinite. For a deterministic shift $x$, smooth the prefix
barriers, deficit restrictions and mixed future barriers
$Y_m^{v_i}+V_j^i\le B+x$ at width
$\delta=t^{-P_0}$, with $P_0>3$. Mollify Lipschitz marks at the same
width if needed. For fixed $H$, the aggregate second derivative bound
is at most $Ct^2\delta^{-2}$. Thus interpolation costs
$Ct^2\delta^{-2}e^{-2g}=o(t^{-1})$ when $c_{\mathrm{gap}}>P_0+4$.

At both endpoint laws, the prefix coordinates have positive variances
bounded below by \eqref{gaussian:gpm:positive-variance}. The mixed barrier forms
$Y_m+V_j$ have variance $m+j+O(1)$, and deficit boundaries have variance
$m+O(1)$. Their unconditional density bounds show that removing the smoothing
introduces an error of at most
$Ct\delta=o(t^{-1})$. Redundant coordinates and singular full covariance
matrices cause no difficulty, since no conditional covariance is inverted.

\emph{Gap restoration and endpoint shifts.}
Set $q_t=C(1+\log t)^2$. Under the actual law, a finite Gaussian
union bound gives a super-polynomially small exception outside which,
at each site, $|Y_k-Y_m|\le q_t$ at every integer $m\le k\le a$,
and $|Y_t-Y_{\lfloor t\rfloor}|\le q_t$.
The first increments have variance at most $g+C$ and the final
increment has bounded variance; their rooted means are bounded.
On $\mathcal C$, $d_-\gg q_t$, so all gap barriers are automatic.
The future barriers are bracketed by
\[
 K_i^- =\{\max_j V_j^i\le d_i-Cq_t\},\qquad
 K_i^+ =\{\max_j V_j^i\le d_i+Cq_t\}.
\]
The endpoint marks differ from $(d_i-V_N^i)/\sqrt t$ by at most
$Cq_t/\sqrt t$. Apply this deterministic bracketing under the actual law
before comparing its bounding observables to the product law. The gap
event itself is not transported and is never asserted independent.

Under the product law, conditional on the prefix, a future survives with
probability at most $Cd_+/\sqrt t$. Its shell of width $Cq_t$ has
probability at most $Cq_t/\sqrt t$, by Lemma~\ref{gaussian:gpm:killed-cdf} and
reflection. A two-path shell thus has probability at most $Cq_td_+/t$.
After prefix integration and critical multiplication, its cost is
\begin{equation}
 t\,Q(\hbox{prefix and pair shell})
 \le C\frac{\log(2+m)}m q_td_+=o(1).
 \label{gaussian:gpm:shell-budget}
\end{equation}
The mark-shift error similarly costs at most
\begin{equation}
 C\frac{\log(2+m)}m d_+^2\frac{q_t}{\sqrt t}
 =C\log(2+m)(\log t)^2\frac{q_t}{\sqrt t}=o(1).
 \label{gaussian:gpm:mark-shift}
\end{equation}
Bounded signed $H$ and signed marks are handled by absolute values in
these estimates.

\emph{Identification of the independent future marks.}
It remains to identify each future mark under the product law. Write
\[
 g_i(d)=\PP(\max_jV_j^i\le d),\qquad
 f_i(d)=\E[\ind_{\max_jV_j^i\le d}
                   h_i((d-V_N^i)/\sqrt t)].
\]
For Brownian motion killed at $d>0$, the endpoint deficit
$u=d-W_N>0$ has subprobability density
\begin{equation}
 p_{N,d}(u)=\frac1{\sqrt{2\pi N}}
  \left[e^{-(u-d)^2/(2N)}-e^{-(u+d)^2/(2N)}\right].
 \label{gaussian:gpm:reflection-density}
\end{equation}
Write $\phi(x)=(2\pi)^{-1/2}e^{-x^2/2}$ for the standard normal
density. For $\varepsilon=d/\sqrt N\le1/2$, Taylor's formula and the
integrable third derivative of $\phi$ give the uniform $L^1$ expansion
\[
 \phi(v-\varepsilon)-\phi(v+\varepsilon)
      =2\varepsilon v\phi(v)+O_{L^1}(\varepsilon^3).
\]
Replacing $h_i(\sqrt{N/t}\,v)$ by $h_i(v)$ costs
$C_h(t-N)/t$ after integration against $v\phi(v)$.
Lemma~\ref{gaussian:gpm:killed-cdf} transfers this calculation to the actual
future marginal, yielding
\begin{equation}
 |f_i(d)-\rhoR(h_i)g_i(d)|
 \le C_h\left[\frac{q_t}{\sqrt t}
          +\frac{d^3}{t^{3/2}}+\frac{dm}{t^{3/2}}\right]+C_Dt^{-D}.
 \label{gaussian:gpm:one-future-mark}
\end{equation}
The same leading Brownian coefficient
$\sqrt{2/\pi}\,d/\sqrt N$ occurs in marked and unmarked expectations;
it cancels exactly. No persistence constant for the exact Gaussian path
has been assumed or computed.

Since $|f_i|\le\|h_i\|_\infty g_i$ and $g_i\le Cd_+/\sqrt t$, the
product-factorization error is at most
\[
 C_h\left[\frac{q_td_+}{t}+\frac{d_+^4}{t^2}
                         +\frac{d_+^2m}{t^2}\right]+C_Dt^{-D}.
\]
We integrate against the joint law of the initial history and coarse
observations, with the factor $\ind_{\mathcal C}H(U)$, and multiply by $t$. The factor
\eqref{gaussian:gpm:prefix-survival} is retained before this multiplication:
\begin{align}
 &C_h\|H\|_\infty\frac{\log(2+m)}m
       \left[q_td_++\frac{d_+^4}t+\frac{d_+^2m}t\right]\notag\\
 &\quad\le C_h\|H\|_\infty\bigg[
 \frac{q_t\log t\log(2+m)}{\sqrt m}
 +\frac mt(\log t)^4\log(2+m)\notag\\
 &\hspace{43mm}+\frac mt(\log t)^2\log(2+m)\bigg]
 \longrightarrow0.
 \label{gaussian:gpm:product-budget}
\end{align}
Every term vanishes for $m=t^{3/4}$. Without the prefix survival factor,
the term $m^2/t$ would grow. Combining the discarded prefixes,
interpolation, gap, shell, mark shift and \eqref{gaussian:gpm:product-budget}
proves \eqref{gaussian:gpm:pair-marks}. The bounded scalar in front of $Q$ gives
the unnormalized-root version.

\emph{The one-root assertion.}
For one root, choose instead $m=t^{1/2}$ with the same logarithmic gap,
$d_-=t^{1/16}$ and $d_+=\sqrt m\log t$. The small-deficit discard from
\eqref{gaussian:gpm:one-strip} is $O(t^{-5/8}\sqrt{\log t})=o(t^{-1/2})$.
The sufficient one-root prefix survival bound
$O((\log m)^{3/2}/\sqrt m)$ follows from the strip estimate and a Gaussian
endpoint tail, splitting the deficit at $C_D\sqrt{m\log m}$.
Apply the same separation of the future from the earlier observations
and the killed-endpoint comparison, now with one future block.
After multiplication by $\sqrt t$, the reflection residual is
bounded by
\[
 q_tQ(A_m)+d_+^3/t+d_+m/t=o(1).
\]
The corresponding shell and mark-shift estimates are smaller. This proves
the one-root assertion, with the same truncation for bounded Lipschitz
marks as in Lemma~\ref{gaussian:gpm:killed-cdf}.
\end{proof}

\subsection{Fixed-barrier marked spatial measures}
\label{gaussian:gpm:marked-spatial}

In this subsection we combine separated mark estimates with control of nearby points to obtain the Gaussian marked measure at a fixed barrier.
We first pass from separated endpoint tests to a spatial measure on the
same Gaussian space. The factorization theorem controls the centered
marks at fixed positive separation; the two-root envelope controls the
integral near the diagonal. These two estimates are used in this order,
with the barrier $B$ fixed throughout.

For bounded Lipschitz $h$ and $f\in C(K)$ define
\[
 \mu_t^{B,h}(f)=\sqrt t\int_K f(z)W_t(z)\ind_{A_n^B(z)}
                    h((2t-G_t(z))/\sqrt t)\dd A(z).
\]
These are the product tests of the nonnegative finite measure obtained
by attaching the mark $(2t-G_t(z))/\sqrt t$ to $\mu_t^B$ at each $z$.
The tested quantity may be signed when $f$ or $h$ is signed. Its centered
version is $\mu_t^{B,h}(f)-\rhoR(h)\mu_t^B(f)$.

\begin{corollary}[Rayleigh product at fixed barrier]
\label{gaussian:gpm:marked-L2}
For every fixed $B$, every $f\in C(K)$ and bounded Lipschitz $h$,
\begin{equation}
 \mu_t^{B,h}(f)-\rhoR(h)\mu_t^B(f)\longrightarrow0
       \quad\hbox{in }L^2,
 \qquad
 \mu_t^{B,h}(f)\longrightarrow\rhoR(h)Z_B(f)
       \quad\hbox{in }L^2.
 \label{gaussian:gpm:marked-L2-bound}
\end{equation}
Equivalently the marked finite measures converge in probability, weakly
on $K\times\R$, to $Z_B\otimes\rhoR$. The convergence is joint with
any finite list of exact coarse observations on this same Gaussian space.
\end{corollary}
\begin{proof}
\emph{Centered marks in $L^2$.}
Expand the square of the first difference. On $|z-w|\ge r>0$, apply
Theorem~\ref{gaussian:gpm:factorization} to
$h-\rhoR(h)$ at the two roots. Its Rayleigh mean is zero, so the
integrated contribution tends to zero uniformly at each fixed $r$.

For near-diagonal domination, we use \eqref{gaussian:gpm:external-bound} with both
exponential weights evaluated at scale $t$, keeping the terminal core separate from the annuli.
On $j$-annuli with $j\le n$, area cancels $e^{2j}$,
and the critical factor is $t$. Splitting at $j=t/2$ gives
\[
 \limsup_{t\to\infty}t\iint_{|z-w|\le e^{-L}}
       \E[F_t^B(z)F_t^B(w)]\dd A(z)\dd A(w)
 \le C\sum_{j\ge L-O(1)}(1+j)^{-3/2}.
\]
The portion $t/2<j\le n$ is $O(t^{-1/2}\log t)$.
For the remaining ball $|z-w|\le e^{-n}$, we have $q\le t$ and
$m=n$ in \eqref{gaussian:gpm:external-bound}; since $t-n<1$, its area times the
critical pair bound is $O(t^{-1/2})$. Thus the displayed right side
tends to zero as $L\to\infty$. Boundedness of the centered marks makes
this an absolute domination of the near-diagonal square integrand.
First take $t\to\infty$ at fixed separation, then separation to zero.
The first assertion follows. The second follows from
Theorem~\ref{gaussian:gp2:gaussian-main}.

\emph{Weak convergence of the marked measures.}
The $L^2$ statements give product-test convergence, and the following
mark-tail estimate supplies tightness in the unbounded coordinate.
Finite sums of products $f(z)h(u)$ form a convergence-determining class
on compact sets. Their convergence in probability follows from
\eqref{gaussian:gpm:marked-L2-bound}. To prevent escape in the mark coordinate,
choose bounded Lipschitz cutoffs $0\le h_R\le1$ equal to one on $[-R,R]$
and zero outside $[-R-1,R+1]$.
For each fixed $R$, the expected discarded mass converges to
$(1-\rhoR(h_R))\E Z_B(K)$, which tends to zero as $R\to\infty$.
The finite first moment here follows from the proved fixed-barrier third
moment; it is not asserted for the untruncated canonical mass.
Markov's inequality gives tightness of the marked finite measures.
Their total masses are $\mu_t^B(K)\to Z_B(K)$, proving weak convergence
in probability. The coarse observations are already defined on the same
space, so joint convergence with them is immediate.
\end{proof}

\paragraph{Scope of the Gaussian conclusions.}
All moment and marked statements in this section keep $B$ fixed.
The canonical measure and its exact-family fixed-barrier restrictions are
identified on their original Gaussian space. The random restriction
$Z_B$ is measurable with respect to the completed sigma-field generated
by the integer observations of this Gaussian family; measurability with
respect to the limiting generalized field alone is not needed. The maximum estimate permits
$Z_B\to M_\kappa$ as $B\to\infty$, but does not assert a moment bound
uniform in $B$. The matrix transfer of any statement here requires its
own first/second moment and cutoff proofs. The separated endpoint
factorization is not a merging-root Rayleigh theorem.

Section~\ref{chap:fixed} uses the separated centered-mark estimate,
together with its own matrix comparison and diagonal bound, to prove
the matrix marked concentration. Its joint identification uses
Proposition~\ref{p2:matrix} and Corollary~\ref{sp:fb:identification}.
Those later applications do not enter the Gaussian arguments above;
the same-space Gaussian $L^2$ limits and the matrix convergence in
distribution have their respective hypotheses and proofs.

\endgroup

\begingroup
\newcommand{\E}{\mathbb E}
\newcommand{\PP}{\mathbb P}
\newcommand{\dd}{\,\mathrm d}
\newcommand{\ind}{\mathbf1}
\section{Scalar regularization and exact normalization}\label{chap:scalar}

In this section we shall compare regularized and original determinant weights while keeping their expected normalizing constants precise enough for the critical limit.

The scalar estimates control the normalizers used in the local-factor,
signed comparison, endpoint and compact-window arguments. We begin
with the sharp positive one-point moment formula, then compare a
regularized determinant with an unregularized determinant after an
auxiliary Gaussian perturbation. The exact deterministic centering cancels
the cost of that perturbation.

The positive one-point theorem of Cipolloni--Landon
\cite[Theorem 9.7]{CL} is the external asymptotic input. The
minor expansion and conditional Jensen comparison below are elementary;
the potential identity follows from the scalar form of the matrix Dyson
equation (MDE). We retain their
short proofs to record the Gaussian normalization, the restriction
$\lambda\ge2$, and the exact deterministic cancellation. Together they
give the quantitative positive-weight $L^1$ comparison in
Theorem~\ref{scalar:thm:sr-regularization}; its total-variation conclusion
is a consequence of this comparison and the normalizer bounds.
These scalar arguments are independent of the later local-factor,
two-root, signed-source and measure-identification arguments.

\subsection{Model, density centering, and scalar asymptotics}

 We  first fix the field centering and recall the sharp determinant moment formula that supplies the expected normalizing constants.
Let
\[
 X=N^{-1/2}(\xi_{ab})_{a,b\le N},\qquad
 \E\xi=\E\xi^2=0,\quad \E|\xi|^2=1,\quad
 \sup_N\E|\xi|^p\le C_p<\infty
\]
for every fixed positive $p$. The entries are independent and
identically distributed within each $N$; their law may depend on $N$.
By \emph{common moment bounds}, we mean that the same prescribed bounds
$(C_p)_{p>0}$ hold for every law in the family. Put
$\kappa=\E|\xi|^4-2\ge-1$; this is the cumulant of the entry law
at the current $N$. Fix a compact bulk disk $|z|\le r_0<1$.

Let $\rho^z$ denote the symmetric, probability-normalized
deterministic Hermitization density. Its Stieltjes transform is
\[
 m^z(w)=\int_{\mathbb R}\frac{\rho^z(x)}{x-w}\,\dd x,\qquad
 -\frac1m=w+m-\frac{|z|^2}{w+m}.
\]
The branch has $\Im m>0$ in the upper half-plane and
$m(w)\sim-1/w$ at infinity. Define
\begin{align}
 D_\eta(A)&=\det(AA^*+\eta^2I_N),\nonumber\\
 I_z(\eta)&=\int_{\mathbb R}\log(x^2+\eta^2)\rho^z(x)\,\dd x,\nonumber\\
 \Phi_N(z,\eta;X)&=\frac12\log D_\eta(X-zI)
                             -\frac N2 I_z(\eta).              \label{scalar:eq:sr-field}
\end{align}
At $\eta=0$ the logarithm may equal $-\infty$; its positive
exponential is assigned value zero. The deterministic identity proved
below gives $I_z(0)=|z|^2-1$, so
\begin{equation}
 e^{\lambda\Phi_N(z,0;X)}
 =|\det(X-zI)|^\lambda
       e^{-\lambda N(|z|^2-1)/2}.                              \label{scalar:eq:sr-raw}
\end{equation}
For $\lambda>0$ set
\[
 W_{\lambda,\eta}=e^{\lambda\Phi_N(z,\eta;X)},\qquad
 Z_{\lambda,\eta}(z;X)=\E W_{\lambda,\eta},\qquad
 \frac{\dd Q_{\lambda,\eta}}{\dd\PP}
       =\frac{W_{\lambda,\eta}}{Z_{\lambda,\eta}}.
\]
At $\lambda=0$ set $W_{0,\eta}=Z_{0,\eta}=1$, including on determinant
zeros. The normalizers are ordinary expectations, not random masses.

The following estimate applies on every fixed compact interval of
nonnegative charges. Its quantitative uniformity will allow the entry
law and the charge to vary with $N$ within the stated class.

\begin{theorem}[Sharp positive one-point moments; Cipolloni--Landon]\label{scalar:thm:sr-sharp}
Fix a finite $L$. Uniformly for $\lambda\in[0,L]$, $|z|\le r_0$,
and complex iid laws with prescribed common moment bounds,
there are $\alpha>0$ and $C<\infty$ such that
\begin{equation}
 Z_{\lambda,0}(z;X)
 =N^{\lambda^2/8}c_\lambda
   \exp\left\{\frac{\lambda(\lambda-2)}8
                   \kappa(1-|z|^2)^2\right\}
         [1+O(N^{-\alpha})].                                 \label{scalar:eq:sr-sharp}
\end{equation}
Here
\[
 c_\lambda=\frac{(2\pi)^{\lambda/4}}{G(1+\lambda/2)}
\]
is positive and continuous on the stated interval, and $G$ is the
Barnes $G$ function. In addition,
\begin{equation}
 Z_{\lambda,0}(z;X)\le C_L N^{\lambda(\lambda+2)/8},
       \qquad 0\le\lambda\le L .                              \label{scalar:eq:sr-raw-polynomial}
\end{equation}
\end{theorem}

This is the complex case of \cite[Theorem 9.7, equation (9.29)]{CL};
Theorem 1.4, equation (1.14), at one point gives the same sharp
formula. Its separation condition is then empty. The polynomial
bound is \cite[Corollary 2.9, equation (2.22)]{CL}.
When an estimate is stated for $L\ge10$, we enlarge the chosen finite
interval to $[0,\max\{10,L\}]$. No subcritical restriction occurs in
this scalar theorem.

The Ginibre Barnes-$G$ factor has the direct predecessor
\cite[Theorem 1.1]{WW}. Webb--Wong state uniform spatial errors on
fixed bulk annuli away from the origin, with the fixed origin treated
separately. We use \cite{CL} for uniformity on the full closed bulk
disk and over the prescribed iid moment class.

The uniform moment-class interpretation is quantitative.
Section 1.4 of \cite{CL} specifies dependence on the common fixed
moment bounds. In the proof of its Theorem 9.7 one can fix the
matching Gaussian variance $T=N^{-1/10}$, then choose the
submicroscopic exponent $\delta_w>0$ and the homogenization
parameters once for these moment bounds and the exponent interval.
After division by the one-point scale in that paper,
$D(z,\lambda)=N^{\lambda^2/8}$ in the complex class,
the relative comparison error in (9.24) is bounded by
\[
 N^{6\delta_w}(N^{-1/10}+N^{-1/2}),
\]
the scalar removal error in (9.31) contains
$N^{-\delta_w^{3/2}/4}$, and Corollary 4.7 supplies another fixed
positive power after those choices. Choosing $\delta_w$ small gives
one $\alpha>0$ and one collection of constants for the family.
The proof itself uses an $N$-dependent moment-matched ensemble.
Thus \eqref{scalar:eq:sr-sharp} is not obtained by extrapolating a
qualitative fixed-law limit to arbitrary triangular arrays.
Common moment bounds, the fixed bulk margin, and the fixed exponent
interval remain essential.
In applications below, we first enlarge the common moment bounds to include
all normalized Gaussian convolutions in \eqref{scalar:eq:sr-convolution},
and then fix the rate, constants, and sufficiently-large-$N$
threshold for these bounds and the whole charge interval. Thus the
same assertion applies to deterministic choices $\lambda=\lambda_N$
within that interval and to the corresponding triangular families.
No convergence of their fourth cumulants is needed for the scalar
estimate, whose leading factor uses the cumulant at the current $N$.

\subsection{Gaussian minors and conditional Jensen}

The purpose of this subsection is to compare a regularized determinant with a determinant after Gaussian perturbation by expanding both into the same minors.
The minor coefficients compare the regularized determinant with the
second moment of a Gaussian perturbation. The complementary-minor
expansion and vanishing-cross-term calculation already appear in
Feldheim--Paquette--Zeitouni \cite[proof of Lemma 4.2, (18)--(20)]{FPZ}.
Their real-Gaussian calculation has the circular-complex counterpart
recorded below. We use that algebraic calculation, rather than their
block-matrix theorem. Conditional Jensen then extends the comparison
to real charges $\lambda\ge2$.

Every Gaussian matrix $G$ in this section has independent standard
circular complex entries with
$\E G_{ab}= \E G_{ab}^2=0$ and $\E|G_{ab}|^2=1$.
Thus $N^{-1/2}G$, rather than $G$, is normalized Ginibre.

\begin{lemma}[Elementary minor comparison]\label{scalar:lem:sr-minor}
For a deterministic complex $N$ by $N$ matrix $A$, put
\[
 b_k(A)=\sum_{\substack{I,J\subset\{1,\ldots,N\}\\|I|=|J|=N-k}}
                      |\det A_{I,J}|^2,\qquad 0\le k\le N,
\]
with the empty determinant equal to one. Then, for $\eta,h\ge0$,
\begin{align}
 D_\eta(A)&=\sum_{k=0}^N\eta^{2k}b_k(A),                       \label{scalar:eq:sr-minors}\\
 \E_G|\det(A+hG)|^2
          &=\sum_{k=0}^N k!h^{2k}b_k(A).                     \label{scalar:eq:sr-gaussian-minors}
\end{align}
Consequently for $s>0$ and every real $\lambda\ge2$,
\begin{equation}
 D_\eta(A)^{\lambda/2}
 \le \exp\!\left(\frac{\lambda N\eta^2}{2s}\right)
          \E_G|\det(A+\sqrt{s/N}\,G)|^\lambda .               \label{scalar:eq:sr-jensen}
\end{equation}
\end{lemma}

\begin{proof}
Expanding $\det(AA^*+\eta^2I)$ by principal minors and applying
Cauchy--Binet to each principal minor of $AA^*$ gives
\eqref{scalar:eq:sr-minors}.

Expand $\det(A+hG)$ by the subsets of rows and columns occupied
by Gaussian entries. The complementary deterministic minor is
$A_{I,J}$, and the Gaussian minor has size $k$. Holomorphic
Gaussian monomials with different entry exponent vectors are
orthogonal in $L^2$. Different row/column sets or total degrees
therefore give zero cross terms. Within a fixed $k$ by $k$
Gaussian minor, only identical permutations survive in its squared
determinant expectation. Each has expectation one and there are
$k!$ permutations. This proves \eqref{scalar:eq:sr-gaussian-minors},
including $k=0$ and $k=N$.

Put $x=N\eta^2/s$. The elementary inequality $x^k/k!\le e^x$
and positivity of the $b_k$ give
\[
 D_\eta(A)\le e^{N\eta^2/s}
                  \E_G|\det(A+\sqrt{s/N}\,G)|^2.
\]
The map $v\mapsto v^{\lambda/2}$ is convex for every real
$\lambda\ge2$. Conditional Jensen proves \eqref{scalar:eq:sr-jensen}.
This is valid at noninteger $\lambda=2\sqrt2$, and at $\lambda=2$
Jensen is an equality. No claim with this direction is made for
$0<\lambda<2$.
\end{proof}

No unitary invariance of the target iid matrix was used. The
auxiliary Gaussian is the only invariant ensemble in this
calculation, and nonsingularity of $A$ is unnecessary.

\subsection{The exact deterministic cancellation}

In this subsection we shall show that the deterministic centering cancels the exponential cost of the Gaussian perturbation when both are chosen consistently.
The preceding comparison has an exponential cost proportional to $N$.
To cancel it, we keep the deterministic potential exactly; an
asymptotic replacement could leave an exponential error. The MDE
selects both the perturbation size and the corresponding change of
spatial argument.

\begin{lemma}[MDE potential identity]\label{scalar:lem:sr-potential}
For $|z|\le r_0$ and $\eta>0$ let
\[
 y=\Im m^z(i\eta)>0,\qquad
 u=\frac{y}{\eta+y},\qquad s=\frac{\eta}{y}.
\]
Then
\begin{equation}
 I_z(\eta)=\frac{\eta^2}{s}+\log(1+s)
                       +\frac{|z|^2}{1+s}-1.                 \label{scalar:eq:sr-potential}
\end{equation}
In particular $I_z(0)=|z|^2-1$, $I_z'(\eta)=2y$,
and $s=O_{r_0}(\eta)$ as $\eta\downarrow0$, uniformly in $z$.
\end{lemma}

\begin{proof}
The imaginary-axis MDE gives
\[
 \frac1y=\eta+y+\frac{|z|^2}{\eta+y},\qquad
 u=\frac1{1+s},\qquad y^2=u-|z|^2u^2 .
\]
Direct differentiation under the deterministic integral yields
\[
 I_z'(\eta)=2\eta\int\frac{\rho^z(x)}{x^2+\eta^2}\,\dd x
           =2y.
\]
For $J_z(s)=\log(1+s)+|z|^2/(1+s)-1$,
\[
 J_z'(s)=u-|z|^2u^2=y^2=\eta^2/s^2.
\]
Therefore
\[
 \frac{\dd}{\dd\eta}\{\eta^2/s(\eta)+J_z(s(\eta))\}
 =2\eta/s+s'(-\eta^2/s^2+J_z'(s))=2y=I_z'(\eta).
\]
The integration constant is zero. Indeed, as $\eta\to\infty$,
$y=\eta^{-1}+O(\eta^{-3})$, $s=\eta^2+O(1)$, and both
sides of \eqref{scalar:eq:sr-potential} equal $2\log\eta+o(1)$.
As $\eta\downarrow0$, $y\to\sqrt{1-|z|^2}$, uniformly positive
in the fixed bulk. Taking the limit gives $I_z(0)=|z|^2-1$ and
the claimed bound on $s$.
\end{proof}

We normalize the perturbed matrix to restore entry variance $1/N$.
The resulting change of spatial argument is included in the same
cancellation. Define
\begin{equation}
 X^{(s)}=\frac{X+\sqrt{s/N}\,G}{\sqrt{1+s}},\qquad
 \xi_s=\frac{\xi+\sqrt s\,g}{\sqrt{1+s}},\qquad
 w=\frac z{\sqrt{1+s}},                                     \label{scalar:eq:sr-convolution}
\end{equation}
where $g$ is standard circular complex Gaussian independent of $\xi$.

\begin{proposition}[Exact upper comparison]\label{scalar:prop:sr-exact-upper}
For every deterministic $\eta>0$, choose
$s=\eta/\Im m^z(i\eta)$. For every real $\lambda\ge2$,
\begin{equation}
 Z_{\lambda,\eta}(z;X)
       \le Z_{\lambda,0}(w;X^{(s)}).                          \label{scalar:eq:sr-exact-upper}
\end{equation}
This has no additional deterministic exponential factor.
\end{proposition}

\begin{proof}
Apply \eqref{scalar:eq:sr-jensen} conditionally on $X$ to $A=X-zI$,
multiply by $\exp(-\lambda N I_z(\eta)/2)$, and average.
The determinant scaling is
\[
 \det(X-zI+\sqrt{s/N}G)
       =(1+s)^{N/2}\det(X^{(s)}-wI).
\]
Using \eqref{scalar:eq:sr-raw}, the entire remaining prefactor is
\[
 \exp\left\{\frac{\lambda N}{2}
 [\,\eta^2/s+\log(1+s)+|w|^2-1-I_z(\eta)\,]\right\}=1
\]
by \eqref{scalar:eq:sr-potential}. This proves the assertion.
All expectations can first be interpreted by Tonelli for
nonnegative variables; the positive scalar moment bounds give
finiteness.
\end{proof}

\subsection{Common moments and positive regularized bounds}

In this subsection we bound positive moments of the regularized weights uniformly over the entry laws produced by Gaussian perturbation.
To apply the sharp raw estimate after Gaussian perturbation, all
normalized entry laws must satisfy common moment bounds. This includes
perturbation sizes depending deterministically on $N$, $z$ and $\eta$.
The normalized convolution satisfies exactly
\begin{equation}
 \E\xi_s=\E\xi_s^2=0,\quad \E|\xi_s|^2=1,\quad
 \kappa_s=\frac{\kappa}{(1+s)^2}.                              \label{scalar:eq:sr-cumulant}
\end{equation}
The cumulant identity follows from addition of independent
fourth mixed cumulants and degree-four scaling; the Gaussian
mixed fourth cumulant is zero. For every fixed $p\ge1$,
\begin{equation}
 \sup_{s\ge0}\E|\xi_s|^p
       \le2^{p-1}(C_p+\E|g|^p).                              \label{scalar:eq:sr-common}
\end{equation}
For $0<p<1$ the analogous bound follows by subadditivity.
Thus all these laws, including deterministic choices
$s=s(N,z,\eta)$, satisfy common moment bounds. Their entries remain
independent. Also $|w|\le|z|\le r_0$.

\begin{proposition}[Positive regularized moments]\label{scalar:prop:sr-moments}
For every fixed compact charge interval $[2,L]$,
\begin{equation}
 \sup_{\lambda\in[2,L],\,|z|\le r_0,\,\eta\ge0}
 \frac{Z_{\lambda,\eta}(z;X)}
      {N^{\lambda(\lambda+2)/8}}\le C_L.                      \label{scalar:eq:sr-allheight}
\end{equation}
For each fixed bounded height interval $0\le\eta\le\eta_0$ one
also has the sharper bound
\begin{equation}
 Z_{\lambda,\eta}(z;X)\le C_{L,\eta_0}N^{\lambda^2/8}.
                                                                    \label{scalar:eq:sr-sharp-upper}
\end{equation}
Both bounds are uniform for the prescribed common moment class.
\end{proposition}

\begin{proof}
For $\eta>0$, use \eqref{scalar:eq:sr-exact-upper} and
\eqref{scalar:eq:sr-common}. Apply the raw polynomial bound
\eqref{scalar:eq:sr-raw-polynomial}; this works even for unbounded
deterministic $s$, proving \eqref{scalar:eq:sr-allheight}.
For bounded $\eta$, the deterministic $s$ is bounded uniformly
on the fixed bulk, and the sharp formula
\eqref{scalar:eq:sr-sharp}, applied to the same common moment class,
gives \eqref{scalar:eq:sr-sharp-upper}. The case $\eta=0$ follows
directly from Theorem~\ref{scalar:thm:sr-sharp}.
\end{proof}

The polynomial estimate alone suffices for all the weighted
expectations on exceptional events used below. For example, at $\chi=2\sqrt2$,
\[
 \E e^{2\chi\Phi(z,\eta)}
       \le CN^{4+\sqrt2},\qquad
 \E e^{4\chi\Phi(z,\eta)}
       \le CN^{16+2\sqrt2}.
\]
Since $Z_{\chi,0}\asymp N$ by \eqref{scalar:eq:sr-sharp},
the raw critical exponential weight, normalized by its expectation,
has fourth moment at most $CN^{12+2\sqrt2}$. Stronger bounded-height moment exponents are
available from \eqref{scalar:eq:sr-sharp-upper} when convenient.

\subsection{One-sided upper ratios and microscopic total variation}

The purpose of this subsection is to quantify how removing microscopic regularization changes the determinant weights, their expectations, and the probability laws they define.
The exact upper comparison controls the normalizer even when
$N\eta$ is large. To compare the weights in $L^1$ and the rooted
laws in total variation, we also need a pointwise lower bound;
its error has the larger scale $N\eta$.

Cipolloni--Landon already establish capped positive-weight $L^1$
regularization at $\eta=N^{-1-\delta}$ for fixed small $\delta>0$,
including the critical charge: see \cite[Lemma 4.3, (4.11), and
Proposition 6.1, (6.5)--(6.7)]{CL}. Their results also cover real
matrices and smaller positive charges. In the overlapping complex,
$\lambda\ge2$ setting, the exact cancellation above gives the uncapped
error $O(N\eta+N^{-\alpha})$ for $0\le\eta\le c/N$, including
$\eta=N^{-1}(\log N)^{-q}$, and a one-sided normalizer bound at all
small deterministic heights. These height and error estimates require
the argument below; a theorem with fixed $\delta$ cannot be used with
$\delta$ tending to zero without uniform parameter control.

\begin{theorem}[Quantitative scalar regularization]\label{scalar:thm:sr-regularization}
Fix a compact charge interval $\lambda\in[2,L]$.
There are a small constant $\eta_0>0$, $\alpha>0$, and $C<\infty$
such that uniformly for $|z|\le r_0$ and $0\le\eta\le\eta_0$,
\begin{align}
 \frac{Z_{\lambda,\eta}}{Z_{\lambda,0}}
       &\le1+C(\eta+N^{-\alpha}),                            \label{scalar:eq:sr-upper-ratio}\\
 W_{\lambda,\eta}
       &\ge e^{-d_{\lambda,\eta}}W_{\lambda,0},\qquad
 d_{\lambda,\eta}=\lambda N\int_0^\eta
                         \Im m^z(iv)\,\dd v\le CN\eta,        \label{scalar:eq:sr-pointwise-lower}\\
 e^{-CN\eta}\le\frac{Z_{\lambda,\eta}}{Z_{\lambda,0}}
       &\le1+C(\eta+N^{-\alpha}).                            \label{scalar:eq:sr-ratio}
\end{align}
For $0\le\eta\le c/N$, with any fixed $c<\infty$,
\begin{align}
 \left|\frac{Z_{\lambda,\eta}}{Z_{\lambda,0}}-1\right|
       &\le C(N\eta+N^{-\alpha}),                            \label{scalar:eq:sr-two-sided}\\
 \frac{\E|W_{\lambda,\eta}-W_{\lambda,0}|}{Z_{\lambda,0}}
       &\le C(N\eta+N^{-\alpha}),                            \label{scalar:eq:sr-L1}\\
 \|Q_{\lambda,\eta}-Q_{\lambda,0}\|_{\mathrm{TV}}
       &\le C(N\eta+N^{-\alpha}).                            \label{scalar:eq:sr-TV}
\end{align}
Here $\|\sigma\|_{\mathrm{TV}}=|\sigma|(\Omega)$ is the full
variation norm, as in Section~\ref{sp:common}. For probability measures,
$\|P-Q\|_{\mathrm{TV}}=2\sup_A|P(A)-Q(A)|$.
In particular for $\eta_{\rm m}=N^{-1}t^{-q}$,
$t=\frac12\log N$ and fixed $q>0$, the last three errors are
$O(t^{-q}+N^{-\alpha})$. At $\lambda=\chi=2\sqrt2$,
both corresponding normalizers are comparable with $N$.
\end{theorem}

\begin{proof}
Let $h=1-|z|^2$ and use \eqref{scalar:eq:sr-convolution}. The exact
cumulant and spatial transformations give
\[
 \kappa_s(1-|w|^2)^2
       =\kappa\,\frac{(h+s)^2}{(1+s)^4}
       =\kappa h^2+O(s),
\]
uniformly in the moment class and the fixed bulk, for small $s$.
Apply the sharp formula \eqref{scalar:eq:sr-sharp} to the right side
of \eqref{scalar:eq:sr-exact-upper} and to $Z_{\lambda,0}(z;X)$.
Their factors $N^{\lambda^2/8}c_\lambda$ cancel exactly. The
leading factors are bounded above and below uniformly. Since
$s=O(\eta)$, their quotient is at most
\[
 e^{C\eta}\frac{1+CN^{-\alpha}}{1-CN^{-\alpha}}
       \le1+C'(\eta+N^{-\alpha})
\]
for large $N$ and sufficiently small fixed $\eta_0$.
The cancellation uses the same charge on both sides and therefore
also holds for a moving deterministic charge in the fixed compact
interval; no derivative estimate for the scalar asymptotic error
is required. The height and the resulting $s$ are deterministic
parameters throughout this comparison.
This proves \eqref{scalar:eq:sr-upper-ratio} for all those small heights,
including heights with $N\eta$ tending to infinity.

For the lower comparison, $D_\eta(X-zI)\ge D_0(X-zI)$
pointwise. The deterministic centering increment is
\[
 \frac{\lambda N}{2}(I_z(\eta)-I_z(0))
       =\lambda N\int_0^\eta \Im m^z(iv)\,\dd v
       =d_{\lambda,\eta}\le CN\eta.
\]
This proves \eqref{scalar:eq:sr-pointwise-lower} and its integrated
version. It does not reverse Jensen. Combining it with the
upper ratio gives \eqref{scalar:eq:sr-two-sided}.

For the unnormalized $L^1$ estimate use the exact identity
\[
 \E|W_{\lambda,\eta}-W_{\lambda,0}|
 =Z_{\lambda,\eta}-Z_{\lambda,0}
        +2\E(W_{\lambda,0}-W_{\lambda,\eta})_+ .
\]
Its right side is at most
\[
 \{C(\eta+N^{-\alpha})+2(1-e^{-d_{\lambda,\eta}})\}
                                           Z_{\lambda,0},
\]
which proves \eqref{scalar:eq:sr-L1}.
Furthermore, the normalized laws satisfy
\[
 Q_{\lambda,\eta}\ge
 \frac{e^{-d_{\lambda,\eta}}}{1+C(\eta+N^{-\alpha})}
                         Q_{\lambda,0}
\]
as nonnegative measures. Subtract the displayed common
subprobability measure from both laws. Each residual has mass at
most $C(N\eta+N^{-\alpha})$, so their difference has full variation
at most twice this amount. Absorbing the factor two in $C$ proves
\eqref{scalar:eq:sr-TV}. For any bounded signed test $F$,
$|\E_{Q_{\lambda,\eta}}F-\E_{Q_{\lambda,0}}F|$ is bounded by
$\|F\|_\infty\|Q_{\lambda,\eta}-Q_{\lambda,0}\|_{\mathrm{TV}}$.

At $\eta_{\rm m}=N^{-1}t^{-q}$ the claimed rate follows.
The sharp raw formula gives $Z_{\chi,0}\asymp N$ since
$\chi^2/8=1$, and the two-sided ratio gives the same for
$Z_{\chi,\eta_{\rm m}}$. Determinant-zero events create no
exception: $W_{\lambda,0}=0$ there.
\end{proof}

\begin{remark}[The two different scales]
The upper ratio has error $O(\eta+N^{-\alpha})$ at every sufficiently
small deterministic height. The pointwise lower comparison and
two-sided calibration have the larger scale $O(N\eta+N^{-\alpha})$.
The upper estimate is therefore available at terminal integer heights
$\eta_J=t^C/N\to0$, while a small two-sided or total-variation error
is used only at submicroscopic heights. The proof does not assert a
two-sided $O(\eta)$ comparison.
\end{remark}

The $L^1$ statement may be used with the same bounded random test on an
enlarged probability space carrying the matrix, including initial
history or independent reference variables. It uses only the
scalar expectation of the absolute difference; the test need not
be independent of the matrix. If the observable itself changes
with the regularization, a separate estimate for the difference
of the two observables is required. The scalar statement does \emph{not}
allow multiplication by a second unbounded exponential weight;
comparisons with two such weights require their separate pointwise
reconstruction and error bounds for the pair.

\subsection{Exact variance rescaling and the ensembles in the proof}

{We shall now rescale each matrix together with its deterministic density so that the scalar estimates apply at its exact variance.}
The actual aged matrix and the Gaussian reference have different
entry variances. Each uses its own deterministic Hermitization
density. Rescaling both the determinant and this density makes the
scalar transfer exact.

\begin{lemma}[Variance scaling]\label{scalar:lem:sr-variance}
Suppose $X$ has entry variance $a/N$, with $a$ in a fixed compact
subset of $(0,\infty)$. Let
\[
 \widehat X=X/\sqrt a,\qquad
 \widehat z=z/\sqrt a,\qquad \widehat\eta=\eta/\sqrt a,
\]
and suppose $\widehat z$ remains in a fixed unit-variance bulk disk.
For transfer of the scalar estimates assume that $\widehat X$ belongs
to the preceding complex iid class with its prescribed common moment
bounds. The displayed scaling identities themselves are deterministic.
Write $\rho_a^z$ for the exact variance-$a$ deterministic density.
Define $I_{a,z}$ and $\Phi_{a,N}$ as in \eqref{scalar:eq:sr-field},
using $\rho_a^z$ for the deterministic density. Then
\begin{align}
 \rho_a^z(x)&=a^{-1/2}\rho_1^{\widehat z}(x/\sqrt a),\nonumber\\
 I_{a,z}(\eta)&=\log a+I_{\widehat z}(\widehat\eta),\nonumber\\
 D_\eta(X-zI)&=a^N D_{\widehat\eta}(\widehat X-\widehat zI),\nonumber\\
 \Phi_{a,N}(z,\eta;X)
       &=\Phi_N(\widehat z,\widehat\eta;\widehat X).            \label{scalar:eq:sr-variance}
\end{align}
All scalar statements above transfer with uniform constants under
this exact change of variables.
\end{lemma}
\begin{proof}
The first identity follows by scaling the MDE, equivalently the
deterministic Hermitization spectral variable. In the logarithmic
integral the change of variable contributes $\log a$ because
the density has total mass one. The determinant contributes
$N\log a$ before its factor $1/2$. Thus the two deterministic
$N\log a/2$ terms cancel exactly. Bounded scaling preserves the
common moment bounds and the fixed bulk condition.
\end{proof}

For the aged additive matrix
\[
 X_s=\sqrt{1-\sigma}\,\widehat X
              +\sqrt{\sigma}\,G^{\mathrm{norm}}_0
              +N^{-1/2}B_s,\qquad \sigma=N^{-1/200},
\]
$G^{\mathrm{norm}}_0$ is normalized Ginibre and the entries of
$B_s$ have variance $s$. Its entry variance is $(1+s)/N$.
The normalized entries satisfy common moment bounds, even when
$\widehat X$ is an $N$-dependent moment-matched law. Their fourth
cumulant is exactly
\[
 \kappa_{\mathrm{norm},s}
       =\left(\frac{1-\sigma}{1+s}\right)^2\widehat\kappa .
\]
Here $\widehat\kappa$ is the fourth cumulant of the variance-one entry
variables $\sqrt N\,\widehat X_{ab}$.
The local-factor centering uses $\rho_{1+s}^z$, so
Lemma~\ref{scalar:lem:sr-variance} applies at every deterministic time
in the bounded intervals used by the coupling, including $\tau=t^b/N$.

The independent Gaussian reference initialized from
$G^{\mathrm{norm}}/B_v$, with
$B_v=\sqrt{1-|v|^2}$, has at additive time $s$ the same marginal
law as a normalized Ginibre matrix scaled to variance
$a=B_v^{-2}+s$, at spatial point zero. On a fixed bulk compact and
the bounded additive-time intervals used below, $a$ lies in a fixed
compact subset of $(0,\infty)$. The same lemma
therefore applies with the reference system's own centering.
Only its marginal scalar law is used here.

No factor $e^{CN(a-1)}$ or $e^{CN\tau}$ is introduced in either
case. The cancellation in \eqref{scalar:eq:sr-variance} is exact.

\subsection{Applications of the scalar estimates}

The purpose of this subsection is to specify how the scalar estimates control normalizers, exceptional events, and removal of microscopic regularization in the later comparisons.

The local-factor calibration in Section~\ref{chap:lf} uses
$\chi=2\sqrt2$ and $\eta_{\rm m}=t^{-q}/N$. The preceding bounds give
\[
 Z_{\chi,0}\asymp N,\quad Z_{\chi,\eta_{\rm m}}\asymp N,\quad
 \frac{Z_{\chi,\eta_{\rm m}}}{Z_{\chi,0}}
       =1+O(t^{-q}+N^{-\alpha}),
\]
together with \eqref{scalar:eq:sr-L1}--\eqref{scalar:eq:sr-TV}.
Thus an independently proved reconstruction
$(\E U_v)(\E V_v)=Z_{\chi,\eta_{\rm m}}(1+o(t^{-P}))$
has the same form with $Z_{\chi,0}$ when $q>P$, without a pair
normalizer. The endpoint comparisons in Sections~\ref{chap:oneroot}
and~\ref{chap:separated} use the upper ratio
\eqref{scalar:eq:sr-upper-ratio} at larger vanishing heights.
Section~\ref{chap:window} fixes all one-root constants and rates on
$[2,4]$. Positive moments control weighted exceptional events and
cap removal; any larger fixed charges needed for H\"older estimates
are included by enlarging $L$ before the required moment orders and
powers in the failure probability bounds are selected.

The bounded mean correction
$-\kappa m^z(i\eta)^4/4$ and the factor $\sqrt2$ in the field definition
change the scalar exponent or multiply an unnormalized weight by a
bounded deterministic quantity. For an exponential weight normalized by its expectation, the
deterministic factor cancels exactly. Negative coefficients of capped coordinates are treated
by the signed-source theorem; no negative raw determinant moment
is claimed or used in this section.

\endgroup

\begingroup
\newcommand{\E}{\mathbb E}
\newcommand{\PP}{\mathbb P}
\newcommand{\Tr}{\operatorname{Tr}}
\newcommand{\HS}{\mathrm{HS}}
\newcommand{\dd}{\mathrm d}
\newcommand{\one}{\mathbf1}
\newcommand{\cS}{\mathcal S}
\newcommand{\cU}{\mathcal U}
\newcommand{\cM}{\mathcal M}
\newcommand{\br}[1]{\langle #1\rangle}
\newcommand{\norm}[1]{\lVert #1\rVert}
\section{Full block resolvent insertions at unequal heights}\label{chap:i2}

In this section we shall control resolvent products with block insertions at unequal heights and nearby spatial points for the later critical comparisons.

The logarithmic comparison arguments use products of resolvents at
different heights with constant block insertions. Unequal heights and
arbitrary bounded deterministic insertions are already covered by
Cipolloni--Landon \cite[Proposition B.2]{CL}. Their logarithmic estimate
applies at coincident or fixed-separated spatial points. Here
Theorem~\ref{i2:thm:I2} gives logarithmic control for the stated range of
nearby spatial points, without a lower bound on their separation, for all four constant
block directions, including off-diagonal insertions. The matrix has an
independent Gaussian component whose age is comparable to $N^{-1/200}$.

The proof begins with the exact stochastic equations and their
deterministic propagator. An initial estimate at the Gaussian age,
followed by drift and martingale bounds, controls the whole
characteristic. Coarse and large-height extensions then give the
integrated insertion estimate used in signed comparison. Its
coordinate heights are at most one, while the integration variable
extends to $N^{10}$. Each use of an external local law is stated
with the required domain and probability bound.

\subsection{Statement, conventions, and local-law inputs}

 {We first state the resolvent estimate together with the ranges and existing local laws needed to prove it.}
Fix compact disks $K\Subset K'\Subset\{z:|z|<1\}$. Let
$\widehat X=N^{-1/2}(\widehat\xi_{ab})$ have independent identically
distributed complex entries, possibly depending on $N$, satisfying
\[
 \E\widehat\xi=\E\widehat\xi^2=0,
 \qquad \E|\widehat\xi|^2=1,
 \qquad \sup_N\E|\widehat\xi|^r<\infty\quad(r<\infty).
\]
All constants below may depend on the fixed moment bounds. Put
\[
 t=\tfrac12\log N,\quad \sigma=N^{-1/200},\quad
 c_*\sigma\le q_N\le C_*\sigma,\quad
 \theta=-\log(1-q_N),
\]
where $0<c_*<C_*<\infty$ are fixed, and define the aged matrix
\[
 X=\sqrt{1-q_N}\,\widehat X+\sqrt{q_N}\,G.
\]
Here $G$ is an independent normalized complex Ginibre matrix. This
form allows the slightly different ages produced by moment matching
and exact variance normalization.

Write
\[
 H^z=\begin{pmatrix}0&X-zI\\(X-zI)^*&0\end{pmatrix},\qquad
 G_z(i\eta)=(H^z-i\eta)^{-1},\qquad
 \br R=(2N)^{-1}\Tr R.
\]
The symbols $E_{ab}$ denote the constant $2$ by $2$ matrix units tensored
with $I_N$. All insertion matrices in this section belong to this
four-dimensional block space. For such a matrix $C$, put
\[
 \norm C_{\rm bl}=\sum_{a,b=1}^2|C_{ab}|.
\]
This norm is equivalent, with absolute constants, to the operator norm
on the constant block space. It suffices to estimate the sixteen
pairs of matrix units; expansion in this basis then gives the bound
for arbitrary constant block insertions.

Let $M_z(i\eta)$ be the deterministic solution of the Hermitized
matrix Dyson equation, with positive imaginary part. The deterministic
two-resolvent solution is defined by
\begin{equation}
 \cM_{z,x;w,y}[C]
 =M_z(ix)C M_w(iy)
   +M_z(ix)\cS[\cM_{z,x;w,y}[C]]M_w(iy),
 \label{i2:eq:M12}
\end{equation}
where
\[
 \cS[R]=2\br{RE_{22}}E_{11}+2\br{RE_{11}}E_{22}.
\]
For deterministic $C,D$ define
\[
 \mathcal E_{z,w}(x,y;C,D)
 =\br{G_z(ix)C G_w(iy)D}
       -\br{\cM_{z,x;w,y}[C]D}.
\]

\begin{theorem}[Full block insertion estimate]\label{i2:thm:I2}
Take $p=16$ and $b_0=24$. There is a fixed $c>0$ such that, for every
$D_*>0$, outside an event of probability at most $C_{D_*}N^{-D_*}$,
\begin{equation}
 |\mathcal E_{z,w}(x,y;C,D)|
 \le C_{D_*}\norm C_{\rm bl}\norm D_{\rm bl}
       \left(\frac{t^{16}}{Nxy}+N^{-c}\right)
 \label{i2:eq:I2}
\end{equation}
simultaneously for $z,w\in K$, $|z-w|\le C_0N^{-1/4}$, and
$t^{24}/N\le x,y\le1$. The constant $C_0$ is fixed. In particular,
no lower separation is required.

If $x\le\sigma$, the additive term may be omitted, and the resulting
bound $Ct^{16}/(Nxy)$ extends to every $y\ge x$.
For arbitrary $z,w\in K$, the bound extends when
$y\ge\eta_c>0$ is bounded below, in the form
\begin{equation}
 |\mathcal E_{z,w}(x,y;C,D)|
 \le C_{D_*,\eta_c}\norm C_{\rm bl}\norm D_{\rm bl}
 \left(\frac{t^{16}}{Nxy}+\frac{N^{-c}}{1+y}\right),
 \qquad t^{24}/N\le x\le2.
 \label{i2:eq:coarse-I2}
\end{equation}
All these exponents are fixed before the exponent
$b$ used to choose the duration of the later matrix flow. They
are independent of $D_*$; the constants and the threshold for sufficiently large
$N$ may depend on $D_*$.
\end{theorem}

The proof uses the following external statements, within their stated
ranges.
\begin{enumerate}[label=\textup{(E\arabic*)}]
\item\label{i2:input:coincident} Cipolloni--Landon \cite{CL}, arXiv:2605.29962v3,
Proposition B.2(ii). At coincident
spatial points, it permits arbitrary bounded deterministic insertions
and unequal heights, with error $(\log N)^{10}/(Nxy)$.
The lower heights are at least $(\log N)^{10}/N$; the first height is
at most $N^{-\epsilon}$ and the second may be as large as $N^C$.
Thus the initial logarithmic exponent is $p_0=10$.
\item\label{i2:input:polynomial} The same paper, Proposition B.2(i), gives the error
$N^\epsilon/(Nxy)$ for arbitrary spatial points and insertions in the
fixed-power range. Lemma B.1 covers
$N^{\epsilon-1}\le x\le C$ and $C^{-1}\le y\le N^C$ with the
same error. Taking a fixed $C>2$ patches the overlap near $y=1$.
\item Cipolloni--Landon \cite{DEL}, arXiv:2509.15189v1,
Theorem 3.1, gives
\[
 |\br{(G-M)B}|\le\frac{C_D\log N}{N\eta}
\]
for bounded deterministic $B$ in the domain
$N\eta\rho\ge100C^2\log N$, $\eta/\rho\le N^{-\xi}$,
with overwhelming probability. In the bulk,
$\rho\asymp1$ on the fine characteristics below. To cover both
heights up to $2a$, where $a=N^{-10\xi}$ and $\xi=1/10000$,
choose the upper-domain exponent in that theorem strictly less than
$10\xi$. This includes the first characteristic of height
$O(\sigma)$ and the second one between $\sigma$ and $a$. This yields
$p_1=1$ independently of the eventual lower-height exponent.
At larger heights use the usual fixed-power one-resolvent law,
with error $N^\xi/(N\eta)$.
\end{enumerate}
The hypotheses concern the unconditional law at deterministic flow
times. They are not applied after conditioning on a future event.

\subsection{Exact stochastic system}

Next we shall move the spatial points and regularization heights with the matrix flow to cancel the leading drift.
The Ornstein--Uhlenbeck (OU) flow preserves entry variance. Moving the
spatial points and heights along deterministic characteristics cancels
the leading resolvent drift and exposes the scalar local-law errors.
This characteristic-flow approach is used for two-resolvent laws in
\cite[Appendix B.2, equation (B.28)]{CL}. The equations below record it in all sixteen
block coordinates, with the exact forcing and quadratic variations
needed here.

Realize the aged matrix as the endpoint of
\[
 \dd X_s=-\tfrac12X_s\,\dd s+N^{-1/2}\dd B_s,
 \qquad 0\le s\le\theta,
\]
where the complex Brownian entries satisfy
\[
 \E[\dd B_{ab}\dd\overline B_{cd}]
 =\delta_{ac}\delta_{bd}\,\dd s,
 \qquad \E[\dd B_{ab}\dd B_{cd}]=0.
\]
For terminal data $(z_i,\eta_i)$, solve backwards the deterministic
characteristics
\begin{equation}
 \dot z_{i,s}=-z_{i,s}/2,\qquad
 \dot\eta_{i,s}=-v_{i,s}-\eta_{i,s}/2,
 \qquad m_{i,s}=iv_{i,s}.
 \label{i2:eq:characteristic}
\end{equation}
The term $-\eta/2$ is needed for this variance-preserving convention.
The characteristics remain in $K'$ for all sufficiently large $N$.

\begin{lemma}[Exact resolvent and insertion equations]\label{i2:lem:SDE}
Use the filtration generated by the initial matrix and the Brownian
path, with Brownian increments independent of the initial matrix.
Let $G_i=(H^{z_{i,s}}_s-i\eta_{i,s})^{-1}$ and
$\epsilon_i=\br{G_i}-m_i$. On every interval on which both
$\eta_{i,s}$ are positive,
\begin{equation}
 \dd G_i=-G_i\dd W G_i+
       \left(\tfrac12G_i+\epsilon_iG_i^2\right)\dd s,
 \label{i2:eq:dG}
\end{equation}
where
\[
 \dd W=N^{-1/2}
 \begin{pmatrix}0&\dd B\\\dd B^*&0\end{pmatrix}.
\]
For $C_1=E_{11},C_2=E_{22},C_3=E_{12},C_4=E_{21}$, set
$F_{ab}=\br{G_1C_aG_2C_b}$. If $J_{12}=J_{21}=2$ and all other
entries of $J$ vanish, then
\begin{equation}
 \dd F=(F+FJF+H)\dd s+\dd M,
 \label{i2:eq:16SDE}
\end{equation}
where
\[
 H_{ab}=\epsilon_1\br{G_1^2C_aG_2C_b}
             +\epsilon_2\br{G_1C_aG_2^2C_b}.
\]
The martingale term is given explicitly by
\begin{align}
 \dd M_{ab}
 &=-\br{G_1\dd W G_1C_aG_2C_b}
   -\br{G_1C_aG_2\dd W G_2C_b}
   =-\br{\dd W Q_{ab}},\label{i2:eq:exact-martingale}\\
 Q_{ab}&=G_1C_aG_2C_bG_1+G_2C_bG_1C_aG_2.
 \label{i2:eq:martingale-coefficient}
\end{align}
Thus the quadratic feedback closes on the sixteen coordinates of $F$.
The remaining terms are the displayed forcing $H$ and martingale $M$;
their coefficients still depend on the underlying resolvents.

For coordinate indices $\alpha=(a,b)$ and $\beta=(c,d)$, the exact
quadratic covariations are
\begin{align}
 \dd[M_\alpha,M_\beta]
 &=\frac{\Tr\bigl((Q_\alpha)_{21}(Q_\beta)_{12}\bigr)
        +\Tr\bigl((Q_\alpha)_{12}(Q_\beta)_{21}\bigr)}{4N^3}
        \,\dd s,\label{i2:eq:bilinear-bracket}\\
 \dd[M_\alpha,\overline{M_\beta}]
 &=\frac{\Tr\bigl((Q_\alpha)_{21}(Q_\beta)_{21}^*\bigr)
        +\Tr\bigl((Q_\alpha)_{12}(Q_\beta)_{12}^*\bigr)}{4N^3}
        \,\dd s.\label{i2:eq:sesquilinear-bracket}
\end{align}
Here $(Q_\alpha)_{ij}$ is an $N$ by $N$ block, and the traces in
these two formulas are ordinary block traces.
\end{lemma}
\begin{proof}
\emph{The Brownian contraction and resolvent drift.}
For every adapted matrix $R$ evaluated before the Brownian
increment, direct contraction of its blocks gives
\begin{equation}
 \dd W R\dd W
 =\left(\frac{\Tr R_{22}}{N}E_{11}
        +\frac{\Tr R_{11}}{N}E_{22}\right)\dd s
 =\cS[R]\dd s.
 \label{i2:eq:exact-contraction}
\end{equation}
For example, the upper-left block is
$N^{-1}\dd B R_{22}\dd B^*$. Its $(i,l)$ entry contracts to
\[
 \delta_{il}N^{-1}\sum_j(R_{22})_{jj}\dd s.
\]
The lower-right block is analogous. The off-diagonal blocks contain two unconjugated or
two conjugated Brownian increments and vanish. In particular, the
factor two in the definition of $\cS$ follows from the normalized
trace $\br R=(2N)^{-1}\Tr R$.

For any square matrix $A$, put
$H=\left(\begin{smallmatrix}0&A\\A^*&0\end{smallmatrix}\right)$.
The identity
\[
 (H-\zeta I)^{-1}=(H+\zeta I)(H^2-\zeta^2I)^{-1}
\]
shows that its diagonal resolvent blocks are
$\zeta(AA^*-\zeta^2I)^{-1}$ and
$\zeta(A^*A-\zeta^2I)^{-1}$. Their traces agree because $AA^*$
and $A^*A$ have identical eigenvalues with multiplicity, including
zero. Applying this with $A=X-z_iI$ and $\zeta=i\eta_i$ yields
\begin{equation}
 \cS[G_i]=\br{G_i}I.
 \label{i2:eq:equal-block-traces}
\end{equation}

The spatial characteristic $\dot z_i=-z_i/2$ gives
$\dd H_i=-H_i\dd s/2+\dd W$. It\^o inversion therefore gives
\begin{align*}
 \dd G_i
 &=-G_i\dd W G_i
   +\left(\tfrac12G_iH_iG_i
          +i\dot\eta_iG_i^2+G_i\cS[G_i]G_i\right)\dd s\\
 &=-G_i\dd W G_i
   +\left(\tfrac12G_i
       +(i\eta_i/2+i\dot\eta_i+\br{G_i})G_i^2\right)\dd s.
\end{align*}
Here $G_iH_iG_i=G_i+i\eta_iG_i^2$. By
\eqref{i2:eq:characteristic}, the coefficient in parentheses is
$\br{G_i}-iv_i=\epsilon_i$, proving \eqref{i2:eq:dG}.

\emph{The insertion equation and its quadratic covariations.}
Apply the product rule to $\br{G_1C_aG_2C_b}$. The two linear
drifts contribute $F_{ab}\dd s$, and the two scalar-error drifts
contribute $H_{ab}\dd s$. The cross variation is
\begin{align*}
 \br{G_1\dd W G_1C_aG_2\dd W G_2C_b}
 &=\br{G_1\cS[G_1C_aG_2]G_2C_b}\dd s\\
 &=\left(2F_{a2}F_{1b}+2F_{a1}F_{2b}\right)\dd s
   =(FJF)_{ab}\dd s.
\end{align*}
There is no transpose or adjoint on either factor $F$. The calculation
does not require either insertion to be diagonal or self-adjoint, so
it covers $E_{12}$ and $E_{21}$ as well. The stochastic terms give
the first expression in \eqref{i2:eq:exact-martingale}; cyclicity of the
trace gives the second. This proves \eqref{i2:eq:16SDE}.

Finally, for a block matrix $Q$,
\[
 \br{\dd W Q}
 =\frac1{2N^{3/2}}\sum_{i,j=1}^N
   \left(\dd B_{ij}(Q_{21})_{ji}
        +\dd\overline B_{ij}(Q_{12})_{ij}\right).
\]
Contracting this formula for $Q_\alpha$ and $Q_\beta$ proves
\eqref{i2:eq:bilinear-bracket}. Contracting it with the complex
conjugate of the second expression proves
\eqref{i2:eq:sesquilinear-bracket}. In particular,
\[
 \dd[M_\alpha,\overline{M_\alpha}]
 =\frac{\norm{(Q_\alpha)_{21}}_{\mathrm{HS}}^2
        +\norm{(Q_\alpha)_{12}}_{\mathrm{HS}}^2}{4N^3}\dd s.
\]
All computations are valid after localization. On a compact time
interval with a positive minimum of both heights, the resolvent
bound $\norm{G_i}\le\eta_i^{-1}$ makes the displayed martingale
coefficients bounded, so no localization is needed there.
\end{proof}

\subsection{Deterministic stability and the propagator}

In this subsection we bound the growth of errors under the deterministic evolution, keeping the dependence on both regularization heights explicit.

The two-body stability operator and its finite-dimensional reduction
were analyzed by Cipolloni--Erd\H{o}s--Schr\"oder
\cite[Lemma 6.1]{CES}. We give the coordinate calculation for both
denominator inverses, including the coarse-height bound without an
upper height restriction. Their explicit form gives the linear
propagator, including the order of its left and right factors.
The backward characteristics then express its growth in terms of
the terminal heights.

Let $K=(K_{ab})_{a,b=1}^4$, where
$K_{ab}=\br{M_1C_aM_2C_b}$. Along
\eqref{i2:eq:characteristic} one has $M_i'=M_i/2$, hence $K'=K$;
the exact backward formulas below verify this directly.
Set $\overline F_{ab}
=\br{\cM_{z_1,\eta_1;z_2,\eta_2}[C_a]C_b}$.
Equation \eqref{i2:eq:M12} implies
\begin{equation}
 \overline F=K+\overline FJK
 =K(I-JK)^{-1}=(I-KJ)^{-1}K.
 \label{i2:eq:Fbar}
\end{equation}
In particular, $\overline F'=\overline F+\overline FJ\overline F$.

\begin{lemma}\label{i2:lem:stability}
For $|z_i|\le r_*<1$ and $0<\eta_i\le\Lambda<\infty$, both full
four-dimensional denominator inverses satisfy
\[
 \norm{(I-JK)^{-1}}+\norm{(I-KJ)^{-1}}
 \le\frac C\gamma,
 \qquad \gamma=|z_1-z_2|^2+\eta_1+\eta_2.
\]
Here $C$ may depend on $r_*$ and $\Lambda$.
If one height is at least a fixed positive constant, these inverses
are bounded uniformly over all positive values of the other height,
and over both spatial points in the fixed bulk disk.
\end{lemma}
\begin{proof}
Write
\[
 M_i=\begin{pmatrix}iv_i&-z_i u_i\\-\overline z_i u_i&iv_i\end{pmatrix},
 \qquad u_i=\frac{v_i}{\eta_i+v_i},
 \qquad v_i^2=u_i-|z_i|^2u_i^2.
\]
To cover all block insertions, compute the full coordinate matrix
before estimating its active two-dimensional denominator. Set
$a=iv_1$, $b=-z_1u_1$, $c=-\overline z_1u_1$,
$d=iv_2$, $e=-z_2u_2$, and $f=-\overline z_2u_2$.
The normalized trace and the order
$(E_{11},E_{22},E_{12},E_{21})$ give the full matrix
\[
 2K=\begin{pmatrix}
 ad&ce&cd&ae\\
 bf&ad&af&bd\\
 af&cd&cf&ad\\
 bd&ae&ad&be
 \end{pmatrix}.
\]
Thus, with $P=\left(\begin{smallmatrix}0&1\\1&0\end{smallmatrix}\right)$,
$R_0=\left(\begin{smallmatrix}af&bd\\cd&ae\end{smallmatrix}\right)$,
and $L_0=\left(\begin{smallmatrix}cd&af\\ae&bd\end{smallmatrix}\right)$,
\begin{equation}
 I-JK=\begin{pmatrix}D&-R_0\\0&I_2\end{pmatrix},
 \qquad
 I-KJ=\begin{pmatrix}PDP&0\\-L_0&I_2\end{pmatrix},
 \label{i2:eq:full-denominators}
\end{equation}
where the active two-dimensional denominator is
\[
 D=\begin{pmatrix}
 1-z_1\overline z_2u_1u_2&v_1v_2\\
 v_1v_2&1-\overline z_1z_2u_1u_2
 \end{pmatrix}.
\]
Its determinant has the positive decomposition
\begin{align}
 \det D={}&(1-u_1)(1-|z_1|^2u_1u_2)
 +(1-u_2)u_1(1-|z_2|^2u_2)\nonumber\\
 &+u_1u_2|z_1-z_2|^2.
 \label{i2:eq:positive-det}
\end{align}
To verify this identity, expand the determinant and substitute
$v_i^2=u_i-|z_i|^2u_i^2$; the quartic product cancels, leaving
\[
 \det D=1-u_1u_2+u_1u_2
 \bigl(|z_1-z_2|^2-|z_1|^2(1-u_1)-|z_2|^2(1-u_2)\bigr).
\]
Let $\delta=1-r_*^2$. Since $0<u_i<1$ and $0<v_i\le1$,
the scalar MDE gives $v_i(\eta_i+v_i)\ge\delta$.
For $\eta_i\le\Lambda$, it follows that
$v_i\ge\delta/(\Lambda+1)$ and
$u_i\ge\delta/(\Lambda+1)^2$.
Also $1-u_i=\eta_i/(\eta_i+v_i)\asymp\eta_i$.
All coefficients of the first two terms in
\eqref{i2:eq:positive-det} are therefore bounded below by positive
constants depending only on the displayed compact ranges. The upper
bounds follow from the same formulas, so $\det D\asymp\gamma$.
The entries of $D$ are bounded, whence $\norm{D^{-1}}\le C/\gamma$.
The full inverses, with no omission of the off-diagonal insertions, are
\begin{equation}
 \begin{split}
 (I-JK)^{-1}&=\begin{pmatrix}D^{-1}&D^{-1}R_0\\0&I_2\end{pmatrix},\\
 (I-KJ)^{-1}&=\begin{pmatrix}PD^{-1}P&0\\L_0PD^{-1}P&I_2\end{pmatrix}.
 \end{split}
 \label{i2:eq:full-denominator-inverses}
\end{equation}
The entries of $L_0,R_0$ are uniformly bounded for all positive
heights. This proves the first assertion.
If $\eta_1\ge a_*>0$, then
$\det D\ge\delta(1-u_1)\ge\delta a_* /(a_*+1)$.
If instead $\eta_2\ge a_*$, use the same positive decomposition with
the indices interchanged. The full inverse formulas prove the second
assertion, without an upper restriction on either height.
\end{proof}

The stability estimate now gives the fundamental solution for the
linearized error equation. Set $D_L=I-KJ$, $D_R=I-JK$ and
$E=F-\overline F$. Then
\begin{equation}
 \dd E=(E+\overline FJE+EJ\overline F+EJE+H)\dd s+\dd M.
 \label{i2:eq:error-SDE}
\end{equation}
Its deterministic linear propagator is
\begin{equation}
 \cU(v,u)[A]=e^{v-u}D_L(v)^{-1}D_L(u)
       A D_R(u)D_R(v)^{-1},\qquad u\le v.
 \label{i2:eq:propagator}
\end{equation}
Indeed, write $L(v,u)=D_L(v)^{-1}D_L(u)$ and
$R(v,u)=D_R(u)D_R(v)^{-1}$. Since $K'=K$,
\[
 \partial_v L(v,u)=\overline F(v)J L(v,u),\qquad
 \partial_v R(v,u)=R(v,u)J\overline F(v).
\]
Together with $L(u,u)=R(u,u)=I$, these identities prove that
\eqref{i2:eq:propagator} is exactly the fundamental solution of the
linear part of \eqref{i2:eq:error-SDE}. There is no transpose on its
right factor. Since $K(u)=e^{u-v}K(v)$,
\[
 D_L(v)^{-1}D_L(u)
 =e^{u-v}I+(1-e^{u-v})D_L(v)^{-1},
\]
and the analogous right identity holds. For $0\le v-u\le\theta$,
the stability lemma consequently gives
\[
 \norm{\cU(v,u)}
 \le C\left(1+\frac{v-u}{\gamma_v}\right)^2.
\]
The characteristic formulas below identify this bound with
\begin{equation}
 \norm{\cU(v,u)}\le C(\gamma_u/\gamma_v)^2.
 \label{i2:eq:prop-bound}
\end{equation}
To express this estimate in the actual height variables, solve the
characteristics explicitly. Denote the terminal MDE parameters by
$(z_i,\eta_i,v_i,u_i)$.
In backward time $r=\theta-s$, the exact formulas are
\begin{equation}
 \begin{aligned}
 z_{i,r}&=e^{r/2}z_i,&
 v_{i,r}&=e^{-r/2}v_i,&
 u_{i,r}&=e^{-r}u_i,\\
 \eta_{i,r}&=e^{r/2}\eta_i+2\sinh(r/2)v_i,&
 M_{i,r}&=e^{-r/2}M_i.
 \end{aligned}
 \label{i2:eq:exact-backward-characteristics}
\end{equation}
They preserve both scalar MDE identities, have positive height, and
differentiate to \eqref{i2:eq:characteristic} after reversing time.
The uniqueness of the positive-imaginary MDE solution identifies
them with the given characteristic. Since $\theta\asymp\sigma\to0$
and $K\Subset K'$, the spatial points remain in $K'$.
The lower and upper bounds on $v_i$ at bounded terminal heights give
\begin{equation}
 x_r\asymp x+r,\qquad y_r\asymp y+r,\qquad
 \gamma_r\asymp\gamma_0+r
 \label{i2:eq:back-heights}
\end{equation}
for bounded terminal heights. More explicitly,
\[
 \gamma_r=e^r|z_1-z_2|^2
       +e^{r/2}(\eta_1+\eta_2)
       +2\sinh(r/2)(v_1+v_2).
\]
Applying this identity between any two times shows
$\gamma_u\asymp\gamma_v+(v-u)$ and proves
\eqref{i2:eq:prop-bound}. In particular, if $x\le y$, then
$\gamma_0\ge y$ and hence
\begin{equation}
 \norm{\cU(\theta,\theta-r)}
 \le C((y+r)/y)^2.
 \label{i2:eq:unequal-prop}
\end{equation}
When $y$ is bounded below, the propagator is bounded instead.

\subsection{The initial estimate at the Gaussian age}

Next we are able to obtain the starting estimate at the Gaussian smoothing scale and extend it from coincident to nearby spatial points.
In the fine range $x\le\sigma$, the smaller height at the initial
end of a backward characteristic is comparable to $\sigma$. We first use the
coincident-point estimate there, then interpolate in the real spatial
variables to cover nearby points. The interpolation error has a
negative power of $N$; after transport, the initial bound is a
vanishing fraction of the required logarithmic estimate.

\begin{lemma}\label{i2:lem:seed}
Suppose $t^{24}/N\le x\le\sigma$, $y\ge x$, and
$|z-w|\le C_0N^{-1/4}$. At the initial end of the characteristics,
\begin{equation}
 \norm{E_0}\le\frac{C_D t^{10}}{N\sigma(y+\sigma)}
 \label{i2:eq:seed}
\end{equation}
with arbitrarily fast polynomial decay of the failure probability,
uniformly over each fixed polynomial parameter range needed below.
\end{lemma}
\begin{proof}
The initial heights satisfy $x_0\asymp\sigma$ and
$y_0\asymp y+\sigma$. The initial points remain in $K'$ and have
distance at most $CN^{-1/4}$. At coincident points, external input
item~\ref{i2:input:coincident} yields $Ct^{10}/(Nx_0y_0)$, since a fixed
$\epsilon<1/200$ puts $x_0$ in its allowed upper-height range.
Both initial heights exceed the logarithmic lower cutoff. If the
two initial heights need to be interchanged, cyclicity interchanges
the two resolvents and their insertions at the coincident reference
point; the smaller height still satisfies the required upper bound.

Keep $z_0,x_0,y_0$ fixed and vary the second point. For one coordinate
of the matrix error define
\[
 f(w)=Nx_0y_0\,\mathcal E_{z_0,w}(x_0,y_0;C_a,C_b).
\]
Work in a fixed small bulk neighborhood with a fixed larger
neighborhood still contained in $K'$. On a common fixed-power
local-law event, item~\ref{i2:input:polynomial} gives $\norm f_\infty\le CN^\epsilon$ on
the larger neighborhood. Here B.2(i) covers the small-height
range, and its overlap with B.1 covers $y_0$ near or above one,
up to any polynomial upper bound fixed in advance.

All spatial derivatives below are real derivatives in either
$\operatorname{Re}w$ or $\operatorname{Im}w$; no holomorphy in
$w$ is asserted. For every fixed integer $m\ge1$, deterministic
resolvent differentiation and Lemma~\ref{i2:lem:stability} give
\[
 \norm{\partial^m f}_\infty\le C_mN\sigma^{-m}.
\]
Here are details for both terms. Since the real spatial derivative
of the Hermitization is a constant norm-one matrix,
$\norm{\partial^m G_w(iy_0)}\le m!y_0^{-m-1}$. Thus
\[
 \left|\partial^m\bigl(Nx_0y_0
       \br{G_{z_0}(ix_0)C_aG_w(iy_0)C_b}\bigr)\right|
 \le C_mN y_0^{-m}\le C_mN\sigma^{-m}.
\]
For the deterministic term, $M$ and each of its fixed-order real
spatial derivatives are bounded uniformly on bounded-height bulk
domains, including at height zero. Indeed, its scalar parameter
$u$ solves
\[
 \eta=(1-u)\sqrt{\frac{1-|z|^2u}{u}},\qquad 0<u\le1.
\]
On such a domain the derivative of the right side in $u$ is
strictly negative and bounded away from zero, so implicit
differentiation is uniform. Differentiating the denominator
inverse in \eqref{i2:eq:Fbar} therefore gives
\[
 |\partial^m\overline F_{ab}|
 \le C_m\gamma^{-m-1},\qquad
 \gamma=x_0+y_0+|z_0-w|^2.
\]
The inequalities $x_0\asymp\sigma$ and
$y_0\ge c\sigma$ now imply
$Nx_0y_0|\partial^m\overline F_{ab}|
 \le C_mN\sigma^{-m}$.
For $y_0$ bounded below and possibly polynomially large, the
denominator inverse is uniformly bounded instead. The explicit
MDE formula gives $M_w(iy_0)=O(y_0^{-1})$ and the same bound for
each fixed-order real spatial derivative. Differentiation of
\eqref{i2:eq:Fbar} then gives
$|\partial^m\overline F_{ab}|\le C_m/y_0$, which is also
sufficient after multiplication by $Nx_0y_0$. These derivative
estimates are deterministic.

The spatial interpolation uses the finite-difference inequality
\[
 \norm{\partial f}_\infty
 \le C_m\left(h^{-1}\norm f_\infty
           +h^{m-1}\norm{\partial^m f}_\infty\right).
\]
It follows by a fixed $m$-point difference formula whose Taylor
coefficients through order $m-1$ cancel except for the first
derivative. The larger neighborhood provides all its evaluation
points for sufficiently small $h$. Taking
$h=\sigma N^{-(1-\epsilon)/m}$ gives
\[
 \norm{\partial f}_\infty
 \le C_m\sigma^{-1}N^{\epsilon+(1-\epsilon)/m}.
\]
Apply this separately in the two real spatial directions. With
$\epsilon=1/1000$ and $m=32$,
\begin{equation}
 |f(w_0)-f(z_0)|
 \le CN^{-1/4+1/200+\epsilon+(1-\epsilon)/m}
 =CN^{-0.21278125}.
 \label{i2:eq:seed-fd}
\end{equation}
This is negligible relative to $t^{10}$. To obtain the common
event used above, take a sufficiently fine fixed polynomial mesh
of the spatial variables and of the other parameters in the fixed
polynomial range. Their deterministic resolvent derivatives are
polynomially bounded because the initial heights are at least
$c\sigma$. Choose the external failure exponent larger than the
mesh exponent and the desired final failure exponent. A union
bound and deterministic interpolation then give the claimed
continuum event. The same argument applies to the coincident-point
input. This proves \eqref{i2:eq:seed}.
\end{proof}

The transported initial estimate is absorbed at the desired scale, because
\[
 \norm{\cU(\theta,0)E_0}
 \le\frac{C_Dt^{10}(y+\sigma)}{N\sigma y^2},
\]
and its ratio to $t^{16}/(Nxy)$ is at most
\[
 C_Dt^{-6}\frac{x(y+\sigma)}{\sigma y}\le2C_Dt^{-6}.
\]
Here $x\le\min\{\sigma,y\}$.

\subsection{Drift, bracket, and nonlinear stopping}

Next we control drift and random fluctuations until a possible first failure, then show that this failure is unlikely.
We first estimate the stopped equation at a fixed deterministic
terminal time. A common event on a polynomial mesh and continuity
then exclude a random first exit. Throughout, the local-law estimates
are applied to unconditioned matrices at deterministic times.

First consider the range where both one-resolvent estimates are
logarithmic. Fix a deterministic terminal time $\theta$ and its
deterministic characteristics. Write $s$ for forward time and
$r=\theta-s$ only when evaluating deterministic integrals. As above,
$x_r=\eta_{1,s}\asymp x+r$ and
$y_r=\eta_{2,s}\asymp y+r$, with $x\le y$.
Stop at the first failure of the adapted one-resolvent bounds or of
the proposed error bound; denote this stopping time by $\tau_*$. The
estimates below use only the inequalities
\[
 |\epsilon_i|\le C_Dt/(N\eta_i),\qquad
 \br{\operatorname{Im}G_i}\le C.
\]
They are imposed before the stop, rather than by conditioning on a
future event. Their simultaneous validity along the required
characteristics is addressed in the uniformization argument below.
Multiplication by $\one_{\{s\le\tau_*\}}$ gives predictable
stopped integrands.

\emph{Forcing and quadratic variation.}
Normality and Ward's identity give
\[
 G_i^*G_i=\frac{\operatorname{Im}G_i}{\eta_i},\qquad
 \norm{G_i}_{\HS}^2=\norm{G_i^2}_1
   =\Tr|G_i|^2\le CN/\eta_i,
 \qquad \norm{G_i}\le\eta_i^{-1}.
\]
Since every basis insertion has operator norm one, trace-norm
duality yields
\[
 |\br{G_1^2C_aG_2C_b}|
 \le\frac1{2N}\norm{G_1^2}_1\norm{C_aG_2C_b}
 \le\frac C{x_ry_r}.
\]
By cyclicity, the same bound holds with $G_2$ squared instead.
Consequently, before the stop,
\begin{equation}
 \norm{H_s}\le\frac{C_Dt}{N}
       (x_r^{-2}y_r^{-1}+x_r^{-1}y_r^{-2}).
 \label{i2:eq:forcing}
\end{equation}
The exact martingale coefficient is the sum of chains in
\eqref{i2:eq:martingale-coefficient}. For its first chain,
\[
 \norm{G_1C_aG_2C_bG_1}_{\HS}^2
 \le \norm{G_1}_{\HS}^2\norm{G_2}^2\norm{G_1}^2
 \le\frac{CN}{x_r^3y_r^2}.
\]
The second chain has squared Hilbert--Schmidt norm at most
$CN/(x_r^2y_r^3)$. Squaring their sum costs an absolute constant,
and retaining only its off-diagonal blocks decreases its squared
Hilbert--Schmidt norm. Thus \eqref{i2:eq:sesquilinear-bracket} gives
\begin{equation}
 \dd[M_{ab},\overline M_{ab}]_s
 \le\frac C{N^2}(x_r^{-3}y_r^{-2}+x_r^{-2}y_r^{-3})\dd s,
 \qquad r=\theta-s,
 \label{i2:eq:bracket}
\end{equation}
before the stop. The exact covariance factor is $1/(4N^3)$;
the chain estimates supply one factor $N$, leaving $N^{-2}$.
Cross-coordinate brackets are bounded by Cauchy--Schwarz. There are
only sixteen complex coordinates.

\emph{Transport to a fixed terminal time.}
Using \eqref{i2:eq:unequal-prop}, the stopped transported drift is at most
\begin{align}
 \int_0^\theta\norm{\cU(\theta,s)}
       \one_{\{s\le\tau_*\}}\norm{H_s}\dd s
 &\le\frac{C_Dt}{Ny^2}\int_0^\theta
  \left(\frac{y+r}{(x+r)^2}+\frac1{x+r}\right)\dd r\nonumber\\
 &\le\frac{C_Dt^2}{Nxy}.
 \label{i2:eq:drift-integral}
\end{align}
For clarity, put $d=y-x\ge0$ and $L=\log(1+\theta/x)$.
The integral in the first line is exactly
\[
 2L+d\left(\frac1x-\frac1{x+\theta}\right)
 \le2L+\frac yx.
\]
The lower-height floor and $\theta\le C$ give
$1+L\le C\log N\le Ct$. Since $y/x\ge1$, this proves the
second line of \eqref{i2:eq:drift-integral}. Its ratio to
$t^{16}/(Nxy)$ is $O_D(t^{-14})$.

For this fixed deterministic $\theta$, define the forward-time
martingale
\[
 \mathcal N_\theta(u)=\int_0^u\cU(\theta,s)
       \one_{\{s\le\tau_*\}}\dd M_s,
 \qquad 0\le u\le\theta.
\]
The kernel's terminal argument is fixed. We do not assert that
$\theta\mapsto\mathcal N_\theta(\theta)$ is a martingale, or
replace that terminal argument by a random time in the following
calculation. Coordinate Cauchy--Schwarz, the fixed dimension,
\eqref{i2:eq:bracket}, and the squared propagator bound show that each
coordinate's sesquilinear bracket at time $\theta$ is at most
\begin{align}
 &\frac C{N^2y^4}\int_0^\theta
 \left(\frac{(y+r)^2}{(x+r)^3}+
       \frac{y+r}{(x+r)^2}\right)\dd r\nonumber\\
 &\hspace{25mm}\le\frac{C\log N}{N^2x^2y^2}.
 \label{i2:eq:bracket-integral}
\end{align}
Indeed, setting $q=x+r$ turns the integrand into
$2/q+3d/q^2+d^2/q^3$, whose integral is
\begin{align*}
 2L+3d\left(\frac1x-\frac1{x+\theta}\right)
 +\frac{d^2}{2}\left(\frac1{x^2}-\frac1{(x+\theta)^2}\right)
 &\le 2L+3\frac yx+\frac{y^2}{2x^2}\\
 &\le C(1+L)\frac{y^2}{x^2}.
\end{align*}
Only deterministic integrals have been rewritten in backward time.
Each real or imaginary coordinate of $\mathcal N_\theta$ has
quadratic variation bounded by \eqref{i2:eq:bracket-integral}.
The continuous exponential martingale inequality, followed by a
finite union over signs and coordinates, therefore gives, for every
fixed $\delta>0$,
\begin{equation}
 \PP\left(\norm{\mathcal N_\theta(\theta)}
      >\frac{\delta t^{16}}{Nxy}\right)
 \le C\exp(-c_\delta t^{31}).
 \label{i2:eq:fixed-terminal-martingale-tail}
\end{equation}
The exponent follows from $t^{32}/\log N\asymp t^{31}$.
This bound is smaller than every prescribed $N^{-D_*}$ for
sufficiently large $N$. It is an unconditional estimate for the
stopped process; the probability of reaching the local-law stop
still has to be controlled separately.

\emph{The nonlinear improvement.}
Let the nonlinear stopping threshold be
$K t^{16}/(Nx_ry_r)$ with fixed $K$. Since
$\norm{EJE}\le C\norm E^2$, its contribution before the stop is at most
\begin{equation}
 \int_0^\theta\norm{\cU(\theta,s)}
       \one_{\{s\le\tau_*\}}\norm{E_sJE_s}\dd s
 \le\frac{CK^2t^{32}}{N^2y^2}\int_0^\theta\frac{\dd r}{(x+r)^2}
 \le\frac{CK^2t^{32}}{N^2xy^2}.
 \label{i2:eq:quadratic}
\end{equation}
Its ratio to $t^{16}/(Nxy)$ is
$CK^2t^{16}/(Ny)\le CK^2t^{-8}$ under the lower-height floor.
For every fixed $K$, it is therefore a vanishing fraction of the
stopping threshold. Together with the transported initial estimate, the
forcing estimate and \eqref{i2:eq:fixed-terminal-martingale-tail},
these bounds give a strict improvement of the stopped extension
at each deterministic terminal time, outside the stated martingale
exception. The same calculation applies with $\theta$ replaced
by any fixed earlier terminal time whose characteristic obeys the
same bounds. Exclusion of a random first exit requires the
uniformization argument that follows. The logarithmic exponents have been
chosen in the finite order
\[
 p_0=10,\quad p_1=1,\quad p=16,\quad b_0=24.
\]
No later history exponent enters this choice.

\begin{lemma}[Uniform estimates and the first exit]\label{i2:lem:uniformization}
The fine stopped bounds above and the coarse stopped bounds established
below, together with the external local laws, hold on one event of
overwhelming probability for all the
spatial points, heights, and observation times required in
Theorem~\ref{i2:thm:I2}. The strict stopped improvement rules out the
error exit on that event.
\end{lemma}
\begin{proof}
We give explicit polynomial meshes. Realize the OU process on
$[0,1]$ only to construct a common modulus event. Observation ages
are restricted to $[c_{\rm age}\sigma,C_{\rm age}\sigma]$; the
lemma for the initial estimate is used only in that age range. Let the upper
spectral height be $N^C$, with fixed $C\ge0$. It suffices to use
$C=10$ for the integrated application. Prove the pointwise estimates
on a fixed-factor enlargement of the desired domain: lower floor
$t^{24}/(2N)$, upper height $2N^C$, a slightly enlarged bulk disk,
and twice the desired separation bound.

\emph{The stopped extension and intermediate observation times.}
For each deterministic characteristic put
$T_v=t^{16}/(N\eta_{1,v}\eta_{2,v})$. Use the stopping rule above,
with nonlinear threshold $T_v$ (that is, choose $K=1$ in the
preceding estimates) and with fixed slack in the local-law bounds.
Its first failure time is $\tau_*$. The stopped extension is
\begin{align*}
 \widehat E_v={}&\cU(v,0)E_0
 +\int_0^v\cU(v,u)\one_{\{u<\tau_*\}}(E_uJE_u+H_u)\dd u\\
 &+\int_0^v\cU(v,u)\one_{\{u<\tau_*\}}\dd M_u,
\end{align*}
It agrees with $E_v$
up to $\tau_*$ and evolves linearly afterwards. For fixed deterministic
$v$ its stochastic convolution has the deterministic terminal kernel
needed for \eqref{i2:eq:bracket-integral}; it is not asserted to be a
martingale as a function of $v$. Nor is this expression a representation
of $E_{v\wedge\tau_*}$.

The initial estimate also gives the required improvement at intermediate
observation times $v$, including $v$ arbitrarily close to zero.
Indeed, write $x_v=\eta_{1,v}$ and $y_v=\eta_{2,v}$, with
$x_v\le C y_v$. The characteristic formulas give
$x_0\asymp x_v+v$ and $y_0\asymp y_v+v$. Combining the original
initial estimate at time zero with the propagator to time $v$ gives
\[
 \frac{\norm{\cU(v,0)E_0}}{T_v}
 \le C_Dt^{-6}\frac{x_vy_0}{x_0y_v}
 \le C_D't^{-6}.
\]
The same calculation holds at $v=0$. The stopped forcing,
quadratic, and bracket integrals have the same bounds with $v$
as their terminal observation time and $x_v,y_v$ as their terminal
heights. This uses the original initial estimate; it does not apply an
aged-matrix theorem at an intermediate age or require
$v\asymp\sigma$.

\emph{A common event on the deterministic mesh.}
The resolvent identity, for heights at least $N^{-1}$, gives
\begin{align}
 &|\mathcal E(s,z,w;x,y)-\mathcal E(s',z',w';x',y')|\nonumber\\
 &\quad\le CN^3\big(\norm{X_s-X_{s'}}+|z-z'|+|w-w'|
                         +|x-x'|+|y-y'|\big).
 \label{i2:eq:det-modulus}
\end{align}
This holds for every basis insertion. For the random trace it follows
by replacing one resolvent at a time. For the deterministic reference,
the explicit denominator formula gives first derivatives $O(N^2)$;
the stated $O(N^3)$ envelope is more than sufficient.

For a time mesh $\delta=N^{-L}$, there is an event of probability
at least $1-C_{D,L}N^{-D}$ on which
\begin{equation}
 \sup_{s\le1}\norm{X_s}_F\le C_DN,\qquad
 \sup_{|s-r|\le\delta}\norm{X_s-X_r}
       \le C_{D,L}N^{1-L/2}.
 \label{i2:eq:OU-modulus}
\end{equation}
Indeed $\E\norm{X_0}_F^{2q}\le C_qN^q$, so Markov's inequality
at level $N$ gives the stated bound on the failure probability by choosing
$q$ large. The explicit OU solution and the entrywise Brownian
supremum bound on $[0,1]$ give the same $O(N)$ Frobenius bound
uniformly in time. The drift increment is then at most $CN\delta$. The Brownian
reflection bound, followed by a union over $N^2$ entries and
$O(N^L)$ time cells, bounds each unnormalized Brownian increment by
$C_{D,L}\sqrt{\delta\log N}$. The normalized Frobenius increment is
therefore at most $CN^{1/2}\sqrt{\delta\log N}$, which is covered
by \eqref{i2:eq:OU-modulus}. Pairs of nearby times lie in at most two
adjacent cells.

Choose spatial and spectral mesh spacing $h=N^{-K_m}$ and time
mesh spacing $\delta=N^{-L_m}$, where
\[
 K_m=12+4C,\qquad L_m=32+8C.
\]
Include both observation-age endpoints and zero and the endpoint of
each characteristic. The four real spatial coordinates, two heights,
observation time and characteristic terminal time require at most
$CN^Q$ tests, with
\[
 Q=6K_m+2C+2L_m=136+42C.
\]
This is a deliberate upper bound. Choose every external local-law
failure exponent larger than $D+Q+20$ before taking this finite union.
Their logarithmic exponents remain $p_0=10$ and $p_1=1$.
The martingale failure probability $C\exp(-ct^{31})$ also survives
this union. Thus, after assigning a sufficiently small fixed fraction
of the threshold to the martingale, the strict deterministic
improvements already proved imply
\begin{equation}
 \norm{\widehat E_{v_j}}\le T_{v_j}/4
 \label{i2:eq:mesh-improvement}
\end{equation}
at every deterministic time mesh point, for sufficiently large $N$.
The one-point local laws here are applied only at deterministic
flow times and deterministic characteristic parameters. At those
times the unconditioned OU matrix has the required uniform moment
bounds. Take each local-law stopping threshold to be twice the
corresponding bound on the common event. The one-resolvent identity
gives the deterministic modulus $CN^2$ for its normalized trace
when heights are at least $N^{-1}$; the deterministic $m$ has
bounded first derivatives in the bulk. Equations
\eqref{i2:eq:OU-modulus} and the chosen meshes extend these laws to
the continuum with errors smaller than the fixed slack. The
imaginary-trace bound follows as well: in the logarithmic range,
$t/(N\eta)\le Ct^{-23}$ at the floor, and in the coarse range
the fixed-power one-point error tends to zero uniformly at and
above its cutoff. Thus the local-law part of
the adapted stop is not reached on the common event. No law is
applied after conditioning on that event.

\emph{Excluding the first exit and interpolating the parameters.}
Suppose a nonlinear first crossing occurs at $\tau_*$. Let $v_j$
be the last mesh time before it. Since $v_j\le\tau_*$,
$\widehat E_{v_j}=E_{v_j}$. Along a characteristic,
\[
 |z'_v|+|w'_v|\le C,\qquad
 |\eta'_{1,v}|+|\eta'_{2,v}|\le C(1+N^C).
\]
Equations \eqref{i2:eq:det-modulus} and \eqref{i2:eq:OU-modulus} give
\[
 \norm{E_{\tau_*}-E_{v_j}}
 \le C\big(N^{4-L_m/2}+N^{C+3-L_m}\big)
 =o(t^{16}N^{-1-2C}).
\]
Moreover $|\partial_v\log T_v|\le CN$, and hence
$T_{v_j}/T_{\tau_*}=1+o(1)$. Since
$T_{\tau_*}\ge ct^{16}N^{-1-2C}$, the last display and
\eqref{i2:eq:mesh-improvement} contradict
$\norm{E_{\tau_*}}=T_{\tau_*}$. If the crossing lies in the first
time cell, use the improvement from the initial estimate at $v_j=0$. This proves that
the nonlinear stop is not reached.

Finally, nearest-mesh interpolation in endpoint parameters has
absolute cost at most
\[
 C\big(N^{3-K_m}+N^{4-L_m/2}\big)
 =o(t^{16}N^{-1-2C}).
\]
The neighboring heights differ relatively by $1+o(1)$ because they
are at least $N^{-1}$. The enlarged bulk, height and separation
ranges cover all boundary choices. This proves the continuous
parameter conclusion on the same event.

\emph{The coarse threshold.}
For the arbitrary-spatial-point coarse extension, use the pointwise
stopped improvement proved below with
\[
 L_v=T_v+\frac{N^{-c}}{1+\eta_{2,v}},\qquad c=1/10,
\]
in place of $T_v$ throughout the stopping argument. The second
height is bounded below along the whole backward characteristic.
Both summands are positive, $L_v\ge T_v$, and
$|\partial_v\log L_v|\le CN$. Thus every modulus and threshold
comparison above remains valid. In this coarse case the ordinary
fixed-power initial estimate at the original time zero is
$O(N^\epsilon/(N\sigma y_0))$. The propagator is bounded, and,
for the fixed sufficiently small $\epsilon$ used below, this is
$o(N^{-c}/(1+y_v))$ uniformly in $v$; hence the intermediate-time
improvement from the initial estimate remains valid too. The proof below supplies the
other coarse stopped contributions before this uniformization is
invoked. No uncountable union, random
selection of a parameter before applying a pointwise law, or
terminal-time martingale assertion is used.
\end{proof}

\begin{corollary}[Logarithmic error bound along the whole characteristic]
\label{i2:cor:full-characteristic}
Fix the age constants and the bulk and near-pair separation
constants above. For every $D>0$, there is one event with complement
of probability at most $C_DN^{-D}$ such that, for every terminal
age $\theta\asymp\sigma$, every admissible terminal pair, and
$t^{24}/N\le x,y\le\sigma$, the whole characteristic
$0\le v\le\theta$ satisfies
\[
 \norm{E_v}\le\frac{C_Dt^{16}}{N\eta_{1,v}\eta_{2,v}},\qquad
 |\epsilon_{i,v}|\le\frac{C_Dt}{N\eta_{i,v}},\qquad
 \br{\operatorname{Im}G_{i,v}}\le C.
\]
The same statement holds for any fixed positive multiple of the
lower-height floor, with constants changed but with the same
logarithmic exponents.
\end{corollary}
\begin{proof}
Along these characteristics both heights are at most $C\sigma$,
so the logarithmic one-resolvent input applies throughout. The
two labels may first be interchanged to order their terminal heights.
The fine first-exit argument just proved controls every intermediate
time, using its inherited initial estimate at time zero. Its continuum
interpolation is in the original terminal parameters and the
intermediate observation time, so the conclusion includes the
initial segment of each characteristic. A fixed multiple of the
floor changes only fixed constants in the stopped estimates.
\end{proof}

\subsection{Coarse heights, arbitrary spatial points, and spectral tails}

The purpose of this subsection is to extend the resolvent bounds to larger heights and arbitrary bulk points, retaining the decay needed for integration.
Note that the remaining ranges have a bounded propagator. For nearby spatial
points with smaller height at most $\sigma$, we retain the
logarithmic error without an additive polynomial remainder; for arbitrary points we allow the stated
additive remainder. At large heights the error
must also retain its inverse-height decay, which will make the
integrated insertion estimate finite.

\emph{Nearby spatial points with a coarse second height.}
Keep $|z-w|\le C_0N^{-1/4}$. Choose $\xi=1/10000$ and $a=N^{-10\xi}$. Then $\sigma\ll a$.
If $x\le\sigma$ and $y\ge a$, the propagator is bounded over the
aging interval. Indeed, for $a\le y\le1$ the ratio bound gives
$\norm{\cU}\le C(1+\theta/y)^2\le C$, and for $y\ge1$ the
uniform coarse-height denominator bound gives the same conclusion.
The exact backward characteristic also gives $y_r\asymp y$ uniformly
over $0\le r\le\theta$, even on a fixed polynomial range of $y$.
The fine first-height law remains logarithmic. Its forcing contributes
at most $C_Dt/(Nx_r^2y)$, whose integral is at most
$C_Dt/(Nxy)$. The coarse one-point error contributes at most
\[
 \frac{CN^\xi}{Nx_ry^2}
\]
to the forcing. Its integrated size relative to $1/(Nxy)$ is
\begin{align}
 N^\xi\frac{x}{y}\log(1+\theta/x)
 &\le N^\xi\frac\theta y\nonumber\\
 &\le CN^{11\xi-1/200}=CN^{-0.0039}.
 \label{i2:eq:coarse-relative}
\end{align}
This is the needed relative bound; an absolute polynomial error
alone would not complete the stopping argument with only the logarithmic error. The initial estimate
is at most $C_Dt^{10}/(N\sigma y)$, whose ratio to
$t^{16}/(Nxy)$ is at most $C_Dt^{-6}x/\sigma$.
The transported bracket is bounded by
\[
 \frac C{N^2}\left(
 y^{-2}\int_0^\theta x_r^{-3}\dd r
 +y^{-3}\int_0^\theta x_r^{-2}\dd r\right)
 \le\frac C{N^2x^2y^2},
\]
using $x\le y$. Thus the same martingale and nonlinear stopping
argument applies, with a strict improvement for all sufficiently
large $N$. This proves the logarithmic bound without an additive polynomial remainder in the remaining
range $y\ge a$ on every fixed polynomial height range.
If the smaller terminal height exceeds $\sigma$, the direct
fixed-power two-resolvent law gives an absolute error at most
$N^{-1+2/200+\epsilon}$ for $x,y\le1$. This is the source of the
additive $N^{-c}$ term in \eqref{i2:eq:I2}. The ordering of the two
heights can be exchanged by cyclicity of the trace and the
deterministic two-resolvent equation, with $C,D$ exchanged as well.

\emph{Arbitrary spatial points.}
For arbitrary spatial points with $y\ge\eta_c>0$, the propagator is
bounded and a fixed-power initial estimate gives
\[
 \norm{E_0}\le CN^{-1+1/200+\epsilon}/y.
\]
Use the threshold
\[
 L_r=\frac{t^{16}}{Nx_ry_r}+\frac{N^{-c}}{y_r},
\]
where $c>0$ is fixed and small. The square of the first term is
controlled by \eqref{i2:eq:quadratic}. The cross term, relative to the
additive threshold, is at most $Ct^{16}\log N/(Ny)$, while the
square of the additive term has relative size $C\theta N^{-c}/y$.
More explicitly, writing $A=t^{16}$ and $\delta_N=N^{-c}$, the
three integrated nonlinear contributions are bounded by
\[
 C\left(\frac{A^2}{N^2xy^2}
       +\frac{A\delta_N\log(1+\theta/x)}{Ny^2}
       +\frac{\theta\delta_N^2}{y^2}\right).
\]
Divide the first term by $A/(Nxy)$ and the remaining terms by
$\delta_N/y$. All three ratios tend to zero since $y\ge\eta_c$.
For definiteness, fix the external fixed-power loss
$\epsilon_0=1/1000$ and take $c=1/10$ in these extensions.
The exponent in the initial estimate, $-1+1/200+\epsilon_0=-0.994$, and the preceding
coarse/coarse exponent $-1+2/200+\epsilon_0=-0.989$ are both
strictly smaller than $-c$. The initial estimate is therefore a vanishing
fraction of the additive threshold, uniformly in the stated ranges.
The forcing calculation \eqref{i2:eq:coarse-relative} is
even stronger here. This proves \eqref{i2:eq:coarse-I2}, first on fixed
polynomial height ranges. It gives the estimate directly in the Gaussian-perturbed model for
arbitrary spatial points when the second spectral height is bounded
below.
For $x>\sigma$, apply B.2(i) or B.1 directly at the final matrix:
$N^\epsilon/(Nxy)\le N^{-1+1/200+\epsilon}/y$, which is absorbed
by the additive $N^{-c}/y$ term. This also handles bounded coarse
heights just above one: choose the fixed constant in B.1 larger
than $2$ and $\eta_c^{-1}$. Since $y\ge\eta_c$, the factors
$y^{-1}$ and $(1+y)^{-1}$ are comparable, with constants depending
only on $\eta_c$, which gives the stated form of
\eqref{i2:eq:coarse-I2}.

\emph{Large spectral heights.}
For $y\ge1$, Ward improves to $\norm{G_2^2}_1\le CN/y^2$.
The coarse forcing is $CN^\xi/(Nx_ry^3)$ and its relative
integrated cost is at most $CN^\xi\theta/y^2=o(1)$. The bracket is
bounded by
\[
 \frac C{N^2}(x_r^{-3}y^{-2}+x_r^{-2}y^{-4})\dd r.
\]
The large-height initial estimate has the required $y^{-1}$ decay:
use B.1 at arbitrary points, or Lemma~\ref{i2:lem:seed} for nearby
points, where the latter was derived from the coincident-point
input B.2(ii). Thus the
propagated error retains $y^{-1}$ decay, rather than merely an
absolute $N^{-c}$ remainder.

There is no need to extrapolate an external theorem to infinite
spectral height. Propagate first to $1\le y\le N^3$. With overwhelming
probability, $\norm{H^w}\le N$ uniformly
on the fixed bulk set; entry moment bounds and a union bound suffice.
For $y\ge N^3$,
\[
 G_w(iy)=\frac{i}{y}I+O(N/y^2),\qquad
 \cM_{z,x;w,y}[C]=\frac{i}{y}M_z(ix)C+O(1/y^2).
\]
The first remainder follows from the norm-convergent Neumann series,
since $\norm{H^w}/y\le N^{-2}$. For the second, the scalar MDE gives
$M_w(iy)=iI/y+O(y^{-2})$, while $M_z(ix)$ is uniformly bounded in
the bulk. Insert this into \eqref{i2:eq:M12} and use its bounded
large-$y$ denominator inverse. These remainder constants are uniform
for $t^{24}/N\le x\le2$.
Using the averaged one-resolvent law with insertion $CD$ gives
\[
 |\mathcal E_{z,w}(x,y;C,D)|
 \le\frac{Ct}{Nxy}+\frac{CN}{xy^2}+\frac C{y^2}
 \le\frac{Ct}{Nxy}
\]
for $x$ in the fine range. At coarse $x$, the corresponding
fixed-power one-resolvent law produces the permissible additive
$N^{-c}/y$. This completes the large-height extension.
Combining the initial estimate, stopped improvement, uniformization, and the
separate coarse-height cases proves Theorem~\ref{i2:thm:I2}.

\subsection{Integrated coefficient bound}

Next we shall integrate the resolvent error to obtain the coefficient bound needed for the weighted Gaussian comparison.
The signed comparison integrates one resolvent height. Its
deterministic insertion decays as $(1+u)^{-1}$, so the inverse-height
error just proved has an integrable majorant after insertion. We
keep the endpoint-height restrictions separate from the integration
range and use one common event for all coordinates.

Let a deterministic coordinate list $(z_i,\eta_i)$ consist of
fine coordinates at two points whose distance is at most
$C_0N^{-1/4}$, together with finitely many coarse coordinates at
arbitrary points of the fixed bulk compact set. Assume
\[
 \eta_{\rm mes}\le\eta_i\le1,\qquad
 \eta_{\rm mes}\ge t^{24}/N,
\]
and that the coarse heights belong to $[a_0/2,1]$, where
$0<a_0\le1$ is fixed. These are the endpoint-height ranges used by
the signed comparison. The integration variable below still runs
up to $N^{10}$.
Put
\[
 \mathcal A_i(u)=\frac{M_{z_i}(iu)}{1-\br{M_{z_i}(iu)^2}},
 \qquad \norm{\mathcal A_i(u)}_{\rm bl}\le\frac C{1+u}.
\]
The last deterministic bound follows from the bulk MDE and remains
valid at all positive spectral heights. To see this explicitly, write
$m=iv$, $\omega=v/(h+v)$, and $\varrho=|z|^2$ at height $h>0$.
Then $0<\omega<1$, $v^2=\omega-\varrho\omega^2$, and
\[
 1-\br{M^2}=1+\omega-2\varrho\omega^2
 \ge\min\{1,2(1-\varrho)\}>c_K>0.
\]
The inequality follows by concavity in $\omega$ and the fixed
bulk margin. Also $v\le\min\{1,h^{-1}\}$ and
$\omega\le\min\{1,h^{-2}\}$: the latter bounds follow from
$v^2\le\omega\le v/h$. The explicit entries of $M$ therefore
give $\norm M_{\rm bl}\le C/(1+h)$, proving the asserted bound
for $\mathcal A_i$ without any small-height singularity.

\begin{corollary}[Integrated insertion error]\label{i2:cor:I6}
For each orientation $p,q\in\{1,2\}$, define
\begin{align*}
 \mathfrak V_{ir}^{pq}
 &=\int_{\eta_i}^{N^{10}}
   \br{G_{z_i}(iu)\mathcal A_i(u)E_{pp}
             G_{z_r}(i\eta_r)E_{qq}}\dd u,\\
 V_{ir}^{pq}
 &=\int_{\eta_i}^{N^{10}}
   \br{\cM_{z_i,u;z_r,\eta_r}[\mathcal A_i(u)E_{pp}]
             E_{qq}}\dd u.
\end{align*}
With overwhelming probability,
\begin{equation}
 \sup_{i,r,p,q}|\mathfrak V_{ir}^{pq}-V_{ir}^{pq}|
 \le C_D\frac{t^{17}}{N\eta_{\rm mes}}+C_DN^{-c}.
 \label{i2:eq:I6}
\end{equation}
The same conclusion holds for the reversed insertion orientations
and fixed finite linear combinations used in the mixed
second-cumulant coefficient. No derivatives of $\mathcal A_i$
along the matrix flow are required.
\end{corollary}
\begin{proof}
We first record reciprocity of the deterministic reference, since
the integration may put the larger height in either position.
Write $K^{12}_{ab}=\br{M_1C_aM_2C_b}$ and define $K^{21}$ with
the two resolvents exchanged. Cyclicity gives
$K^{21}=(K^{12})^{\mathsf T}$, where this is the ordinary
transpose of the coordinate matrix. Since $J=J^{\mathsf T}$,
\eqref{i2:eq:Fbar} implies
\[
 \overline F^{21}=(\overline F^{12})^{\mathsf T}.
\]
The random coordinate matrices obey the same identity by trace
cyclicity. By bilinearity, the full error thus satisfies
\[
 \mathcal E_{z,w}(x,y;C,D)
 =\mathcal E_{w,z}(y,x;D,C)
\]
for all constant block insertions. No adjoint or complex conjugation
is involved in this exchange.

Expand the deterministic insertion $\mathcal A_i(u)E_{pp}$ in the
constant block basis after applying Theorem~\ref{i2:thm:I2}. Bilinearity
pays only its block norm. The portion $u\le1$ is empty when
$\eta_i=1$; otherwise its integrated error is at most
\[
 C_D\int_{\eta_i}^1
 \left(\frac{t^{16}}{Nu\eta_r}+N^{-c}\right)\dd u
 \le C_D\frac{t^{17}}{N\eta_r}+C_DN^{-c}.
\]
If one coordinate is coarse, the arbitrary-spatial extension applies;
trace cyclicity gives the required orientation when the larger
height occurs in the first position. For $u\ge1$, the norm bound
$\norm{\mathcal A_i(u)}_{\rm bl}\le C/u$ and the large-height
error give an integrable majorant
\[
 C_D\left(\frac{t^{16}}{N\eta_r u^2}+
                 \frac{N^{-c}}{u^2}\right).
\]
Its integral is at most $C_Dt^{16}/(N\eta_r)+C_DN^{-c}$.
Taking the supremum proves \eqref{i2:eq:I6}. The fixed finite basis and
the common event in Lemma~\ref{i2:lem:uniformization} justify
simultaneity, including for $O(t)$ endpoint coordinates: the bound
is a maximum, so no factor equal to the number of coordinates is
introduced here. Fixed finite linear combinations have coefficients
bounded independently of $N$. The argument uses only insertion values at each
spectral height, so no time derivative of the insertion is taken.
\end{proof}

\subsection{Exact second-cumulant normalization}

The purpose of this subsection is to identify the deterministic covariance coefficient and check its constant in the chosen trace and derivative conventions.

The integrated estimate reduces the second-cumulant term to a
deterministic coefficient. We first identify its primitive, then
check the entry derivative and both block orientations to fix the
normalization exactly. The fourth-cumulant term is a separate
calculation. The logarithmic covariance formula appears in
\cite[(2.33)]{CL}; see also the resolvent calculation leading to
\cite[(8.21)]{CL}.
We record the primitive in the present block and derivative conventions.

Fix spatial points and the second spectral parameter. Write
\[
 \Sigma_{ir}(u)=
 \sum_{p\ne q}
 \br{\cM_{z_i,u;z_r,\eta_r}[\mathcal A_i(u)E_{pp}]E_{qq}},
\]
where the sum contains exactly $(p,q)=(1,2),(2,1)$, with no weight.
Let $D_{ir}(u,\eta_r)$ be the active two-dimensional determinant from
Lemma~\ref{i2:lem:stability}. In endpoint notation it is
\[
 D_{ir}
 =1+(u_i u_r|z_i z_r|)^2-(m_i m_r)^2
                  -2u_i u_r\operatorname{Re}(z_i\overline z_r).
\]
Here $u_j=v_j/(\eta_j+v_j)$ and $m_j=iv_j$. This $u_j$ is a
deterministic MDE parameter, not the integration variable $u$.

\begin{lemma}\label{i2:lem:primitive}
For arbitrary spatial points in the fixed bulk and positive
cutoffs,
\begin{equation}
 -\frac i2\int_{\eta_i}^{\infty}\Sigma_{ir}(u)\dd u
       =-\frac14\log D_{ir}(\eta_i,\eta_r).
 \label{i2:eq:primitive}
\end{equation}
The logarithm is real, and the answer is symmetric in the two
spatial-height pairs.
\end{lemma}
\begin{proof}
Fix the spatial point, put $Z=zE_{12}+\overline zE_{21}$, and
differentiate with respect to the complex spectral variable $w$.
Since $M^{-1}=-(w+m)I-Z$ and $m=\br M$, differentiation gives
\[
 \partial_wM=(1+\partial_wm)M^2,
 \qquad
 1+\partial_wm=\frac1{1-\br{M^2}},
 \qquad
 \partial_wM=M\mathcal A.
\]
Only deterministic matrices are multiplied here. In the
four-coordinate notation of \eqref{i2:eq:Fbar}, write
$K'=\partial_{w_1}K$ and set
\[
 K'_{ab}=\br{M_1\mathcal A_1C_aM_2C_b},\qquad
 Q_{ab}=\br{\cM_{12}[\mathcal A_1C_a]C_b}.
\]
Applying \eqref{i2:eq:M12} to this input, including its self-energy
term, gives the matrix identity
\[
 Q=K'+QJK,
 \qquad Q=K'(I-JK)^{-1}.
\]
Here the inverse multiplies on the right. No interchange
of $K'$ with that inverse is asserted. The determinant of $I-JK$
is the active two-dimensional determinant $D$, because the other
two rows of $J$ vanish. Jacobi's determinant identity therefore gives
\begin{align*}
 \partial_{w_1}\log D
 &=-\operatorname{Tr}_4\big((I-JK)^{-1}JK'\big)\\
 &=-\operatorname{Tr}_4\big(JK'(I-JK)^{-1}\big)
   =-\operatorname{Tr}_4(JQ)\\
 &=-2(Q_{12}+Q_{21})=-2\Sigma.
\end{align*}
The middle equality uses only cyclicity of the finite-dimensional
trace. The factor two is precisely $J_{12}=J_{21}=2$; the definition
of $\Sigma$ itself has two unweighted orientations.
Along $w_1=iu$ this becomes
\[
 \partial_u\log D=-2i\Sigma(u).
\]
The determinant is positive by \eqref{i2:eq:positive-det}, with the
same conclusion at large height by continuation. It tends to one
as the first height tends to infinity. Integrating proves
\eqref{i2:eq:primitive}. The explicit determinant formula proves
symmetry, although the integral representation singles out one
coordinate.
\end{proof}

We next check that the random coefficient has the same factor. Use the
finite-integral field, with $R=N^{10}$,
\[
 Y_j=iN\int_{\eta_j}^{R}\br{G_j(iv)-M_j(iv)}\dd v.
\]
Its deterministic centering plays no role in entry derivatives.
For an oriented Hermitized entry derivative,
\begin{equation}
 \partial_{ba}Y_j
 =-\frac i2\int_{\eta_j}^{R}(G_j(iv)^2)_{ab}\dd v
 =\frac12\big(G_j(i\eta_j)-G_j(iR)\big)_{ab}.
 \label{i2:eq:endpoint-derivative}
\end{equation}
This uses $\partial_vG(iv)=iG(iv)^2$. The factor $1/2$ is the
product of $N$ with the normalized trace factor $(2N)^{-1}$.
Here the entry derivative treats one oriented Hermitized matrix
entry as the variable, without an additional $N^{-1/2}$ factor;
its variance is already $1/N$.

To check the remaining normalization explicitly, put $W=H^0$,
$P_i=G_i\mathcal A_i$ and
$\Delta G_r=G_r(i\eta_r)-G_r(iR)$. For an upper-block entry
$h=W_{\alpha\beta}$, with $\alpha$ in block $1$ and $\beta$ in
block $2$, the mixed covariance term differentiates the field
with respect to $W_{\beta\alpha}$. It contributes
$(P_i)_{\beta\alpha}(\Delta G_r)_{\alpha\beta}/2$.
There is also the reverse outer orientation. The outer field
identity, normalized trace, entry covariance and field derivative
therefore give, in this order, the scalar product
\[
 (-iN)\frac1{2N}\frac1N\frac12=-\frac{i}{4N}
\]
before summing the two oriented block index sets. Their sum satisfies
\[
 \frac1{2N}\sum_{\substack{\alpha,\beta\\
                 \text{in opposite blocks}}}
 (P_i)_{\beta\alpha}(\Delta G_r)_{\alpha\beta}
 =\sum_{p\ne q}\br{P_iE_{pp}\Delta G_rE_{qq}}.
\]
Thus the coefficient, for real field tests, is
\begin{align}
 \mathcal V_{ir}^{(2)}(X)
 =\operatorname{Re}\bigg[
 -\frac i2\int_{\eta_i}^{R}\sum_{p\ne q}
 \br{G_i(iu)\mathcal A_i(u)E_{pp}
       \big(G_r(i\eta_r)-G_r(iR)\big)E_{qq}}\dd u
 \bigg].
 \label{i2:eq:raw-second-coefficient}
\end{align}
Taking the real part introduces no factor two. One must not commute
$\mathcal A_i$ through the random resolvent; its displayed position
is the insertion covered by the theorem.

\begin{corollary}[Normalized second-cumulant coefficient]\label{i2:cor:second}
For the coordinate family of Corollary~\ref{i2:cor:I6},
\[
 \max_{i,r}\left|
 \mathcal V_{ir}^{(2)}(X)+\frac14\log D_{ir}
 \right|
 \le C_D\frac{t^{17}}{N\eta_{\rm mes}}+C_DN^{-c}
\]
on one event of overwhelming probability.
\end{corollary}
\begin{proof}
Apply Corollary~\ref{i2:cor:I6} to the two orientations in
\eqref{i2:eq:raw-second-coefficient}. The discarded $G_r(iR)$ endpoint
has absolute size at most
\[
 \frac C R\int_{\eta_i}^{R}\frac{\dd u}{u(1+u)}
       \le \frac{C\log N}{R}.
\]
The deterministic tail above $R$ is $O(R^{-1})$, since $M_i(iu)$
and $\mathcal A_i(u)$ are $O(u^{-1})$ and the deterministic
two-resolvent stability is bounded there. Lemma~\ref{i2:lem:primitive}
then identifies the coefficient exactly.
\end{proof}

The direct calculation \eqref{i2:eq:endpoint-derivative} fixes the
normalization in our trace and entry-derivative conventions.
It agrees with \cite[equations (8.21) and (2.33)]{CL}.
Multiplying each field by $\sqrt2$ doubles the second-cumulant
covariance contribution to $-\tfrac12\log D_{ir}$. A deterministic
mean subtraction does not change its entry derivatives. The value
of that mean correction and the fourth-cumulant addition require
their separate contraction calculations.

\paragraph{Use in the signed comparison.}
The insertion estimates and exact second coefficient above concern
the aged model. Section~\ref{chap:signed} combines them with the
other cumulant contractions. Relative replacement and signed
normalizer calibration are established separately in
Sections~\ref{chap:replacement} and~\ref{chap:calibration}; they
supply the corresponding conclusions for the target law and its
normalized weights.

\endgroup

\begingroup
\newcommand{\E}{\mathbb E}
\newcommand{\PP}{\mathbb P}
\newcommand{\Tr}{\operatorname{Tr}}
\newcommand{\HS}{\mathrm{HS}}
\newcommand{\dd}{\mathrm d}
\newcommand{\one}{\mathbf 1}
\newcommand{\br}[1]{\langle #1\rangle}
\newcommand{\norm}[1]{\lVert #1\rVert}
\section{Positive overlap and independent reference noises}\label{chap:noise}

In this section we shall use overlap bounds for singular vectors to construct independent reference noises close enough to the matrix noises for later comparisons.

Indeed, this section constructs independent reference Brownian motions for each
prescribed pair of spatial points satisfying the separation conditions
below. We first consider the noises driving the small positive singular
values. A positive resolvent overlap controls the off-diagonal block of
their covariance-density matrix. We then multiply the noise increments
by the inverse square root of this matrix, an operation called whitening.
The reference noises have exact Brownian laws, and their differences
from the original noises are controlled up to an adapted stopping time.
Up to that time, the matrix-valued quadratic variation of the difference
is of quadratic order in the covariance error. We introduce the signed
noises only after this positive-index construction, because the full
signed covariance matrix is singular.

Positive two-resolvent estimates and singular-vector overlap decay were
proved by Cipolloni--Erd\H{o}s--Xu
\cite[Theorem 3.5 and Corollary 3.6]{CEX}. In the bulk, their two-resolvent
theorem requires heights above the microscopic scale by a fixed power,
and their overlap corollary retains a fixed-power loss. Here the relative
overlap estimate is proved on a common time event down to the logarithmic
height cutoff for the stated nearby points and matrices with an
independent Gaussian component.

The overlap argument uses the full-characteristic estimate under the
unweighted matrix law in Corollary~\ref{i2:cor:full-characteristic}.
The reference-process construction then supplies the noises used in the
relaxation and local-factor sections. We also give an extension at fixed
positive separation, using the additional estimates under the unweighted
matrix law stated at the end of this section. No coverage of the
intervening separation range is asserted here.

\subsection{Model, geometry, and estimates under the unweighted law}

The purpose of this subsection is to specify the matrix flow and spatial separation range in which the reference noises can be constructed.
We first fix compact disks $\mathcal D\Subset\mathcal D'\Subset\{|z|<1\}$.
Let $\widehat X=N^{-1/2}(\widehat\xi_{ab})$ have complex i.i.d. entries
with $\E\widehat\xi=\E\widehat\xi^2=0$, $\E|\widehat\xi|^2=1$, and
uniform bounds on every fixed moment. Define
\[
 t=\tfrac12\log N,\qquad d=t^{1/3},\qquad
 \sigma\asymp N^{-1/200},\qquad \tau=H/N,\qquad H=t^b,
\]
where $b$ is any fixed exponent. More precisely, $\sigma$ is deterministic
and lies in $[c_{\rm age}N^{-1/200},\allowbreak C_{\rm age}N^{-1/200}]$ for fixed
positive constants; all conclusions are uniform over this interval.
The theorem for full block insertions already permits this range of
Gaussian variance fractions.
The proof below uses $-\log(1-\sigma)\asymp N^{-1/200}$, initial
characteristic heights comparable to $\sigma$, and $\tau=o(\sigma)$,
all uniformly on this interval. Let $G^{\rm age}$ be an independent
normalized complex Ginibre matrix and set
\begin{equation}
 X_0=\sqrt{1-\sigma}\,\widehat X+\sqrt\sigma\,G^{\rm age},
 \qquad X_s=X_0+N^{-1/2}B_s,\qquad 0\le s\le\tau.
 \label{noise:eq:lfn-additive-model}
\end{equation}
The complex Brownian entries obey
$\E[\dd B_{ab}\dd\overline B_{cd}]=\delta_{ac}\delta_{bd}\dd s$
and $\E[\dd B_{ab}\dd B_{cd}]=0$.
The required deterministic spatial points satisfy
\begin{equation}
 z,w\in\mathcal D,\qquad
 c_0N^{-1}e^d\le |z-w|^2,\qquad
 |z-w|\le C_0N^{-1/4}.
 \label{noise:eq:lfn-geometry}
\end{equation}
Constants $c_0,C_0>0$ are fixed. The variable $d$ in this section is
always $t^{1/3}$. If $j=\lfloor\log|z-w|^{-1}\rfloor$, then
$N|z-w|^2\asymp e^{2(t-j)}$; the lower bound $e^d$ is the uniform,
weaker estimate used below, not a redefinition $d=t-j$.

For any matrix $Y$, write
\[
 H^v(Y)=\begin{pmatrix}0&Y-vI\\(Y-vI)^*&0\end{pmatrix},
 \quad G_v(Y,i\eta)=(H^v(Y)-i\eta)^{-1},
 \quad \br R=(2N)^{-1}\Tr R.
\]
The symbols $E_{ab}$ are the constant two-by-two matrix units tensored
with $I_N$.

\paragraph{Exact conversion to the aged OU flow.}
Put $a_\sigma=-\log(1-\sigma)$ and $u=\log(1+s)$. An OU process
$Y_r$ with initial value $Y_0=\widehat X$ can be coupled so that
\begin{equation}
 \dd Y_r=N^{-1/2}\dd\mathcal B_r-\tfrac12Y_r\dd r,
 \quad Y_{a_\sigma}=X_0,
 \quad Y_{a_\sigma+u}=\frac{X_s}{\sqrt{1+s}}.
 \label{noise:eq:lfn-ou-conversion}
\end{equation}
Indeed the first OU interval produces exactly the Gaussian mixture
in \eqref{noise:eq:lfn-additive-model}; afterwards
$\dd\mathcal B_u=(1+s)^{-1/2}\dd B_s$ has entry variance $\dd u$.
At an observation time $s$, set $r_*=a_\sigma+\log(1+s)$ and
\[
 \widetilde v=\frac{v}{\sqrt{1+s}},\qquad
 \widetilde\eta=\frac\eta{\sqrt{1+s}}.
\]
Then
\begin{equation}
 G_v(X_s,i\eta)=(1+s)^{-1/2}
       G_{\widetilde v}(Y_{r_*},i\widetilde\eta).
 \label{noise:eq:lfn-resolvent-conversion}
\end{equation}
Positive singular vectors are unchanged by this positive scalar
normalization. The Brownian noises used in the whitening construction
below are in the additive time $s$; the OU time is used only to prove
the overlap bound.
All martingale estimates on the aging interval use its natural
filtration generated by $\widehat X$ and the OU Brownian motion up to
the current aging time. This filtration is distinct from the later
additive-time filtration, whose time-zero sigma-field contains $X_0$.
In particular, the aging Brownian motion is never treated as a
martingale after its endpoint $X_0$ has been revealed.

\paragraph{Estimates along the entire characteristic.}
For each terminal pair and terminal heights of common order $a$,
$t^{b_0}/N\le a\le\sigma$, use the backward characteristics on
$[0,r_*]$ given by
\begin{equation}
 \dot z_i=-z_i/2,\qquad
 \dot\eta_i=-v_i-\eta_i/2,\qquad m_i=iv_i.
 \label{noise:eq:lfn-characteristic}
\end{equation}
Their terminal spatial points and heights are the normalized ones
in \eqref{noise:eq:lfn-resolvent-conversion}. Since $r_*\asymp\sigma$,
the deterministic bulk bounds imply
\begin{equation}
 \eta_{i,r}\asymp a+r_*-r=:a_r,
 \qquad \eta_{i,0}\asymp\sigma,
 \qquad |z_1(r)-z_2(r)|\le C N^{-1/4}.
 \label{noise:eq:lfn-characteristic-scales}
\end{equation}
The stopped propagation argument for full block insertions starts
with the coincident-point logarithmic estimate at height comparable
to $\sigma$. It gives one event on which, along the entire characteristic,
\begin{equation}
 |o_r-\overline o_r|\le\frac{C_Dt^p}{Na_r^2},\qquad
 |\epsilon_{i,r}|\le\frac{C_Dt^{p_1}}{Na_r},\qquad
 \br{\operatorname{Im}G_{i,r}}\le C,
 \label{noise:eq:lfn-ordinary-input}
\end{equation}
where $\epsilon_i=\br{G_i}-m_i$, $p_1\le p$, and $o,\overline o$
are defined below. One may take $p=16$, $b_0=24$ and $p_1=1$.
The same proof permits any fixed positive multiple of the floor
$t^{b_0}/N$: all stopping ratios change only by fixed constants.
We use this fixed-factor extension for the normalized terminal height
$a/\sqrt{1+s}$ and for meshes at admissible boundaries.
The probability of the complement is at most $C_DN^{-D}$, for every
fixed $D>0$, after prescribed polynomial meshes are included.
This estimate has no additive polynomial remainder along the whole
aging characteristic; a terminal estimate for already-aged snapshots
alone would not be a sufficient input. In particular, the initial segment is supplied by
the same stopped proof, rather than by applying an aged-matrix theorem
to an unaged intermediate matrix. Corollary~\ref{i2:cor:full-characteristic} gives precisely this assertion,
including intermediate times, fixed multiples of the lower-height floor,
and ages comparable to $\sigma$. Its coordinate bound applies to
$o=2F_{21}$ with only a fixed factor in the constant.

\subsection{The positive scalar equation}

In this section we shall see that positivity lets us control resolvent overlap relative to its deterministic value, even when an absolute error estimate would be too weak.
The useful scalar is positive on every sample path. Dividing it by its
positive deterministic counterpart allows the forcing and quadratic
variation to be controlled relatively, even when an absolute two-resolvent
error would be too large. We work first in the OU convention along one
fixed deterministic characteristic. Suppress its time variable, and put
\[
 A_i=\operatorname{Im}G_i,\qquad
 F_{ab}=\br{G_1E_{aa}G_2E_{bb}}\quad(a,b=1,2),\qquad
 J=\begin{pmatrix}0&2\\2&0\end{pmatrix}.
\]
Each $A_i$ is strictly positive and block diagonal. Spectral calculus
gives $G_i=iA_i+B_i$, with $B_i$ Hermitian and block off-diagonal.
Consequently
\begin{equation}
 Q:=\br{A_1A_2}=-\Tr F>0,\qquad
 F_{12}=\overline{F_{21}},\qquad o=2F_{21}.
 \label{noise:eq:lfn-positive-identities}
\end{equation}
Here $\Tr F$ is the ordinary trace of the two-by-two coordinate matrix.

The exact resolvent equation from the OU characteristic is
\[
 \dd G_i=-G_i\dd W G_i+\bigl(\tfrac12G_i+\epsilon_iG_i^2\bigr)\dd r,
 \qquad
 \dd W=N^{-1/2}\begin{pmatrix}0&\dd\mathcal B\\
                      \dd\mathcal B^*&0\end{pmatrix}.
\]
The covariance map is
$\mathcal S[R]=2\br{RE_{22}}E_{11}+2\br{RE_{11}}E_{22}$.
Thus product It\^o gives the exact forced equation
\begin{equation}
 \dd F=(F+FJF+\mathcal H)\dd r+\dd\mathcal M,
 \quad
 \mathcal H_{ab}=\epsilon_1\br{G_1^2E_{aa}G_2E_{bb}}
             +\epsilon_2\br{G_1E_{aa}G_2^2E_{bb}}.
 \label{noise:eq:lfn-four-coordinate}
\end{equation}
Its quadratic feedback closes on these four coordinates; the
displayed forcing and martingale coefficients still depend on the
underlying resolvents.
For example the cross variation is
$2F_{a2}F_{1b}+2F_{a1}F_{2b}$, which fixes the orientation of $J$.
Since
$\Tr(FJF)=4\operatorname{Re}(F_{21})\Tr F$, taking minus the trace
in \eqref{noise:eq:lfn-four-coordinate} yields
\begin{equation}
 \dd Q=(1+2\operatorname{Re}o)Q\dd r+f_Q\dd r+\dd M_Q,
 \qquad f_Q=-\Tr\mathcal H.
 \label{noise:eq:lfn-Q-SDE}
\end{equation}

For completeness, write the deterministic one-point solution as
\[
 M_i=\begin{pmatrix}iv_i&-z_i u_i\\-\overline z_i u_i&iv_i\end{pmatrix},
 \quad u_i=\frac{v_i}{\eta_i+v_i},\quad
 v_i^2=u_i-|z_i|^2u_i^2.
\]
Let $K_{ab}=\br{M_1E_{aa}M_2E_{bb}}$ and
$\overline F=K(I-JK)^{-1}$. With $Y=v_1v_2$, the active denominator is
\[
 I-JK=\begin{pmatrix}
 1-z_1\overline z_2u_1u_2&Y\\
 Y&1-\overline z_1z_2u_1u_2
 \end{pmatrix}.
\]
Its determinant has the positive decomposition
\begin{align}
 \Delta={}&(1-u_1)(1-|z_1|^2u_1u_2)
 +(1-u_2)u_1(1-|z_2|^2u_2)
 +u_1u_2|z_1-z_2|^2.
 \label{noise:eq:lfn-Delta}
\end{align}
Therefore, uniformly in the stated bulk and fine-height range,
\begin{equation}
 \overline Q:=-\Tr\overline F=\frac{Y}{\Delta}>0,
 \qquad \overline Q\asymp\frac1{|z_1-z_2|^2+a_r}.
 \label{noise:eq:lfn-Qbar}
\end{equation}
Along \eqref{noise:eq:lfn-characteristic}, $M_i'=M_i/2$, hence $K'=K$ and
$\overline F'=\overline F+\overline FJ\overline F$. In particular,
\begin{equation}
 \overline Q'=(1+2\operatorname{Re}\overline o)\overline Q,
 \qquad \overline o=2\overline F_{21}.
 \label{noise:eq:lfn-Qbar-ODE}
\end{equation}

\begin{lemma}[Relative forcing and quadratic-variation bounds]
\label{noise:lem:lfn-relative}
The forcing and real martingale in \eqref{noise:eq:lfn-Q-SDE} obey
\begin{align}
 |f_Q|&\le
 \left(\frac{|\epsilon_1|}{\eta_1}
              +\frac{|\epsilon_2|}{\eta_2}\right)Q,
 \label{noise:eq:lfn-relative-forcing}\\
 \dd[M_Q]&\le\frac{CQ^2}{N}
                 (\eta_1^{-2}+\eta_2^{-2})\dd r.
 \label{noise:eq:lfn-relative-bracket}
\end{align}
\end{lemma}
\begin{proof}
The diagonal-block part of $G_i^2$ is
\[
 D_i=(H_i^2-\eta_i^2)(H_i^2+\eta_i^2)^{-2},\qquad
 -A_i/\eta_i\le D_i\le A_i/\eta_i.
\]
The scalar $\epsilon_i$ is purely imaginary; put
$s_i=\operatorname{Im}\epsilon_i$. The diagonal-block trace in
$f_Q$ gives exactly
$f_Q=s_1\br{D_1A_2}+s_2\br{A_1D_2}$.
The Loewner inequalities, tested against the positive matrix $A_2$
or $A_1$, prove \eqref{noise:eq:lfn-relative-forcing}. This trace comparison
does not require commutation: for example, sandwich
$D_1+A_1/\eta_1\ge0$ by $A_2^{1/2}$ and take its trace.

The martingale part of $\dd A_i$ is
$-A_i\dd W G_i-G_i^*\dd W A_i$. Hence
$\dd M_Q=-\br{\dd W T}$ with the Hermitian matrix
\[
 T=G_1A_2A_1+A_1A_2G_1^*
             +G_2A_1A_2+A_2A_1G_2^*.
\]
The entry-noise normalization and normalized trace give
\[
 \dd[M_Q]=\frac{\norm{T_{12}}_{\HS}^2}{2N^3}\dd r
          \le\frac{\norm T_{\HS}^2}{4N^3}\dd r.
\]
Here $T_{12}$ is the upper-right $N$ by $N$ block; Hermitian
symmetry gives $T_{21}=T_{12}^*$.
Set $P=A_1^{1/2}A_2A_1^{1/2}\ge0$. Then
$\Tr P=2NQ$ and
\[
 G_1A_2A_1=(G_1A_1^{-1/2})P A_1^{1/2},\qquad
 \norm P_{\HS}\le\Tr P.
\]
Normality gives $G_1^*G_1=A_1/\eta_1$, so
$\norm{G_1A_1^{-1/2}}=\eta_1^{-1/2}$ and
$\norm{A_1^{1/2}}\le\eta_1^{-1/2}$. It follows that
$\norm{G_1A_2A_1}_{\HS}\le2NQ/\eta_1$.
Apply the same bound to its adjoint and interchange the two indices
for the remaining terms. Substitution proves
\eqref{noise:eq:lfn-relative-bracket}. Both estimates retain a factor of
$Q$ at the necessary multiplicative order.
\end{proof}

\begin{proposition}[Positive overlap on a common time event]
\label{noise:prop:lfn-overlap}
For every $D>0$, with probability at least $1-C_DN^{-D}$,
\begin{equation}
 \br{\operatorname{Im}G_z(X_s,ia)
               \operatorname{Im}G_w(X_s,ia)}
 \le \frac{C}{|z-w|^2+a}
 \label{noise:eq:lfn-overlap}
\end{equation}
simultaneously for $0\le s\le\tau$, $t^{b_0}/N\le a\le\sigma$,
and $z,w\in\mathcal D$ with $|z-w|\le C_0N^{-1/4}$.
The exponent $b_0=24$ is fixed before $b$ and before $D$.
\end{proposition}
\begin{proof}
First fix a deterministic observation time and terminal parameters,
and use the OU characteristic. At its initial point,
\eqref{noise:eq:lfn-characteristic-scales} and Ward's identity give
\[
 Q_0\le\br{A_{1,0}^2}^{1/2}\br{A_{2,0}^2}^{1/2}
 \le C/\sigma.
\]
Here $\br{A_i^2}\le\br{A_i}/\eta_i$.
The squared initial separation is at most $CN^{-1/2}=o(\sigma)$,
so \eqref{noise:eq:lfn-Qbar} gives $\overline Q_0\asymp\sigma^{-1}$.
Thus $Q_0/\overline Q_0\le C$ on the initial local-law event.
Only the upper separation bound is used in this initialization.
The positive-overlap conclusion therefore includes coincident
points. A fixed-separated extension instead requires its own
initial estimate, since its deterministic $\overline Q_0$ is of
order one rather than $\sigma^{-1}$.

Stop at the first adapted failure of the bounds in
\eqref{noise:eq:lfn-ordinary-input}, with fixed slack. Up to this stop,
logarithmic It\^o applied to $\mathcal R=Q/\overline Q$ gives
\begin{equation}
 \dd\log\mathcal R=\alpha_r\dd r+\dd\mathcal N_r
                         -\tfrac12\dd[\mathcal N]_r,
 \quad
 \alpha_r=2\operatorname{Re}(o-\overline o)+f_Q/Q,
 \quad \dd\mathcal N=\dd M_Q/Q.
 \label{noise:eq:lfn-log-ratio}
\end{equation}
One may initially localize away from $Q=0$. At positive heights
$Q$ is continuous and strictly positive on the compact time interval,
so this localization can then be removed. The preceding lemma and
\eqref{noise:eq:lfn-ordinary-input} imply the deterministic stopped bounds
\begin{equation}
 \int_0^{r_*}|\alpha_r|\dd r\le\frac{C_Dt^p}{Na}=o(1),
 \qquad [\mathcal N]_{r_*}\le\frac C{Na},
 \label{noise:eq:lfn-integrated-bounds}
\end{equation}
since $b_0>p+2$ and
$\int_0^{r_*}a_r^{-2}\dd r\le C/a$.
The exponential inequality for a continuous martingale of bracket
at most $C/(Na)$ yields
\[
 \PP\left(\sup_{r\le r_*}\mathcal N_r\ge1\right)
       \le e^{-cNa}\le e^{-ct^{b_0}}.
\]
We have used the stopped martingale, not conditioned on a future
local-law event. On that event the stop is absent; we add its
failure probability afterwards. The negative bracket term in
\eqref{noise:eq:lfn-log-ratio} may be discarded for this upper bound.
Equations \eqref{noise:eq:lfn-Qbar} and \eqref{noise:eq:lfn-integrated-bounds}
therefore give $Q\le C\overline Q$. Since $b_0>1$,
$e^{-ct^{b_0}}$ is smaller than every fixed power $N^{-D}$.
The overlap constant can be fixed independently of $D$: its initial
Ward bound is fixed, and for each $D$ the drift is eventually less
than one. The sufficiently large $N$ threshold may depend on $D$.

Here are explicit details for the common time event. Take spatial
and height mesh spacing $N^{-12}$ and observation-time mesh spacing
$N^{-32}$. There are at most $CN^{100}$ tests, including boundaries
and a harmless enlargement of the compact set and parameter bounds.
At heights at least $N^{-1}$, the resolvent identity gives the
deterministic modulus
\[
 |Q(X,z,w,a)-Q(X',z',w',a')|
 \le CN^3\bigl(\norm{X-X'}+|z-z'|+|w-w'|+|a-a'|\bigr).
\]
With overwhelming probability, the additive Brownian path satisfies
\[
 \sup_{|s-s'|\le N^{-32}}\norm{X_s-X_{s'}}\le C_DN^{-15}.
\]
This follows from the entrywise reflection inequality and a union
over $N^2$ entries and time cells. Thus the interpolation error is
$O_D(N^{-9})$. The right side of \eqref{noise:eq:lfn-overlap} has uniformly
polynomial first derivatives and is bounded below by a positive
constant on the enlarged parameter set. Mesh interpolation therefore
preserves the estimate with increased fixed slack. Prescribe the
failure exponent under the unweighted matrix law larger than $D+110$ before taking the
union; the exponential bound above also survives this union. The
local-law event on each entire characteristic follows from the
pathwise insertion estimate, as specified in the opening section.

Finally \eqref{noise:eq:lfn-resolvent-conversion} transfers the bound to the
additive model. Its two prefactors and the normalized separation and
height are all comparable to their additive counterparts, uniformly
for $s\le\tau=o(1)$. This proves \eqref{noise:eq:lfn-overlap}.
\end{proof}

\subsection{Low-index overlap and the exact noise covariance}

In this subsection we shall turn positive resolvent overlap bounds into control of the correlations between the noises driving the smallest singular values.

We first pass from the resolvent trace to individual singular vectors.
Positivity prevents cancellation in the spectral sum. A one-point counting
bound ensures that all of the first $K$ singular values receive sufficient
weight, so no separate rigidity estimate is needed for this passage.

Choose any fixed $k>b_0+1$ and put $K=\lceil t^k\rceil$.
The exponent can subsequently be enlarged to meet the needs of the
other local-factor arguments. All claims below hold for any such fixed choice.

\paragraph{The needed counting bound from the one-point law.}
There is a fixed $A>0$ such that, with $a=AK/N$, all times and spatial points
under consideration satisfy
\begin{equation}
 \lambda_K^v(s)\le C a,
 \label{noise:eq:lfn-low-counting}
\end{equation}
on the common local-law event. Here the positive singular values of
$X_s-v$ are ordered increasingly. To see this without an additional
rigidity hypothesis, write $h_a(x)=a/(x^2+a^2)$.
The deterministic bulk density and the one-point local law give
\[
 \frac1N\sum_i h_a(\lambda_i^v)\ge c,
 \qquad \frac1N\sum_i h_{Ca}(\lambda_i^v)\le C_1.
\]
For $x>Ca$ one has $h_a(x)\le(2/C)h_{Ca}(x)$.
Choose the fixed $C$ so that $2C_1/C\le c/2$. Since $h_a\le1/a$,
there are at least $cNa/2$ singular values below $Ca$.
Taking $A\ge2/c$ proves \eqref{noise:eq:lfn-low-counting}.
Both heights lie in the logarithmic fine range for large $N$.

Let $u_i^v,v_i^v$ denote the normalized left and right singular
vectors. Spectral decomposition gives
\begin{equation}
 2N\br{A_zA_w}
 =\sum_{i,l}h_a(\lambda_i^z)h_a(\lambda_l^w)
 \left(|\br{u_i^z,u_l^w}|^2+|\br{v_i^z,v_l^w}|^2\right),
 \label{noise:eq:lfn-positive-spectral}
\end{equation}
where in vector inner products $\br{x,y}=x^*y$.
Every term is nonnegative. By \eqref{noise:eq:lfn-low-counting}, the
two weights in a selected term with $i,l\le K$ are at least $c/a$.
Proposition~\ref{noise:prop:lfn-overlap} therefore yields
\begin{equation}
 |\br{u_i^z,u_l^w}|^2+|\br{v_i^z,v_l^w}|^2
 \le\frac{CK^2}{N|z-w|^2+K}
 \le\frac{CK^2}{1+N|z-w|^2},\qquad i,l\le K.
 \label{noise:eq:lfn-vector-overlap}
\end{equation}

The positive-index additive singular-value noises have the precise
normalization
\begin{equation}
 \dd b_i^v=\sqrt2\operatorname{Re}\bigl((u_i^v)^*\dd B\,v_i^v\bigr),
 \qquad
 \dd\lambda_i^v=\frac{\dd b_i^v}{\sqrt{2N}}+\text{drift}\,\dd s.
 \label{noise:eq:lfn-positive-noise}
\end{equation}
Adapted versions of singular vectors may be used away from collisions;
the noise is unchanged by the simultaneous phase choice. The complex
entry covariance gives exactly
\begin{equation}
 \frac{\dd[b_i^z,b_l^w]_s}{\dd s}
 =\operatorname{Re}\left(\br{u_i^z,u_l^w}\br{v_l^w,v_i^z}\right),
 \qquad \frac{\dd[b_i^v,b_l^v]_s}{\dd s}=\delta_{il}.
 \label{noise:eq:lfn-cross-covariance}
\end{equation}
The full positive vector $b=(b_1^z,\ldots,b_K^z,b_1^w,\ldots,b_K^w)$
has covariance density
\begin{equation}
 \mathsf C_s=\begin{pmatrix}I&R_s\\R_s^{\mathsf T}&I\end{pmatrix},
 \qquad
 \norm{R_s}_{\rm op}\le K\max_{i,l}|(R_s)_{il}|
 \le\frac{CK^3}{1+N|z-w|^2}\le CK^3e^{-d}.
 \label{noise:eq:lfn-cross-matrix}
\end{equation}
All bounds hold for the whole additive interval on the same event.

Thus the positive covariance has identity diagonal blocks and an
operator-norm error of order $K^3/(1+N|z-w|^2)$. These are the quantities
passed to the next section: the common adapted-time control, the identity
marginal brackets, and the small off-diagonal block $R_s$. Squaring this
last error in the whitening identity gives the later
$K^6/(1+N|z-w|^2)^2$ bound for the noise-difference bracket. The factor
$e^{-2d}$ therefore comes from quadratic covariance control, not from an
additional separation hypothesis.

\subsection{Whitening, filtrations, and independent reference processes}

Now we construct Brownian references independent of the original initial information and each other, with a controlled difference from the matrix noises.

The covariance-square-root coupling and its quadratic bracket error
appear in Cipolloni--Erd\H{o}s--Schr\"oder
\cite[equations~(7.25)--(7.30), and equations
(7.53)--(7.54)]{CES}, also used in \cite[Appendix C]{CL}.
Here the covariance estimate from the preceding section holds under
the unweighted matrix law throughout the time interval. To obtain an
exact Brownian reference on every sample path in the filtration specified
below, we whiten the positive-noise increments until an adapted stopping
time and continue with fresh noise afterwards. The vector L\'evy
characterization below proves that the Brownian vector is independent of
the enlarged initial sigma-field. The reference processes also use their
own initial spectra; their independence is from the original initial
information, as required for the local-factor argument.

Let $\mathcal F_0^{\rm mat}$ contain the entire initial matrix $X_0$
and any original initial information, and assume that future matrix
Brownian increments are independent of it. Let
$(\mathcal F_s^{\rm mat})$ be the corresponding usual filtration.
Sample independent normalized Ginibre matrices $G_z^{\rm ref},G_w^{\rm ref}$
and fresh Brownian families independently of the sigma-field generated
jointly by $\mathcal F_0^{\rm mat}$ and the \emph{entire} actual matrix
process. Set
\begin{equation}
 \mathcal G_0=\mathcal F_0^{\rm mat}\vee
               \sigma(G_z^{\rm ref},G_w^{\rm ref}),
 \label{noise:eq:lfn-initial-sigma-fields}
\end{equation}
and let $(\mathcal G_s)$ be the usual filtration additionally containing
the actual process and the fresh Brownian families up to time $s$.
The actual positive noises remain continuous martingales in this
enlargement.

\begin{theorem}[Exact independent reference noises]
\label{noise:thm:lfn-whitening}
For each prescribed deterministic pair satisfying
\eqref{noise:eq:lfn-geometry}, there is a construction of two full positive
Brownian families $(\beta_i^v)_{1\le i\le N}$, $v=z,w$, such that:
\begin{enumerate}
\item All $2N$ components are independent standard Brownian motions
relative to $(\mathcal G_s)$, and the entire vector is independent of
$\mathcal G_0$.
\item There is an adapted stopping time $\vartheta$ with
$\PP(\vartheta\le\tau)\le C_DN^{-D}$ for every fixed $D>0$, such that
before $\vartheta$ the low positive-index vectors obey
\begin{equation}
 \left\|\frac{\dd[b-\beta]_s}{\dd s}\right\|_{\rm op}
 \le\delta_N,
 \qquad
 \delta_N=\frac{CK^6}{(1+N|z-w|^2)^2}
 \le Ct^{6k}e^{-2d}.
 \label{noise:eq:lfn-delta}
\end{equation}
Here brackets of vectors are their matrix-valued quadratic variations.
\item Extend only afterwards by $\beta_{-i}^v=-\beta_i^v$.
Reference signed singular-value processes initialized from
$G_v^{\rm ref}/\sqrt{1-|v|^2}$ can be chosen as measurable strong
solutions of their own initial spectrum and their own Brownian family.
The two reference processes are mutually independent and their joint
family is independent of $\mathcal F_0^{\rm mat}$.
\end{enumerate}
\end{theorem}
\begin{proof}
\textit{Stopped construction on the positive indices.}
Use the first failure of the current additive-time one-point and
overlap bounds with fixed slack, together with
$\inf\{s\ge0:\norm{R_s}_{\rm op}\ge1/2\}$, as $\vartheta$.
Initial failure gives $\vartheta=0$, and the infimum of the empty
set is infinity. On the common event proved above this stopping
time exceeds $\tau$, since $K^3e^{-d}\to0$.
Choose a predictable covariance-density version. Before the stop,
the spectrum of $\mathsf C_s$ lies in $[1/2,3/2]$.
To define the integrand on every sample point, including after an
initial failure, set
\[
 f(x)=\bigl(\max\{x,1/2\}\bigr)^{-1/2},\qquad
 f(\mathsf C_s)\ \text{by symmetric functional calculus}.
\]
This is a bounded predictable matrix function and agrees with
$\mathsf C_s^{-1/2}$ whenever its spectrum lies in $[1/2,3/2]$.
For the low vector define
\begin{equation}
 \dd\beta_s=\one_{\{s\le\vartheta\}}f(\mathsf C_s)\dd b_s
                +\one_{\{s>\vartheta\}}\dd W_s^{\rm fresh}.
 \label{noise:eq:lfn-whitening-definition}
\end{equation}
The stopped indicators are predictable, and their values at the
single time $\vartheta$ do not affect these stochastic integrals.
Thus no inverse of a possibly singular covariance after the stop
is evaluated. On the active interval this is exactly positive
covariance whitening.
\textit{Exact Brownian law and the squared coupling error.}
Complete all indices above $K$ by distinct fresh Brownian motions,
independent also of $W^{\rm fresh}$. Every cross bracket with these
new components vanishes. The full vector starts from zero, is a
continuous $(\mathcal G_s)$-martingale, and has bracket exactly
$I_{2N}\dd s$. Vector L\'evy characterization makes it a standard
$2N$-dimensional Brownian motion relative to this filtration.
In particular its components are mutually independent and its entire
path is independent of the initial sigma-field $\mathcal G_0$.

Before $\vartheta$ the low-vector cross bracket equals
$\dd[b,\beta]_s=\mathsf C_s^{1/2}\dd s$. Hence, exactly,
\[
 \frac{\dd[b-\beta]_s}{\dd s}
  =\mathsf C_s+I-2\mathsf C_s^{1/2}
  =(\mathsf C_s^{1/2}-I)^2.
\]
Functional calculus on $[1/2,3/2]$ bounds its norm by
$C\norm{\mathsf C_s-I}_{\rm op}^2=C\norm{R_s}_{\rm op}^2$.
Equation \eqref{noise:eq:lfn-cross-matrix} proves \eqref{noise:eq:lfn-delta}.
The high-probability statement concerns closeness to the actual
noises; Brownian independence itself is exact, also after the stop.

\textit{Reference initial spectra and measurable solutions.}
Now impose the signed extension. The reference stochastic differential equation (SDE) in additive time
is, for $\mathcal I_N=\{-N,\ldots,-1,1,\ldots,N\}$,
\begin{equation}
 \dd\mu_i^v=\frac{\dd\beta_i^v}{\sqrt{2N}}
 +\frac1{2N}\sum_{\substack{j\in\mathcal I_N\\j\ne i}}
                  \frac{\dd s}{\mu_i^v-\mu_j^v},
 \qquad \mu_{-i}^v=-\mu_i^v.
 \label{noise:eq:lfn-reference-SDE}
\end{equation}
The sum includes $j=-i$; for positive $i$ this supplies the
self-mirror term $1/(4N\mu_i^v)$. There is no OU drift in this
additive-time equation. The Ginibre initial spectra are almost surely
strictly positive and simple. Existence and strong uniqueness for this
signed SDE, more generally for continuous martingale drivers of
diagonal bracket at most one, are given by
Bourgade--Cipolloni--Huang, Appendix B, Proposition B.1
\cite{BCH}. The localization proof away from collisions and explosion
gives the usual measurable strong solution; hence each $\mu^v$ is a
measurable function solely of its own initial spectrum and its own
positive Brownian family.

The inputs $(G_z^{\rm ref},\beta^z)$ and
$(G_w^{\rm ref},\beta^w)$ are mutually independent, and together
independent of $\mathcal F_0^{\rm mat}$. Their measurable solutions
have the asserted independence. Only the positive covariance was
inverted: the full signed covariance is singular and must never be
used in \eqref{noise:eq:lfn-whitening-definition}.
\end{proof}

\begin{remark}[Independence needed for local factors]
The Brownian vector is independent of $\mathcal G_0$, but the reference
processes are not: their initial spectra belong to $\mathcal G_0$.
Their relevant independence is from $\mathcal F_0^{\rm mat}$.
The references are coupled to the future matrix noise; the independence
used here concerns the original initial information.
Thus any test $F$ in \eqref{lf:eq:lfl-lf3} must be
$\mathcal F_0^{\rm mat}$-measurable. Local caps and local factors may
depend on the corresponding reference inputs only. A joint coupling
good event is used to compare actual and reference objects and is
never inserted into those local factors.

The construction is for prescribed deterministic roots, with constants
uniform in the stated geometry. Finite or polynomially many prescribed
pair constructions may be placed on a common event, with their own
enlargements if necessary. No simultaneous family of mutually
independent references indexed by every point of a continuum is
asserted or needed.
\end{remark}

\subsection{Fixed positive separation}

Next we extend the construction of reference noises to spatial points that remain a fixed positive distance apart.
For a pair with $|z-w|\ge r_0>0$, the same conclusion holds with
$\delta_N\le C_{D,r_0}K^6/N^2$, provided the following two additional
estimates under the unweighted matrix law are used on its whole characteristic:
\begin{equation}
 Q_0=\overline Q_0+o(1),\qquad
 |o_r-\overline o_r|\le\frac{C_{D,r_0}t^p}{Na_r^2}.
 \label{noise:eq:lfn-fixed-separated-premises}
\end{equation}
These are the fixed-separated logarithmic two-resolvent law, applied
at the polynomial initial height and along the entire fine-height
characteristic, respectively. The fixed-separated case of
Cipolloni--Landon, Proposition B.2(ii), gives the required estimates
\cite{CL}. For any fixed $0<\varepsilon<1/200$, the characteristic heights are
between a fixed positive multiple of $t^{24}/N$ and $C\sigma\le
N^{-\varepsilon}$ for all sufficiently large $N$, and the spatial
separation stays bounded below by a fixed multiple of $r_0$.
At each deterministic OU time the matrix is again normalized i.i.d.
with the same uniform fixed-moment bounds. The cited estimate is
applied under this unweighted matrix law before the mesh union is taken;
there is no conditioning on the future common event.
The one-point local-law bounds and the construction of a common time event are as above. The initial error tends to zero because
the initial height is comparable to $\sigma$ and
$t^p/(N\sigma^2)=o(1)$.

In this regime $\overline Q_0\asymp_{r_0}1$, so the first assertion
of \eqref{noise:eq:lfn-fixed-separated-premises} gives
$Q_0/\overline Q_0\le C_{r_0}$; the Ward estimate $Q_0\le C/\sigma$
alone would not do so. The proof of
Proposition~\ref{noise:prop:lfn-overlap} then gives $Q\le C_{r_0}$.
The positive spectral calculation yields
$\norm R\le C_{r_0}K^3/N$, and exactly the same whitening proof gives
the stated $K^6/N^2$ bound. This clause does not interpolate between
fixed separation and the near-pair upper bound $C_0N^{-1/4}$.

\endgroup

\begingroup

\section{Bounded-bracket weak laws and same-noise relaxation}\label{chap:relax}

The purpose of this section is  first to control flows with bounded noise quadratic variation, and then to show how common Brownian drivers reduce the effect of different initial data.

The independent reference noises of the preceding section will be used in
flows whose driving martingales may be correlated. This section first
proves a weak law that uses only a bound on their diagonal brackets.
The later sections verify its initial hypothesis, obtain particle counts,
and establish same-noise relaxation. The weak law is conditional on the
initial information, so it can be applied to the initial configurations
used in those comparisons.

\subsection{A logarithmic weak law for bounded-bracket singular-value flows}
\label{relax:sec:lfwl}

Next we prove a local spectral estimate for correlated driving noises with bounded individual quadratic variations. In particular, this section separates the external well-posedness results from the
logarithmic estimate proved here. Throughout,
\[
 \mathcal I_N=\{-N,\ldots,-1,1,\ldots,N\},\qquad
 \phi=(\log N)^M,\qquad h_N=\frac{\phi^6}{N},
 \qquad M>4.
\]
Constants may depend on fixed bulk bounds and on $M$, but not on the
size of the index set on which the noises are compared. We use the additive, rather than
Ornstein--Uhlenbeck, convention
\begin{equation}
 dx_i(s)=\frac{db_i(s)}{\sqrt{2N}}
   +\frac1{2N}\sum_{j\in\mathcal I_N\setminus\{i\}}
                   \frac{ds}{x_i(s)-x_j(s)}.
 \label{relax:eq:lfwl-sde}
\end{equation}
In particular the sum includes $j=-i$. The signed conventions are
$x_{-i}=-x_i$ and $b_{-i}=-b_i$.

Let $(\mathcal G_s)$ be a filtration whose initial sigma-field contains
the initial configuration. The drivers are continuous local martingales
relative to this filtration, starting at zero, with absolutely
continuous brackets and
\begin{equation}
 q_i(s):=\frac{d\langle b_i\rangle_s}{ds}\le1,
 \qquad i>0.
 \label{relax:eq:lfwl-bracket-hyp}
\end{equation}
For $i,j\in\mathcal I_N$, write
$q_{ij}(s)=d\langle b_i,b_j\rangle_s/ds$ for the signed bracket densities.
No independence assumption is imposed. Positivity of the bracket
matrix implies $|q_{ij}(s)|\le1$ for its signed entries as well.

For strictly ordered positive initial particles, strong existence,
uniqueness, noncollision and nonexplosion follow from
Bourgade--Cipolloni--Huang, \emph{Fluctuations for non-Hermitian
dynamics}, arXiv:2409.02902v2, Appendix B, Proposition B.1.
For the weak law we use the characteristic continuation argument of
\cite[Proposition~4.15]{BCH}. Both arguments allow correlated drivers
with diagonal brackets at most one. Here the regularity scale is a fixed
power of $\log N$, and conditioning on the initial information fixes the
possibly random reference measure. The stopping argument gives the
conditional probability and the uniform time--height range in
\eqref{relax:eq:lfwl-prob}--\eqref{relax:eq:lfwl-conclusion}.

\subsubsection{The conditional deterministic reference}

Next we shall condition on the initial data to fix a deterministic reference and control the flow uniformly over the allowed driving noises.
The reference may depend on the initial configuration, but not on its
future noise. Conditioning on the initial sigma-field fixes the reference
transform and every characteristic used below. It does not require the
martingale drivers to become independent Brownian motions. The resulting
estimate will be uniform over their allowed bracket processes.

Let $\mu_0$ be a symmetric probability measure, measurable with
respect to $\mathcal G_0$, and let $m_0$ be its Stieltjes transform.
Assume that fixed $r,c,C>0$ satisfy
\begin{equation}
 c\le\Im m_0(E+i\eta)\le C,\qquad
 |m_0(E+i\eta)|+|m_0'(E+i\eta)|\le C
 \label{relax:eq:lfwl-reference-initial}
\end{equation}
for $|E|\le8r$ and $0<\eta\le8r$. Conditional on $\mathcal G_0$,
this is a deterministic reference. Its free convolution with the
semicircle law of variance $s$ has transform $m_s$, characterized by
\begin{equation}
 m_s(w)=m_0(w+s m_s(w)),\qquad
 \partial_s m_s=m_s\partial_wm_s.
 \label{relax:eq:lfwl-burgers}
\end{equation}
These are the deterministic characteristic identities used in BCH
(4.39)--(4.41). By reducing a fixed $T_*>0$, the bounds
\begin{equation}
 c_1\le\Im m_s\le C_1,
 \qquad |m_s|+|\partial_wm_s|\le C_1
 \label{relax:eq:lfwl-reference-evolved}
\end{equation}
hold for $s\le T_*$ on $|E|\le4r$, $0<\eta\le4r$.
Indeed the subordination point stays in the larger initial rectangle,
and differentiation of \eqref{relax:eq:lfwl-burgers} gives the denominator
$1-s m_0'$, bounded away from zero for $T_*C$ sufficiently small.
The corresponding backward characteristics obey
\begin{equation}
 \dot w_u=-m_u(w_u),\qquad
 m_u(w_u)=m_s(w_s),\qquad
 \eta_u=\eta_s+(s-u)\Im m_s(w_s),
 \label{relax:eq:lfwl-characteristic}
\end{equation}
where $0\le u\le s$ and $\eta_u=\Im w_u$. Their spatial derivatives
are bounded by fixed constants on the smaller rectangles used below.
All their coefficients are initial-data measurable, so their use in
stochastic integrals is nonanticipative.

Define
\[
 m_{N,s}(w)=\frac1{2N}\sum_{i\in\mathcal I_N}
                       \frac1{x_i(s)-w}.
\]
The initial regularity hypothesis for the next proposition is
\begin{equation}
 |m_{N,0}(E+i\eta)-m_0(E+i\eta)|
       \le A\frac{\phi^{1/4}}{\sqrt{N\eta}},
 \qquad |E|\le4r,\quad h_N\le\eta\le4r,
 \label{relax:eq:lfwl-initial}
\end{equation}
where $A$ is fixed. It is a hypothesis about the actual initial
configuration, not its conditional expectation.

\begin{proposition}[Conditional bounded-bracket continuation]
\label{relax:prop:lfwl-continuation}
Assume \eqref{relax:eq:lfwl-reference-initial} and
\eqref{relax:eq:lfwl-initial}, with strictly ordered positive initial
particles, and let \eqref{relax:eq:lfwl-sde} satisfy
\eqref{relax:eq:lfwl-bracket-hyp}. Put
$T_N=\min\{T_*,\phi^{-2}\}$. There are fixed constants
$C_2,C_3,c_2>0$ such that, except on an event of conditional
probability at most
\begin{equation}
 C_2N^{C_3}\exp(-c_2\phi^{1/2}),
 \label{relax:eq:lfwl-prob}
\end{equation}
given $\mathcal G_0$, one has simultaneously
\begin{equation}
 |m_{N,s}(E+i\eta)-m_s(E+i\eta)|
       \le C_2(A+1)\frac{\phi^{1/4}}{\sqrt{N\eta}},
 \qquad
 0\le s\le T_N,\quad |E|\le r,\quad h_N\le\eta\le r.
 \label{relax:eq:lfwl-conclusion}
\end{equation}
The constants are uniform over the allowed bracket processes. If the
initial hypotheses hold only on a $\mathcal G_0$-measurable event,
their failure probability is added to \eqref{relax:eq:lfwl-prob}.
\end{proposition}

\begin{proof}
\textit{The empirical equation and its two error terms.}
We first derive the exact empirical equation. Write
$g_i=(x_i-w)^{-1}$. It\^o's formula and symmetrization give
\[
 -\frac1{4N^2}\sum_{i\ne j}\frac{g_i^2}{x_i-x_j}
       =m_N\partial_wm_N-\frac1{4N^2}\sum_i g_i^3.
\]
The diagonal quadratic variation in It\^o's formula adds
$(4N^2)^{-1}\sum_i q_i g_i^3ds$. Consequently
\begin{align}
 dm_{N,s}(w)
   &=m_{N,s}(w)\partial_wm_{N,s}(w)\,ds+dM_s(w)+dR_s(w),
                                                        \label{relax:eq:lfwl-ito}\\
 dM_s(w)&=-\frac1{2N\sqrt{2N}}
                       \sum_i\frac{db_i(s)}{(x_i(s)-w)^2},
                                                        \label{relax:eq:lfwl-mart}\\
 dR_s(w)&=\frac1{4N^2}\sum_i
                       \frac{q_i(s)-1}{(x_i(s)-w)^3}\,ds.
                                                        \label{relax:eq:lfwl-cubic}
\end{align}
This is also BCH (4.48). In particular the cubic term is retained
when the diagonal brackets are not one. Off-diagonal correlations enter
the martingale bracket, while the cubic correction uses the individual
$q_i$. These are distinct error terms and will be bounded separately.

For $w=E+i\eta$, the identity
\[
 \sum_i|x_i-w|^{-2}=2N\frac{\Im m_N(w)}{\eta}
\]
and $|q_{ij}|\le1$ yield
\begin{equation}
 d\langle M,\overline M\rangle_s
     \le\frac C N\left(\frac{\Im m_{N,s}(w)}\eta\right)^2ds,
 \qquad
 |dR_s(w)|\le\frac{C\Im m_{N,s}(w)}{N\eta^2}\,ds.
 \label{relax:eq:lfwl-noise-bounds}
\end{equation}
The first bracket is the positive total quadratic variation; separate
real and imaginary martingales obey the same bound. All the constants
are independent of the covariance rank and of the number of drivers
on which the covariance construction may depend.

Set $\Delta_s=m_{N,s}(w_s)-m_s(w_s)$ along
\eqref{relax:eq:lfwl-characteristic}. Equations
\eqref{relax:eq:lfwl-burgers} and \eqref{relax:eq:lfwl-ito} imply
\begin{equation}
 d\Delta_s=\Delta_s\partial_wm_{N,s}(w_s)\,ds
                                     +dM_s(w_s)+dR_s(w_s).
 \label{relax:eq:lfwl-delta}
\end{equation}
\textit{The coarse estimate and the common stopping domain.}
First, on a coarse region $\eta\ge r/2$, the deterministic bounds
$|\partial_wm_N|\le\eta^{-2}$ and $\Im m_N\le\eta^{-1}$
bound the coefficient in \eqref{relax:eq:lfwl-delta} by a constant, the
accumulated martingale bracket by $CT_N/N$, and the remainder by
$CT_N/N$. The initial hypothesis and the exponential martingale
inequality therefore give
\begin{equation}
 |m_{N,s}-m_s|\le C(A+1)\phi^{1/4}/\sqrt N
 \label{relax:eq:lfwl-coarse}
\end{equation}
uniformly in the required coarse region and times, with an exception
of the form \eqref{relax:eq:lfwl-prob}. Parameter uniformity follows from
the mesh argument below and does not use a matrix local law at
positive times. The larger initial and reference rectangles provide
fixed lateral and upper buffers for these characteristics.

For the finer region choose $C_0>2C_1$ and consider
\[
 \mathcal D_s=\{E+i\eta:|E|\le2r-C_0s,\quad
                           h_N\le\eta\le2r\}.
\]
For all sufficiently large $N$ this contains the target rectangle.
Every backward characteristic from its fine part stays inside its
lower and lateral boundaries; if it enters $\eta\ge r/2$ it is
already covered by the coarse estimate. Stop at the first failure
of \eqref{relax:eq:lfwl-coarse} or of
\begin{equation}
 |m_{N,s}(w)-m_s(w)|<\frac{\phi}{\sqrt{N\Im w}}
 \label{relax:eq:lfwl-stop}
\end{equation}
on this region. A countable dense parameter set and continuity define
an adapted stopping time. Initial failure is excluded for large $N$
by \eqref{relax:eq:lfwl-initial}. Before the stop,
$c'\le\Im m_{N,s}\le C'$ because $N\eta\ge\phi^6$.

\textit{Maximal control along a conditional characteristic.}
For a fixed deterministic terminal pair $(s,w_s)$, conditional on
$\mathcal G_0$, the characteristic is deterministic. By
\eqref{relax:eq:lfwl-noise-bounds} and
$-\dot\eta_u\asymp1$, its stopped martingale and remainder satisfy,
for each physical time $v\le s$,
\begin{equation}
 \langle M,\overline M\rangle_v\le\frac C{N\eta_v},
 \qquad \operatorname{Var}_{[0,v]}R\le\frac C{N\eta_v}.
 \label{relax:eq:lfwl-integrated-noise}
\end{equation}
Here $\operatorname{Var}_{[0,v]}R$ denotes the total variation of the
path $R$ on $[0,v]$.
On dyadic ranges of the deterministic decreasing height $\eta_v$,
the exponential inequality, separately for real and imaginary parts,
gives
\begin{equation}
 \sup_{u\le v}|M_u|+\operatorname{Var}_{[0,v]}R
       \le C\frac{\phi^{1/4}}{\sqrt{N\eta_v}}
 \label{relax:eq:lfwl-mart-envelope}
\end{equation}
simultaneously in physical time, except with probability
$C(1+\log N)e^{-c\phi^{1/2}}$. The deterministic remainder is
absorbed because $N\eta_v\ge\phi^6$.

\textit{Derivative control above the cutoff and strict improvement.}
We next control the drift coefficient without using values below
$h_N$. At a fine point $w=E+i\eta$ put
\[
 e_\eta=\frac{\phi}{\sqrt{N\eta}},\qquad
 w^+=w+i\eta\sqrt{e_\eta}.
\]
Here $e_\eta\le\phi^{-2}$, so $w^+$ stays above $w$ and is at a
distance $o(\eta)$ from it. Both points are in the stopped region
or its controlled coarse part. Taylor's formula along the upward
segment gives
\[
 |\partial_wm_N(w)|
 \le C+\frac{C e_\eta}{|w^+-w|}
           +C\frac{|w^+-w|}{\eta^2}
 \le C+\frac{C\phi^{1/2}}{\eta(N\eta)^{1/4}}.
\]
Indeed $|m_N''|\le C\Im m_N/\eta^2$ on the segment, and the
reference derivative is bounded. This is the upward version of BCH
Lemma A.2. If a derivative is subsequently estimated by Cauchy's
formula, disks centered on the upward segment may be chosen with
their lower edge at or above the original height; no control below
the stopped lower boundary is used. Integrating along a
characteristic yields
\begin{equation}
 \int_0^v |\partial_wm_{N,u}(w_u)|\,du
       \le CT_N+\frac{C\phi^{1/2}}{(N\eta_v)^{1/4}}
       \le C.
 \label{relax:eq:lfwl-derivative-integral}
\end{equation}

Integral Gronwall applied to \eqref{relax:eq:lfwl-delta}, retaining its
initial term, now gives
\begin{equation}
 |\Delta_v|\le C\left(|\Delta_0|
                  +\frac{\phi^{1/4}}{\sqrt{N\eta_v}}\right)
       \le C(A+1)\frac{\phi^{1/4}}{\sqrt{N\eta_v}}.
 \label{relax:eq:lfwl-improvement}
\end{equation}
We used $\eta_0\ge\eta_v$. This strictly improves
\eqref{relax:eq:lfwl-stop} for large $N$.

\textit{Simultaneous parameters and the first-exit conclusion.}
The preceding bound is first proved on fixed terminal characteristics.
To use it at a first-exit point, we now construct one event on which it
holds simultaneously. The maximal estimate already covers every physical
time on each grid characteristic; the grid below discretizes terminal
parameters, not the martingale's physical time.
Use grids of mesh
$\delta=N^{-C_4}$ in terminal time, real coordinate and imaginary
coordinate, including the final time. Fix an arbitrary target
physical time $v$ and target point $w\in\mathcal D_v$ in the fine
region. Choose a grid terminal time $S\ge v$ with $S-v\le\delta$.
Move the terminal real coordinate inward by $O(\delta)$ and the
terminal imaginary coordinate upward by $O(\delta)$ before
rounding to the grid. The constants in these moves exceed
$C_0+2C_1+2$. Thus the grid terminal point lies inside
$\mathcal D_S$; its backward characteristic stays above $h_N$
and within the shrinking lateral boundary. At the physical time
$v$, that characteristic differs from $w$ by $O(\delta)$.
The upper boundary is handled by the already controlled coarse
region and its fixed buffer. For each grid characteristic the
maximal martingale estimate holds at every physical time up to
$S$, so it applies at this same $v$ without a union over physical
time.

For interpolation, the physical integration endpoint is always
held fixed. More explicitly, compare the target and grid
characteristics on $[0,v\wedge\sigma]$, where $\sigma$ is the
common stopping time for the whole parameter domain. Their
distance there is $O(\delta)$ by the deterministic characteristic
bounds. On $\eta\ge h_N/2$, the estimates
$|\partial_w^j m_N|\le j!\eta^{-j-1}$ are polynomial in $N$.
The integrated identity \eqref{relax:eq:lfwl-delta}, evaluated at this
same $v\wedge\sigma$ for both paths, expresses their stopped
martingale difference as endpoint differences minus the drift and
cubic-remainder differences. These are bounded by $CN^4\delta$
using the resolvent derivative bounds and $|q_i-1|\le1$.
This is Lipschitz dependence on terminal parameters at a fixed
physical endpoint; no Lipschitz dependence on $v$ itself is
asserted or needed. Taking, for example, a sufficiently large
fixed $C_4>10$ makes the interpolation error negligible compared
with the smallest martingale threshold, of order at least
$\phi^{1/4}N^{-1/2}$. The polynomial grid size and dyadic-height
union produce \eqref{relax:eq:lfwl-prob}. All grid exponents are fixed
before any later noise-window exponents.

To make the simultaneous interpretation precise, first apply the
stochastic identities on a countable dense family of terminal
parameters. Define their extension by the integrated identity
\eqref{relax:eq:lfwl-delta}, with the same physical endpoint and common
stop. The empirical resolvents are continuous, and the drift and
cubic-remainder integrals are continuous in terminal parameters by
the deterministic derivative bounds just used. This defines a
jointly parameter-continuous version of the stopped martingale
family, agreeing with each original stochastic integral on the
dense family. The $CN^4\delta$ comparison and hence the improvement
then extend to every terminal pair on a single event. In particular,
evaluating that continuous family at a first-exit parameter does
not invoke a stochastic estimate with an anticipatively selected
terminal characteristic.

Thus the strict improvement holds throughout the stopped domain,
so its fine stopping boundary is not reached. Intersecting with the
coarse event proves \eqref{relax:eq:lfwl-conclusion}. Since $M>4$,
\eqref{relax:eq:lfwl-prob} is smaller than $N^{-D}$ for every fixed $D$
once $N$ is sufficiently large. All estimates remain valid
conditionally on $\mathcal G_0$, completing the proof.
\end{proof}

\subsubsection{Verification of the initial input after polynomial Gaussian age}

Next we show that a short Gaussian smoothing step supplies the initial spectral regularity needed by the conditional weak law.

The continuation proposition assumes a logarithmic initial estimate for
the actual configuration. Gaussian aging supplies that estimate from a
fixed-power matrix local law at larger imaginary height. The backward
characteristic carries the terminal observation back to height comparable
to the terminal height plus the age, where the initial error is small
enough to close the same
stopping argument. All spectral parameters here include their real part;
this off-axis control will be used for interval counts.

\begin{proposition}[Aged iid initial regularity]
\label{relax:prop:lfwl-aged}
Let $\sigma=c_NN^{-a_{\rm age}}$ with fixed $a_{\rm age}\in(0,1)$,
where $c_N$ is deterministic and lies in a fixed compact subinterval
of $(0,\infty)$, and
\[
 A_s=\sqrt{1-\sigma}\,\widehat X+N^{-1/2}B_s,
 \qquad 0\le s\le\sigma,
\]
where $\widehat X$ is a complex iid matrix with centered normalized
entries, vanishing complex second moment, and a fixed common bound
on every finite entry moment. Fix the spatial parameter $z$ in a
compact subset of the unit disk. The variance-$1-\sigma+s$
deterministic Hermitization density is denoted by $\mu_s^z$.
For a sufficiently small fixed $r>0$, with overwhelming probability,
\begin{equation}
 |m_{N,\sigma}^z(E+i\eta)-m_\sigma^z(E+i\eta)|
 \le C\frac{\phi^{1/4}}{\sqrt{N\eta}}
       +\frac{C N^\xi}{N(\eta+\sigma)},
 \quad |E|\le4r,\quad h_N\le\eta\le4r,
 \label{relax:eq:lfwl-aged-estimate}
\end{equation}
where $0<\xi<(1-a_{\rm age})/4$ is fixed. In particular the aged
configuration satisfies \eqref{relax:eq:lfwl-initial}, after increasing
its fixed constant. The statement is uniform over the allowed age
factors, the fixed bulk set and the prescribed common moment class,
and admits a union over
any fixed polynomial-size family of spatial parameters.
\end{proposition}

\begin{proof}
The variance-dependent MDE, or its scalar equation differentiated in
the variance, gives \eqref{relax:eq:lfwl-burgers}. Its branch is uniformly
stable near zero in the fixed spatial bulk. Thus it supplies
\eqref{relax:eq:lfwl-reference-evolved}, and any fixed derivatives needed
on smaller rectangles, throughout the aging interval.

The only random-matrix local law used at the beginning of this
interval is the fixed-power averaged law: CL
arXiv:2605.29962v3, Theorem 2.3, (2.8). Its spectral variable is the
general $E+i\eta$. After the harmless variance rescaling it gives
\[
 |m_{N,0}^z(w)-m_0^z(w)|\le\frac{N^\xi}{N\Im w}
 \quad\hbox{for }\Im w\ge c_0\sigma
\]
on the required fixed rectangle, with overwhelming
probability. The rescaled spatial point remains in a fixed enlarged
bulk disk and the entry moment bounds remain uniform.

Repeat the proof of Proposition~\ref{relax:prop:lfwl-continuation} on the
sloping domain
\[
 h_N+c_0(\sigma-s)\le\eta,
\]
where $0<c_0<c_1/2$. The polynomial age $\sigma$ is smaller than the
fixed logarithmic horizon $\phi^{-2}$ for all sufficiently large $N$.
Backward characteristics increase in height at
speed at least $c_1$, so they stay above this boundary and start at
height at least $c_0\sigma$. The upward derivative segment stays in
the domain as before. The coarse region follows from the same
coarse resolvent argument, with initial error $O(N^{-1+\xi})$;
there is no logarithmic local-law assumption at the unaged time.

The integrated martingale, cubic remainder and derivative estimates
are unchanged. Keeping the initial term in
\eqref{relax:eq:lfwl-improvement}, instead of imposing
\eqref{relax:eq:lfwl-initial}, gives
\[
 |\Delta_s|\le C\frac{N^\xi}{N\eta_0}
                +C\frac{\phi^{1/4}}{\sqrt{N\eta_s}}.
\]
Its initial contribution is strictly below the stopping threshold:
\[
 \frac{N^\xi/(N\eta_0)}{\phi/\sqrt{N\eta_s}}
 \le\frac{N^\xi}{\phi\sqrt{N\eta_0}}
 \le C\phi^{-1}N^{\xi-(1-a_{\rm age})/2}=o(1).
\]
At time $\sigma$, the characteristic identity gives
$\eta_0\ge\eta_\sigma+c_1\sigma$, proving
\eqref{relax:eq:lfwl-aged-estimate}. Its second term is uniformly
$o((N\eta)^{-1/2})$, since
\[
 \frac{N^\xi\sqrt{N\eta}}{N(\eta+\sigma)}
       \le C N^{\xi-(1-a_{\rm age})/2}=o(1).
\]

The actual matrix process exists even if the unaged spectrum has a
collision. One may apply It\^o directly to its resolvent, or begin
the singular-value calculation at a positive time and then let that
time decrease to zero at positive imaginary height. For each
prescribed deterministic spatial point, Gaussian smoothing makes
the terminal spectrum strictly ordered and positive almost surely;
the same holds for any prescribed finite polynomial-size list.
No simultaneous simplicity or positivity assertion over a
continuum of spatial points is used. The transform regularity
estimate itself can be made spatially uniform. A spatial mesh costs only a fixed
polynomial factor, absorbed by the initial local-law failure
choice and the stretched-logarithmic martingale bound.
\end{proof}

\subsubsection{Interval counts and the order of logarithmic parameters}

Next we turn the weak law into counts on shifted intervals and then into spacing bounds for sufficiently separated particle indices.

The full real-variable range of the weak law controls masses of shifted
intervals. We first compare their Poisson smoothings, then account for the
sharp endpoints, and only then convert counts into rank gaps. An estimate
restricted to $E=0$ would not supply this interval-by-interval argument.

\begin{corollary}[Uniform central interval counting]
\label{relax:cor:lfwl-counts}
On the event of Proposition~\ref{relax:prop:lfwl-continuation}, after
shrinking the real interval once more, the empirical measure
$\mu_{N,s}=(2N)^{-1}\sum_i\delta_{x_i(s)}$ satisfies
\begin{equation}
 |\mu_{N,s}(I)-\mu_s(I)|
   \le C\phi\sqrt{\frac{|I|}{N}}+C h_N(1+\log N)
 \label{relax:eq:lfwl-counts}
\end{equation}
for every central interval $I$, simultaneously in $0\le s\le T_N$.
The reference density is bounded above and below there.

Put $t=\tfrac12\log N$. For any fixed $a>6M+3$, set $Q_0=t^a$.
In a fixed smaller central index block, with harmless fixed changes
to its endpoints,
\begin{equation}
 c\frac{|i-j|}{N}\le|x_i(s)-x_j(s)|
       \le C\frac{|i-j|}{N},\qquad |i-j|\ge Q_0,
 \label{relax:eq:lfwl-gaps}
\end{equation}
and $|x_i(s)|\asymp|i|/N$ for $|i|\ge Q_0$ in that block.
Particles outside a fixed smaller central index block are separated
from zero by a fixed constant.
\end{corollary}

\begin{proof}
Use the weaker consequence
$|m_{N,s}-m_s|\le\phi/\sqrt{N\eta}$ and suppress $s$.
For $I=[a,b]$ of length $h\ge h_N$, analyticity gives
\[
 \frac d{d\eta}\int_a^b\Im(m_N-m)(E+i\eta)\,dE
   =\Re(m_N-m)(b+i\eta)-\Re(m_N-m)(a+i\eta).
\]
Integrate this identity from $h_N$ to $h$. The integral at height
$h$ is bounded by $C\phi\sqrt{h/N}$, and integrating the endpoint
errors gives the same bound. This proves the corresponding
Poisson-smoothed mass discrepancy at height $h_N$.

Bounded imaginary transforms imply that every central interval of
length $v\ge h_N$ has mass at most $Cv$, for either measure.
The absolute error between a sharp interval and its Poisson
smoothing is $O(h_N)$ within distance $h_N$ of each endpoint.
On a dyadic band at distance $v\ge h_N$, the difference of the two
indicator kernels is at most $C h_N/v$ and the band has mass at
most $Cv$. Each such band costs $C h_N$. There are $O(1+\log N)$
bands before reaching the fixed exterior of the bulk region.
Beyond that exterior the kernel is $O(h_N)$ and the total mass is
at most one. Applying this to both measures proves
\eqref{relax:eq:lfwl-counts}. When $h<h_N$, the same imaginary-transform
bound gives mass $O(h_N)$, already included in the additive term.

For intervals with length at least $t^{6M+2}/N$, the error in
\eqref{relax:eq:lfwl-counts}, divided by their length, tends to zero.
The fixed factor between $t^M$ and $(\log N)^M$ is harmless.
Thus such intervals have between $cN|I|$ and $CN|I|$ particles.
Since $Q_0\gg t^{6M+2}$, the standard rank-to-count inversion,
allowing two endpoint particles, gives \eqref{relax:eq:lfwl-gaps}.
Counting from zero gives the individual location statement;
counting to a fixed small real endpoint gives separation of the
remaining ranks from zero. All estimates were uniform in $s$, so
these consequences are uniform in time as well.
\end{proof}

\paragraph{Application to the noise-interpolating family.}
Let the initial spectrum be that of the aged matrix from
Proposition~\ref{relax:prop:lfwl-aged}, and let $b_i,\beta_i$ be the
actual and reference drivers in a common filtration containing this
matrix. For each fixed $\alpha\in[0,1]$, define
\[
 b_i^\alpha=\alpha b_i+(1-\alpha)\beta_i.
\]
If each endpoint has diagonal bracket at most one, then
\[
 \frac{d\langle b_i^\alpha\rangle_s}{ds}
 \le\left(\alpha\sqrt{\frac{d\langle b_i\rangle_s}{ds}}
 +(1-\alpha)\sqrt{\frac{d\langle\beta_i\rangle_s}{ds}}\right)^2
 \le1.
\]
The initial spectrum is the same for every $\alpha$, so the
deterministic density at variance $1+s$ is also a common reference.
Proposition~\ref{relax:prop:lfwl-continuation} and
Corollary~\ref{relax:cor:lfwl-counts} apply to each fixed $\alpha$ with
uniform failure constants. They make no claim of simultaneous
noncollision or a uniform tangent event over all $\alpha$.
The interpolated particle system need not be the spectrum of an iid
matrix. No iid local law is applied to its positive-time marginals;
the matrix local law is used only to initialize the earlier Gaussian
smoothing argument.
The later tangent-integration argument uses Fubini and its separate
energy estimate.

The order of parameter choices is explicit. Fix $M$, and then fix
$a>6M+3$. The theorem has already been proved up to
$T_N=\phi^{-2}$ for large $N$. For every subsequently chosen fixed
$b$, the interval $[0,\tau]$, $\tau=t^b/N$, is contained in that
horizon for sufficiently large $N$. Finally choose the fixed noise
windows $L=t^\ell$ and $K=t^k$ as required by the different-noise
energy proof. Their enlargement changes neither $M$ nor $a$,
because only the dimension-independent bounds
$q_i^\alpha\le1$ and $|q_{ij}^\alpha|\le1$ enter this section.
The eventual threshold in $N$ may depend on all ultimately fixed
exponents. The argument retains the cubic correction when the diagonal
brackets differ from one.

\subsection{Logarithmic same-noise relaxation}
\label{relax:sec:lfsn}

The purpose of this subsection is to show that common Brownian drivers reduce the effect of different initial data on the smallest singular values.
In the following, $u$ and $t$ denote particle-flow times.
We use signed indices
$\mathcal I_N=\{-N,\ldots,-1,1,\ldots,N\}$ and the additive complex
singular-value convention
\begin{equation}
 dx_i(u)=\frac{db_i(u)}{\sqrt{2N}}
   +\frac1{2N}\sum_{j\in\mathcal I_N\setminus\{i\}}
                         \frac{du}{x_i(u)-x_j(u)},
 \qquad x_{-i}=-x_i,\quad b_{-i}=-b_i.
 \label{relax:eq:lfsn-sde}
\end{equation}
For $i>0$ the $b_i$ are independent standard Brownian motions relative
to a filtration containing the full initial configuration. The sum
includes $j=-i$. Strictly ordered positive initial data give a
noncolliding strong solution by Bourgade--Cipolloni--Huang,
\emph{Fluctuations for non-Hermitian dynamics},
arXiv:2409.02902v2, Appendix B, Proposition B.1.
The logarithmic relaxation estimate is proved below using the initial
counts and the common noise. The external well-posedness theorem supplies
the solutions, while the quantitative logarithmic losses and time range
are established in this section.

\subsubsection{Initial configurations and choice of parameters}

Now we are able to connect the two initial configurations by an interpolation whose spectral regularity remains strong enough for the relaxation argument.
The two endpoint configurations $s_i,r_i$ are symmetric, strictly
ordered, and have fixed spectral norms on an event of overwhelming
probability. For $e=s,r$, write
$\mu_N^e=(2N)^{-1}\sum_{i\in\mathcal I_N}\delta_{e_i}$ and let
$\mu^e$ denote the corresponding reference measure. On a fixed
neighborhood of zero, the reference densities $\rho_s,\rho_r$ are
positive with uniformly bounded first two derivatives. Assume that
\begin{equation}
 \rho_s(0)=\rho_r(0)=A>0,\qquad
 |\mu_N^e(I)-\mu^e(I)|
 \le C\left(P\sqrt{\frac{|I|}{N}}+\frac{P^2}{N}\right),
 \quad e=s,r.
 \label{relax:eq:lfsn-endpoint-count}
\end{equation}
The constants are uniform over the endpoint class. Here $P$ is a
fixed power of $\log N$.
The endpoints are included in the initial sigma-field of the common
Brownian family. In the matrix application this is the enlarged initial
field of \eqref{noise:eq:lfn-initial-sigma-fields};
Theorem~\ref{noise:thm:lfn-whitening} supplies a positive Brownian family
independent of that entire field. This permits conditioning on both
endpoint configurations in the arguments below.

For the aged iid endpoint, these assumptions follow from the
logarithmic off-axis weak law obtained by Brownian smoothing, the
sharp-counting consequence of
Proposition~\ref{relax:prop:lfwl-continuation}, and deterministic bulk
density regularity. The reference endpoint is the appropriately scaled complex Ginibre
singular-value configuration, sampled independently of
$\mathcal F_0^{\rm mat}$ in \eqref{noise:eq:lfn-initial-sigma-fields}.
The fixed endpoint norm bound follows, for example, from
Cipolloni--Landon, arXiv:2605.29962v3, Theorem 2.5, equation (2.11),
at the last singular value, together with Lemma 2.1(ii).
The edge quantile is bounded and its
$N^\delta N^{-2/3}$ error is $o(1)$ when $\delta<2/3$.

The central quantile interpolation with empirical exterior atoms follows
the construction in \cite[Section~7.6, equations~(7.57)--(7.60)]{CES}.
The proof below verifies the buffered construction and the logarithmic
counting error needed here.

\begin{lemma}[Quantile interpolation and initial tangent domination]
\label{relax:lem:lfsn-quantile}
Under \eqref{relax:eq:lfsn-endpoint-count}, let
$x_i^\alpha=\alpha s_i+(1-\alpha)r_i$, $0\le\alpha\le1$.
Let $\mu_{N,\alpha}=(2N)^{-1}
\sum_{i\in\mathcal I_N}\delta_{x_i^\alpha}$, and write
$m_{N,\alpha}$ for its Stieltjes transform.
There are symmetric probability measures $\widetilde\mu_\alpha$,
measurable with respect to the endpoint configurations, with density
$A$ at zero. Write $\widetilde m_\alpha$ for their Stieltjes transforms.
For $w=E+i\eta$ in a fixed smaller upper-half-plane rectangle adjacent
to zero, we have
\begin{equation}
 c\le\Im\widetilde m_\alpha\le C,\qquad
 |\widetilde m_\alpha|+|\widetilde m_\alpha'|\le C,
 \qquad
 |m_{N,\alpha}-\widetilde m_\alpha|
 \le C\left(\frac P{\sqrt{N\eta}}+\frac{P^2}{N\eta}\right).
 \label{relax:eq:lfsn-quantile-law}
\end{equation}
Moreover,
\begin{equation}
 |s_i-r_i|\le CP^2
       \left(\frac{\sqrt i}{N}+\frac{i^2}{N^2}\right),
 \qquad 1\le i\le N.
 \label{relax:eq:lfsn-initial-tangent}
\end{equation}
\end{lemma}

\begin{proof}
The reference keeps the empirical exterior atoms and replaces only a
central block by a smooth quantile interpolation. We use a smaller real
interval with a fixed buffer from the artificial cutoff. This buffer is
needed because the empirical and continuous cutoff locations need not
coincide.

We may first symmetrize each endpoint reference measure. Symmetry of
the empirical measure preserves the interval bound under this operation;
the local $C^2$ bounds and the density $A$ at zero are also preserved.
We use these symmetric reference densities for the quantiles below.

Let $\Gamma_e$ be the signed local quantile determined by
$u=\int_0^{\Gamma_e(u)}\rho_e(v)\,dv$. Choose $u_*>0$ small enough
that these quantiles are defined on a fixed neighborhood of
$[-2u_*,2u_*]$. Put $K=\lfloor2Nu_*\rfloor$, $u_K=K/(2N)$, and set
\[
 \Gamma_\alpha=\alpha\Gamma_s+(1-\alpha)\Gamma_r,\qquad
 \widetilde\mu_\alpha
 =\int_{-u_K}^{u_K}\delta_{\Gamma_\alpha(u)}\,du
    +\frac1{2N}\sum_{K<|i|\le N}\delta_{x_i^\alpha}.
\]
The two masses add to one. The central quantile derivative is
bounded above and below, and its first three derivatives are bounded.
Counting from zero gives, uniformly for $|i|\le K$,
\[
 |e_i-\Gamma_e(i/(2N))|\le C D_N,
 \qquad D_N=\frac P{\sqrt N}+\frac{P^2}{N}.
\]
The same estimate holds after interpolation. In particular, the
retained exterior atoms remain a fixed distance from zero, uniformly
in $\alpha$. The continuous central density is uniformly $C^2$ and
equals $A$ at zero, since $\Gamma_\alpha'(0)=1/A$.
Work on a fixed smaller rectangle whose real interval has positive
distance from all central cutoff points
$\pm\Gamma_\alpha(u_K)$ and all retained atoms. The atoms' local
transform and derivative are bounded. For the continuous part,
integration by parts and subtraction of the value of the density
derivative at the real part of the spectral parameter give bounded
transform and derivative. Positivity of the central density gives a
uniform positive lower bound for the imaginary part. These prove
the deterministic reference bounds in
\eqref{relax:eq:lfsn-quantile-law}.

We next justify the local increment estimate. For central positive
indices $i<j$, put $h=(j-i)/(2N)$ and $L_e=e_j-e_i$.
The endpoint interval estimate, including the two endpoint atoms,
and the lower density bound imply
\[
 L_e\le C\left(h+P\sqrt{L_e/N}+P^2/N\right)
       \le C(h+P^2/N).
\]
It follows that the mass error on this particle interval is at most
\[
 \varepsilon_h=C\left(P\sqrt{h/N}+P^2/N\right).
\]
The starting quantile mass differs from $i/(2N)$ by $O(D_N)$.
The bounded second derivative of $\Gamma_e$ makes the effect of
that common starting displacement on a quantile increment of
length $h$ at most $ChD_N\le C\varepsilon_h$. Consequently
\[
 |(e_j-e_i)-\{\Gamma_e(j/(2N))-\Gamma_e(i/(2N))\}|
 \le C\varepsilon_h.
\]
Symmetry gives the corresponding estimate on the negative side
and across zero; the missing zero index changes counts by at most
one atom. This increment estimate passes linearly to the interpolated
positions. Its one-step case bounds central neighboring gaps by
$CP^2/N$.

For an interval $I$ in a fixed buffered interior real neighborhood,
use its first and last particles and the preceding increment bound.
Comparing the reference mass at their displaced starting quantile
costs only $ChD_N$ again. The two boundary gaps cost $CP^2/N$.
The same argument applied to an empty interval uses just one such
gap. Thus
\[
 |\mu_{N,\alpha}(I)-\widetilde\mu_\alpha(I)|
 \le C\left(P\sqrt{|I|/N}+P^2/N\right)
\]
for these interior intervals. We do not assert this estimate near
the artificial central cutoff: its empirical and continuous endpoints
can differ by $O(D_N)$.

Here is the required control of that exterior when passing to
transforms. The empirical exterior atoms cancel exactly in
$\nu_\alpha=\mu_{N,\alpha}-\widetilde\mu_\alpha$, and its total
mass is zero. Couple each central empirical atom with its adjacent
quantile cell of mass $1/(2N)$. The global central quantile error
above and the bounded density give
\[
 F_\alpha(x):=\nu_\alpha(( -\infty,x]),\qquad
 \sup_x|F_\alpha(x)|\le CD_N.
\]
For $w=E+i\eta$ in a still smaller fixed rectangle, integration by
parts, followed by
$\int_{\mathbb R}(x-w)^{-2}\,dx=0$, gives
\[
 m_{N,\alpha}(w)-\widetilde m_\alpha(w)
 =\int_{\mathbb R}\frac{F_\alpha(x)-F_\alpha(E)}{(x-w)^2}\,dx.
\]
Near $E$, the interior interval estimate bounds the numerator by
$C(P\sqrt{|x-E|/N}+P^2/N)$. Outside a fixed neighborhood of $E$
use the global $CD_N$ bound; its contribution is $O(D_N)$.
The elementary integrals
\[
 \int_0^c\frac{\sqrt v\,dv}{v^2+\eta^2}\le C\eta^{-1/2},
 \qquad
 \int_0^c\frac{dv}{v^2+\eta^2}\le C\eta^{-1}
\]
now prove the claimed transform error. The exterior $O(D_N)$ term
is absorbed because $\eta$ is at most a fixed constant. This argument
also accounts for any mismatch at the artificial central endpoints.

Finally, counting from zero more precisely gives
\[
 |e_i-\Gamma_e(i/(2N))|\le C(P\sqrt i+P^2)/N.
\]
Equality of the two densities at zero implies
$|\Gamma_s(i/(2N))-\Gamma_r(i/(2N))|\le Ci^2/N^2$.
Since $P\ge1$ and $i\ge1$, these bounds imply
\eqref{relax:eq:lfsn-initial-tangent} in the central range. Outside it,
$i^2/N^2$ is bounded below by a fixed positive constant, and the
endpoint norm bound proves the same assertion.
\end{proof}

We now estimate how the particle flow changes when its initial
configuration is interpolated. We use the regular conditional reference,
the propagated off-axis counts, and the bound for the initial derivative
proved above. We fix their logarithmic scales before choosing the later
noise windows:
\begin{equation}
 \varphi_0=(\log N)^L,\quad L>4,\qquad
 P=C\varphi_0^4,\qquad \Phi=P^8,\qquad
 \varphi=\Phi^2=P^{16},\quad \ell_N=1+\log N.
 \label{relax:eq:lfsn-ledger}
\end{equation}
Constant factors in a logarithmic parameter may be enlarged.
The endpoint weak law at height $\varphi_0^6/N$ gives
\eqref{relax:eq:lfsn-endpoint-count}: its sharp interval error is
$C\varphi_0\sqrt{|I|/N}+C\varphi_0^6\log N/N$.
The stated $P$ absorbs the second term.

Evolve each fixed quantile interpolation by
\eqref{relax:eq:lfsn-sde}, and evolve its initial reference by free
convolution. Applying Proposition~\ref{relax:prop:lfwl-continuation}
with parameter $\Phi$ is legitimate because
\eqref{relax:eq:lfsn-quantile-law} at height $\Phi^6/N$ has prefactor
$CP$, which is smaller than $\Phi^{1/4}=P^2$.
The propagated error is
\begin{equation}
 |m_{N,\alpha,u}-\widetilde m_{\alpha,u}|
 \le\frac{C(P+\Phi^{1/4})}{\sqrt{N\eta}}
 \ll\frac{\varphi}{\sqrt{N\eta}},
 \quad
 \eta\ge\Phi^6/N=\varphi^3/N,\quad u\le\varphi^{-2}.
 \label{relax:eq:lfsn-propagated}
\end{equation}
The conditional failure probability is at most
$CN^C\exp(-c\Phi^{1/2})$. The continuation proposition supplies the
horizon $\Phi^{-2}$; the displayed estimate uses the smaller horizon
$\varphi^{-2}=\Phi^{-4}$ needed in the next argument.
The reference measure may have remote atoms and depend on the
initial configurations; the conditional proof of that proposition
allows both. Its coarse region uses the particle resolvent SDE
directly, and does not require an iid law for the interpolated
particles.

\subsubsection{A conditional estimate for the tangent transform}

Now in this subsection we use a Stieltjes transform to control how the particle flow responds to changes in its initial configuration.

We keep the interpolation parameter fixed in this subsection. The
tangent-weighted Stieltjes observable and its stochastic advection
equation go back to Bourgade \cite[(1.9), Lemma~2.1 and (2.3)]{BourgadeGaps}.
We adapt the hard-edge argument of
\cite[Lemmas~4.20--4.23]{BCH}, with logarithmic losses and a bound uniform
over the terminal time and energy. The propagated weak law supplies
particle counts and a regular reference; positivity then turns transform
control into bounds for individual tangent weights. Pairing the signed
indices in the imaginary martingale retains the small real-energy
factor. The final interpolation over the initial parameter is treated
in the next subsection.

Fix an interpolation parameter and suppress it. Set
\[
 y=\varphi^5/N,\qquad \eta_*=\varphi^3/N,\qquad
 T\le\varphi^{-2}.
\]
Let $\mathcal W$ be the event that
\begin{equation}
 |m_{N,u}(E+i\eta)-\widetilde m_u(E+i\eta)|
 \le\frac{\varphi}{\sqrt{N\eta}},
 \qquad |E|\le R,\quad \eta_*\le\eta\le R,\quad u\le T.
 \label{relax:eq:lfsn-weak}
\end{equation}
The free-convolution reference obeys
\begin{equation}
 c\le\Im\widetilde m_u\le C,\qquad
 |\widetilde m_u|+|\partial_w\widetilde m_u|\le C
 \label{relax:eq:lfsn-reference}
\end{equation}
on a fixed larger neighborhood. Define the auxiliary tangent flows
by
\begin{equation}
 \partial_u v_i=\frac1{2N}\sum_{j\ne i}
                  \frac{v_j-v_i}{(x_i-x_j)^2},\qquad
 v_i^{(1)}(0)=\frac{\operatorname{sgn}(i)\sqrt{|i|}}N,\quad
 v_i^{(2)}(0)=\frac{\operatorname{sgn}(i)i^2}{N^2}.
 \label{relax:eq:lfsn-tangent}
\end{equation}
For either profile or their sum write
\[
 f_u(w)=\sum_{i\in\mathcal I_N}\frac{v_i(u)}{x_i(u)-w},
 \qquad
 B(E,t)=E\left(\sqrt{\frac N{E\vee t}}+N(E\vee t)\right).
\]
The characteristic $w_{u;t,E}$ has terminal value $E+iy$ at time
$t$ and satisfies $\partial_u w=-\widetilde m_u(w)$.

\begin{proposition}[Estimate for the tangent transform]
\label{relax:prop:lfsn-block}
There are fixed $K,c,C>0$ such that, outside an event of probability
\begin{equation}
 \mathbb P(\mathcal W^c)+CN^{20}\exp(-c\varphi^3),
 \label{relax:eq:lfsn-block-prob}
\end{equation}
one has simultaneously for $0\le u\le t\le T$ and $y\le E\le r$,
with fixed small $r>0$,
\begin{equation}
 |\Im f_u(w_{u;t,E})-\Im f_0(w_{0;t,E})|
       \le K\ell_N B(E,t).
 \label{relax:eq:lfsn-block}
\end{equation}
The constants are uniform over admissible initial configurations
and over the fixed interpolation parameter.
\end{proposition}

\begin{proof}
\textit{Counts, positivity and the initial transform.}
First, the sharp-counting argument applied to
\eqref{relax:eq:lfsn-weak} gives for an interval of length $y$ the error
\[
 O\left(\varphi\sqrt{y/N}+\eta_*\ell_N\right).
\]
Dividing by $y$ gives
$O(\varphi^{-3/2}+\ell_N\varphi^{-2})=o(1)$.
Every positive-side $y$-block in a fixed interior neighborhood
therefore contains between $cNy$ and $CNy$ particles.

The tangent evolution preserves oddness, positivity on the positive
indices, and monotonicity, and does not increase its supremum norm. For positivity, the
positive-coordinate off-diagonal coefficients are
$((x_i-x_j)^{-2}-(x_i+x_j)^{-2})/(2N)\ge0$.
At a first contact $v_i=v_{i+1}$, the terms with $j>i+1$
have two nonnegative factors in $\dot v_{i+1}-\dot v_i$,
and those with $j<i$ have two nonpositive factors. This proves
monotonicity. The maximum principle gives $|v_i|\le2$.

For $w=e+i\eta$, $e\ge0$, signed pairing gives
\begin{equation}
 f_u(w)=2\sum_{i>0}\frac{x_iv_i}{x_i^2-w^2},\qquad
 \Im f_u(w)=4e\eta\sum_{i>0}
       \frac{x_iv_i}{|x_i-w|^2|x_i+w|^2}\ge0.
 \label{relax:eq:lfsn-positive-transform}
\end{equation}
The initial rank bound $i\le CN(x_i(0)+y)$ and upper block counts
reduce the initial transform to the two integrals
\begin{align}
 \int_0^c\frac{x\sqrt{x+y}}{|x-w|^2|x+w|^2}\,dx
       &\le\frac C{\eta\sqrt H}, \nonumber\\
 \int_0^c\frac{x(x+y)^2}{|x-w|^2|x+w|^2}\,dx
       &\le C(H/\eta+\ell_N),\qquad H=e\vee\eta .
 \label{relax:eq:lfsn-initial-integrals}
\end{align}
For $e\le\eta$, scale using denominator comparable to
$(x^2+\eta^2)^2$; for $e>\eta$, split at $e/2,2e$.
The second integral has a logarithmic far tail. The first block is
bounded by its maximum since $y\le\eta$; distant particles contribute
$CNe\eta$. Consequently
\begin{equation}
 \Im f_0(e+i\eta)\le
 C\ell_N e\left(\sqrt{N/H}+NH\right).
 \label{relax:eq:lfsn-initial-transform}
\end{equation}
Symmetry and \eqref{relax:eq:lfsn-reference} imply
$|\Re\widetilde m_u(e+i\eta)|\le C|e|$. Along the characteristic,
\begin{equation}
 e_u\asymp E,\quad \eta_u\asymp y+t-u,\quad
 \max(e_u,\eta_u,u)\asymp S:=E\vee t.
 \label{relax:eq:lfsn-char-bounds}
\end{equation}
These statements include $t<y$, since $E\ge y$.

\textit{Stopping and recovery of the tangent profile.}
Let $\tau$ be the first time at which \eqref{relax:eq:lfsn-weak} or
\eqref{relax:eq:lfsn-block} fails, with $K$ a sufficiently large fixed
constant.
Conditional on the initial data every terminal characteristic is
deterministic; taking a countable dense supremum over terminal pairs
therefore defines an adapted stopping time. The fixed outer energy
strip cannot cause a tangent stop: if $E\ge r/4$, positivity,
$v_i\le2$ and the weak law give $|\Im f_u|+|\Im f_0|\le CN$,
whereas $B(E,t)\ge NE^2$. It suffices to improve the stop on an
inner strip, using blocks in a fixed slightly larger strip.

Before the stop, set the terminal time equal to the current time.
For a block $I_j=\{i>0:jy\le x_i<(j+1)y\}$, $j\ge1$,
\eqref{relax:eq:lfsn-positive-transform} gives
\begin{equation}
 \sum_{i\in I_j}v_i(u)
 \le Cy\,\Im f_u(jy+iy)
 \le CK\ell_N\,y B(jy,u).
 \label{relax:eq:lfsn-blocksum}
\end{equation}
For an individual particle, use a block two blocks later. It has
at least $cNy$ particles and all their weights are at least $v_i$.
This also handles the first block without comparing its cardinality
to that of its neighbor. Thus
\begin{equation}
 v_i(u)\le CK\ell_N\left[
 \frac{x_i+y}{\sqrt{N\max(x_i+y,u)}}+
 (x_i+y)\max(x_i+y,u)\right].
 \label{relax:eq:lfsn-profile}
\end{equation}
Small particles use later blocks inside the inner energy domain.
For particles a fixed positive distance from zero the same bound
follows directly from $v_i\le2$. Hence no circular domain enlargement
has been used.

\textit{Gradient and paired-kernel bounds.}
Differentiating the paired transform and using the block counts
gives, with $H=\max(e,\eta,u)$,
\begin{equation}
 |f'_u(w)|\le CK\ell_N^2
       \left(\sqrt{N/H}+NH\right)\frac{|w|}{\eta}.
 \label{relax:eq:lfsn-gradient}
\end{equation}
The gradient estimate uses the following two kernel integrals, whose
verification is given below:
\begin{align}
 \int_0^c\frac{x^2\,dx}
   {\sqrt{\max(x,u)}\,|x-w|^2|x+w|^2}
       &\le\frac C{\eta\sqrt H},\nonumber\\
 \int_0^c\frac{x^2\max(x,u)\,dx}{|x-w|^2|x+w|^2}
       &\le C(H/\eta+\ell_N).
 \label{relax:eq:lfsn-gradient-integrals}
\end{align}
Split at $e/2,2e$ when $e>\eta$ and additionally at $u$ on the
upper tail; when $e\le\eta$ use $(x^2+\eta^2)^2$.
The second integrand is $O(x^{-1})$ for $x>2H$, which accounts
for the extra logarithm. First-block maxima and distant particles
are absorbed. The estimates also hold on discs of radius $c\eta$
with comparable parameters.

For the imaginary martingale, the positive and negative indices must be
paired before absolute values are taken. This preserves the factor $e$
and improves the far-tail decay. The required kernel is
\[
 k_w(x)=\Im\{(w-x)^{-2}+(w+x)^{-2}\}.
\]
If $e\le2\eta$, then
$|k_w(x)|\le Ce\eta/(x^2+\eta^2)^2$. If $e>2\eta$, use
\[
 \frac{C\eta|x-e|}{((x-e)^2+\eta^2)^2}+\frac{C\eta}{e^3}
\]
on $e/2\le x\le2e$, $C\eta/e^3$ on the lower tail, and
$Ce\eta/x^4$ on the upper tail. The improved upper-tail power
is a consequence of the signed pairing. With $a=x+y$, these
bounds yield
\[
 \int_0^c\frac{a}{\sqrt{\max(a,u)}}|k_w(x)|\,dx
 \le\frac{Ce}{\eta\sqrt H},\qquad
 \int_0^c a\max(a,u)|k_w(x)|\,dx
 \le\frac{CeH}{\eta}.
\]
The same positive envelopes control blockwise suprema at zeros
of the kernel. There is no additional logarithmic loss:
\begin{equation}
 \sum_{i>0}v_i(u)|k_w(x_i(u))|
 \le CK\ell_N\frac e\eta
                  \left(\sqrt{N/H}+NH\right).
 \label{relax:eq:lfsn-paired-kernel}
\end{equation}

\textit{Verification of the kernel estimates.}
Put
\[
 \Delta_w(x)=((x-e)^2+\eta^2)((x+e)^2+\eta^2),\qquad
 R=e\vee\eta,\qquad H=e\vee\eta\vee u.
\]
Throughout, $0<y\le\eta$ and the upper integration limit $c$ is
fixed. If $e\le\eta$, then
$\Delta_w(x)\asymp(x^2+\eta^2)^2$. Scaling $x=\eta v$ gives
\[
 \int_0^c\frac{x\sqrt{x+y}}{\Delta_w(x)}\,dx
 \le C\eta^{-3/2}\int_0^\infty
       \frac{v\sqrt{v+1}}{(1+v^2)^2}\,dv
 \le C\eta^{-3/2}.
\]
For the quadratic initial profile, the same scaling gives an
integrand bounded by $Cv(v+1)^2/(1+v^2)^2$, so the integral is
at most $C(1+\log_+(c/\eta))$. If $e>\eta$, the denominator on
the lower, middle and upper intervals is respectively comparable to
$e^4$, $e^2((x-e)^2+\eta^2)$ and $x^4$. The first initial
integral has tail contributions $O(e^{-3/2})$ and middle
contribution $O((\eta\sqrt e)^{-1})$. The second has lower
contribution $O(1)$, middle contribution $O(e/\eta)$, and upper
contribution $O(\log_+(c/e))$. These prove
\eqref{relax:eq:lfsn-initial-integrals}, since
$\log_+(c/R)\le C\ell_N$ under $\eta\ge y=\varphi^5/N$.
In particular, the logarithm in the second integral comes from its
$x^{-1}$ upper tail.

For the gradient estimates, first note the elementary bound
\[
 J_0:=\int_0^c\frac{x^2}{\Delta_w(x)}\,dx
 \le\int_{\mathbb R}\frac{dx}{(x-e)^2+\eta^2}
 =\frac\pi\eta.
\]
If $u\le R$, then
$x^2/\sqrt{\max(x,u)}\le x^{3/2}$, so the first gradient
integral is bounded by the first initial integral with $y=0$.
If $u\ge R$, it is instead bounded by $J_0/\sqrt u$.
These two cases give $C/(\eta\sqrt H)$. Since
$\max(x,u)\le x+u$, the second gradient integral is bounded by
the second initial integral with $y=0$, plus $uJ_0$.
This is at most $C(H/\eta+\ell_N)$ and proves
\eqref{relax:eq:lfsn-gradient-integrals}. Its far tail is again
$O(x^{-1})$ when $x>2H$.

For the imaginary martingale kernel, direct algebra gives
\[
 k_w(x)=-4e\eta\,
 \frac{(e^2+\eta^2)^2+2(e^2-\eta^2)x^2-3x^4}
      {\Delta_w(x)^2}.
\]
This formula gives the stated $Ce\eta/(x^2+\eta^2)^2$
bound when $e\le2\eta$ and the $Ce\eta/x^4$ bound when
$x>2e>4\eta$. On $e/2\le x\le2e$, use separately the two
terms in
\[
 k_w(x)=-2\eta\left[
 \frac{e-x}{((e-x)^2+\eta^2)^2}
 +\frac{e+x}{((e+x)^2+\eta^2)^2}\right].
\]
Their absolute integral is at most $C/\eta$, because
\[
 \int_{\mathbb R}\frac{\eta|v|}{(v^2+\eta^2)^2}\,dv
 =\frac1\eta,
\]
and the second term contributes only $O(\eta/e^2)$ there.
The lower interval has kernel bound $C\eta/e^3$.

Set $a=x+y$. When $e\le2\eta$, use
$a/\sqrt{\max(a,u)}\le\sqrt{x+\eta}$ if $u\le\eta$,
and $a/\sqrt{\max(a,u)}\le(x+\eta)/\sqrt u$ otherwise.
Scaling against $(x^2+\eta^2)^{-2}$ gives the bound
$Ce/(\eta\sqrt H)$ for the first weighted kernel integral.
For the second, use
$a\max(a,u)\le(x+\eta)^2+u(x+\eta)$; scaling gives
$Ce(1+u/\eta)\le CeH/\eta$.

When $e>2\eta$, one has $a\asymp e$ in the middle interval.
Its contributions are consequently at most
$Ce/(\eta\sqrt H)$ and $CeH/\eta$.
On the lower interval, the integrated weights are bounded by
$Ce^2/\sqrt H$ and $Ce^2H$. After multiplication by
$C\eta/e^3$, these contributions are at most
$C\eta/(e\sqrt H)$ and $C\eta H/e$.
On the upper interval, $a\asymp x$ and
$\max(a,u)\ge H$, so its contributions are bounded by
\[
 \frac{Ce\eta}{\sqrt H}\int_{2e}^\infty x^{-3}\,dx
 \le\frac{C\eta}{e\sqrt H},\qquad
 Ce\eta\int_{2e}^\infty(x^2+ux)x^{-4}\,dx
 \le\frac{C\eta H}{e}.
\]
Both are dominated by the middle-interval bounds. This proves the
two weighted paired-kernel estimates. Unlike the initial and
gradient quadratic profiles, the paired quadratic integrand has an
integrable $x^{-2}$ far tail; no logarithmic loss occurs.

Finally, these integral estimates also control the block maxima
used above. Outside the first block, $x+y\asymp x$, and
$\max(x+y,u)\asymp\max(x,u)$. On $0\le x\le y$,
$\Delta_w(x)\ge cR^4$. Multiplying the relevant maxima by the
block length $y$ gives, for the initial profiles, bounds
$Cy^{5/2}/R^4$ and $Cy^4/R^4$, and, for the gradient profiles,
\[
 \frac{Cy^3}{\sqrt{\max(y,u)}R^4}
      \le\frac C{\eta\sqrt H},\qquad
 \frac{Cy^3\max(y,u)}{R^4}\le\frac{CH}{\eta}.
\]
These follow from $y\le\eta\le R$, splitting $u\ge R$ and
$u<R$ in the first inequality. For paired-kernel block maxima,
replace the middle envelope's first term by the positive bound
\[
 \frac{C\eta}{((x-e)^2+\eta^2)^{3/2}}.
\]
It still has integral $O(\eta^{-1})$ and varies on scale at least
$\eta$. The other regional envelopes agree within constants at
their boundaries. Their supremum on a block of length $y\le\eta$
is therefore controlled by a constant times their block average.
The factors $x+y$ and $\max(x+y,u)$ vary by bounded ratios on
such a block, including the first one. The upper count $CNy$ per
block now proves the stated sum estimates even at zeros of the
actual kernel. Particles outside a fixed local neighborhood
contribute at most $CNe\eta$ to the paired imaginary sums and
$CN|w|$ to the differentiated transform; these are absorbed by
the displayed bounds since $H\ge\eta$.

\textit{Exact tangent equation and the imaginary drift.}
We next derive the stochastic identity. It uses the unit diagonal
brackets of the common Brownian family; the later martingale bound uses
only the coarser bound on its covariance entries. Write $a_i=x_i-w$.
Pairing the tangent and particle interaction drifts in
$\sum_i v_i/a_i$ gives
$(2N)^{-1}\sum_{i\ne j}v_i/(a_i^2a_j)$.
The It\^o contribution is
$(2N)^{-1}\sum_i v_i/a_i^3$ because
$d\langle b_i\rangle=du$. Together they equal $m_N f'$.
Thus along the deterministic characteristic,
\begin{equation}
 df_u(w_u)=(m_N-\widetilde m_u)(w_u)f'_u(w_u)\,du
 -\frac1{\sqrt{2N}}\sum_i
      \frac{v_i(u)}{(w_u-x_i(u))^2}\,db_i(u).
 \label{relax:eq:lfsn-exact-sde}
\end{equation}
This also follows from BCH Lemma 4.21 or equation (4.74); its general cubic
correction is proportional to
$d\langle b_i\rangle/du-1$ and is exactly zero here.
Negative-index correlations add no mixed It\^o term because
the observable is a sum of single-coordinate functions.

The imaginary part of the drift retains a necessary factor $E$.
For a symmetric configuration, both
$(m_N-\widetilde m)(i\eta)$ and $f'(i\eta)$ are purely imaginary,
so their product is real. If $e\le\eta/4$, Cauchy's derivative
estimate on discs of radius $c\eta$, followed by integration from
the imaginary axis to $e+i\eta$, supplies a factor $e/\eta$.
These discs stay above $y$: the characteristic has
$e\ge(1-o(1))E\ge(1-o(1))y$ and $\eta\ge4e$.
The discs can cross the imaginary axis, where the two-sided weak
law and symmetry apply. If $e>\eta/4$, use $|w|\le Ce$ directly
in \eqref{relax:eq:lfsn-gradient}. Combining the two regimes gives
\[
 |\Im[(m_N-\widetilde m_u)f'_u](w_u)|
 \le CK\ell_N^2\frac{\varphi E}{\sqrt N\,\eta_u^{3/2}}
                    \left(\sqrt{N/S}+NS\right).
\]
Since $-\dot\eta_u\asymp1$, the total stopped drift is at most
\begin{equation}
 CK\ell_N^2 B(E,t)\frac{\varphi}{\sqrt{Ny}}.
 \label{relax:eq:lfsn-drift}
\end{equation}
Relative to $K\ell_N B$, this is
$O(\ell_N\varphi^{-3/2})=o(1)$.

\textit{Bracket, maximal estimate and simultaneous terminal parameters.}
For each fixed terminal pair $(t,E)$, let $M$ denote the imaginary
part of the stopped stochastic integral in
\eqref{relax:eq:lfsn-exact-sde}. Pairing $\pm i$ in that integral gives
the coefficient
$v_i[(w-x_i)^{-2}+(w+x_i)^{-2}]/\sqrt{2N}$.
Even the coarse covariance bound
$|d\langle b_i,b_j\rangle|\le du$ and
\eqref{relax:eq:lfsn-paired-kernel} imply
\begin{equation}
 \langle M\rangle_t
 \le\frac C N\int_0^{t\wedge\tau}
       \left(\sum_{i>0}v_i|k_{w_u}(x_i)|\right)^2du
 \le\frac{CK^2\ell_N^2 B(E,t)^2}{Ny}.
 \label{relax:eq:lfsn-bracket}
\end{equation}
The continuous-martingale maximal inequality, for a fixed terminal
pair $(t,E)$, therefore gives
\begin{equation}
 \mathbb P\left\{\sup_{s\le t}|M_s|>
            K\ell_N B(E,t)\varphi^{-1}\right\}
       \le2e^{-cNy\varphi^{-2}}=2e^{-c\varphi^3}.
 \label{relax:eq:lfsn-martingale-tail}
\end{equation}
Here the terminal pair is held fixed; the stochastic integral is
a martingale in its physical integration time.

The preceding maximal estimate holds in physical integration time for a
fixed terminal pair. We now put all terminal pairs on one event without
introducing a further union over physical time. Use mesh $N^{-10}$ in the
terminal time and energy, including endpoints. Characteristic maps
are uniformly Lipschitz by \eqref{relax:eq:lfsn-reference};
$|f'|\le CN/y^2\le CN^3$ controls interpolation at every physical
time. The function $B$ is Lipschitz with constant $CN$, as seen
from
\[
 B(E,t)=
 \begin{cases}
  \sqrt{NE}+NE^2,&t\le E,\\
  E\sqrt{N/t}+NEt,&t>E.
 \end{cases}
\]
Also $B\ge E\sqrt N\ge N^{-1/2}$.
Hence the interpolation error $CN^{-7}$ is negligible relative
to the stopping threshold. Round terminal time upwards so that
the current physical time lies in the approximating interval.
For boundary energies work first on a fixed slightly larger
neighborhood. There are at most $CN^{20}$ terminal pairs, and
\eqref{relax:eq:lfsn-martingale-tail} already controls all physical times.
No additional time union is needed.

Combining drift, martingale and mesh bounds gives before the stop
\[
 |\Im f_s(w_s)-\Im f_0(w_0)|
 \le CK\ell_N B(E,t)
     (\ell_N\varphi^{-3/2}+\varphi^{-1})+CN^{-7}
 <\tfrac12K\ell_N B(E,t)
\]
for all sufficiently large $N$. The tangent boundary cannot be
reached first. Intersecting with $\mathcal W$ removes the other
boundary and proves the proposition.
\end{proof}

\subsubsection{Integration over the initial interpolation parameter}

The purpose of this subsubsection is to show that integrating the response to changes in the initial data gives the required comparison between the two flows.

The particle bound below has the index and time profile of
\cite[Proposition~4.6, equation~(4.10)]{BCH}. Their estimate has an
$N^\varepsilon$ loss and starts at a fixed-power time
$N^{-1+\omega_t}$. Here the loss is $\varphi^{15}$ on
$\varphi^{100}/N\le t\le\varphi^{-2}$, with the logarithmic parameter
specified in \eqref{relax:eq:lfsn-ledger}. The present tangent argument
uses independent unit Brownian drivers; BCH also treats correlated
continuous martingales with diagonal brackets at most one.

The block estimate holds for each fixed initial interpolation parameter,
with uniform constants. We obtain the tangent bound, construct a
jointly measurable interpolation, and integrate in that parameter.
The exceptional parameters are controlled by their Lebesgue measure;
the proof does not require a common good event for every real parameter.

\begin{theorem}[Same-noise relaxation at logarithmic time]
\label{relax:thm:lfsn-relaxation}
Assume the endpoint count, density, norm and initial independence
conditions above, with the parameter choices
\eqref{relax:eq:lfsn-ledger}. Couple the two endpoint flows using the
same positive-index Brownian motions. For every fixed $D>0$,
with probability at least $1-N^{-D}$,
\begin{equation}
 |s_i(t)-r_i(t)|
 \le C_D\varphi^{15}\frac iN
 \left[\frac1{\sqrt{Nt}}+\max\{i/N,t\}\right],
 \quad
 1\le i\le N,\qquad
 \frac{\varphi^{100}}N\le t\le\varphi^{-2}.
 \label{relax:eq:lfsn-relaxation}
\end{equation}
Constants are uniform over the fixed spatial bulk and the admissible
moment class. In the matrix application with this independent reference
endpoint, the comparison path is independent of
$\mathcal F_0^{\rm mat}$.
\end{theorem}

\begin{proof}
\textit{Particle bounds and domination of the actual tangent.}
For each deterministic $\alpha$, apply
Proposition~\ref{relax:prop:lfsn-block} and
\eqref{relax:eq:lfsn-propagated}.
Within the fixed central index range, uniform counts imply
$x_i(t)\asymp i/N$ for $i\ge K_0Ny$,
where $K_0$ is fixed and sufficiently large. A later block with at
least $cNy$ particles, monotonicity and
\eqref{relax:eq:lfsn-blocksum} yield
\[
 v_i^{(1)}(t)+v_i^{(2)}(t)
 \le C\ell_N\frac iN
       [(Nt)^{-1/2}+\max\{i/N,t\}].
\]
For $i<K_0Ny$, compare with
$i_0=\lceil K_0Ny\rceil=O(\varphi^5)$.
Since $t\ge\varphi^{100}/N\gg i_0/N$, this costs at most
$C\varphi^5$, including $i=1$. Outside the fixed central index
range the maximum principle gives the same conclusion.
For completeness, initial tangent domination is applied on the
positive coordinates. Put
$c_{ij}=(2N(x_i-x_j)^2)^{-1}$. An odd tangent $q$ solves there
\[
 \dot q_i=\sum_{\substack{j>0\\j\ne i}}
 (c_{ij}-c_{i,-j})(q_j-q_i)
       -2\sum_{j>0}c_{i,-j}q_i.
\]
All off-diagonal coefficients are nonnegative, since
$|x_i-x_j|\le x_i+x_j$. This is a positive semigroup with a
nonnegative killing term. Apply comparison separately to $q$ and
$-q$, using \eqref{relax:eq:lfsn-initial-tangent}. Since
$CP^2\le\varphi$ for large $N$, the absolute actual tangent on
positive indices is at most $\varphi(v^{(1)}+v^{(2)})$.
Its fixed-$\alpha$ estimate therefore has prefactor
$C\ell_N\varphi^6$, smaller than the conservative $\varphi^{15}$.

\textit{A common measurable version.}
We specify the interpolation version before integrating the estimates. Construct the flows first for rational
$\alpha$, including the endpoints. On a common probability-one
event, all these flows exist and same-noise maximum contraction gives
\[
 \sup_{t\le\varphi^{-2}}\sup_i
 |x_i(t,\alpha)-x_i(t,\beta)|
 \le|\alpha-\beta|\sup_i|s_i(0)-r_i(0)|.
\]
Extend by uniform continuity in $\alpha$ to every parameter as a
continuous path. On the endpoint norm event this is a uniformly
Lipschitz map. For each fixed $\alpha$, it agrees almost surely
with that parameter's noncolliding strong solution: local continuous
dependence applies on its compact path, whose minimum gaps are
positive. The same localization gives differentiability in the
initial parameter and identifies the derivative with
\eqref{relax:eq:lfsn-tangent}. Fubini makes these identifications valid
for almost every $\alpha$ on one probability-one event. No claim
of simultaneous noncollision for every real $\alpha$ is needed.

\textit{Exceptional parameters and scalar integration.}
Let $\mathcal G_\alpha$ denote the fixed-parameter estimate,
simultaneously in the required indices and times, using this joint
measurable version. For every fixed $D_0$,
$\mathbb P(\mathcal G_\alpha^c)\le N^{-D_0}$ uniformly in
$\alpha$. Fubini and Markov give
\[
 \mathbb P\{|\{\alpha:\mathcal G_\alpha^c\}|>N^{-4}\}
 \le N^{4-D_0}.
\]
For each rational time, the scalar Lipschitz endpoint map is
absolutely continuous, so its fundamental theorem of calculus gives
\[
 s_i(t)-r_i(t)=\int_0^1\partial_\alpha x_i(t,\alpha)\,d\alpha.
\]
Take the countable intersection over rational times and all indices.
The identified tangent paths are continuous in time for almost every
$\alpha$ and bounded by the fixed endpoint norm difference through
maximum contraction. Dominated convergence therefore extends the
identity to every time simultaneously. This uses scalar absolute
continuity, without requiring a differentiability theorem for maps
into a space of continuous paths.
The bad set contributes at most $CN^{-4}$, uniformly in $i,t$.
It is absorbed because the target without the logarithmic factor
is at least $N^{-3/2}$ for $i\ge1$, $t\le1$.
Choosing $D_0>D+5$ and including the initial exceptional events
proves \eqref{relax:eq:lfsn-relaxation}.

\textit{Independence of the comparison path.}
For independence, let $\mathcal F_0^{\rm mat}$ contain the original
initial matrix and let the actual singular-value Brownian vector be
Brownian relative to a filtration containing it. The vector is
therefore independent of $\mathcal F_0^{\rm mat}$.
Add a fresh scaled Ginibre initial configuration independent of
that filtration. The comparison path is its strong solution with
the common Brownian vector, hence is independent of
$\mathcal F_0^{\rm mat}$. Conditional estimates use the enlarged
initial field
\[
 \mathcal G_0=\mathcal F_0^{\rm mat}
        \vee\sigma(\text{fresh reference Ginibre configuration}).
\]
The Brownian vector is independent of $\mathcal G_0$, but the
comparison path is not, since its own initial value belongs to
$\mathcal G_0$. For the application at two spatial points, take the common Brownian family to be the
corresponding reference family from Theorem~\ref{noise:thm:lfn-whitening}.
Its established independence from this enlarged initial field supplies
the same conclusion. The reference process still uses its own initial
spectrum, so the distinction between the two initial fields is unchanged.
\end{proof}

\paragraph{Time range and spatial families.}
The chosen time lies in the required interval. By
\eqref{relax:eq:lfsn-ledger}, we have $P=C(\log N)^{4L}$ and
$\varphi=P^{16}=C^{16}(\log N)^{64L}$.
If $\tau=(\log N)^b/N$, choose the fixed exponent $b$ after $L$
with $b>6400L$. The strict inequality absorbs every fixed
multiplicative constant for sufficiently large $N$.
Then $\tau\ge\varphi^{100}/N$, whereas asymptotically
$\tau\ll\varphi^{-2}$ and $\tau\ll N^{-a}$ for every fixed
$a<1$. Thus the logarithmic relaxation fits inside the available
polynomial Gaussian age. The theorem asserts a coupling at one spatial point,
uniform in $i,t$; a polynomial family of deterministic spatial
points is covered by choosing a larger failure exponent before taking
the union bound. It does not
assert a simultaneous auxiliary weighted-flow estimate for every
interpolation parameter, nor a continuum of independently coupled
spatial points.

\endgroup

\begingroup
\newcommand{\E}{\mathbb E}
\newcommand{\Pp}{\mathbb P}
\newcommand{\cE}{\mathcal E}
\newcommand{\cB}{\mathcal B}
\newcommand{\cS}{\mathcal S}
\newcommand{\cL}{\mathcal L}
\newcommand{\dd}{\,\mathrm d}
\newcommand{\tr}{\operatorname{Tr}}
\section{Different-noise relaxation with adapted moving weights}\label{chap:different}

The purpose of this section is to show that a small change in the driving noise causes only a small displacement of the central singular values.

\paragraph{Noise interpolation.}
We use the noise interpolation and short-interaction comparison of
\cite[Proposition~4.5]{BCH} for signed singular-value systems with a
common initial configuration and correlated martingale drivers. BCH
uses a fixed-power index window and a fixed-power covariance bound.
Here \eqref{different:eq:noise} permits a covariance-operator bound
$Ct^{6k}e^{-2t^{1/3}}$ on a logarithmic window, with
$t=\tfrac12\log N$. A moving spatial weight controls the central
particles, while global dissipative estimates bound the remaining
indices. Its mixed covariation with the difference noise enters the
energy identity. Stopping and parameter integration yield the
$N^{-1}t^{-R}$ comparison in \eqref{different:eq:conclusion}, uniformly
up to time $\tau=t^b/N$, with overwhelming
probability.

\subsection{Hypotheses and statement}

We first state the counting and covariance assumptions that make two singular-value flows with the same initial spectrum comparable.
Let
\[
 \mathcal I_N=\{-N,\ldots,-1,1,\ldots,N\},\qquad
 t=\tfrac12\log N,\qquad H=t^b,\qquad \tau=H/N,
\]
where $b>0$ is fixed. There is no zero index. All processes are adapted to one filtration satisfying the usual conditions. Let $b_i,\beta_i$, $1\le i\le N$, be real continuous local martingales starting from zero, with absolutely continuous brackets whose diagonal densities are at most one. Extend them antisymmetrically:
\[
 b_{-i}=-b_i,\qquad \beta_{-i}=-\beta_i.
\]
No independence between coordinates or between the two driver families is assumed in this theorem. In the intended application, each positive-index family separately is a standard Brownian vector.

Let $x(0)$ be an initial configuration measurable at time zero, with
\[
 0<x_1(0)<\cdots<x_N(0),\qquad x_{-i}(0)=-x_i(0).
\]
For $\alpha\in[0,1]$ set
\[
 b_i^\alpha=\alpha b_i+(1-\alpha)\beta_i,
 \qquad
 \chi_i^\alpha=\frac{b_i^\alpha}{\sqrt{2N}},
 \qquad
 \xi_i=\frac{b_i-\beta_i}{\sqrt{2N}}.
\]
The systems to be compared solve
\begin{equation}
 \dd x_i^\alpha
 =\dd\chi_i^\alpha+
 \frac1{2N}\sum_{j\ne i}\frac{\dd s}{x_i^\alpha-x_j^\alpha},
 \qquad x_i^\alpha(0)=x_i(0).
 \label{different:eq:system}
\end{equation}
We impose three hypotheses: existence and noncollision for each deterministic interpolation parameter, central counting estimates, and a small adapted covariance for the driver difference.

\paragraph{(WP): bounded-bracket well-posedness.}
For each deterministic $\alpha$, equation \eqref{different:eq:system} has a unique nonexplosive strong solution that remains strictly ordered and positive on positive indices throughout $[0,\tau]$. The solution is obtained by localization of the locally Lipschitz drift in the open signed Weyl chamber. This input is supplied by \cite[Appendix B, Proposition B.1]{BCH}. Indeed, Cauchy--Schwarz for predictable covariations gives
\begin{equation}
 \frac{\dd[b_i^\alpha]}{\dd s}
 \le\left(\alpha\sqrt{\frac{\dd[b_i]}{\dd s}}
 +(1-\alpha)\sqrt{\frac{\dd[\beta_i]}{\dd s}}\right)^2\le1.
 \label{different:eq:diagonal}
\end{equation}
The cited statement permits correlated continuous martingale drivers; it does not require the interpolated diagonal bracket to equal one.

\paragraph{(REG): off-axis counting and nested central blocks.}
There are fixed $a>0$, $0<\kappa_1<\kappa_2<1$, and positive constants $c_0,C_0,r_0$, such that, with $Q_0=t^a$ and
\[
 I_1=\{i:|i|\le\kappa_1N\},\qquad
 I_2=\{i:|i|\le\kappa_2N\},
\]
the following properties hold up to an adapted stopping time $\rho_\alpha$:
\begin{align}
 c_0\frac{|i-j|}{N}
 &\le |x_i^\alpha(s)-x_j^\alpha(s)|
 \le C_0\frac{|i-j|}{N},
 &&i,j\in I_2,\quad |i-j|\ge Q_0, \label{different:eq:reg-gap}\\
 |x_i^\alpha(s)-x_j^\alpha(s)|&\ge c_0,
 &&i\in I_1,\quad j\notin I_2, \label{different:eq:reg-far}\\
 |x_i^\alpha(s)|&\ge r_0,
 &&i\notin I_1. \label{different:eq:reg-outside}
\end{align}
For every fixed $D>0$,
\begin{equation}
 \sup_{\alpha\in[0,1]}\Pp(\rho_\alpha\le\tau)
 \le C_DN^{-D}.
 \label{different:eq:reg-prob}
\end{equation}
The stopping inequalities may be taken with fixed slack relative to the counting theorem. The constants and exponents in \eqref{different:eq:reg-gap}--\eqref{different:eq:reg-outside} are fixed before the later windows $L,K$ are selected. These are the conclusions of the bounded-bracket weak-law/counting results, namely Propositions~\ref{relax:prop:lfwl-aged} and~\ref{relax:prop:lfwl-continuation}, Corollary~\ref{relax:cor:lfwl-counts}, and its following noise-interpolation argument. An imaginary-axis trace estimate alone is not asserted to imply this input.

\paragraph{(NOISE): an adapted low-index covariance estimate.}
For $K=t^k$, there is a stopping time $\rho_{\rm n}$ such that the positive-index covariance of $b-\beta$ on $1\le i\le K$ satisfies
\begin{equation}
 \left\|\frac{\dd[(b-\beta)_{1:K}]}{\dd s}\right\|_{\rm op}
 \le\delta_N,
 \qquad
 \delta_N\le Ct^{6k}e^{-2d_N},\qquad d_N\ge t^{1/3},
 \quad s<\rho_{\rm n}.
 \label{different:eq:noise}
\end{equation}
For every fixed $D>0$, $\Pp(\rho_{\rm n}\le\tau)\le C_DN^{-D}$. The signed extension of \eqref{different:eq:noise} has operator norm at most $2\delta_N$; no inverse of its singular signed covariance is used.

\begin{theorem}[Comparison for different driving noises]\label{different:thm:main}
Assume \textup{(WP)}, \textup{(REG)}, and \textup{(NOISE)}. Fix $R>0$, and choose fixed exponents in the order
\begin{equation}
 \ell>\max\{a+1,\tfrac32b+R+5\},\qquad
 L=t^\ell,\qquad k>\ell+2,\qquad K=t^k.
 \label{different:eq:order}
\end{equation}
Then, for every fixed $D>0$ and all sufficiently large $N$,
\begin{equation}
 \Pp\left(
 \sup_{0\le s\le\tau}\max_{|i|\le L}
 |x_i^1(s)-x_i^0(s)|>N^{-1}t^{-R}
 \right)\le C_DN^{-D}.
 \label{different:eq:conclusion}
\end{equation}
Constants and the sufficiently large $N$ threshold may depend on the fixed exponents and the probability exponent $D$. No subsequent enlargement of $k$ or $\ell$ changes the earlier counting exponent $a$.
\end{theorem}

In the application, $x^1$ is the actual singular-value flow and $x^0$ uses the reference Brownian family constructed in Section~\ref{chap:noise}, with the same initial singular values. Changing this common initial spectrum to an independent Gaussian reference is the separate same-noise comparison in Section~\ref{chap:relax}.

\subsection{Tangents and global dissipative estimates}

In this subsection we first bound the response to changing the noise across all indices, so distant interactions can be controlled.
For a fixed deterministic $\alpha$, localization away from collisions gives the tangent $u_i=\partial_\alpha x_i^\alpha$:
\begin{equation}
 \dd u_i=\dd\xi_i+\sum_{j\ne i}c_{ij}(u_j-u_i)\dd s,
 \qquad
 c_{ij}=\frac1{2N(x_i^\alpha-x_j^\alpha)^2},\qquad u(0)=0.
 \label{different:eq:tangent}
\end{equation}
The rigorous identification and integration in $\alpha$ are given in Section~\ref{different:sec:alpha}. Write $\cB$ for the generator in \eqref{different:eq:tangent}. Its conductances are symmetric and nonnegative. The diagonal covariance bounds imply, for every $i$,
\begin{equation}
 \frac{\dd[\chi_i^\alpha]}{\dd s}\le\frac1{2N},\qquad
 \frac{\dd[\xi_i]}{\dd s}\le\frac2N,\qquad
 \left|\frac{\dd[\chi_i^\alpha,\xi_i]}{\dd s}\right|\le\frac CN.
 \label{different:eq:global-brackets}
\end{equation}
Before $\rho_{\rm n}$, for $|i|\le K$ the last two bounds improve to
\begin{equation}
 \frac{\dd[\xi_i]}{\dd s}\le\frac{C\delta_N}{N},\qquad
 \left|\frac{\dd[\chi_i^\alpha,\xi_i]}{\dd s}\right|
 \le\frac{C\sqrt{\delta_N}}N.
 \label{different:eq:low-mixed}
\end{equation}

The global estimates below use only dissipativity and diagonal bracket bounds. For the local comparison we also separate interactions across index distances at least $L$: counting makes their total rate small enough to control by the global supremum estimate. Remove every unordered pair with $|i-j|\ge L$ and at least one endpoint in $I_1$. Let $\cS$ retain all other pairs, and put $\cL=\cB-\cS$. Define the short-interaction tangent by
\begin{equation}
 \dd w=\cS w\dd s+\dd\xi,\qquad w(0)=0.
 \label{different:eq:short-tangent}
\end{equation}
Both $\cB$ and $\cS$ are symmetric Markov generators. The supremum tail will control the removed interactions, while the $L^2$ moment will later control exceptional interpolation parameters.

The supremum estimate uses the dissipative argument of
\cite[Lemma~4.12]{BCH}. The proof below gives the stopped maximal form and the
energy moment needed for the later parameter integration.

\begin{lemma}\label{different:lem:global}
For $y=u$ or $y=w$,
\begin{align}
 \Pp\left(\sup_{s\le\tau}\|y(s)\|_\infty>A\right)
 &\le4N\exp(-cA^2N/\tau), \label{different:eq:cosh-tail}\\
 \E\sup_{s\le\tau}\|y(s)\|_2^2&\le C\tau.
 \label{different:eq:l2-moment}
\end{align}
In particular, the bound $\sup_{s\le\tau}\|w(s)\|_\infty\le t^2\sqrt H/N$ fails with probability at most $4Ne^{-ct^4}$.
\end{lemma}
\begin{proof}
For $\theta>0$, the interaction part of It\^o's formula for $\sum_i\cosh(\theta y_i)$ is nonpositive, since each pair contributes a negative multiple of
\[
 (y_i-y_j)\bigl(\sinh(\theta y_i)-\sinh(\theta y_j)\bigr)\ge0.
\]
The diagonal It\^o correction is at most $C\theta^2N^{-1}\sum_i\cosh(\theta y_i)$. Consequently
\[
 e^{-C\theta^2s/N}\sum_i\cosh(\theta y_i(s))
\]
is a nonnegative supermartingale, first under collision/size localization and then by Fatou. The maximal inequality gives $4N\exp(-\theta A+C\theta^2\tau/N)$; optimizing $\theta$ proves \eqref{different:eq:cosh-tail}. Cross-coordinate correlations do not enter this separable function.

For \eqref{different:eq:l2-moment}, the interaction in $\dd\|y\|_2^2$ is nonpositive. The covariance matrix $Q_\xi=\dd[\xi]/\dd s$ is positive semidefinite and has $\tr Q_\xi\le C$. The scalar energy martingale has bracket at most $C\|y\|_2^2\dd s$, because $\|Q_\xi\|_{\rm op}\le\tr Q_\xi$. The Burkholder--Davis--Gundy (BDG) inequality and Young's inequality yield
\[
 \E\sup_{s\le\tau}\|y(s)\|_2^2
 \le C\tau+\tfrac12\E\sup_{s\le\tau}\|y(s)\|_2^2.
\]
Localization and Fatou complete the proof.
\end{proof}

The following is the short/long-interaction Duhamel comparison in
\cite[(4.24)--(4.26)]{BCH}, with the present logarithmic parameters.

\begin{lemma}[Long-interaction error]\label{different:lem:long}
On the counting event, except with probability $4Ne^{-ct^4}$,
\begin{equation}
 N\sup_{s\le\tau}\|u(s)-w(s)\|_\infty
 \le C\frac{t^2H^{3/2}}L.
 \label{different:eq:long-error}
\end{equation}
\end{lemma}
\begin{proof}
Equations \eqref{different:eq:reg-gap}--\eqref{different:eq:reg-far} imply, on every removed pair, $c_{ij}\le CN/|i-j|^2$. Hence
\[
 \sup_i\sum_{j:\,\text{removed}}c_{ij}\le CN/L,
 \qquad \|\cL w\|_\infty\le CN\|w\|_\infty/L.
\]
The difference $e=u-w$ satisfies the ordinary equation $\dot e=\cB e+\cL w$ with zero initial data. The homogeneous propagator of $\cB$ is an $\ell^\infty$ contraction. Its ordinary Duhamel integral and Lemma~\ref{different:lem:global} give
\[
 \sup_{s\le\tau}\|e(s)\|_\infty
 \le C\tau\frac NL\frac{t^2\sqrt H}{N}.
\]
This is \eqref{different:eq:long-error}. The propagator is used only in this ordinary integral, after the common additive noise has canceled.
\end{proof}

\subsection{The moving weight and its exact energy identity}

In this subsection it is shown that a weight that follows the particles isolates the central region while keeping every stochastic correction in the energy identity.
The tangent for the short-interaction system still receives noise at every index.
Following the moving-weight argument in
\cite[Lemma~4.13, equations~(4.27)--(4.35)]{BCH}, we use a weight that
stays bounded below on the target block $|i|\le L$ and decays
exponentially before the boundary $|i|=K$ of the covariance estimate.
The weight follows particle positions; its stochastic motion and mixed
covariation with the tangent forcing enter the product rule. We retain
both terms when estimating the drift and martingale bracket at the
logarithmic scales.

Put $\zeta=N/L$. Choose a smooth nondecreasing $f:[0,\infty)\to[0,\infty)$ that equals its argument on $[0,r_0/2]$ and is constant for arguments at least $r_0$. Set
\[
 \psi_\zeta(x)=f\bigl(\sqrt{x^2+\zeta^{-2}}\bigr),
 \qquad \phi_i=e^{-\zeta\psi_\zeta(x_i^\alpha)},
 \qquad v_i=\phi_iw_i,
 \qquad \cE=\sum_i v_i^2.
\]
Then $\psi_\zeta$ is even and smooth,
\begin{equation}
 |\psi_\zeta'|\le C,\qquad |\psi_\zeta''|\le C\zeta.
 \label{different:eq:psi-derivatives}
\end{equation}
Near zero it is exactly $\sqrt{x^2+\zeta^{-2}}$, which smooths the absolute value on the scale $L/N$. It is constant outside the fixed central neighborhood, where its derivatives vanish.

By symmetry, counting, and monotonicity,
\begin{equation}
 \phi_i\ge c\quad(|i|\le L),\qquad
 \phi_i\le e^{-cK/L}\quad(|i|\ge K).
 \label{different:eq:weight-support}
\end{equation}
For retained pairs on which the weights differ, at least one endpoint is in $I_1$ and $|i-j|<L$. Both endpoints are then in $I_2$ for large $N$. At distances below $Q_0$, monotonicity compares the gap with one of index distance at least $Q_0$. Thus $|x_i-x_j|\le CL/N$ on every such pair, and
\begin{equation}
 |\log(\phi_i/\phi_j)|\le C,
 \qquad
 \frac{(\phi_i-\phi_j)^2}{\phi_i\phi_j}
 \le C\zeta^2(x_i-x_j)^2.
 \label{different:eq:weight-ratio}
\end{equation}
For retained pairs whose endpoints both lie outside $I_1$, the two weights are equal and both derivatives $\psi_i',\psi_j'$ vanish.

Suppress $\alpha$, write $\Delta_{ij}=x_i-x_j$, $a_i=\phi_i^2$, and $\psi_i'=\psi_\zeta'(x_i)$. The derivatives are
\[
 a_i'=-2\zeta\psi_i'a_i,
 \qquad a_i''=\bigl(4\zeta^2(\psi_i')^2-2\zeta\psi_i''\bigr)a_i.
\]
The product rule for $a_iw_i^2$ separates the short-interaction drift, the drift of the particle-dependent weight, the two pure quadratic-variation terms, the mixed covariation, and the martingale:
\begin{align}
 \dd\cE={}&D_{\rm sh}\dd s+D_{\rm mov}\dd s+
 \sum_i a_i\dd[\xi_i]
 +\sum_i\bigl(2\zeta^2(\psi_i')^2-\zeta\psi_i''\bigr)v_i^2\dd[\chi_i]
 \nonumber\\
 &\quad-4\zeta\sum_i\psi_i'\phi_i v_i\dd[\chi_i,\xi_i]
 +\dd M,\label{different:eq:exact-energy}\\
 D_{\rm sh}={}&-2\sum_{i<j:\,\mathrm{ret}}c_{ij}
 \bigl[(v_i-v_j)^2-w_iw_j(\phi_i-\phi_j)^2\bigr],\label{different:eq:short-drift}\\
 D_{\rm mov}={}&-\frac\zeta N\sum_{i<j}
 \frac{\psi_i'v_i^2-\psi_j'v_j^2}{\Delta_{ij}},\label{different:eq:moving-drift}\\
 \dd M={}&2\sum_i\phi_i v_i\dd\xi_i
 -2\zeta\sum_i\psi_i'v_i^2\dd\chi_i.\label{different:eq:energy-martingale}
\end{align}
The factor $-4$ in the mixed drift comes from $2a_i'w_i\dd[x_i,w_i]$. It records the correlation between the interpolated driver and the driver difference; independence of reference noises across sites does not remove it.

\subsection{Bounds for the terms in the energy identity}

In this subsection we bound the energy drift and its random fluctuations while retaining the correlations created by the moving weight. We let
\begin{equation}
 \nu_N=K\delta_N+Ne^{-cK/L}.
 \label{different:eq:nu}
\end{equation}
Constants denoted by $c$ in exponential tails may decrease in subsequent estimates.

\begin{proposition}[Moving-weight energy bound]\label{different:prop:energy}
Before $\rho_\alpha\wedge\rho_{\rm n}$, the energy satisfies
\begin{align}
 \dd\cE&\le C\frac{N\log N}{L}\cE\dd s+
 C\frac{\nu_N}{N}\dd s+\dd M,\label{different:eq:LF5-drift}\\
 \dd[M]&\le C\left(\frac{N}{L^2}\cE^2+
 \frac{\nu_N}{N}\cE\right)\dd s.\label{different:eq:LF5-bracket}
\end{align}
\end{proposition}
\begin{proof}
The short-interaction dissipation absorbs the singular part of the moving drift. We bound the remaining drift terms and then the bracket of the complete energy martingale, retaining all mixed covariations.

\paragraph{Short-interaction commutator.}
The identity underlying \eqref{different:eq:short-drift} is
\[
 (\phi_i^2w_i-\phi_j^2w_j)(w_i-w_j)
 =(v_i-v_j)^2-w_iw_j(\phi_i-\phi_j)^2.
\]
By \eqref{different:eq:weight-ratio}, the commutator has absolute value at most
\[
 C\frac{\zeta^2}{N}
 \sum_{\substack{i<j:\,\mathrm{ret}\\\phi_i\ne\phi_j}}(v_i^2+v_j^2)
 \le C\frac NL\cE,
\]
because at most $2L$ active neighbors occur for each index. The negative dissipation retains the size
\[
 \sum_{i<j:\,\mathrm{ret}}\frac{(v_i-v_j)^2}{N\Delta_{ij}^2}.
\]

\paragraph{Short part of the moving-weight drift.}
Use the exact decomposition
\[
 \psi_i'v_i^2-\psi_j'v_j^2
 =\tfrac12(\psi_i'-\psi_j')(v_i^2+v_j^2)
 +\tfrac12(\psi_i'+\psi_j')(v_i-v_j)(v_i+v_j).
\]
The first part costs $C\zeta^2(v_i^2+v_j^2)/N$ per pair, by \eqref{different:eq:psi-derivatives}. Young's inequality absorbs the second part into any fixed small portion of the preceding dissipation, at the same remainder cost. Summation gives $CN\cE/L$.

\paragraph{Long part of the moving-weight drift.}
On removed pairs, \eqref{different:eq:reg-gap}--\eqref{different:eq:reg-far} yield $|\Delta_{ij}|^{-1}\le CN/|i-j|$. Hence
\[
 |D_{\rm mov,long}|
 \le C\frac\zeta N\sum_i v_i^2
       \sum_{j:\,\mathrm{removed}}\frac1{|\Delta_{ij}|}
 \le C\zeta\log N\,\cE.
\]
Retained pairs outside $I_1$ contribute zero, so this covers the complete moving drift.

\paragraph{Pure weight It\^o term and additive noise.}
Equations \eqref{different:eq:global-brackets} and \eqref{different:eq:psi-derivatives} bound the pure weight correction by $C\zeta^2\cE/N=CN\cE/L^2$. The additive trace term obeys
\[
 \sum_i\phi_i^2\frac{\dd[\xi_i]}{\dd s}
 \le\frac C N\left(K\delta_N+Ne^{-cK/L}\right)
 =C\nu_N/N.
\]

\paragraph{The mixed weight/noise drift.}
Inside $|i|\le K$, \eqref{different:eq:low-mixed} and Cauchy--Schwarz give
\[
 4\zeta\sum_{|i|\le K}|\psi_i'\phi_i v_i|
 \left|\frac{\dd[\chi_i,\xi_i]}{\dd s}\right|
 \le C\frac{\sqrt{K\delta_N\cE}}L.
\]
Outside this block, use \eqref{different:eq:global-brackets} and
$\sum_{|i|>K}\phi_i^2\le CNe^{-cK/L}$. The combined bound is
\begin{equation}
 C\frac{\sqrt{\nu_N\cE}}L
 \le C\frac NL\cE+C\frac{\nu_N}{NL}
 \le C\frac NL\cE+C\frac{\nu_N}{N}.
 \label{different:eq:cross-payment}
\end{equation}
The low-block factor $\sqrt{\delta_N}$ ensures that this mixed contribution has the same small forcing scale $\nu_N/N$ as the additive noise.

\paragraph{Martingale bracket.}
Write $M=M_{\rm n}+M_{\rm w}$ according to \eqref{different:eq:energy-martingale}. Positive semidefiniteness and the diagonal bounds give $|\dd[\chi_i,\chi_j]/\dd s|\le C/N$, even for arbitrary coordinate correlations. Thus
\[
 \frac{\dd[M_{\rm w}]}{\dd s}
 \le C\frac{\zeta^2}{N}\left(\sum_i v_i^2\right)^2
 =C\frac N{L^2}\cE^2.
\]
For the low-index part of $M_{\rm n}$ use the operator covariance estimate \eqref{different:eq:noise}; it costs $C\delta_N\cE/N$. For its complement,
\[
 \frac{\dd[M_{\rm n,high}]}{\dd s}
 \le\frac C N\left(\sum_{|i|>K}\phi_i|v_i|\right)^2
 \le Ce^{-cK/L}\cE.
\]
The inequality $\dd[A+B]\le2\dd[A]+2\dd[B]$, applied twice, absorbs all cross brackets and proves \eqref{different:eq:LF5-bracket}. The preceding drift bounds prove \eqref{different:eq:LF5-drift}.
\end{proof}

\subsection{Stopping and the failure probability}

In this subsection we stop before the energy becomes large and show that this stopping event has a very small probability. We fix a deterministic $\alpha$. We first obtain its tangent estimate using adapted stops, then integrate the exceptional parameters in the next section. Set $R_1=R+2$ and
\[
 \cE_*=N^{-2}t^{-2R_1}.
\]
Stop at the first energy crossing, at $\rho_\alpha\wedge\rho_{\rm n}$, or at $\tau$. Before this stop the accumulated drift divided by $\cE_*$ is at most
\begin{equation}
 C\frac{Ht}{L}+CH\nu_Nt^{2R_1},
 \label{different:eq:drift-ratio}
\end{equation}
and the stopped bracket divided by $\cE_*^2$ is at most
\begin{equation}
 \Lambda_N=C\left(\frac H{L^2}+H\nu_Nt^{2R_1}\right).
 \label{different:eq:Lambda}
\end{equation}
Both bounds follow by integrating Proposition~\ref{different:prop:energy} over the deterministic horizon $\tau=H/N$, while the stopped energy is at most $\cE_*$. They therefore apply to the stopped martingale without conditioning on a future event.

Under \eqref{different:eq:order},
\[
 \nu_N\le Ct^{7k}e^{-2t^{1/3}}
       +\exp\{2t-ct^{k-\ell}\}.
\]
Consequently \eqref{different:eq:drift-ratio} is $o(1)$ and
\[
 \Lambda_N=o(t^{-1}),\qquad
 \Lambda_N^{-1}/\log N\longrightarrow\infty.
\]
Since $\cE(0)=0$, an energy exit before $\rho_\alpha\wedge\rho_{\rm n}\wedge\tau$ requires the stopped continuous martingale to increase by at least a fixed fraction of $\cE_*$. The exponential martingale inequality gives
\begin{equation}
 \Pp(\text{energy exit before }\rho_\alpha\wedge\rho_{\rm n}\wedge\tau)
 \le2e^{-c/\Lambda_N}.
 \label{different:eq:exit-prob}
\end{equation}
Since $\Lambda_N^{-1}/\log N\to\infty$, this is smaller than every fixed power $N^{-D}$.

Equations \eqref{different:eq:weight-support} and \eqref{different:eq:exit-prob} give
\[
 \sup_{s\le\tau,\,|i|\le L}|w_i(s)|\le CN^{-1}t^{-R-2}.
\]
The long-interaction correction is at most
\[
 CN^{-1}t^{2+3b/2-\ell}=o(N^{-1}t^{-R-2}).
\]
Together with \eqref{different:eq:reg-prob} and \eqref{different:eq:noise}, we have proved for every fixed deterministic $\alpha$, uniformly in $\alpha$,
\begin{equation}
 \Pp\left(\sup_{s\le\tau,\,|i|\le L}|u_i(s,\alpha)|
              >CN^{-1}t^{-R-2}\right)
 \le C_DN^{-D}.
 \label{different:eq:fixed-alpha}
\end{equation}
The estimate is uniform in the deterministic parameter, but does not assert one event on which it holds for every $\alpha$. The failure bound can be made strong enough to absorb the loss when the global tangent moment is used to integrate exceptional parameters.

\subsection{Integration in the noise parameter}\label{different:sec:alpha}

In this subsection we integrate the noise derivative to compare the endpoint flows, controlling the contribution of exceptional interpolation parameters.
To compare the endpoint flows, we need the parameter integral of the tangent. The proof has two parts: identify that integral with the endpoint difference, and bound the contribution from parameters where \eqref{different:eq:fixed-alpha} fails. Joint measurability and a global moment bound make Fubini applicable, without an intersection of uncountably many good events.

\begin{lemma}\label{different:lem:increments}
For deterministic $\alpha,\gamma\in[0,1]$,
\begin{equation}
 \E\sup_{s\le\tau}\|x^\alpha(s)-x^\gamma(s)\|_2^2
 \le C|\alpha-\gamma|^2\tau.
 \label{different:eq:alpha-increments}
\end{equation}
\end{lemma}
\begin{proof}
Set $D_i=x_i^\alpha-x_i^\gamma$. The common initial configuration cancels, so $D(0)=0$ regardless of its moments. The exact difference equation is
\[
 \dd D_i=(\alpha-\gamma)\dd\xi_i
 +\sum_{j\ne i}c_{ij}^{\alpha,\gamma}(D_j-D_i)\dd s,
 \qquad
 c_{ij}^{\alpha,\gamma}
 =\frac1{2N(x_i^\alpha-x_j^\alpha)(x_i^\gamma-x_j^\gamma)}>0.
\]
Positivity uses common ordering. The interaction in $\dd\|D\|_2^2$ is nonpositive. Its scalar martingale bracket is at most $C|\alpha-\gamma|^2\|D\|_2^2\dd s$, and the additive trace is at most $C|\alpha-\gamma|^2\dd s$. BDG and Young prove \eqref{different:eq:alpha-increments} under localization; Fatou removes localization.
\end{proof}

\begin{lemma}[Tangent identification and integration]\label{different:lem:FTC}
The family can be chosen jointly measurably in $(\alpha,\omega)$. Its tangent satisfies
\begin{equation}
 \E\int_0^1\sup_{s\le\tau}\|u(s,\alpha)\|_2^2\dd\alpha\le C\tau,
 \label{different:eq:integrated-moment}
\end{equation}
and, almost surely, for every $i$ and every $s\in[0,\tau]$,
\begin{equation}
 x_i^1(s)-x_i^0(s)=\int_0^1u_i(s,\alpha)\dd\alpha.
 \label{different:eq:FTC}
\end{equation}
\end{lemma}
\begin{proof}
\emph{Localization and measurability.}
Replace the singular drift by smooth globally Lipschitz drifts agreeing with it when all relevant gaps exceed $1/m$ and all coordinates have absolute value below $m$. After subtracting the continuous forcing $\chi^\alpha$, each localized equation is an ordinary integral equation depending affinely on $\alpha$ through that forcing. Pathwise differentiation gives a continuously differentiable parameter map and the localized tangent equation.

For each fixed $\alpha$, (WP) implies that the original path stays in a compact subset of the open chamber on $[0,\tau]$. One sufficiently large localization therefore agrees along the whole path. Continuous dependence of the localized equation on $\alpha$ supplies a random neighborhood in which the same original equation is valid. Thus the local pathwise derivative is precisely \eqref{different:eq:tangent}. These localized constructions provide a canonical jointly measurable family. By (WP) and Fubini, for almost every sample the noncolliding parameters have full Lebesgue measure.

\emph{Identification at a fixed time.}
Fix a rational time $s\in[0,\tau]$ and put
\[
 Y^\alpha(s)=x^\alpha(s)-x^0(s).
\]
Lemma~\ref{different:lem:increments} puts $\alpha\mapsto Y^\alpha(s)$ in
$H^1((0,1);L^2(\Omega;\mathbb R^{2N}))$. Subtracting $x^0$ avoids any unnecessary moment assumption on the common initial configuration. Localized difference quotients converge pointwise to $u(s,\alpha)$. Their global $L^2(\Omega\times(0,1))$ bounds from \eqref{different:eq:alpha-increments} give uniform integrability in $L^1$; every weakly convergent $L^2$ subsequence therefore has the same limit $u$. The moment estimate in Lemma~\ref{different:lem:global} ensures that this limit belongs to $L^2$. This identifies the Hilbert weak derivative.

By Fubini, the one-dimensional Sobolev fundamental theorem applies for almost every sample at this fixed rational time. Its endpoint traces agree with the actual endpoint processes by the $L^2$ continuity in \eqref{different:eq:alpha-increments}. Intersecting over countably many rational times gives \eqref{different:eq:FTC} on that dense set.

\emph{Extension to every physical time.}
Integrating Lemma~\ref{different:lem:global} in $\alpha$ yields \eqref{different:eq:integrated-moment}. Hence $\int_0^1\sup_{s\le\tau}|u_i(s,\alpha)|\dd\alpha$ is finite almost surely. The tangent is continuous in time for almost every $\alpha$, so dominated convergence makes its parameter integral continuous in time. The two endpoint processes are continuous as well. This extends the rational-time identity to all $s\in[0,\tau]$ simultaneously.
\end{proof}

\begin{proof}[Proof of Theorem~\ref{different:thm:main}]
Using the jointly measurable tangent version from Lemma~\ref{different:lem:FTC}, define
\[
 B_\alpha=\left\{\sup_{s\le\tau,\,|i|\le L}|u_i(s,\alpha)|
                    >CN^{-1}t^{-R-2}\right\},
\]
with the fixed threshold constant from \eqref{different:eq:fixed-alpha}.
This event is jointly measurable in $(\alpha,\omega)$, by continuity
in physical time and a countable time supremum. Its probability is
bounded by \eqref{different:eq:fixed-alpha}, including the failure
probabilities of the unweighted regularity and noise stopping conditions. We therefore need no joint
measurability assumption on a separately chosen family of stopping
times $\rho_\alpha$. Choose the probability exponent $D_0$ in that
bound. By Cauchy--Schwarz in probability and Lemma~\ref{different:lem:global},
\begin{align*}
 &\E\int_0^1\mathbf1_{B_\alpha}
       \sup_{s\le\tau}\|u(s,\alpha)\|_\infty\dd\alpha\\
 &\hspace{2em}\le
 \int_0^1\Pp(B_\alpha)^{1/2}
 \left(\E\sup_{s\le\tau}\|u(s,\alpha)\|_2^2\right)^{1/2}\dd\alpha
 \le C\sqrt\tau\,N^{-D_0/2}.
\end{align*}
Markov's inequality at a fixed small multiple of $N^{-1}t^{-R}$ gives failure at most
\[
 C N^{(1-D_0)/2}t^{b/2+R}.
\]
Taking $D_0>2D+3$ makes this $O(N^{-D})$. On the remaining parameter set, the integral of the bound in \eqref{different:eq:fixed-alpha} is $CN^{-1}t^{-R-2}$, a vanishing fraction of $N^{-1}t^{-R}$. Equation \eqref{different:eq:FTC} now proves \eqref{different:eq:conclusion}.
\end{proof}

\subsection{Parameter order and the local-factor comparison}

In this subsection we shall choose the comparison scales in a fixed order and pass the resulting particle bounds to the local-factor construction.

We choose the parameters in the following order: first the counting scale $a$ and the available weak-law horizon, then the Brownian horizon exponent $b$ and the particle accuracy $R$, and finally $\ell$ and $k$ by \eqref{different:eq:order}. The companion counting proof is established on a horizon such as $t^{-2M}$ with fixed $M$, which contains $t^b/N$ for every ultimately fixed $b$. Its only covariance inputs are diagonal bracket bounds and their positive-semidefinite consequence. Thus enlarging the later noise window does not alter $a$.

In the energy proof, the probability exponent affects constants and sufficiently large $N$ thresholds. In the final parameter integration, one chooses the fixed-alpha failure strength $D_0$ after the desired final strength $D$. This order leaves the previously chosen window exponents unchanged.

For Section~\ref{chap:lf}, take $x^1$ to be the actual flow and $x^0$ the flow driven by Section~\ref{chap:noise}'s reference Brownian family from the same initial singular values. Theorem~\ref{different:thm:main} gives precisely the uniform particle displacement in \eqref{lf:eq:lfl-different}, with $\varepsilon=t^{-R}$ on $|i|\le L$ throughout $[0,H/N]$. This estimate controls the endpoint kernels and the particle-displacement part of the stochastic integral. The separate driver-difference covariance estimate controls the change of integrator.

Section~\ref{chap:relax} compares $x^0$ with the independently initialized Gaussian reference driven by the same reference noises. Independence of that reference from the original initial matrix comes from Section~\ref{chap:noise}'s construction. The local-factor proof combines these two comparisons with exact reconstruction, then applies local caps and weighted exceptional-event estimates.

\endgroup

\begingroup
\newcommand{\E}{\mathbb E}
\newcommand{\PP}{\mathbb P}
\newcommand{\one}{\mathbf1}
\newcommand{\dd}{\,\mathrm d}
\section{Local factors and exact normalizing expectations}\label{chap:lf}

The purpose of this section is to approximate determinant weights using the initial matrix and independent local factors, with errors relative to exact normalizing expectations.

\paragraph{Local factors and normalization.}
The construction follows the independent local/global decomposition of
\cite[Theorem~3.8]{CL} and the comparison of the resulting normalizing
expectations in \cite[Lemma~4.4 and Theorem~4.6]{CL}. In that decomposition, spatial points are
separated by at least $N^{-1/2+\delta}$ and the elapsed time is
$N^{-1+\omega}$, for fixed positive exponents. Here the logarithmic
relaxation estimates approximate the field at time $t^b/N$ for the
prescribed pair in \eqref{lf:eq:lfl-geometry}, whose separation can be
as small as a constant times $N^{-1/2}e^{t^{1/3}/2}$.

For the complex Gaussian-aged model below, we first replace the actual
drivers while keeping the initial spectrum, then use the same reference
drivers to compare with an independently initialized Gaussian spectrum.
At each point, the local variables use the full reference spectrum and
its driving Brownian motions. They are independent of the original
initial matrix. Local caps and weighted
exceptional-event bounds give
\eqref{lf:eq:lfl-lf1}--\eqref{lf:eq:lfl-lf3}. We use exact
one-point expectations for normalization, and the upper comparison for
the normalized unregularized product at the two points
holds for every $0\le F\le1$ measurable with respect to the initial
matrix.

\subsection{Model, normalization, and estimates used below}

{We shall fix the model and the three coupled spectra, keeping the reference independence and exact normalizing expectations explicit.}
We first fix a compact disk $\mathcal D\Subset\{|v|<1\}$ and prescribed
deterministic points $z,w\in\mathcal D$ satisfying
\begin{equation}
 c_0N^{-1}e^d\le |z-w|^2,\qquad
 |z-w|\le C_0N^{-1/4},\qquad
 t=\tfrac12\log N,\quad d=t^{1/3}.
 \label{lf:eq:lfl-geometry}
\end{equation}
All constants are uniform in this geometry. The reference construction
is made for each prescribed deterministic pair of points.

Let $\widehat X=N^{-1/2}(\widehat\xi_{ab})$ have complex iid entries
with $\E\widehat\xi=\E\widehat\xi^2=0$,
$\E|\widehat\xi|^2=1$, and common bounds on all fixed moments.
Use independent normalized complex Ginibre $G^{\rm age}$ and set
\[
 X_0=\sqrt{1-\sigma}\,\widehat X+\sqrt\sigma\,G^{\rm age},
 \quad \sigma\asymp N^{-1/200},\qquad
 X_s=X_0+N^{-1/2}B_s,\quad 0\le s\le\tau.
\]
The age $\sigma$ is deterministic in
$[c_{\rm age}N^{-1/200},C_{\rm age}N^{-1/200}]$, for fixed positive
constants. All conclusions are uniform over this interval.
The entries of $B$ are standard complex Brownian motions with
$\E\dd B_{ab}\dd\overline B_{cd}
 =\delta_{ac}\delta_{bd}\dd s$. Put
\[
 H=N\tau=t^b,\quad S=N\eta_*=t^c,
 \quad \eta_{\rm m}=r_NN^{-1}t^{-q},\quad
 B_v=\sqrt{1-|v|^2},\quad \chi=2\sqrt2.
\]
Here $r_N$ is deterministic in a fixed compact subinterval of
$(0,\infty)$, and every conclusion is uniform in $r_N$.
In particular, the exact physical height
$\sqrt{1+\tau}\,N^{-1}t^{-q}$ is covered. The same microscopic
height is used in each actual and reference endpoint comparison.
For the exact variance-$(1+s)$ Hermitization density $\rho_s^v$,
define
\begin{equation}
 \Phi_s(v,\eta)
 =\frac12\log\det\big((X_s-v)(X_s-v)^*+\eta^2I\big)
 -\frac N2\int_{\mathbb R}\log(x^2+\eta^2)\rho_s^v(x)\dd x.
 \label{lf:eq:lfl-field}
\end{equation}
The value at zero may be $-\infty$; its positive exponential is
then zero. All normalizers below are exact expectations.

Write $\mathcal F_0^{\rm mat}$ for the original initial-matrix
sigma field. On an enlarged space use fresh mutually independent
Ginibre matrices $G_z^{\rm ref},G_w^{\rm ref}$. The enlarged initial
sigma field is
\[
 \mathcal G_0=\mathcal F_0^{\rm mat}
       \vee\sigma(G_z^{\rm ref},G_w^{\rm ref}).
\]
At each point, we use three signed singular-value processes. The intermediate
process changes only the drivers; the final reference also changes the
initial spectrum:
\begin{align*}
 \lambda^v&:\ \hbox{actual spectrum and actual drivers }b^v,\\
 x^v&:\ \hbox{same initial spectrum and reference drivers }\beta^v,\\
 \mu^v&:\ \hbox{initial singular values of }G_v^{\rm ref}/B_v
                 \hbox{ and drivers }\beta^v.
\end{align*}
Their signed index set is $\mathcal I_N=\{\pm1,\ldots,\pm N\}$,
with no zero index, and their additive SDE is
\begin{equation}
 dy_i=\frac{d\gamma_i}{\sqrt{2N}}
       +\frac1{2N}\sum_{j\in\mathcal I_N\setminus\{i\}}
                         \frac{\dd s}{y_i-y_j},
 \qquad y_{-i}=-y_i,\quad\gamma_{-i}=-\gamma_i.
 \label{lf:eq:lfl-sde}
\end{equation}
The sum includes $j=-i$.

\paragraph{Noise and relaxation inputs.}
The preceding sections supply the following estimates and independence
statements for these three processes.

\begin{enumerate}
\item \emph{Positive overlap and independent reference noises}
{\hyperref[chap:noise]{Section~\ref*{chap:noise}}} constructs the full positive-index
Brownian vector $\beta$, independent of $\mathcal G_0$, with independent
components across both points. Each $\mu^v$ is the strong measurable
solution determined by its own initial spectrum and its own Brownian
family. Thus $\mu^z,\mu^w$ are mutually independent and their joint
family is independent of $\mathcal F_0^{\rm mat}$, although it is not
independent of $\mathcal G_0$.
For $K=t^k$, before an adapted exceptional stop,
\begin{equation}
 \left\|\frac{d\langle b-\beta\rangle_s}{ds}\right\|_{\rm op}
       \le\delta_N,
 \qquad \delta_N\le Ct^{6k}e^{-2d}
 \label{lf:eq:lfl-noise}
\end{equation}
on the positive low-index vector of length $2K$.

\item \emph{Bounded-bracket off-axis weak law and counting}
{\hyperref[chap:relax]{Section~\ref*{chap:relax}}} fixes a logarithmic counting
exponent $a$ before $b,\ell,k$. Its fixed-interpolation-parameter
counting bounds hold with overwhelming probability, even
when the interpolated drivers have diagonal bracket less than one.

\item \emph{Comparison for different driving noises}
{\hyperref[chap:different]{Section~\ref*{chap:different}}} gives, for every fixed $R,D>0$,
\begin{equation}
 \sup_{s\le\tau}\max_{|i|\le L}
       |\lambda_i^v(s)-x_i^v(s)|\le\frac\varepsilon N,
 \qquad \varepsilon=t^{-R},\quad L=t^\ell,
 \label{lf:eq:lfl-different}
\end{equation}
outside an event of probability at most $C_DN^{-D}$, with the parameter order stated
in Section~\ref{lf:sec:lfl-parameters}.

\item \emph{Logarithmic same-noise relaxation} is proved in
{\hyperref[chap:relax]{Section~\ref*{chap:relax}}}. Its theorem, labeled
{\hyperref[relax:thm:lfsn-relaxation]{\ref*{relax:thm:lfsn-relaxation}}}, supplies on the actual coupling
$x^v,\mu^v$,
\begin{equation}
 |x_i^v(s)-\mu_i^v(s)|
 \le C_D\varphi^{15}\frac iN
       \left[(Ns)^{-1/2}+\max\{i/N,s\}\right],
 \quad s_0\le s\le\tau,\quad 1\le i\le c_*N,
 \label{lf:eq:lfl-same}
\end{equation}
where $\varphi=(\log N)^M$ and $s_0=\varphi^{100}/N$.
Here $M>64L_{\rm sn}$ is fixed after that section's initial
weak-law exponent $L_{\rm sn}$, so this $\varphi$ eventually
dominates its $C^{16}(\log N)^{64L_{\rm sn}}$; we use the resulting
weaker prefactor, later starting time, and narrower horizon.
It holds with overwhelming probability. Section~\ref{chap:relax}
proves this logarithmic range using a random smooth reference and
buffered initial-data quantile interpolation. Its same-noise tangent
calculation uses unit diagonal brackets. The different-noise counting
proof retains the cubic correction needed for diagonal brackets at most
one.
\end{enumerate}

\paragraph{Uniformity in the deterministic age and height factors.}
Theorem~\ref{i2:thm:I2} permits comparable Gaussian ages. The overlap proof
uses the OU age $-\log(1-\sigma)$ and initial characteristic heights
comparable to $\sigma$. The aged weak-law proof starts from height
comparable to $\sigma$ and retains the initial error
$CN^\xi/[N(\eta+\sigma)]$. Replacing the age by any allowed factor
therefore changes only fixed constants; its polynomial margins and
the later counting and relaxation exponents stay the same.
For the microscopic height, the endpoint proofs below use only
$\eta_{\rm m}^{-1}\le CNt^q$ and $N\eta_{\rm m}\le Ct^{-q}$.
The scalar comparison is uniform at such heights, and the CL cutoff
lemma permits both cutoffs to change by bounded positive factors.
All other reconstruction bounds involve $\eta_*$ rather than
$\eta_{\rm m}$. Thus the error bounds and parameter order
apply uniformly to the stated $r_N$; no approximate substitution for
the physical height or its deterministic centering is being made.

We also use the following marginal estimates under the matrix
probability law, before any exponential reweighting: the
imaginary-axis bound $|m_N-m|\le Ct/(N\eta)$ on the relevant
characteristics, bounded imaginary transforms, and lower rigidity
$\lambda_i^v,x_i^v,\mu_i^v\ge ci/N$ for $i\ge I_0=t^{m_0}$ in
a central index block. Outside that block the spectra are bounded
away from zero. These statements are uniform in the required times,
with overwhelming probability. The marginal law of
$x^v$, conditional on its initial singular values, equals that of
$\lambda^v$: its positive drivers are Brownian and independent of
the initial matrix, and the signed SDE is strongly unique. This
does not assert joint independence of $x^z,x^w$.

\paragraph{Scalar inputs.}
The scalar positive-moment and regularization estimates are separate
from the coupling. We use the following uniform estimates at
the finitely many heights and charges used below:
\begin{align}
 m_v^{\rm m}:=\E e^{\chi\Phi_\tau(v,\eta_{\rm m})}
       &\asymp N,
 &m_v^0:=\E e^{\chi\Phi_\tau(v,0)}&\asymp N,
                                                        \label{lf:eq:lfl-scalar1}\\
 \frac{m_v^{\rm m}}{m_v^0}
       &=1+O(t^{-q}+N^{-\alpha}),
                                                        \label{lf:eq:lfl-scalar2}\\
 \E e^{r\Phi_s(v,\eta)}&\le C_rN^{r(r+2)/8},
 &r&\in\{2\chi,4\chi\}.
                                                        \label{lf:eq:lfl-moments}
\end{align}
Here $\alpha>0$ is fixed. The last inequality also holds for the
initial height $B_v\tau$. These are proved, independently of the local-factor construction, in
{\hyperref[chap:scalar]{Section~\ref*{chap:scalar}}}:
Theorem {\hyperref[scalar:thm:sr-sharp]{\ref*{scalar:thm:sr-sharp}}}, Theorem
{\hyperref[scalar:thm:sr-regularization]{\ref*{scalar:thm:sr-regularization}}}, and Proposition
{\hyperref[scalar:prop:sr-moments]{\ref*{scalar:prop:sr-moments}}}, with exact variance scaling in
Lemma {\hyperref[scalar:lem:sr-variance]{\ref*{scalar:lem:sr-variance}}}. They are not inferred from a
comparison at two spatial points. The local
one-sided cutoff estimate of CL Lemma 2.12 is used explicitly below.
All variance and deterministic-centering conversions are exact.

\subsection{The local factors and the theorem}

In this subsection we shall define the local factors and state the field approximation and expectation bounds they must satisfy.

Each reference process has its own exact deterministic density
$\rho_s^{\mu,v}$, corresponding to its initial scale $1/B_v$ and
additive variance $s$. Let $\Phi_s^{\mu,v}(\eta)$ be the field
defined by the signed reference spectrum and this density, with
the same normalization as \eqref{lf:eq:lfl-field}. The first local
variable is the reference field increment from the mesoscopic height
$\eta_*$ down to $\eta_{\rm m}$. The second is its stochastic integral
along the linearized characteristic after the early interval $[0,s_0]$.
Both use the full reference spectrum and only the initial data and
Brownian family at the same point. Put
\[
 \nu_{v,s}=\eta_*+B_v(\tau-s),
\]
\begin{align}
 L_{1,v}&=\Phi_\tau^{\mu,v}(\eta_{\rm m})
                          -\Phi_\tau^{\mu,v}(\eta_*),
                                                        \label{lf:eq:lfl-L1}\\
 L_{2,v}&=\frac12\Re\frac1{\sqrt{2N}}
       \sum_{i\in\mathcal I_N}\int_{s_0}^\tau
                      \frac{d\beta_i^v(s)}{\mu_i^v(s)-i\nu_{v,s}}.
                                                        \label{lf:eq:lfl-L2}
\end{align}
For fixed $A_{\rm loc}>0$, cap these two local variables and separate
the reference factor from the initial-matrix factor:
\begin{align}
 \mathcal K_v&=\{L_{1,v}\le(\log N)^{3/4}\}
                      \cap\{L_{2,v}\le A_{\rm loc}\log N\},
                                                        \label{lf:eq:lfl-cap}\\
 U_v&=e^{\chi(L_{1,v}+L_{2,v})}\one_{\mathcal K_v},
 &V_v&=e^{\chi\Phi_0(v,B_v\tau)}.
                                                        \label{lf:eq:lfl-UV}
\end{align}
The event $\mathcal K_v$ depends only on the reference inputs at $v$.
The reconstruction event will be used to estimate errors, separately
from the definition of $U_v$.

The theorem first approximates the field with an arbitrarily small
inverse-logarithmic error and compares the exact normalizing
expectations. Its last conclusion bounds the expectation of the
normalized unregularized product at the two points, tested against a
permitted function measurable with respect to the full original initial
information, by an expectation over that initial information.

\begin{theorem}[Fixed-time local factors]
\label{lf:thm:lfl-main}
Assume the noise, relaxation, and scalar inputs stated above. For every fixed
$P,D,A_0>0$, the logarithmic parameters may be chosen so that the
following hold uniformly in \eqref{lf:eq:lfl-geometry}.
There is a deterministic $\rho_N=o(t^{-P})$ and an event of
probability at least $1-C_DN^{-D}$ on which, for $v=z,w$,
\begin{equation}
 \left|\Phi_\tau(v,\eta_{\rm m})
       -\Phi_0(v,B_v\tau)-L_{1,v}-L_{2,v}\right|
       \le\rho_N.
 \tag{LF1}\label{lf:eq:lfl-lf1}
\end{equation}
Moreover
\begin{equation}
 \frac{(\E U_v)(\E V_v)}{m_v^{\rm m}}=1+o(t^{-P}),
 \qquad
 \frac{(\E U_v)(\E V_v)}{m_v^0}=1+o(t^{-P}).
 \tag{LF2}\label{lf:eq:lfl-lf2}
\end{equation}
For every $\mathcal F_0^{\rm mat}$-measurable $0\le F\le1$,
\begin{equation}
 \frac{\E[e^{\chi\Phi_\tau(z,0)+\chi\Phi_\tau(w,0)}F]}
            {m_z^0m_w^0}
 \le(1+o(t^{-P}))
       \E\left[\frac{V_zV_w}{(\E V_z)(\E V_w)}F\right]
       +O(N^{-A_0}).
 \tag{LF3}\label{lf:eq:lfl-lf3}
\end{equation}
The two reference processes, and consequently $U_z,U_w$, are
mutually independent and jointly independent of
$\mathcal F_0^{\rm mat}$.
\end{theorem}

\subsection{Exact dynamics and the field approximation at one point}

Now we follow the logarithmic field along its characteristic to compare two spectra driven by the same noise.
We first derive the exact field dynamics. The identity below holds for
both the actual flow $\lambda$ and the flow
$x$ with the same initial spectrum and the reference drivers, using
their respective drivers and the same deterministic density.
We then compare $x$ with $\mu$ using their common drivers, before
adding the $\lambda$-to-$x$ error.

We suppress the spatial label. Let $\zeta_s=i\eta_s$ be the deterministic
characteristic ending at $i\eta_*$ at time $\tau$, with
$\zeta_s'=-m_s(\zeta_s)$. For the analytic centered logarithm
\[
 \mathcal A_s=\frac12\sum_{i\in\mathcal I_N}
                \log(\lambda_i(s)-\zeta_s)
       -N\int\log(x-\zeta_s)\rho_s(x)\dd x,
\]
one has $\Re\mathcal A_s=\Phi_s(v,\eta_s)$. Branch constants do
not affect this equality or its differential. Direct It\^o
calculation gives
\begin{equation}
 d\mathcal A_s
   =\frac1{2\sqrt{2N}}\sum_i
                \frac{db_i(s)}{\lambda_i(s)-\zeta_s}
       -\frac N2(m_{N,s}(\zeta_s)-m_s(\zeta_s))^2\dd s.
 \label{lf:eq:lfl-log-ito}
\end{equation}
For completeness, if $g_i=(\lambda_i-\zeta)^{-1}$, symmetrization
of the interaction in the half-logarithm gives
$-(8N)^{-1}[(\sum_i g_i)^2-\sum_i g_i^2]$. The It\^o correction
is $-(8N)^{-1}\sum_i g_i^2$, and the characteristic term is
$(m/2)\sum_i g_i$. The deterministic logarithmic integral has
characteristic derivative $m^2/2$. Their combination is precisely
$-N(m_N-m)^2/2$.

On the imaginary axis this real drift is nonnegative, because
$m_N-m$ is purely imaginary. We use its absolute value. The
height derivative, directly from \eqref{lf:eq:lfl-field}, is
\begin{equation}
 \partial_\eta\Phi_s(v,\eta)=N\Im(m_{N,s}(i\eta)-m_s(i\eta)).
 \label{lf:eq:lfl-height-derivative}
\end{equation}
The preceding identities are exact. The remaining comparisons use
deterministic characteristic estimates and stopped local laws. Uniform
deterministic bulk smoothness gives
\begin{equation}
 \eta_s\asymp\nu_s,
 \quad |\eta_s-\nu_s|\le C\eta_s(\tau+\eta_*),
 \quad |\eta_0-B_v\tau|\le C(\tau^2+\eta_*).
 \label{lf:eq:lfl-characteristic-errors}
\end{equation}
The local law and $-\dot\eta_s\asymp1$ therefore imply
\begin{equation}
 \frac N2\int_0^\tau|m_{N,s}-m_s|^2\dd s
       \le C\frac{t^2}{S}.
 \label{lf:eq:lfl-drift-budget}
\end{equation}
Integrating \eqref{lf:eq:lfl-height-derivative} over the initial
height change costs $C(t\tau+tS/H)$.

Pairing signed indices reduces the real martingale in
\eqref{lf:eq:lfl-log-ito} to
\[
 \frac1{\sqrt{2N}}\sum_{i=1}^N\int
              f_\eta(\lambda_i)\,db_i,
 \qquad f_\eta(x)=\frac{x}{x^2+\eta^2},
 \qquad |f_\eta'(x)|\le\frac1{x^2+\eta^2}.
\]
On $[0,s_0]$ the characteristic height is comparable to $\tau$;
bounded imaginary transforms give martingale bracket
$C\varphi^{100}/H$. The height replacement
$\eta_s\to\nu_s$ has bracket at most
$C(\tau+\eta_*)^2\log(1+H/S)$, by
\eqref{lf:eq:lfl-characteristic-errors}. This is polynomially small
in $N$ for fixed logarithmic exponents.

For the same-noise comparison of $x$ with $\mu$, use $u=Ns$ and
$v=N\nu_s=S+B_v(H-u)$. Lower rigidity above $I_0$ and
\eqref{lf:eq:lfl-same} bound its bracket by
\begin{align}
 C\varphi^{30}\int_{\varphi^{100}}^H
   \sum_{i=I_0}^{c_*N}\frac{i^2}{(i^2+v^2)^2}
                 \left[\frac1u+\frac{i^2+u^2}{N^2}\right]\dd u
 &\le C\varphi^{30}
   \left[\frac{\log H}{H}+\frac HN
           +\frac{H^2}{N^2}\log(1+H/S)\right].
 \label{lf:eq:lfl-same-bracket}
\end{align}
Here $\sum_i i^2/(i^2+v^2)^2\le C/v$, and splitting at $H/2$
gives $\int du/(uv)=O(\log H/H)$. For the shape term use
$\sum_{i\le c_*N}i^4/(i^2+v^2)^2\le CN$.
For $i<I_0$ use only the regularizer in the denominator. Their
leading contribution is
$C\varphi^{30}I_0^3\log H/(HS^3)$, absorbed by the first term
when $S\ge I_0$. Their remaining terms carry $N^{-2}$.
Outside the central block, bounding both martingales separately
costs bracket $O(\tau)$.

Each integral is stopped at the first pointwise adapted failure of
the estimates just used. At $t$ times the square root of a
deterministic stopped bracket, the exponential martingale inequality
has tail $\exp(-ct^2)$, smaller than every fixed $N^{-D}$.
The probability of the stop is then added. No future good event is
conditioned upon. Thus the early and same-noise martingale errors
are bounded by
\begin{equation}
 C\left[t\sqrt{\frac{\varphi^{100}}H}
       +t\varphi^{15}\sqrt{\frac{\log H}{H}}
       +N^{-c_1}t^{C_1}\right].
 \label{lf:eq:lfl-one-mart}
\end{equation}

\subsubsection{Endpoint logarithms and exact density centerings}

Next we compare endpoint logarithms with their exact density centerings, so the approximation preserves the field's normalization.
Define
\[
 h_*(x)=\frac12\log\frac{x^2+\eta_*^2}{x^2+\eta_{\rm m}^2},
 \qquad
 |h_*'(x)|\le C\min\{\eta_{\rm m}^{-1},|x|^{-1},
                                      \eta_*^2|x|^{-3}\}.
\]
The actual field difference between $\eta_{\rm m}$ and $\eta_*$ is
$-\sum_{i>0}h_*(\lambda_i)+N\int h_*\rho$; the same formula
with $x_i$ holds for the intermediate flow. First compare that flow
with $\mu$, using each process's exact deterministic density. Applying
\eqref{lf:eq:lfl-same} at time $\tau$, summing below and above
$I_0$, gives the random endpoint error
\begin{equation}
 C\varphi^{15}\left[
  \frac{I_0^2t^q+S}{\sqrt H}
  +\frac{I_0^3t^q+S^2\log N}{N}
  +\tau(I_0^2t^q+S)\right]+CN\eta_*^2.
 \label{lf:eq:lfl-one-endpoint}
\end{equation}
For $i<I_0$, the first term uses
$\eta_{\rm m}^{-1}\le CNt^q$ and the sum of $i$; it assumes no
lower bound on the smallest singular value. For larger indices
the derivative bound uses $\min\{N/i,S^2N/i^3\}$.
Indices outside the central block are bounded by the kernels
separately, giving $CN\eta_*^2$.

The actual and reference initial densities agree at zero, with
value $B_v/\pi$. They are even and smooth near zero, so at time
$\tau$ their difference is $O(\tau+x^2)$ there. The bounds
$\int h_*\dd x=O(\eta_*)$ and
$\int_{|x|\le c}x^2h_*\dd x=O(\eta_*^2)$ give a deterministic
centering error $O((HS+S^2)/N)$. The exterior costs
$O(S^2/N)$ because its distance from zero is fixed. In particular
there is no order-one density-centering discrepancy.

Combining the characteristic, martingale, endpoint and density-centering
estimates gives the error in approximating the field of the flow with
the same initial spectrum and the reference drivers by the field of its Gaussian reference:
\begin{equation}
 \rho_{\rm one}\le C\left[
  \frac{t^2}{S}+\frac{tS}{H}
  +t\sqrt{\frac{\varphi^{100}}H}
  +t\varphi^{15}\sqrt{\frac{\log H}{H}}
  +\frac{\varphi^{15}(I_0^2t^q+S)}{\sqrt H}
  +N^{-c_1}t^{C_1}\right].
 \label{lf:eq:lfl-one-budget}
\end{equation}

\subsection{Adding the different-noise comparison}

The purpose of this subsection is to add the errors caused by changing both the particle paths and the driving noise.

It remains to compare $\lambda$ with $x$. Their initial spectra and
deterministic centerings agree. Particle displacement changes the
endpoint kernels and the stochastic integrands; the separate covariance
estimate controls the change of integrator. Choose endpoint and
martingale windows $I_0\le J=t^j\le K_1=t^{k_1}\le L\le K$.
Below $I_0$, \eqref{lf:eq:lfl-different} and the regularizer give
cost $C\varepsilon I_0t^q$. Between $I_0$ and $J$, rigidity and
the derivative envelope give
$C\varepsilon\sum_i\min\{i^{-1},S^2i^{-3}\}$.
Above $J$, bounding the two kernels separately gives
$CS^2\sum_{i>J}i^{-2}$. Thus the extra endpoint error is
\begin{equation}
 C\varepsilon[I_0t^q+\log(2+S)]+C\frac{S^2}{J}.
 \label{lf:eq:lfl-diff-endpoint}
\end{equation}
Only $J\le L$ is required; no condition $J\ll\sqrt H$ is used.
There is no additional deterministic-centering term in this comparison.

For martingales on $i\le K_1$, the particle displacement has
bracket $C\varepsilon^2/S^2$; the indices below the rigidity cutoff contribute
$C\varepsilon^2I_0/S^3$, which is absorbed by $S\ge I_0$.
The driver-difference bracket \eqref{lf:eq:lfl-noise} gives
$C\delta_N\log(1+H/S)$. The two high-index martingales are
bounded separately, giving $CH/K_1$. Their influence on the
low-index Dyson drift is already contained in
\eqref{lf:eq:lfl-different}; the tail bracket alone would not prove
that comparison. The cross quadratic variations are bounded by a fixed
multiple of the sum of these three bounds. The additional stopped
martingale error is therefore
\begin{equation}
 Ct\left[\frac\varepsilon S
       +\sqrt{\delta_N\log(1+H/S)}
       +\sqrt{\frac H{K_1}}\right].
 \label{lf:eq:lfl-diff-mart}
\end{equation}
All integrands are adapted, and the stopped martingale estimate applies
to each of these three contributions.

We have proved \eqref{lf:eq:lfl-lf1} with
\begin{equation}
 \rho_N\le C\left[
  \rho_{\rm one}
  +\varepsilon(I_0t^q+\log(2+S)) +\frac{S^2}{J}
  +\frac{t\varepsilon}{S}
  +t\sqrt{\delta_N\log(1+H/S)}
  +t\sqrt{\frac H{K_1}}\right].
 \label{lf:eq:lfl-total-budget}
\end{equation}

\subsection{Order of parameter choices}
\label{lf:sec:lfl-parameters}

The purpose of this subsection is to choose the parameters in order so every error in approximating the field meets the prescribed logarithmic accuracy.
The common order of parameter choices is given in
Section~\ref{main:parameter-schedule}; the estimates below justify
these choices.

Fix the charge range, $P,D,A_0,A_{\rm loc}$ and a failure order $D_0$
large enough to control the explicit moment contributions in
Section~\ref{lf:sec:lfl-weights}. These moment requirements depend on the charges
and the cap envelope, so they can be fixed before the reconstruction
parameters. The preceding noise and relaxation estimates fix
$M,m_0,a$ and their local-law floors before the terminal time and noise
windows. The weak-law estimate, for example, permits $a>6M+3$ and has horizon
$(\log N)^{-2M}$, which contains every later fixed $t^b/N$.
Let $C_*$ dominate the fixed lower-height exponents required by the one-point estimates used above. Choose
\begin{align*}
 q&>P+2,\qquad c>\max\{m_0,P+3,C_*\},\\
 b&>\max\left\{
  100M+2P+4,\ c+P+2,
  2[15M+\max\{2m_0+q,c\}+P+2]\right\},\\
 j&>2c+P+2,\qquad
 k_1>\max\{j,b+2P+4\},\qquad R>P+q+m_0+3,\\
 \ell&>\max\{a+1,k_1+1,3b/2+R+5\},\\
 k&>\max\{\ell+2,b_{\rm ins}+2,m_0+2\}.
\end{align*}
Here $b_{\rm ins}$ is the fixed insertion/overlap height floor from
the independent-noise comparison. All exponents are fixed, and
$e^{-2t^{1/3}}$ dominates every fixed power of $t$.
These inequalities make every term in \eqref{lf:eq:lfl-total-budget}
$o(t^{-P})$: $c$ controls the characteristic drift, $b$ the same-noise
and early-time errors, $j,k_1$ the endpoint and martingale tails, and
$R$ the particle displacement. The final windows $\ell,k$ then meet
the different-noise and covariance requirements. Thus the order is
\[
 \text{accuracy and moment bounds}
 \ \longrightarrow\ (M,m_0,a)
 \ \longrightarrow\ (q,c)
 \ \longrightarrow\ b
 \ \longrightarrow\ (j,k_1,R)
 \ \longrightarrow\ \ell\ \longrightarrow\ k.
\]
Increasing an arbitrary polynomial failure exponent changes constants
and the threshold in $N$, not these dependencies. For large $N$,
$\tau=o(\sigma)$, $s_0<\tau$, and all required particle windows
are smaller than every fixed power of $N$.

\subsection{Local caps and weighted bounds on exceptional events}
\label{lf:sec:lfl-weights}

Next we cap the local factors and bound exceptional weighted contributions before comparing their exact expectations.
The cap on $L_1$ follows from CL Lemma 2.12, (2.31), applied to
the terminal reference spectrum. In law this is a scaled complex
Ginibre matrix with variance $B_v^{-2}+\tau$, in a fixed compact
positive interval. Dividing by its scale moves both cutoffs only
by bounded factors and leaves the centered field unchanged.
Choose the fixed cutoff exponent in that lemma larger than $c$ and its
parameter $\delta\in(0,1/4)$, so the logarithmic error exponent
$1/2+\delta$ is strictly below $3/4$.
The resulting cap failure is smaller than any prescribed $N^{-D_0}$.

After signed pairing, $L_2$ has bracket at most
$C\log(1+H/S)=O(\log t)$ on an adapted marginal local-law stop.
Its cap failure is consequently bounded by
\[
 C\exp\{-c(\log N)^2/\log\log N\}+C_{D_0}N^{-D_0}.
\]
The cap itself uses the full reference path. The adapted stop is used
only to prove its probability bound and does not enter the local factor.
Because the cap uses only the reference spectrum and its driving
Brownian motions at the same point, the local factors remain
independent, and
\begin{equation}
 0\le U_v\le N^{A_*},\qquad A_*:=\chi(A_{\rm loc}+1).
 \label{lf:eq:lfl-cap-envelope}
\end{equation}
No lower cap is needed for a positive exponential weight.

We next compare the relevant expectations by controlling the weighted
exceptional set, before using any estimate of the reference expectation. Put
$a_v=e^{\chi\Phi_\tau(v,\eta_{\rm m})}$ and $b_v=U_vV_v$.
The scalar moment bound at $2\chi$ gives
\begin{equation}
 \E a_v^2+\E V_v^2\le CN^{4+\sqrt2},\qquad
 \E b_v^2\le CN^{2A_*+4+\sqrt2}.
 \label{lf:eq:lfl-second-moments}
\end{equation}
Let $\Omega$ intersect the two reconstruction events and both
reference caps, with $\PP(\Omega^c)\le C_{D_0}N^{-D_0}$. This joint
event is used only in the comparison of expectations; it is not part
of either $U_v$.
On $\Omega$, $b_v/a_v=\exp(O(\rho_N))$, so, for large $N$,
\[
 |a_v-b_v|\one_\Omega\le C\rho_N a_v\one_\Omega.
\]
Its expectation is at most $C\rho_N m_v^{\rm m}$ directly; this
uses no estimate of $\E b_v$. On the complement,
Cauchy--Schwarz and $m_v^{\rm m}\ge cN$ imply
\begin{equation}
 \frac{\E[(a_v+b_v)\one_{\Omega^c}]}{m_v^{\rm m}}
 \le CN^{A_*+1+1/\sqrt2-D_0/2}.
 \label{lf:eq:lfl-one-bad}
\end{equation}
It follows without any prior asymptotic formula for $\E b_v$ that
\begin{equation}
 \frac{\E|a_v-b_v|}{m_v^{\rm m}}
       \le C\rho_N+CN^{A_*+1+1/\sqrt2-D_0/2}.
 \label{lf:eq:lfl-L1-weight}
\end{equation}
Since $U_v$ and $V_v$ are independent,
$\E b_v=(\E U_v)(\E V_v)$. Taking $D_0$ sufficiently large
proves the first part of \eqref{lf:eq:lfl-lf2};
\eqref{lf:eq:lfl-scalar2} proves the version at height zero. This calibrates the
product of the exact expectations directly, with no conditional or
survival normalization.

For the pair comparison we need a weighted two-root exceptional bound.
Write $W_v^0=e^{\chi\Phi_\tau(v,0)}/m_v^0$. At charge $4\chi$,
\eqref{lf:eq:lfl-moments} and \eqref{lf:eq:lfl-scalar1} give
\[
 \E(W_v^0)^4\le CN^{12+2\sqrt2}.
\]
For any event $E$ of probability $O(N^{-D_0})$, H\"older gives
\begin{equation}
 \E[W_z^0W_w^0\one_E]
 \le\big(\E(W_z^0)^4\E(W_w^0)^4\big)^{1/4}\PP(E)^{1/2}
 \le CN^{6+\sqrt2-D_0/2}.
 \label{lf:eq:lfl-pair-bad}
\end{equation}
If needed, the capped reference product has the analogous bound
$CN^{2A_*+6+\sqrt2-D_0/2}$ after division by $m_z^0m_w^0$.
Thus one may choose, for example,
\[
 D_0>\max\left\{D+2,
  2(A_*+1+1/\sqrt2+A_0+2),
  2(2A_*+6+\sqrt2+A_0+2)\right\}.
\]
All required moment orders and cap sizes are fixed before applying the
estimates that hold with overwhelming probability.

\subsection{Proof of the upper comparison at height zero}

In this subsection we use positivity to extend the comparison to unregularized determinant products while retaining bounded observations of the initial matrix.
Monotonicity of the regularized random determinant and bounded
deterministic density give the pointwise inequality
\begin{equation}
 \Phi_\tau(v,0)\le\Phi_\tau(v,\eta_{\rm m})+CN\eta_{\rm m}
       \le\Phi_\tau(v,\eta_{\rm m})+Ct^{-q}.
 \label{lf:eq:lfl-monotonicity}
\end{equation}
It remains valid at a zero determinant in the positive-weight sense.
For $0\le F\le1$ measurable with respect to
$\mathcal F_0^{\rm mat}$, the two pointwise reconstructions yield,
on $\Omega$,
\[
 e^{\chi\Phi_\tau(z,0)+\chi\Phi_\tau(w,0)}F
       \le e^{C(\rho_N+t^{-q})}U_zU_wV_zV_wF.
\]
Use \eqref{lf:eq:lfl-pair-bad} on the complement and drop
$\one_\Omega$ from the positive reference integral. Exact
independence then gives
\[
 \E[U_zU_wV_zV_wF]
       =(\E U_z)(\E U_w)\E[V_zV_wF].
\]
Divide by the exact normalizers $m_z^0m_w^0$ and apply
\eqref{lf:eq:lfl-lf2}. The chosen $D_0$ gives $O(N^{-A_0})$ on
the exceptional set, proving \eqref{lf:eq:lfl-lf3} and completing
Theorem~\ref{lf:thm:lfl-main}.

We can factor the two local reference expectations only after bounding
the weighted contribution from the exceptional event in the joint field
approximation.
The error in \eqref{lf:eq:lfl-lf3} is a multiplicative $1+o(t^{-P})$ on the
initial pair expectation, which may grow polynomially. Later
compensation and barrier estimates control that expectation.
All exceptional-set bounds above use positive determinant moments.

Uniformity on fixed bounded positive-charge intervals also enters the
product formula of \cite[Lemma~4.4]{CL}. For the compact deterministic
charge range below, the reconstruction event and local caps are
independent of the charge. Moment bounds on one common positive interval
control the exceptional-set contributions, so the same parameters
compare the exact expectations and give the unequal-charge pair bound.

\begin{proposition}[Uniformity in the positive exponents]
\label{lf:prop:lfl-uniform-charges}
Fix a nonempty compact set $J\subset[2,\infty)$ and write
$\lambda_+=\max J$. For $\lambda\in J$, define
\begin{align*}
 U_v(\lambda)&=e^{\lambda(L_{1,v}+L_{2,v})}\one_{\mathcal K_v},
 &V_v(\lambda)&=e^{\lambda\Phi_0(v,B_v\tau)},\\
 m_v^{\rm m}(\lambda)&=\E e^{\lambda\Phi_\tau(v,\eta_{\rm m})},
 &m_v^0(\lambda)&=\E e^{\lambda\Phi_\tau(v,0)}.
\end{align*}
Under the hypotheses of Theorem~\ref{lf:thm:lfl-main}, for every fixed
$P,D,A_0>0$ one common choice of logarithmic parameters, reconstruction
event and local caps gives, uniformly for $\lambda\in J$,
\begin{equation}
 m_v^{\rm m}(\lambda)\asymp m_v^0(\lambda)
       \asymp N^{\lambda^2/8},\qquad
 \frac{(\E U_v(\lambda))(\E V_v(\lambda))}
      {m_v^a(\lambda)}=1+o(t^{-P}),\quad a\in\{{\rm m},0\}.
 \label{lf:eq:lfl-uniform-calibration}
\end{equation}
For arbitrary deterministic charges $\chi_z,\chi_w\in J$ and every
$\mathcal F_0^{\rm mat}$-measurable $0\le F\le1$, the same choice gives
\begin{align}
 &\frac{\E[e^{\chi_z\Phi_\tau(z,0)+\chi_w\Phi_\tau(w,0)}F]}
              {m_z^0(\chi_z)m_w^0(\chi_w)}\notag\\
 &\qquad\le(1+o(t^{-P}))
   \E\left[\frac{V_z(\chi_z)V_w(\chi_w)}
      {(\E V_z(\chi_z))(\E V_w(\chi_w))}F\right]
       +O(N^{-A_0}).
 \label{lf:eq:lfl-uniform-pair}
\end{align}
All constants and remainders are uniform in these deterministic
charges, including charges depending on $N$, and in the root geometry
\eqref{lf:eq:lfl-geometry}. The local factors retain the independence
stated in Theorem~\ref{lf:thm:lfl-main}.
\end{proposition}

\begin{proof}
First enlarge the prescribed common moment bounds to include the
normalized Gaussian convolutions used by the initial and terminal
fields. Fix the scalar results cited above on the single charge
interval $[2,4\lambda_+]$. Their constants and positive rate $\alpha$
are then common to all these laws and charges. In particular,
\[
 m_v^0(\lambda)\asymp N^{\lambda^2/8},\qquad
 \frac{m_v^{\rm m}(\lambda)}{m_v^0(\lambda)}
       =1+O(t^{-q}+N^{-\alpha}),\qquad \lambda\in J,
\]
and their all-height moment bound applies also at $B_v\tau$.
Exact variance scaling preserves these centered-field statements.

The variables $L_{1,v},L_{2,v}$, the caps $\mathcal K_v$ and the
reconstruction event do not depend on the charge. On the same event,
exponentiating \eqref{lf:eq:lfl-lf1} at any $\lambda\in J$ incurs a
relative error $O(\lambda_+\rho_N)$, uniformly. Also
\[
 0\le U_v(\lambda)\le N^{A_J},\qquad
 A_J:=\lambda_+(A_{\rm loc}+1).
\]
For $a_v(\lambda)=e^{\lambda\Phi_\tau(v,\eta_{\rm m})}$ and
$b_v(\lambda)=U_v(\lambda)V_v(\lambda)$, the scalar bound at
$2\lambda$ gives
\[
 \E a_v(\lambda)^2+\E V_v(\lambda)^2
       \le C N^{\lambda^2/2+\lambda/2},\qquad
 \E b_v(\lambda)^2
       \le C N^{2A_J+\lambda^2/2+\lambda/2}.
\]
Thus the normalized exceptional contribution is at most
\[
 \frac{\E[(a_v(\lambda)+b_v(\lambda))\one_{\Omega^c}]}
      {m_v^{\rm m}(\lambda)}
 \le C N^{A_J+\lambda^2/8+\lambda/4-D_0/2}
 \le C N^{A_J+\lambda_+^2/8+\lambda_+/4-D_0/2}.
\]
On $\Omega$, $|a_v(\lambda)-b_v(\lambda)|\le
C\lambda_+\rho_N a_v(\lambda)$. Integration and exact independence
therefore prove \eqref{lf:eq:lfl-uniform-calibration}, uniformly, without
a prior asymptotic formula for $\E b_v(\lambda)$.

For the pair estimate, we put
$W_v(\lambda)=e^{\lambda\Phi_\tau(v,0)}/m_v^0(\lambda)$.
The scalar bound at $4\lambda$ yields
\[
 \E W_v(\lambda)^4\le C N^{3\lambda^2/2+\lambda}.
\]
For any event $E$ with $\PP(E)\le C_{D_0}N^{-D_0}$,
H\"older with exponents $4,4,2$ consequently gives
\begin{align*}
 \E[W_z(\chi_z)W_w(\chi_w)\one_E]
 &\le C N^{\frac38(\chi_z^2+\chi_w^2)
                 +\frac14(\chi_z+\chi_w)-D_0/2}\\
 &\le C N^{3\lambda_+^2/4+\lambda_+/2-D_0/2}.
\end{align*}
The capped reference product, normalized by
$m_z^0(\chi_z)m_w^0(\chi_w)$, has the same bound with the additional
power $2A_J$. One finite choice that makes all these exceptional-set
contributions sufficiently small is
\[
 D_0>\max\left\{D+2,\;
 2(A_J+\lambda_+^2/8+\lambda_+/4+A_0+2),\;
 2(2A_J+3\lambda_+^2/4+\lambda_+/2+A_0+2)\right\}.
\]
Fix this failure order before choosing the common reconstruction and
local-law parameters. The event is charge-independent, so no union
over a continuum of charges is needed.

Finally, \eqref{lf:eq:lfl-monotonicity} and the two reconstructions imply
on $\Omega$ the uniform pointwise estimate
\[
 e^{\chi_z\Phi_\tau(z,0)+\chi_w\Phi_\tau(w,0)}F
 \le e^{C\lambda_+(\rho_N+t^{-q})}
 U_z(\chi_z)U_w(\chi_w)V_z(\chi_z)V_w(\chi_w)F.
\]
We bound the unregularized product on the exceptional event using the
preceding pair estimate, drop the
indicator of $\Omega$ from the nonnegative reference integral, and
factor the two independent reference expectations. Inserting
\eqref{lf:eq:lfl-uniform-calibration}, with $q>P+2$, proves
\eqref{lf:eq:lfl-uniform-pair}. The relative error still multiplies the initial pair expectation.
The uniformity is over deterministic charges, including deterministic
sequences depending on $N$.
\end{proof}

\subsection{Exact variance scaling and sign conventions}

Now we check that variance rescaling preserves the centered field and that the stochastic identities use consistent signs.
Under the variance conversion
$u=\log(1+s)$, $\widetilde X=X_s/\sqrt{1+s}$,
$\widetilde v=v/\sqrt{1+s}$ and
$\widetilde\eta=\eta/\sqrt{1+s}$, the random and deterministic
terms of \eqref{lf:eq:lfl-field} both acquire
$N\log(1+s)/2$, which cancels. The same exact cancellation
applies to the scaled reference Ginibre field. Thus the centered
field, its physical exponent and its exact expectations are preserved
under the simultaneous matrix, spatial-point and height rescaling.

The direct calculation \eqref{lf:eq:lfl-log-ito} agrees with CL
arXiv:2605.29962v3, (3.12); the later display (3.36) prints the
opposite quadratic-drift sign. The sentence preceding (3.35) also
reverses the height derivative obtained from the defining field.
Here both signs follow from \eqref{lf:eq:lfl-field} and the calculation
leading to \eqref{lf:eq:lfl-log-ito}. The drift is estimated in absolute value, so these
sign discrepancies in those displays do not affect the error bounds or the
definition of $L_1+L_2$.

\subsection{Global parameter order and limit conventions}
\label{main:parameter-schedule}

Now we place all parameter choices and limit passages in one order, so later accuracy requirements do not alter earlier assumptions.
Indeed, this section  fixes a common order of parameter choices for the matrix
comparisons and limit passages. The local parameter choices are justified by the estimates in
their respective sections; the order below combines their sufficient bounds.
In particular, use Section~\ref{lf:sec:lfl-parameters},
Section~\ref{oneroot:complete-parameter-order} and the parameter
arguments in Sections~\ref{chap:separated} and~\ref{chap:annuli}.

All limits have a fixed compact spatial set and a fixed finite
original-field test family. Each desired logarithmic accuracy is
fixed before selecting logarithmic cutoff and flow exponents. Coefficient
supports, charge bounds, derivative orders, and cap sizes are fixed
before those exponents. We first fix the bounds needed for the caps, positive moments and
polynomial mesh. We then fix the failure and high-moment orders needed
to control the corresponding exceptional-set contributions, before
choosing the logarithmic field-approximation parameters. All these orders
and exponents are fixed independently of $N$.

Here is one explicit finite order that simultaneously meets the
parameter bounds for the matrix comparisons. We fix an error precision $P>10$, use the
available derivative order $m=8$, and choose $P_*>P+m+10$.
First fix the finite coarse list, smooth cylinder tests, positive
charge intervals and required cap-shell error bounds, then the cap sizes.
For bounded Lipschitz endpoint test functions, we use mollification at
scale $t^{-s}$,
with derivative cost covered by $t^{m(1+s)}$ and restoration weighted
by the surviving mass. The comparison for a compact range of exponents
uses determinant exponents $[2,4]$ and the local-factor moment interval $[2,16]$;
other cap and exceptional-event estimates may require further fixed
entry moments, all supplied by \eqref{main:model}. For the signed tests
use the fixed cap and order-eight construction of
Corollary~\ref{signed:cor:ss-eight}; its shell and scalar-remainder
costs are those of \eqref{signed:eq:ss-shell-budget} and
Lemma~\ref{signed:lem:ss-sixth-remainder}. The weighted local-factor exceptional-set
bounds are fixed in Section~\ref{lf:sec:lfl-weights}, including its
uniform-charge extension. Fix also the replacement parameter
$0<\delta_0<1/1200$ from Theorem~\ref{replacement:thm:rgfc-main}
and the sufficiently small derivative exponent $\varepsilon$ in
Lemma~\ref{replacement:lem:rgfc-input} before these bounds.
The age is $\sigma=N^{-1/200}$; every ultimately fixed $b,q$ has
$t^b/N=o(\sigma)$ and $N^{-1}t^{-q}\ge N^{-1-\delta_0}$ eventually.
Thus these domain requirements do not impose an upper bound on $b$.
Choose the
weak-law and same-noise constants $M,m_0,a_{\rm count}$, with
$a_{\rm count}>6M+3$, and let $C_*$ dominate the fixed comparison and
local-law height floors. Next choose
\[
 s>P_*+4,\qquad q>P_*+4,\qquad
 c_{\rm LF}>\max\{m_0,P_*+3,C_*\}.
\]
The barrier smoothing width is $t^{-s}$ and the microscopic proxy
has height $N^{-1}t^{-q}$, up to a fixed bounded positive factor.
The compensated annular argument retains its width-one smoothing
and $O(t^m)$ derivative cost. It uses $P_*>P+m+10$ as stated in
Section~\ref{chap:annuli}; it does not replace that cost by the
finer retained-barrier smoothing cost or choose $P_*$ again after $s$.
For the main retained-history flow choose
\begin{equation}
\begin{split}
 b_{\rm main}>\max\big\{&100M+2P_*+4,\ c_{\rm LF}+P_*+2,\\
 &2[15M+\max\{2m_0+q,c_{\rm LF}\}+P_*+2],\\
 &2[m(1+s)+P_*+4],\ m(1+s)+P_*+19,\\
 &400,\ C_*+1,\ 24\big\}.
\end{split}
\label{main:parameter-b}
\end{equation}
Then choose $a_{\rm hist}>\max\{5,s+P_*+5\}$ and set
\[
 \tau_{\rm main}=t^{b_{\rm main}}/N,\qquad
 J_{\rm main}=\left\lfloor t-
 (b_{\rm main}/2+a_{\rm hist})\log t-C_K\right\rfloor.
\]
Thus $\eta_{J_{\rm main}}\asymp
N^{-1}t^{b_{\rm main}+2a_{\rm hist}}$.
The logarithmic-floor history corollary in
Section~\ref{chap:calibration}, with $L=b_{\rm main}+2a_{\rm hist}>1$,
gives $O(t^{1-a_{\rm hist}}+t^{1-2a_{\rm hist}}+
t^{2-b_{\rm main}-4a_{\rm hist}})$ before the
$O(t^{1+s})$ first-derivative cost. This is $o(t^{-P_*})$.
The two error powers in the exponential-weight comparison, after the full
derivative cost, are
$t^{2-b_{\rm main}/2+m(1+s)}$ and
$t^{17-b_{\rm main}+m(1+s)}$, both $o(t^{-P_*})$ by
\eqref{main:parameter-b}. The annular errors, after summation and
critical normalization, are $O(t^{2-P})=o(1)$.

We choose the noise windows for the main local-factor comparison in the order already justified:
\begin{align*}
 j&>2c_{\rm LF}+P_*+2,&
 k_1&>\max\{j,b_{\rm main}+2P_*+4\},\\
 R_{\rm noise}&>P_*+q+m_0+3,&
 \ell&>\max\{a_{\rm count}+1,k_1+1,
                         3b_{\rm main}/2+R_{\rm noise}+5\},\\
 k&>\max\{\ell+2,b_{\rm ins}+2,m_0+2\}.&&
\end{align*}
Here $b_{\rm ins}$ is the fixed insertion/overlap floor; none of
these later windows changes $J_{\rm main}$. The already fixed
failure exponents and moment bounds change constants and the threshold in $N$,
not this dependence order.

The separated variable-height scalar upper bound uses a
\emph{second auxiliary coupling}. After fixing $J_{\rm main}$, choose
\[
 c_{\rm scalar}>\max\{c_{\rm LF},
                  b_{\rm main}+2a_{\rm hist}+3\},
 \qquad \eta_{*,\rm scalar}=t^{c_{\rm scalar}}/N.
\]
Choose $b_{\rm scalar}$ by the same finite local-factor bounds with
$c_{\rm scalar}$ replacing $c_{\rm LF}$, and then choose its own
later noise windows. This gives
$\eta_{J_{\rm main}}<\eta_{*,\rm scalar}\ll
\tau_{\rm scalar}=t^{b_{\rm scalar}}/N$ for large $N$.
That calculation has no retained-history test and estimates the same
original scalar quantity through a separate moment-matched comparison matrix.
In particular, the choice of $b_{\rm scalar}$ does not alter
$b_{\rm main},a_{\rm hist}$ or $J_{\rm main}$.
The fixed-fraction prefix $J_\varepsilon=\lfloor(1-\varepsilon)t\rfloor$
used for the global event and nearby masses is a third, distinct
cutoff. Gaussian bridge and interpolation parameters are chosen in
their own finite parameter choices before the matrix limits; they do not
alter these flow choices.

The error in replacing the retained-history density is multiplied by
the test supremum norm. The error from transporting the heights is
multiplied by the aggregate first-derivative bound
$O(t^{1+s})$. The signed comparison and Gaussian interpolation involve
their own displayed derivative norms. In particular, the full derivative
bound for the fine barrier test $t^{m(1+s)}$ does not enter the error from replacing the density by its local-factor approximation.
For the separated comparison take its smoothing exponent $A_s=s$;
the stronger bounds above cover its fixed $P>3$ requirements. The
compact-window desmoothing bound \eqref{window:eq:cw-desmoothed}
is applied with $Q\le P$ and $Q_*=P_*$. The annular width-one cost
remains $O(t^m)$, with $P_*$ fixed before $s$.

The aggregate comparison uses a separate cutoff that is a fixed power of the matrix size.
Fix its $\epsilon>0$, polynomial-list exponent $L$ and common moment bounds,
then use the sixth-order entry expansion and local-law slack
$\xi=\epsilon/100$. The common-event failure exponents and finite entry-moment bounds
of Section~\ref{aggregate:sec:event} are fixed after this polynomial
test count, before $N\to\infty$. They impose no condition on
$b_{\rm main},a_{\rm hist}$ or the auxiliary scalar flow.
The $L=9$ choice is the global-prefix application only; it does not
restrict the standalone theorem. Proposition~\ref{p34:global} uses the fixed original entry law.
Its finite-$N$ comparison may depend on fixed $\epsilon$; the tail in
Corollary~\ref{aggregate:cor:au-maximum} is independent of $\epsilon$
only after the matrix limit.

For the change of Gaussian weighting scale, Section~\ref{barrier:sec:root} fixes its
buffer coefficient $A_G>2$ and strip remainder order $D_G>1/2$
before the matrix limit. The additive comparison and strip estimates
\eqref{gaussian:gpm:additive-comparison} and
Lemma~\ref{gaussian:gpm:strips} account for their intermediate threshold
lifts within their proofs; the final strip and root-return conclusions
concern the prescribed barrier $B$.
These Gaussian auxiliary choices do not redefine the retained matrix
cutoff. They are independent of the aggregate list exponent $L$.

The Gaussian mark proof of Theorem~\ref{gaussian:gpm:factorization} instead
uses its own prefix of order $t^{3/4}$ for a pair and $t^{1/2}$ for
one root. Choose its smoothing power $P_0>3$, then the gap coefficient
$K_G>P_0+4$ and $g=\lceil K_G\log t\rceil$, before the large-time
limit. The gap choice gives the interpolation error
$t^2t^{2P_0}e^{-2g}=o(t^{-1})$; the later shell, mark-shift and
killed-law product errors retain their prefix-survival factors.
The multiplicative
external-root estimate of Theorem~\ref{gaussian:gpm:external-theorem}
is used in annular geometry; an additive error is not substituted
for its relative bound.

At fixed barrier, separated spatial matching takes the matrix limit
before the diagonal-depth cutoff tends to infinity, keeping the outer
tail, terminal annular sum, summed comparison error and cut core
separate. Finite-cylinder approximation likewise fixes its finite list
at each approximation precision before the matrix limit. These are
the independent approximation orders used in Section~\ref{chap:fixed}.

The canonical Gaussian construction and fixed-$B$ estimates precede
identification of the matrix limit. The order for removing the barrier is $N\to\infty$, then
$\varepsilon\downarrow0$ at fixed $B$, then $B\to\infty$.
For original-field approximation the order is fixed coarse scale,
$N\to\infty$, quadrature refinement, then coarse scale to infinity.
For a compact window of radius $A_{\rm win}$, first fix an auxiliary
radius $b_{\rm tail}>A_{\rm win}$ for \eqref{cw:tail-domination};
the requirement $t\ge4b_{\rm tail}^2$ is eventual, not a choice depending
on $N$. This radius is distinct from both flow exponents. The mark
cutoff is removed after the matrix limit, as in
Section~\ref{window:sec:tails}. None of these steps claims rates
uniform in the parameters that are sent to their limits later.

The external estimates and sign conventions are from \cite{CL}, in particular
Lemma~2.12 and equations~(3.12), (3.35)--(3.36); the sign qualification
above is retained. Strong well-posedness is supplied by
\cite[Appendix B, Proposition B.1]{BCH}. The logarithmic coupling
estimates are proved in Sections~\ref{chap:noise}--\ref{chap:different}.

\endgroup

\begingroup
\newcommand{\E}{\mathbb E}
\newcommand{\PP}{\mathbb P}
\newcommand{\Tr}{\operatorname{Tr}}
\newcommand{\br}[1]{\langle #1\rangle}
\newcommand{\D}{\mathcal D}
\newcommand{\HS}{\mathrm{HS}}
\newcommand{\ind}{\mathbf1}
\section{Sparse exponential weights and Stein comparison}\label{chap:signed}

The purpose of this section is to compare exponentially weighted field histories with Gaussian vectors, using capped coordinates to control errors relative to exact normalizing constants.

We compare the logarithmic field under exponential weights involving
at most three coordinates. The coefficients of the exponent may have
either sign; smooth coordinate caps make every such weight bounded
and positive. The term signed source below refers to the coefficient
vector. The comparison uses the full-block insertion estimates of
Section~\ref{chap:i2}, the exact complex cumulant expansion, and
relative error bounds for the weighted expectations. It gives the
normalizing constant for the weight and a Gaussian test-function
comparison. This normalizing constant is derived from the Stein
identity after its relative error has been controlled.

Signed Laplace parameters and capped auxiliary fields already occur in
\cite[Propositions~2.15 and~7.5]{CL}. The extension proved here is a bound for the
relative error in the Stein identity at logarithmic cutoffs, uniformly for the stated
growing history lists and coefficient vectors supported on at most
three coordinates, with the summed derivative
cost in \eqref{signed:eq:ss-norm}.

\subsection{Coordinates, exponential weights, and statement}

Now we specify the capped coordinates, positive and negative exponent coefficients, and derivative bounds used in the weighted Gaussian comparison.
Let $\widehat X=N^{-1/2}(\widehat\xi_{ab})$ be complex iid with
\[
 \E\widehat\xi=\E\widehat\xi^2=0,\qquad
 \E|\widehat\xi|^2=1,\qquad
 \sup_N\E|\widehat\xi|^p<\infty\quad(p<\infty).
\]
Put $\sigma=N^{-1/200}$ and
\[
 X=\sqrt{1-q_N}\widehat X+\sqrt{q_N}G,\qquad
 c\sigma\le q_N\le C\sigma,
\]
where $G$ is independent normalized complex Ginibre.
The fourth cumulant $\kappa$ below is that of the normalized
entries of $X$, so $\kappa=\E|\xi|^4-2\ge-1$.
Write $t=\frac12\log N$.

There are $n\le C_0t$ deterministic coordinates $(z_i,\eta_i)$
in a fixed bulk compact. The fine coordinates lie at two spatial
points $z,w$ with $|z-w|\le C N^{-1/4}$, and obey
\begin{equation}
 \eta_i\in[\eta_{\rm mes},1],\qquad
 \eta_{\rm mes}\ge c_b t^b/N.
 \label{signed:eq:ss-heights}
\end{equation}
This list includes both histories and their external roots at
their exact real clocks. A fixed finite additional list has
arbitrary spatial points in the bulk and heights in $[a_0,1]$,
where $a_0>0$ is fixed. Every capped coordinate has height at most one, as required by
the scalar tail estimate below.

With normalized trace $(2N)^{-1}\Tr$, set
\[
 H^z=\begin{pmatrix}0&X-z\\(X-z)^*&0\end{pmatrix},
 \quad G_z(i\eta)=(H^z-i\eta)^{-1},\quad
 M_z(i\eta)=
 \begin{pmatrix}m&-zu\\-\overline z u&m\end{pmatrix}.
\]
The displayed $2$ by $2$ matrices are tensored with $I_N$.
Write $m_i=m^{z_i}(i\eta_i)$ and
\begin{equation}
 E_i=-\frac\kappa4m_i^4,\qquad
 Y_i=\sqrt2(\Phi_N(z_i,\eta_i)-E_i),\qquad
 \Sigma_{ir}=2\left(-\frac14\log D_{ir}
                    +\frac\kappa4m_i^2m_r^2\right).
 \label{signed:eq:ss-centering}
\end{equation}
Here $D_{ir}$ is the exact stability determinant
\[
 D_{ir}=1+(u_i u_r|z_i z_r|)^2-(m_i m_r)^2
                   -2u_i u_r\Re(z_i\overline z_r).
\]
This is the restriction of the common covariance
\eqref{gaussian:gp2:covariance} to the exact clocks
$s_i=\operatorname{arctanh}u_i$: the MDE relation gives
$-m_i^2=u_i(1-u_i|z_i|^2)$. Thus the finite-height and Gaussian
conventions agree, including the fourth-cumulant term. We recall
the feature expansion in these endpoint variables. Put
$p_i=-m_i^2>0$ and $q_{ir}=u_i u_rz_i\overline z_r$.
Then $D_{ir}=|1-q_{ir}|^2-p_ip_r>0$, and
\[
 -\tfrac12\log D_{ir}
 =\Re\sum_{k\ge1}\frac{q_{ir}^k}{k}
  +\frac12\sum_{a\ge1}\frac{(p_ip_r)^a}{a}
       (1-q_{ir})^{-a}(1-\overline q_{ir})^{-a}.
\]
As in \eqref{gaussian:gp2:features}, expanding the two binomial factors
writes every term as a nonnegative multiple of a Gram kernel. The term $a=1$ with
both binomial indices zero has coefficient $p_ip_r/2$.
Adding the cumulant term changes it to
$(1+\kappa)p_ip_r/2\ge0$. The series converge at the positive
heights in question. A singular covariance is allowed.

Fix a smooth odd nondecreasing saturation map $c_A$ satisfying
\[
 c_A(x)=x\ (|x|\le A\log N),\qquad
 |c_A(x)|\le2A\log N,\qquad 0\le c_A'(x)\le1,
\]
with all fixed derivatives bounded uniformly in $N$. In particular
$c_A(x)\ge A\log N$ for $x\ge A\log N$, and the reflected
inequality holds on the negative side.
Put $\mathbf X=c_A(\mathbf Y)$ coordinatewise.
For a deterministic real coefficient vector $a\in\mathbb R^n$ with at
most three nonzero coordinates and $\|a\|_1\le6$, the exponential
weight is normalized by its exact mean to define a new law. Along the coefficient ray
$sa$, define
\[
 Z(sa)=\E e^{sa\cdot\mathbf X},\qquad
 dQ_{sa}=Z(sa)^{-1}e^{sa\cdot\mathbf X}\,d\PP,
 \qquad 0\le s\le1.
\]
The exponential weight is positive for every sign pattern; negative
coefficients do not require negative moments of the unregularized
determinant because the coordinates are capped.

For a bounded smooth scalar function let
\begin{equation}
 \D_k(f)=\sum_{q=0}^k\sum_{r_1,\ldots,r_q}
       \|\partial_{r_1}\cdots\partial_{r_q}f\|_\infty,
 \qquad
 \|\mathbf F\|_{k,1}=\sum_i\D_k(F_i).
 \label{signed:eq:ss-norm}
\end{equation}
The order-zero term is included. Summing the supremum of each
derivative accounts for all coordinate choices, including the
$O(\log N)$ history list, without an additional dimension factor.

For a centered Gaussian vector with covariance $\Sigma$, exponential weighting
by $sa$ shifts its mean to $s\Sigma a$ and leaves its covariance
unchanged, even when $\Sigma$ is singular. The following theorem
bounds the corresponding integration-by-parts defect for the matrix
field.

\begin{theorem}[Stein defect for sparse exponential weights]
\label{signed:thm:ss-defect}
The cap $A$ can be fixed sufficiently large, and then a finite
integer $m$ can be chosen independently of $b,N,n$, such that
the following holds for all sufficiently large fixed $b$.
Uniformly for the stated coefficient rays and coordinate lists,
\begin{align}
 \left|\E_{Q_{sa}}\left[
  (\mathbf X-s\Sigma a)\cdot\mathbf F(\mathbf X)
  -\sum_{i,r}\Sigma_{ir}\partial_rF_i(\mathbf X)\right]\right|
 &\le r_N\|\mathbf F\|_{m-1,1},                 \label{signed:eq:ss-defect}\\
 r_N&\le C_b\left\{
       \frac{t^2}{\sqrt{N\eta_{\rm mes}}}
       +\frac{t^{17}}{N\eta_{\rm mes}}+N^{-c}\right\}.
                                                        \label{signed:eq:ss-rate}
\end{align}
Here $c>0$ is fixed. The derivative order and displayed logarithmic
exponents are selected before $b$; constants and the eventual
threshold in $N$ may depend on $b$.
\end{theorem}

\subsection{Unweighted inputs and the lower bound for the normalizing constant}

In this subsection we control the fields under the original law before bounding the normalizing constant from below.

We need local laws and scalar tails under the unweighted matrix law
before introducing the weighted law. The local laws control the cumulant
coefficients; the tails give a preliminary lower bound for the exact
normalizing constant.
The logarithmic isotropic one-resolvent law is
Cipolloni--Landon, \emph{Optimal Delocalization for Non-Hermitian
Eigenvectors}, arXiv:2509.15189v1, Theorem 3.1.
On a fixed bulk and a fine height domain it is stronger than
\[
 |v^*(G(iu)-M(iu))w|\le\frac{Ct^2}{\sqrt{Nu}}
\]
for deterministic unit vectors. Its averaged estimate against each deterministic block-basis insertion
gives $Ct^2/(Nu)$. Choose a small fixed $0<\xi<1/100$ and put $a_N=N^{-10\xi}$.
Take the upper-domain exponent in that theorem smaller than
$10\xi$. Above $a_N$, use the ordinary fixed-power local law
in \cite[Theorem 2.3]{CL}.
On a common event the resulting envelope is
\begin{equation}
 \Psi(u)=\frac{Ct^2+N^\xi\ind_{\{u\ge a_N\}}}{\sqrt{Nu}},
 \qquad \eta_{\rm mes}\le u\le1.
 \label{signed:eq:ss-psi}
\end{equation}
The deterministic vector list includes all coordinate pairs,
coordinate/block-constant pairs in both orders, and pairs of
normalized vectors constant on either Hermitization component.
It has polynomial size. Multiplication by a fixed block-constant
matrix reduces to a finite linear combination of this list.
No random vector is used as a test vector of a local law.

Corollaries~\ref{i2:cor:I6} and \ref{i2:cor:second} supply, on a common
event, the mixed coefficient error
\begin{equation}
 \max_{i,r}|\mathcal V_{ir}-V_{ir}|
 \le Ct^{17}/(N\eta_{\rm mes})+N^{-c},
 \qquad V_{ir}=-\tfrac14\log D_{ir}.
 \label{signed:eq:ss-I6}
\end{equation}
It includes arbitrary coarse spatial points, both orientations,
the insertion $M_i/(1-\br{M_i^2})$, and the outer integration
through $N^{10}$. The two-resolvent estimate is needed for this mixed coefficient;
the one-resolvent law controls the other entrywise errors.

The scalar tail bound \cite[Lemma 8.1]{CL}, with the harmless
bounded mean correction and $\sqrt2$ rescaling included, gives
\begin{equation}
 \PP\{|Y_i|>A\log N\}
 \le C N^{-c_0A^2}+N^{-D}
 \label{signed:eq:ss-tail}
\end{equation}
uniformly for heights between its fixed logarithmic floor and
one. The exponent $b$ is chosen above that floor exponent.
The exceptional exponent $D$ can be selected after all fixed
cap and derivative parameters.

Choose a fixed moderate $A_1$ so that a union of
\eqref{signed:eq:ss-tail} over the $n=O(t)$ coordinates gives
$\PP(\max_i|Y_i|\le A_1\log N)\ge1/2$.
Choose $A>A_1$. Before any comparison or normalizer asymptotic,
\begin{equation}
 Z(sa)\ge\tfrac12N^{-6A_1},\qquad 0\le s\le1.
 \label{signed:eq:ss-prelower}
\end{equation}
Indeed the exponent in the weight is at least $-6A_1\log N$ on that
event. Also
\[
 e^{sa\cdot\mathbf X}\le N^{12A},\qquad
 \PP(\max_i|Y_i|\ge A\log N)
       \le CtN^{-c_0A^2}+N^{-D}.
\]
The raw finite-field linear factor is $O(N\log N)$ and the
covariance coefficients are $O(t)$, so $N^2$ bounds all additional
factors on this tail event for large $N$. After division by
\eqref{signed:eq:ss-prelower}, the resulting error is at most
\begin{equation}
 C N^{12A+6A_1+2+o(1)-c_0A^2}
       +C N^{12A+6A_1+2-D}.
 \label{signed:eq:ss-shell-budget}
\end{equation}
Choose $A$ first to make this smaller than any prescribed fixed
polynomial accuracy. For a sufficiently large cap, the negative quadratic
term in the exponent dominates the positive linear terms. This choice does not depend
on the later cumulant truncation order.

\subsection{Finite fields and exact resolvent identities}

In this subsection we write the logarithms as finite resolvent integrals to obtain exact entry derivatives without losing control of the normalizing constants.
To apply entrywise cumulant expansion, we first replace the logarithm
by an integral with a deterministic upper endpoint. This gives a
deterministic field bound and exact entry derivatives. Let $R=N^{10}$
and work initially with
\begin{equation}
 U_i=iN\int_{\eta_i}^{R}
       \br{G_i(iu)-M_i(iu)}\,du-E_i,\qquad
 \widetilde Y_i=\sqrt2\,U_i.
 \label{signed:eq:ss-finite}
\end{equation}
These variables are real and deterministically $O(N\log N)$.
On an overwhelming norm event, $U_i-(\Phi_i-E_i)=O(N^{-5})$.
We define the finite capped fields and their exact normalizing constants by
\[
 \widetilde{\mathbf X}=c_A(\widetilde{\mathbf Y}),\qquad
 \widetilde Z(sa)=\E e^{sa\cdot\widetilde{\mathbf X}},\qquad
 d\widetilde Q_{sa}=\widetilde Z(sa)^{-1}
               e^{sa\cdot\widetilde{\mathbf X}}\,d\PP.
\]
The tildes distinguish these finite-integral fields, normalizing constants
and weighted laws
from the actual $\mathbf X,Z(sa),Q_{sa}$. The following transfer uses
only the preliminary lower bound, so it precedes any asymptotic
formula for the normalizing constant.
For every prescribed fixed $L$, taking the norm-event failure
exponent sufficiently large gives, uniformly over the coefficient rays,
\begin{equation}
 \frac{\widetilde Z(sa)}{Z(sa)}
   =1+O(N^{-5}+N^{-L}),\qquad
 \left\|\frac{e^{sa\cdot\widetilde{\mathbf X}}}{\widetilde Z(sa)}
       -\frac{e^{sa\cdot\mathbf X}}{Z(sa)}\right\|_{L^1(\PP)}
   \le C(N^{-5}+N^{-L}).
 \label{signed:eq:ss-finite-transfer}
\end{equation}
Indeed the cap is Lipschitz, and $\|sa\|_1\le6$, so the two
exponential weights have ratio $1+O(N^{-5})$ on the comparison event.
On its complement both weights are at most $N^{12A}$.
Dividing the resulting absolute complement error by the already
proved actual lower bound \eqref{signed:eq:ss-prelower} costs at most
$CN^{12A+6A_1-D}$. Choose $D>12A+6A_1+L+2$.
This proves the ratio of normalizing constants without an asymptotic
formula for the normalizing constant of the finite-integral field;
the normalized $L^1$ estimate follows from
$\E|u/\E u-v/\E v|\le2\E|u-v|/\E v$ for positive weights.
In particular $\widetilde Z(sa)\ge\frac14N^{-6A_1}$ for large $N$.
The contribution from the finite cap shell
$\{\max_i|\widetilde Y_i|>A\log N\}$ obeys the same bound as
\eqref{signed:eq:ss-shell-budget}, after decreasing $c_0$ once: on the
comparison event, this tail event implies $\max_i|Y_i|\ge(A/2)\log N$ for
large $N$, and the contribution of the complement is bounded as above.
The preceding choice of $A$ is made with this smaller tail constant.
All exact exponential-weight expectations for the finite-integral fields below therefore use
$\widetilde Z(sa)$. Actual normalizers in the theorem and its
corollaries remain $Z(sa)$.

All fixed entry derivatives are endpoint resolvent products; in
particular, for an oriented entry $h=W_{\alpha\beta}$,
\begin{equation}
 \partial_h U_r=
 \tfrac12\{G_r(i\eta_r)-G_r(iR)\}_{\beta\alpha}.
 \label{signed:eq:ss-endpoint}
\end{equation}
Here $W=H^0(X)$ and the spatial shift is deterministic. This is the
entry convention of \eqref{i2:eq:endpoint-derivative}: the oriented
entry already has variance $1/N$, and its derivative carries no
additional $N^{-1/2}$ factor.

For a scalar test $f$, use
\[
 H(y)=e^{sa\cdot c_A(y)}f(c_A(y)).
\]
Product and chain rules, with $\|sa\|_1\le6$, show that every
fixed coordinate derivative sum satisfies
\begin{equation}
 \sum_{q\le k}\sum_{r_1,\ldots,r_q}
 |\partial_{r_1}\cdots\partial_{r_q}H(y)|
 \le C_{k,A}e^{sa\cdot c_A(y)}\D_k(f).
 \label{signed:eq:ss-envelope}
\end{equation}
The coordinate sums are already in $\D_k$; they are not counted
again as $n^k$. Entry derivatives use
\eqref{signed:eq:ss-endpoint} and its higher-order versions.
In the cumulant calculation with the unscaled $\Phi$ coordinates below, use
\[
 \mathcal H(U)=H(\sqrt2\,U),\qquad
 \partial_{r_1}\cdots\partial_{r_k}\mathcal H(U)
 =2^{k/2}H_{r_1\cdots r_k}(\sqrt2\,U).
\]
In the coefficient calculation below, derivatives of $\mathcal H$
are therefore taken with respect to the unscaled coordinates.
The fixed chain factors are included in \eqref{signed:eq:ss-envelope};
the covariance rescaling is performed explicitly after this calculation.

For one output coordinate write
\[
 \Delta=\br G-m,\qquad
 A=\frac M{1-\br{M^2}},\qquad
 P=GA,\qquad
 \underline{WG}=WG+\br G G.
\]
The denominator uses the normalized matrix trace $\br{M^2}$,
including the off-diagonal MDE entries.
Bulk stability gives $\|A(iu)\|\le C/(1+u)$.
Since $M^{-1}=-(iu+Z+m)$ and
$(W-Z-iu)G=I$, one has
$\underline{WG}=I-M^{-1}G+\Delta G$.
Rearranging and taking traces yields the exact identity
\begin{equation}
 \Delta=-\br{\underline{WG}A}
                +\Delta\br{(G-M)A}.
 \label{signed:eq:ss-scalar-identity}
\end{equation}
The second term contributes, after multiplication by $N$ and
integration,
\begin{equation}
 \frac CN\int_{\eta_i}^1
 \frac{(t^2+N^\xi\ind_{\{u\ge a_N\}})^2}{u^2}\,du
 \le Ct^4/(N\eta_{\rm mes})+N^{-c}.
 \label{signed:eq:ss-quadratic}
\end{equation}
Here the averaged input controls both $\Delta$ and
$\br{(G-M)A}$, by decomposing the deterministic block matrix $A$
in that basis. Above one, the averaged estimate of
Lemma~\ref{signed:lem:ss-tail}, together with $A=O(u^{-1})$,
gives the tail $O(N^{-1+2\xi}\log N)$ computed after that lemma.

\subsubsection{The complex cumulant expansion}

Now we record the cumulant expansion for complex entries and their conjugates, fixing the normalization of the mean and covariance coefficients.
The multivariate cumulant algebra and its complex formulation are
standard; see \cite[Proposition~3.2 and Appendix~C]{EKS}.
We retain the expansion below to fix the entry-derivative conventions.
The relative weighted remainder is estimated separately in
Lemma~\ref{signed:lem:ss-sixth-remainder}.

For an upper-block entry $h=W_{\alpha\beta}$, write $x=h$ and
$y=\overline h$. Treat $x,y$ as the two formal variables of the
complex entry derivatives, condition on the other independent entries,
and use the convention
\begin{equation}
 \E[hF]=\sum_{p+q\le K}
 \frac{\kappa_{p+1,q}(h)}{p!\,q!}
       \E[\partial_x^p\partial_y^q F]
       +\operatorname{Rem}_K.
 \label{signed:eq:ss-cumulant}
\end{equation}
Thus the covariance term is $N^{-1}\E\partial_yF$, and the
balanced fourth term is
$\kappa(2N^2)^{-1}\E\partial_x\partial_y^2F$.
The reverse outer orientation has the reflected formula.
Writing $\kappa_{p,q}$ for cumulants of the unnormalized entry,
the complete outer-$x$ operators at the next two orders are
\begin{align*}
 N^{-3/2}\left(
  \tfrac{\kappa_{3,0}}2D_x^2
       +\kappa_{2,1}D_xD_y
       +\tfrac{\kappa_{1,2}}2D_y^2\right),\\
 N^{-2}\left(
  \tfrac{\kappa_{4,0}}6D_x^3
  +\tfrac{\kappa_{3,1}}2D_x^2D_y
  +\tfrac{\kappa_{2,2}}2D_xD_y^2
  +\tfrac{\kappa_{1,3}}6D_y^3\right).
\end{align*}
For the outer $y$, the coefficient at $D_x^pD_y^q$ is instead
$N^{-(p+q+1)/2}\kappa_{p,q+1}/(p!q!)$.
The operators retain all unbalanced third and fourth cumulants.
Only the first two entry moments are fixed by the model. The absolute
remainder of this fixed-order expansion is estimated later.
The expansion is applied to the raw product $WP$ together with its
multiplier, so every resulting coefficient and test derivative stays
inside the same expectation.

Since $\partial_yP_{\beta\alpha}
=-G_{\beta\beta}P_{\alpha\alpha}$, the pure covariance term,
after both outer orientations and the normalized trace, is
\[
 -\frac1{2N^2}\left[
  \Big(\sum_{\beta\in2}G_{\beta\beta}\Big)
  \Big(\sum_{\alpha\in1}P_{\alpha\alpha}\Big)
 +\Big(\sum_{\alpha\in1}G_{\alpha\alpha}\Big)
  \Big(\sum_{\beta\in2}P_{\beta\beta}\Big)\right].
\]
The two diagonal-block traces of a square Hermitization resolvent
are equal, so this expression is exactly
$-\br G\br{GA}$. It cancels the renormalization in
\eqref{signed:eq:ss-scalar-identity}, with the scalar multiplier still
present. There is no leading Gaussian mean term.

\subsubsection{The mixed second-cumulant coefficient}

Now we identify the logarithmic covariance term while keeping its random error inside the weighted expectation.
In the raw $\Phi$ convention, the remaining second-cumulant
coefficient is
\begin{equation}
 \mathcal V_{ir}=
 \Re\left[-\frac i2\int_{\eta_i}^{R}
 \sum_{p\ne q}
 \br{G_i(iu)A_i(u)E_{pp}
      \{G_r(i\eta_r)-G_r(iR)\}E_{qq}}\,du\right].
 \label{signed:eq:ss-mixed-second}
\end{equation}
The factor $1/2$ comes from
\eqref{signed:eq:ss-endpoint}. There are exactly two unweighted block
orientations; taking real parts adds no further factor.
The matrix order is retained, so $A_i$ is not commuted through
a random resolvent.

Lemma~\ref{i2:lem:primitive} identifies the deterministic primitive
in these same trace and insertion conventions. It follows from
\[
 \partial_{w_i}\log D_{ir}
 =-2\sum_{p\ne q}
      \br{M_{ir}[A_iE_{pp}]E_{qq}},
\]
and equals $-\log D_{ir}/4$ after integration.
Equation \eqref{signed:eq:ss-I6} bounds its random error.
All test derivatives remain inside the expectation; on the
common event their multiplier is bounded by
\eqref{signed:eq:ss-envelope}. Thus the error is proportional to the
exact expectation of the exponential weight, without any independence assumption.
This supplies the logarithmic part of the raw covariance. The next
calculation identifies the third-order error and the fourth-cumulant
mean and covariance corrections before converting to the $\sqrt2$-scaled
real field.

\subsection{Complete third- and fourth-order estimates}

This subsection separates the fourth-cumulant mean and covariance corrections from the terms that vanish in the comparison.
The third cumulants contribute only errors. At fourth order, two
balanced contractions survive: one gives the deterministic leading
mean and the other the fourth-cumulant part of the covariance. We
separate these terms using physical sites, then estimate their
replacement errors and the remaining complex orientations.

Group Hermitization indices by physical site, each carrying two
internal components. The deterministic MDE and $MA$ are constant
within each physical site and vanish between distinct sites.
Different spatial parameters change their coefficients but not
this property. Separate the $O(N)$ equal-site pairs first.
For outer heights at most one, bounded resolvent entries give costs
$O(N^{-1/2})$ at third cumulant order and $O(N^{-1})$ at fourth
order. The large-height estimates below supply integrable decay
for the same contributions above one.

For distinct sites, call a factor joining them off-site. Its
deterministic MDE value vanishes, so its error is bounded by $\Psi$ from
\eqref{signed:eq:ss-psi} on the ordinary event. All other entries
are bounded. The following estimates concern outer heights at most one.
Ward gives
\[
 \|G_r\|_{\HS}\le C\sqrt{N/\eta_r},\qquad
 \|P_i(u)\|_{\HS}\le C\sqrt{N/u}.
\]
Let $\Psi_r=\Psi(\eta_r)$ and
$q_N^{\rm loc}=Ct^2/\sqrt{N\eta_{\rm mes}}+CN^{-1/2+6\xi}$.
For sufficiently large $b$, $q_N^{\rm loc}=o(1)$.
Splitting at $a_N$ gives
\begin{align}
 I_1:=\int_{\eta_i}^1\Psi(u)\,du&\le CN^{-1/2+2\xi},
 &I_2:=\int_{\eta_i}^1\Psi(u)^2\,du&\le CN^{-1+3\xi},
 \nonumber\\
 I_3:=\int_{\eta_i}^1\frac{\Psi(u)}{\sqrt u}\,du
       &\le CN^{-1/2+2\xi},
 &I_4:=\frac1{\sqrt N}\int_{\eta_i}^1\frac{\Psi(u)}u\,du
       &\le\frac{Ct^2}{N\sqrt{\eta_{\rm mes}}}
                                      +CN^{-1+6\xi}.
 \label{signed:eq:ss-integrals}
\end{align}
Logarithms are absorbed into $N^\xi$ only in terms already
polynomially decreasing, never next to an exclusively logarithmic
fine-scale denominator.

\subsubsection{Third cumulants}

Now we show that every third-cumulant term is small enough for the weighted comparison at logarithmic regularization scales.
For labelled entry derivatives $D_1,\ldots,D_k\in\{\partial_x,\partial_y\}$,
write $[k]=\{1,\ldots,k\}$ and let $\mathcal P(T)$ be the set of
partitions of $T$. Put $D_T=\prod_{\ell\in T}D_\ell$, with
$D_\varnothing$ the identity. Each $r_B$ below ranges from $1$ to $n$.
All chain multiplicities are fixed by the labelled rule
\[
 D_1\cdots D_k(P\mathcal H(U))
 =\sum_{S\subseteq[k]}(D_SP)
   \sum_{\pi\in\mathcal P([k]\setminus S)}
   \sum_{(r_B)_{B\in\pi}}
     \mathcal H_{(r_B)_{B\in\pi}}
       \prod_{B\in\pi}D_BU_{r_B}.
\]
The empty partition contributes $\mathcal H$. This counts each
labelled allocation and partition once, including repeated entry
directions. At orders two and three it gives exactly the four
and seven classes below.

At third cumulant order there are two entry derivatives. After the
trace and field factors, their common prefactor is $N^{-3/2}$.
The labelled rule gives the four classes
\begin{equation}
 (D^2P_i)\mathcal H,\quad
 (DP_i)\mathcal H_rDU_r,\quad
 P_i\mathcal H_rD^2U_r,\quad
 P_i\mathcal H_{rs}DU_rDU_s.
 \label{signed:eq:ss-third-classes}
\end{equation}
Every product in the first three classes has one or three off-site
factors; the last class has three. This follows by counting
physical-site switches along the derivative edges, independently
of their complex orientations. A block partner from $A_i$ does
not change the site.

For a product with exactly one off-site factor, replace both diagonal factors by their
deterministic MDE values plus errors. If both are deterministic,
the remaining total off-site sum is $O(N)$ by the local law for
two deterministic block-constant vectors. Removing the equal-site
pairs changes it by $O(N)$. This contribution is $O(N^{-1/2})$.

With exactly one diagonal error, retain the deterministic-vector
contraction by summing its free index first. For example,
\[
 \sup_\beta\left|\sum_\alpha(G_r)_{\alpha\beta}\right|
       \le C(1+\sqrt N\,\Psi_r),
\]
so
\begin{equation}
 N^{-3/2}\left|\sum_{\alpha,\beta}
       \delta_\beta(u)(G_r)_{\alpha\beta}\right|
 \le C\{\Psi(u)/\sqrt N+\Psi(u)\Psi_r\}.
 \label{signed:eq:ss-ordered-contraction}
\end{equation}
The diagonal error depends on the remaining index $\beta$, so it is
bounded only after this contraction. With two diagonal errors,
absolute values and Ward give:
\[
 N^{-3/2}\sum_{\alpha,\beta}
 |\delta_\alpha(u)\delta_\beta(u)(G_r)_{\alpha\beta}|
 \le C\Psi(u)^2/\sqrt{\eta_r}.
\]
The scalar random derivatives $\mathcal H,\mathcal H_r,\ldots$ stay outside these
entry-index sums but inside expectation.

These choices give the complete integrated third-order estimates:
\[
\begin{array}{c|c|c}
\text{class}&\text{diagonal errors}&\text{three off-site factors}\\ \hline
D^2P_i&
 CI_3&CI_4\\
(DP_i)DU_r&
 C(I_1/\sqrt N+I_1\Psi_r+I_2/\sqrt{\eta_r})&
 CN^{-1/2}\Psi_r\log N\\
P_iD^2U_r&
 C\Psi_r&C\Psi_r/\sqrt{N\eta_r}\\
P_iDU_rDU_s&
 \text{none}&
 C\min\{\Psi_r/\sqrt{N\eta_s},\Psi_s/\sqrt{N\eta_r}\}.
\end{array}
\]
For the first row all factors occur at $u$; one maximum and two
Hilbert--Schmidt bounds give $I_4$.
For the second row use \eqref{signed:eq:ss-ordered-contraction};
for its product with three off-site factors take the endpoint as the maximum.
For the third row the endpoint diagonal errors multiply
$\int_{\eta_i}^1u^{-1/2}du\le2$.
For the last two rows, if an endpoint is below $a_N$, choose
that fine endpoint as the maximum. Its contribution is at most
$Ct^2/(N\eta_{\rm mes})$. If both endpoints are above $a_N$,
the contribution is polynomially small.
Thus the fine denominators carry only the stated logarithmic losses;
$N^\xi$ occurs in terms with a remaining polynomial saving. The
entire third-order tensor is bounded coefficientwise by
\begin{equation}
 C\left\{\frac{t^2}{\sqrt{N\eta_{\rm mes}}}
        +\frac{t^2}{N\eta_{\rm mes}}+N^{-c}\right\}.
 \label{signed:eq:ss-third-bound}
\end{equation}

\subsubsection{Fourth cumulants and the exact factor two}

Now we compute the surviving fourth-cumulant corrections and check the numerical factors that determine the limiting mean and covariance.
At fourth cumulant order there are three entry derivatives. The
labelled rule gives seven classes:
\begin{align}
 &(D^3P_i)\mathcal H,\quad (D^2P_i)\mathcal H_rDU_r,\quad
 (DP_i)\mathcal H_rD^2U_r,\quad (DP_i)\mathcal H_{rs}DU_rDU_s,
 \nonumber\\
 &P_i\mathcal H_rD^3U_r,\quad
 P_i\mathcal H_{rs}(D^2U_r)DU_s,\quad
 P_i\mathcal H_{rsu}DU_rDU_sDU_u.
 \label{signed:eq:ss-fourth-classes}
\end{align}
Only the first and third classes have products whose factors are all diagonal in their physical-site indices.
We compute these two coefficients before taking expectations, then
bound every remaining class.

\paragraph{The mean correction.}
Write $x=W_{\alpha\beta}$, $y=W_{\beta\alpha}$ and
$s_i(u)=(M_iA_i)_{\alpha\alpha}=\partial_{iu}m_i(u)$.
For distinct sites the only all-diagonal third derivative of
$P_{\beta\alpha}$ is
\begin{equation}
 (\partial_x\partial_y^2P_{\beta\alpha})_{\rm diag}
       =-2m_i^3s_i.
 \label{signed:eq:ss-pure-fourth}
\end{equation}
The derivative edges must alternate $y,x,y$. The two labels
attached to the identical $y$ derivatives give the factor two.
The balanced cumulant coefficient is $\kappa/(2N^2)$.
Including both outer orientations, the trace/field factor
$N/(2N)=1/2$, and the outer $-i$ integration gives
\begin{equation}
 i\kappa\int_{\eta_i}^{\infty}
             m_i(u)^3\partial_{iu}m_i(u)\,du
       =-\frac\kappa4m_i(\eta_i)^4=E_i.
 \label{signed:eq:ss-mean}
\end{equation}
Here $dm_i/du=i\partial_{iu}m_i$ and the upper endpoint vanishes.
The finite endpoint at $N^{10}$ has negligible error.
This mean is removed once by \eqref{signed:eq:ss-finite}.
There are $N(N-1)$ distinct-site pairs; replacing their count by
$N^2$ changes the leading expression by $O(N^{-1})$, already
included with the separately bounded equal-site contribution.
Thus $E_i$ is the deterministic leading mean, not an assertion
that it equals the exact finite-$N$ expectation.

\paragraph{The covariance correction.}
The mixed surviving class is $(DP_i)\mathcal H_rD^2U_r$. The
following display records its leading coefficients after discarding
the upper endpoint. At finite
$R$, the second expression has $m_r^2$ replaced by
$m_r(\eta_r)^2-m_r(R)^2$; this costs $O(R^{-2})$ uniformly.
\[
 (\partial_yP_{\beta\alpha})_{\rm diag}=-m_i s_i,
 \qquad
 (\partial_x\partial_yU_r)_{\rm diag}=-m_r^2/2.
\]
There are two choices of which $y$ derivative hits $P$.
The chain-rule product is therefore
\[
 2(-m_i s_i)(-m_r^2/2)=m_i s_i m_r^2.
\]
After both outer orientations and the cumulant denominator,
the integrated coefficient is
\begin{equation}
 -\frac{i\kappa}2m_r^2
       \int_{\eta_i}^{\infty}m_i(u)\partial_{iu}m_i(u)\,du
       =\frac\kappa4m_i^2m_r^2.
 \label{signed:eq:ss-fourth-covariance}
\end{equation}
Before discarding the upper endpoint, the same calculation gives
\[
 \frac\kappa4
 \bigl(m_i(\eta_i)^2-m_i(R)^2\bigr)
 \bigl(m_r(\eta_r)^2-m_r(R)^2\bigr).
\]
Thus the field derivative's $1/2$, the chain-rule factor two,
and the two outer orientations have all been included.
Passing from the raw $\Phi$ coordinates to
$Y=\sqrt2(\Phi-E)$ multiplies the covariance by two, as in
\eqref{signed:eq:ss-centering}; the mean correction has already
been subtracted once.

\paragraph{The remaining fourth-order terms.}
Unbalanced complex fourth cumulants cannot produce either
product with all factors diagonal in their physical-site indices: including the outer entry, alternation forces
two entries of each orientation. Each such product contains a factor that improves the bound.
The total crude fourth-order pair count is already $O(1)$.
An error in either leading product is bounded by $CI_1+C\Psi_r$.
The second class of \eqref{signed:eq:ss-fourth-classes} has an outer
off-site factor and costs $CI_1$; the fourth has endpoint
off-site factors and costs $C\Psi_r+C\Psi_s$; each final class
has the outer off-site $P_i$ and costs $CI_1$.
Equal-site pairs cost $O(N^{-1})$.
Consequently all fourth-order errors are bounded by
\begin{equation}
 Ct^2/\sqrt{N\eta_{\rm mes}}+CN^{-c}.
 \label{signed:eq:ss-fourth-bound}
\end{equation}

\paragraph{Large-height contributions.}
The preceding estimates treated the outer integral up to height one.
For its remaining interval, we need decay of the random replacement
errors as well as decay of the resolvents. The following lemma and
contractions provide both and also justify the quadratic tail used
in \eqref{signed:eq:ss-quadratic}.

\begin{lemma}
\label{signed:lem:ss-tail}
On one overwhelming event, uniformly for $1\le u\le N^{10}$,
the required entrywise and deterministic-vector tests satisfy
\[
 |(G_i(iu)-M_i(iu))_{ab}|
 +|v^*(G_i(iu)-M_i(iu))w|
 \le CN^{-1/2+\xi}/u .
\]
Consequently the corresponding error for $P_i=G_iA_i$ is
$CN^{-1/2+\xi}/u^2$. For a bounded deterministic block insertion $B$,
\[
 |\br{(G_i(iu)-M_i(iu))B}|\le CN^\xi/N .
\]
\end{lemma}
\begin{proof}
The rigidity theorem \cite[Theorem 2.5 and Lemma 2.1(ii)]{CL}
gives a fixed spectral norm bound with overwhelming
probability.
The deterministic support assertion also holds for the matrix spectral
measure of $M$. Indeed,
\[
 -M(w)^{-1}=(w+m(w))I+Z,\qquad
 \Im M(w)=M(w)^*\bigl(\Im w+\Im m(w)\bigr)M(w)\ge0.
\]
Together with $M(w)=-I/w+O(w^{-2})$, this gives a positive
matrix spectral measure of total mass $I$. Its normalized trace
is $\rho(x)\,dx$. Positivity therefore forces every matrix entry
of that measure to vanish outside the scalar support, which is
uniformly bounded by the cited deterministic support estimate.
Choose a fixed rectangle containing both the random
Hermitization spectra and the supports of the deterministic matrix
spectral measures. On its horizontal edges at imaginary parts
$\pm1$, use the ordinary anisotropic local law. Place the vertical
edges a fixed distance outside both supports. There the deterministic
density is $O(|\Im w|)$, so for $|\Im w|\ge N^{-1/2}$ the ordinary
isotropic error is $O(N^{-1/2+\xi})$. Below that height the fixed
spectral gap bounds the resolvent derivative, and continuation from
height $N^{-1/2}$ costs $O(N^{-1/2})$. A polynomial contour mesh
makes the bounds simultaneous. Cauchy's formula for the two matrix
spectral measures gives the stated $1/u$ decay outside the rectangle;
the remaining bounded interval is covered by the ordinary law.
Since $\|A_i(iu)\|\le C/u$, the assertion for $P_i$ follows.

For the averaged estimate use
\[
 G(iu)=G(i)-i(1-u)G(i)G(iu).
\]
The deterministic two-resolvent solution with identity insertion is
$(M(i)-M(iu))/(i(1-u))$, for $u\ne1$. Subtract the corresponding
averaged identity with $B$ and use the averaged law at height one and
\cite[Lemma B.1]{CL}, whose error is $CN^\xi/(Nu)$.
Here that lemma is used with $\eta_1=1$, $\eta_2=u$, identity
insertion, and its fixed upper-domain exponent at least $10$;
thus its full domain includes $1\le u\le N^{10}$.
Multiplication by $|1-u|$ gives $CN^\xi/N$; continuity covers $u=1$.
\end{proof}

Put $\varepsilon_{\rm out}=N^{-1/2+\xi}$. For $u\ge1$,
\[
 \|G_i(iu)\|_{\HS}\le C\sqrt N/u,
 \qquad \|P_i(u)\|_{\HS}\le C\sqrt N/u^2.
\]
An outer diagonal replacement error is bounded by
$C\varepsilon_{\rm out}/u$ for $G_i$, and by
$C\varepsilon_{\rm out}/u^2$ for $P_i$.
Consider the products with exactly one off-site factor in $(DP_i)DU_r$.
With exactly one outer diagonal error, first contract the free
endpoint index as in \eqref{signed:eq:ss-ordered-contraction}; with
two errors, use the endpoint Hilbert--Schmidt bound. Including
the common $N^{-3/2}$ factor, the resulting integrand is at most
\[
 \frac C{u^3}\left(
  \frac{\varepsilon_{\rm out}}{\sqrt N}
  +\varepsilon_{\rm out}\Psi_r
  +\frac{\varepsilon_{\rm out}^2}{\sqrt{\eta_r}}\right).
\]
The zero-error term is $O(N^{-1/2}u^{-3})$.
For products with three off-site factors, putting the endpoint in the supremum
norm and the two outer factors in Hilbert--Schmidt norm gives
$CN^{-1/2}\Psi_r u^{-3}$. These are integrable, and
\[
 \frac{\varepsilon_{\rm out}^2}{\sqrt{\eta_r}}
 \le C_b N^{-1/2+2\xi}t^{-b/2}.
\]
The ordered contraction preserves the needed saving for the
exactly-one-error term throughout the large-height interval.

For $P_iD^2U_r$, the endpoint diagonal errors cost
$C\Psi_r u^{-2}$, and the products with three off-site factors cost
$C\Psi_r u^{-2}/\sqrt{N\eta_r}$.
For $P_iDU_rDU_s$, the two choices of endpoint supremum give
\[
 \frac C{u^2}\min\left\{
  \frac{\Psi_r}{\sqrt{N\eta_s}},
  \frac{\Psi_s}{\sqrt{N\eta_r}}\right\}.
\]
If either endpoint is below $a_N$, use that endpoint for the
supremum. After integration the bound is
$Ct^2/(N\sqrt{\eta_r\eta_s})\le Ct^2/(N\eta_{\rm mes})$;
if both are above $a_N$, it is $O(N^{-1+11\xi})$.
The third-cumulant terms with both entry derivatives on the outer resolvent expression have $u^{-4}$ decay and
integrated replacement error $O(\varepsilon_{\rm out})$.
At fourth order the pair-count prefactor is bounded; every
replacement error or off-site factor therefore gives
$C(\varepsilon_{\rm out}+\sum_{r\in\mathcal I}\Psi_r)$ after integration,
where $\mathcal I$ contains only the at most three endpoint labels
of that product,
using outer decay at least $u^{-2}$. Equal-site pairs retain
their $O(N^{-1/2})$ and $O(N^{-1})$ savings.
Thus all these tails satisfy
\eqref{signed:eq:ss-third-bound}--\eqref{signed:eq:ss-fourth-bound}, including
unequal fine and coarse endpoints.
Finally, $\|G_i(iu)\|\le u^{-1}$ and $\|A_i(iu)\|\le C/u$
give $u^{-q-2}$ decay for every fixed derivative $D^qP_i$.
The averaged
quadratic tail is at most
\[
 N\int_1^{N^{10}}\frac{CN^\xi}{N}
                  \frac{CN^\xi}{Nu}\,du
 \le CN^{-1+2\xi}\log N .
\]
The third-order errors therefore satisfy
\eqref{signed:eq:ss-third-bound} over the full outer interval; the
fourth-order calculation gives $E_i$, $\kappa m_i^2m_r^2/4$ and
\eqref{signed:eq:ss-fourth-bound}. The averaged quadratic tail is
polynomially small. These are the low-order inputs for the relative
remainder estimate below, with their logarithmic powers unchanged.

\subsection{High orders, caps, and completion of the defect}

In this subsection we control the remaining cumulant terms and cap errors before dividing by the exact normalizing constant. Indeed, the remaining terms must be bounded relative to the exact finite
normalizing constant before the weighted identity is normalized. We first
complete the argument with a fixed, cap-dependent truncation order,
then use a common replacement segment to fix the order at six. Both
arguments use the preliminary lower bound from
\eqref{signed:eq:ss-finite-transfer}; neither uses the normalizer
asymptotic proved in the next section.

\paragraph{Retained higher orders and an initial remainder bound.}
At retained entry derivative order $k$, the cumulant factor is
$O_k(N^{-(k+1)/2})$. Including the trace normalization, the
factor $N$ in the field, and $N^2$ entry pairs gives
$O_k(N^{(3-k)/2})$. For $k\ge4$, bounded good-event endpoint
products and the integrable outer decay make the total
$O_k(N^{-1/2})$ times the exact normalizing constant and the test-function derivative norm.
The number of retained orders is fixed.

For the final cumulant remainder use the absolute cap envelope,
not a factorization of a weighted entry moment.
On the event $|h|\le N^{-2/5}$, a rank-two resolvent identity
shows that all required modified-entry resolvent entries stay
bounded for every replacement $|\theta|\le N^{-2/5}$.
Indeed its finite-dimensional denominator is $I+O(N^{-2/5})$
on the original bounded-entry event. All fixed endpoint
derivative products remain bounded there.
The truncated remainder is therefore bounded by
\begin{equation}
 C_KN^2\E|h|^{K+2}N^{12A}\D_{K+1}(f)
 \le C_KN^{1-K/2+12A}\D_{K+1}(f).
 \label{signed:eq:ss-remainder}
\end{equation}
Choose $K>2(12A+6A_1+D_*+3)$. The preliminary lower bound for the normalizing constant of the finite-integral fields
converts this absolute estimate into
$C N^{-D_*}\widetilde Z(sa)\D_{K+1}(f)$, so the remainder is
relative before division by the normalizing constant. Bad modified-entry events have deterministic
resolvent derivative bounds polynomial in $N$, since
$\eta_{\rm mes}\ge N^{-1}$; choose their failure exponent
after $K,A$. Larger replacements are bounded using sufficiently
high fixed entry moments, also chosen afterwards. No statement
of the form $\E[W_s|h|^{K+2}]=Z(sa)\E|h|^{K+2}$ is used.

\paragraph{The relative defect and removal of cap errors.}
Combining \eqref{signed:eq:ss-quadratic}, \eqref{signed:eq:ss-I6},
\eqref{signed:eq:ss-third-bound}, \eqref{signed:eq:ss-fourth-bound}, and
the higher-order estimates first gives the raw centered identity
\[
 \E[U_i\mathcal H(U)]
 =\sum_r C^\Phi_{ir}\E[\mathcal H_r(U)]
       +\mathcal R_i^\Phi,\qquad
 C^\Phi_{ir}=-\tfrac14\log D_{ir}
                 +\tfrac\kappa4m_i^2m_r^2 .
\]
Multiply by $\sqrt2$ and use
$\mathcal H_r(U)=\sqrt2H_r(\sqrt2U)$.
This gives the identity for $\widetilde Y$, with $\Sigma=2C^\Phi$:
\[
 \E[\widetilde Y_i H(\widetilde Y)]
 =\sum_r\Sigma_{ir}\E[H_r(\widetilde Y)]+\mathcal R_i,
 \qquad
 |\mathcal R_i|\le r_N\widetilde Z(sa)\D_{m-1}(f).
\]
The mean \eqref{signed:eq:ss-mean} is already subtracted; it is not
inserted a second time. Each coefficient estimate is pathwise
on a common good event, before applying
\eqref{signed:eq:ss-envelope}. The complement is bounded using the same preliminary lower bound for the normalizing constant and an
arbitrarily large failure exponent. Thus the error displayed above
is already proportional to the exact normalizing constant.

The difference $\widetilde Y_i-c_A(\widetilde Y_i)$ is supported
on the cap shell and is deterministically $O(N\log N)$.
Derivatives of $c_A$ different from those of the identity are
also shell-supported. Equation \eqref{signed:eq:ss-shell-budget}
controls these errors, including covariance factors $O(t)$.
Off the shell,
\[
 H_r=e^{sa\cdot y}(\partial_rf+sa_rf).
\]
The bound on the finite cap shell is relative to the preliminary
lower estimate for the normalizing constant, so these terms can be added to the same error.
Only now divide by the exact $\widetilde Z(sa)$ to obtain the
finite-field scalar weighted identity under $\widetilde Q_{sa}$.
Summing it with $f=F_i$ uses the norm
\eqref{signed:eq:ss-norm}, with no extra dimension factor.

\paragraph{Return to the actual fields.}
Restore the actual fields using
\eqref{signed:eq:ss-finite-transfer}. For $\|x\|_\infty\le2A\log N$, the scalar defect
integrand
\[
 d_i(x)=(x_i-s(\Sigma a)_i)f(x)
                  -\sum_r\Sigma_{ir}\partial_rf(x)
\]
has supremum at most $Ct\D_1(f)$ and summed first-derivative
supremum at most $Ct\D_2(f)$, because
$\max_{ir}|\Sigma_{ir}|\le Ct$ and $\|a\|_1\le6$.
The good-event coordinate displacement is $O(N^{-5})$.
The relative weight estimate and the bound on the normalized contribution of the complement
therefore transfer each defect expectation with error at most
$Ct(N^{-5}+N^{-L})\D_2(f)$. This is absorbed into $r_N\D_{m-1}(f)$
after fixing $L$ large and decreasing $c$ if needed; here $m-1\ge2$.
Summing over outputs introduces only the existing
$\|\mathbf F\|_{m-1,1}$. Thus every final expectation and
normalizer is the actual one, without the lossy absolute error
$N^{12A}N^{-5}$.
This proves Theorem~\ref{signed:thm:ss-defect}.

The order of choices is explicit: first the fixed bulk, common moment bounds
and coefficient bound; then $A_1$, then $A$, then $K$ and $m$;
then all failure and fixed-moment orders; finally $b$.
The one-resolvent and integrated insertion estimates have fixed logarithmic exponents.
Thus \eqref{signed:eq:ss-rate} has exponents independent of $b$.

\subsubsection{A fixed sixth-order expansion}

Next we keep the expansion order fixed by showing that changing one small entry barely changes the exponential weight.
The initial remainder argument allowed the derivative order to
depend on the cap size. We now replace that estimate alone.
The exponential weight varies by a bounded relative factor along one common
entry-replacement disk, which permits truncation at derivative
order six and gives $m=8$.

\begin{lemma}[Stability of the exponential weight under entry replacement]
\label{signed:lem:ss-segment}
Let $\rho_N=N^{-2/5}$ and include
$\max_{ab}|X_{ab}|\le\rho_N$ in the common overwhelming event.
For one independent complex entry, replace its value $h$ by any
$\theta$ with $|\theta|\le\rho_N$, keeping the Hermitian paired entry.
For the finite-integral fields $\widetilde{\mathbf Y}$, uniformly on this
whole replacement disk,
\[
 \max_i|\widetilde Y_i(X^{ab,\theta})-\widetilde Y_i(X)|
 \le C\rho_N,\qquad
 e^{-C\rho_N}\le
 \frac{W_{sa}(X^{ab,\theta})}{W_{sa}(X)}
 \le e^{C\rho_N},
\]
where $W_{sa}=\exp(sa\cdot c_A(\widetilde{\mathbf Y}))$.
All fixed-order endpoint entry derivatives remain bounded.
\end{lemma}
\begin{proof}
The rank-two Woodbury inverse is a $2$ by $2$ matrix
$I+O(\rho_N)$, since the four original entries it uses are bounded.
It bounds all modified resolvent entries uniformly on the disk.
The endpoint derivative formula \eqref{signed:eq:ss-endpoint}, with the
fixed $\sqrt2$ factor, gives a bounded gradient of each finite
field along this disk; integrate from $h$ to $\theta$.
The cap is Lipschitz and $\|sa\|_1\le6$, so the logarithm of the
weight changes by at most $C\rho_N$. This also holds on cap shells.
\end{proof}

\begin{lemma}[Relative sixth-order cumulant remainder]
\label{signed:lem:ss-sixth-remainder}
Apply the entry cumulant formula through derivative order six to the
entire per-entry scalar function, summed over outputs,
\[
 \mathcal T_{ab}(X)=
 \sum_i\int_{\eta_i}^{N^{10}}(P_i(u))_{ba}
       W_{sa}(X)F_i(c_A(\widetilde{\mathbf Y}(X)))\,du,
\]
including both outer orientations and their fixed raw field factors.
The total remainder has absolute value at most
\[
 C(N^{-6/5}+N^{-10})\,\widetilde Z(sa)\|\mathbf F\|_{7,1}.
\]
The exponent ten can be replaced by any prescribed fixed exponent
by increasing the failure exponent and the required entry-moment orders.
\end{lemma}
\begin{proof}
Condition on the other entries. Taylor-expand the entire scalar
function about the zero entry through degree six. The conditioning
justifies the entry moment--cumulant identity; it does not separate
the exponential weight from the entry being expanded. Expand each
derivative of order $k\le6$ on the right-hand side of the cumulant
identity through degree $6-k$. Polynomial terms cancel by the
moment--cumulant identity. Suppressing finitely many real directional
multiplicities, the remainder is bounded by
\[
 C\E[|h|^8\mathcal S_7\mathcal T_{ab}(h)]
 +C\sum_{k=1}^6N^{-(k+1)/2}
       \E[|h|^{7-k}\mathcal S_7\mathcal T_{ab}(h)],
\]
where $\mathcal S_7$ is the supremum of the sum of absolute seventh
directional derivatives along the segment from zero to $h$.
Every field in a derivative is evaluated at the same modified matrix,
with the Hermitian paired entry changed simultaneously. The common
disk in Lemma~\ref{signed:lem:ss-segment} contains this whole segment.

On the common event, Lemma~\ref{signed:lem:ss-segment}, the chain rule,
bounded endpoint entries, and integrable outer factors give
\[
 \mathcal S_7\mathcal T_{ab}(h)
 \le C W_{sa}(X)\|\mathbf F\|_{7,1}.
\]
The cap derivatives are uniformly bounded, including on shells.
Since $|h|\le\rho_N$ and $N^{-1/2}\le\rho_N$, every entry factor
preceding this derivative bound is at most $C\rho_N^8$.
The physical pair sum is therefore at most
\[
 CN^2\rho_N^8\E W_{sa}\|\mathbf F\|_{7,1}
 =CN^{-6/5}\widetilde Z(sa)\|\mathbf F\|_{7,1}.
\]
This is a deterministic truncation inside the weighted expectation.
The comparison of exponential weights along the segment leaves the original
$W_{sa}(X)$ inside that expectation, whose mean is precisely
$\widetilde Z(sa)$. It does not factor a weighted entry moment.

On the complement, at every replacement point the cap bounds the
exponential weight by $N^{12A}$, and the deterministic resolvent estimate gives
$\mathcal S_7\mathcal T_{ab}\le
 CN^{12A+7}(1+\log N)\|\mathbf F\|_{7,1}$.
Cauchy--Schwarz and sufficiently high fixed entry moments bound all
exceptional terms after the pair sum by
$CN^{12A+14-D/2}\|\mathbf F\|_{7,1}$, when the common event fails with
probability at most $N^{-D}$. Taking
$D>2(12A+6A_1+60)$ and applying the preliminary lower bound for the normalizing constant
only to this exceptional part makes it smaller than
$CN^{-10}\widetilde Z(sa)\|\mathbf F\|_{7,1}$.
The event includes large entries; no derivative of its indicator is
taken. All moment and failure choices are fixed before $N\to\infty$.
Finally \eqref{signed:eq:ss-finite-transfer} compares the normalizing constant for the finite-integral fields
with the actual normalizing constant, with a relative error.
\end{proof}

\begin{corollary}
\label{signed:cor:ss-eight}
Theorem~\ref{signed:thm:ss-defect} and all its normalization and Gaussian
comparison consequences hold with $m=8$, with the same rate
\eqref{signed:eq:ss-rate}, including its proved $t^{17}$ insertion term
and arbitrary finite coarse family.
\end{corollary}
\begin{proof}
Retain entry derivative orders four, five and six, which cost
$O(N^{-1/2})$ by the bounds already proved. Replace only the last
absolute-envelope estimate by Lemma~\ref{signed:lem:ss-sixth-remainder}.
The low-order coefficients and cap-shell bounds are unchanged.
The extra $N^{-6/5}$ is absorbed into $N^{-c}$.
Thus seven vector derivatives suffice; the OU argument below
uses eight scalar derivatives. The order of choices is now the
bulk domain and coefficient bounds, $A_1,A$, the fixed order six, the
failure exponent and required moment orders, and finally the logarithmic exponent $b$.
\end{proof}

\subsection{Normalizing constants and Gaussian comparison}

This section uses the weighted Stein identity to estimate normalizing constants and compare smooth tests with the corresponding Gaussian law.

The defect now holds under the actual weighted law, with its error
controlled independently of any Gaussian formula for the normalizing
constant. We use constant vector tests to integrate along each permitted
coefficient ray, then gradient tests to compare with the Gaussian law.

The Gaussian conversion uses the standard Ornstein--Uhlenbeck solution
of the Stein equation; see \cite[Lemmas~1(3) and~2(1), and (10)]{Meckes}.
We record the derivative contraction in the present summed norm;
positive semidefinite, possibly singular covariance is allowed.

\begin{corollary}[Normalizing constants for sparse exponential weights]
\label{signed:cor:ss-normalizer}
Under Theorem~\ref{signed:thm:ss-defect},
\[
 \log Z(a)=\tfrac12a^{\mathsf T}\Sigma a+O(r_N)
\]
uniformly for the stated sparse coefficient vectors.
\end{corollary}
\begin{proof}
All coordinates are capped, so differentiation under expectation
is legitimate. For a fixed permitted $a$, every $sa$ with $0\le s\le1$
has support contained in that of $a$ and $\|sa\|_1\le6$, hence remains in the
theorem's scope. Apply \eqref{signed:eq:ss-defect} with the constant
vector $\mathbf F=a$. Its norm is $\|a\|_1\le6$, since the
order-zero term in \eqref{signed:eq:ss-norm} is retained:
\[
 \frac d{ds}\log Z(sa)
 =\E_{Q_{sa}}(a\cdot\mathbf X)
 =s\,a^{\mathsf T}\Sigma a+O(r_N).
\]
Integrate from $s=0$, where $Z(0)=1$.
\end{proof}

\begin{corollary}[Smooth comparison, including singular covariance]
\label{signed:cor:ss-ou}
If $\mathbf G_a\sim N(\Sigma a,\Sigma)$ and $g$ is bounded and
$C^m$, then
\[
 |\E_{Q_a}g(\mathbf X)-\E g(\mathbf G_a)|
       \le Cr_N\D_m(g).
\]
\end{corollary}
\begin{proof}
Use the Gaussian Ornstein--Uhlenbeck semigroup
\[
 P_ug(x)=\E g\bigl(\Sigma a+e^{-u}(x-\Sigma a)
       +\sqrt{1-e^{-2u}}(\mathbf G_a-\Sigma a)\bigr)
\]
and $h_L=-\int_0^L(P_ug-\E g(\mathbf G_a))\,du$.
Its generator
$\mathcal Lh=\sum_{ir}\Sigma_{ir}\partial_{ir}h
 -(x-\Sigma a)\cdot\nabla h$
satisfies $\mathcal Lh_L=g-P_Lg$.
For every $q\ge1$,
\[
 \sum_{r_1,\ldots,r_q}
 \|\partial_{r_1}\cdots\partial_{r_q}h_L\|_\infty
 \le\frac1q\sum_{r_1,\ldots,r_q}
 \|\partial_{r_1}\cdots\partial_{r_q}g\|_\infty.
\]
In particular
$\sum_i\D_{m-1}(\partial_i h_L)\le\D_m(g)$.
Apply Theorem~\ref{signed:thm:ss-defect} with $\mathbf F=\nabla h_L$.
Boundedness of $g$ allows $L\to\infty$ by dominated convergence.
The argument uses the Gaussian vector directly and is valid for a
positive semidefinite singular $\Sigma$. No inverse covariance or
lower eigenvalue bound is used.
\end{proof}

Taking $a=0$ gives $Z(0)=1$ and $Q_0=\PP$, so the same comparison
applies to the capped field under the unweighted matrix law. For bounded tests,
\eqref{signed:eq:ss-tail} also removes the caps at the
chosen polynomial accuracy. The zero-charge consequence therefore
remains available without a separate normalization argument.

\paragraph{Exact cancellation of exponential factors for history tests.}
For later use let $p=2e_{r_z}$, $q=2e_{r_w}$ and
$a_j=p+q-2e_c$. Each coefficient vector has at most three nonzero
coordinates and $\ell^1$ norm at most six, including any coordinate
coincidences. Applying Corollary~\ref{signed:cor:ss-normalizer}
separately to $p$, $q$, and $a_j$ gives
\[
 \log\frac{e^{2j}Z(a_j)}{Z(p)Z(q)}
 =2j+4\Sigma_{r_zr_w}-4\Sigma_{r_zc}
       -4\Sigma_{r_wc}+2\Sigma_{cc}+O(r_N).
\]
If the separate exact-clock covariance calculation makes all four
entries $j+O(1)$, their leading terms cancel and this ratio stays
bounded above and below. The covariance geometry is not replaced
by a signed moment assumption.

For this algebraic application, fix a smooth nonincreasing
$0\le\theta\le1$ that equals one on $(-\infty,0]$ and zero
on $[1,\infty)$. Let $0\le H_0\le1$ be a smooth function of
a fixed finite coarse list, with bounded derivatives through
order $m$. For an observed list $\mathcal I$ of $O(t)$ history
coordinates, including $c=(j,z)$, set
\[
 \mathcal B(x)=H_0(x_{\rm coarse})
       \prod_{(k,v)\in\mathcal I}
          \theta(x_{k,v}-2k-B),\qquad
 g_j(x)=e^{2(x_c-2j)}\mathcal B(x).
\]
The transition width is one. The factor at $c$ gives
$0\le g_j\le e^{2(B+1)}$ and
$\D_m(g_j)\le C_{B,H_0,m}(1+t)^m$.
The exponential factors cancel exactly:
\[
 \frac{e^{2j}Z(a_j)}{Z(p)Z(q)}
       \E_{Q_{a_j}}g_j(\mathbf X)
 =e^{-2j}
   \frac{\E[e^{2X_{r_z}+2X_{r_w}}\mathcal B(\mathbf X)]}
        {Z(p)Z(q)}.
\]
The identical capped negative coordinate occurs on both sides
of this cancellation. It remains valid when that cap is active.
Require every full transition interval $[2k+B,2k+B+1]$ to lie
inside $[-A\log N,A\log N]$. For the intended histories
$1\le k\le t+O(1)$, fixed $B$ and $A>2$, this holds for all
large $N$. Monotonicity and the identity region of $c_A$ then
give equality of the capped and uncapped barrier factors
pointwise, including on both saturated sides.
Caps on bounded smooth coarse tests can instead be removed using
their tail bound.
Positive-root caps may be removed after cancellation using
positive moments and scalar tails. The coordinate with the negative coefficient
remains capped under the tilted law; no inverse determinant moment is invoked.
The displayed cancellation is an algebraic application; it
does not establish a barrier probability estimate or a
comparison of normalizing constants for different ensembles.

\endgroup

\begingroup
\newcommand{\E}{\mathbb E}
\newcommand{\Pp}{\mathbb P}
\newcommand{\C}{\mathbb C}
\newcommand{\R}{\mathbb R}
\newcommand{\1}{\mathbf 1}
\newcommand{\D}{\mathcal D}
\newcommand{\cX}{\mathsf X}
\newcommand{\cY}{\mathsf Y}
\newcommand{\cE}{\mathcal E}
\newcommand{\norm}[1]{\left\lVert #1\right\rVert}
\newcommand{\abs}[1]{\left\lvert #1\right\rvert}
\newcommand{\Tr}{\operatorname{Tr}}
\newcommand{\supp}{\operatorname{supp}}
\section{Entry replacement for exponential weights}\label{chap:replacement}

In this section we change the entry law under the stated moment conditions while controlling weighted expectations relative to their exact normalizing constants.
We compare capped exponential weights between complex iid matrix
ensembles. The coefficients in these weights may have either sign.
The normalized mixed entry moments agree through degree three and differ
by \(O(N^{-1/200})\) at degree four. We first compare the unnormalized
weighted expectations relative to their exact normalizing constants,
then divide to compare the tilted laws. An event determined by the
matrix with the chosen entry set to zero allows us to bound the
fifth-moment remainder relative to the normalizing constant. We retain
the complete fourth-order Taylor polynomial and control its large-entry
tails separately. We also construct explicitly an admissible matched
Gaussian-divisible endpoint with initial Gaussian age exactly \(N^{-1/200}\).
The later comparison across regularization scales uses these separate
ratios of normalizing constants together with local-factor reconstruction.

This argument adapts the Green-function comparison approach in
\cite[Propositions~9.6 and~9.11]{CL}. The exact normalizing constants are
retained through the relative replacement estimate, so the conclusion
also controls the normalized laws for the permitted capped weights.

\subsection{Fields, common moment hypotheses, and the conclusion}

In this subsection we state the moment and regularization conditions under which capped exponential weights can be compared between two ensembles. We first fix a closed disk \(K=\{z:\abs z\le r\}\), where \(r<1\).
For an \(N\times N\) matrix \(M\), let
\[
 H^z(M)=
 \begin{pmatrix}0&M-z\\(M-z)^*&0\end{pmatrix},
 \qquad G^z(w;M)=(H^z(M)-w)^{-1}.
\]
Let \(\rho_z\) denote the deterministic variance-one Hermitized density and
write
\[
 I_z(\eta)=\int_{\R}\log(x^2+\eta^2)\rho_z(x)\,dx .
\]
Our real regularized field is
\begin{equation}\label{replacement:eq:rgfc-field}
 \Phi_N(z,\eta;M)
 =\frac12\log\det\big((M-z)(M-z)^*+\eta^2\big)
       -\frac N2 I_z(\eta),\qquad \eta>0.
\end{equation}
This is the convention of \cite[(2.13)]{CL}; in particular at height zero,
when the determinant is nonzero, it is
\(\log\abs{\det(M-z)}-\frac N2(\abs z^2-1)\).
No height-zero determinant occurs in the comparison theorem below.

The entry laws may depend on \(N\).
Write \(x=\xi/\sqrt N\), \(y=\upsilon/\sqrt N\) for their actual
matrix-entry normalizations.  We assume
\begin{equation}\label{replacement:eq:rgfc-profile}
 \E\xi=\E\upsilon=0,\qquad
 \E\abs{\xi}^2=\E\abs{\upsilon}^2=1,\qquad
 \E\xi^2=\E\upsilon^2=0,\qquad
 \sup_N\big(\E\abs{\xi}^k+\E\abs{\upsilon}^k\big)\le C_k
\end{equation}
for every fixed positive integer \(k\).
The bounds \((C_k)\) are common to both entry laws and uniform in \(N\).
Independence of the real and imaginary parts is not assumed.
Set
\begin{equation}\label{replacement:eq:rgfc-matching}
 \begin{aligned}
 \sigma&=N^{-1/200},\\
 \E x^p\overline{x}^{\,q}&=\E y^p\overline{y}^{\,q}
                                  &&(p+q\le3),\\
 \abs{\E x^p\overline{x}^{\,q}-\E y^p\overline{y}^{\,q}}
       &\le C\sigma N^{-2}         &&(p+q=4).
 \end{aligned}
\end{equation}
All mixed complex moments in these formulas are required.  In particular,
matching only the holomorphic third moment would be insufficient.
Enlarge \(C_4\) if necessary so that it also bounds the fixed constant in the
fourth-moment discrepancy.

Fix
\[
 0<\delta<\frac1{1200},\qquad n=n_N\le C_{\rm list}\log N,
 \qquad z_j\in K,\quad N^{-1-\delta}\le\eta_j\le1 .
\]
The list and real deterministic centers \(e_j\) may depend on \(N\), but
are fixed during the comparison.  Define
\[
 \cY_j(M)=\sqrt2\,\big[\Phi_N(z_j,\eta_j;M)-e_j\big].
\]
Let \(A>0\) be fixed.  Choose a smooth cap \(c_A=c_{A,N}\) that is odd,
nondecreasing, and \(1\)-Lipschitz, equals the identity on
\([-A\log N,A\log N]\), and satisfies
\begin{equation}\label{replacement:eq:rgfc-cap}
 \abs{c_A(u)}\le2A\log N,\qquad
 \sup_{N\ge3}\norm{c_A^{(k)}}_\infty\le C_{A,k},\quad 1\le k\le5 .
\end{equation}
Such caps exist by integrating a smooth nonnegative cutoff of the
derivative.  Put \(\cX_j(M)=c_A(\cY_j(M))\).

For \(a\in\R^n\) satisfying
\begin{equation}\label{replacement:eq:rgfc-charge}
 \abs{\supp a}\le3,\qquad \norm a_1\le6,
\end{equation}
define the positive weight and its ensemble normalizer by
\[
 W_a(M)=\exp(a\cdot\cX(M)),\qquad
 Z_{a,L}=\E_L W_a(M),\qquad
 \E_{Q_{a,L}}f=\frac{\E_L[W_a(M)f(\cX(M))]}{Z_{a,L}} .
\]
Here \(L\) specifies the matrix law.  For deterministic \(M\), \(W_a(M)\)
is a value, not an expectation.  Uniformly in all matrices and the
coefficient family,
\begin{equation}\label{replacement:eq:rgfc-deterministic-cap}
 N^{-C_A}\le W_a(M)\le N^{C_A},\qquad C_A=12A .
\end{equation}
In particular, the weight remains positive when some coefficients are negative.

For \(f\in C_b^k(\R^n)\), use the derivative norm
\begin{equation}\label{replacement:eq:rgfc-testnorm}
 \D_k(f)=\sum_{\ell=0}^k
   \sum_{j_1,\ldots,j_\ell=1}^n
       \norm{\partial_{j_1}\cdots\partial_{j_\ell}f}_\infty .
\end{equation}
The term \(\ell=0\) is \(\norm f_\infty\), which controls the
exceptional-event and quotient terms. The sum of derivative suprema
already accounts for every coordinate choice in the growing list.

\begin{theorem}[Relative entry replacement under exponential weights]\label{replacement:thm:rgfc-main}
Let \(X,Y\) be iid matrices satisfying
\eqref{replacement:eq:rgfc-profile}--\eqref{replacement:eq:rgfc-matching}.
With the deterministic list, fixed centers, cap, and source family above,
put
\[
 \epsilon_N=N^{6\delta}(\sigma+N^{-1/2}),\qquad
 c_0=\frac1{200}-6\delta>0 .
\]
Then, uniformly over this family,
\begin{equation}\label{replacement:eq:rgfc-main}
 \frac{Z_{a,Y}}{Z_{a,X}}=1+O(\epsilon_N),\qquad
 \abs{\E_{Q_{a,Y}}f-\E_{Q_{a,X}}f}
       \le C\epsilon_N\D_5(f),\qquad f\in C_b^5(\R^n).
\end{equation}
In particular \(\epsilon_N\le2N^{-c_0}\).
The same conclusions hold uniformly along each coefficient path \(sa\),
\(0\le s\le1\).
The constants depend only on the fixed bulk disk, the cap, the bounds
on the coefficient vector and list length, \(\delta\), and the common
moment bounds. They do not depend on spatial separation.
\end{theorem}

The theorem compares the stated positive-height exponential weights
with capped coordinates. It does not remove caps from coordinates with
negative coefficients or compare height-zero fields.
The relative differential inequality controls the change in the
normalizing constant before any quotient is taken. Consequently the
first assertion remains informative when that constant is polynomially small.

\subsection{Derivative estimates and an event determined by the remaining entries}

In this subsection we obtain derivative bounds on an event independent of the replaced entry, making the Taylor remainder compatible with weighted expectations.
We use \cite[Proposition 9.1, (9.5)--(9.7), and Lemma 9.4,
(9.18)--(9.19)]{CL}, in version 3.  Their height domain is exactly
\([N^{-1-\delta},1]\), which is why we have imposed this upper endpoint.
We use the magnitude estimates in their actual-entry convention. They
apply to an ordinary
iid matrix from a common moment class and uniformly over replacements
of one entry by a complex number of modulus at most \(N^{-2/5}\).

\begin{lemma}[Derivative input from Cipolloni--Landon]\label{replacement:lem:rgfc-input}
Let \(U\) be any iid matrix satisfying \eqref{replacement:eq:rgfc-profile}.
For an entry position \(e=(a,b)\), let \(U(\zeta;e)\) be \(U\) with
that entry replaced by \(\zeta\), and set \(h_N=N^{-2/5}\).
For every fixed \(\varepsilon>0\) in the derivative bound and failure exponent
\(L>0\), with probability at least \(1-C_LN^{-L}\),
\begin{equation}\label{replacement:eq:rgfc-input}
 \sup_{\abs\zeta\le h_N}
 \abs{\partial_\zeta^p\partial_{\bar\zeta}^{\,q}
              \cY_j(U(\zeta;e))}
       \le C_{p+q}N^{(p+q)\delta+\varepsilon},
 \quad 1\le p+q\le5,\quad 1\le j\le n .
\end{equation}
The constants are uniform over iid laws with the common moment bounds,
entry positions, and admissible deterministic lists.
\end{lemma}

\begin{proof}
The uniform entrywise resolvent bound in \cite[(9.5)]{CL} is
\(N^{\delta+\varepsilon'}\).
For completeness the relevant differentiation has fixed finite rank.
With \(J_e=E_{a,b+N}\) in the Hermitized space,
\[
 \partial_\zeta H^z=J_e,\quad
 \partial_{\bar\zeta}H^z=J_e^*,\quad
 \partial_\zeta G=-GJ_eG,\quad
 \partial_{\bar\zeta}G=-GJ_e^*G .
\]
The convention \eqref{replacement:eq:rgfc-field} gives
\[
 \partial_\zeta\Phi_N
       =\tfrac12 G_{b+N,a},\qquad
 \partial_{\bar\zeta}\Phi_N
       =\tfrac12 G_{a,b+N}.
\]
One may verify these directly from the positive determinant in
\eqref{replacement:eq:rgfc-field}.  Repeated differentiation yields finitely many
products of \(p+q\) resolvent entries.  Thus their absolute values are
bounded by \(C_{p+q}N^{(p+q)(\delta+\varepsilon')}\).
Take \(\varepsilon'\) sufficiently small.  The deterministic centers
differentiate to zero; the factor \(\sqrt2\) is harmless.
This also derives precisely the part of \cite[(9.18)]{CL} needed here,
without using any printed sign convention for the resolvent derivative.
The \(O(\log N)\) coordinate union is absorbed by increasing the
failure exponent in the input estimate.

Uniformity under these moment bounds is part of the moment dependence
in the input estimates.  Equivalently, if uniformity failed, a violating
sequence of laws, points, and entry positions would itself be an
admissible triangular array to which the same estimates apply.
\end{proof}

Let \(\mu_X,\mu_Y\) be the laws of the actual entries \(x,y\), and define
\begin{equation}\label{replacement:eq:rgfc-mixture}
 \mu_\theta=(1-\theta)\mu_X+\theta\mu_Y,\qquad 0\le\theta\le1 .
\end{equation}
Every matrix with iid entries from \(\mu_\theta\) satisfies the same
moment bounds and the same first two moments.  We apply the preceding input
separately at each deterministic \(\theta\).  No event simultaneous in
an uncountable \(\theta\)-family is needed.

Fix the exponent in the final derivative bound with \(0<\varepsilon<\delta\),
decreasing the exponent in the input bound first for the finitely many chain-rule
products. The event used for Taylor expansion must be independent of
the entry that will be integrated. We therefore define it from the matrix
with that entry set to zero, using the whole replacement disk.
For a fixed \(e,\theta\), let \(V=V_\theta^e\) have entry \(e\) equal to
zero and all its other entries independent with law \(\mu_\theta\).
Define \(\cE_e(V)\) to be the property that every inequality
\eqref{replacement:eq:rgfc-input} holds for \(V+\zeta E_e\), over the entire closed
disk \(\abs\zeta\le h_N\).
Here \(E_e\) is the \(N\times N\) matrix unit.
The event is measurable in \(V\): all functions are continuous for
positive height, so the suprema can equally be taken over a countable
dense subset of the disk.

\begin{lemma}[Setting one entry to zero]\label{replacement:lem:rgfc-deleted}
Uniformly in \(\theta,e\),
\begin{equation}\label{replacement:eq:rgfc-deleted-prob}
 \Pp_V(\cE_e^c)\le C_LN^{-L}
\end{equation}
for any prescribed \(L\), after fixing the input constants.
On \(\cE_e\), for all sufficiently large \(N\),
\begin{equation}\label{replacement:eq:rgfc-source-local}
 \tfrac12 W_a(V)\le W_a(V+\zeta E_e)\le2W_a(V),
                         \qquad \abs\zeta\le h_N .
\end{equation}
If \(U_\theta\) is iid with law \(\mu_\theta\) and
\(Z(\theta)=\E W_a(U_\theta)\), then
\begin{equation}\label{replacement:eq:rgfc-deleted-relative}
 \E_V[W_a(V)\1_{\cE_e}]\le4Z(\theta).
\end{equation}
\end{lemma}

\begin{proof}
Adjoin an independent \(u\sim\mu_\theta\).
The completed matrix \(V+uE_e\) is iid.
Its event from Lemma~\ref{replacement:lem:rgfc-input} implies the property
\(\cE_e(V)\), since the conclusion already holds at every replacement
value, including zero.  Hence
\(\Pp_V(\cE_e^c)=\Pp_{V,u}(\cE_e^c)\le C_LN^{-L}\).
The iid input is applied to the completed matrix, while the resulting
property depends only on the matrix \(V\) with that entry set to zero.

By the cap chain rule and the charge bound, on \(\cE_e\),
\begin{equation}\label{replacement:eq:rgfc-log-derivative}
 \abs{\partial_\zeta^p\partial_{\bar\zeta}^{\,q}
                    \log W_a(V+\zeta E_e)}
       \le C_{p+q}N^{(p+q)\delta+\varepsilon},
                   \qquad 1\le p+q\le5,
\end{equation}
where the exponent in the input bound was decreased first to absorb the finitely many
products in this formula.
Integration of the first derivatives along the straight segment from
zero to \(\zeta\) gives
\[
 \abs{\log W_a(V+\zeta E_e)-\log W_a(V)}
       \le C N^{-2/5+\delta+\varepsilon}=o(1).
\]
This proves \eqref{replacement:eq:rgfc-source-local}.
The assertion holds with either sign of a coefficient and when
the cap is active.

Finally, positivity and independence of \(u,V\) imply
\[
 \begin{split}
 Z(\theta)
 &=\E_{V,u}W_a(V+uE_e)\\
 &\ge \tfrac12\,\E_V[W_a(V)\1_{\cE_e}]
                         \Pp(\abs u\le h_N)
 \ge\tfrac14\,\E_V[W_a(V)\1_{\cE_e}].
 \end{split}
\]
The last inequality follows from the common moment bound, which gives
\(\Pp(\abs u>h_N)\to0\) uniformly in \(\theta\).
The lower bound uses a nonnegative integral over \(\cE_e\), so
\eqref{replacement:eq:rgfc-deleted-relative} has no additive exceptional
error. This estimate retains the exact normalizing constant in
the Taylor remainder.
\end{proof}

\subsection{Changing the entry law and estimating the Taylor remainder}

This subsection uses moment matching and tail bounds to control how weighted expectations change along a mixture of entry laws. For this purpose we first fix \(a,f\) and write
\[
 F_f(M)=W_a(M)f(\cX(M)),\qquad
 g_f(\theta)=\E F_f(U_\theta),\qquad Z(\theta)=g_1(\theta).
\]
The integrand is bounded in absolute value by
\(N^{C_A}\norm f_\infty\).
The finite product measure
\(\mu_\theta^{\otimes N^2}\) is polynomial in \(\theta\); consequently
differentiation of this expectation is exact.  With independent
\(x\sim\mu_X\), \(y\sim\mu_Y\), both independent of \(V_\theta^e\),
\begin{equation}\label{replacement:eq:rgfc-exact-derivative}
 g_f'(\theta)=\sum_{e=1}^{N^2}
       \E\big[F_f(V_\theta^e+yE_e)-F_f(V_\theta^e+xE_e)\big].
\end{equation}
Thus each derivative summand has an iid mixture background before
deleting one entry, within the scope of the preceding input.

\begin{lemma}\label{replacement:lem:rgfc-chain}
For \(0\le k\le5\), on \(\cE_e\),
\begin{equation}\label{replacement:eq:rgfc-chain}
 \sup_{\abs\zeta\le h_N}
 \abs{\partial_\zeta^p\partial_{\bar\zeta}^{\,q}
               F_f(V+\zeta E_e)}
       \le C_kN^{k\delta+\varepsilon}\D_k(f)W_a(V),
             \qquad p+q=k .
\end{equation}
For \(k=0\) the factor \(N^\varepsilon\) may be omitted.
\end{lemma}

\begin{proof}
Apply the product rule to the exponential weight and the test.
Each term in the derivative of \(f(\cX)\) consists of a derivative
\(\partial_{j_1}\cdots\partial_{j_\ell}f\), summed over its coordinate
indices, and a product of field derivatives whose total entry order
is \(k\).  Each cap derivative is bounded.
Their powers are therefore at most \(N^{k\delta+\varepsilon}\) after
choosing the input slack sufficiently small.
For the weight use \eqref{replacement:eq:rgfc-log-derivative}, then
\eqref{replacement:eq:rgfc-source-local}.  The number of partitions of an integer
at most five is fixed.  All coordinate sums are already contained in
\(\D_k(f)\), so no extra \(n^k\) appears.
\end{proof}

\begin{proposition}[Relative differential inequality]\label{replacement:prop:rgfc-R1}
For every prescribed \(D>0\), the input failure and entry moment orders
can be fixed so that, uniformly in \(0\le\theta\le1\),
\begin{equation}\label{replacement:eq:rgfc-R1}
 \abs{g_f'(\theta)}
       \le C\epsilon_N Z(\theta)\D_5(f)+N^{-D}\D_5(f).
\end{equation}
These fixed orders are chosen after the cap and \(D\), and do not depend
on later fixed logarithmic cutoff exponents.
\end{proposition}

\begin{proof}
Fix an entry position and abbreviate \(V=V_\theta^e\).
On \(\cE_e\), use the complete polynomial
\begin{equation}\label{replacement:eq:rgfc-taylor}
 T_4(v;V)=\sum_{p+q\le4}\frac{v^p\bar v^{\,q}}{p!\,q!}
       \partial_\zeta^p\partial_{\bar\zeta}^{\,q}F_f(V).
\end{equation}
The real two-dimensional Taylor formula, equivalently expressed with
Wirtinger derivatives, and Lemma~\ref{replacement:lem:rgfc-chain} give
\begin{equation}\label{replacement:eq:rgfc-small-remainder}
 \abs{F_f(V+vE_e)-T_4(v;V)}
    \le C\abs v^5N^{5\delta+\varepsilon}\D_5(f)W_a(V),
                    \qquad \abs v\le h_N .
\end{equation}

\paragraph{The complete polynomial and small-entry remainder.}
First integrate the full polynomial \eqref{replacement:eq:rgfc-taylor}, with no
truncation of the entry law.  Because its coefficients and \(\cE_e\)
are measurable in \(V\), all degrees through three cancel exactly
between \(x\) and \(y\).  The fourth-degree part is bounded by
\begin{equation}\label{replacement:eq:rgfc-degree4}
 \abs{\E[\1_{\cE_e}(T_4(y;V)-T_4(x;V))]}
 \le C\sigma N^{-2+4\delta+\varepsilon}\D_4(f)
                      \E[W_a(V)\1_{\cE_e}] .
\end{equation}
For the small-entry remainders, independence gives
\[
 \E[W_a(V)\1_{\cE_e}(\abs x^5+\abs y^5)]
   =\E[W_a(V)\1_{\cE_e}]\,\E(\abs x^5+\abs y^5)
   \le CN^{-5/2}\E[W_a(V)\1_{\cE_e}].
\]
Inserting small-entry indicators only decreases this upper bound.
Thus their total contribution is at most
\begin{equation}\label{replacement:eq:rgfc-degree5}
 C N^{-5/2+5\delta+\varepsilon}\D_5(f)
                               \E[W_a(V)\1_{\cE_e}] .
\end{equation}
We have not factored an exponential weight evaluated at \(V+xE_e\) from a function
of \(x\).

\paragraph{Large entries and the event \(\cE_e^c\).}
The use of full moments also requires the tails of every Taylor
monomial, as well as the tail of the original observable.  On \(\cE_e\), for arbitrary \(v\),
\[
 \abs{F_f(V+vE_e)-T_4(v;V)}
 \le N^{C_A}\norm f_\infty
    +C N^{C_A+4\delta+\varepsilon}\D_4(f)
                                  \sum_{j=0}^4\abs v^j .
\]
For every fixed \(j,L\), the common \((j+L)\)-th moment gives
\begin{equation}\label{replacement:eq:rgfc-tail}
 \E[\abs v^j\1_{\{\abs v>N^{-2/5}\}}]
 \le N^{2L/5}\E\abs v^{j+L}
 \le C_{j,L}N^{-j/2-L/10},
                         \qquad v\in\{x,y\}.
\end{equation}
Choose \(L\) large after \(A,D\).  Then both the original observable
and every term of the polynomial on these tails contribute at most
\(N^{-D-3}\D_5(f)\) per position.
In particular no cancellation of truncated moments has been invoked.

On \(\cE_e^c\), do not Taylor-expand.  The exact summand in
\eqref{replacement:eq:rgfc-exact-derivative} is bounded by
\(2N^{C_A}\norm f_\infty\).  Taking the failure exponent in
\eqref{replacement:eq:rgfc-deleted-prob} larger than \(C_A+D+3\) makes this part
at most \(N^{-D-3}\D_5(f)\) per position, after increasing constants.
No derivative estimate on a bad deleted background is required.

\paragraph{The relative bound after summing entries.}
Apply \eqref{replacement:eq:rgfc-deleted-relative} to
\eqref{replacement:eq:rgfc-degree4} and \eqref{replacement:eq:rgfc-degree5}, and sum all \(N^2\)
positions in \eqref{replacement:eq:rgfc-exact-derivative}.  We obtain
\[
 \abs{g_f'(\theta)}
 \le C\big(\sigma N^{4\delta+\varepsilon}
               +N^{-1/2+5\delta+\varepsilon}\big)
                         Z(\theta)\D_5(f)
       +N^{-D}\D_5(f).
\]
Choose the exponent in the final derivative bound so that \(0<\varepsilon<\delta\).
Both powers are bounded by the error term \(\epsilon_N\) in
\eqref{replacement:eq:rgfc-R1}.
\end{proof}

\subsection{Normalizers, normalized expectations, and center changes}

In this subsection we compare normalizing constants before taking quotients, then restore each ensemble's natural deterministic centering.

\begin{proof}[Proof of Theorem~\ref{replacement:thm:rgfc-main}]
We first compare the normalizing constants directly from the
unnormalized inequality. Taking \(f=1\) in
\eqref{replacement:eq:rgfc-R1} gives the required scalar inequality. The deterministic bound
\eqref{replacement:eq:rgfc-deterministic-cap} is available throughout the mixture,
before any comparison:
\[
 Z(\theta)\ge N^{-C_A}.
\]
Therefore
\begin{equation}\label{replacement:eq:rgfc-logZ}
 \abs{\frac{d}{d\theta}\log Z(\theta)}
       \le C\epsilon_N+N^{-D+C_A}.
\end{equation}
Choose \(D>C_A+c_0+2\) and integrate.  Uniformly for \(0\le\theta\le1\),
\begin{equation}\label{replacement:eq:rgfc-relativeZ}
 \abs{\log\frac{Z(\theta)}{Z(0)}}\le C\epsilon_N,\qquad
 \frac{Z(\theta)}{Z(0)}=1+O(\epsilon_N).
\end{equation}
This proves the relative comparison by estimating the normalizing
constant itself. We can now normalize the general weighted expectation.

For a general test, let \(q_f(\theta)=g_f(\theta)/Z(\theta)\).
Positivity implies \(\abs{q_f(\theta)}\le\norm f_\infty\).
Differentiating the quotient and using
\eqref{replacement:eq:rgfc-R1} for \(f\) and \(1\), we find
\[
 \begin{split}
 \abs{q_f'(\theta)}
 &\le \frac{\abs{g_f'(\theta)}}{Z(\theta)}
       +\norm f_\infty\frac{\abs{Z'(\theta)}}{Z(\theta)}\\
 &\le C\big(\epsilon_N+N^{-D+C_A}\big)\D_5(f).
 \end{split}
\]
Integration over \([0,1]\) gives the second assertion.
All estimates depend on the source only through
\eqref{replacement:eq:rgfc-charge}, so the same constants apply along its ray.
\end{proof}

\begin{corollary}[Compensated microscopic normalizer]\label{replacement:cor:rgfc-comp}
Let \(p=2e_{r_z}\), \(q=2e_{r_w}\), and \(a_j=p+q-2e_c\) for
coordinates in the admissible list, and use identical caps, centers,
and regularizers in both ensembles.  Define
\[
 \mathcal C_{j,L}^{\rm mic}
       =e^{2j}\frac{Z_{a_j,L}}{Z_{p,L}Z_{q,L}} .
\]
Then
\[
 \frac{\mathcal C_{j,Y}^{\rm mic}}
      {\mathcal C_{j,X}^{\rm mic}}=1+O(\epsilon_N).
\]
The assertion does not assume that either compensated normalizer is
bounded.
\end{corollary}

\begin{proof}
Apply \eqref{replacement:eq:rgfc-main} separately to the three sources and take
their ratio.  The deterministic factor \(e^{2j}\) cancels.
If some coordinates coincide, combine the source coefficients;
the same \(\ell^1\) and support bounds still hold.
\end{proof}

\begin{remark}[Scope of the replacement comparison]\label{replacement:rem:rgfc-boundary}
Corollary~\ref{replacement:cor:rgfc-comp} is a relative comparison at the same
microscopic regularizers.  It does not identify
\(\mathcal C_{j,Y}^{\rm mic}\) with any mesoscopic
normalizer.  That requires a separate local-factor reconstruction and
normalization argument.  Likewise raw positive determinant normalizers
require separate scalar regularization and positive cap removal.
No inverse determinant moment or negative cap removal has entered this
section.
\end{remark}

\begin{lemma}[Restoring natural deterministic centers]\label{replacement:lem:rgfc-centers}
Suppose each endpoint center is changed by a vector of
\(\ell^\infty\) norm at most \(s_N\), with the other conventions fixed.
This changes each source normalizer by a factor \(1+O(s_N)\) and each
normalized bounded smooth expectation by at most
\(Cs_N\D_1(f)\), provided \(s_N=o(1)\).
In particular, natural fourth-cumulant corrections can be restored
after comparison with an additional \(O(\sigma)\) error.
\end{lemma}

\begin{proof}
The change in the \(\sqrt2\)-scaled field is at most \(\sqrt2s_N\) per coordinate.
By the cap Lipschitz property and the bound on the sum of absolute coefficients, the two
log-weights differ pointwise by at most \(6\sqrt2s_N\).
This bounds the ratio of their expectations and, after normalization,
the ratio of their densities by \(1+O(s_N)\).
The fundamental theorem of calculus bounds the change in the test
by \(\sqrt2s_N\sum_j\norm{\partial_jf}_\infty\).
There is no extra factor \(n\) outside \(\D_1\).

For the natural correction
\[
 e_j=E_\kappa(z_j,\eta_j)
       =-\frac{\kappa}{4}m(z_j,i\eta_j)^4,\qquad
 \kappa=\E\abs{\xi}^4-2,
\]
the balanced fourth part of \eqref{replacement:eq:rgfc-matching} implies
\(\abs{\kappa_X-\kappa_Y}=O(\sigma)\).
The deterministic \(m\) is bounded on the fixed bulk and stated
height interval.  Thus \(s_N=O(\sigma)\).
During entry replacement the centers remain frozen. The bound just
proved then returns each endpoint to its own natural center.
\end{proof}

\begin{corollary}[Unweighted bounded-test comparison]\label{replacement:cor:rgfc-unweighted}
Under the same moment and coordinate hypotheses, the uncapped vector
\(\cY\) satisfies
\begin{equation}\label{replacement:eq:rgfc-unweighted}
 \abs{\E f(\cY(Y))-\E f(\cY(X))}
       \le C\epsilon_N\D_5(f),\qquad f\in C_b^5(\R^n).
\end{equation}
No cap and no source normalization are needed in this statement.
It applies in particular to every fixed finite list of coarse
positive heights at most one.
Natural endpoint centers may be used with an additional
\(C\sigma\D_1(f)\) error.
\end{corollary}

\begin{proof}
Apply the same mixture argument directly to the uncapped field, with
\(W\equiv1\), \(Z\equiv1\), and
\(F_f(M)=f(\cY(M))\).
The coordinate derivative input is unchanged.
The chain-rule bound uses \(\D_k(f)\), the observable is globally
bounded by \(\norm f_\infty\), and the deleted good event remains
independent of the auxiliary entry.  The complete Taylor cancellation,
fifth moment, and polynomial tails give \eqref{replacement:eq:rgfc-R1} with \(Z=1\).
Integrate directly. The field need not be bounded because the test is
bounded globally; no cap-removal tail estimate is required.
Lemma~\ref{replacement:lem:rgfc-centers} treats the centers.
\end{proof}

\subsection{An explicit matched law with the prescribed initial age}

This subsection constructs a law with a prescribed Gaussian component at the start and the required moment comparison at the endpoint.
Gaussian-divisible laws matching the first three moments with a small
fourth-moment discrepancy have a real-variable antecedent in
\cite[Lemma~3.4]{EYY}; see also the remark after Definition~9.5
in \cite{CL}. We give an explicit complex construction, checking all mixed
moments and the endpoint $\kappa=-1$.

We construct the final variance-one law to satisfy the preceding
replacement theorem. Bernoulli thinning, followed by variance
normalization, amplifies its mixed third moments by the factor that
Gaussian addition will remove. The choice of thinning probability
makes this cancellation exact. The fourth moments are allowed the
stated small discrepancy, including at the endpoint cumulant
\(\kappa=-1\).

\begin{lemma}[Elementary Gaussian-divisible matching]\label{replacement:lem:rgfc-matchlaw}
Let \(\xi\) obey the normalized moment hypotheses
\eqref{replacement:eq:rgfc-profile}.  Fix \(0<v<1/2\), put \(p=(1-v)^3\), and let
\(I\) be Bernoulli with success probability \(p\), independent of \(\xi\).
Let \(g\) be an independent circular complex Gaussian with
\(\E\abs g^2=1\).  Define
\begin{equation}\label{replacement:eq:rgfc-matchlaw}
 \widehat\xi=\frac{I\xi}{\sqrt p},\qquad
 \upsilon=\sqrt{1-v}\,\widehat\xi+\sqrt v\,g .
\end{equation}
Then \(\widehat\xi,\upsilon\) have mean zero, variance one, and vanishing
complex second moment; they satisfy common bounds for every fixed
moment, uniformly for \(v\le1/2\).
Every mixed moment of \(\upsilon\) of total degree at most three equals
the corresponding moment of \(\xi\).
At degree four,
\begin{equation}\label{replacement:eq:rgfc-fourth-exact}
 \E\upsilon^a\overline{\upsilon}^{\,b}
 =\frac{\E\xi^a\bar\xi^{\,b}}{1-v}
       +\1_{\{a=b=2\}}\big(4v(1-v)+2v^2\big),
                         \qquad a+b=4 .
\end{equation}
Thus all fourth discrepancies are \(O(v)\), with constants controlled
by the common fourth moment.
\end{lemma}

\begin{proof}
For each positive integer \(k\),
\[
 \E\abs{\widehat\xi}^k=p^{1-k/2}\E\abs\xi^k .
\]
Since \(p\ge1/8\), these moments are uniformly bounded.
The same follows for \(\upsilon\) by
\(\abs{u+w}^k\le2^{k-1}(\abs u^k+\abs w^k)\).
The first and second moment assertions follow directly by independence.

For total degree three, a Gaussian term with odd Gaussian degree is
zero.  A term with Gaussian degree two either has zero circular
Gaussian moment or leaves a first moment of \(\widehat\xi\), also zero.
Only the pure non-Gaussian term survives, with coefficient
\[
 (1-v)^{3/2}p^{-1/2}=1 .
\]
This proves all mixed third identities.

At total degree four the pure non-Gaussian coefficient is
\((1-v)^2p^{-1}=(1-v)^{-1}\).
The only nonzero terms containing \(g\) have balanced degree
\((a,b)=(2,2)\): two Gaussian factors yield
\(4v(1-v)\E\abs{\widehat\xi}^2\E\abs g^2=4v(1-v)\),
and four Gaussian factors yield
\(v^2\E\abs g^4=2v^2\).
For the other fourth monomials a possible Gaussian pairing leaves
\(\E\widehat\xi^2\) or its conjugate, which is zero.
This proves \eqref{replacement:eq:rgfc-fourth-exact}.
\end{proof}

In particular
\begin{equation}\label{replacement:eq:rgfc-kurtosis}
 \E\abs{\upsilon}^4
       =\frac{\E\abs\xi^4}{1-v}+4v-2v^2,\qquad
 \kappa_\upsilon
       =\frac{\E\abs\xi^4}{1-v}-2(1-v)^2 .
\end{equation}
\paragraph{The endpoint fourth cumulant.}
At the extremal value \(\kappa_\xi=-1\), this gives
\(\kappa_\upsilon=-1+5v+O(v^2)\), an admissible fourth cumulant.
The atom at zero in \(\widehat\xi\) is permitted by the iid hypotheses;
the endpoint itself has a Gaussian component. This explicit calculation
requires no positive distance of the target cumulant from \(-1\).

\paragraph{Matching after the full additive evolution.}
The thinning parameter must use the total Gaussian fraction after
normalizing the final matrix. We choose it so that the initial
Gaussian variance is exactly \(\sigma\), while the endpoint has the
moment matching needed for relative replacement.

\begin{proposition}[Exact initial age and matched final normalization]
\label{replacement:prop:rgfc-age}
Let \(t=\frac12\log N\), fix \(b>0\), and set
\[
 \sigma=N^{-1/200},\qquad \tau=\frac{t^b}{N},\qquad
 v=\frac{\sigma+\tau}{1+\tau}.
\]
Apply Lemma~\ref{replacement:lem:rgfc-matchlaw} at this value of \(v\) to the target
normalized entry law, and let \(\widehat X\) have independent entries
\(\widehat\xi/\sqrt N\).
Let \(G_0,G_1\) be independent normalized Ginibre matrices, independent
of \(\widehat X\), and define
\begin{equation}\label{replacement:eq:rgfc-age}
 X_0=\sqrt{1-\sigma}\,\widehat X+\sqrt\sigma\,G_0,\qquad
 X_\tau=X_0+\sqrt\tau\,G_1,\qquad
 Y=\frac{X_\tau}{\sqrt{1+\tau}} .
\end{equation}
The initial matrix \(X_0\) has Gaussian variance exactly \(\sigma\).
The final matrix \(Y\) is variance-one iid and satisfies
\eqref{replacement:eq:rgfc-matching} with the target \(X\).
All base, initial, and endpoint normalized laws have common moment
bounds uniform in \(N\).
\end{proposition}

\begin{proof}
Every fixed \(b\) gives \(\tau=o(\sigma)\), so \(v\asymp\sigma\) and
\(v<1/2\) for large \(N\).
The independent Gaussian combination
\[
 G_{\rm comb}
       =\frac{\sqrt\sigma\,G_0+\sqrt\tau\,G_1}{\sqrt{\sigma+\tau}}
\]
is normalized Ginibre and is independent of \(\widehat X\).
Since \(1-v=(1-\sigma)/(1+\tau)\), we have exactly
\[
 \frac{X_\tau}{\sqrt{1+\tau}}
       =\sqrt{1-v}\,\widehat X+\sqrt v\,G_{\rm comb}.
\]
Its entry law is therefore \(\upsilon/\sqrt N\) from
Lemma~\ref{replacement:lem:rgfc-matchlaw}, giving degree-three exact matching and
fourth discrepancy \(O(vN^{-2})=O(\sigma N^{-2})\).
The common moment bounds were proved in that lemma and persist under
the bounded Gaussian coefficients in \eqref{replacement:eq:rgfc-age}.
\end{proof}

Matching is imposed on the final normalized endpoint; the choice of
base law already includes the later Brownian increment. Along the
additive evolution, \(X_s/\sqrt{1+s}\) has Gaussian fraction
\((\sigma+s)/(1+s)\asymp\sigma\), uniformly for \(0\le s\le\tau\).
Thus the exact initial age and the comparable normalized ages required
by Sections~\ref{chap:i2}, \ref{chap:noise}, \ref{chap:lf} and
\ref{chap:signed} are all available. The common moment bounds make
the \(N\)-dependence of the base law admissible for these inputs.

\paragraph{Returning coordinates and deterministic centers.}
Variance normalization changes the matrix, spatial point and height
together. The density-centered field is exactly invariant under this
joint change, as the next lemma records.

\begin{lemma}[Exact variance and coordinate conversion]
\label{replacement:lem:rgfc-variance}
Let \(s>0\), \(M=\sqrt s\,\widetilde M\),
\(z=\sqrt s\,\widetilde z\), and
\(\eta=\sqrt s\,\widetilde\eta\).
Define the variance-\(s\) deterministic potential by
\[
 I_{s,z}(\eta)=\log s+
             I_{z/\sqrt s}(\eta/\sqrt s)
\]
and use it in \eqref{replacement:eq:rgfc-field} to define
\(\Phi_{N,s}(z,\eta;M)\).
Then
\begin{equation}\label{replacement:eq:rgfc-variance}
 \Phi_{N,s}(z,\eta;M)
       =\Phi_N(\widetilde z,\widetilde\eta;\widetilde M)
\end{equation}
exactly.
For Proposition~\ref{replacement:prop:rgfc-age}, take \(s=1+\tau\).
\end{lemma}

\begin{proof}
The positive matrix inside the determinant is \(s\) times the
corresponding matrix for
\((\widetilde M,\widetilde z,\widetilde\eta)\).
Its half log determinant therefore increases by
\(\frac N2\log s\), precisely the increase in the deterministic
centering.  They cancel.
\end{proof}

The natural fourth-cumulant center obeys the same exact scaling.
If $\kappa_{\widetilde M}$ is the fourth cumulant of the normalized
entry $\sqrt N\,\widetilde M_{ab}$, then the physical variable
$\sqrt N M_{ab}$ has cumulant
$\kappa_{\rm phys}=s^2\kappa_{\widetilde M}$, while
\[
 m_{s,z}(i\eta)=s^{-1/2}m_{1,\widetilde z}(i\widetilde\eta).
\]
Hence
$\kappa_{\rm phys}m_{s,z}(i\eta)^4
 =\kappa_{\widetilde M}m_{1,\widetilde z}(i\widetilde\eta)^4$.
Thus carrying the frozen center numerically under this conversion is
also the natural variance-corrected convention; a change of endpoint
entry law is handled separately by Lemma~\ref{replacement:lem:rgfc-centers}.

A normalized coordinate \((z,\eta)\) corresponds at the additive
endpoint to
\[
 (z_{\rm phys},\eta_{\rm phys})
       =(\sqrt{1+\tau}\,z,\sqrt{1+\tau}\,\eta),
\]
with the variance-adjusted deterministic centering.
There is no \(O(N\tau)\) error.
The source replacement theorem is invoked only in its normalized
height domain \(\eta\le1\).  When an earlier flow theorem also needs
physical heights at most one, the actual coarse list can be chosen
inside a fixed smaller interval, as for integer clocks \(k\ge1\).

\subsection{Parameter order and uses in flow and comparison}

Now we fix the replacement parameters and explain how the resulting estimates return the flow comparisons to the original ensemble. We fix the cap, coefficient bound, list-length bound, \(\delta<1/1200\), and final
accuracy first. Choose the exponent in the input derivative bound sufficiently small
and the failure and high-moment orders after these constants.
Only afterwards choose any fixed microscopic logarithmic exponent
\(q\) and mesoscopic exponent \(b\).
For such fixed exponents,
\[
 \eta_{\rm mic}=\frac{t^{-q}}{N}\ge N^{-1-\delta},
 \qquad \tau=\frac{t^b}{N}=o(\sigma)
\]
for all sufficiently large \(N\).
Each choice depends only on parameters already fixed.

We have proved the unnormalized relative inequality
\eqref{replacement:eq:rgfc-R1}, the relative comparison of each capped
normalizing constant, the normalized bounded-test comparison, and the
unweighted uncapped bounded-test comparison \eqref{replacement:eq:rgfc-unweighted}.
The construction of the matched endpoint also supplies the Gaussian-aged
initial law for the local-factor reconstruction in Section~\ref{chap:lf}.
Section~\ref{chap:calibration}
combines that reconstruction with the separate ratios of normalizing
constants to compare microscopic and mesoscopic weights. The unweighted
comparison compares the coarse field observations in the two ensembles,
while Section~\ref{chap:signed} gives the Gaussian approximation under
that initial law at its permitted heights. Positive cap removal and
terminal annular estimates retain their separate proofs.

\endgroup

\begingroup
\newcommand{\E}{\mathbb E}
\newcommand{\PP}{\mathbb P}
\newcommand{\Tr}{\operatorname{Tr}}
\newcommand{\ind}{\mathbf1}
\newcommand{\dd}{\,\mathrm d}
\newcommand{\br}[1]{\langle #1\rangle}
\section{Normalizing constants and comparison of field histories}\label{chap:calibration}

In this section we compare the exact normalizing constants while keeping the earlier field values that control pairs of nearby points.
We first compare the observed history along the additive matrix flow,
keeping its density centering and integer scale labels. The sparse
Stein comparison then gives the initial size of the normalizing constant
and an exponential moment with a negative coefficient on the capped
history coordinate. Local-factor
reconstruction and removal of the positive caps give both sides of the
normalizing-constant comparison; relative entry replacement transfers it
to the original law. These steps precede any rare-barrier test. The
coordinate with the negative coefficient stays capped throughout.

The use of a negative coarse coefficient follows the near-diagonal
arguments in \cite[proof of Lemma~2.10]{CFLW} and
\cite[proof of Lemma~7.8, especially (7.44)]{CL}.
The mixed-Laplace asymptotics in \cite[Proposition~7.5]{CL} form part
of that subcritical framework. The comparison below treats critical
root charges in the terminal annuli of
\eqref{calibration:eq:sc-domain}, retaining relative error control and
the negative history cap. The positive integrability gap required in
\cite[Lemma~2.8]{CFLW} closes at criticality; the retained-history
estimates in Section~\ref{chap:annuli} supply the summability needed here.

\subsection{Coordinates and the relative assertion}

We first specify the fields and caps needed to compare microscopic and mesoscopic normalizing constants with small relative error.
Let $X_0=N^{-1/2}(\xi^{(0)}_{ab})$ be complex iid, with
\[
 \E\xi^{(0)}=\E(\xi^{(0)})^2=0,\qquad
 \E|\xi^{(0)}|^2=1.
\]
Assume Gaussian age comparable to $\sigma=N^{-1/200}$ and
common uniform bounds for all fixed entry moments. Let
\[
 X_s=X_0+N^{-1/2}B_s,\quad 0\le s\le\tau,\qquad
 M=\frac{X_\tau}{\sqrt{1+\tau}},
\]
where $B$ is standard complex matrix Brownian motion independent
of $X_0$. The additive matrix has entry variance $(1+s)/N$; its
variance-one form is $X_s/\sqrt{1+s}$. By
Lemma~\ref{replacement:lem:rgfc-variance}, passing between these forms
rescales the spatial point and height together, with exact cancellation
of the density-centering change. Set
\[
 t=\tfrac12\log N,\quad d=t^{1/3},\quad H=t^b,\quad
 \tau=H/N,\quad \eta_{\rm m}=N^{-1}t^{-q},
 \quad T=\lfloor t-d/2-C_0\rfloor .
\]
All points belong to a fixed compact subset of the unit disk, with
a fixed slightly larger compact used after rescaling. We consider
the domain
\begin{equation}
 t/2<j\le T,\qquad c_1e^{-j}\le |z-w|\le c_2e^{-j},
 \qquad j\in\mathbb Z.
 \label{calibration:eq:sc-domain}
\end{equation}
The constants $c_1,c_2,C_0$ are fixed. The local-factor theorem is
used only in this domain.

Let $\rho_s^v$ be the deterministic Hermitization density of variance
$1+s$, and write
\[
 \Phi_s(v,\eta)=\frac12\log\det((X_s-v)(X_s-v)^*+\eta^2)
  -\frac N2\int\log(x^2+\eta^2)\rho_s^v(x)\dd x .
\]
For a variance-one matrix with fourth cumulant $\kappa$, use the
bounded mean correction
\[
 E_\kappa(v,\eta)=-\frac{\kappa}{4}m^v(i\eta)^4,\qquad
 Y(v,\eta)=\sqrt2\,[\Phi(v,\eta)-E_\kappa(v,\eta)].
\]
At physical variance $1+s$, the mean correction is the one for
the normalized matrix, point, and height, using its actual normalized
fourth cumulant. This is the deterministic fourth-cumulant correction,
not an identity with the finite-$N$ expectation. The same convention
is kept in every scalar normalizer below.

The exact clock is specified by
\begin{equation}
 u_k=\tanh k,\qquad
 \eta_k(v)=(1-u_k)\sqrt{u_k^{-1}-|v|^2}.
 \label{calibration:eq:sc-clock}
\end{equation}
Define the physical points
$v_z=\sqrt{1+\tau}z$, $v_w=\sqrt{1+\tau}w$.
The endpoint root coordinates are
\[
 y_z=Y^M(z,\eta_{\rm m}),\qquad
 y_w=Y^M(w,\eta_{\rm m}),\qquad
 y_c=Y^M(z,\eta_j(z)).
\]
The initial coordinates are
\[
 x_z=Y^{X_0}(v_z,B_{v_z}\tau),\quad
 x_w=Y^{X_0}(v_w,B_{v_w}\tau),\quad
 x_c=Y^{X_0}(v_z,\eta_j(v_z)),\quad B_v=\sqrt{1-|v|^2}.
\]
The positive root coordinates move from microscopic to mesoscopic
height through the local-factor argument. The observed history
coordinate instead keeps its original integer label $j$ and exact
clock height at the rescaled point; its small transport error is
proved below.

Fix an odd smooth nondecreasing cap $c_A$ with
\begin{equation}
 c_A(x)=x\quad (|x|\le A\log N),\qquad
 |c_A(x)|\le2A\log N,\qquad 0\le c_A'(x)\le1.
 \label{calibration:eq:sc-cap}
\end{equation}
All fixed positive-order derivatives are uniformly bounded. The
same cap is used in all three coordinates appearing in the exponential weight. Write
\begin{align}
 Z_{\rm m}(a_j)&=\E e^{2c_A(y_z)+2c_A(y_w)-2c_A(y_c)},&
 Z_{\rm m}(p)&=\E e^{2c_A(y_z)},&
 Z_{\rm m}(q_*)&=\E e^{2c_A(y_w)}, \label{calibration:eq:sc-Zm}\\
 Z_0(a_j)&=\E e^{2c_A(x_z)+2c_A(x_w)-2c_A(x_c)},&
 Z_0(p)&=\E e^{2c_A(x_z)},&
 Z_0(q_*)&=\E e^{2c_A(x_w)}. \label{calibration:eq:sc-Z0}
\end{align}
Here $a_j=2e_{r_z}+2e_{r_w}-2e_c$,
$p=2e_{r_z}$, and $q_*=2e_{r_w}$.
For any variance-one endpoint law $L$ with fourth cumulant
$\kappa_L$, the uncapped positive normalizers in this section are
\begin{equation}
 m_{L,v}^{\eta}
   :=\E_L\exp\{\chi[\Phi^L(v,\eta)-E_{\kappa_L}(v,\eta)]\},
 \qquad \chi=2\sqrt2,\qquad \eta\ge0.
 \label{calibration:eq:sc-positive-normalizers}
\end{equation}
In particular $m_{L,v}^0$ retains the same mean-corrected
convention at the raw endpoint.

\begin{theorem}[Comparison of microscopic and mesoscopic normalizing constants]
\label{calibration:thm:sc-main}
For every fixed $P>0$, the fixed parameters can be selected in the
order specified below so that, uniformly on \eqref{calibration:eq:sc-domain},
\begin{equation}
 \frac{e^{2j}Z_{\rm m}(a_j)}
             {Z_{\rm m}(p)Z_{\rm m}(q_*)}
 =
 \frac{e^{2j}Z_0(a_j)}{Z_0(p)Z_0(q_*)}
       [1+o(t^{-P})].
 \label{calibration:eq:sc-bridge}
\end{equation}
Both ratios are bounded above and below by positive constants.
Moreover $Z_{\rm m}(a_j)\ge cN$ for large $N$.

If an original complex iid ensemble is compared with $M$ by the
relative iid-mixture theorem, matching mixed complex moments through
degree three and with degree-four mismatch $O(\sigma)$, then
\eqref{calibration:eq:sc-bridge} also holds with the microscopic normalizers of
that original ensemble on the left. For either endpoint law $L$,
each positive microscopic denominator can be replaced first by
$m_{L,v}^{\eta_{\rm m}}$ and then by $m_{L,v}^0$, with relative
error $o(t^{-P})$.
\end{theorem}

\subsection{The history estimate along the matrix flow}

Now we show that the retained field history changes little during the matrix flow when the deterministic centering is kept exact.
The three principal error terms below come from the martingale along
a characteristic, the return to the prescribed initial height or
spatial point, and the quadratic resolvent drift. Exact density
centering is needed before any of these errors is estimated.

\begin{lemma}[History transport with exact density centering]
\label{calibration:lem:sc-history}
Let $n\le C t$ deterministic points and heights lie in the enlarged
bulk and satisfy $\eta_r\in[\eta_{\min},C]$, where
$\eta_{\min}\ge cN^{-1}e^d$. For every fixed $D>0$, outside
an event of probability at most $C_DN^{-D}$,
\begin{equation}
 \max_r|\Phi_\tau(v_r,\eta_r)-\Phi_0(v_r,\eta_r)|
 \le C_D\left[
 t\sqrt{\frac{\tau}{\eta_{\min}}}
 +\frac{t\tau}{\eta_{\min}}
 +\frac{t^2\tau}{N\eta_{\min}^2}+N^{-c}\right].
 \label{calibration:eq:sc-history}
\end{equation}
The same bound includes the changes in variance, clock, and deterministic mean correction
needed to compare $y_c$ with $x_c$, and all histories through $T$.
In particular their maximum coordinate error, denoted
$\rho_{\rm hist}$, satisfies
\begin{equation}
 \rho_{\rm hist}\le C_D\left[
 t^{1+b/2}e^{-d/2}+t^{b+1}e^{-d}
       +t^{b+2}e^{-2d}+N^{-c}\right]=o(t^{-P})
 \label{calibration:eq:sc-history-terminal}
\end{equation}
for every fixed $P$.
\end{lemma}

\begin{proof}
\emph{A common event under the unweighted matrix law.}
We first obtain uniform deterministic-parameter bounds, then use
their pointwise thresholds to stop the martingale. The averaged
one-resolvent law \cite[Theorem 3.1]{DEL}, including fixed $2$ by $2$
block insertions, gives in a fine range
\[
 |\br{(G_s-M_s)A}|\le C_Dt/(N\eta),\qquad \|A\|\le1.
\]
Only a fixed deterministic block basis is needed; linearity then
covers all bounded block insertions used below.
We use the fine domain exponent $1/200$ in the logarithmic local law
just cited and split at $a=N^{-1/100}$. Above $a/4$, the ordinary fixed-power law
\cite[Theorem 2.3]{CL} gives $C_DN^{1/100}/(N\eta)$.
The split is buffered: use the fine
law below $a/2$ and the coarse law above $a/4$. Normalize $X_s$
by $\sqrt{1+s}$ when invoking these statements; its Gaussian age
remains comparable to $\sigma$ and its moment bounds are uniform.

Here is a quantitative common-event argument. Use a time mesh $N^{-40}$ and spatial/height meshes $N^{-12}$
on the enlarged bulk and the height interval
$[\eta_{\mathrm{low}},C+1]$, where
$\eta_{\mathrm{low}}=c_*N^{-1}e^d$ and the fixed $c_*>0$ is
chosen below the heights of all normalized endpoints and return
intervals. On the fine part of this mesh, bulk positivity gives
$N\eta\Im m\gg C_D^2\log N$, so the lower-domain condition of
\cite[Theorem 3.1]{DEL} is satisfied. No logarithmic local law
at a fixed multiple of $N^{-1}$ is invoked.
Uniform high moments give a polynomial
operator-norm bound. The Brownian maximal inequality, followed by
an entry union and the Frobenius bound, gives
\[
 \sup_{|s-r|\le N^{-40}}\|X_s-X_r\|\le C_DN^{-19}
\]
outside an event whose probability decays faster than any fixed inverse power of the matrix size. Resolvent Lipschitz
constants are at most $CN^2$. Interpolation thus costs at most
$CN^{-10}$ in the deterministic parameters and $CN^{-17}$ in time.
A mesh of at most $N^{80}$ points suffices; select the pointwise
failure exponent as $D+100$. These errors are smaller than the
required one-resolvent bounds. In particular the fine event gives
\[
 |m_{N,s}-m_s|\le C_Dt/(N\eta),\qquad \Im m_{N,s}\le C.
\]
No random vector is inserted into an isotropic local law.

\emph{The exact characteristic identity.}
For each endpoint solve the deterministic backward characteristic
\[
 \dot\eta_s=-\Im m_s^v(i\eta_s),\qquad \eta_\tau=\eta_r.
\]
It satisfies $\eta_s\ge\eta_r$ and $|\eta_0-\eta_r|\le C\tau$.
At a fixed spatial point the positive singular-value drivers obey
$\dd[b_i,b_k]_s=\delta_{ik}\dd s$. Independence across different
spatial points is not used. Direct It\^o calculus gives
\begin{align}
 \dd\Phi_s(v,\eta_s)
 &=\dd\mathcal M_s
       -\frac N2(m_{N,s}(i\eta_s)-m_s(i\eta_s))^2\dd s,
                                                        \label{calibration:eq:sc-Ito}\\
 \dd\mathcal M_s
 &=\frac1{\sqrt{2N}}\sum_{i=1}^N
        \frac{\lambda_i(s)}{\lambda_i(s)^2+\eta_s^2}
                  \dd b_i(s).                         \label{calibration:eq:sc-mart}
\end{align}
For the signed indices $i\in\mathcal I_N$, put
$g_i=(\lambda_i(s)-i\eta_s)^{-1}$. Symmetrizing the signed Dyson
interaction in the half-logarithm gives
\[
 -(8N)^{-1}\big[(\textstyle\sum g_i)^2-\textstyle\sum g_i^2\big],
\]
and the diagonal It\^o correction is $-(8N)^{-1}\sum g_i^2$.
The characteristic contributes $Nmm_N$ to the random logarithm.
The deterministic logarithmic center has characteristic derivative
$Nm^2/2$ and is subtracted. Thus the total drift is
\[
 -\frac N2m_N^2+Nmm_N-\frac N2m^2
       =-\frac N2(m_N-m)^2,
\]
exactly as in \eqref{calibration:eq:sc-Ito}.
On the imaginary
axis the displayed real drift is nonnegative, and we bound it in
absolute value.

\emph{Stopping and the three error costs.}
Stop at the first time the current resolvents leave slightly enlarged
deterministic envelopes from the common local-law bounds. This
buffered exit depends only on the current and past matrices, so it is
an adapted stopping time; initial failure is included in the exceptional
event. The common event prevents exit, but we do not condition the
martingale on that event.
For the stopped martingale, Ward's identity gives
\[
 [\mathcal M^\theta]_\tau
 \le\frac1{2N}\int_0^{\tau\wedge\theta}
             \sum_i\frac1{\lambda_i^2+\eta_s^2}\dd s
 =\frac12\int_0^{\tau\wedge\theta}
              \frac{\Im m_{N,s}(i\eta_s)}{\eta_s}\dd s
 \le C\tau/\eta_r.
\]
The exponential martingale inequality at threshold
$C_Dt\sqrt{\tau/\eta_r}$ has probability $2e^{-c_Dt^2}$.
A union over $Ct$ endpoints is smaller than each fixed $N^{-D}$.
The absolute drift costs $C_Dt^2\tau/(N\eta_r^2)$.
The exact centered derivative
\[
 \partial_\eta\Phi_0(v,\eta)
       =N\Im(m_{N,0}^v(i\eta)-m_0^v(i\eta))
\]
bounds the error from the initial height adjustment by $C_Dt\tau/\eta_r$.
The coarse range costs at most
$C_D[t\sqrt{\tau/a}+N^{1/100}\tau/a
             +N^{2/100}\tau/(Na^2)]$, which is $N^{-c}$
for a fixed $c>0$ and all sufficiently large $N$. This proves
\eqref{calibration:eq:sc-history}.

\emph{Variance rescaling, height adjustment, and mean correction.}
For the variance rescaling, the deterministic density scales by
\[
 \rho_{\tau}^{\sqrt{1+\tau}z}(x)
 =(1+\tau)^{-1/2}\rho_0^z(x/\sqrt{1+\tau}).
\]
Hence
\begin{equation}
 \Phi^M(z,\eta)=
   \Phi_\tau(\sqrt{1+\tau}z,\sqrt{1+\tau}\eta).
 \label{calibration:eq:sc-variance}
\end{equation}
The two terms each acquire $N\log(1+\tau)/2$, which cancels.
At $v=\sqrt{1+\tau}z$ and integer $k\ge1$, formula
\eqref{calibration:eq:sc-clock} gives
\[
 |\sqrt{1+\tau}\eta_k(z)-\eta_k(v)|\le C\tau\eta_k(z).
\]
The height adjustment contributes at most $Ct\tau$ in the fine region.
If a return to the original spatial point is wanted, for a real
direction $e$ and its fixed block insertion $J_e=\partial_eH^z$,
\[
 \partial_e\Phi_0(z,\eta)
       =\tfrac12\Tr[(G_0^z-M_0^z)J_e].
\]
The deterministic identity follows by differentiating the MDE
potential; its two spatial derivatives are $2u\Re z$ and
$2u\Im z$. Thus the same averaged block law bounds this gradient
by $Ct/\eta$, paying an $O(\tau)$ spatial displacement by
$Ct\tau/\eta$. There is no uncentered $O(N\tau)$ error.

Gaussian addition leaves the physical fourth cumulant unchanged;
variance normalization then gives the actual normalized cumulant
$\kappa_\tau=\kappa_0/(1+\tau)^2$. The bounded deterministic
mean correction changes by $O(\tau)$ in these returns, using bulk
MDE smoothness up to height zero. Thus the estimate includes the
mean-corrected coordinates. Substitution of
$\eta_{\min}\ge cN^{-1}e^d$ proves
\eqref{calibration:eq:sc-history-terminal}. Keep each original integer label:
recomputing a floor of the rescaled separation is unnecessary.
\end{proof}

\begin{corollary}[History transport at a logarithmic floor]
\label{calibration:cor:sc-history-log}
Fix $L>1$ and $b>0$. Under the other assumptions of
Lemma~\ref{calibration:lem:sc-history}, its estimate \eqref{calibration:eq:sc-history}
and its exact variance, clock, spatial, and mean returns remain valid
when the height assumption is replaced by
$\eta_{\min}\ge c t^L/N$. The constants may depend on these fixed
exponents, the bulk compact, and the prescribed failure exponent.
In particular, with $\tau=t^b/N$,
\begin{equation}
 \rho_{\rm hist}\le C_D\left[
 t^{1+(b-L)/2}+t^{1+b-L}+t^{2+b-2L}+N^{-c_0}\right].
 \label{calibration:eq:sc-history-log}
\end{equation}
This estimate is not asserted to vanish for every $L>1$ and $b>0$.
If $L=b+2a>1$, it becomes
\begin{equation}
 \rho_{\rm hist}\le C_D\left[
 t^{1-a}+t^{1-2a}+t^{2-b-4a}+N^{-c_0}\right].
 \label{calibration:eq:sc-history-buffer}
\end{equation}
For a test with first derivative norm $O(t^\nu)$, where $\nu\ge0$,
for each fixed $P>0$ the transport error is $o(t^{-P})$ whenever
$a>P+\nu+1$; all exponents are fixed before $N\to\infty$.
\end{corollary}
\begin{proof}
Repeat the preceding common-event construction with
$\eta_{\mathrm{low}}=c_*t^L/N$, choosing $c_*>0$ below the
normalized endpoints and all return intervals. Bulk positivity and
$L>1$ give
\[
 N\eta_{\mathrm{low}}\Im m\ge c't^L
       \gg C_D^2\log N.
\]
Thus the lower-domain condition of \cite[Theorem 3.1]{DEL} still holds
for every fixed failure exponent. The buffered split at
$a_0=N^{-1/100}$ (the quantity denoted $a$ in that proof) ensures
$\eta/\Im m\le N^{-1/200}$ on the fine portion; the ordinary
fixed-power law supplies the overlapping coarse portion.
The normalized evolving matrices have the same uniform moment bounds
and comparable Gaussian age as before. Since
$\eta\ge c_*t^L/N\ge N^{-1}$ eventually, the $CN^2$ resolvent
Lipschitz bound, the time and spatial/height meshes, and their
$N^{80}$ union bound are unchanged.

The same buffered adapted stopping time therefore gives bracket
$C\tau/\eta_r$, martingale threshold
$C_Dt\sqrt{\tau/\eta_r}$, absolute drift
$C_Dt^2\tau/(N\eta_r^2)$, and centered initial height or spatial
return $C_Dt\tau/\eta_r$. The coarse contribution is still
$N^{-c_0}$, because $b,L$ are fixed and the split is a fixed power
of $N$. The exact variance identity and the relative clock displacement
$O(\tau\eta_k)$ are unchanged, and bounded mean corrections cost
$O(\tau)$. This proves the same estimate on the enlarged height domain.
Substitution gives \eqref{calibration:eq:sc-history-log} and
\eqref{calibration:eq:sc-history-buffer}. For a test with the stated
summed first-derivative bound, the coordinate error is multiplied by
$t^\nu$ once. The condition $a>P+\nu+1$ makes the martingale term
$o(t^{-P})$ and also bounds the smaller errors from adjusting the height or spatial point and from the drift at the required rate. Finally $e^{t^{1/3}}/t^L\to\infty$ for every
fixed $L$, so all original terminal-annulus applications and
\eqref{calibration:eq:sc-history-terminal} retain their stronger
$o(t^{-P})$ conclusion for every fixed $P$.
\end{proof}

\subsection{Initial geometry and the negative-cap moment}

We first obtain bounds on the initial normalizing constants and on the
negative-cap moment from the sparse-weight theorem. Positive-cap removal
then identifies the uncapped initial root factors. The resulting lower
bound allows us to turn local-factor reconstruction into a relative
microscopic comparison.

The initial sparse-weight theorem applies to every coefficient vector supported
on at most three coordinates with $\ell^1$ norm at most six,
uniformly along its ray. It supplies
\begin{equation}
 \log Z_0(a)=\tfrac12a^{\mathsf T}\Sigma a+O(r_N),
 \qquad r_N=o(t^{-P})
 \label{calibration:eq:sc-source}
\end{equation}
after enlarging the fixed exponent $b$. All initial capped heights
are between its logarithmic floor and one.

\begin{lemma}
\label{calibration:lem:sc-initial}
Uniformly on \eqref{calibration:eq:sc-domain},
\begin{align}
 Z_0(p),Z_0(q_*)&\asymp N/H, \label{calibration:eq:sc-one-scale}\\
 Z_0(a_j)&\asymp (N/H)^2e^{-2j}, \label{calibration:eq:sc-signed-scale}\\
 \E e^{-4c_A(x_c)}&\le Ce^{8j}\le CN^4. \label{calibration:eq:sc-negative0}
\end{align}
\end{lemma}

\begin{proof}
At exact clocks $u=\tanh s$, $v=\tanh r$, the stability determinant is
\[
 D=1-uv-uv\{(1-u)|z|^2+(1-v)|w|^2\}+uv|z-w|^2.
\]
On a fixed bulk, with $\delta>0$ depending only on the bulk margin,
\[
 \delta(1-uv)+uv|z-w|^2
       \le D\le(1-uv)+uv|z-w|^2.
\]
Moreover $\eta_s(z)=2B_ze^{-2s}[1+O(e^{-2s})]$. Thus the
external heights $B_vH/N$ have clocks
$\ell_v=t-(b/2)\log t+O(1)$, and
$\ell_v-T\longrightarrow\infty$. The exact covariance is
\[
 \Sigma_{ir}=-\tfrac12\log D_{ir}
                +\tfrac{\kappa_0}{2}m_i^2m_r^2.
\]
The cumulant term is uniformly bounded. Consequently
\[
 \Sigma_{r_vr_v}=\ell_v+O(1),\qquad
 \Sigma_{r_zr_w},\Sigma_{r_zc},\Sigma_{r_wc},\Sigma_{cc}=j+O(1).
\]
These bounds also hold at $\kappa_0=-1$; no inverse covariance is
used. Apply \eqref{calibration:eq:sc-source} to $p,q_*,a_j$. In the compensated
ratio the leading terms cancel exactly:
$2j+4j-4j-4j+2j=0$. This proves
\eqref{calibration:eq:sc-one-scale}--\eqref{calibration:eq:sc-signed-scale}.
Finally $-4e_c$ has one nonzero coordinate and $\ell^1$ norm four, so the same theorem
applies and gives
$\E e^{-4c_A(x_c)}=\exp(8\Sigma_{cc}+O(r_N))$, proving
\eqref{calibration:eq:sc-negative0}.
\end{proof}

Put
\[
 H_0=e^{-2c_A(x_c)},\qquad H_\tau=e^{-2c_A(y_c)}.
\]
By Lemma~\ref{calibration:lem:sc-history} and cap Lipschitz continuity,
$H_\tau/H_0=\exp(O(\rho_{\rm hist}))$ on its event, including
when the negative cap is active. On its complement
$H_\tau^2\le N^{8A}$. Selecting its probability exponent after
$A$ gives
\begin{equation}
 \E H_\tau^2\le Ce^{8j}+N^{8A-D}\le CN^4.
 \label{calibration:eq:sc-negative-tau}
\end{equation}
This estimate keeps the same cap on the mesoscopic history coordinate
at both endpoints. It supplies exactly the negative moment needed
for positive-cap removal below.

\subsection{Removing positive caps without a microscopic lower tail}

We next show that the positive caps can be removed because large values have controlled moments and small values carry little weight.
Indeed, the positive caps are removed using scalar moments and an upper tail.
The contribution from their lower sides will instead be bounded by a deterministic weight estimate.
First, Proposition {\hyperref[scalar:prop:sr-exact-upper]{\ref*{scalar:prop:sr-exact-upper}}}
in Section~\ref{chap:scalar} gives for fixed $r\ge2$
\[
 \E e^{r\Phi(v,\eta;X)}
       \le \E e^{r\Phi(w,0;\widetilde X_s)},\qquad
 s=\eta/\Im m^v(i\eta),\quad w=v/\sqrt{1+s}.
\]
The deterministic MDE potential cancels exactly in this inequality.
The sharp one-root theorem \cite[Theorem 9.7]{CL}, uniformly for laws satisfying the
common moment bounds, applies to the normalized Gaussian-convolved
law. In particular at every height used in this section,
\begin{equation}
 \E e^{4Y}\le CN^4,\qquad \E e^{8Y}\le CN^{16},\qquad
 \E e^{2y_v}\asymp N.
 \label{calibration:eq:sc-positive}
\end{equation}
The last asymptotic concerns the microscopic endpoint roots $y_v$:
it uses Theorem~\ref{scalar:thm:sr-regularization} together with the
sharp raw one-root moment. The initial roots have the different scale
given in \eqref{calibration:eq:sc-one-scale}. Bounded mean corrections
affect only the constants. The weaker
positive polynomial moments already in the local-factor section also suffice
with enlarged fixed error exponents.

The upper-tail estimate follows from the same scalar input: for each fixed
$r\ge\sqrt2$, uniformly for heights in $[0,1]$,
\[
 \E e^{rY}\le C_r N^{r^2/4}.
\]
Indeed apply the sharp positive raw moment and the exact
regularization upper comparison with raw charge $\sqrt2r\ge2$;
the bounded mean correction changes only $C_r$.
Taking $r=2A$ in Markov's inequality proves directly
\begin{equation}
 \PP(Y>A\log N)\le C_A N^{-A^2}
                    \le C_A N^{-c_1A^2}
 \label{calibration:eq:sc-upper-tail}
\end{equation}
for fixed $A\ge1$ and, for example, $c_1=1/2$.
The scalar tail bound \cite[Lemma 8.1]{CL} supplies
a two-sided bound only at mesoscopic heights. We do not extend
that statement to the microscopic positive coordinates.

\begin{lemma}[Positive-cap removal with a capped history weight]
\label{calibration:lem:sc-cap-removal}
Suppose $R_v=e^{2Y_v}$, $0\le H_c\le N^{4A}$,
$\E H_c^2\le CN^4$, and \eqref{calibration:eq:sc-positive}--\eqref{calibration:eq:sc-upper-tail}
hold for $Y_z,Y_w$. For every fixed $L>0$, a sufficiently large
fixed $A$ gives
\[
 \left|\E H_c e^{2c_A(Y_z)+2c_A(Y_w)}
           -\E H_c R_zR_w\right|\le CN^{-L}.
\]
The same assertion holds for each single positive normalizer with
$H_c$ omitted.
\end{lemma}

\begin{proof}
\emph{The lower sides of the positive caps.}
Let $\epsilon_A=N^{-2A}$ and $R_A(y)=e^{2c_A(y)}$. The cap
properties imply
\[
 R_A(y)\le e^{2y}+\epsilon_A.
\]
Indeed for $y\ge0$ one has $c_A(y)\le y$; for
$y<-A\log N$ one has $c_A(y)\le-A\log N$.
Writing $R_+(y)=\min\{R_A(y),e^{2y}\}$ gives
$0\le R_A-R_+\le\epsilon_A$, whereas $e^{2y}-R_+(y)$ vanishes
unless $y>A\log N$. Thus the absolute difference in the lemma
is at most
\begin{align*}
 &\epsilon_A\E[H_c(R_z+R_w)]+\epsilon_A^2\E H_c\\
 &\hspace{15mm}
 +\E[H_cR_zR_w\ind_{\{\max(Y_z,Y_w)>A\log N\}}].
\end{align*}
Cauchy--Schwarz bounds the first line by
$CN^{-2A+4}+CN^{-4A+2}$, without any lower-tail probability.
\emph{The upper sides of the positive caps.}
For the second line, use $H_c\le N^{4A}$,
\[
 \E(R_zR_w)^2
     \le(\E e^{8Y_z}\E e^{8Y_w})^{1/2}\le CN^{16},
\]
and \eqref{calibration:eq:sc-upper-tail}. Its bound is
$CN^{4A+8-c_1A^2/2}$. Choose $A$ large enough for the requested
fixed $L$. The single-root argument drops $H_c$ and one root.
\end{proof}

Apply Lemma~\ref{calibration:lem:sc-cap-removal} initially using
\eqref{calibration:eq:sc-negative0}, and at the endpoint using
\eqref{calibration:eq:sc-negative-tau}. If
\[
 V_v=e^{2x_v},\quad \ell_v^{\,*}=\E V_v,\quad
 J_0=\E[V_zV_wH_0],\qquad m_v=\E e^{2y_v},
\]
then, for any required fixed $L$,
\begin{align}
 \ell_v^{\,*}&=Z_0(2e_{r_v})+O(N^{-L})\asymp N/H,
                                                       \label{calibration:eq:sc-lscale}\\
 J_0&=Z_0(a_j)+O(N^{-L})\asymp (N/H)^2e^{-2j}.
                                                       \label{calibration:eq:sc-Jscale}
\end{align}
These are initial-law estimates, obtained before microscopic signed
calibration. In particular, their use below introduces no assumption
about the unknown microscopic signed normalizer.

\subsection{Two-sided local-factor calibration}

In this subsection we approximate microscopic weights using independent local factors and the initial matrix, then compare their exact joint normalizing constants.
Let $U_z,U_w$ be the capped reference local factors supplied by
the local-factor theorem, with $\chi=2\sqrt2$. They are mutually independent
and their joint family is independent of the original initial-matrix
sigma field $\mathcal F_0^{\rm mat}$. They obey $U_v\le N^{A_*}$
for a fixed $A_*$. Their defining individual cap indicators remain
inside the separate factors, so this independence is exact.

\eqref{lf:eq:lfl-lf1} is written for $\Phi$; we denote its reconstruction
error by $\rho_{\rm LF}$. The normalized microscopic height
$\eta_{\rm m}$ corresponds exactly to the physical height
$\sqrt{1+\tau}\,\eta_{\rm m}$ at $v=\sqrt{1+\tau}z$.
The local-factor theorem permits this deterministic bounded
positive multiplier of $N^{-1}t^{-q}$, with the same height in
its actual and reference endpoints. We use that version here.
Let $E_v^0$ and $E_v^{\rm m}$ be,
respectively, the deterministic mean corrections in $x_v$ and $y_v$.
Define the exact deterministic factors
\begin{equation}
 k_v=e^{\chi(E_v^0-E_v^{\rm m})},\qquad U'_v=k_vU_v.
 \label{calibration:eq:sc-mean-factor}
\end{equation}
They are bounded above and below and preserve independence. In the
following identity, the local-factor event means the reconstruction event
intersected with the individual reference-cap event, so that the
capped $U_v$ agrees there with its uncapped reference factor. \eqref{lf:eq:lfl-lf1}
and \eqref{lf:eq:lfl-lf2} become
\begin{align}
 e^{2y_v}&=U'_vV_v\exp(O(\rho_{\rm LF}))
       \quad\hbox{on the local-factor event},                  \label{calibration:eq:sc-LF1}\\
 \frac{(\E U'_v)\ell_v^{\,*}}{m_v}
 &=1+o(t^{-P}).                                      \label{calibration:eq:sc-LF2}
\end{align}
The second formula is exact in its deterministic factors: multiplying
out $k_v$, $V_v=e^{\chi(\Phi_0-E_v^0)}$, and
$m_v=e^{-\chi E_v^{\rm m}}\E e^{\chi\Phi_\tau}$
cancels both mean corrections. They are not replaced by approximate
normalizers.

Let $\Omega$ intersect both local-factor reconstruction events, the two
individual reference-cap events, and the history-flow event.
It has complement probability $C_DN^{-D}$ for any fixed $D$
selected after the cap and moment costs. On this event,
\begin{equation}
 e^{2y_z+2y_w}H_\tau
    =U'_zU'_wV_zV_wH_0
                  \exp(O(\rho_{\rm LF}+\rho_{\rm hist})).
 \label{calibration:eq:sc-product}
\end{equation}
This is the two-sided reconstruction \eqref{calibration:eq:sc-LF1}. The bounded-test comparison \eqref{lf:eq:lfl-lf3}
cannot be used directly with $F=H_0$, because $H_0$ may be $N^{4A}$.

The weighted contributions from the exceptional event satisfy
\begin{align}
 \E[e^{2y_z+2y_w}H_\tau\ind_{\Omega^c}]
       &\le CN^{4A+8-D/2},                          \label{calibration:eq:sc-bad1}\\
 \E[U'_zU'_wV_zV_wH_0\ind_{\Omega^c}]
       &\le CN^{4A+2A_*+8-D/2}.                    \label{calibration:eq:sc-bad2}
\end{align}
These follow from H\"older with fourth moments of the positive
roots, the probability of $\Omega^c$, and the deterministic cap
envelopes. No independence from $\Omega$ is used.
After removing the event with these estimates, exact independence
gives
\[
 \E[U'_zU'_wV_zV_wH_0]
     =(\E U'_z)(\E U'_w)J_0 .
\]
It follows for any required fixed $L$ that
\begin{equation}
 S_\tau:=\E[e^{2y_z+2y_w}H_\tau]
   =(\E U'_z)(\E U'_w)J_0
      [1+O(\rho_{\rm LF}+\rho_{\rm hist})]+O(N^{-L}).
 \label{calibration:eq:sc-two-sided}
\end{equation}
By \eqref{calibration:eq:sc-positive}, \eqref{calibration:eq:sc-lscale}, and
\eqref{calibration:eq:sc-LF2}, $\E U'_v\asymp H$. Thus the main term in
\eqref{calibration:eq:sc-two-sided} satisfies
\begin{equation}
 (\E U'_z)(\E U'_w)J_0\asymp N^2e^{-2j}\ge cN.
 \label{calibration:eq:sc-micro-lower}
\end{equation}
These bounds already follow from the initial sparse-weight theorem,
scalar one-root normalization and \eqref{calibration:eq:sc-LF2}. The lower
bound therefore makes the absolute exceptional errors relative before the
microscopic signed normalizer is identified.

Apply Lemma~\ref{calibration:lem:sc-cap-removal} to the endpoint numerator
and both individual endpoint denominators. Equations
\eqref{calibration:eq:sc-two-sided}--\eqref{calibration:eq:sc-micro-lower},
\eqref{calibration:eq:sc-LF2}, and the analogous initial cap removals give
\[
 \frac{e^{2j}Z_{\rm m}(a_j)}
           {Z_{\rm m}(p)Z_{\rm m}(q_*)}
 =\frac{e^{2j}J_0}{\ell_z^{\,*}\ell_w^{\,*}}[1+o(t^{-P})]
 =\frac{e^{2j}Z_0(a_j)}{Z_0(p)Z_0(q_*)}[1+o(t^{-P})].
\]
This proves the first assertion of Theorem~\ref{calibration:thm:sc-main},
including its natural positive lower bound. The negative cap has
remained present throughout.

\subsection{Original-ensemble return and the parameter order}

In this subsection we transfer the normalizing-constant comparison to the original entry law and choose the parameters to preserve its relative precision.
Theorem~\ref{replacement:thm:rgfc-main} supplies, separately for
$a=a_j,p,q_*$,
\begin{equation}
 \frac{Z_{\rm target}(a)}{Z_{\rm m}(a)}
       =1+O(N^{-c}).
 \label{calibration:eq:sc-replacement}
\end{equation}
This is the consequence of applying its unnormalized differential
inequality first to $f=1$, using the deterministic two-sided cap
lower bound. It is not inferred from comparison of already
normalized laws. The coordinates and deterministic centers are
held fixed during interpolation. The natural mean corrections of the
original ensemble and the comparison matrix $M=X_\tau/\sqrt{1+\tau}$
differ by $O(\sigma)$; cap Lipschitz continuity
therefore restores them with a multiplicative $1+O(\sigma)$ error,
as in Lemma~\ref{replacement:lem:rgfc-centers}.
The exact variance conversion \eqref{calibration:eq:sc-variance} and exact
mean factors \eqref{calibration:eq:sc-mean-factor} have already been included.
Only after obtaining each of these relative ratios of normalizing constants do we
combine \eqref{calibration:eq:sc-replacement} for the three sources
to prove the original-ensemble claim.

For each original individual positive normalizer,
Lemma~\ref{calibration:lem:sc-cap-removal} with the negative factor omitted
uses only the scalar upper tail and positive moments. Theorem~\ref{scalar:thm:sr-regularization} applies to the exact
density-centered normalizers
\[
 Z_{L,v}^{\eta,\Phi}:=\E_L e^{\chi\Phi^L(v,\eta)}.
\]
The bounded deterministic mean correction is smooth at height
zero in the fixed bulk, so
$E_{\kappa_L}(v,\eta)-E_{\kappa_L}(v,0)=O(\eta)$ uniformly
over the common moment class. Consequently, for
$\eta=\eta_{\rm m}$,
\[
 \frac{m_{L,v}^{\eta}}{m_{L,v}^0}
 =e^{-\chi[E_{\kappa_L}(v,\eta)-E_{\kappa_L}(v,0)]}
      \frac{Z_{L,v}^{\eta,\Phi}}{Z_{L,v}^{0,\Phi}}
 =1+O(t^{-q}+N^{-\alpha}+\eta).
\]
The same argument applies to the endpoint law of the comparison matrix $M$.
At a common height $u\ge0$, mean corrections cancel exactly:
\[
 \frac{e^{\chi[\Phi^L(v,u)-E_{\kappa_L}(v,u)]}}{m_{L,v}^{u}}
 =\frac{e^{\chi\Phi^L(v,u)}}{Z_{L,v}^{u,\Phi}}.
\]
If the numerator remains at height $\eta$ while its denominator
has been replaced by $m_{L,v}^{0}$, simultaneous conversion to
density-centered notation retains the deterministic factor
\[
 e^{-\chi[E_{\kappa_L}(v,\eta)-E_{\kappa_L}(v,0)]}=1+O(\eta).
\]
Thus normalization at a common height is unchanged by the deterministic mean correction, and a change
of height leaves precisely the deterministic factor displayed above.
No isolated mean-corrected denominator is identified with a
density-centered denominator.
Zero determinants in a discrete ensemble have positive-root weight
zero and cause no problem. No claim about their logarithmic lower
tails is made.

Fix the desired precision $P$ first, allowing extra fixed powers
for later smooth-test derivative costs. Choose $A$ sufficiently
large for the scalar cap errors and fix the local-factor cap envelope
$A_*$. After these caps and all fixed charge and moment costs are
known, fix the required failure and high-moment orders. Then choose
the reconstruction exponents in the local-factor theorem's proved order, with
$q>P+2$, choosing $b$ large enough for both the precision of the local-factor reconstruction and \eqref{signed:eq:ss-rate}
before fixing the later windows. All logarithmic exponents are thus
chosen with the failure exponents already fixed. A sufficient bound for the total relative error is
\[
 O(\rho_{\rm LF}+\rho_{\rm hist}+t^{-q}+r_N+N^{-c})
       +\hbox{the specified cap and exceptional remainders},
\]
which is $o(t^{-P})$.

\begin{remark}[Scope of the normalizing-constant comparison]
The normalizing constant for the signed coefficient vector is calibrated before any rare
barrier test is inserted. An upper comparison for a compensated
barrier still needs the bounded test, its summed derivative norm,
the same negative cap in source and compensation, and weighted
exceptional estimates. Gaussian barrier bounds and the terminal
core estimate are additional inputs. This theorem neither
replaces them nor asserts unrestricted two-root convergence.
\end{remark}

\endgroup

\begingroup
\section{One-root comparison and the endpoint cutoff}\label{chap:oneroot}

In this section, at a fixed barrier, we compare determinant-weighted histories with Gaussian histories and show that the terminal cutoff loses negligible expected mass.
Indeed, we compare the actual determinant-weighted law with the corresponding
Gaussian law. The comparison includes the full integer barrier history,
the prescribed coarse observations, and bounded Lipschitz terminal marks. The proof first separates the local
endpoint increment from the initial matrix. Calibration at three nearby
positive charges controls both signs of this increment; capped signed
sources control the initial increment. Together they give an endpoint
displacement of order \(\log^2\log N\) with error under the determinant-weighted law smaller than any prescribed
fixed inverse power of the logarithm of the matrix size. Independent Gaussian estimates then
complete the integer history and prove the cutoff at
\(d=(\tfrac12\log N)^{1/3}\). The mean-corrected endpoint is retained
throughout, including in the cutoff; exact determinant zeros carry
zero positive-root weight.

\subsection{Coordinates, measures, and Gaussian input}\label{oneroot:coordinates-measures-and-gaussian-input}

We next fix the corrected fields, weighted laws, and Gaussian barrier estimates needed to compare one spatial point and its observations.
Let \(A_N=N^{-1/2}(\xi_{ab})_{a,b=1}^N\), where the complex entries
are independent and identically distributed, with
\[
 \mathbb E\xi=\mathbb E\xi^2=0,\qquad
 \mathbb E|\xi|^2=1,\qquad
 \mathbb E|\xi|^r\le C_r\quad(r\ge1),\qquad
 \kappa=\mathbb E|\xi|^4-2.
\]
The entry law may vary with \(N\) under common moment bounds;
all conclusions are uniform when the parameters and bulk compact are fixed.
Let \(K\Subset\{|z|<1\}\) be contained in a fixed compact disk and
let \(\rho^z\) be the exact deterministic variance-one Hermitization density.
Define
\[
 \Phi_N(z,\eta)=\tfrac12\log\det\big((A_N-z)(A_N-z)^*+\eta^2I\big)
 -\tfrac N2\int_{\mathbb R}\log(x^2+\eta^2)\rho^z(x)\,dx.
\]
At zero height its positive exponential is defined to be zero on exact
determinant zeros. The scalar transform \(m^z\) is the Stieltjes transform
of \(\rho^z\). Fix an integer \(k_0\ge1\) and a real barrier height \(B\).
For real \(s>0\), put
\[
 u_s=\tanh s,\qquad p_z(u)=u(1-u|z|^2),\qquad
 \eta_s(z)=(1-u_s)\sqrt{u_s^{-1}-|z|^2}.
\]
The clock satisfies \(m^z(i\eta_s)^4=p_z(u_s)^2\).
Put
\[
 t=\tfrac12\log N,\qquad n=\lfloor t\rfloor,\qquad
 \chi=2\sqrt2,
\]
\[
 \begin{aligned}
 Y_N(z)&=\sqrt2\Phi_N(z,0),\\
 Y_N^0(z)&=Y_N(z)+\frac{\sqrt2\kappa}{4}(1-|z|^2)^2,\\
 X_{N,k}(z)&=\sqrt2[\Phi_N(z,\eta_k(z))-E_\kappa(z,\eta_k(z))].
 \end{aligned}
\]
where \(E_\kappa=-\kappa m^4/4\) is the deterministic finite correction
and \(\eta_k\) is the exact clock. Thus \(Y_N\) is the raw uncorrected
endpoint, while \(Y_N^0\) is the canonically corrected endpoint used in
the comparison and cutoff sections. Since \(m^z(i0)^4=(1-|z|^2)^2\), this is the
continuous zero-height version of the correction in \(X_{N,k}\). The
difference \(Y_N^0-Y_N\) is uniformly bounded on the bulk compact; it is
kept in the barrier/cutoff conventions and cancels exactly from the
normalized root. The raw endpoint mark in the main theorem is recovered when we pass to the marked measures: its displacement is \(O(t^{-1/2})\), and the bounded-Lipschitz error vanishes after total-mass tightness has been proved. 
Write
\[
 W_N^0(z)=\frac{e^{2Y_N^0(z)}}{\mathbb Ee^{2Y_N^0(z)}},
 \qquad dQ_{N,z}^{\rm det}=W_N^0(z)dP,
\]
\[
 A_{N,j}^B(z)=\{X_{N,k}(z)\le2k+B:\ k_0\le k\le j\}.
\]
The probability measure just defined is the determinant-weighted, or rooted, law. At a zero determinant its weight is zero, so the value assigned
to the unregularized field inside a rooted test is immaterial.

The centered Gaussian comparator has covariance
\begin{align*}
 C_\kappa(s,z;q,w)&=-\tfrac12\log D(u_s,z;u_q,w)
       +\tfrac\kappa2p_z(u_s)p_w(u_q),\\
 D(u,z;v,w)&=1-uv-uv\{(1-u)|z|^2+(1-v)|w|^2\}+uv|z-w|^2.
\end{align*}
The positivity of this complete covariance, including \(\kappa=-1\), is
proved by the feature expansion \eqref{gaussian:gp2:features}; no negative-variance
independent summand is introduced. The parameter \(\kappa\) is the original
entry-law fourth cumulant. The Gaussian root is at the real time \(t\),
while its barrier is observed at integer times:
\[
 dQ_{t,z}^G=e^{2G_t(z)-2\operatorname{Var}G_t(z)}dP,
 \qquad A_j^B=\{G_k(z)\le2k+B:\ k_0\le k\le j\}.
\]

The needed Gaussian input is, for \(j\asymp t\), \(1\le R\le j\), and
any fixed \(D\),
\[
 Q_{t,z}^G\{A_j^B,\ 0\le B+2j-G_j\le R\}
 \le \frac{C_D R(1+\sqrt{\log(2+j)})(1+R+\sqrt{\log(2+j)})}
              {j^{3/2}}+C_Dj^{-D}.
 \tag{OR1}\label{oneroot:eq:or1}
\]
Also, if \(t-j=O(\log t)\),
\[
 Q_{t,z}^G(A_j^B\setminus A_n^B)=o(t^{-1/2}).
 \tag{OR2}\label{oneroot:eq:or2}
\]

Lemma~\ref{gaussian:gpm:strips} and its following completion argument
in Section~\ref{gaussian:gpm:endpoint}
verify \eqref{oneroot:eq:or1} by endpoint regression, a bounded covariance comparison
with a Brownian bridge, an \(O(\sqrt{\log j})\) auxiliary-noise maximum,
and bridge reflection. \eqref{oneroot:eq:or2} follows by splitting at a terminal increment
of size \(C\log t\), applying \eqref{oneroot:eq:or1} to the small-deficit part, and using
Gaussian increment tails. Neither statement uses a Markov property of
the exact Gaussian scale process. These inputs are not conclusions
inferred from the matrix full-history comparison in Theorem~\ref{oneroot:thm:one-root}.

\begin{theorem}[One-root comparison and the endpoint cutoff]\label{oneroot:thm:one-root}
Under the preceding hypotheses, the following hold uniformly for \(z\in K\).
For every prescribed fixed \(P>1/2\), there are fixed exponents \(b,a>0\)
and a fixed \(C_K\) such that
\[
 J=\lfloor t-(b/2+a)\log t-C_K\rfloor,
 \qquad
 Q_{N,z}^{\rm det}(|Y_N^0-X_{N,J}|>C(1+\log t)^2)=o(t^{-P}).
\]
Let \(S_N\) be a fixed finite prefix of the root's integer history and
let \(U_N\) be a fixed finite list of centered coarse coordinates
\(\Phi_N(w_\ell,a_\ell)-E_\kappa(w_\ell,a_\ell)\), with
\(a_\ell\in[a_0,1]\), \(a_0>0\), and \(w_\ell\) in a fixed bulk compact.
There is no separation assumption between these points and \(z\).
Let \(S_G,U_G\) be the corresponding Gaussian coordinates, with
\(U_{G,\ell}=G_{s_\ell}(w_\ell)/\sqrt2\) and
\(\eta_{s_\ell}(w_\ell)=a_\ell\).
For every fixed \(H\in C_b^{M_H}\), where \(M_H\ge\max\{m,5\}\) and \(m\) is the fixed Stein derivative order from Corollary~\ref{signed:cor:ss-eight}, and every bounded Lipschitz \(h\),
\begin{align*}
\sqrt t\,\big|&\mathbb E_{Q_{N,z}^{\rm det}}
 [H(S_N,U_N)\mathbf1_{A_{N,n}^B}
 h((2t-Y_N^0)/\sqrt t)]\\
 &-\mathbb E_{Q_{t,z}^G}
 [H(S_G,U_G)\mathbf1_{A_n^B}
 h((2t-G_t)/\sqrt t)]\big|\longrightarrow0.
\end{align*}
Finally, set \(d=t^{1/3}\) and
\[
 \mu_N^B(dz)=\sqrt t\,W_N^0(z)\mathbf1_{A_{N,n}^B(z)}\,d^2z,
 \quad
 \widetilde\mu_N^B(dz)=\mu_N^B(dz)\mathbf1_{\{Y_N^0(z)\le2t-d\}}.
\]
Then \(\mathbb E\|\mu_N^B-\widetilde\mu_N^B\|_{\rm TV}\to0\).
For the cut density \(U_N^{\rm cut}(z)=W_N^0(z)\mathbf1_{A_{N,n}^B(z)}
\mathbf1_{\{Y_N^0(z)\le2t-d\}}\), each fixed \(C>0\) satisfies
\[
 t\int_K\int_{K\cap\{|w-z|^2\le CN^{-1}e^d\}}
   \mathbb E[U_N^{\rm cut}(z)U_N^{\rm cut}(w)]\,d^2w\,d^2z\le C'te^{-d}.
\]
\end{theorem}

The coarse sites in the theorem may lie outside \(K\), within their
prescribed fixed bulk compact, and may repeat or coincide with the root.
Their heights include one. Their matrix coordinates have no leading
\(\sqrt2\), which explains the Gaussian factor \(1/\sqrt2\).
These finite coarse observations are retained in every comparison below.

\paragraph{Proof structure.}
Sections~\ref{oneroot:scalar-facts-including-the-previously-delicate-variable-height}
through \ref{oneroot:completing-the-final-integer-barrier} establish the endpoint
and history assertions, including the mark;
Section~\ref{oneroot:strong-cutoff-repair-at-dt13} proves the cutoff assertions.
Section~\ref{oneroot:complete-parameter-order} records the
order of parameter choices. The local-factor, signed-source, replacement,
and history-flow theorems are used in their quantitative forms from Sections~\ref{chap:lf}, \ref{chap:signed}, \ref{chap:replacement}
and~\ref{chap:calibration}. In particular, the local-factor input includes
Theorem~\ref{lf:thm:lfl-main}, Proposition~\ref{lf:prop:lfl-uniform-charges}
and the full-space weighted proof in Section~\ref{lf:sec:lfl-weights};
the history input includes Corollary~\ref{calibration:cor:sc-history-log}
and its proof. The Gaussian estimates
\eqref{oneroot:eq:or1}--\eqref{oneroot:eq:or2} are supplied independently by Lemma~\ref{gaussian:gpm:strips} and its following terminal-completion proof in Section~\ref{gaussian:gpm:endpoint}; they do not
assume the matrix full-history theorem. The fixed-barrier spatial measure limit, mixed-time Cauchy statement,
removal of the barrier, and the main joint convergence conclusion are proved
separately.

\subsection{Scalar moments and variable heights}\label{oneroot:scalar-facts-including-the-previously-delicate-variable-height}

This subsection controls positive determinant moments across the regularization heights needed to connect the retained history with the terminal field.
Theorem~\ref{scalar:thm:sr-sharp}, Proposition~\ref{scalar:prop:sr-moments}
and Theorem~\ref{scalar:thm:sr-regularization}, using the
one-point theorem \cite{CL}, prove, uniformly for every
fixed charge \(r\) in a compact subinterval of \([2,\infty)\),
\[
 Z_r(0):=\mathbb Ee^{r\Phi(z,0)}\asymp N^{r^2/8},
 \qquad \mathbb Ee^{r\Phi(z,\eta)}\le C_rN^{r(r+2)/8},
 \tag{OR3}\label{oneroot:eq:or3}
\]
and, for \(\eta_m=t^{-q_m}/N\),
\[
 Z_r(\eta_m)/Z_r(0)=1+O(t^{-q_m}+N^{-\alpha}),
 \qquad
 \left\|Q_{r,\eta_m}-Q_{r,0}\right\|_{\rm TV}
 \le C(t^{-q_m}+N^{-\alpha}).
 \tag{OR4}\label{oneroot:eq:or4}
\]

For endpoint displacement we also need regularizers above the
microscopic scale. The scalar proof gives the following upper ratio
at every deterministic height tending to zero:
\[
 Z_r(\eta)/Z_r(0)\le1+C_r(\eta+N^{-\alpha})\le C_r.
 \tag{OR5}\label{oneroot:eq:or5}
\]

Here is the essential derivation, so the range is explicit.
Gaussian-minor expansion and Jensen, valid for \(r\ge2\), yield
\[
 \mathbb Ee^{r\Phi(z,\eta;X)}
 \le\mathbb Ee^{r\Phi(w,0;\widetilde X_s)},
 \quad s=\eta/\Im m^z(i\eta),\quad w=z/\sqrt{1+s},
\]
where \(G\) has independent standard circular complex entries of
variance one and \(\widetilde X_s=(X+\sqrt{s/N}G)/\sqrt{1+s}\).
The full minor and Jensen proof is Lemma {\hyperref[scalar:lem:sr-minor]{\ref*{scalar:lem:sr-minor}}}
of Section~\ref{chap:scalar}. Write
\(I_z(\eta)=\int\log(x^2+\eta^2)\rho^z(x)\,dx\).
The deterministic density-centering prefactor cancels exactly because
\[
 I_z(\eta)=\eta^2/s+\log(1+s)+|z|^2/(1+s)-1.
\]

The normalized convolution has common moment bounds, and its fourth
cumulant is \(\kappa/(1+s)^2\). The sharp raw one-root formula, uniform
under these moment bounds and in the fixed charge interval, differs from the target
formula by a bounded factor, and by \(1+O(s+N^{-\alpha})\) when
\(s\to0\). Since \(s=O(\eta)\) in the bulk, \eqref{oneroot:eq:or5} follows. The restriction
\(\eta\le c/N\) is needed for the small total-variation error in \eqref{oneroot:eq:or4},
not for the constant upper ratio \eqref{oneroot:eq:or5}.

Thus \eqref{oneroot:eq:or5} applies at every \(\eta_k\), \(J\le k\le n\), since
\(\eta_J\asymp t^{b+2a}/N\to0\). This is a direct scalar comparison, with all logarithmic exponents
fixed before \(N\to\infty\). It does not require moving the lower
endpoint in a local-factor reconstruction theorem.

Positive two-sided cap removal is uniform over these regularizers. For
the upper tail use a higher fixed positive moment in \eqref{oneroot:eq:or3}; for the
lower side of a positive-root cap use its exponentially small weight. For a bounded
monotone cap which is not assumed contractive, the capped upper-shell
term is bounded by its polynomial cap envelope times a higher
positive-moment Markov bound. Prescribe a cap-shell error exponent
\(D_{\rm cap}\) first, then choose the fixed positive moment charges and
cap size \(A\) so this error is at most \(N^{-D_{\rm cap}}\). Choose
the independent overwhelming-event failure exponent afterward, beyond
all fixed cap and moment costs. This is a positive-weight argument and does not integrate an additive
exceptional-probability tail against an unbounded exponential.
For every signed-source or replacement invocation below, choose the
common cap to be odd, smooth, nondecreasing and \(1\)-Lipschitz,
equal to the identity on \([-A\log N,A\log N]\) and bounded by
\(2A\log N\). The preceding broader scalar observation does not
enlarge the cap class in either of those comparison theorems.

\subsection{One-root factorization on the full coupling space}\label{oneroot:arbitrarily-accurate-one-root-full-space-lf-factorization}

This subsection separates the local increment from the initial matrix on a common probability space before normalizing the determinant weight.
We first choose a matched Gaussian-divisible comparison matrix, called the surrogate below, by the explicit
construction in Proposition~\ref{replacement:prop:rgfc-age} and its proof using Lemma~\ref{replacement:lem:rgfc-matchlaw}, so
\[
 Y=X_\tau/\sqrt{1+\tau},\qquad
 X_\tau=X_0+N^{-1/2}B_\tau,
 \qquad \tau=t^b/N,
\]
and \(X_0\) has Gaussian age exactly \(\sigma=N^{-1/200}\). Matching is
imposed on \(Y\), not on \(X_0\).

To use the stated two-root local-factor theorem for this single root, we adjoin a
deterministic auxiliary point at distance comparable to \(N^{-1/4}\)
inside the enlarged bulk compact, apply that theorem, and retain the
desired marginal. Its local reference family remains independent of the
entire initial matrix. This construction is for each deterministic root;
no simultaneous coupling over a continuum of spatial points is asserted.

At the physical site \(v=\sqrt{1+\tau}\,z\), the one-root reconstruction is
\[
 \Phi_\tau(v,\eta_m^{\rm phys})=C_0+S+o(t^{-P_*}),
 \qquad C_0=\Phi_0(v,B_v\tau),\quad S=L_1+L_2,
 \tag{OR6}\label{oneroot:eq:or6}
\]
with overwhelming probability. The variables \(L_1,L_2\) are the reference
increments in \eqref{lf:eq:lfl-L1}--\eqref{lf:eq:lfl-L2} at this root. Here
\(\eta_m^{\rm phys}=\sqrt{1+\tau}\,\eta_m\): the point and height
are rescaled together. The local-factor theorem permits this bounded
deterministic height multiplier, with the same height at its actual and
reference endpoints. The reference variable \(S\) is independent of
the entire initial matrix.

Fix once and for all \(0<p<\min\{\chi-2,1/2\}\), for example \(p=1/4\).
Use one common cap event \(\mathcal K\), depending only on the local
reference variables, for the three charges
\[
 r\in\{\chi-p,\chi,\chi+p\}\subset(2,\infty),
 \qquad U_r=e^{rS}\mathbf1_{\mathcal K}.
\]
Proposition~\ref{lf:prop:lfl-uniform-charges} provides a common local
cap and calibration over a fixed positive charge interval. Its
reconstruction is charge-independent, and multiplication by these
three charges changes only fixed constants. The proof extends the
weighted estimate \eqref{lf:eq:lfl-L1-weight} to these charges and gives
\[
 \mathbb E\left|e^{r\Phi_\tau(v,\eta_m^{\rm phys})}
                    -U_r e^{rC_0}\right|
 =o(t^{-P_*})\mathbb Ee^{r\Phi_\tau(v,\eta_m^{\rm phys})},
 \tag{OR7}\label{oneroot:eq:or7}
\]
with the failure exponent chosen after the three charges and the local
cap. The event \(\mathcal K\) is never enlarged by the coupling event.
Integrating \eqref{oneroot:eq:or7} first compares the two unnormalized masses; exact
independence then factors the reference mass and calibrates each
normalizer. Thus the denominator in \eqref{oneroot:eq:or8} is justified before division.

At \(r=\chi\), normalize both sides of \eqref{oneroot:eq:or7}. On the \textbf{full coupling
space}, not only after conditional expectation,
\[
 \left\|
 \frac{e^{\chi\Phi_\tau(v,\eta_m^{\rm phys})}}
      {\mathbb Ee^{\chi\Phi_\tau(v,\eta_m^{\rm phys})}}
 -\frac{U_\chi e^{\chi C_0}}
       {\mathbb EU_\chi\,\mathbb Ee^{\chi C_0}}
 \right\|_{L^1}=o(t^{-P_*}).
 \tag{OR8}\label{oneroot:eq:or8}
\]
\eqref{oneroot:eq:or4} changes the first density to the actual determinant density at an
additional \(O(t^{-q_m}+N^{-\alpha})\). The second density gives exactly
the product law \(Q_U\otimes Q_V\) to the local reference family and the
initial matrix. Here \(Q_U\) has density \(U_\chi/\mathbb E U_\chi\)
on the local reference family, and \(Q_V\) has density
\(e^{\chi C_0}/\mathbb E e^{\chi C_0}\) on the initial matrix.
This statement on the full coupling space follows from the already
proved pointwise local-factor reconstruction and weighted exceptional estimate; it is stronger
than, and supplies a proof of, the restriction to initial-matrix tests.

\subsection{Local and initial increment tails}\label{oneroot:two-sided-local-and-initial-increment-tails}

The purpose of this subsection is to control both tails of the local and initial increments so that neither can move the terminal field too far.
The initial root height is \(B_v\tau\), corresponding to the exact clock
\[
 \ell_v=t-(b/2)\log t+O_K(1).
\]

Corollary~\ref{signed:cor:ss-normalizer}, applied at the three positive
charges and followed by positive cap removal, gives:
\[
 \mathbb Ee^{rC_0}\asymp(B_v\tau)^{-r^2/8}.
\]
Together with \eqref{oneroot:eq:or3}--\eqref{oneroot:eq:or4} and \eqref{oneroot:eq:or7}, this proves
\[
 \mathbb EU_r\asymp(N\tau)^{r^2/8}.
 \tag{OR9}\label{oneroot:eq:or9}
\]

Only a fixed compact interval of positive charges is needed. The initial
coefficients \(r/\sqrt2\) satisfy the bound of six in the signed-source
estimate; passing from the mean-corrected \(\sqrt2\)-scaled field to \(C_0=\Phi_0\)
introduces only a bounded deterministic factor, using the initial law's own fourth cumulant. 
Consequently
\[
 \mathbb E_{Q_U}e^{\pm pS}
 =\frac{\mathbb EU_{\chi\pm p}}{\mathbb EU_\chi}
 \le Ct^{C_p},
 \qquad Q_U(|S|>R)\le Ct^{C_p}e^{-pR}.
 \tag{OR10}\label{oneroot:eq:or10}
\]
The equalities require the same cap at all three charges. Both signs
of the rooted increment are controlled by positive charges
\(\chi-p\) and \(\chi+p\); no uncapped negative local exponential
moment is used.

Set
\[
 J=\left\lfloor t-(b/2+a)\log t-C_K\right\rfloor,
 \qquad C_J=\Phi_0(v,\eta_J(v)),
\]
with \(C_K\) chosen so that \(\ell_v-J\ge a\log t\) uniformly. The
initial same-site clock difference is \(a\log t+O(1)\).

The initial increment requires a different argument because both
coordinates belong to the same matrix. To control the two-sided tail
of \(C_0-C_J\), first keep both coordinates capped. Under the initial critical root the relevant source rays in
the unit coordinates are
\[
 (2\pm p/\sqrt2)e_{\ell_v}\mp(p/\sqrt2)e_J.
 \tag{OR11}\label{oneroot:eq:or11}
\]
They have support at most two and absolute charge at most
\(2+\sqrt2p<6\). Corollary~\ref{signed:cor:ss-normalizer} computes their normalizer ratios. The exact
same-site covariance gives
\[
 \operatorname{Var}[(G_{\ell_v}-G_J)/\sqrt2]
       =O(1+\log t),
\]
and the mean of this difference under the root at \(\ell_v\) is also
\(O(1+\log t)\). Therefore the two signed exponential moments of the
\textbf{capped difference} are at most \(Ct^{C_p}\).

The positive initial root cap can be removed in total variation using
\eqref{oneroot:eq:or3} and the source normalization. Under that positive root, the event
that either of these two mesoscopic coordinates lies outside the cap
identity interval is polynomially negligible: combine the scalar cap
tail with Cauchy--Schwarz and the positive second moment. 
Thus
\[
 Q_V(|C_0-C_J|>R)\le Ct^{C_p}e^{-p(R-O(1))}+N^{-D}.
 \tag{OR12}\label{oneroot:eq:or12}
\]

This argument never removes a negative cap inside an exponential moment.
It uses the capped signed ray to bound a bounded tail event, and only
then removes the coordinate caps on an event of controlled positive-root
probability. Directly asserting an uncapped moment with coefficient
\(-p\) at \(C_J\) would be an unjustified inverse-determinant step.

\subsection{Endpoint displacement under the surrogate and target laws}\label{oneroot:endpoint-displacement-under-the-surrogate-then-the-target}

In this subsection we bound the gap between the terminal field and its last retained observation, first for the comparison matrix and then for the original law.
Indeed, we first combine the independent local increment with the initial
capped-increment tail under the surrogate root. We then pass from the
regularized weight at microscopic height to the actual determinant root. Return to the target
law uses a smooth test of positive-height coordinates; the raw endpoint
estimate is applied separately in that law.

Corollary~\ref{calibration:cor:sc-history-log}, with its proof, with
\(L=b+2a>1\) and \(\eta_J\asymp\tau t^{2a}\), gives
the actual coordinate error
\[
 \rho_{\rm hist}\le C_D[t^{1-a}+t^{1-2a}
                         +t^{2-b-4a}+N^{-c}],
 \tag{OR13}\label{oneroot:eq:or13}
\]
with overwhelming probability. The first term is the
martingale square-root term and must be retained. For sufficiently large
fixed \(a\), this tends to zero. Combining \eqref{oneroot:eq:or6} and \eqref{oneroot:eq:or13},
\[
 \Phi_\tau(v,\eta_m^{\rm phys})
       -\Phi_\tau(v,\eta_J^{\rm phys})
   =S+C_0-C_J+o(1).
\]
Positive second moments and the prior one-root normalizer turn the
exceptional probability into arbitrarily small polynomial actual-root
mass, after the failure order is fixed sufficiently large. The
full-space comparison \eqref{oneroot:eq:or8} then permits the local and initial tail
bounds \eqref{oneroot:eq:or10} and \eqref{oneroot:eq:or12} to be combined. For
\(R=C(1+\log t)^2\),
\[
 Q_Y^{\rm det}\left(
 |\Phi_Y(z,\eta_m)-\Phi_Y(z,\eta_J)|>R\right)
 =o(t^{-P})
 \tag{OR14}\label{oneroot:eq:or14}
\]
for every prescribed fixed \(P\), once the local-factor and regularization accuracies have been
selected accordingly. The exponentially small terms
\(t^Ce^{-c(\log t)^2}\) are smaller than every fixed inverse power of \(t\).

To handle the raw endpoint, put \(d_m=\Phi(z,\eta_m)-\Phi(z,0)\) on the
nonzero-determinant support. Determinant monotonicity gives
\(d_m\ge-Ct^{-q_m}\), while the rooted moment bound gives
\[
 \mathbb E_{Q^{\rm det}}e^{\chi d_m}
 =\frac{\mathbb E[e^{\chi\Phi(z,\eta_m)}\mathbf1_{\{\det(X-z)\ne0\}}]}{Z_\chi(0)}
 \le Z_\chi(\eta_m)/Z_\chi(0)\le C.
\]

The root assigns zero mass to exact determinant zeros, which can have
positive unrooted probability for discrete entry laws. Thus the
displayed bound is an inequality in that generality; no expression
\(0\cdot\infty\) or absolute continuity in the reverse direction is
used. Therefore \(Q^{\rm det}(d_m>R)\le Ce^{-\chi R}\). Together with
\eqref{oneroot:eq:or14} this yields the two-sided raw displacement for the surrogate. The
bounded deterministic correction gives the same displacement estimate
for the surrogate law's corresponding corrected endpoint.

For the original ensemble, transfer only the bounded smooth tail test of
the two \textbf{regularized} coordinates \(\Phi(\eta_m)-\Phi(\eta_J)\),
using transition width one. Its \(\mathcal D_5\) is bounded
independently of the translating threshold. Theorem~\ref{replacement:thm:rgfc-main}, positive cap removal, and \eqref{oneroot:eq:or4}
give a polynomial replacement error plus \(O(t^{-q_m})\).
Lemma~\ref{replacement:lem:rgfc-centers} restores the natural
deterministic centers after comparison at frozen centers. To remove
the coordinate caps from this bounded test, the microscopic root's
lower cap side is controlled by its small positive weight, and its upper
side by a higher positive moment. The other coordinate is mesoscopic;
its cap event is controlled by the scalar tail under the unweighted law
and a positive-root Cauchy--Schwarz estimate. No microscopic lower-tail probability
estimate is required.
Both cutoffs exceed \(N^{-1-\delta}\) eventually. Apply the same
comparison between the unregularized and regularized weights directly in the target law. The bounded finite
mean correction between \(Y_N^0\) and \(X_{N,J}\) is absorbed into the threshold. 
Thus
\[
 Q_{N,z}^{\rm det}\{|Y_N^0-X_{N,J}|>C(1+\log t)^2\}
       =o(t^{-P}).
\tag{OR15}\label{oneroot:eq:or15}
\]

This proves the arbitrary fixed endpoint precision in
Theorem~\ref{oneroot:thm:one-root}, while every differentiated
coordinate remains at positive height.

\subsection{Retained history with fixed coarse observations}\label{oneroot:retained-history-comparison-with-fixed-coarse-observations}

In this subsection we compare the retained barrier history while keeping the prescribed coarse observations in the same expectation.
First, we let \(H\in C_b^{M_H}\), with bounded derivatives through \(M_H\ge\max\{m,5\}\), be a fixed function of a fixed early prefix and
finitely many coarse fields \(\Phi(w_\ell,a_\ell)-E_\kappa\), with
\(a_\ell\in[a_0,1]\), all spatial points in fixed bulk compacts. No
separation from the variable root is required; the prescribed list may
include repeats and sites outside \(K\). Smooth each retained barrier
at width \(h=t^{-s}\). There are \(O(t)\) factors, so summing over
all derivative coordinates gives the following bound for a product
test \(F\), including \(H\), in the summed derivative norm \eqref{signed:eq:ss-norm}:
\[
 \|F\|_\infty\le C_H,
 \qquad \mathcal D_r(F)\le C_{H,r}t^{r(1+s)}.
 \tag{OR16}\label{oneroot:eq:or16}
\]

The following five steps compare these bounded smooth expectations.
Every stated error is absolute under a normalized rooted law. The
bounds account for the full sum of derivatives in \eqref{oneroot:eq:or16}
before the smoothing is removed.

\begin{enumerate}
\setlength{\itemsep}{0pt}\setlength{\parskip}{0pt}
\item
  \textbf{Target to surrogate.} Theorem~\ref{replacement:thm:rgfc-main} for the positive-root source, together with the
  uniform scalar proxy and cap removal, gives error
  \(O(N^{-c}\mathcal D_5(F))+o(t^{-P_*})\|F\|_\infty\).
\item
  \textbf{Exact normalization and history flow.} Write the surrogate as
  \(X_\tau/\sqrt{1+\tau}\). The centered field identity sends every
  coordinate to its physical point/height exactly, with no \(N\tau\)
  error. Keep the root history at the physical site
  \(v=\sqrt{1+\tau}z\), and move its clocks to \(\eta_k(v)\). The common
  history error is \eqref{oneroot:eq:or13}, so its test cost is
  \(C\rho_{\rm hist}\mathcal D_1(F)\). Each fixed coarse coordinate has
  a positive height and incurs a polynomially small flow error. The
  history-flow estimate permits any fixed upper height, so it includes
  the physical coarse heights slightly above one. At the initial
  matrix, return each such height by \(O(\tau)\) to its original
  value before applying Corollary~\ref{signed:cor:ss-ou}. At a cutoff bounded below by
  \(a_0/2\), the centered height derivative and bounded mean return
  give a polynomially small error. This preserves height-one
  observations while keeping the signed list within height one.
\item
  \textbf{Restriction to the initial matrix.} The flowed test is
  measurable with respect to the entire initial matrix, including its
  coarse observations. \eqref{oneroot:eq:or8} and independence replace the actual root by
  \(e^{\chi C_0}/\mathbb Ee^{\chi C_0}\), with error
  \(o(t^{-P_*})\|F\|_\infty\), and no derivative loss.
\item
  \textbf{Initial matrix to Gaussian.} Apply
  Corollary~\ref{signed:cor:ss-ou}
  with the external root at \(\ell_v\), the fine history at the same
  physical site \(v\), and the fixed coarse list. Positive cap removal
  is controlled as above. Its error is at most 
\[
  C_b[t^{2-b/2}+t^{17-b}+N^{-c}]\mathcal D_m(F),
  \tag{OR17}\label{oneroot:eq:or17}
\] 
where the fixed derivative order \(m\) precedes \(b\).
\item
  \textbf{Gaussian parameter return.} Under the exact covariance, moving
  the external root from \(\ell_v\) to \(t\) changes every retained
  history mean by \(O(e^{-2(\ell_v-J)})=O(t^{-2a})\), as in
  \eqref{gaussian:gpm:prefix-decay}. For a fixed coarse
  observation the corresponding change is polynomially small. Returning
  the physical spatial points to the normalized ones, and restoring the
  target fourth cumulant, changes covariance and mean entries by
  \(O(\tau+\sigma)\). Gaussian interpolation costs
  \(O(\tau+\sigma)(\mathcal D_1(F)+\mathcal D_2(F))\); no inverse
  covariance is used. At the same fine site the exact kernel's radial
  derivative is uniformly bounded, and at a positive coarse height its
  stability denominator is bounded below.
\end{enumerate}

The physical root and its fine history stay at the same point until
the Gaussian stage. Their common spatial return is then controlled by
the same-site covariance kernel. This also keeps the lowest-height
root inside the exact normalizer throughout, instead of paying its
spatial displacement through a random field gradient.

The Gaussian upper and lower barrier sandwiches differ on coordinate
strips. For exact integer clock \(k\ge k_0\), the variance is comparable
to \(k\), so their union has probability at most
\[
 C h\sum_{k=k_0}^J k^{-1/2}\le C h\sqrt t.
\]
This bound is uniform in the strip centers and rooted means. 
For a signed insertion \(H\), the error from removing the smoothing is bounded by
\(\|H\|_\infty\) times the same strip probability. Smooth upper
and lower barrier products therefore transfer the hard barriers
without a matrix anti-concentration input. With the parameter order
in Section~\ref{oneroot:complete-parameter-order}, the retained
comparison has absolute error \(o(t^{-P})\).

\subsection{Translating endpoint strips}\label{oneroot:translating-endpoint-strips}

Now we transfer moving endpoint thresholds to the Gaussian field so that values close to the barrier can be estimated.
We append to the retained barrier product a smooth factor for
\[
 X_{N,J}\ge2J+B-r_N^{\rm strip},
\]
with any deterministic threshold \(r_N^{\rm strip}\). This is an existing
coordinate. Its derivative bounds have the same form \eqref{oneroot:eq:or16} and are
independent of the translating threshold. The Gaussian endpoint strip
has probability \(O(h/\sqrt t)\); the other strips have the preceding
\(O(h\sqrt t)\) bound. 
Therefore
\[
 \left|Q_{N,z}^{\rm det}(A_{N,J}^B,\ 0\le D_{N,J}\le r_N^{\rm strip})
       -Q_{t,z}^G(A_J^B,\ 0\le D_J\le r_N^{\rm strip})\right|
 =o(t^{-P}),
 \tag{OR18}\label{oneroot:eq:or18}
\]
where \(D_{N,J}=B+2J-X_{N,J}\), \(D_J=B+2J-G_J\). The same fixed smooth coarse insertion \(H\) can be included with
the derivative bound \eqref{oneroot:eq:or16}. Thus the estimate is uniform in the
translating endpoint threshold, which may depend on \(N\); it uses
the retained endpoint coordinate, not merely a fixed early prefix.

\subsection{Completing the integer history and terminal mark}\label{oneroot:completing-the-final-integer-barrier}

Now we restore the final integer barriers and terminal mark with errors that vanish after the critical normalization.
To restore the final integer barriers, we first control an upper
overshoot under the actual one-root law. For any fixed
\(0<p<\chi\), H\"older and the charge-\(\chi\) case of \eqref{oneroot:eq:or5} give
\[
\begin{aligned}
 \mathbb E_{Q^{\rm det}}e^{p[\Phi(z,\eta)-\Phi(z,0)]}
 &=\frac{\mathbb E e^{(\chi-p)\Phi(z,0)+p\Phi(z,\eta)}}{Z_\chi(0)}\\
 &\le[Z_\chi(\eta)/Z_\chi(0)]^{p/\chi}\le C.
\end{aligned}
 \tag{OR19}\label{oneroot:eq:or19}
\]
The power \(\chi-p\) is positive, so the numerator in \eqref{oneroot:eq:or19} is zero
at an exact determinant zero; the rooted identity is well defined.
Since \(n-J=O(\log t)\),
\[
 Q^{\rm det}\left(\max_{J\le k\le n}
 [\Phi(z,\eta_k)-\Phi(z,0)]>R\right)
 \le C(1+\log t)e^{-pR}.
\]
Combine this with \eqref{oneroot:eq:or15}, bounded finite mean corrections, and the
nonpositive subtraction \(-2(k-J)\). For \(R=C(1+\log t)^2\),
\[
 Q^{\rm det}\left(\max_{J\le k\le n}
 [X_{N,k}-X_{N,J}-2(k-J)]>R\right)=o(t^{-P}).
 \tag{OR20}\label{oneroot:eq:or20}
\]
On its complement, survival through \(J\) followed by a later crossing
forces \(0\le D_{N,J}\le R\). Applying \eqref{oneroot:eq:or18} and the Gaussian strip bound
\eqref{oneroot:eq:or1} gives
\[
 Q^{\rm det}(A_{N,J}^B\setminus A_{N,n}^B)
 \le C\frac{(1+\log t)^{9/2}}{t^{3/2}}+o(t^{-P}).
 \tag{OR21}\label{oneroot:eq:or21}
\]
The right side is \(o(t^{-1/2})\) for \(P>1/2\). This loss is not
claimed to be arbitrarily small in powers of \(t\); it has the natural
endpoint-strip contribution shown.

Insert retained barriers between the matrix and Gaussian full barriers.
The retained comparison from
Section~\ref{oneroot:retained-history-comparison-with-fixed-coarse-observations},
the matrix completion \eqref{oneroot:eq:or21}, and
the independent Gaussian completion \eqref{oneroot:eq:or2} then prove
\[
 \sqrt t\left|
 \mathbb E_{Q_{N,z}^{\rm det}}
    [H(S_N,U_N)\mathbf1_{A_{N,n}^B}]
 -\mathbb E_{Q_{t,z}^G}
    [H(S_G,U_G)\mathbf1_{A_n^B}]
 \right|\longrightarrow0.
\tag{OR22}\label{oneroot:eq:or22}
\]
The bounded coarse insertion costs only \(\|H\|_\infty\) when discarding
a terminal barrier loss. The original finite fourth-cumulant correction
is retained throughout. This proves the full-history comparison under the determinant-weighted
law, including the prescribed coarse observations.

The same comparison permits a bounded Lipschitz terminal mark. More
precisely, for such an \(h\), replace the two integrands in \eqref{oneroot:eq:or22}
respectively by
\[
 H(S_N,U_N)\mathbf1_{A_{N,n}^B}
       h((2t-Y_N^0)/\sqrt t),\qquad
 H(S_G,U_G)\mathbf1_{A_n^B}
       h((2t-G_t)/\sqrt t).
 \tag{OR22a}\label{oneroot:eq:or22a}
\]

To justify this at critical precision, first replace either terminal
mark by the mark of its retained deficit \((2J-X_{N,J})/\sqrt t\) or
\((2J-G_J)/\sqrt t\). On the event where the endpoint displacement bound holds, the argument changes by
at most \(C(\log t)^2/\sqrt t\); the Gaussian analogue follows from its
final \(O(\log t)\)-length increment. The additive rectangle comparison
\eqref{gaussian:gpm:additive-comparison}, followed by Brownian
reflection as in the Gaussian section, gives
\[
 Q_{t,z}^G(A_J^B)\le C\sqrt{\log(2+t)}/\sqrt t+O(t^{-D}),
\]
and the retained comparison gives the same bound for the matrix law.
Multiplying the mark error by this survival bound before multiplying by
\(\sqrt t\) gives \(O((\log t)^{5/2}/\sqrt t)=o(1)\); the endpoint
exceptions cost \(o(t^{1/2-P})\). Mollify \(h\) at width \(t^{-s}\) and apply the retained smooth
comparison to its retained-coordinate mark. Each derivative of this
factor has the mollification cost together with the coordinate scale
\(t^{-1/2}\), and its fixed-order summed cost fits \eqref{oneroot:eq:or16}. The uniform
mollification error is negligible even after multiplication by
\(\sqrt t\). Complete the barriers using boundedness as in \eqref{oneroot:eq:or21}.
This proves \eqref{oneroot:eq:or22} with \eqref{oneroot:eq:or22a}. A crude mark sup-norm error multiplied by
\(\sqrt t\) without the survival factor would not suffice.

\subsection{Endpoint cutoff and the remaining near-diagonal disk}
\label{oneroot:strong-cutoff-repair-at-dt13}

In this subsection we show that, at a fixed barrier, the terminal cutoff removes negligible expected mass and bounds the density enough to control the closest pairs.
The endpoint strip estimate now removes a terminal layer in expected
total variation. The resulting cutoff also gives a deterministic
density bound, which controls the remaining spatial core without a pair comparison. 
Next we set \(d=t^{1/3}\) and define
\[
 \mu_N^B(dz)=\sqrt t\,W_N^0(z)\mathbf1_{A_{N,n}^B(z)}\,d^2z,
 \quad
 \widetilde\mu_N^B(dz)=\mu_N^B(dz)
                   \mathbf1_{\{Y_N^0(z)\le2t-d\}}.
\]
On the complement of the endpoint exception in \eqref{oneroot:eq:or15}, full survival and
\(Y_N^0>2t-d\) imply
\[
 A_{N,J}^B,\qquad
 0\le D_{N,J}\le d+C(1+\log t)^2+C_B=:R_t.
\]
The term \(-2(t-J)\) actually improves this upper bound and has simply been discarded. 
Since \(R_t\asymp t^{1/3}\), \eqref{oneroot:eq:or18} and \eqref{oneroot:eq:or1} give
\[
 \sqrt t\,Q_{N,z}^{\rm det}(A_{N,n}^B,Y_N^0>2t-d)
 \le C t^{-1/3}\sqrt{\log(2+t)}+o(1).
 \tag{OR23}\label{oneroot:eq:cutoff-prob}
\]
The difference of these measures is positive. Tonelli's theorem
therefore identifies its expected total variation with the integral
of the uniform rooted loss in \eqref{oneroot:eq:cutoff-prob}, giving
\[
\mathbb E\|\mu_N^B-\widetilde\mu_N^B\|_{\rm TV}\to0.
\tag{OR24}\label{oneroot:eq:cutoff-tv}
\]
Any bounded mark and bounded coarse insertion multiply this loss by their supremum norms. 
The rate in \eqref{oneroot:eq:cutoff-prob}
comes directly from the endpoint strip with \(R_t\asymp t^{1/3}\).

Theorem~\ref{scalar:thm:sr-sharp}, with the bounded deterministic
mean correction, gives \(\mathbb Ee^{2Y_N^0}\asymp N\). 
Therefore the cut root density
\[
 U_N^{\rm cut}(z)=W_N^0(z)\mathbf1_{A_{N,n}^B(z)}
                    \mathbf1_{\{Y_N^0(z)\le2t-d\}}
\]
satisfies pointwise
\[
 0\le U_N^{\rm cut}(z)\le CN e^{-2d},\qquad \mathbb EU_N^{\rm cut}(z)\le1.
\]
On the core \(|z-w|^2\le C N^{-1}e^d\),
\[
 t\int_K\int_{|w-z|^2\le C N^{-1}e^d}
       \mathbb E[U_N^{\rm cut}(z)U_N^{\rm cut}(w)]\,d^2w\,d^2z
 \le Ct\,(N^{-1}e^d)\,(Ne^{-2d})=Cte^{-d}.
 \tag{OR25}\label{oneroot:eq:terminal-core}
\]

The pointwise density bound is used at one point and the expectation bound at the other; 
the disk has area \(O(N^{-1}e^d)\). 
Thus the core estimate uses no pair moment theorem. 
The expected total-variation conclusion \eqref{oneroot:eq:cutoff-tv} controls bounded tests, while
the later second-moment argument keeps the cut measures and its own pair estimates.

\subsection{Order of parameter choices}\label{oneroot:complete-parameter-order}

Now we choose the parameters in an order that preserves the precision required for the full history and terminal mark.
The order of parameter choices for the main applications is in
Section~\ref{main:parameter-schedule}; the local estimates below justify
its corresponding choices.

Fix the precision and all failure exponents before choosing the
logarithmic exponents. The letters \(a\) (history buffer),
\(c_{\rm int}\) (internal height in the local-factor construction), and \(q_m\) (microscopic
proxy) denote different parameters. The counting exponent for the local-factor construction, written
\(a_{\rm count}\) here, is also separate from the history buffer;
it is already fixed by the weak-law and noise estimates.

\begin{enumerate}
\setlength{\itemsep}{0pt}\setlength{\parskip}{0pt}
\item
  Fix the desired endpoint/retained precision \(P>1/2\), and choose
  \(P_*>P+3\). Fix the finite coarse list and smooth insertion. Fix
  \(p=1/4\) and the finite charge interval containing
  \(\chi-p,\chi,\chi+p\) and all higher positive moment charges used to
  bound cap tails.
\item
  Prescribe the cap-shell error exponent \(D_{\rm cap}\), then
  choose the cap size \(A\) using the fixed higher positive charges.
  Fix also the common local-only cap for \eqref{oneroot:eq:or7}--\eqref{oneroot:eq:or10}. The higher
  charges control the scalar tails; they do not enlarge the coefficient
  range allowed by the signed-source estimate. Choose the fixed derivative order from Corollary~\ref{signed:cor:ss-eight} \(m\), and
  require \(H\in C_b^{M_H}\), \(M_H\ge\max\{m,5\}\).
  Fix the failure and high-moment orders after both cap envelopes
  and every fixed moment cost; the failure order need not equal
  \(D_{\rm cap}\). Choose a barrier smoothing exponent
  \(s>P_*+2\), so \(h=t^{-s}\) controls each union of Gaussian strips at
  precision \(o(t^{-P_*})\).
\item
  The weak-law and noise parameters \(M,m_0,a_{\rm count}\) are
  fixed independently of the later logarithmic windows. Let
  \(C_{\rm source}\) dominate the resulting fixed signed-source,
  local-law and one-site height floors. 
  Choose 
  \[
  q_m>P_*+2,\qquad
  c_{\rm int}>\max\{m_0,P_*+3,C_{\rm source}\}.
  \]
\item
  Choose \(b\) larger than all of 
  \[
  \begin{gathered}
  100M+2P_*+4,\quad c_{\rm int}+P_*+2,\\
  2[15M+\max\{2m_0+q_m,c_{\rm int}\}+P_*+2],\\
  2[m(1+s)+P_*+4],\quad m(1+s)+P_*+19,\\
  400,\quad C_{\rm source}+1,\quad24.
  \end{gathered}
  \tag{OR26}\label{oneroot:eq:reconstruction-order}
  \] 
  The first two lines ensure the local-factor reconstruction bounds in
  Section~\ref{lf:sec:lfl-parameters}. The following two bounds control \eqref{oneroot:eq:or17} after including the full derivative
  factor \(t^{m(1+s)}\). The remaining fixed thresholds preserve the
  source/local-law domains. For the one-root proof some are deliberately
  redundant.
\item
  Choose the history buffer after \(b\), with 
  \[
  a>\max\{5,s+P_*+3\}.
  \tag{OR27}\label{oneroot:eq:history-buffer}
  \] 
  Then \eqref{oneroot:eq:or13} multiplied by \(\mathcal D_1(F)\le Ct^{1+s}\) is
  \(o(t^{-P_*})\). The Gaussian root-shift cost
  \(t^{-2a}\mathcal D_1(F)\) is smaller still. Enlarging \(a\) changes
  neither the smallest initial root cutoff \(B_v\tau\) nor the requirements in
  \eqref{signed:eq:ss-rate}; it only moves the retained history to larger
  heights.
\item
  Choose the auxiliary noise windows for the local-factor construction in the order of
  Section~\ref{lf:sec:lfl-parameters}, using its fixed counting
  parameters and the failure and moment orders already chosen above.
  These choices of noise windows impose no new lower bound on \(b\). Set \(J\)
  with the uniform bulk constant \(C_K\). For every
  ultimately fixed set of exponents, \(\tau\ll\sigma\),
  \(\eta_m\ge N^{-1-\delta}\), and \(J/t\to1\) hold for large \(N\).
\end{enumerate}

This order of choices permits both the arbitrary \(o(t^{-P})\) endpoint tail and
the \(o(t^{-1/2})\) full-history conclusion. It does not choose a new
exponent depending on \(N\), and it does not require the history buffer
to affect the previously chosen reconstruction time exponent.

\subsection{Scope of the one-root comparison}\label{oneroot:audit-of-remaining-premises}

Now we clarify which observations and marks the one-root theorem retains and which later conclusions require estimates for pairs.
The comparison uses positive determinant moments and the same local
cap at the three nearby charges. Negative coefficients occur only
inside the capped initial increment calculation; removing coordinate
caps afterward is a statement about bounded tail events. The
variable-height upper ratio \eqref{oneroot:eq:or5} is a direct scalar estimate.

The Gaussian endpoint strips, terminal completion and retained-survival
bound used for the mark are proved independently of the matrix
comparison. They require no Markov property of the exact Gaussian
scale process. With these inputs, the argument above proves the full integer-history
comparison under the determinant-weighted law and both cutoff conclusions
of Theorem~\ref{oneroot:thm:one-root}.

The one-root theorem retains arbitrary prescribed coarse bulk sites,
repeated observations and heights through one. Its bounded Lipschitz
marks are comparison tests; the stronger bounded-test scope of
\eqref{oneroot:eq:cutoff-tv} follows from total variation. Spatial
fixed-barrier moment limits, Gaussian marked factorization and final
measure identification use the corresponding pair estimates and measure
convergence results.

\endgroup

\begingroup
\newcommand{\E}{\mathbb E}
\newcommand{\PP}{\mathbb P}
\newcommand{\ind}{\mathbf 1}
\newcommand{\D}{\mathcal D}
\newcommand{\br}[1]{\langle #1\rangle}
\newcommand{\norm}[1]{\lVert #1\rVert}
\newcommand{\abs}[1]{\lvert #1\rvert}
\section{Separated marked two-root comparison}\label{chap:separated}

The purpose of this section is to compare two determinant-weighted histories at fixed positive separation while retaining their terminal marks and prescribed coarse observations.
Indeed, this section proves a finite-time matrix-to-Gaussian comparison
for two roots at fixed positive separation, bounded Lipschitz marks, and
a fixed finite list of coarse field observations. We use the local-factor, signed-source, and relative replacement
estimates. Below we prove the extensions needed for points at fixed
positive separation and for variable endpoint heights.
Joint microscopic calibration comes first: one-point comparisons do
not control a test weighted by both roots. The argument then transports
the retained pair history, restores both terminal barriers and marks,
and finally restores the exact pair prefactors. Its Gaussian inputs
come from finite rectangle comparison and Brownian reflection; they do
not use the Gaussian fixed-barrier limit.

\subsection{Conventions and statement}

We first state the pair comparison with its exact normalizing factors so that the matrix and Gaussian expectations remain correctly scaled.
Let \(K\Subset\{|z|<1\}\), \(k_0\ge1\), \(B\in\mathbb R\), and
\[
 t=\tfrac12\log N,\qquad n=\lfloor t\rfloor,\qquad\chi=2\sqrt2.
\]
The matrix \(X=N^{-1/2}(\xi_{ab})\) has complex iid entries satisfying
\[
 \E\xi=\E\xi^2=0,\qquad \E|\xi|^2=1,\qquad
 \sup_N\E|\xi|^p<\infty\quad\hbox{for every fixed }p.
\]
Put \(\kappa=\E|\xi|^4-2\). With the exact variance-one singular-value
density \(\rho_z\), write
\[
 \Phi_N(z,\eta)
 =\tfrac12\log\det\big((X-z)(X-z)^*+\eta^2\big)
       -\tfrac N2\int_{\mathbb R}\log(x^2+\eta^2)\rho_z(x)\,dx .
\]
At zero a vanishing determinant gives \(\Phi_N=-\infty\), and its
positive exponential is zero. Define the exact clocks and fields by
\[
 u_s=\tanh s,\quad
 \eta_s(z)=(1-u_s)\sqrt{u_s^{-1}-|z|^2},\quad p_z(u)=u(1-u|z|^2),
\]
\[
 X_{N,s}(z)=\sqrt2\{\Phi_N(z,\eta_s(z))+\kappa p_z(u_s)^2/4\},
 \qquad
 Y_N^0(z)=\sqrt2\{\Phi_N(z,0)+\kappa(1-|z|^2)^2/4\}.
\]
Let
\[
 W_N^0(z)=\frac{e^{2Y_N^0(z)}}{\E e^{2Y_N^0(z)}},\qquad
 R_N^0(z)=\frac{2t-Y_N^0(z)}{\sqrt t},\qquad
 A_{N,j}^B(z)=\bigcap_{k=k_0}^j\{X_{N,k}(z)\le2k+B\}.
\]

The centered Gaussian family \(G_s(z)\) has the exact covariance
\[
 C_\kappa(s,z;u,w)
 =-\tfrac12\log D(\tanh s,z;\tanh u,w)
       +\tfrac\kappa2 p_z(\tanh s)p_w(\tanh u),
\]
where
\[
 D(a,z;b,w)
 =1-ab-ab\{(1-a)|z|^2+(1-b)|w|^2\}+ab|z-w|^2.
\]
Set
\[
 W_t^G(z)=e^{2G_t(z)-2C_\kappa(t,z;t,z)},\qquad
 R_t^G(z)=\frac{2t-G_t(z)}{\sqrt t},\qquad
 A_{G,j}^B(z)=\bigcap_{k=k_0}^j\{G_k(z)\le2k+B\}.
\]
Coarse vectors \(C_N,C_G\) use a fixed finite list of coordinates
with heights in a fixed compact subinterval of \((0,1)\).
The matrix coordinates use \(X_{N,s}/\sqrt2=\Phi_N+\kappa p_z(u_s)^2/4\)
and the Gaussian coordinates use \(G_s/\sqrt2\), at the corresponding
exact clocks. Repeated observations and fixed coarse spatial sites
are allowed. Integer clocks \(k\ge1\), at a dense collection of
bulk spatial points, give an exhausting coarse family. This theorem
keeps its stated compact height interval inside \((0,1)\).

\begin{theorem}[Comparison at fixed positive separation]
\label{separated:thm:sm-comparison}
Fix \(r>0\), bounded Lipschitz marks \(h_1,h_2\), and a bounded smooth
test \(\Psi\) of the fixed coarse observations, with bounded derivatives through the
fixed orders required by Corollary~\ref{signed:cor:ss-eight} and
Theorem~\ref{replacement:thm:rgfc-main}.
Uniformly for \(z,w\in K\) with \(|z-w|\ge r\),
\begin{align}
t\Big|&
 \E\!\left[W_N^0(z)W_N^0(w)\Psi(C_N)
 h_1(R_N^0(z))h_2(R_N^0(w))
 \ind_{A_{N,n}^B(z)\cap A_{N,n}^B(w)}\right]\nonumber\\
&-\E\!\left[W_t^G(z)W_t^G(w)\Psi(C_G)
 h_1(R_t^G(z))h_2(R_t^G(w))
 \ind_{A_{G,n}^B(z)\cap A_{G,n}^B(w)}\right]\Big|
 \longrightarrow0 .                                      \label{separated:eq:sm-main}
\end{align}
The reference retains the original \(B\), exact covariance, exact
normalizers, and terminal mark convention.
\end{theorem}

Write
\[
 a_N(z,w)=\E W_N^0(z)W_N^0(w),\qquad
 a_t^G(z,w)=e^{4C_\kappa(t,z;t,w)}
\]
for the two pair prefactors, and \(Q_N^{(2)},Q_t^{G,(2)}\) for the
corresponding normalized pair laws. Their densities with respect to
the unweighted matrix and Gaussian laws are, respectively,
\(W_N^0(z)W_N^0(w)/a_N(z,w)\) and
\(W_t^G(z)W_t^G(w)/a_t^G(z,w)\). A normalized comparison must
therefore be multiplied by its own prefactor on each side. Neither
prefactor is replaced by one. All constants below may depend on the
fixed separation \(r\).

\subsection{Local factors and sources at fixed separation}

Now we shall extend the local-factor estimates to fixed positive separation while keeping the reference factors independent of the original initial matrix.
The local-factor construction uses Sections~\ref{chap:noise}--\ref{chap:lf};
the signed-source theorem is Theorem~\ref{signed:thm:ss-defect}, and
relative replacement is Theorem~\ref{replacement:thm:rgfc-main}.
The estimates below at fixed positive separation verify the extra geometry needed
to use these results for this pair.

The corrected history-flow estimate is proved in
{\hyperref[chap:calibration]{Section~\ref*{chap:calibration}}},
Lemma {\hyperref[calibration:lem:sc-history]{\ref*{calibration:lem:sc-history}}} and its logarithmic-floor
Corollary {\hyperref[calibration:cor:sc-history-log]{\ref*{calibration:cor:sc-history-log}}}, used with $L=b+2a>1$.
We state the quantitative bound used below, including its stochastic term.

\paragraph{Full characteristic and independent local factors.}
CL \cite[Proposition B.2(ii)]{CL} explicitly permits
\(|z-w|\ge r\), arbitrary deterministic insertions of norm at most one,
and unequal heights:
\begin{equation}\label{separated:eq:sm-B2}
 \left|\br{G_z(ix)A G_w(iy)B}
       -\br{\mathcal M_{z,x;w,y}[A]B}\right|
 \le\frac{C_{D,r}(\log N)^{10}}{N|xy|}.
\end{equation}
Its ranges are
\((\log N)^{10}/N\le|x|\le N^{-\varepsilon}\) and
\((\log N)^{10}/N\le|y|\le N^C\).
This estimate holds under the unweighted iid matrix law.

On the aging characteristic of the local-factor construction, both
heights are at most \(C\sigma\), where \(\sigma=N^{-1/200}\), and
the separation is at least \(r/2\). Choose
\(\varepsilon<1/200\). Equation \eqref{separated:eq:sm-B2} therefore holds on
the entire characteristic, with overwhelming probability
after the polynomial time/parameter mesh already proved in the
insertion argument. Intermediate matrices are iid in one common
moment class. At the initial height comparable to \(\sigma\), the
positive overlap obeys
\[
 Q_0=\bar Q_0+O_r(t^{10}/(N\sigma^2)),\qquad
 \bar Q_0\asymp_r1.
\]
Thus \(Q_0/\bar Q_0\le C_r\). This bounded initial ratio
is the input to the stopped overlap argument; the Ward bound
\(Q_0\le C/\sigma\) alone would not supply it.

Recall $Q=\br{\Im G_z\Im G_w}$ and $o$ from
\eqref{noise:eq:lfn-positive-identities}, and write
$\epsilon_i=\br{G_i}-m_i$ for the marginal resolvent error.
The exact positive-overlap equation and its relative estimates are
\[
 dQ=(1+2\operatorname{Re}o)Q\,ds+f_Q\,ds+dM_Q,\qquad
 |f_Q|/Q\le\sum_i|\epsilon_i|/\eta_i,\qquad
 d[M_Q]/Q^2\le C(N\eta^2)^{-1}ds .
\]
Here \(|o-\bar o|\le C_rt^{10}/(N\eta^2)\) by
\eqref{separated:eq:sm-B2}, and \(|\epsilon_i|\le C_Dt/(N\eta)\) by the
logarithmic one-point law. At the floor \(t^{24}/N\), the stopped
logarithmic-ratio drift integrates to \(o(1)\), and its bracket is
at most \(C/(N\eta_{\min})\). Exponential martingale concentration
therefore proves \(Q\le C_r\) on a common time event with arbitrary
polynomial failure strength. The projector/counting and positive
whitening calculation gives the following bounds for \(K_0=t^k\).
Here \(R_s\) is the off-diagonal covariance block in
\eqref{noise:eq:lfn-cross-matrix}, with index window \(K=K_0\):
\begin{equation}\label{separated:eq:sm-noise}
 \|R_s\|_{\mathrm{op}}\le C_rK_0^3/N,\qquad
 \left\|\frac{d\langle b-\beta\rangle}{ds}\right\|_{\mathrm{op}}
       \le C_rK_0^6/N^2.
\end{equation}
This is smaller than the noise error already allowed in the
different-noise energy proof. Every such occurrence in its drift,
bracket, and long-range tail is absorbed. The same-noise and
one-site logarithm reconstructions do not depend on the other root.

Consequently the two reference local families are mutually
independent and independent of the original initial-matrix
sigma field. They are not independent of the enlarged field
containing their own reference Ginibre initial spectra.
Only the independence from the original initial matrix is used when
the local factors are integrated out of the pair tilt.

\paragraph{Extension of the comparison at fixed separation.}
Equation \eqref{separated:eq:sm-B2} covers all sixteen constant-block
insertions. If either height is at most \(N^{-\varepsilon}\), put
it first by cyclicity. In the integral used below, the endpoint height
is at most one and the integration variable runs only to \(R=N^{10}\).
If both heights exceed \(N^{-\varepsilon}\), CL B.2(i) below one,
patched with B.1 above height \(1/2\), gives an integrated absolute
\(N^{-c}\) error. The discarded resolvent at the upper endpoint
\(R\) costs \(O(\log N/R)\) directly, as in Corollary~\ref{i2:cor:I6}; no assertion with
both heights unbounded is required. The deterministic insertion
\(A_i(u)=M_i(u)/(1-\br{M_i(u)^2})\) has norm \(O((1+u)^{-1})\);
the contribution above one is integrable. Arbitrary coarse spatial
points use the coarse bound in Theorem~\ref{i2:thm:I2}. Thus the integrated
coefficient bound remains
\[
 \max_{i,j}|\mathcal V_{ij}-V_{ij}|
 \le C_rt^{17}/(N\eta_{\rm mes})+N^{-c}.
\]
The other signed cumulant, cap, source-ray, and OU-Stein estimates
do not use the separation. Hence, for \(O(t)\) fine coordinates at
these two roots and a fixed coarse list, the same source theorem
holds, with derivative order \(m\) fixed before \(b\), coefficient
bound \(\|a\|_1\le6\), at most three nonzero coefficients, and error
\begin{equation}\label{separated:eq:sm-S1}
 r_N\le C_b\left\{\frac{t^2}{\sqrt{N\eta_{\rm mes}}}
                +\frac{t^{17}}{N\eta_{\rm mes}}+N^{-c}\right\}.
\end{equation}
The source-ray normalizer identity in
Corollary~\ref{signed:cor:ss-normalizer} holds with the same extension.
It supplies the exact initial pair-to-single ratio before any pair
law is normalized.

\subsection{Joint microscopic calibration}

The purpose of this subsection is to compare joint weights before dividing to show that microscopic regularization changes the law weighted at both points by a vanishing amount.

We first use the matched surrogate of
Proposition~\ref{replacement:prop:rgfc-age}. Its initial age is
exactly \(\sigma\), its additive duration is \(\tau=t^b/N\), and the
final variance-one normalization is matched through all mixed third
moments with fourth discrepancy \(O(\sigma)\). Spatial points,
heights, and deterministic potentials are scaled together by the
exact variance conversion of Lemma~\ref{replacement:lem:rgfc-variance}.
The physical endpoint corresponding to normalized \((v,\eta)\) is
\((\sqrt{1+\tau}\,v,\sqrt{1+\tau}\,\eta)\). In the flow formulas
below the variables denote the corresponding physical coordinates;
the stated uniform estimates apply to these endpoints. The exact
conversion is retained when returning to the target coordinates.

The bounded natural mean corrections cancel from the normalized
laws and from pair-to-single ratios, so suppress them in this scalar
calculation. The raw and regularized pair weights are
\[
 A=e^{\chi\Phi_\tau(z,0)+\chi\Phi_\tau(w,0)},\qquad
 B_m=e^{\chi\Phi_\tau(z,\eta_m)+\chi\Phi_\tau(w,\eta_m)},\qquad
 \eta_m=N^{-1}t^{-q},
\]
and \(m_v^0=\E e^{\chi\Phi_\tau(v,0)}\).
Here \(\eta_m\) is first a physical cutoff. Every estimate below is
uniform for cutoff multipliers in \([1/2,2]\). Thus a
target cutoff \(N^{-1}t^{-q}\) is handled at the physical cutoff
\(\sqrt{1+\tau}\,N^{-1}t^{-q}\) on return. The ratio
\(\E A/(m_z^0m_w^0)\) is the corresponding pair prefactor, since
the deterministic endpoint corrections cancel.
Write \(L_v=L_{1,v}+L_{2,v}\) for the sum of the reference increments
in \eqref{lf:eq:lfl-L1}--\eqref{lf:eq:lfl-L2}.
For any prescribed \(P\), the reconstruction at fixed positive separation gives
\[
 \Phi_\tau(v,\eta_m)=C_{0,v}+L_v+o(t^{-P}),\qquad
 C_{0,v}=\Phi_0(v,B_v\tau),\quad B_v=\sqrt{1-|v|^2},
\]
outside an event of probability \(O_D(N^{-D})\). With the local-only caps,
\[
 U_v=e^{\chi L_v}\ind_{\mathcal K_v},\qquad
 V_v=e^{\chi C_{0,v}},
\]
the three groups \(U_z,U_w,(V_z,V_w)\) are independent. Each
\(U_v\) contains only its own reference local cap. The joint coupling
good event is removed in weighted \(L^1\), rather than inserted into
these independent factors.

Fixed positive moments and deterministic local caps bound all
second moments used here by fixed powers of \(N\). Choosing the
exceptional exponent after these fixed costs gives
\begin{equation}\label{separated:eq:sm-proxy-rawineq}
 \E|B_m-U_zU_wV_zV_w|\le C\rho_N\E B_m+N^{-D},
 \qquad \rho_N=o(t^{-P}),
\end{equation}
and the equivalent reversed two-sided expectation inequalities.
At one root, the already proved scalar calibration gives
\[
 (\E U_v)(\E V_v)=m_v^0(1+o(t^{-P})),\qquad m_v^0\asymp N .
\]

At the initial heights, apply the source-ray theorem separately
to the pair source \((2,2)\) and its two single sources. Positive
cap removal is elementary at every required height. For the
monotone \(1\)-Lipschitz cap,
\[
 e^{2c_A(Y)}\le e^{2Y}+N^{-2A},\qquad
 |e^{2Y}-e^{2c_A(Y)}|
 \le e^{2Y}\ind_{\{Y>A\log N\}}+N^{-2A}.
\]
A higher fixed positive moment and H\"older bound each pair shell;
choose the cap large enough. This gives
\begin{equation}\label{separated:eq:sm-initialratio}
 \frac{\E[V_zV_w]}{\E V_z\,\E V_w}
 =e^{4C_{\kappa_0}(\ell_z,z;\ell_w,w)}(1+o(t^{-P}))
 =e^{4K_\kappa(z,w)}(1+o(t^{-P})),
\end{equation}
where \(\ell_v\) is the exact clock of \(B_v\tau\) and
\[
 K_\kappa(z,w)=-\log|z-w|
             +\tfrac\kappa2(1-|z|^2)(1-|w|^2).
\]
The last covariance error is \(O_r(\tau+\sigma)\), including
the exact variance/spatial conversion and actual endpoint cumulants.
Independence and \eqref{separated:eq:sm-initialratio} show that the proxy
expectation is comparable to \(N^2\). Only now does
\eqref{separated:eq:sm-proxy-rawineq} imply
\begin{equation}\label{separated:eq:sm-proxyL1}
 \E|B_m-U_zU_wV_zV_w|/N^2=o(t^{-P}).
\end{equation}
For example, the first inequality gives
\((1-C\rho_N)\E B_m\le\E[U_zU_wV_zV_w]+N^{-D}\);
its reversed inequality supplies the other direction. Thus the
regularized pair mass is established by absorption, before division
and without presuming a regularized pair bound.

Independently, CL \cite[Theorem 1.4, (1.14)]{CL} gives
\begin{equation}\label{separated:eq:sm-CLpair}
 \frac{\E A}{m_z^0m_w^0}
 =|z-w|^{-\chi^2/2}
   \exp\left\{\frac{\kappa\chi^2}{4}
                    (1-|z|^2)(1-|w|^2)\right\}
       (1+O_r(N^{-c}))
 =e^{4K_\kappa(z,w)}(1+O_r(N^{-c})).
\end{equation}
The normalized one-point constants cancel. Its fixed-power
separation parameter is distinct from our logarithmic flow exponent
\(b\); fixed positive separation satisfies the source condition.
The estimates are uniform in the common moment class, including
the triangular matched laws.

It follows that \(\E B_m/\E A=1+o(t^{-P})\).
Centered determinant monotonicity gives \(B_m\ge e^{-Ct^{-q}}A\).
For nonnegative \(A,B\) with \(B\ge e^{-\varepsilon}A\),
\[
 \E|A-B|\le|\E B-\E A|+2(1-e^{-\varepsilon})\E A .
\]
Choose \(q>P+2\). The one-sided determinant comparison and the
calibrated total masses now control the full absolute difference of
pair weights, giving the joint conclusions
\begin{equation}\label{separated:eq:sm-pairTV}
 \E|A-B_m|/\E A=o(t^{-P}),\qquad
 \|Q_{\rm raw}^{(2)}-Q_m^{(2)}\|_{\rm TV}=o(t^{-P}).
\end{equation}

For the target matrix, relative replacement transfers the capped
regularized pair and single normalizers, with centers fixed during
interpolation. Remove the positive caps by the moment argument
above. This argument works at \(\eta_m\) below the range of the
mesoscopic two-sided tail law; that tail law is not invoked there.
CL supplies \eqref{separated:eq:sm-CLpair} for the target itself.
The same monotonicity proof gives \eqref{separated:eq:sm-pairTV} for the
target. The source normalizers are compared before division; replacing
one-point tilted laws separately would not establish this pair estimate.

\subsection{A variable-height joint scalar upper bound}

Now we vary one regularization height to control large increases on the final scales under the weight at both points.

\begin{lemma}\label{separated:lem:sm-jointscalar}
Uniformly at fixed separation, the following bound holds for every
deterministic height \(\eta\in[\eta_m,\eta_J(z)]\):
\[
 \frac{\E e^{\chi\Phi(z,\eta)+\chi\Phi(w,0)}}
      {\E e^{\chi\Phi(z,0)+\chi\Phi(w,0)}}\le C_r,
       \qquad \eta_m\le\eta\le\eta_J(z),
\]
where \(J\) is selected in the next section. In particular, the bound holds
at the final integer-clock heights \(\eta_k(z)\), \(J\le k\le n\).
\end{lemma}
\begin{proof}
Theorem~\ref{scalar:thm:sr-regularization}, in the one-sided form
\eqref{oneroot:eq:or5}, gives the following comparison at every sufficiently
small height. Set \(s=\eta/\Im m^z(i\eta)=O(\eta)\). For any fixed real
\(\lambda\ge2\),
\begin{equation}\label{separated:eq:sm-scalarupper}
 \E e^{\lambda\Phi(z,\eta)}
 \le \E e^{\lambda\Phi(z/\sqrt{1+s},0;X^{(s)})}
 \le(1+C(\eta+N^{-\alpha}))\E e^{\lambda\Phi(z,0)}
 \le C N^{\lambda^2/8}.
\end{equation}
Here \(X^{(s)}\) is the normalized Gaussian convolution; its common
moment bounds and the last comparison follow from the scalar proof
and the sharp one-point theorem of Cipolloni--Landon, as explained in that scalar
section. No lower regularization comparison
is used, so \(N\eta\) may be a large logarithmic power.

Treat the retained-history parameters \(b,a\), and hence \(J\),
as already fixed for the following scalar calculation. Use a second
matched auxiliary flow
\(\tau'=t^{b'}/N\), \(\eta_*'=t^{c'}/N\), where
\(c'>b+2a+3\) and then \(b'\) is sufficiently large for the local-factor
inequalities with these auxiliary parameters. Since \(\eta_J\asymp t^{b+2a}/N\), all endpoints lie below
\(\eta_*'\ll\tau'\). This scalar calculation has no history test.
The auxiliary exponents \(c',b'\) do not replace \(b,a\) in the
retained integer \(J\).

Fix one such deterministic height \(\eta\). The following construction
may be chosen separately for each \(\eta\); its constants are uniform over
the stated interval. Replace the lower endpoint at \(z\) by \(\eta\), keeping \(\eta_m\)
at \(w\). This is an upper-bound reconstruction on the same
coupling space. The absolute derivative of
\[
 \tfrac12\log\frac{x^2+\eta_*'^2}{x^2+\eta^2}
\]
decreases as \(\eta\) increases in this interval. The nonnegative
centering error majorants also decrease, and the characteristic
martingale is independent of the lower endpoint. Thus every error in local-factor
reconstruction is dominated by the original \(\eta_m\) error.
At the first root choose a cap depending only on its reference local
family, common to the endpoint list \(\{\eta_m,\eta\}\); at the other root
use its microscopic endpoint. The two roots use separate such caps.
CL Lemma 2.12 applies to each increment with its cutoff exponent
chosen above \(c'\), uniformly over these deterministic endpoint choices.
For the later application to the \(O(\log t)\) final integer heights, use
a cap common to that finite list; the union bound retains arbitrary
polynomial failure strength. No simultaneous construction for a continuum
of heights is needed. Each cap remains reference-only, so independence from the entire
auxiliary initial matrix and between the two local families is preserved.

Let \(V'_v=e^{\chi\Phi_0(v,B_v\tau')}\).
The weighted reconstruction inequalities give
\[
 \E e^{\chi\Phi(z,\eta)+\chi\Phi(w,\eta_m)}
 \le(1+C\rho'_N)(\E U'_{z,\eta})(\E U'_{w,m})
                         \E[V'_zV'_w]+N^{-D},
\]
while the one-site inequalities give
\[
 (\E U'_{z,\eta})(\E V'_z)
 \le(1+C\rho'_N)m_\eta(z)+N^{-D},
\]
and the analogous upper bound at \(w\).
Here \(m_\eta(z)=\E e^{\chi\Phi(z,\eta)}\).
Only these upper inequalities are used; no lower bound on
\(m_\eta(z)\) is needed. The scalar upper estimate applies at each
variable endpoint of this auxiliary calculation.
The initial pair ratio from Corollary~\ref{signed:cor:ss-normalizer}, with the fixed-separated extension above is bounded by \(C_r\), so
\eqref{separated:eq:sm-scalarupper} makes the joint expectation at most \(C_rN^2\).
The raw \(w\) factor is bounded by \(e^{Ct^{-q}}\) times its
microscopic proxy. Divide by \eqref{separated:eq:sm-CLpair}.
Relative positive-source replacement and the same positive cap
removal transfer this auxiliary-surrogate bound to the target.
\end{proof}

For fixed \(0<p<\chi\), H\"older under the joint raw pair law gives
\[
 \E_{Q_N^{(2)}}e^{p[\Phi(z,\eta)-\Phi(z,0)]}
 \le\left(
 \frac{\E e^{\chi\Phi(z,\eta)+\chi\Phi(w,0)}}{\E A}
                 \right)^{p/\chi}\le C_r .
\]
Set \(R=A_0(1+\log t)^2\).
A union over the last \(O(\log t)\) heights gives, for every fixed \(P\),
\begin{equation}\label{separated:eq:sm-overshoot}
 Q_N^{(2)}\left\{
 \max_{J\le k\le n}[\Phi(z,\eta_k)-\Phi(z,0)]>R\right\}
 =o(t^{-P}),
\end{equation}
and similarly at \(w\). These are estimates under the joint raw
pair law, supplied to the retained-history argument together with
\eqref{separated:eq:sm-pairTV} and the exact pair-to-single ratios.
The auxiliary scalar flow is discarded after these estimates are
obtained; it does not move the retained integer \(J\).

\subsection{Retained comparison and two-sided terminal displacement}

The purpose of this subsection is to compare both retained histories and bound the terminal displacement under the joint weight, with the coarse observations preserved. We fix \(P>3\), smoothing width \(\varepsilon=t^{-A_s}\) with
\(A_s>P+3\), and \(M=\max(m,5)\), where \(m\) is the fixed derivative order from
Corollary~\ref{signed:cor:ss-eight}. A product of \(O(t)\) smoothed barrier factors,
the coarse test, and boundedly many endpoint factors satisfies, in the
norm \eqref{signed:eq:ss-norm},
\[
 \D_M(F)\le Ct^{M(1+A_s)} .
\]
Threshold translations do not affect this norm.
Choose \(q>P+3\), the internal local-factor cutoffs, and then \(b\)
large enough for the local-factor inequalities and
\[
 [t^{2-b/2}+t^{17-b}+N^{-c}]t^{M(1+A_s)}=o(t^{-P}).
\]
Afterwards choose \(a>P+A_s+4\) and
\begin{equation}\label{separated:eq:sm-J}
 J=\left\lfloor t-(b/2+a)\log t-C_K\right\rfloor,\qquad
 \eta_J(v)\asymp\tau t^{2a}.
\end{equation}
The corrected stopped history-flow bound is
\begin{align}
 \rho_{\rm hist}
 &\le C_D\left[
 t\sqrt{\frac{\tau}{\eta_J}}
       +t\frac{\tau}{\eta_J}
       +t^2\frac{\tau}{N\eta_J^2}+N^{-c}\right]\nonumber\\
 &\le C_D[t^{1-a}+t^{1-2a}+t^{2-b-4a}+N^{-c}].          \label{separated:eq:sm-history}
\end{align}
The square-root term is necessary. Since
\(\D_1(F)\le Ct^{1+A_s}\), this error after testing is \(o(t^{-P})\).
Weighted exceptional events are controlled by fixed positive moments
or cap envelopes before normalization.

Relative replacement transports the smoothed test under the
microscopic pair source. Equations \eqref{separated:eq:sm-pairTV},
\eqref{separated:eq:sm-proxyL1}, and \eqref{separated:eq:sm-history} turn it into an
initial-matrix-measurable test. Independence removes the two local
factors exactly, leaving the initial pair tilt
\(V_zV_w/\E[V_zV_w]\), to which the comparison of Corollary~\ref{signed:cor:ss-ou} with the fixed-separated extension and rate \eqref{separated:eq:sm-S1}
applies.

To apply the replacement theorem to capped coordinates, we first
cap the history and coarse test coordinates. Their heights are at least
\(\eta_J\), within the mesoscopic tail range. The union over \(O(t)\)
coordinates, H\"older with fixed higher positive pair moments, and
the polynomial lower bound on the pair normalizer make the weighted
active-cap event arbitrarily small polynomially in \(N\).
This removes caps from bounded marks and coarse tests.
The monotone cap preserves barrier indicators pointwise once their
entire transition interval lies in \([-A\log N,A\log N]\),
which holds for all \(k_0\le k\le J\) for the fixed sufficiently
large cap and all sufficiently large \(N\).
The microscopic positive source caps were removed separately by
moments, not by a mesoscopic tail estimate.

The initial external clocks satisfy
\(\ell_v=t-(b/2)\log t+O_K(1)\).
Changing them to \(t\) alters retained same-root means by
\(O(e^{-2(\ell_v-J)})=O(t^{-2a})\), and separated cross means by
\(O_r(\tau)\). Coarse cross terms are \(O(\tau)\).
Exact variance/spatial conversion and fourth-cumulant changes
give \(O_r(\tau+\sigma)\) covariance and mean errors.
Gaussian interpolation bounds these errors using \(\D_1+\D_2\), without
an inverse covariance. The \(N\log(1+\tau)/2\) potential terms
cancel exactly.

Finally, Gaussian marginal densities at all \(k\ge k_0\) are
uniformly bounded, including \(\kappa=-1\): the positive feature
expansion retains nonzero terms of degree at least two, giving a
uniform positive variance floor. Upper and lower barrier sandwiches
therefore differ by at most \(Ct\varepsilon=o(t^{-P})\).
The same holds for a translated endpoint strip. Sandwiching also transfers the bound to matrix barriers. For a
signed coarse insertion, bound the error by its supremum norm times
the strip probability; no ordering of signed tests is used.
We have proved
\begin{equation}\label{separated:eq:sm-retained}
 \left|\E_{Q_N^{(2)}}F(X_{\le J},C_N)
       -\E_{Q_t^{G,(2)}}F(G_{\le J},C_G)\right|=o(t^{-P})
\end{equation}
for the retained tests, including one-path barriers and endpoint
strips, at the original \(B\).

\begin{lemma}[Terminal displacement]\label{separated:lem:sm-endpoint}
With \(R=A_0(1+\log t)^2\), for \(v=z,w\),
\[
 Q_N^{(2)}\{|Y_N^0(v)-X_{N,J}(v)|>CR\}=o(t^{-P}).
\]
\end{lemma}
\begin{proof}
We prove the estimate for \(v=z\); interchanging the roots gives the case \(v=w\).
Choose \(0<p<\min(\chi-2,\sqrt2)\), and one common local-only cap
at powers \(\chi-p,\chi,\chi+p\).
One-site reconstruction, sharp positive moments, and the initial
one-source normalizer from Corollary~\ref{signed:cor:ss-normalizer} give
\[
 \E[e^{qL_v}\ind_{\mathcal K_v}]
       \asymp(N\tau)^{q^2/8},
       \qquad q\in\{\chi-p,\chi,\chi+p\}.
\]
Indeed the raw moment is \(N^{q^2/8}\) times a bounded positive
factor, while its initial moment is
\(\exp(q^2\ell_v/4)\) times a bounded factor.
Write \(Q_{U_v}\) for the law of the local reference family with density
\(U_v/\E U_v\). Hence \(\E_{Q_{U_v}}e^{\pm pL_v}\le t^{C_{p,b}}\).

Put \(C_{J,v}=\Phi_0(v,\eta_J(v))\), and write \(Q_{V_zV_w}\) for the
initial-matrix law with density \(V_zV_w/\E[V_zV_w]\).
Under this pair tilt use the three-coordinate source
\[
 (\chi\pm p)C_{0,v}\mp pC_{J,v}+\chi C_{0,w}.
\]
Its unit-coordinate absolute charge is at most
\(4+\sqrt2p<6\). The Gaussian exponent difference from the pair
source is \(O(\log t)\): the same-root increment has variance and
covariance of that order, and its covariance with the separated
second root is bounded. Thus the exponential moments of the
capped increment are at most \(t^{C_{p,b,a}}\).
The negative coefficient stays capped in this normalizer calculation.
For the active-cap event, combine the ordinary mesoscopic two-sided
coordinate tail \eqref{signed:eq:ss-tail} with higher positive pair
moments and H\"older. All three initial heights are in that tail
range. This removes the event with arbitrary polynomial error. We transfer only the
resulting tail estimate, not an uncapped negative moment:
\[
 Q_{V_zV_w}\{|C_{0,v}-C_{J,v}|>R\}
       \le t^{C_{p,b,a}}e^{-c_pR}+N^{-D}.
\]
Independence, reconstruction, and \eqref{separated:eq:sm-history} now control
the regularized endpoint.

For the raw endpoint put
\(d_v=\Phi(v,\eta_m)-\Phi(v,0)\ge-Ct^{-q}\).
On the raw-pair support \(A>0\), if \(d_z>s\), then
\(B_m/A\ge e^{\chi(s-Ct^{-q})}\), so
\[
 Q_{\rm raw}^{(2)}(d_z>s)
 \le e^{-\chi(s-Ct^{-q})}\E B_m/\E A\le Ce^{-\chi s}.
\]
Use \eqref{separated:eq:sm-pairTV} and the same estimate at \(w\).
Bounded natural mean corrections are absorbed into \(CR\).
The target endpoint tests are smoothed at positive \(\eta_m\)
before replacement; the raw logarithm is never differentiated.
\end{proof}

\subsection{Gaussian survival and killed strips}

At fixed separation, the exact covariance implies
\[
 C_\kappa(k,v;\ell,v)=\min(k,\ell)+O(1),\qquad
 C_\kappa(k,z;\ell,w)=O_r(1).
\]
Under the Gaussian double tilt \(G_k(v)-2k\) has bounded mean
for \(k\le t\).

\begin{lemma}[Pair survival and a one-path strip]
\label{separated:lem:sm-Gaussian}
Uniformly for the present parameters,
\begin{equation}\label{separated:eq:sm-survival}
 Q_t^{G,(2)}(A_{G,J}^B(z)\cap A_{G,J}^B(w))
       \le C_{r,B}(1+\log t)/t+O(t^{-D}).
\end{equation}
For \(1\le R\le J\),
\begin{align}
 &Q_t^{G,(2)}
 \{A_{G,J}^B(v),\ 0\le B+2J-G_J(v)\le R\}\nonumber\\
 &\qquad\le
 \frac{C_D R(1+\sqrt{\log J})(1+R+\sqrt{\log J})}{J^{3/2}}
       +C_DJ^{-D}.                                      \label{separated:eq:sm-strip}
\end{align}
\end{lemma}
\begin{proof}
Apply Lemma~\ref{gaussian:gpm:strips} with observation horizon
\(j=J\) and both external root clocks equal to the real time \(t\).
Its normalized double-root law is exactly \(Q_t^{G,(2)}\).
Fixed separation gives the bounded rooted means and covariance
bounds above. Equation~\eqref{gaussian:gpm:pair-survival}, with
\(J\asymp t\), proves \eqref{separated:eq:sm-survival}.
The marginal double-root clause of
\eqref{gaussian:gpm:one-strip} proves \eqref{separated:eq:sm-strip}.
The factors involving \(\log(2+J)\) and \(1+\log J\) are
comparable in this range, and the error power is arbitrary.
The independent proof there uses
\eqref{gaussian:gpm:additive-comparison}, regression and Brownian
reflection. It retains the original \(B\), handles singular arrays
by regularization, and uses no Markov property of the exact scale
process. The positive variance bound
\eqref{gaussian:gpm:positive-variance} includes \(\kappa=-1\).
\end{proof}

The intermediate pair bound has a logarithmic loss. This loss
does not weaken Theorem~\ref{separated:thm:sm-comparison}; the terminal mark
calculation below has enough margin. No sharper Gaussian pair bound is needed here.

\subsection{Late barriers, marks, and exact prefactors}

In this subsection we restore the late barriers and terminal marks, then recover the exact normalizing factors to complete the comparison at separated points.

On the good events of \eqref{separated:eq:sm-overshoot} and
Lemma~\ref{separated:lem:sm-endpoint},
\[
 X_{N,k}(v)-X_{N,J}(v)\le CR,\qquad J\le k\le n .
\]
A later crossing after retained survival therefore implies
\[
 A_{N,J}^B(v),\qquad
 0\le B+2J-X_{N,J}(v)\le CR .
\]
Drop the other barrier at this step, use the retained comparison
\eqref{separated:eq:sm-retained}, then \eqref{separated:eq:sm-strip}, and sum the two
possible crossing roots. We obtain
\begin{align}
 Q_N^{(2)}\big[
 &(A_{N,J}^B(z)\cap A_{N,J}^B(w))
       \setminus(A_{N,n}^B(z)\cap A_{N,n}^B(w))\big]\nonumber\\
 &\le C(1+\log t)^{9/2}t^{-3/2}+o(t^{-P}).               \label{separated:eq:sm-late}
\end{align}
The exponent \(9/2\) comes from \(R^2\sqrt{\log t}\).
This is already \(o(t^{-1})\), so dropping the second barrier is
legitimate for this strip estimate.

Replace the mark by
\(R_{N,J}^*(v)=(2J-X_{N,J}(v))/\sqrt t\).
Its good-event displacement is at most \(CR/\sqrt t\), since
\(t-J=O(\log t)\). Retained comparison and
\eqref{separated:eq:sm-survival} bound the probability of both retained
barriers by \(C(1+\log t)/t+o(t^{-P})\).
The survival-weighted Lipschitz error is therefore
\begin{equation}\label{separated:eq:sm-mark}
 C\frac{R}{\sqrt t}\frac{1+\log t}{t}+o(t^{-P})
       \le C(1+\log t)^3t^{-3/2}+o(t^{-P}).
\end{equation}
An unweighted \(R/\sqrt t\) bound would be insufficient.
A bounded coarse insertion costs only its supremum norm here.

For the Gaussian field, increments from \(J\) to the last integer
heights and to \(t\) have mean \(O(t-J)\) and variance
\(O(t-J+1)\) under the double tilt. Their maximum absolute value exceeds
\(R=(1+\log t)^2\) with probability \(o(t^{-P})\), for every fixed
\(P\), by a Gaussian tail and an \(O(\log t)\) union.
The same strip and survival proofs give Gaussian versions of
\eqref{separated:eq:sm-late} and \eqref{separated:eq:sm-mark}.

For bounded Lipschitz marks, convolve each mark with a smooth
probability mollifier at width \(\varepsilon=t^{-A_s}\).
Their supremum and Lipschitz constants stay bounded, the uniform
approximation error is \(O(\varepsilon)\), and
\(\|h_{i,\varepsilon}^{(q)}\|_\infty\le C_q\varepsilon^{-q}\)
for every \(1\le q\le M\).
Composing with the retained deficit divided by \(\sqrt t\) only
reduces derivative costs, so the displayed bound
\(\D_M(F)\le Ct^{M(1+A_s)}\) still applies. The extra approximation
error after multiplication by \(t\) is \(O(t\varepsilon)=o(1)\).
Thus this mollification fits the same parameter choice and extends
the preceding estimates to the stated Lipschitz marks.

Apply \eqref{separated:eq:sm-retained} to the retained marks and the coarse
test, then restore both terminal barriers and marks on both sides.
The normalized pair expectations differ by
\[
 O((1+\log t)^{9/2}t^{-3/2})+o(t^{-P}).
\]
Finally \eqref{separated:eq:sm-CLpair} and the exact covariance imply
\[
 a_N(z,w)=e^{4K_\kappa(z,w)}(1+O_r(N^{-c})),\qquad
 a_t^G(z,w)=e^{4K_\kappa(z,w)}(1+O_r(N^{-1})).
\]
Both prefactors are bounded above and below at fixed separation.
Restoring them adds \(O_r(N^{-c})\) times a bounded test.
Multiplication by \(t\) proves \eqref{separated:eq:sm-main}.

The following positive cutoff-loss bound permits arbitrary bounded
measurable marks and coarse insertions. The matrix-to-Gaussian
comparison itself retains the theorem's bounded Lipschitz marks,
smooth coarse test and compact coarse-height interval inside \((0,1)\).

\begin{corollary}[Separated endpoint-cut removal]
\label{separated:cor:sm-cut}
Set \(d=t^{1/3}\) and
\(\mathcal C_N(v)=\{Y_N^0(v)\le2t-d\}\).
For each fixed \(r>0\), uniformly on \(|z-w|\ge r\),
\begin{align}
 t\E\big[&W_N^0(z)W_N^0(w)
 \ind_{A_{N,n}^B(z)\cap A_{N,n}^B(w)}
 (1-\ind_{\mathcal C_N(z)}\ind_{\mathcal C_N(w)})\big]
 \le C_{r,B}t^{-1/3}(1+\log t)+o(1).\label{separated:eq:sm-cut}
\end{align}
Consequently the theorem also holds with both endpoint cuts on the
matrix side and the uncut Gaussian comparator. Bounded marks and
bounded coarse insertions multiply this loss by their supremum norms.
\end{corollary}
\begin{proof}
Outside the two endpoint exceptions of Lemma~\ref{separated:lem:sm-endpoint},
both full barriers and \(Y_N^0(z)>2t-d\) imply both retained
barriers and
\[
 0\le B+2J-X_{N,J}(z)\le
 R_t:=d+C(1+\log t)^2+C_B\asymp t^{1/3}.
\]
The omitted term \(-2(t-J)\) improves the upper bound.
Use \eqref{separated:eq:sm-retained} with this translating endpoint strip,
retaining the second barrier. The joint strip bound in
Lemma~\ref{gaussian:gpm:strips}, specifically \eqref{gaussian:gpm:pair-strip} gives
\begin{align*}
 &Q_N^{(2)}(A_{N,n}^B(z)\cap A_{N,n}^B(w),Y_N^0(z)>2t-d)\\
 &\quad\le
 \frac{C_{r,B}R_t(1+\log J)(R_t+1+\sqrt{\log J})}{J^2}
       +o(t^{-P})
 \le C_{r,B}t^{-4/3}(1+\log t)+o(t^{-P}).
\end{align*}
Sum the identical estimate at \(w\), and multiply by \(t\) and
the bounded pair prefactor \(a_N(z,w)\). Since \(P>3\), this
proves the claim. A one-path strip would not suffice at this cutoff;
the two retained barriers are essential here.
\end{proof}

\subsection{Order of parameter choices}

Now we choose the parameters so that the separate height estimates leave the retained history and final comparison precision intact.
All auxiliary parameters are fixed before the matrix limit.
Fix \(P>3\), the supremum and Lipschitz bounds for the marks, the coarse derivative bounds, the coefficient bound,
and a cap sufficient for the required positive shells.
Fix the source derivative order \(m\), \(M=\max(m,5)\), and
\(A_s>P+3\). Choose \(q>P+3\), the internal cutoffs in the local-factor construction, and then \(b\)
large enough for its established inequalities and the smoothed
source error. Choose \(a>P+A_s+4\) afterwards and \(J\) from
\eqref{separated:eq:sm-J}. For the separate scalar calculation only, choose
\(c'>b+2a+3\), then \(b'\) to satisfy the corresponding local-factor
inequalities. These latter choices do not alter the retained history.
The exceptional and high-moment orders are fixed after all source
and sitewise local cap envelopes and positive charges, before the
logarithmic exponents. Here the internal cutoffs in the local-factor construction mean its preliminary
counting and height parameters; its later particle and noise windows
are chosen after \(b\), in the order of
Section~\ref{lf:sec:lfl-parameters}. The auxiliary flow repeats that
order with \(b'\). These later windows use the fixed failure
exponents and impose no new condition on the retained \(b,a,J\). All logarithmic exponents remain fixed before
\(N\to\infty\).

This proves the stated finite-time comparison at fixed positive separation.
The Gaussian spatial construction, mixed-time Cauchy convergence,
removal of the matrix barrier, and integration near the spatial diagonal
are treated separately.

\endgroup

\begingroup
\newcommand{\E}{\mathbb E}
\newcommand{\PP}{\mathbb P}
\newcommand{\Tr}{\operatorname{Tr}}
\newcommand{\D}{\mathcal D}
\newcommand{\ind}{\mathbf1}
\newcommand{\br}[1]{\langle #1\rangle}
\newcommand{\dd}{\,\mathrm d}
\section{Bounds for pairs of nearby points}\label{chap:annuli}

The purpose of this section is to control nearby spatial pairs on successive annuli, leaving only the smallest disk to the separate terminal-cutoff estimate.

We compare the entire bounded compensated barrier observable for the
original complex iid ensemble and the exact Gaussian field. The comparison
controls both its weighted numerator and its normalizing constants. It proves the residual terminal-annulus estimate and the
outer-annulus domination required for near-diagonal uniform integrability.
The outer estimate uses an extension of the aged insertion theorem,
proved below. Its integrated error is small even where its extra pointwise
error is not. The terminal core and removal of the endpoint cutoff are treated
in the one-root section; the precise identification at fixed positive separation and the main joint
convergence theorem are proved separately.

\subsection{Conventions, inputs, and conclusions}

We first state the two annular bounds whose survival factors make the nearby-pair contribution summable outside the terminal core.
We let $X=N^{-1/2}(\xi_{ab})$ be complex iid, with
\[
 \E\xi=\E\xi^2=0,\qquad \E|\xi|^2=1,\qquad
 \sup_N\E|\xi|^r<\infty\quad(r<\infty).
\]
All statements are uniform for laws satisfying fixed common entry-moment bounds and
spatial points in a fixed compact $K\Subset\{|z|<1\}$. A fixed enlarged
compact is used after deterministic variance rescaling. Write
\[
 t=\tfrac12\log N,\qquad d=t^{1/3},\qquad
 \sigma=N^{-1/200},\qquad H=t^b,\qquad \tau=H/N,
 \qquad \eta_{\rm m}=N^{-1}t^{-q}.
\]
The parameters $b,q$ are fixed after a requested logarithmic precision;
they never depend on $N$. We use the density-centered potential of
Section~\ref{chap:scalar} and the exact clocks of
Section~\ref{chap:oneroot}. To specify the signed-source convention, put
\[
 E_\kappa(z,\eta)=-\frac\kappa4m^z(i\eta)^4,\qquad
 Y(z,\eta)=\sqrt2\,[\Phi(z,\eta)-E_\kappa(z,\eta)],
 \qquad \chi=2\sqrt2.
\]
Here $\kappa=\E|\xi|^4-2$. The raw normalized positive root is
\[
 W_N^0(z)=\frac{e^{\chi\Phi(z,0)}}{\E e^{\chi\Phi(z,0)}}
         =\frac{e^{2Y(z,0)}}{\E e^{2Y(z,0)}}.
\]
The bounded deterministic correction cancels exactly from this ratio.
At determinant zeros the numerator is zero. Set
\[
 u_k=\tanh k,\qquad
 \eta_k(z)=(1-u_k)\sqrt{u_k^{-1}-|z|^2},\qquad
 X_{N,k}(z)=Y(z,\eta_k(z)),
\]
for $k\ge k_0\ge1$, and define
\[
 A_{N,n}^B(z)=\{X_{N,k}(z)\le2k+B:\ k_0\le k\le n\}.
\]

\paragraph{Estimates used below.}
The full-block resolvent insertion estimates give the exact Riccati equation, its
propagator, the logarithmic law for nearby points and the integrated coefficient estimate.
The local-factor sections give pointwise reconstruction, weighted scalar
calibration, and two mutually independent local reference factors that are
also independent of the entire initial matrix. The sparse-weight comparison gives its explicit cumulant identity, normalizing constants along coefficient rays and smooth Stein
comparison. The relative replacement section gives both ratios of unnormalized
exponential moments and comparison under the normalized exponential weight. The scalar
section {\hyperref[chap:scalar]{Section~\ref*{chap:scalar}}}, specifically
{\hyperref[scalar:thm:sr-sharp]{\ref*{scalar:thm:sr-sharp}}}, {\hyperref[scalar:thm:sr-regularization]{\ref*{scalar:thm:sr-regularization}}},
{\hyperref[scalar:prop:sr-moments]{\ref*{scalar:prop:sr-moments}}}, and {\hyperref[scalar:lem:sr-variance]{\ref*{scalar:lem:sr-variance}}}, supplies sharp raw
normalization, positive moments, regularization calibration, and exact
variance conversion without using any two-root result.
Their extension to outer spatial pairs is proved in Section~\ref{annuli:sec:outer}.
Theorem~\ref{gaussian:gpm:external-theorem} supplies the independent
Gaussian envelope restated in Lemma~\ref{annuli:lem:gaussian}, including
external roots strictly deeper than the final retained barrier.

The two annular bounds use different horizons. Outer annuli stop at
$s_*=\lfloor3t/4\rfloor$, leaving a macroscopic gap after every
compensating depth $j\le t/2$. Terminal annuli retain the history up
to $n_T$, a distance of order $t^{1/3}$ before the real terminal root.
Their area cancels the factor $e^{2j}$ in the pair density.
The retained histories then supply polynomial survival factors in
place of the subcritical exponent gain discussed in
Section~\ref{chap:calibration}. These are the two estimates summed
below; the smaller remaining disk uses
\eqref{oneroot:eq:terminal-core} separately.

\begin{theorem}[Terminal and outer annular bounds]
\label{annuli:thm:annuli}
Fix $B$ and any $P>2$. The fixed parameters can be chosen so that the
following estimates hold uniformly in the indicated domains.
Let
\[
 n_T=\lfloor t-d/2-C_0\rfloor,\qquad s_*=\lfloor3t/4\rfloor.
\]
For integers $t/2<j\le n_T$ and
$c_1e^{-j}\le|z-w|\le c_2e^{-j}$,
\begin{equation}
 e^{-2j}\E[W_N^0(z)W_N^0(w)
       \ind_{A_{N,n_T}^B(z)\cap A_{N,n_T}^B(w)}]
 \le\frac{C_B}{(1+j)^{3/2}(1+n_T-j)}+C_Bt^{-P}.
 \label{annuli:eq:terminal}
\end{equation}
For integers $k_0\le j\le t/2$ in the same annular geometry,
\begin{equation}
 e^{-2j}\E[W_N^0(z)W_N^0(w)
       \ind_{A_{N,s_*}^B(z)\cap A_{N,s_*}^B(w)}]
 \le\frac{C_B}{(1+j)^{3/2}(1+s_*-j)}+C_Bt^{-P}.
 \label{annuli:eq:outer}
\end{equation}
Fixed-width changes in annular endpoints affect only the constants.
The negative source coordinate is always an observed depth $j\ge k_0$;
no source is placed at the invalid clock depth zero.
\end{theorem}

\begin{corollary}[Near-diagonal integrability outside the core]
\label{annuli:cor:ui}
Suppose
\[
 0\le U_N(z)\le W_N^0(z)\ind_{A_{N,\lfloor t\rfloor}^B(z)}.
\]
For bounded compactly supported $f$, the outer annuli satisfy, for
$L\ge k_0+2$,
\begin{align}
 \limsup_{N\to\infty}t\iint_{\{e^{-t/2}\lesssim|z-w|\le e^{-L}\}}
 |f(z)f(w)|\E[U_N(z)U_N(w)]\dd A(z)\dd A(w)
 \le C_{B,f}(1+L)^{-1/2}. \label{annuli:eq:outer-ui}
\end{align}
The residual annuli $t/2<j\le n_T$ contribute
$O(t^{-1/2}\log t)+O(t^{2-P})=o(1)$. Thus the pair integrals are uniformly
integrable near the diagonal outside the terminal core
$|z-w|^2\lesssim N^{-1}e^d$.
\end{corollary}

\subsection{Two-resolvent estimates for the outer annuli}
\label{annuli:sec:outer}

This section concerns the unweighted Gaussian-divisible matrix law at positive age, before exponential weighting.
The separation may shrink to order $N^{-1/4}$, so the estimate at fixed positive
separation in Section~\ref{chap:separated} alone is insufficient.
The proof first controls a mixed pointwise error and then integrates
it; only the integrated coefficient is required by the signed-source
comparison.
Take an OU observation age $S\in[c_{\rm age}\sigma,C_{\rm age}\sigma]$,
with fixed $0<c_{\rm age}<C_{\rm age}<\infty$, and fine heights
at least $t^{b_0}/N$, with $p=16$, $b_0=24$. The spectral domain
used below is
\[
 t^{b_0}/N\le\min(x,y)\le1,\qquad
 \max(x,y)\le R=N^{10},
\]
with fixed-factor buffers for normalization and meshes. The observation
endpoints in the integrated statement are at most one; only the
integration variable may run to $R$. Put
\[
 \delta=|z-w|^2,\qquad cN^{-1/2}\le\delta\le C_K,
 \qquad \epsilon=1/1000,\qquad
 A_\delta=N^{-1+2\epsilon}(1+\delta/\sigma)^2.
\]
For norm-one constant block insertions, let $F$ denote the sixteen
averaged two-resolvent coordinates and $\overline F$ their deterministic
reference. The exact already established equation is
\begin{equation}
 \dd E=(E+\overline FJE+EJ\overline F+EJE+\mathcal H)\dd s+\dd\mathcal M,
 \qquad E=F-\overline F.
 \label{annuli:eq:riccati}
\end{equation}
The OU characteristic is $z'=-z/2$, $\eta'=-\Im m-\eta/2$.
Writing $D_L=I-KJ$, $D_R=I-JK$, with $K'=K$, its deterministic
linear propagator is exactly
\begin{equation}
 U(v,u)[B]=e^{v-u}D_L(v)^{-1}D_L(u)
                       B D_R(u)D_R(v)^{-1}.
 \label{annuli:eq:propagator}
\end{equation}
Bulk stability gives
\begin{equation}
 \|U(v,u)\|\le C(\Gamma_u/\Gamma_v)^2,
 \qquad \Gamma_s=|z_1(s)-z_2(s)|^2+\eta_1(s)+\eta_2(s).
 \label{annuli:eq:gamma-prop}
\end{equation}
The algebraic propagator \eqref{i2:eq:propagator} has no near-pair
restriction; its stability bound therefore applies in this outer domain. At a terminal
point $x\le y$, backward time $r=v-u$ gives
$x_r\asymp x+r$, $y_r\asymp y+r$,
$\Gamma_r\asymp\delta+x+y+r$. In addition one may bound the propagator
by $C((y+r)/y)^2$.

\begin{proposition}[Pointwise and integrated bounds]
\label{annuli:prop:outer-i2}
With overwhelming probability, uniformly over the
stated ages, bulk spatial points and spectral heights,
\begin{align}
 \|E(x,y)\|&\le C_D\left[
       \frac{t^{16}}{Nxy}+\frac{A_\delta}{(\delta+x+y)^2}\right],
                    &&x,y\le\sigma,\label{annuli:eq:mixed-fine}\\
 \|E(x,y)\|&\le C_D\left[
       \frac{t^{16}}{Nxy}+\frac{N^{-.49}}{y}\right],
                    &&x\le\sigma\le y.\label{annuli:eq:mixed-coarse}
\end{align}
The direct polynomial law covers $x,y\ge\sigma$ and is smaller than
the displayed coarse bound. Consequently the integrated insertion
coefficient at a fixed height $\eta_{\rm mes}=t^b/N\le\eta\le1$
has error
\begin{equation}
 C_D\left[\frac{t^{17}}{N\eta_{\rm mes}}+N^{-.48}\right].
 \label{annuli:eq:outer-i6}
\end{equation}
\end{proposition}

\begin{proof}
First suppose $x,y\le\sigma$. The initial characteristic heights are
comparable to $\sigma$. The polynomial local law at arbitrary spatial points,
CL Proposition B.2(i), gives initial error
$C_DN^{-1+\epsilon}/\sigma^2$. Propagation by
\eqref{annuli:eq:gamma-prop} puts it below
$C_DN^{-\epsilon}A_\delta/\Gamma_v^2$. The factor
$(1+\delta/\sigma)^2$ is necessary when $\delta\gg\sigma$.
Set
\[
 P_r=\frac{t^{16}}{Nx_ry_r},\qquad S_r=\frac{A_\delta}{\Gamma_r^2}.
\]
Stop at a large fixed multiple of $P+S$ and the ordinary one-resolvent
stop. Here $P_r$ is the pointwise logarithmic envelope, distinct from
the requested final precision $P$ in Theorem~\ref{annuli:thm:annuli}. The logarithmic marginal laws apply throughout these paths.
The forcing and bracket estimates from the full-block insertion argument are
\begin{align}
 \|\mathcal H_r\|&\le\frac{Ct}{N}
       (x_r^{-2}y_r^{-1}+x_r^{-1}y_r^{-2}),\notag\\
 \dd[\mathcal M]_r&\le\frac{C}{N^2}
       (x_r^{-3}y_r^{-2}+x_r^{-2}y_r^{-3})\dd r.
 \label{annuli:eq:forcing-bracket}
\end{align}
Their transported drift is $o(P_v)$ and their transported bracket is at
most $Ct/(N^2x^2y^2)$. A fixed fraction of $P_v$ is exceeded with
probability at most $\exp(-ct^{31})$ for the stopped martingale.

The quadratic term splits into three estimates:
\begin{align}
 \int\|U\|P_r^2\dd r&\le CP_v\frac{t^{16}}{Ny},\label{annuli:eq:PP}\\
 \int\|U\|P_rS_r\dd r
 &\le CS_v\frac{t^{16}}N
       \int\frac{\dd r}{(x+r)(y+r)}
 \le CS_v\frac{t^{17}}{Ny},\label{annuli:eq:PS}\\
 \int\|U\|S_r^2\dd r
 &\le CS_v A_\delta\int_0^\infty(\Gamma_v+r)^{-2}\dd r
 \le CS_v\frac{A_\delta}{\Gamma_v}.\label{annuli:eq:SS}
\end{align}
The integral in \eqref{annuli:eq:PS} is at most $C\log N/y$, including equal
heights. Also
\begin{equation}
 \frac{A_\delta}{\delta}
 =N^{-1+2\epsilon}\left(\delta^{-1}+2\sigma^{-1}
                                  +\delta\sigma^{-2}\right)
 \le CN^{-.498}.
 \label{annuli:eq:A-small}
\end{equation}
All relative coefficients tend to zero, since $Ny\ge t^{24}$.
The propagated initial error is smaller than the stopping threshold by a fixed power
of the matrix size. The estimate therefore improves the threshold,
proving \eqref{annuli:eq:mixed-fine}.

Next suppose $x\le\sigma\le y$. The propagator is bounded because the
aging length is $O(\sigma)$. The initial polynomial bound is at most
$CN^{-1+1/200+\epsilon}/y$. Set $B_r=N^{-.49}/y$ and stop at $P+B$.
The three nonlinear relative coefficients are
\[
 \frac{Ct^{16}}{Ny},\qquad \frac{Ct^{17}}{Ny},\qquad
 \frac{C\sigma N^{-.49}}y\le CN^{-.49}.
\]
For marginal forcing, choose $\xi=10^{-4}$ and $a=N^{-10\xi}$.
If $\sigma\le y\le a$, use the logarithmic marginal law with upper
source-domain exponent strictly below $10\xi$. The second marginal
forcing divided by $P_v$ is at most
\[
 Ct^{1-p}\frac{x}{y}\log(1+C\sigma/x)
 \le Ct^{1-p}\frac\sigma y=o(1).
\]
For $a\le y\le1$, the polynomial marginal tolerance instead gives
\[
 Ct^{-p}N^\xi\frac\sigma y
 \le Ct^{-p}N^{-1/200+11\xi}=o(1).
\]
Using the polynomial tolerance down to $y=\sigma$ would not close this
relative estimate; the preceding logarithmic range is necessary.
For $y\ge1$, the stronger Ward decay $y^{-2}$ supplies an additional
factor $y^{-1}$ in the forcing. The bracket is unchanged. This proves
\eqref{annuli:eq:mixed-coarse}, with its written $y$ dependence. If both
heights exceed $\sigma$, B.2(i) directly gives the stronger bound
when both are at most one. B.1 patches the part where one height is
above $1/2$ and up to $R$, with the other bounded by one (or its
fixed-factor buffer). No bound with both heights unbounded is asserted.

For uniformity, use the stopped linear extension from
Lemma~\ref{i2:lem:uniformization}: it agrees with $E$ before the first exit and
continues linearly afterward. At a fixed deterministic terminal point,
the deterministic-kernel stochastic integral is a martingale in its
integration variable. No such integral is treated as a martingale in
terminal time. The same initial estimate works at every intermediate
terminal point $v$: since $\Gamma_0\asymp\delta+\sigma$,
\[
 \|U(v,0)E_0\|\le
 C N^{-\epsilon}A_\delta/\Gamma_v^2.
\]
In the mixed $P+B$ regime its bound is
$CN^{-.994}/y_0=o(N^{-.49}/y_v)$, because $y_0\asymp y_v$.
Thus the initial error remains strictly below the threshold even as $v\downarrow0$;
no aged input is reapplied at an intermediate age.
The existing polynomial spatial/height/time meshes apply:
resolvent and deterministic-reference derivatives have polynomial bounds,
ordinary OU matrix increments give the time modulus, and $P+S$, $P+B$
have positive thresholds with polynomial logarithmic derivatives.
The last deterministic mesh point before an alleged first exit preserves
the strict improvement. Enlarge the separation boundary by a fixed factor
and choose pointwise failure strength after the mesh exponent. The ages
remain comparable to $\sigma$; the larger domain is used only to estimate increments. The local-law
estimate is still applied within its stated range.

Finally the deterministic insertion is bounded below height one and
$O(1/u)$ above one. For $u,\eta\le\sigma$,
\[
 \int_0^\sigma\frac{A_\delta}{(\delta+u+\eta)^2}\dd u
 \le\frac{A_\delta}{\delta+\eta}\le CN^{-.498}.
\]
Where $\max(u,\eta)\ge\sigma$, the additional envelope is
$N^{-.49}/\max(u,\eta)$, whose integral up to one is
$O(N^{-.49}\log N)$ in either height ordering. The insertion decay makes
the upper tail integrable. The logarithmic part integrates to
$Ct^{17}/(N\eta)$. The discarded upper-endpoint resolvent
$G(iR)$ costs $O(\log N/R)$ directly, and the deterministic tail
is $O(R^{-1})$, as in Corollary~\ref{i2:cor:I6}. For fine--coarse pairs with a fixed positive lower bound on the coarse
height, use \eqref{i2:eq:coarse-I2} at arbitrary spatial points.
The shrinking interface above \(\sigma\) instead uses
\eqref{annuli:eq:mixed-coarse} and its marginal-law split;
trace cyclicity supplies the opposite height orientation.
This proves \eqref{annuli:eq:outer-i6} on its stated endpoint domain.
\end{proof}

\begin{remark}
At $\delta=N^{-1/2}$ the extra pointwise term in
\eqref{annuli:eq:mixed-fine} can be $N^{.002}$. The conclusion being used is
\eqref{annuli:eq:outer-i6}; the pointwise error need not satisfy a logarithmic bound uniformly in space.
\end{remark}

\begin{lemma}[Positive overlap and local-reference extension]
\label{annuli:lem:outer-lf}
The local-factor and signed-source comparisons at two roots extend to the outer range
$|z-w|\ge cN^{-1/4}$ with their same logarithmic precision and parameter
order.
\end{lemma}
\begin{proof}
For the positive overlap $Q=\br{\Im G_z\Im G_w}$, its deterministic
counterpart at comparable height $a$ is $\overline Q\asymp(\delta+a)^{-1}$.
The initial polynomial estimate under the unweighted law gives
\[
 \frac{|Q_0-\overline Q_0|}{\overline Q_0}
 \le CN^{-1+\epsilon}\frac{\delta+\sigma}{\sigma^2}
 \le CN^{-.989}.
\]
This replaces the near-pair Ward initialization, which is insufficient
relatively when $\delta\gg\sigma$. In the exact stopped logarithmic
ratio equation, the additional drift from the coefficient error has integral
\[
 \int|o-\overline o|\dd r
 \le\frac{Ct^{16}}{Na}+\frac{A_\delta}{\delta+a}
 \le\frac{Ct^{16}}{Na}+CN^{-.498}.
\]
The integrated relative marginal forcing is bounded by $C_Dt/(Na)$,
while the relative martingale bracket is bounded by $C/(Na)$,
as in the positive-overlap argument.
The adapted positive-overlap argument and its deterministic terminal
mesh yield $Q\le C/(\delta+a)$ with overwhelming probability.
The counting and whitening proof consequently gives, on the low-index
window $K=t^k$, noise discrepancy
\[
 \left\|\frac{\dd[b-\beta]}{\dd s}\right\|
 \le\frac{CK^6}{(1+N\delta)^2}\le\frac{CK^6}{N}.
\]
This is smaller than the discrepancy already handled by the moving-weight
energy estimate used in the local-factor construction. Its one-site reconstruction and weighted moments
are unchanged. The independent references remain independent of the
entire initial matrix. In the signed-source proof the only changed
second-cumulant input is \eqref{annuli:eq:outer-i6}; all third- and fourth-cumulant product bounds,
cap estimates and derivative bounds are already spatially uniform.
Thus the local-factor construction and signed-source comparison apply uniformly
in the outer range. The available inputs for compensated transport are
the integrated coefficient \eqref{annuli:eq:outer-i6}, the same weighted
local-factor calibration and the same source-ray normalizers. Their constants remain uniform as the separation approaches the stated
$N^{-1/4}$ boundary. The reference construction is made for each
prescribed deterministic pair; no independence across different pair
constructions is asserted.
\end{proof}

\subsection{The compensated bounded test and its caps}

Next we use the barrier to bound the compensating exponential while keeping the weighted pair expectation unchanged.

The pair weight grows as the roots approach. We insert a negative coefficient
at the observed depth $j$ and cancel its exponential factor inside the test. The barrier
at that same coordinate bounds the compensating exponential and all
its derivatives. This gives a bounded observable while retaining the
exact unnormalized pair density.

Let $n$ equal either $n_T$ or $s_*$ in the respective domain of
Theorem~\ref{annuli:thm:annuli}. Introduce two formal positive-root labels
$r_z,r_w$ and the observed negative label $c=(j,z)$, with $k_0\le j\le n$.
Fix one smooth odd nondecreasing cap $c_A$, with
\[
 c_A(x)=x\ (|x|\le A\log N),\qquad |c_A(x)|\le2A\log N,
 \qquad 0\le c_A'\le1.
\]
All fixed derivatives are uniformly bounded. At every layer set
$\mathbf X=c_A(\mathbf Y)$ coordinatewise and
\[
 p=2e_{r_z},\quad q_*=2e_{r_w},\quad a=p+q_*-2e_c,
 \qquad Z_M(a)=\E_Me^{a\cdot\mathbf X},
\]
\begin{equation}
 \mathcal C_M=\frac{e^{2j}Z_M(a)}{Z_M(p)Z_M(q_*)},\qquad
 \dd Q_{M,a}=Z_M(a)^{-1}e^{a\cdot\mathbf X}\dd\PP_M.
 \label{annuli:eq:source}
\end{equation}
All exponential weights are positive. Let $\theta$ be smooth decreasing,
identically one on $(-\infty,0]$ and zero on $[1,\infty)$. For $0<h\le1$
and a fixed smooth coarse factor $0\le F\le1$, whose bounded derivatives
are the ones required below and whose finite coordinate list has heights in
$[a_0,1]$ for a fixed $a_0>0$, define
\begin{align}
 \mathcal B_{B,h}(x)&=F(x_{\rm coarse})
   \prod_{v=z,w}\prod_{k=k_0}^{n}
       \theta\left(\frac{x_{k,v}-2k-B}{h}\right),\notag\\
 g(x)&=e^{2(x_c-2j)}\mathcal B_{B,h}(x).\label{annuli:eq:g}
\end{align}
The test does not use the two positive-root coordinates. Its factor at
$c$ gives, including on derivative supports,
\begin{equation}
 0\le g\le e^{2(B+h)},\qquad
 \D_m(g)\le C_{B,F,m}(1+t/h)^m.
 \label{annuli:eq:g-derivatives}
\end{equation}
Here $\D_m$ is the summed derivative-supremum norm in
\eqref{signed:eq:ss-norm}; no second dimension factor is inserted outside it.
Pointwise cancellation gives the exact identity
\begin{equation}
 \mathcal C_M\E_{Q_{M,a}}g
 =e^{-2j}\frac{\E_M[e^{2X_{r_z}+2X_{r_w}}\mathcal B_{B,h}(\mathbf X)]}
                    {Z_M(p)Z_M(q_*)}.
 \label{annuli:eq:compensation}
\end{equation}
The identical capped coordinate $X_c=c_A(Y_c)$ occurs in the source
and in $g$, so the cancellation remains exact even when the negative
cap is active. That cap is retained through every comparison. For $A>2$ and large $N$,
all barrier transition intervals lie inside the cap identity region,
and monotonicity gives
\[
 \theta((c_A(y)-2k-B)/h)=\theta((y-2k-B)/h)
\]
for every real $y$. Thus the barrier itself is exactly the uncapped one.

\begin{lemma}[Removal of the microscopic positive caps]
\label{annuli:lem:cap}
On the aged surrogate let $R_v=e^{2Y_v^{\rm m}}$,
$R_v^A=e^{2c_A(Y_v^{\rm m})}$ and $S=e^{-2c_A(Y_c)}$.
For every fixed $L$, selecting $A$ and then the failure exponent
large enough gives
\begin{equation}
 \E[S|R_z^AR_w^A-R_zR_w|]\le C_LN^{-L}.
 \label{annuli:eq:cap-remove}
\end{equation}
The same holds with a bounded test and for individual positive scalar
normalizers. No microscopic lower-tail theorem is required.
\end{lemma}
\begin{proof}
The initial mesoscopic ray $-4e_c$ has absolute charge four, within the
proved charge-six range. Corollary~\ref{signed:cor:ss-normalizer}
applied to this capped ray gives
\[
 \E S_0^2=\exp(8\Sigma_{cc}+O(r_N))\le Ce^{8j}\le CN^4.
\]
The history-flow event gives $S_\tau/S_0=e^{O(\rho_H)}$; on its complement
$S_\tau^2\le N^{8A}$. Choosing failure strength $D>8A+10$ yields
$\E S_\tau^2\le CN^4$. This is a capped mesoscopic moment, not an inverse
raw determinant moment.

Let $\varepsilon_A=N^{-2A}$. Since $R_v^A\le R_v+\varepsilon_A$,
\begin{align}
 \E[S|R_z^AR_w^A-R_zR_w|]
 &\le\varepsilon_A\E[S(R_z+R_w)]+\varepsilon_A^2\E S\notag\\
 &\quad+\E[SR_zR_w\ind_{\max(Y_z^{\rm m},Y_w^{\rm m})>A\log N}].
 \label{annuli:eq:cap-split}
\end{align}
The verified positive moment bound
\[
 \E e^{uY_v^{\rm m}}\le C_uN^{u^2/4+u/(2\sqrt2)},\qquad u\ge\sqrt2,
\]
implies, by Cauchy--Schwarz, that the first two terms are at most
$CN^{-2A+4+\sqrt2/2}+CN^{-4A+2}$. For fixed $u$ a higher fixed positive
moment gives
\[
 \E[e^{uY}\ind_{Y>A\log N}]\le C_AN^{-A^2+O_u(A)}.
\]
Apply H\"older with $S^2$, $R_z^4$ on the upper tail, and $R_w^4$.
The last term of \eqref{annuli:eq:cap-split} is $N^{-cA^2+O(A)}$.
Choose $A$ after $L$, then $D$ after $A$. The same proof at the initial
mesoscopic layer uses its already proved cap tails.
\end{proof}

\subsection{Transport of the numerator and the normalizers}

In this subsection we compare the weighted numerator and its normalizing constants together, removing exceptional events before using independence of the local factors.

First we transport the surrogate numerator and its scalar normalizers.
Then we apply relative replacement to return to the original ensemble.
Proposition~\ref{replacement:prop:rgfc-age} represents the matched
surrogate exactly as
$M=X_\tau/\sqrt{1+\tau}$, where the initial matrix $X_0$ has Gaussian age
exactly $\sigma$ (and comparable normalized ages along the flow).
Matching is imposed on $M$. The physical root
points are $v_z=\sqrt{1+\tau}z$, $v_w=\sqrt{1+\tau}w$.
Choose a fixed compact $K_*$ strictly inside the unit disk that
contains $K$ in its interior. All physical points belong to $K_*$ for
large $N$, and the estimates use its uniform bulk constants.
The microscopic physical height is $\sqrt{1+\tau}\eta_{\rm m}$;
the local-factor theorem's bounded deterministic multiplier version applies.
Let the initial
root height at $v$ be the height $B_v\tau$ used in the local-factor construction, where $B_v$ is the positive
bulk characteristic speed fixed in that section. Its exact real clock is
\begin{equation}
 \ell_v=t-\tfrac b2\log t+O_K(1).
 \label{annuli:eq:ell}
\end{equation}
In both domains under consideration $\ell_v>n$ for large $N$.

Write $Y_v^0$ for the initial corrected root field at the positive
height $B_v\tau$. This superscript denotes the initial layer, not
evaluation at zero height. Put
\[
 V_v=e^{2Y_v^0},\quad
 U'_v=e^{\chi(E_v^0-E_v^{\rm m})}U_v,\quad
 \mu_v=\E U'_v,\quad m_v=\E R_v,\quad v_v=\E V_v.
\]
The local factors $U_v$ are those in the density-centered local-factor theorem;
the written factor converts the density-centered local reconstruction
to the corrected $\sqrt2$-scaled field. Write $\rho_L$ for the reconstruction
error in \eqref{lf:eq:lfl-lf1} and $\epsilon_L$ for the relative scalar-calibration
error in \eqref{lf:eq:lfl-lf2}, with precision chosen below. These estimates give
\begin{equation}
 R_v=U'_vV_ve^{O(\rho_L)},\qquad
 \frac{\mu_vv_v}{m_v}=1+O(\epsilon_L).
 \label{annuli:eq:LF}
\end{equation}
The first relation holds on the reconstruction event intersected with
the root's own local-cap event; the scalar calibration is unconditional.
The two $U'_v$ are independent of each other and of the entire initial
matrix. Their own local cap indicators remain in their individual
factors; a joint coupling indicator is never included in their definition.

Initial signed ray normalizers, their exact covariance geometry and
positive cap removal yield
\begin{equation}
 v_z,v_w\asymp N/H,\qquad
 J_0:=\E[V_zV_wS_0]\asymp(N/H)^2e^{-2j}.
 \label{annuli:eq:initial-scales}
\end{equation}
Sharp scalar normalization gives $m_v\asymp N$, so
\begin{equation}
 \mu_v\asymp H,\qquad
 \mu_z\mu_wJ_0\asymp N^2e^{-2j}\ge cN.
 \label{annuli:eq:lower-scale}
\end{equation}
In the outer range the last lower bound improves to $cN^{3/2}$.
These bounds are obtained before any relative division. The initial
signed ray has two positive roots and the same capped negative
coordinate; its exact covariance exponent supplies $e^{-2j}$.
No terminal signed-normalizer estimate is assumed to establish this
initial lower bound.

\begin{lemma}[Comparison of history coordinates along the flow]
\label{annuli:lem:history}
On an event of overwhelming probability, the capped history/coarse
vectors $\mathbf H_\tau,\mathbf H_0$, including $c$, satisfy
\[
 \|\mathbf H_\tau-\mathbf H_0\|_\infty\le\rho_H,
 \qquad S_\tau/S_0=e^{O(\rho_H)}.
\]
For histories with minimum height $\eta_*$ one may take
\begin{equation}
 \rho_H\le C_D\left[
 t\sqrt{\tau/\eta_*}+t\tau/\eta_*
       +\frac{t^2\tau}{N\eta_*^2}+\sigma+N^{-c}\right].
 \label{annuli:eq:history-error}
\end{equation}
For $n=n_T$, $\eta_*\asymp N^{-1}e^d$; for $n=s_*$,
$\eta_*\asymp N^{-3/4}$. Thus $\rho_H$ beats every fixed inverse power
of $t$ in either regime.
\end{lemma}
\begin{proof}
Apply Lemma~\ref{calibration:lem:sc-history} to the common history
and coarse list, followed by the exact variance identity of
Lemma~\ref{replacement:lem:rgfc-variance}. The characteristic equation is
\[
 \dd\Phi=\dd M-\frac N2(m_N-m)^2\dd s.
\]
At each site the positive-index singular-value Brownian drivers have
$\dd[b_i,b_j]=\delta_{ij}\dd s$; cross-site independence is not needed.
The martingale bracket is at most $C\tau/\eta_*$, giving the first term
of \eqref{annuli:eq:history-error} after a concentration bound with overwhelming
probability. The local-law square gives the third term. Returning the cutoff
along the centered deterministic characteristic gives the second term.
All stochastic integrals are adapted and stopped before local-law failure.
A polynomial parameter mesh uses ordinary OU matrix modulus and resolvent
Lipschitz estimates. The exact variance identity
\[
 \Phi^{X_\tau/\sqrt{1+\tau}}(z,\eta)
 =\Phi_\tau^{X}(\sqrt{1+\tau}z,\sqrt{1+\tau}\eta)
\]
cancels the two deterministic $N\log(1+\tau)/2$ contributions exactly.
The remaining clock, bounded mean and fourth-cumulant changes fit the
written error. In particular the stochastic square-root term cannot be
omitted. At the terminal floor the three explicit terms are
$t^{1+b/2}e^{-d/2}$, $t^{b+1}e^{-d}$ and $t^{b+2}e^{-2d}$.
\end{proof}

\begin{proposition}
\label{annuli:prop:observable}
For a bounded $C^1$ function $g$ of the history/coarse vector,
\begin{align}
 &|\mathcal C_{\rm sur}\E_{Q_{{\rm sur},a}}g(\mathbf H_\tau)
       -\mathcal C_0\E_{Q_{0,a}}g(\mathbf H_0)|\notag\\
 &\quad\le C\left[(\rho_L+\rho_H+\epsilon_L)\|g\|_\infty
                        +\rho_H\D_1(g)\right]
       +O(N^{-L'}\|g\|_\infty).
 \label{annuli:eq:observable}
\end{align}
\end{proposition}
\begin{proof}
On the common event $\Omega$ of reconstruction, local caps and history
flow, \eqref{annuli:eq:LF} gives
\begin{align}
 &|R_zR_wS_\tau g(\mathbf H_\tau)
       -U'_zU'_wV_zV_wS_0g(\mathbf H_0)|\notag\\
 &\quad\le CU'_zU'_wV_zV_wS_0
       [(\rho_L+\rho_H)\|g\|_\infty+\rho_H\D_1(g)].
 \label{annuli:eq:numerator-pathwise}
\end{align}
This is a relative pointwise reconstruction, not a one-root $L^1$ error
multiplied by another unbounded root. If $U'_v\le CN^{A_*}$ and
$\PP(\Omega^c)\le C_DN^{-D}$, the positive moments give
\begin{align}
 \E[R_zR_wS_\tau\ind_{\Omega^c}]
   &\le CN^{4A+8+\sqrt2-D/2},\notag\\
 \E[U'_zU'_wV_zV_wS_0\ind_{\Omega^c}]
   &\le CN^{2A_*+4A+8+\sqrt2-D/2}.
 \label{annuli:eq:bad-budget}
\end{align}
Choose $D$ after all cap constants and the desired absolute accuracy.
We bound these absolute errors before dividing by
\eqref{annuli:eq:lower-scale}. We then remove the indicator of $\Omega$
before using independence, since this joint event need not factor. Then
\[
 \E[U'_zU'_wV_zV_wS_0g(\mathbf H_0)]
 =\mu_z\mu_w\E[V_zV_wS_0g(\mathbf H_0)].
\]
Taking expectations in \eqref{annuli:eq:numerator-pathwise}, using
\eqref{annuli:eq:bad-budget}, and repeating with $g=1$ gives both numerator
transport and comparison of normalizing constants without the barrier restriction. Restore positive caps by
Lemma~\ref{annuli:lem:cap}, use the two individual ratios in \eqref{annuli:eq:LF}, and
divide using \eqref{annuli:eq:lower-scale}. The initial compensated normalizer
is bounded above and below by Corollary~\ref{signed:cor:ss-normalizer} and the covariance geometry in \eqref{annuli:eq:initial-scales}. This proves
\eqref{annuli:eq:observable}. No division by a rare barrier probability occurs.
\end{proof}

\subsection{Replacement and restoration of the original parameters}

Now we return the compensated expectation to the original entry law and Gaussian parameters while preserving the retained history observations.
Theorem~\ref{replacement:thm:rgfc-main} applies to the same cap
and separately to the rays $a,p,q_*$. Their supports have at most
three coordinates and total absolute charges at most six. With $\epsilon_G=N^{-c}$ it gives
\[
 Z_{\rm tar}(a)/Z_{\rm sur}(a)=1+O(\epsilon_G),\qquad
 |\E_{Q_{{\rm tar},a}}g-\E_{Q_{{\rm sur},a}}g|
       \le C\epsilon_G\D_5(g),
\]
and the analogous two positive scalar ratios. Since the surrogate
compensated normalizer is already bounded,
\begin{equation}
 |\mathcal C_{\rm tar}\E_{Q_{{\rm tar},a}}g
       -\mathcal C_{\rm sur}\E_{Q_{{\rm sur},a}}g|
 \le C\epsilon_G(\D_5(g)+\|g\|_\infty).
 \label{annuli:eq:gfc}
\end{equation}
The factor $e^{2j}$ cancels inside the relative normalizer ratio. The
comparison of exponential moments is proved by the unnormalized
Green-function comparison inequality,
not by setting $g=1$ in an already normalized probability comparison.
Lemma~\ref{replacement:lem:rgfc-centers} returns the frozen
deterministic centers to each law's own fourth cumulant, adding only $O(\sigma)(\|g\|_\infty+\D_1(g))$.

If exact variance rescaling moves a fixed coarse height slightly above
one, return it to its original height at the same physical point before
invoking Corollary~\ref{signed:cor:ss-ou}. The centered height derivative at that positive cutoff is
$O(N^\xi)$ on the common event, so this $O(\tau)$ change is polynomially
small. Thus every capped coordinate remains in the height domain of Theorem~\ref{signed:thm:ss-defect}.

Corollaries~\ref{signed:cor:ss-ou} and
\ref{signed:cor:ss-normalizer} compare the initial normalized law
and its separate normalizing expectations, giving
\begin{equation}
 |\mathcal C_0\E_{Q_{0,a}}g-\mathcal C_{G,0}\E g(G_{a,0})|
 \le Cr_N\D_m(g),\qquad
 r_N\le C_b(t^{2-b/2}+t^{17-b}+N^{-c}).
 \label{annuli:eq:S1}
\end{equation}
Here $G_{a,0}$ has mean $\Sigma_0a$ and covariance $\Sigma_0$, where
$\Sigma_0$ is the covariance matrix of the initial Gaussian coordinate list;
$\mathcal C_{G,0}$ is the same compensated ratio as in
\eqref{annuli:eq:source}, with each matrix normalizer replaced by the
uncapped Gaussian exponential moment
$Z_{G,0}(u)=\exp(\tfrac12u^{\mathsf T}\Sigma_0u)$.
After restoring the parameters, we denote the corresponding objects by
$G_a$ and $\mathcal C_G$.
The fixed derivative order $m$ is chosen before $b$. Section~\ref{annuli:sec:outer}
justifies the same assertion for outer pairs.
Keep the fine root and its history at their common physical spatial site
until this Gaussian stage. Return the physical points, exact clocks and
fourth cumulant there. A sufficient bound for the covariance error is
\begin{equation}
 \delta_C\le C[\tau/|z-w|^2+\sigma+N^{-c}],
 \label{annuli:eq:deltaC}
\end{equation}
with the stronger $O(\tau+\sigma)$ bound available for the actual common
radial rescaling at exact clocks. It is super-polynomially small in $t$
in both annular domains. The coefficient vector has sum of absolute values six, so
the tilted mean moves by at most $6\delta_C$. The compensated Gaussian
normalizer changes relatively by $O(\delta_C)$. Convex Gaussian covariance
interpolation, valid also for singular covariance, therefore costs
\[
 C\delta_C(\|g\|_\infty+\D_1(g)+\D_2(g)).
\]
Keep the original integer $j$ through this return; recomputing a floor
may change it by one.

Combining these bounds with Proposition~\ref{annuli:prop:observable} gives
\begin{align}
 |\mathcal C_{\rm tar}\E_{Q_{{\rm tar},a}}g-\mathcal C_G\E g(G_a)|
 \le C\big[& (\rho_L+\epsilon_L)\|g\|_\infty
       +\rho_H(\|g\|_\infty+\D_1(g))\notag\\
 &+\epsilon_G\D_5(g)+r_N\D_m(g)\notag\\
 &+\delta_C(\|g\|_\infty+\D_1(g)+\D_2(g))\big]+N^{-L'}.
 \label{annuli:eq:master}
\end{align}
This transports the whole bounded compensated observable, not just its
normalizer.

To prove the unweighted annular bound, now take $F\equiv1$ in
\eqref{annuli:eq:g}. The preceding transport remains valid for the
stated bounded coarse insertions. For the original raw root,
determinant monotonicity gives
$\Phi(z,0)\le\Phi(z,\eta_{\rm m})+CN\eta_{\rm m}$.
The sharp scalar calibration gives
\[
 Z_{\rm reg}(z)/Z_{\rm raw}(z)
       =1+O(t^{-q}+N^{-\alpha}).
\]
Consequently, pointwise,
\begin{equation}
 W_N^0(z)\le[1+C(t^{-q}+N^{-\alpha})]
           \frac{e^{2Y(z,\eta_{\rm m})}}{\E e^{2Y(z,\eta_{\rm m})}}.
 \label{annuli:eq:raw-upper}
\end{equation}
The deterministic correction cancels in the individually normalized
weight at each point.
Apply this pointwise upper inequality at both
sites and dominate the hard barrier by the smooth product. Positive
upper cap tails are controlled by moments under the unweighted law, while their lower cap
increment only increases the capped upper comparison. Equivalently one
may use the established positive cap estimates and exact compensation.
Calibrate the two individual scalar normalizers. It follows that
\begin{align}
 e^{-2j}\E[W_N^0(z)W_N^0(w)\ind_{A_{N,n}^B(z)\cap A_{N,n}^B(w)}]
 \le[1+O(t^{-q}+N^{-\alpha})]
       \mathcal C_{\rm tar}\E_{Q_{{\rm tar},a}}g+O(N^{-L'}).
 \label{annuli:eq:raw-to-g}
\end{align}
No one-site total-variation estimate is multiplied by a second root.
This is an upper comparison for the raw pair observable. It uses
two-sided calibration only for the individual scalar normalizers.
Width $h=1$ suffices: the smooth product is supported in the hard barrier
at height $B+1$. If $h=t^{-s}$ is desired, all derivative costs remain
explicit in \eqref{annuli:eq:g-derivatives}; an upper bound requires only raising the barrier,
without an identity after removing the smoothing.

\subsection{The exact Gaussian envelope with deeper external roots}

In this subsection we bound pairs of Gaussian weights at scales finer than the last observed barrier scale, keeping the additional correlations explicit.
The independent proof is Theorem~\ref{gaussian:gpm:external-theorem}
in Section~\ref{gaussian:gpm:external}. We restate its envelope with
the present notation, including roots deeper than the final barrier.
The covariance is
\begin{align}
 C(s,z;r,w)&=-\tfrac12\log D(u_s,z;u_r,w)
              +\tfrac\kappa2p_z(u_s)p_w(u_r),\notag\\
 p_z(u)&=u(1-u|z|^2),\notag\\
 D(u,z;v,w)&=1-uv-uv[(1-u)|z|^2+(1-v)|w|^2]
                             +uv|z-w|^2.
 \label{annuli:eq:Gaussian-cov}
\end{align}
The family is positive semidefinite for $\kappa\ge-1$, including the
endpoint $\kappa=-1$. Each root is normalized using its own exact variance:
$W_\ell^G(v)=\exp(2G_\ell(v)-2C(\ell,v;\ell,v))$.
We write $A_n^{B,G}(v)=\bigcap_{k=k_0}^n\{G_k(v)\le2k+B\}$ for the
Gaussian barrier event.

\begin{lemma}[General two-root Gaussian envelope]
\label{annuli:lem:gaussian}
Let the external clocks satisfy $\ell_z,\ell_w\ge n$. Set
\[
 q_\infty=\log_+(1/|z-w|),\qquad
 q=\min(\ell_z,\ell_w,q_\infty),\qquad m=\min(n,\lfloor q\rfloor).
\]
Then
\begin{equation}
 \E[W_{\ell_z}^G(z)W_{\ell_w}^G(w)
                  \ind_{A_n^{B,G}(z)\cap A_n^{B,G}(w)}]
 \le\frac{C_B e^{4q-2m}}{(1+m)^{3/2}(1+n-m)}.
 \label{annuli:eq:Gaussian-general}
\end{equation}
If $j\le n$, $|q_\infty-j|\le C_1$ and $q_\infty\le n+C_1$, then
\begin{equation}
 e^{-2j}\E[W_{\ell_z}^G(z)W_{\ell_w}^G(w)
                  \ind_{A_n^{B,G}(z)\cap A_n^{B,G}(w)}]
 \le\frac{C_{B,C_1}}{(1+j)^{3/2}(1+n-j)}.
 \label{annuli:eq:Gaussian-annulus}
\end{equation}
The constants do not depend on the gaps $\ell_v-n$.
\end{lemma}

\paragraph{Normalization of the Gaussian bound.}
The exact kernel satisfies the uniform logarithmic profile
\eqref{gaussian:gp2:log-profile}:
\[
 |C(s,z;r,w)-\min(s,r,q_\infty)|\le C.
\]
The proof of Theorem~\ref{gaussian:gpm:external-theorem} retains
the exact double-root factor in \eqref{gaussian:gpm:double-tilt},
then uses a multiplicative two-leaf comparison and the killed-tree
estimate \eqref{gaussian:gpm:two-tree}. It includes the delayed
barrier start and zero-length continuations. Thus the root depths
in \eqref{annuli:eq:ell} are covered without turning a bounded
covariance discrepancy into an additive rare-event error. The
correlation depth $q$ and retained depth $m$ remain distinct in
\eqref{annuli:eq:Gaussian-general}; only the stated annular hypotheses
permit \eqref{annuli:eq:Gaussian-annulus}.

For the Gaussian vector containing the external roots and histories,
write
\[Z_G(a)=\exp(\tfrac12a^T\Sigma a).\]
Its exact compensation gives
\begin{equation}
 \mathcal C_G\E g(G_a)
 =e^{-2j}\E[W_{\ell_z}^G(z)W_{\ell_w}^G(w)\mathcal B_{B,h}(G)]
 \le\frac{C_B}{(1+j)^{3/2}(1+n-j)}.
 \label{annuli:eq:Gaussian-compensated}
\end{equation}
The inequality uses the barrier $B+1$ and discards a coarse factor bounded
by one. The logarithm of $\mathcal C_G$ is exactly
\begin{equation}
 2j+4\Sigma_{r_zr_w}-4\Sigma_{r_zc}-4\Sigma_{r_wc}+2\Sigma_{cc}.
 \label{annuli:eq:exact-normalizer}
\end{equation}
Every displayed covariance is $j+O(1)$ in the annular domain, so this
normalizer is bounded above and below. Root self-variances cancel; they
must not be replaced by a common approximate variance before this step.

\subsubsection*{A separate check for the outer annuli}

Next we check the direction and normalization of the outer-annulus bound by moving the Gaussian weighting scales back. 
The preceding theorem is the proof input in both annular domains.
The following separate calculation checks its direction and scalar
normalization in the outer range by returning both external clocks
to $s_*$ through a one-sided mean comparison. At one spatial site the exact kernel is
\[
 K^z(k,l)=-\tfrac12\log(1-u_ku_l)
 -\tfrac12\log(1-\rho A(u_k,u_l))
 +\tfrac\kappa2p_\rho(u_k)p_\rho(u_l),
 \quad \rho=|z|^2,
\]
where $p_\rho(u)=u(1-\rho u)$ and $A(u,v)=uv(2-u-v)/(1-uv)$. Direct differentiation gives
\[
 \partial_v A=\frac{u(1-v)[2-u(1+v)]}{(1-uv)^2}\ge0.
\]
Thus $K^z(k,\ell)-K^z(k,s_*)\ge-Ce^{-2s_*}$ for $k\le s_*\le\ell$.
This is one-sided: the absolute difference near $k=s_*$ need not be
small. Cross-site changes are at most
$Ce^{-2s_*}/\delta\le CN^{-1/4}$. At the compensating depth
$j\le t/2$, the same-site absolute change is also at most
$Ce^{-2(s_*-j)}\le CN^{-1/4}$. Hence every history mean under the deeper
roots is at least its $s_*$ mean minus $CN^{-1/4}$, while the compensating
coordinate and \eqref{annuli:eq:exact-normalizer} change absolutely by that
amount. Couple the two tilted arrays with the identical centered array.
Their decreasing barriers compare after an $O(N^{-1/4})$ lift, and the
compensating exponential changes only relatively. The same-time case
of \eqref{annuli:eq:Gaussian-general} then gives \eqref{annuli:eq:Gaussian-compensated}
for outer annuli. This separately checks the direction and scalar
normalization of the general external-depth argument.

\subsection{Parameter order and proof of the annular assertions}

We next choose the parameters and sum the annular estimates to control nearby pairs outside the separately bounded terminal core.
We first fix $P>2$, the compact bulk data and barrier height. Take width one
in \eqref{annuli:eq:g}; fix the source and local caps, all higher positive moment
charges used to control their tails, and the finite signed-Stein derivative
order $m$. Fix a precision $P_*>P+m+10$ for local-factor reconstruction and weighted
calibration, and take $q>P_*+3$. Fix the failure exponents
after the cap bounds, positive-moment orders and polynomial mesh sizes; these
orders are included among the weighted-estimate requirements in
Section~\ref{lf:sec:lfl-parameters}. Choose the reconstruction
parameters and then $b$ in that schedule. In particular fix the weak-law and same-noise parameters for the local-factor construction
$(M,m_0,a_{\rm count})$ and its internal-height exponent $c_{\rm LF}$
before $b$. Choose $b$ so that those inequalities and the following bounds hold:
\[
 b>2(m+P+5),\qquad b>m+P+20,\qquad b\ge24.
\]
Then \eqref{annuli:eq:S1} times the full $O(t^m)$ derivative cost is $o(t^{-P})$.
If smoothing width $t^{-s}$ is used instead, replace $m$ in these two
bounds by $m(1+s)$ and request the corresponding higher precision in local-factor reconstruction.
Choose the auxiliary local noise windows afterward in the same order
of parameter choices for the local-factor construction. The fixed replacement-height exponent can be any sufficiently small
$0<\delta_0<1/1200$; eventually
$\eta_{\rm m}\ge N^{-1-\delta_0}$. No exponent depends on $N$.

All history terms and errors from restoring the Gaussian parameters in
\eqref{annuli:eq:master} are smaller than every fixed inverse power of $t$: at the terminal separation,
$\tau/|z-w|^2\le Ct^be^{-d}$; at the outer interface it is at most
$Ct^bN^{-1/2}$. The scalar error in \eqref{annuli:eq:raw-to-g} multiplies an
already bounded compensated quantity. Therefore
\eqref{annuli:eq:master}, \eqref{annuli:eq:raw-to-g} and
\eqref{annuli:eq:Gaussian-compensated} prove both estimates in
Theorem~\ref{annuli:thm:annuli}. A fixed coarse observation can be retained
through all comparisons and then discarded for the upper envelope.

For Corollary~\ref{annuli:cor:ui}, dropping later barriers and cutoffs increases
the nonnegative pair density. A $j$-annulus has $w$-area at most $Ce^{-2j}$
for fixed $z$, canceling the $e^{2j}$ in Theorem~\ref{annuli:thm:annuli}.
Since $s_*-j\ge t/4+O(1)$ in the outer range,
\[
 t\sum_{L-O(1)\le j\le t/2}
  \frac{1}{(1+j)^{3/2}(1+s_*-j)}
 \le C\sum_{j\ge L-O(1)}(1+j)^{-3/2}
 \le C(1+L)^{-1/2}.
\]
There are $O(t)$ annuli, and the integral has its prefactor $t$.
The summed comparison error is therefore at most $Ct^{2-P}$ in
each range. In the terminal range,
\[
 t\sum_{t/2<j\le n_T}
  \frac{1}{(1+j)^{3/2}(1+n_T-j)}
 \le Ct^{-1/2}\log(2+t)\longrightarrow0.
\]
This proves the corollary, first taking $\limsup_N$, then $L\to\infty$.
For the endpoint-cut density used in Section~\ref{chap:fixed},
the remaining core is bounded independently by
\eqref{oneroot:eq:terminal-core}, with contribution $Ct e^{-d}$;
this uses the pointwise cutoff at one root and the one-root expectation
at the other.
\eqref{annuli:eq:Gaussian-annulus} is not asserted at coincident points when
both external roots lie deeper than the final barrier. Likewise the
fixed-distance region $j<k_0$ uses
Theorem~\ref{separated:thm:sm-comparison}. The compensation here is
always at an observed depth $j\ge k_0$.

Proposition~\ref{fixed:prop:fb-diagonal} applies these two annular
sums and the independent core estimate. Its prelimit bound
\eqref{fixed:eq:fb-diagonal-prelimit} retains all four contributions:
the outer tail, the terminal logarithmic sum, the comparison error
and the core. Fixed-distance comparison is used separately in
Section~\ref{fixed:sec:fb-cut-pair}.

\endgroup

\begingroup
\newcommand{\E}{\mathbb E}
\newcommand{\Pp}{\mathbb P}
\newcommand{\R}{\mathbb R}
\newcommand{\C}{\mathbb C}
\newcommand{\dd}{\,\mathrm d}
\newcommand{\ind}{\mathbf 1}
\newcommand{\norm}[1]{\left\|#1\right\|}
\newcommand{\cM}{\mathcal M}
\newcommand{\rhoR}{\rho_{\mathrm R}}
\newcommand{\src}[1]{\texttt{\detokenize{#1}}}
\newcommand{\Var}{\operatorname{Var}}
\section{Fixed-barrier matrix identification and marks}\label{chap:fixed}

In this section we combine the weighted comparisons to identify critical mass, height marks, and their dependence on the field under a fixed barrier.
We prove the fixed-barrier matrix estimates needed for the joint
identification in Corollary~\ref{sp:fb:identification}.
The estimates apply to the original iid ensemble.
The truncated second-moment method also underlies the subcritical
argument in \cite[Section 7.3]{CL}. The critical matrix inputs here
are the one- and two-root comparisons and annular estimates proved
in the preceding sections. The second-moment argument uses an endpoint
cutoff, whose loss tends to zero in expected total variation.
The barrier height is fixed throughout this section.

\subsection{Conventions}
\label{fixed:sec:fb-setting}

Now we fix the barrier, terminal cutoff, and observation lists so that the matrix measures and Gaussian targets use matching conventions.
Let $A_N=N^{-1/2}(\xi_{ab})$ have complex iid entries with
\[
 \E\xi=\E\xi^2=0,\qquad \E|\xi|^2=1,\qquad
 \E|\xi|^q<\infty\quad(q<\infty).
\]
The entry law is fixed, and $\kappa=\E|\xi|^4-2\ge-1$.
We fix a compact disk $K$ strictly inside the unit disk, an integer
$k_0\ge1$, and $B\in\R$. We write $t=\tfrac12\log N$,
$n=\lfloor t\rfloor$, and $d=t^{1/3}$.
The observation lists, test functions, and their derivative bounds
are fixed before $N\to\infty$. Auxiliary comparison parameters
are chosen after these data and the required precision, and remain
fixed in the matrix limit.

We use the deterministic variance-one Hermitized density $\rho^z$
and the corrected fields of \eqref{sp:fb:Phi}--\eqref{sp:fb:endpoints}.
Explicitly,
\begin{align*}
 \Phi_N(z,\eta)
 &=\frac12\log\det\big((A_N-zI)(A_N-zI)^*+\eta^2I\big)
       -\frac N2\int_\R\log(x^2+\eta^2)\rho^z(x)\dd x,\\
 X_{N,s}(z)
 &=\sqrt2\left[\Phi_N(z,\eta_s(z))
                    +\frac\kappa4p_z(u_s)^2\right],\\
 Y_N(z)&=\sqrt2\Phi_N(z,0),\qquad
 Y_N^0(z)=Y_N(z)+\frac{\sqrt2\kappa}{4}(1-|z|^2)^2,
\end{align*}
where
\[
 u_s=\tanh s,\qquad p_z(u)=u(1-u|z|^2),\qquad
 \eta_s(z)=(1-u_s)\sqrt{u_s^{-1}-|z|^2}.
\]
At zero height, $\Phi_N(z,0)=\log|\det(A_N-zI)|
+\tfrac N2(1-|z|^2)$ whenever the determinant is nonzero.
The finite correction does not assert that the fields have zero
expectation at finite $N$.

The critical density, corrected mark, and barrier are
\begin{align*}
 W_N(z)
 &=\frac{e^{2Y_N(z)}}{\E e^{2Y_N(z)}}
   =\frac{e^{2Y_N^0(z)}}{\E e^{2Y_N^0(z)}},&
 R_N^0(z)&=\frac{2t-Y_N^0(z)}{\sqrt t},\\
 A_{N,j}^B(z)
 &=\bigcap_{k=k_0}^j\{X_{N,k}(z)\le2k+B\}.
\end{align*}
The equality of the two densities is exact. We set the density
and the mark to zero at determinant zeros. We then define
\begin{align*}
 W_N^B(z)&=W_N(z)\ind_{A_{N,n}^B(z)},&
 \widehat W_N^B(z)&=W_N^B(z)\ind_{\{Y_N^0(z)\le2t-d\}},\\
 Z_N^B(\dd z)&=\sqrt t\,W_N^B(z)\dd A(z),&
 \widehat Z_N^B(\dd z)&=\sqrt t\,\widehat W_N^B(z)\dd A(z).
\end{align*}
All spatial measures in this section are restricted to $K$, and
$\dd A=\dd x\,\dd y$. The marked measures are
\begin{align}
 \nu_N^B(\dd z,\dd r)=\Xi_N^{B,0}(\dd z,\dd r)
   &:=Z_N^B(\dd z)\,\delta_{R_N^0(z)}(\dd r),\nonumber\\
 \widetilde\nu_N^B(\dd z,\dd r)=\widehat\Xi_N^{B,0}(\dd z,\dd r)
   &:=\widehat Z_N^B(\dd z)\,\delta_{R_N^0(z)}(\dd r).
 \label{fixed:eq:fb-measures}
\end{align}
For a bounded mark $h$, we retain the local notation
\[
 \mu_N^{B,h}(f)=\Xi_N^{B,0}(f\otimes h),\qquad
 \widetilde\mu_N^{B,h}(f)=\widehat\Xi_N^{B,0}(f\otimes h).
\]
In particular, $\mu_N^{B,1}=Z_N^B$ and
$\widetilde\mu_N^{B,1}=\widehat Z_N^B$.
In the later formulas of this section, $U_N=\widehat W_N^B$
and the local abbreviation $R_N$ means the corrected mark $R_N^0$.
The raw mark in the opening section remains the distinct quantity
defined in \eqref{sp:fb:roots}.

We use the exact centered Gaussian family with covariance
\eqref{sp:fb:covariance}, on the realization constructed in
Section~\ref{gaussian:gp2:canonical}. We write $\E_G$ and $\Pp_G$
for its expectation and probability. For real $s>0$, put
\[
 A_j^B(z)=\bigcap_{k=k_0}^j\{G_k(z)\le2k+B\},\qquad
 \mu_s^{G,B}(\dd z)
 =\sqrt s\,e^{2G_s(z)-2C_\kappa(s,z;s,z)}
                   \ind_{A_{\lfloor s\rfloor}^B(z)}\dd A(z).
\]
The root time is real and the barrier times are integers.
Theorem~\ref{gaussian:gp2:gaussian-main} gives, on this same
Gaussian space,
\begin{equation}
 \mu_s^{G,B}\longrightarrow Z_B
 \quad\text{weakly in probability},\qquad
 \mu_s^{G,B}(f)\longrightarrow Z_B(f)\quad\text{in }L^2
 \quad(f\in C(K;\R)).
 \label{fixed:eq:fb-gaussian}
\end{equation}
It also gives $\sup_{s\ge k_0+1}\E_G[\mu_s^{G,B}(K)^3]<\infty$.
The limiting measure is atomless and has the explicit form
\[
 M_\kappa=\sqrt{2/\pi}\,D_\kappa,\qquad
 Z_B(\dd z)=\ind_{A_\infty^B(z)}M_\kappa|_K(\dd z),\qquad
 A_\infty^B(z)=\bigcap_{k=k_0}^\infty\{G_k(z)\le2k+B\}.
\]
Thus the barrier restricts the measure at each spatial point.
We retain $\mu_s^B=\mu_s^{G,B}$ as a local abbreviation.
The normalized one-root laws used below are
\[
 \frac{\dd Q_{N,z}^{\mathrm{det}}}{\dd\Pp}=W_N(z),
 \qquad
 \frac{\dd Q_{s,z}^G}{\dd\Pp_G}
      =e^{2G_s(z)-2C_\kappa(s,z;s,z)}.
\]

We choose a countable dense set $\mathcal D\subset K$ and enumerate
the integer observations $(k,v)$ with $k\ge k_0$ and $v\in\mathcal D$.
Any fixed finite sublist may be enlarged by a fixed finite list
of observations at positive heights in $[a_0,1]$, for some $a_0>0$.
Their sites may lie in a fixed bulk disk $K_*$ whose interior
contains $K$, including sites outside $K$. We write the resulting
finite vectors as follows, decreasing $a_0$ if necessary so that
it bounds every height in the combined list from below:
\[
 C_N^{(m)}=(X_{N,s_i}(v_i))_{i=1}^m,\qquad
 C_G^{(m)}=(G_{s_i}(v_i))_{i=1}^m,
 \qquad \eta_{s_i}(v_i)\in[a_0,1].
\]
Here every $s_i,v_i$ is fixed. Repeated observations are permitted.
The comparisons are applied on $K_*$; the barrier and measures
remain on $K$.

We use the Rayleigh probability measure
\[
 \rhoR(\dd r)=r e^{-r^2/2}\ind_{\{r>0\}}\dd r.
\]
Proposition~\ref{p1:proposition} supplies the Gaussian target.
We now verify the matrix estimates in Proposition~\ref{p2:matrix},
in the order (a), (e), (b), (c), (d).
In particular, marked concentration is \eqref{p2:d}, and the
expected endpoint-cutoff loss is \eqref{p2:e}.
Corollary~\ref{sp:fb:identification} gives the joint identification
from these estimates.

\subsection{Coarse observations and first moments}
\label{fixed:sec:fb-first}

Next we match first moments with retained field observations, providing the link between the critical mass and its Gaussian field.

\begin{lemma}[Fixed positive heights; consequence of Cipolloni--Landon]
\label{fixed:lem:positive-height}
Fix a closed bulk disk $K_*$, $a_0\in(0,1]$, and an integer $m\ge1$.
For a list $(v_i,\eta_i)\in K_*\times[a_0,1]$, let $s_i>0$ be
specified by $\eta_{s_i}(v_i)=\eta_i$, and set
\[
 V_N=(X_{N,s_i}(v_i))_{i=1}^m,\qquad
 V_G=(G_{s_i}(v_i))_{i=1}^m.
\]
For every fixed such list, $V_N\Rightarrow V_G$ under the unweighted
matrix law. For every real polynomial $P$ on $\mathbb R^m$,
\begin{equation}
 \sup_{(v_i,\eta_i)\in(K_*\times[a_0,1])^m}
       |\E P(V_N)-\E_G P(V_G)|\longrightarrow0.
 \label{fixed:eq:positive-height-moments}
\end{equation}
Moreover, for every fixed integer $q\ge1$ the $q$th absolute moments
of $V_N$ are bounded uniformly over these lists and all sufficiently
large $N$. Repeated observations, approaching points, unequal heights
and height one are allowed. The bounds are uniform for entry laws
with common moment bounds; the convergence uses the fixed entry
law and its fourth cumulant $\kappa$. The number of observations and
$a_0$ remain fixed before $N\to\infty$.
\end{lemma}

\begin{proof}
This is a fixed-positive-height consequence of
\cite[Proposition 2.15, Lemmas 8.1--8.2]{CL}. We give the cap and
moment details needed at both applications below.

\paragraph{The unweighted joint law.}
Apply \cite[Proposition 2.15]{CL} with its mesoscopic observation
list equal to the present fixed list, no additional observations,
and any one fixed exponent $\varepsilon\in(0,1)$, for example
$\varepsilon=1/2$. The lower-height condition holds eventually
because $a_0>0$. Write $B_{\rm tail}>0$ for the constant shared by
\cite[Proposition 2.15 and Lemma 8.1]{CL}. Choose a bounded neighborhood
of zero for the real Laplace parameters, and let $L_{\rm CL}>0$ bound
the absolute values of its coordinates. Then choose a single sufficiently
large fixed cap $A$ for that neighborhood and the whole list. The formula (2.38) applies
to
\[
 \mathcal F_N=\bigcap_{i=1}^m
       \{|\Phi_N(v_i,\eta_i)|\le A\log N\}.
\]
The cap can be enlarged so that the explicit error
$N^{mL_{\rm CL}A-A^2/B_{\rm tail}}$ in that formula vanishes; Lemma 8.1 also gives
$\Pp(\mathcal F_N^c)\to0$, uniformly over the lists. Normalize the
capped subprobability law by $\Pp(\mathcal F_N)$. Its Laplace
transforms then converge on that neighborhood to those of the
Gaussian vector with mean and covariance from
\cite[(2.33)--(2.34)]{CL}. Laplace-transform continuity gives the
conditional joint-law limit. Removing the cap changes the expectation
of a bounded test by at most twice its supremum norm times
$\Pp(\mathcal F_N^c)$, and therefore gives the unweighted joint law.
This argument does not require the limiting covariance to be invertible.

For the complex ensemble, the mean in (2.34) is
$-\kappa p_{v_i}(u_{s_i})^2/4$, since
$m^{v_i}(i\eta_i)^2=-p_{v_i}(u_{s_i})$. The covariance in (2.33)
is $C_\kappa(s_i,v_i;s_j,v_j)/2$. The deterministic correction and
factor $\sqrt2$ defining $X_{N,s_i}(v_i)$ hence give exactly the
centered vector $V_G$. This checks the normalization without using
critical weights or fixed-barrier identification.

\paragraph{Actual moments.}
Let $\|\cdot\|_{\rm F}$ denote the unnormalized Frobenius norm.
At heights in $[a_0,1]$, the determinant and deterministic density
bounds give, uniformly in $v\in K_*$,
\[
 |\Phi_N(v,\eta)|\le
 C_{a_0,K_*}N[1+\log(2+\|A_N\|_{\rm F})].
\]
The entry moments consequently imply, for every fixed integer $q\ge1$,
\[
 \sup_{v,\eta}\E|\Phi_N(v,\eta)|^{2q}
       \le C_{a_0,K_*,q}N^{3q}.
\]
Given $M>3q$, choose a sufficiently large fixed $A\ge100$.
Write $B_{\rm Lap}>0$ for the constant in \cite[Lemma 8.2]{CL}.
Lemma 8.1 and Lemma 8.2, equation (8.33), of \cite{CL} then give
\begin{align*}
 \sup_{v,\eta}\Pp\{|\Phi_N(v,\eta)|>A\log N\}
       &\le C N^{-M},\\
 \sup_{v,\eta}\E\!\left[e^{\pm\Phi_N(v,\eta)}
       \ind_{\{|\Phi_N(v,\eta)|\le A\log N\}}\right]&\le C.
\end{align*}
Indeed, at fixed positive heights the Gaussian mean and variance parameters
in the cited formula are bounded; use $L_{\rm CL}=2$ in (8.33), enlarge $A$ until both
$A^2/B_{\rm tail}>M$ and $10+4B_{\rm Lap}+2A-A^2/B_{\rm Lap}<0$,
and choose its residual
accuracy greater than $M$. All these choices precede the matrix
limit. Cauchy--Schwarz bounds the $q$th-moment contribution outside
the cap by $C N^{(3q-M)/2}=o(1)$. Inside the cap, the inequality
$|x|^q\le C_q(e^x+e^{-x})$ gives a uniform bound. The deterministic
correction to $X_N$ is bounded. Repeating at any larger fixed order
therefore gives uniform integrability of every prescribed polynomial
in the fixed vector; H\"older's inequality handles its products.

Finally, if \eqref{fixed:eq:positive-height-moments} failed, choose
a failing sequence of lists and pass to a convergent subsequence
in $(K_*\times[a_0,1])^m$. The uniform capped Laplace formula gives
the Gaussian law along this sequence, including coincident limiting
observations. The preceding uniform integrability gives its polynomial
moments. Continuity of the positive-height covariance
\eqref{sp:fb:covariance} gives the same limiting moments for $V_G$,
a contradiction. This proves the uniform assertion as well as the
fixed-list convergence.
\end{proof}

\begin{lemma}
\label{fixed:lem:fb-first}
For every fixed finite observation vector, real $f\in C(K)$,
and real $\Psi\in C_b^\infty(\R^m)$ with all derivatives bounded,
\begin{align}
 C_N^{(m)}&\Rightarrow C_G^{(m)},\label{fixed:eq:fb-M1}\\
 \E[\widehat Z_N^B(f)\Psi(C_N^{(m)})]
 -\E_G[\mu_t^{G,B}(f)\Psi(C_G^{(m)})]&\longrightarrow0.
 \label{fixed:eq:fb-M2}
\end{align}
Consequently,
\begin{equation}
 \E[\widehat Z_N^B(f)\Psi(C_N^{(m)})]
 \longrightarrow \E_G[Z_B(f)\Psi(C_G^{(m)})].
 \label{fixed:eq:fb-first-limit}
\end{equation}
The target $Z_B$ is measurable with respect to the completed
sigma-field generated by the countable integer observation family.
\end{lemma}

\begin{proof}
We prove the coarse convergence first, then
\eqref{p2:e}, and finally the mixed first moment.
These are clauses (a), (e), and (b) of Proposition~\ref{p2:matrix}.

\paragraph{Coarse convergence.}
Lemma~\ref{fixed:lem:positive-height} applies to the prescribed fixed
list and proves \eqref{fixed:eq:fb-M1}, including repeated observations
and height one. This unweighted input is separate from the critical
mixed first moment proved below.

\paragraph{Endpoint cutoff.}
The density $W_N^0$ in Theorem~\ref{oneroot:thm:one-root}
equals $W_N$ exactly. Its full integer barrier and corrected
endpoint cutoff are also the ones used here.
The uniform estimate \eqref{oneroot:eq:cutoff-prob} therefore gives
\begin{equation}
 \alpha_N:=\sqrt t\,\sup_{z\in K}
 Q_{N,z}^{\mathrm{det}}
   \{A_{N,n}^B(z),\ Y_N^0(z)>2t-d\}
 \le C_B t^{-1/3}\sqrt{\log(2+t)}+o(1)
 \longrightarrow0.
 \label{fixed:eq:fb-cutoff-pointwise}
\end{equation}
The measure identification above and
\eqref{oneroot:eq:cutoff-tv} give
\begin{equation}
 \E\norm{Z_N^B-\widehat Z_N^B}_{\mathrm{TV}}
 \le |K|\alpha_N\longrightarrow0.
 \label{fixed:eq:fb-cutoff-integral}
\end{equation}
This proves \eqref{p2:e} and clause (e).
Its use below concerns bounded first-moment tests.

\paragraph{Mixed first moment.}
We apply Theorem~\ref{oneroot:thm:one-root} with terminal
mark equal to one. All prescribed spatial observations, including
the fixed integer observations, are placed in its additional coarse
list. We let the insertion ignore the root-prefix coordinate.
The additional coordinates in that theorem are
\[
 U_{N,i}
 =\Phi_N(v_i,\eta_{s_i}(v_i))
                  -E_\kappa(v_i,\eta_{s_i}(v_i))
 =\frac{X_{N,s_i}(v_i)}{\sqrt2},\qquad
 U_{G,i}=\frac{G_{s_i}(v_i)}{\sqrt2}.
\]
We therefore use $H(u)=\Psi(\sqrt2\,u)$.
The list is fixed and all required derivatives of $H$ are bounded.
Every height is in a fixed interval $[a_0,1]$; the sites lie in
$K_*$ and need not be separated from the varying root $z\in K$.
Repeated observations are permitted.

The full-history comparison, equation \eqref{oneroot:eq:or22} in
Section~\ref{oneroot:completing-the-final-integer-barrier},
now yields
\begin{align}
 \varepsilon_N^{(1)}:=\sqrt t\,\sup_{z\in K}\Big|
 &\E_{Q_{N,z}^{\mathrm{det}}}
       [\Psi(C_N^{(m)})\ind_{A_{N,n}^B(z)}]\nonumber\\
 &-\E_{Q_{t,z}^G}
       [\Psi(C_G^{(m)})\ind_{A_n^B(z)}]\Big|
 \longrightarrow0.
 \label{fixed:eq:fb-root-comparison}
\end{align}
The Gaussian root is at the real time $t$, and both barriers
retain every integer time through $n$.
The tests are bounded and each normalized root has expectation one,
so spatial integration gives
\begin{equation}
 \left|\E[Z_N^B(f)\Psi(C_N^{(m)})]
       -\E_G[\mu_t^{G,B}(f)\Psi(C_G^{(m)})]\right|
 \le\varepsilon_N^{(1)}\int_K|f|\dd A\longrightarrow0.
 \label{fixed:eq:fb-integrated-first}
\end{equation}
The cutoff estimate already proved gives
\[
 \left|\E[(Z_N^B-\widehat Z_N^B)(f)\Psi(C_N^{(m)})]\right|
 \le\norm f_\infty\norm\Psi_\infty\,
             \E\norm{Z_N^B-\widehat Z_N^B}_{\mathrm{TV}}
 \longrightarrow0.
\]
Together these prove \eqref{fixed:eq:fb-M2}.

Finally, \eqref{fixed:eq:fb-gaussian} holds on the Gaussian space
that also carries $C_G^{(m)}$. Hence
\[
 \left|\E_G[(\mu_t^{G,B}(f)-Z_B(f))\Psi(C_G^{(m)})]\right|
 \le\norm\Psi_\infty
             \norm{\mu_t^{G,B}(f)-Z_B(f)}_{L^2(\Pp_G)}
 \longrightarrow0.
\]
This proves \eqref{fixed:eq:fb-first-limit}, and thus clause (b),
while retaining the dependence between the target and the observations.

The measurability assertion is the integer-history conclusion
proved in Section~\ref{gaussian:gp2:coarse}.
That section also proves the finite-coordinate consequence
\eqref{p1:finite-approximation}. Appending prescribed observations,
including ones outside $K$, does not assert their measurability
with respect to the integer-history sigma-field on $K$.
\end{proof}

\begin{remark}
The independence used here is the conclusion of
Theorem~\ref{noise:thm:lfn-whitening} and the following remark on independence, with
initial sigma-fields defined in \eqref{noise:eq:lfn-initial-sigma-fields}.
The Brownian vector is independent of the enlarged initial field; the
reference processes and their local factors, which use their own initial
spectra, are independent of the original matrix field
$\mathcal F_0^{\mathrm{mat}}$.
Before integrating out a local factor, the history and coarse observations
are returned to that original initial matrix. Arbitrary bounded tests of
those observations are therefore allowed in this independence step.
\end{remark}

\subsection{Pair moments at fixed positive separation}
\label{fixed:sec:fb-cut-pair}

In this subsection we compare weighted pairs a fixed distance apart and show that the terminal cutoffs do not change their limiting contribution.

The square of a spatial integral is an integral over pairs of points.
At a fixed positive distance, we use the separated comparison.
The contribution from nearby points is controlled by a positive
upper bound in the next section.

Theorem~\ref{separated:thm:sm-comparison} gives an absolute error
$o(t^{-1})$ for the uncut, full-barrier pair expectations, uniformly
on $z,w\in K$ with $|z-w|\ge r>0$.
It permits bounded Lipschitz marks and a bounded smooth test of finitely
many coarse observations, with the derivative bounds required there.
The coarse heights lie in a fixed compact subinterval of $(0,1)$.
The exact pair factors are
\[
 a_N(z,w)=\E[W_N(z)W_N(w)],\qquad
 a_t^G(z,w)=e^{4C_\kappa(t,z;t,w)}.
\]
Both are included in the comparison and are bounded above and below
at each fixed separation. We write $Q_N^{(2)}$ for the normalized
matrix pair law, whose density is $W_N(z)W_N(w)/a_N(z,w)$.

\begin{lemma}[The cutoff is negligible for separated pair moments]
\label{fixed:lem:fb-cut-pair}
For each fixed $r>0$,
\begin{align}
 0&\le t\sup_{\substack{z,w\in K\\|z-w|\ge r}}
 \E\!\left[W_N^B(z)W_N^B(w)
       -\widehat W_N^B(z)\widehat W_N^B(w)\right]\nonumber\\
 &\le C_{r,B}t^{-1/3}(1+\log t)+o(1)\longrightarrow0.
 \label{fixed:eq:fb-separated-cut}
\end{align}
Consequently the separated comparison remains valid with both
matrix endpoint cutoffs and no Gaussian endpoint cutoffs, for
the marks and coarse tests stated above.
\end{lemma}

\begin{proof}
The weight $W_N^0$ in Corollary~\ref{separated:cor:sm-cut}
is exactly $W_N$: the deterministic correction cancels from
the normalized density. Its full barriers and cutoff
$Y_N^0\le2t-t^{1/3}$ are also the ones used here.
Thus \eqref{fixed:eq:fb-separated-cut} is precisely
\eqref{separated:eq:sm-cut}.

The proof of that corollary retains both barriers and gives, for
$v=z,w$ and any fixed $P>3$ chosen before its auxiliary parameters,
\begin{align}
 &Q_N^{(2)}\!\left(A_{N,n}^B(z)\cap A_{N,n}^B(w),\
                      Y_N^0(v)>2t-d\right)\nonumber\\
 &\qquad\le C_{r,B}t^{-4/3}(1+\log t)+o(t^{-P}),
 \label{fixed:eq:fb-joint-strip}
\end{align}
uniformly at the fixed separation. Summing over the two possible
failed cutoffs and multiplying by $t\,a_N(z,w)$ explains the
rate in \eqref{fixed:eq:fb-separated-cut}.
Both barriers are needed for this bound.

The difference of the two pair weights in
\eqref{fixed:eq:fb-separated-cut} is nonnegative. Inserting bounded
marks and a bounded coarse test therefore multiplies its absolute
error by at most the product of their supremum norms.
Combining this with Theorem~\ref{separated:thm:sm-comparison}
proves the asserted comparison.
\end{proof}

We record the absolute comparison that will be integrated.
Put
\[
 W_t^{G,B}(z)
 =e^{2G_t(z)-2C_\kappa(t,z;t,z)}\ind_{A_n^B(z)},\qquad
 R_t^G(z)=\frac{2t-G_t(z)}{\sqrt t}.
\]
For bounded Lipschitz $g_1,g_2$ and an allowed coarse test $\Psi$,
\begin{align}
 t\sup_{\substack{z,w\in K\\|z-w|\ge r}}\Big|
 &\E[\widehat W_N^B(z)\widehat W_N^B(w)
       g_1(R_N^0(z))g_2(R_N^0(w))\Psi(C_N)]\nonumber\\
 &-\E_G[W_t^{G,B}(z)W_t^{G,B}(w)
       g_1(R_t^G(z))g_2(R_t^G(w))\Psi(C_G)]
 \Big|\longrightarrow0.
 \label{fixed:eq:fb-separated-comparison}
\end{align}
Here $C_N,C_G$ are the finite coarse vectors in
Theorem~\ref{separated:thm:sm-comparison}, with its stated height domain.
The comparison includes the exact pair factors before spatial
integration. The quadratic arguments below use $\Psi=1$.
The full observation list in \eqref{sp:coarse-vectors},
including height one, is retained through the mixed first
moments in \eqref{fixed:eq:fb-M2}.

\subsection{Near-diagonal uniform integrability}
\label{fixed:sec:fb-diagonal}

With the terminal cutoff retained, we now show that the second-moment contribution vanishes as the diagonal neighborhood shrinks after the matrix limit.

\begin{proposition}[Near-diagonal bound with endpoint cutoffs]
\label{fixed:prop:fb-diagonal}
For integers $L\ge k_0+2$,
\begin{equation}
 \limsup_{N\to\infty}
 t\iint_{\substack{(z,w)\in K^2\\|z-w|\le e^{-L}}}
       \E[U_N(z)U_N(w)]\dd A(z)\dd A(w)
 \le C_B(1+L)^{-1/2}.
 \label{fixed:eq:fb-diagonal}
\end{equation}
Here $U_N=\widehat W_N^B$ is the endpoint-cut density.
\end{proposition}

\begin{proof}
We denote the integral before the limit in
\eqref{fixed:eq:fb-diagonal} by $I_N(L)$.
Fix $P>2$, and then choose the fixed parameters of
Theorem~\ref{annuli:thm:annuli}. We put
\[
 s_*=\lfloor3t/4\rfloor,\qquad
 T=\lfloor t-d/2-C_0\rfloor,
\]
where $C_0$ is the constant in that theorem.
For each fixed $L$ and sufficiently large $N$,
$L\le\lfloor t/2\rfloor<T\le n$ and $s_*\le n$.
We partition the integration region into
\[
 \mathcal A_j
 =\{(z,w)\in K^2:e^{-(j+1)}<|z-w|\le e^{-j}\},
 \qquad L\le j\le T,
\]
and the remaining disk $|z-w|\le e^{-(T+1)}$ in the $w$ variable.

For $L\le j\le\lfloor t/2\rfloor$, the outer estimate
\eqref{annuli:eq:outer} gives, on $\mathcal A_j$,
\begin{equation}
 e^{-2j}\E\!\left[W_N(z)W_N(w)
             \ind_{A_{N,s_*}^B(z)\cap A_{N,s_*}^B(w)}\right]
 \le\frac{C_B}{(1+j)^{3/2}(1+s_*-j)}+C_Bt^{-P}.
 \label{fixed:eq:fb-outer}
\end{equation}
The weight $W_N^0$ in that estimate is again $W_N$.
The full barriers and the endpoint cuts only decrease this
positive integrand. For each $z$, the annulus has area
at most $Ce^{-2j}$, and $t/(1+s_*-j)\le C$.
Its integral can therefore be summed to give the bound
\begin{equation}
 C_B\sum_{j=L}^{\lfloor t/2\rfloor}(1+j)^{-3/2}
       +C_Bt^{2-P}
 \le C_B'(1+L)^{-1/2}+C_Bt^{2-P}
 \label{fixed:eq:fb-outer-sum}
\end{equation}
for the total contribution of these annuli.

For $\lfloor t/2\rfloor<j\le T$, the terminal estimate
\eqref{annuli:eq:terminal} gives \eqref{fixed:eq:fb-outer}
with $s_*$ replaced by $T$. Its total contribution is at most
\begin{equation}
 C_Bt\sum_{\lfloor t/2\rfloor<j\le T}
       \frac{1}{(1+j)^{3/2}(1+T-j)}
       +C_Bt^{2-P}
 \le C_B't^{-1/2}\log(2+t)+C_Bt^{2-P}.
 \label{fixed:eq:fb-terminal-sum}
\end{equation}
In both sums, the annular area cancels the factor $e^{2j}$
in the pointwise bound. The error $t^{2-P}$ accounts for
$O(t)$ annuli and the factor $t$ in $I_N(L)$.

The remaining disk satisfies
\[
 e^{-2(T+1)}\le e^{2C_0}N^{-1}e^d.
\]
Applying \eqref{oneroot:eq:terminal-core} to this remaining disk gives
\begin{align}
 &t\iint_{\substack{(z,w)\in K^2\\|z-w|^2\le CN^{-1}e^d}}
       \E[U_N(z)U_N(w)]\dd A(z)\dd A(w)\nonumber\\
 &\qquad\le Ct\,(N^{-1}e^d)(Ne^{-2d})=Cte^{-d}.
 \label{fixed:eq:fb-core}
\end{align}
At coincidence, the Gaussian estimate \eqref{annuli:eq:Gaussian-general}
has an extra factor when the external roots lie beyond the last
retained barrier. That factor cannot be discarded to extend the
annular bound to coincidence.

Combining \eqref{fixed:eq:fb-outer-sum},
\eqref{fixed:eq:fb-terminal-sum}, and \eqref{fixed:eq:fb-core}, we have
\begin{equation}
 I_N(L)\le C_{B,P}\left[
 (1+L)^{-1/2}+t^{-1/2}\log(2+t)+t^{2-P}+te^{-d}
 \right]
 \label{fixed:eq:fb-diagonal-prelimit}
\end{equation}
for each fixed $L$ and all sufficiently large $N$.
Since $P>2$ and $d=t^{1/3}$, taking $N\to\infty$ proves
\eqref{fixed:eq:fb-diagonal}; we may fix, for example, $P=3$
to obtain the stated constant $C_B$.
The diagonal neighborhood is then allowed to shrink by taking
$L\to\infty$.
\end{proof}

\subsection{Second moments and concentration of the marks}
\label{fixed:sec:fb-quadratic}

 With the terminal cutoff retained, we now match spatial second moments and show that height tests average against the Rayleigh law.
We first combine the pair estimates with the Gaussian limit.
The spatial second moment uses both the matrix and Gaussian
diagonal bounds. For centered marks, the separated Gaussian
factorization and the matrix diagonal bound suffice.

\begin{proposition}
\label{fixed:prop:fb-moments}
For every real $f\in C(K)$,
\begin{equation}
 \E[\widehat Z_N^B(f)^2]
       -\E_G[\mu_t^{G,B}(f)^2]\longrightarrow0,
 \label{fixed:eq:fb-M3}
\end{equation}
and consequently
\begin{equation}
 \E[\widehat Z_N^B(f)^2]\longrightarrow\E_G[Z_B(f)^2].
 \label{fixed:eq:fb-second-limit}
\end{equation}
For real $h\in C_c^\infty(\R)$, \eqref{p2:d} holds.
In particular, these conclusions give clauses (c) and (d) of
Proposition~\ref{p2:matrix}.
\end{proposition}

\begin{proof}
\emph{The Gaussian contribution near the diagonal.}
The convergence in \eqref{fixed:eq:fb-gaussian} takes place on
the same Gaussian probability space. For real $a,b\in C(K)$,
Cauchy--Schwarz therefore gives
\[
 \E_G[\mu_t^{G,B}(a)\mu_t^{G,B}(b)]
       \longrightarrow\E_G[Z_B(a)Z_B(b)].
\]
Finite sums of $a(z)b(w)$ are uniformly dense in $C(K^2;\R)$.
The constant-test case of \eqref{fixed:eq:fb-gaussian} bounds
$\E_G[\mu_t^{G,B}(K)^2]$ for all sufficiently large $t$,
and $Z_B(K)\in L^2$. Uniform approximation now gives
\begin{equation}
 \E_G[(\mu_t^{G,B}\otimes\mu_t^{G,B})(q)]
       \longrightarrow\E_G[(Z_B\otimes Z_B)(q)],
 \qquad q\in C(K^2;\R).
 \label{fixed:eq:fb-gaussian-product}
\end{equation}

Choose a continuous nonincreasing function
$\chi:[0,\infty)\to[0,1]$, equal to one on $[0,1]$ and
zero on $[2,\infty)$, and put
$\chi_r(z,w)=\chi(|z-w|/r)$.
For fixed $r>0$, \eqref{fixed:eq:fb-gaussian-product} applies
to $\chi_r$. Since $Z_B$ is atomless,
$(Z_B\otimes Z_B)(\chi_r)\to0$ almost surely as $r\downarrow0$.
These quantities are bounded by the integrable $Z_B(K)^2$.
Thus
\begin{equation}
 \lim_{r\downarrow0}\lim_{t\to\infty}
 \E_G[(\mu_t^{G,B}\otimes\mu_t^{G,B})(\chi_r)]=0.
 \label{fixed:eq:fb-gaussian-diagonal}
\end{equation}
Since $\chi_r=1$ when $|z-w|\le r$, this also controls the
expected Gaussian product measure of that neighborhood.

\emph{Second moments.}
For all sufficiently large $N$, the density with the endpoint cutoff
is bounded and the Gaussian total mass has a finite second moment.
The pair integrals below are therefore absolutely integrable.
In particular,
\begin{equation}
 \E[\widehat Z_N^B(f)^2]
 =t\iint_{K^2} f(z)f(w)\E[U_N(z)U_N(w)]\dd A(z)\dd A(w).
 \label{fixed:eq:fb-second-square}
\end{equation}
The Gaussian identity has $W_t^{G,B}$ in place of $U_N$ and
$\E_G$ in place of $\E$.

Fix an integer $L\ge k_0+2$ and put $r_L=e^{-L}/2$.
We split both integrals using
$1=(1-\chi_{r_L})+\chi_{r_L}$.
The first part is supported where $|z-w|\ge r_L$.
Equation~\eqref{fixed:eq:fb-separated-comparison}, with both
marks and the coarse insertion equal to one, applies there.
Multiplying its absolute error by $|f(z)f(w)|$ and integrating
gives a quantity tending to zero for each fixed $L$.
The support of $\chi_{r_L}$ lies in $|z-w|\le e^{-L}$.
Using $I_N(L)$ from the preceding proof, we obtain
\begin{align}
 &\left|\E[\widehat Z_N^B(f)^2]
                 -\E_G[\mu_t^{G,B}(f)^2]\right|\nonumber\\
 &\quad\le o_L(1)+\norm f_\infty^2
 \left\{I_N(L)+
 \E_G[(\mu_t^{G,B}\otimes\mu_t^{G,B})(\chi_{r_L})]\right\}.
 \label{fixed:eq:fb-second-difference}
\end{align}
Here $o_L(1)$ tends to zero as $N\to\infty$ with $L$ fixed.
The estimate allows $f$ to change sign.
We first take $N\to\infty$, and then $L\to\infty$.
Equations~\eqref{fixed:eq:fb-diagonal} and
\eqref{fixed:eq:fb-gaussian-diagonal} prove
\eqref{fixed:eq:fb-M3}.
The Gaussian $L^2$ convergence in \eqref{fixed:eq:fb-gaussian}
then gives \eqref{fixed:eq:fb-second-limit}, which is stronger
than clause (c).

\emph{Centered marks.}
We set $\bar h=h-\rhoR(h)$.
This function may lose compact support, but it remains bounded
and Lipschitz, with bounded derivatives, and $\rhoR(\bar h)=0$.
Theorem~\ref{gaussian:gpm:factorization} permits these marks
and the coarse insertion one. Its version for the product
of individually variance-normalized roots gives, for each
fixed $r>0$,
\begin{align}
 t\sup_{\substack{z,w\in K\\|z-w|\ge r}}
 \left|\E_G[W_t^{G,B}(z)W_t^{G,B}(w)
       \bar h(R_t^G(z))\bar h(R_t^G(w))]\right|
 \longrightarrow0.
 \label{fixed:eq:fb-gaussian-centered}
\end{align}
The unmarked term in that factorization is zero because its
coefficient is $\rhoR(\bar h)^2$.
Combining this with \eqref{fixed:eq:fb-separated-comparison},
with $g_1=g_2=\bar h$ and $\Psi=1$, yields
\begin{equation}
 t\sup_{\substack{z,w\in K\\|z-w|\ge r}}
 \left|\E[U_N(z)U_N(w)
       \bar h(R_N^0(z))\bar h(R_N^0(w))]\right|
 \longrightarrow0.
 \label{fixed:eq:fb-matrix-centered}
\end{equation}

The boundedness of the mark and cut density justifies the
square expansion
\begin{align}
 &\E\left|\widehat\Xi_N^{B,0}(f\otimes h)
                 -\rhoR(h)\widehat Z_N^B(f)\right|^2\nonumber\\
 &\quad=t\iint_{K^2}f(z)f(w)
 \E[U_N(z)U_N(w)\bar h(R_N^0(z))\bar h(R_N^0(w))]
                  \dd A(z)\dd A(w).
 \label{fixed:eq:fb-centered-square}
\end{align}
For fixed $L$, we use the same decomposition with $\chi_{r_L}$.
Equation~\eqref{fixed:eq:fb-matrix-centered} makes the
contribution with $1-\chi_{r_L}$ tend to zero.
The absolute contribution with $\chi_{r_L}$ is bounded by
$\norm f_\infty^2\norm{\bar h}_\infty^2I_N(L)$.
Consequently \eqref{fixed:eq:fb-diagonal} gives
\[
 \limsup_{N\to\infty}
 \E\left|\widehat\Xi_N^{B,0}(f\otimes h)
                 -\rhoR(h)\widehat Z_N^B(f)\right|^2
 \le C_B\norm f_\infty^2\norm{\bar h}_\infty^2(1+L)^{-1/2}.
\]
Letting $L\to\infty$ proves \eqref{p2:d}
and clause (d).
\end{proof}

\subsection{Joint identification at a fixed barrier}
\label{fixed:sec:fb-projection}

Now we combine the moment estimates to identify the restricted measure and its marks jointly with the prescribed field observations.
Sections~\ref{fixed:sec:fb-first}--\ref{fixed:sec:fb-quadratic}
prove Proposition~\ref{p2:matrix} for the original ensemble.
Together with Proposition~\ref{p1:proposition}, these estimates
give the joint convergence in Corollary~\ref{sp:fb:identification}.
This includes the spatial, corrected-marked, and coarse limits
for both cut and uncut measures, and the uncut raw-marked limit.
The measure and the observations in the limit use the same
Gaussian family. The observation lists are those in
\eqref{sp:coarse-vectors}, including height one, repeated
observations, and the stated sites outside $K$.
The proof of the corollary establishes weak convergence of
finite positive measures, including total mass, and tightness
in the mark coordinate.

Marked concentration and expected cutoff loss are the separate
assertions \eqref{p2:d} and \eqref{p2:e}, proved in
Sections~\ref{fixed:sec:fb-quadratic} and
\ref{fixed:sec:fb-first}, respectively.
The full finite-time second-moment comparison and its limiting
form are \eqref{fixed:eq:fb-M3} and
\eqref{fixed:eq:fb-second-limit}.
These moment estimates concern the endpoint-cut matrix measures
at fixed $B$.

\endgroup

\begingroup
\newcommand{\E}{\mathbb E}
\newcommand{\Pp}{\mathbb P}
\newcommand{\dd}{\,\mathrm d}
\newcommand{\ii}{\mathrm i}
\newcommand{\eps}{\varepsilon}
\newcommand{\ind}{\mathbf1}
\newcommand{\cA}{\mathcal A}
\newcommand{\cV}{\mathcal V}
\newcommand{\cD}{\mathfrak D}
\newcommand{\cB}{\mathfrak A}
\newcommand{\norm}[1]{\left\|#1\right\|}
\newcommand{\av}[1]{\langle #1\rangle}
\newcommand{\Tr}{\operatorname{Tr}}
\newcommand{\cum}{\operatorname{cum}}
\section{Gaussian comparison for many observations and a spatial barrier bound}\label{chap:aggregate}

The purpose of this section is to compare growing observation lists with Gaussian vectors to control the barrier throughout the spatial disk on an initial fraction of scales.

A barrier imposed at every spatial point requires simultaneous control
of a growing list of observations. This section proves a Gaussian
comparison for lists whose size grows at most as a fixed power of the
matrix size. For each derivative order, all coordinate derivatives are evaluated at
one common point, and the coordinate sum is taken before the supremum.
The error bound has
no additional factor from the number of observations beyond this norm.
The norm itself may grow with the size of the list.
The matrices need no Gaussian divisibility.

We first fix the common event and finite-integral coordinates, derive
the summed raw identity, and identify its second-cumulant covariance.
The third and fourth cumulants and the common-point remainder then
give the summed Stein identity. The Ornstein--Uhlenbeck semigroup
argument and the spatial mesh yield the barrier bound under the unweighted matrix law through
an initial fraction of scales in Proposition~\ref{p34:global}. Throughout, the lower cutoff is a fixed
power above $N^{-1}$.

The external local-law inputs are Cipolloni--Landon
\cite{CL}, Theorems~2.3 and~2.5, Lemma~2.1(ii), Lemma~B.1,
Proposition~B.2(i), equation~(B.15), and Proposition~9.1.
The last may be replaced by the rank-two argument below.
The raw identity and endpoint formulas are derived with the
conventions of their equations~(8.8) and~(8.15); no separately
extracted mean is added to the complete cumulant expansion.
The unweighted two-resolvent estimates for arbitrary spatial points
are used only at fixed-power cutoffs. The last corollary uses the Gaussian global maximum estimate proved
earlier in the paper.

The comparison extends the fixed-list cumulant method of
\cite[Proposition~2.15 and Section~8.1]{CL} to the polynomial lists
specified below. The coefficient bounds are summed at one common
evaluation point, yielding the derivative norm
\eqref{aggregate:eq:norms}; the cutoff remains a fixed power above
$N^{-1}$.
\subsection{Gaussian comparison for polynomially many observations}
\label{aggregate:sec:theorem}

Next we quantify Gaussian comparison for polynomially many observations using a summed derivative norm that avoids an additional factor from the list length.
Let \(A_N=N^{-1/2}(\xi_{ab})\), where the entries have a complex
i.i.d.\ law, possibly depending on $N$, satisfying
\[
 \E\xi=\E\xi^2=0,\qquad \E|\xi|^2=1,\qquad
 \sup_N\E|\xi|^p<\infty\quad\text{for every fixed }p<\infty.
\]
Put \(\kappa=\E|\xi|^4-2\). Fix a compact bulk disk
\(K_+\Subset\{z:|z|<1\}\), \(c_0>0\), and \(0<\eps\le1\).
Fix also $L<\infty$. The deterministic list satisfies
\begin{equation}
 d_N\le N^L,\qquad z_i\in K_+,\qquad
 c_0N^{-1+\eps}\le\eta_i\le1.
 \label{aggregate:eq:list}
\end{equation}
Repeated points and repeated cutoffs are permitted.

Let
\[
 H^z=\begin{pmatrix}0&A_N-zI\\(A_N-zI)^*&0\end{pmatrix},
 \qquad G^z(\ii u)=(H^z-\ii u)^{-1},
 \qquad \av T=(2N)^{-1}\Tr T.
\]
Write
\[
 M^z(w)=\begin{pmatrix}m^z(w)&-zu^z(w)\\
             -\overline z u^z(w)&m^z(w)\end{pmatrix}\otimes I_N,
 \qquad u^z(w)=\frac{m^z(w)}{w+m^z(w)},
\]
where \(m=m^z(w)\) solves
\[
 -\frac1m=w+m-\frac{|z|^2}{w+m}
\]
with the Stieltjes-transform branch. On the positive imaginary axis,
\(m\) is imaginary and \(u^z\) is real.

We use the following deterministic centering,
\[
 \Phi_N(z,\eta)=\frac12\log\det\bigl((A_N-z)(A_N-z)^*+\eta^2I\bigr)
       -\frac N2\int\log(x^2+\eta^2)\rho^z(x)\dd x.
\]
For the list in \eqref{aggregate:eq:list}, set
\begin{equation}
 m_i=m^{z_i}(\ii\eta_i),\quad u_i=u^{z_i}(\ii\eta_i),\quad
 E_i=-\frac\kappa4m_i^4,\quad
 Y_i=\Phi_N(z_i,\eta_i)-E_i,\quad X_i=\sqrt2Y_i.
 \label{aggregate:eq:field}
\end{equation}
The quantity \(E_i\) is an asymptotic deterministic drift, not the
exact finite-\(N\) expectation.
When an argument is displayed, write $m_i(w)=m^{z_i}(w)$ and
$u_i(w)=u^{z_i}(w)$; argument-free $m_i,u_i$ retain the endpoint
values above. A prime denotes differentiation with respect to the
complex spectral variable $w$.

Define
\begin{align}
 D_{ij}&=(1-u_iu_jz_i\overline z_j)
           (1-u_iu_j\overline z_i z_j)-m_i^2m_j^2,\\
 V_{ij}&=-\frac14\log D_{ij},\\
 C^\Phi_{ij}&=V_{ij}+\frac\kappa4m_i^2m_j^2,
 \qquad \Sigma=2C^\Phi.
 \label{aggregate:eq:covariance}
\end{align}
The MDE gives $-m_i^2=u_i(1-u_i|z_i|^2)$. At the exact clocks
$s_i=\operatorname{arctanh}u_i$ from \eqref{sp:fb:scale},
$\Sigma$ is therefore the restriction of
\eqref{gaussian:gp2:covariance}, including its fourth-cumulant term.
Its positivity for every $\kappa\ge-1$ follows from the feature
expansion \eqref{gaussian:gp2:features}. Let $Z$ be a centered Gaussian
vector with this covariance, allowing repeated coordinates and singular
covariance. This is the exact Gaussian comparison kernel, not an
assertion about the finite-$N$ covariance of the matrix field.

For scalar tests and vector fields put
\begin{align}
 \cD_r(F)&=\norm F_\infty+
 \sum_{a=1}^r\sup_x\sum_{j_1,\ldots,j_a}
     |\partial_{j_1}\cdots\partial_{j_a}F(x)|,\\
 \cB_r(f)&=\sum_{a=0}^r\sup_y\sum_{i,j_1,\ldots,j_a}
     |\partial_{j_1}\cdots\partial_{j_a}f_i(y)|.
 \label{aggregate:eq:norms}
\end{align}
For \(a=0\), the inner sum in the second line is \(\sum_i|f_i(y)|\).
The supremum is outside each sum. This differs from the sum of
coordinatewise supremum norms used for signed exponential coefficients. Every derivative in a contraction
below is evaluated at the same random vector, so the pathwise bound
\eqref{aggregate:eq:aggregatebound} uses precisely this weaker norm.
Taking separate suprema first would lose the dimension-independent
bound required for a growing list.

\begin{theorem}[Gaussian comparison for many observations]
\label{aggregate:thm:au-comparison}
For the model and list above, there are constants \(C,N_0\), depending
on the uniform moment bounds, \(K_+,c_0,\eps,L\), but not on the chosen list, such that
for \(N\ge N_0\) and every bounded \(C^8\) test with finite
\(\cD_8\)-norm,
\begin{equation}
 |\E F(X)-\E F(Z)|\le C N^{-\eps/8}\cD_8(F).
 \label{aggregate:eq:comparison}
\end{equation}
There is no uniformity claim as \(\eps\downarrow0\).
The exponent and derivative order are deliberately nonoptimal.
\end{theorem}

The maximum application uses $L=9$. This theorem does not claim
a comparison for arbitrary spatial points at a logarithmic cutoff. All fixed
entry moments needed below are chosen before $N$ tends to infinity.

\subsection{One common event and coordinates defined by finite integrals}
\label{aggregate:sec:event}

Now we define bounded integral coordinates and one common event on which all observations and required resolvent estimates can be compared.
The finite upper endpoint gives coordinates that are deterministically
bounded for every matrix. This makes the cumulant expansion and the
later Gaussian comparison by Stein's method well defined before localization. The actual
field is recovered at the end through a bounded scalar test.
Set \(U=N^{10}\), \(\eta_*=\min_i\eta_i\), and
\begin{equation}
 \widetilde Y_i=
 \ii\int_{\eta_i}^{U}N\av{G^{z_i}(\ii u)-M^{z_i}(\ii u)}\dd u-E_i.
 \label{aggregate:eq:finite}
\end{equation}
The variables are real and satisfy the deterministic bound
\(\max_i|\widetilde Y_i|\le CN\log N\).
They differ from \(Y_i\) by at most \(CN^{-5}\) on an event whose
complement has probability smaller than any prescribed inverse power
of the matrix size. One can use the much stronger
large-\(U\) bound directly: the paired eigenvalues of the Hermitization
cancel the linear term at infinity, and a polynomial bound on
\(\norm{A_N}\) makes the integrated tail negligible.
Only a high-probability comparison is needed because the eventual
scalar test is bounded.

The same event must support every coordinate pair, integration height
and physical-index contraction. Fix a local-law tolerance
\begin{equation}
 \xi_{\rm loc}=\eps/100.
 \label{aggregate:eq:xi}
\end{equation}
On a common event, for \(\eta_i\le u\le1\),
\begin{align}
 |(G^{z_i}(\ii u)-M^{z_i}(\ii u))_{ab}|
    &\le C\frac{N^{\xi_{\rm loc}}}{\sqrt{Nu}}=:\Delta(u),\\
 |\av{(G^{z_i}(\ii u)-M^{z_i}(\ii u))B}|
    &\le C\frac{N^{\xi_{\rm loc}}}{Nu},\\
 \norm{G^{z_i}(\ii u)}_{\mathrm{HS}}
    &\le C\sqrt{N/u}.
 \label{aggregate:eq:singlelaws}
\end{align}
The averaged bound is required for the bounded deterministic block
insertions actually used below. Anisotropic estimates are also imposed
for the normalized vectors constant on one of the two physical components.
Multiplication by a deterministic \(2\times2\) block matrix preserves
this finite class up to fixed linear combinations.

For every needed pair of coordinates and outer integration parameter,
Proposition~B.2(i) of \cite{CL}, with its tolerance chosen
as $\xi_{\rm loc}<\eps$, gives
\begin{equation}
 \left|\av{G^{z_i}(\ii u)B_1G^{z_j}(\ii\eta_j)B_2}
              -\av{\mathcal M_{ij}^{B_1}(u)B_2}\right|
 \le C\frac{N^{\xi_{\rm loc}}}{Nu\eta_j}.
 \label{aggregate:eq:twolaw}
\end{equation}
The notation \(\mathcal M_{ij}^{B_1}\) denotes the deterministic
solution of the two-resolvent matrix Dyson equation. The estimate applies
to arbitrary spatial pairs, including coincident ones. For \(u\ge1\),
Lemma~B.1 applies. The additional decay of the outer insertion makes its
integral finite.

\paragraph{Choosing a common event.}
The list exponent $L$ is fixed first. Failure orders and the finitely
many moment orders may depend on it; they do not grow with $N$.
For \(d_N\le N^9\), a spectral mesh of spacing \(N^{-40}\) on
\([\eta_*,N^{10}]\) has \(O(N^{50})\) points. The relevant one- and
two-resolvent assertions, their physical entries, and their fixed
block-vector tests involve fewer than \(N^{75}\) choices for large \(N\).
One may demand failure exponent \(400\) for the individual local-law
assertions, leaving a common failure probability smaller than
\(N^{-200}\). Deterministic resolvent derivatives and the explicit
finite-dimensional MDE formulas extend the estimates between mesh points.
At \(\eta_*\ge c_0N^{-1+\eps}\), their polynomial derivative bounds
leave a negligible interpolation error at this mesh size.

The deterministic constants in these bounds can depend on \(K_+\) and
\(c_0\). Finitely many enlarged meshes or fixed block insertions do
not change these probability estimates. For general $L$, enlarge the fixed failure exponent after counting
the coordinate pairs and mesh tests and bounding the exceptional
terms deterministically; the numerical choices
above are an example for the maximum application only. This changes fixed
moment requirements, never the local-law tolerance $\xi_{\rm loc}$.
A third coordinate in a chain-rule term does
not require a new three-resolvent local law: the third-order bounds
below use separate Ward estimates.

Also include
\[
 \max_{a,b}|(A_N)_{ab}|\le N^{-2/5}.
\]
For example, the moment \(\E|\xi|^{4096}<\infty\) makes its complement
smaller than \(CN^{-400}\) by Markov's inequality and a union bound.
The high-probability local laws themselves may require additional
fixed moments to realize the chosen failure exponent. These are supplied
by the unchanged all-moments hypothesis. The number 4096 is only a
convenient moment order for truncation, not a claimed optimal or sufficient
standalone moment assumption for all imported local laws.

\paragraph{Single-entry changes.}
On the common event, all relevant resolvent entries are bounded by a
fixed constant, since \(\Delta(\eta_*)=o(1)\). Replacing one matrix
entry by \(\theta\), \(|\theta|\le N^{-2/5}\), changes the
Hermitization by a rank-two matrix. In the Woodbury formula, the inverse
to be taken is a \(2\times2\) matrix of the form \(I+O(N^{-2/5})\).
Thus all modified entries remain bounded uniformly in the full entry
segment, with integrable upper-height decay. This argument preserves
these bounds; it does not claim that the modified matrix retains the
original sharp $N^{-1/2}$ fluctuation accuracy.
This proves the particular mesoscopic stability needed here;
it is also consistent with Proposition~9.1 of \cite{CL}.
No independence of the common event and the changed entry is assumed.

\subsubsection{Large-height control}\label{aggregate:sec:au-large-height}

Next we  control the large-height part of the resolvent integrals so that every cumulant estimate remains integrable.
Also impose a fixed operator-norm bound on $A_N$, whose failure
probability is smaller than any prescribed inverse power of the matrix
size by the rigidity input in
Theorem~2.5 of \cite{CL}. Then we choose a fixed rectangle surrounding the
spectra of all $H^{z_i}$ and the supports of their deterministic matrix
spectral measures, with horizontal sides at imaginary parts $\pm1$
and vertical sides a fixed positive distance outside those spectra.
On the horizontal sides the usual anisotropic local law applies.
On a vertical side, for $|\operatorname{Im}w|\ge N^{-1/2}$, the
fixed gap to the deterministic support gives
$\operatorname{Im}m(w)=O(|\operatorname{Im}w|)$; Theorem~2.3 then gives
an $O(N^{-1/2+\xi_{\rm loc}})$ entrywise and anisotropic error. For smaller imaginary parts, the deterministic gap to the
spectrum and resolvent differentiation extend this same bound.
A polynomial contour mesh makes these bounds simultaneous.

The Cauchy representation of the matrix spectral measures now gives,
for $u\ge2$, entrywise and for the required anisotropic tests,
\[
 |(G^{z_i}(\ii u)-M^{z_i}(\ii u))_{ab}|
 \le C N^{-1/2+\xi_{\rm loc}}/u.
\]
For $1\le u\le2$ the usual local law suffices. Multiplication by
$\cA_i(u)=O(u^{-1})$ gives a corresponding
$O(N^{-1/2+\xi_{\rm loc}}u^{-2})$ error for $P_i$ at large height,
with $\cA_i$ and $P_i$ defined in Section~\ref{aggregate:sec:raw}.
These estimates justify the integrable large-height versions of the
third- and fourth-order contractions below.

For the quadratic averaged remainder one needs the better averaged
power of $N$, not merely the entrywise estimate. Fix a bounded
normalized deterministic block insertion $B$. Apply Lemma~B.1 with
one height equal to one, and use
\[
 G^{z_i}(\ii u)=G^{z_i}(\ii)
       -\ii(1-u)G^{z_i}(\ii)G^{z_i}(\ii u).
\]
For $u>1$, the deterministic two-resolvent solution for insertion $I$ is
$(M^{z_i}(\ii)-M^{z_i}(\ii u))/(\ii(1-u))$, by the two MDEs. At $u=1$ use the usual averaged law.
Subtracting its trace with $B$ gives
\[
 |\av{(G^{z_i}(\ii u)-M^{z_i}(\ii u))B}|
 \le C N^{\xi_{\rm loc}}/N,\qquad 1\le u\le U.
\]
After extracting $\norm{\cA_i(u)}\le C/u$, the quadratic error
integrates to at most $CN^{-1+2\xi_{\rm loc}}\log U$.
These averaged estimates therefore justify the contribution from the
full upper integration interval.

\subsection{The resolvent identity summed over observations}
\label{aggregate:sec:raw}

Next we sum the resolvent identities before estimating test derivatives, keeping the growing observation list under one common bound.
For the \(i\)th coordinate, we write
\[
 \cA_i(u)=\frac{M^{z_i}(\ii u)}{1-\av{M^{z_i}(\ii u)^2}},
 \qquad P_i(u)=G^{z_i}(\ii u)\cA_i(u).
\]
In the bulk,
\[
 \norm{\cA_i(u)}\le\frac C{1+u},\qquad
 \norm{P_i(u)}\le Cu^{-2}\quad(u\ge1).
\]
Equation~(8.15) of \cite{CL} is the algebraic identity
\begin{equation}
 \av{G_i-M_i}=-\av{(WG_i+\av{G_i}G_i)\cA_i}
       +\av{G_i-M_i}\av{(G_i-M_i)\cA_i},
 \label{aggregate:eq:rawres}
\end{equation}
where \(W=H^0\). It is used before extracting any mean. In particular,
we do not replace a product with a random multiplier by the product
of its expectations.

Let \(\sum_{ab}^{\mathrm{or}}\) sum both oriented bipartite blocks
of the Hermitization. Multiplying by \(f_i(\widetilde Y)\), integrating,
and summing \(i\), the raw term is
\begin{equation}
 -\frac\ii2\sum_{ab}^{\mathrm{or}}
 \E\left[W_{ab}\mathcal H_{ab}(W)\right],
 \qquad
 \mathcal H_{ab}(W)=\sum_i\int_{\eta_i}^{U}
       (P_i(u))_{ba}f_i(\widetilde Y(W))\dd u.
 \label{aggregate:eq:summedfunction}
\end{equation}
The renormalization and the deterministic center are retained separately
until the second and fourth cumulants are evaluated.

We use the multivariate cumulant algebra of
\cite[Proposition~3.2 and Appendix~C]{EKS}. The aggregate remainder
bounds below are proved for the whole summed scalar function, so that
the coordinate derivatives retain their common evaluation point.

For each independent
complex entry, apply its two-real-variable cumulant formula to the
single scalar function \(\mathcal H_{ab}\), not separately to its
individual coordinate summands. For \(w=W_{ab}\), the convention is
\[
 \E[w\mathcal H]=
 \sum_{r+s=1}^{q}
 \frac{\cum(w^{[r+1]},\overline w^{[s]})}{r!s!}
 \E[\partial_w^r\partial_{\overline w}^{s}\mathcal H]
 +\mathcal R_q.
\]
Here the sum has $r,s\ge0$, and repeated arguments in a joint cumulant are indicated by brackets;
the two orientations are combined in \eqref{aggregate:eq:summedfunction}.
This is a finite Taylor/cumulant identity with an integral remainder.
The centered first cumulant is zero, and the variance contribution is
\(N^{-1}\partial_{\overline w}\mathcal H\).

For an oriented entry derivative \(D_{ab}\),
\begin{equation}
 D_{ab}G=-GE^{ab}G,\qquad
 D_{ab}\widetilde Y_j
 =\frac12\left(G^{z_j}(\ii\eta_j)-G^{z_j}(\ii U)\right)_{ba}.
 \label{aggregate:eq:endpoint}
\end{equation}
Higher entry derivatives of \(\widetilde Y_j\) are finite sums of
products of endpoint resolvent entries. Their constants depend only
on derivative order. The finite upper-endpoint corrections are
negligible in every retained term.

For each fixed chain-rule type, after contracting physical indices,
the error has the form
\begin{equation}
 \E\sum_{i,j_1,\ldots,j_a}
 Q_{i;j_1\cdots j_a}(W)
 \partial_{j_1}\cdots\partial_{j_a}f_i(\widetilde Y(W)).
 \label{aggregate:eq:commonvector}
\end{equation}
The coefficient estimates in the next sections are pathwise on the
common event, uniformly in every coordinate tuple. Consequently,
\begin{equation}
 \left|\sum_{i,\boldsymbol j}Q_{i;\boldsymbol j}(W)
         \partial_{\boldsymbol j}f_i(\widetilde Y(W))\right|
 \le \max_{i,\boldsymbol j}|Q_{i;\boldsymbol j}(W)|
       \sup_y\sum_{i,\boldsymbol j}|\partial_{\boldsymbol j}f_i(y)|.
 \label{aggregate:eq:aggregatebound}
\end{equation}
The coefficient is bounded uniformly before the coordinate sums are
taken. Thus the bound has no factor $d_N$. Scalar derivatives may
remain random because they share one evaluation point; no independence
between them and the resolvent coefficient is used.

The quadratic term in \eqref{aggregate:eq:rawres}, using the averaged laws, has
absolute integrated coefficient at most
\begin{equation}
 C N\int_{\eta_i}^1\frac{N^{2\xi_{\rm loc}}}{N^2u^2}\dd u
       +\text{integrable tail}
 \le \frac{CN^{2\xi_{\rm loc}}}{N\eta_*}+CN^{-1+2\xi_{\rm loc}}(1+\log N).
 \label{aggregate:eq:quad}
\end{equation}
This is multiplied by the zeroth aggregate derivative norm.

\subsection{Second cumulant and the exact Gaussian covariance}
\label{aggregate:sec:second}

This subsection computes the second-cumulant contribution exactly, identifying the logarithmic covariance that the later fourth-order calculation will complete.
The derivative landing only on \(P_i\) gives
\[
 D_{ba}(P_i)_{ba}=-(G_i)_{bb}(P_i)_{aa}.
\]
For a square Hermitization, the traces of the two diagonal blocks of
\(G_i\) are exactly equal. Summing the two oriented blocks therefore
cancels the \(\av{G_i}G_i\) term in \eqref{aggregate:eq:rawres} exactly.
The assumption \(\E\xi^2=0\) excludes a second-cumulant
same-orientation term. Thus there is no order-one complex Gaussian
mean contribution.

The derivative landing on the multiplier has coefficient
\begin{equation}
 \cV_{ij}(W)=
 -\frac{\ii}{4N}\int_{\eta_i}^{U}
 \sum_{ab}^{\mathrm{or}}(P_i(u))_{ba}
        \left(G^{z_j}(\ii\eta_j)-G^{z_j}(\ii U)\right)_{ab}\dd u.
 \label{aggregate:eq:randomcov}
\end{equation}
It multiplies \(\partial_j f_i(\widetilde Y)\).
The trace is a two-resolvent trace with deterministic block insertions.
The small-height error from \eqref{aggregate:eq:twolaw} is bounded by
\[
 C\int_{\eta_i}^1\frac{N^{\xi_{\rm loc}}}{Nu\eta_j}\dd u
 \le \frac{CN^{\xi_{\rm loc}}\log N}{N\eta_*}.
\]
For \(u\ge1\), factor out \(\norm{\cA_i(u)}\le C/u\) in
Lemma~B.1. This gives an integrable bound
\(CN^{\xi_{\rm loc}}/(N\eta_j u^2)\). The upper endpoint in
\eqref{aggregate:eq:randomcov} contributes a negligible negative power of \(N\).

\paragraph{Deterministic evaluation.}
For clarity, the denominator in the covariance is the determinant of
an explicit \(2\times2\) linear system. If \(M_i,M_j\) denote the
\(2\times2\) MDE blocks, the diagonal entries of the deterministic
two-resolvent solution solve a system with matrix
\begin{equation}
 B_{ij}=\begin{pmatrix}
 1-u_iu_jz_i\overline z_j&-m_im_j\\
 -m_im_j&1-u_iu_j\overline z_i z_j
 \end{pmatrix},\qquad \det B_{ij}=D_{ij}.
 \label{aggregate:eq:stability}
\end{equation}
Indeed the equation is
\(\mathcal M=M_i B M_j+M_i\operatorname{diag}
(\mathcal M_{22},\mathcal M_{11})M_j\), with insertion \(B\).
Let $\Pi_1=E_{11}$ and $\Pi_2=E_{22}$ denote the two block
projections, represented here by their $2\times2$ matrices.
At variable outer height $u$, evaluate the first MDE pair at $\ii u$
and write $D_{ij}(u,\eta_j)$ for the resulting determinant; thus
$D_{ij}=D_{ij}(\eta_i,\eta_j)$.
Taking \(B=\cA_i\Pi_k\), summing the opposite diagonal component,
and differentiating the MDE gives
\begin{equation}
 \sum_{k\ne l}\av{\mathcal M_{ij}^{\cA_i\Pi_k}(u)\Pi_l}
       =\frac\ii2\partial_u\log D_{ij}(u,\eta_j).
 \label{aggregate:eq:primitive}
\end{equation}
The sum has exactly $(k,l)=(1,2),(2,1)$, and the normalized trace
here is one half of the ordinary $2\times2$ trace.
One verifies \eqref{aggregate:eq:primitive} by solving \eqref{aggregate:eq:stability};
the two scalar differentiation identities are
\[
 (M_i\cA_i)_{11}=(M_i\cA_i)_{22}=m_i'(\ii u),
 \qquad
 \partial_w u^z(w)=\frac{w\partial_w m^z(w)-m^z(w)}{(w+m^z(w))^2}.
\]
More explicitly, write $m=m_i(w)$, $n=m_j(i\eta_j)$,
$q=u_i(w)u_jz_i\overline z_j$ and
$\widehat q=u_i(w)u_j\overline z_i z_j$, with primes denoting
$w$-derivatives. In this analytic calculation, $D$ denotes the same
determinant before setting $w=\ii u$. The two opposite diagonal entries
obtained from insertions $\cA_i\Pi_1$ and $\cA_i\Pi_2$ are respectively
\[
 \frac{(1-q)\widehat q'+mm'n^2}{D},\qquad
 \frac{(1-\widehat q)q'+mm'n^2}{D}.
\]
Their half-sum, accounting for normalized trace, is
$-\partial_w\log D/2=(i/2)\partial_u\log D$ at $w=iu$.
Here $\widehat q$ is an analytic expression, not conjugation of the
variable being differentiated.
The same primitive underlies the evaluation (8.20)--(8.21) in the cited paper,
but the two spatial arguments and the two lower cutoffs are kept distinct.
Since \(D_{ij}(u,\eta_j)\to1\) as \(u\to\infty\), integration
of \eqref{aggregate:eq:primitive} with the factor \(-\ii/2\) in
\eqref{aggregate:eq:randomcov} gives \(-\frac14\log D_{ij}\).
Therefore, on the common event,
\begin{equation}
 \max_{i,j}|\cV_{ij}(W)-V_{ij}|
 \le \frac{CN^{\xi_{\rm loc}}\log N}{N\eta_*}+CN^{-8}.
 \label{aggregate:eq:secondbound}
\end{equation}
The coefficient bound is obtained before multiplication by any random
derivative. Together with \eqref{aggregate:eq:quad} and
\eqref{aggregate:eq:aggregatebound}, it supplies the dimension-independent
second-order terms needed in the remaining cumulant calculation.

\subsection{Third cumulants: all four classes}
\label{aggregate:sec:third}

now we are able to bound every third-cumulant term and show that none changes the limiting drift or covariance.
In the schematic cumulant formulas of this and the next subsection,
$Y$ denotes the finite-integral vector $\widetilde Y$, with the already
controlled upper-endpoint corrections suppressed.

The remaining expansion has three tasks. Every third-cumulant term
must be an error; the fourth cumulants must supply exactly the written
drift and covariance correction; the higher orders must retain the
aggregate derivative norm. We keep the physical-entry count separate
from the coordinate sums throughout.

Group the Hermitization indices into \(N\) physical sites, each with
two components. Between distinct physical sites every deterministic
MDE entry is zero. Call a factor joining different sites off-diagonal;
it costs \(\Delta(u)\) at outer height \(u\), or
\(\Delta_j=CN^{\xi_{\rm loc}}/\sqrt{N\eta_j}\) at endpoint \(j\).
All other good-event entries are bounded. The coordinate derivatives
remain at the common vector in \eqref{aggregate:eq:commonvector}, so
each uniform coefficient bound below is inserted into
\eqref{aggregate:eq:aggregatebound} only after the physical contraction.

After trace normalization and the initial factor \(N\), the
third-cumulant prefactor is \(O(N^{-3/2})\) before the physical pair
sum. Ignoring only fixed permutation multiplicities, the four classes are
\begin{equation}
 (D^2P_i)f_i,\quad
 (DP_i)(\partial_jf_i)DY_j,\quad
 P_i(\partial_jf_i)D^2Y_j,\quad
 P_i(\partial_{j\ell}f_i)DY_jDY_\ell.
 \label{aggregate:eq:thirdtypes}
\end{equation}
For distinct physical sites, the first three classes have either one
or three off-diagonal factors; the last has three. This follows by
inserting \(D G=-GEG\) and \eqref{aggregate:eq:endpoint}. It does not depend
on equality of the spatial arguments.

\paragraph{One off-diagonal factor.}
If both diagonal factors have their MDE values, their coefficients
are constant in the physical-site labels. The remaining off-diagonal
resolvent, or \(G\cA_i\), is contracted against deterministic vectors
constant on one component. Its total sum is \(O(N)\): the normalized
anisotropic contraction is bounded, and subtracting the matching-site
terms also costs \(O(N)\). The prefactor makes this \(O(N^{-1/2})\).
The scalar test derivative is common to the physical sum and need not
be independent of it.

If a diagonal factor is replaced by its local-law error, use its
entrywise supremum and an absolute Ward bound for the remaining sum.
For any endpoint,
\[
 \sum_{a,b}|(G_j)_{ab}|\le CN^{3/2}/\sqrt{\eta_j}.
\]
For an outer error and endpoint off-diagonal factor this gives
\begin{equation}
 N^{-3/2}\int_{\eta_i}^1\Delta(u)
       \sum_{a,b}|(G_j)_{ba}|\dd u
 \le \frac{CN^{\xi_{\rm loc}}}{\sqrt{N\eta_j}}.
 \label{aggregate:eq:thirdmixed}
\end{equation}
An endpoint diagonal error with an outer off-diagonal factor instead
gives \(C\Delta_j\int_0^1u^{-1/2}\dd u\).
If all factors have height \(u\), the bound is
\(CN^{-1/2+\xi_{\rm loc}}\int_{\eta_i}^1u^{-1}\dd u\).
More than one diagonal error causes no worse bound because the other
good-event entries, including errors, are bounded. No random vector
is inserted into an anisotropic estimate.

\paragraph{Three off-diagonal factors.}
One entrywise maximum and two Hilbert--Schmidt bounds suffice, even
when all spatial points differ. The following table gives the small-height
integrated coefficient for each class; fixed constants are suppressed.
\begin{center}
\small
\begin{tabular}{p{0.37\textwidth}p{0.52\textwidth}}
\toprule
Class in \eqref{aggregate:eq:thirdtypes} & Three-off-diagonal bound\\
\midrule
\((D^2P_i)f_i\) &
\(N^{-1+\xi_{\rm loc}}\int_{\eta_i}^1u^{-3/2}\dd u
 \le CN^{-1+\xi_{\rm loc}}/\sqrt{\eta_*}\)\\[6pt]
\((DP_i)(\partial_jf_i)DY_j\) &
\(CN^{-1+\xi_{\rm loc}}\log N/\sqrt{\eta_j}\)\\[6pt]
\(P_i(\partial_jf_i)D^2Y_j\) &
\(CN^{-1+\xi_{\rm loc}}/\eta_j\)\\[6pt]
\(P_i(\partial_{j\ell}f_i)DY_jDY_\ell\) &
\(CN^{-1+\xi_{\rm loc}}/\sqrt{\eta_j\eta_\ell}\)\\
\bottomrule
\end{tabular}
\end{center}
For example, the last line follows from
\[
 N^{-3/2}\int_{\eta_i}^1
 \Delta(u)\norm{G_j}_{\mathrm{HS}}\norm{G_\ell}_{\mathrm{HS}}\dd u
 \le\frac{CN^{\xi_{\rm loc}}}{N\sqrt{\eta_j\eta_\ell}}.
\]
Transposes, conjugates, and component partners have the same estimates.

\paragraph{Large heights and matching physical sites.}
For \(u\ge1\), use the probabilistic bounds proved in
Section~\ref{aggregate:sec:au-large-height}, not just $\|G(iu)\|\le u^{-1}$.
For $(DP_i)DY_j$, an outer diagonal replacement error is
$O(N^{-1/2+\xi_{\rm loc}}u^{-3})$ and the absolute endpoint sum is
$O(N^{3/2}/\sqrt{\eta_j})$, giving
$CN^{\xi_{\rm loc}}/\sqrt{N\eta_j}$ after the cumulant factor and integration.
For $P_iD^2Y_j$, an endpoint diagonal error is $\Delta_j$ and
$\sum_{ab}|(P_i)_{ab}|\le CN^{3/2}/u^2$, giving $C\Delta_j$.
For $D^2P_i$ the random outer replacement carries $N^{-1/2+\xi_{\rm loc}}$
and at least three inverse powers of $u$. Three-off-diagonal terms
use the same maximum and Hilbert--Schmidt choices as above, now with
integrable outer decay. All-deterministic diagonal branches use the
flat-vector contraction $O(N/u)$ or $O(N/u^2)$ and the other outer
factors. Thus the integral over large heights satisfies the required
error bound.

Matching physical sites number only \(O(N)\), and bounded entries
give an \(O(N^{-1/2})\) contribution there. Combining the classes,
every third-cumulant coefficient is bounded by
\begin{equation}
 C N^{2\xi_{\rm loc}}\left[
 (N\eta_*)^{-1/2}+\frac{\log N}{N\eta_*}
             +N^{-1/2}\log N\right].
 \label{aggregate:eq:thirdbound}
\end{equation}
In particular, no order-one third-cumulant contribution survives.

\subsection{Fourth cumulants: exact drift and mixed covariance}
\label{aggregate:sec:fourth}

Now we identify the fourth-cumulant terms that determine the drift and the smooth correction to the Gaussian covariance.
There are two surviving coefficients: a deterministic drift and a
mixed covariance term. We compute both from the same local polynomial
before estimating the remaining terms. The seven classes are
\begin{align*}
 &(D^3P_i)f_i;\\
 &(D^2P_i)(\partial_jf_i)DY_j;\\
 &(DP_i)(\partial_jf_i)D^2Y_j;\\
 &(DP_i)(\partial_{j\ell}f_i)DY_jDY_\ell;\\
 &P_i(\partial_jf_i)D^3Y_j;\\
 &P_i(\partial_{j\ell}f_i)(D^2Y_j)DY_\ell;\\
 &P_i(\partial_{j\ell m}f_i)DY_jDY_\ell DY_m.
\end{align*}
The physical pair count times the fourth-cumulant prefactor is \(O(1)\).
At distinct sites, only the first and third classes can contain no
off-diagonal factor. Only the balanced complex cumulant contributes
to those classes.

Here is an explicit calculation of both coefficients. Select an upper
index \(a\) and a lower index \(b\) at different physical sites, and
write \(x=W_{ab}\), \(y=W_{ba}\). Put \(m=m_i(\ii u)\) and
\(p=(M_i\cA_i)_{aa}=(M_i\cA_i)_{bb}\). After replacing every resolvent
factor by its deterministic MDE approximation, the relevant two-site
expansions are
\begin{equation}
 S(x,y)=-\frac{y m p}{1-xy m^2},\qquad
 \delta Y_j(x,y)=\frac12\log(1-xy m_j^2).
 \label{aggregate:eq:twosite}
\end{equation}
They follow from the rank-two inverse and determinant formulas.
The coefficient at derivative order three is therefore
\begin{equation}
 \partial_x\partial_y^2
 \{S(x,y)f_i(Y+\delta Y(x,y))\}\big|_{0}
 =-2m^3p f_i(Y)+mp\sum_jm_j^2\partial_jf_i(Y).
 \label{aggregate:eq:fourthpolynomial}
\end{equation}
All other third-derivative orientations vanish in this deterministic
substitution. Terms with two coordinate derivatives of \(f_i\) cannot
appear at this order, because each \(\delta Y_j\) begins at degree two.
The same expansion also shows that every third-cumulant term with
all resolvents replaced by their MDE approximations has zero coefficient
at entry derivative order two.

The balanced cumulant is
\(\cum(w,w,\overline w,\overline w)=\kappa/N^2\).
Its \(1/(1!2!)\) coefficient, the two orientations, and the raw factor
\(-\ii/2\) yield respectively
\begin{align}
 \ii\kappa\int_{\eta_i}^{\infty}m_i(\ii u)^3
             m_i'(\ii u)\dd u
       &=-\frac\kappa4m_i^4=E_i,\\
 -\frac{\ii\kappa}{2}m_j^2\int_{\eta_i}^{\infty}
             m_i(\ii u)m_i'(\ii u)\dd u
       &=\frac\kappa4m_i^2m_j^2.
 \label{aggregate:eq:fourthprimitives}
\end{align}
The first term cancels the deterministic center in
\eqref{aggregate:eq:finite}; it is not added a second time by a separate
mean extraction. The second is exactly the rank-one correction in
\eqref{aggregate:eq:covariance}. Omitting matching sites changes these
calculations by \(O(N^{-1})\), which is within the error.

Every other term has an off-diagonal factor. Replacing a diagonal
factor in either leading term by its error also gives one local-law
factor. The integrated error is bounded by
\[
 C N^{\xi_{\rm loc}}(N\eta_*)^{-1/2}+CN^{-1}.
\]
The outer integration is harmless, using \(\int_0^1u^{-1/2}\dd u<\infty\)
and the large-height decay. This covers the unbalanced complex fourth
cumulants as well.

Together with \eqref{aggregate:eq:primitive}, these two coefficients
give exactly $C^\Phi$ and cancel the stated deterministic center.
The calculation permits distinct spatial points and unequal cutoffs.

\subsection{Higher orders and the common remainder point}
\label{aggregate:sec:remainder}

In this subsection we bound higher cumulants and Taylor remainders while keeping every test derivative at the same field vector.
The preceding terms determine the covariance; the remaining orders
are errors. For cumulant derivative order \(k\ge4\), the crude physical count is
\[
 N^2N^{-(k+1)/2}=N^{(3-k)/2}.
\]
All good-event endpoint derivatives are bounded at any fixed order.
Every term has an integrated outer factor, whose tail at large height
is integrable. Thus the retained orders \(k=4,5,6\) are bounded by
\(CN^{-1/2}\), times the appropriate aggregate derivative norm.

\paragraph{Retaining the summed derivative norm in the remainder.}
Use \(q=6\) in the scalar expansion of the summed function
\(\mathcal H_{ab}\) in \eqref{aggregate:eq:summedfunction}.
At a Taylor remainder point, all appearances of the field are
\(\widetilde Y(A_N^{(ab,\theta)})\) for one and the same modified
matrix. For each fixed \(\theta\), the chain rule is summed over
\(i,j_1,\ldots,j_a\) before taking any supremum. Outer integration
changes coefficients, not the evaluation point of a test derivative.
On the common event, this gives
\begin{equation}
 \sup_{|\theta|\le N^{-2/5}}
 |D^7\mathcal H_{ab}(A_N^{(ab,\theta)})|
       \le C\cB_7(f),
 \label{aggregate:eq:remainderderivative}
\end{equation}
for every real/complex orientation combination of the seven derivatives.
The estimate follows by bounding the outer products uniformly on
\([\eta_*,1]\), which has length at most one, and by their integrable
large-height bound on \([1,U]\). Different lower endpoints are
implemented by their indicators inside that integral.

For general matrices, the same derivative is bounded by
\(CN^7(1+\log N)\cB_7(f)\). This follows from \(\norm G\le1/u\),
the endpoint formula, and the fact that a term with \(a\) derivatives
on the outer factor integrates to at most \(CN^a(1+\log N)\), while
the other \(7-a\) derivatives contribute at most \(CN^{7-a}\).
A larger fixed polynomial bound would also suffice.

\paragraph{An event determined by the matrix with one entry set to zero.}
To use the small eighth entry moment in the remainder, choose the
segment event as a function of the deleted background alone. For the matrix $V$ with entry $(a,b)$ deleted, let
$\mathcal E_{ab}(V)$ assert that all required resolvent entries remain bounded
for every replacement $|\theta|\le N^{-2/5}$, with $O(u^{-1})$
resolvent decay for $u\ge1$. Only these bounds, and the deterministic
$O(u^{-1})$ decay of $\cA_i(u)$, are part of this event; the sharper
probabilistic bound on $G-M$ is not required for modified matrices. Complete $V$ with an independent iid entry $h'$.
On the common event for this completion and $|h'|\le N^{-2/5}$,
the rank-two Woodbury inverse is $I+O(N^{-2/5})$ uniformly throughout
the replacement disk. Therefore
$\mathbb P_V(\mathcal E_{ab}^{\,c})\le N^{-D}$, after the moment
orders have been chosen. The completion is used only to bound the failure probability of an
event defined by $V$. That event is independent of the actual entry
appearing in the cumulant remainder. Continuity permits a countable rational definition and ensures
measurability. A union bound over physical entries is harmless.
For each replacement matrix the full field vector in every derivative
is the same. In particular, the sum in \eqref{aggregate:eq:remainderderivative}
is taken before the supremum in $\theta$ and before any independent
coordinate mean-value points could be introduced.

The usual integral cumulant remainder involves the eighth entry
moment and finitely many products of lower moments of the same total
order. To see their scale explicitly, condition on the deleted
background and Taylor expand the single function $\mathcal H_{ab}$
and its cumulant-side derivatives about entry zero. If $S_7$ bounds
its seventh directional derivatives along the entry segment, the
remainder is bounded by
\[
 C\E[|w|^8S_7]
 +C\sum_{k=1}^6N^{-(k+1)/2}\E[|w|^{7-k}S_7].
\]
Every moment product has scale $N^{-4}$ per physical entry, since
$N^{-(k+1)/2}\E|w|^{7-k}=O(N^{-4})$. Products of lower moments
are also bounded by a constant times the eighth moment by H\"older.
On this event and for small entry values,
\begin{equation}
 \sum_{a,b}|\mathcal R_{6,ab}|
 \le C N^2\E|N^{-1/2}\xi|^8\cB_7(f)
 =CN^{-2}\cB_7(f).
 \label{aggregate:eq:remaindergood}
\end{equation}
Apply the same split to any independent entry copy in the remainder
formula. The derivative bound is always for the total scalar function
at one replacement vector. Separate coordinate mean-value points
would not give \eqref{aggregate:eq:remainderderivative}.

\paragraph{Bad events and large replacement values.}
The scalar tests are unweighted, so there is no exponential envelope.
On a bad local-law event use the deterministic polynomial bound and
Cauchy--Schwarz with the fixed entry moments. For $L=9$ the displayed event has
failure probability at most \(N^{-200}\); for general $L$ choose the
larger fixed failure exponent of Section~\ref{aggregate:sec:event}, after
counting the tests and bounding the exceptional terms deterministically. These are
choices depending on fixed $L$, not restrictions to the $N^9$ application.
Because the event is determined by the matrix with one entry set to zero,
one can alternatively factor its failure probability from the actual
entry moments. For \(|w|>N^{-2/5}\), use
\[
 \E\bigl[|w|^8\ind_{\{|w|>N^{-2/5}\}}\bigr]
 \le N^{-(p/10)-16/5}\E|\xi|^p,
 \qquad w=N^{-1/2}\xi,\quad p\ge8.
\]
For example \(p=4096\), even after the physical count and the crude
derivative polynomial, makes the contribution negligible.
All retained low-order bad-event terms are controlled by the same
fixed choice of failure exponent. No derivatives of an indicator of the common event are taken;
localization is used only after the exact cumulant expansion.

\begin{proposition}[Summed Stein identity]
\label{aggregate:prop:stein}
For bounded smooth vector fields with finite \(\cB_7\),
\begin{equation}
 \left|\E\left[
 \sum_i\widetilde Y_if_i(\widetilde Y)
 -\sum_{i,j}C^\Phi_{ij}\partial_jf_i(\widetilde Y)
 \right]\right|
 \le C a_N\cB_7(f),
 \label{aggregate:eq:stein}
\end{equation}
where one may take
\begin{equation}
 a_N=N^{2\xi_{\rm loc}}\left[
 (N\eta_*)^{-1/2}
 +\frac{\log N}{N\eta_*}
 +N^{-1/2}(1+\log N)\right]+CN^{-2}+CN^{-8}.
 \label{aggregate:eq:rate}
\end{equation}
After increasing \(N_0\), \(a_N\le C N^{-\eps/8}\).
\end{proposition}

\begin{proof}
Combine \eqref{aggregate:eq:quad}, \eqref{aggregate:eq:secondbound},
\eqref{aggregate:eq:thirdbound}, \eqref{aggregate:eq:fourthprimitives}, the fourth-order
errors, and \eqref{aggregate:eq:remaindergood}. Use
\eqref{aggregate:eq:aggregatebound} once for each of the finitely many chain-rule
types. The finite upper-endpoint corrections are absorbed by
\(CN^{-8}\); all bad-event terms are smaller than this with the
specified failure exponent. Since \(N\eta_*\ge c_0N^\eps\) and
\(\xi_{\rm loc}=\eps/100\), the displayed rate is bounded by the deliberately
weaker power \(N^{-\eps/8}\), absorbing fixed logarithms.
\end{proof}

We make the choices in the following order: fix \(\eps\) and the list exponent
$L$, choose \(q=6\), choose \(\xi_{\rm loc}=\eps/100\), impose the common-event
failure bounds and the needed fixed moments, and only then take
\(N\to\infty\). The required moments may depend on $L$ and the requested failure
order, but are all fixed before the matrix limit. These estimates prove
\eqref{aggregate:eq:stein} with the full $\cB_7$ norm and
rate \eqref{aggregate:eq:rate}; the Gaussian comparison follows next.

\subsection{Ornstein--Uhlenbeck comparison without covariance inversion}
\label{aggregate:sec:OU}

We now turn the summed Stein identity into Gaussian comparison using a semigroup argument that also allows singular covariance.
First we take $F$ smooth. A bounded $C^8$ test is recovered by convolution
with a nonnegative smooth compactly supported kernel of unit mass.
Convolution does not increase $\cD_8$, because every derivative is
evaluated at the same shift; bounded first derivatives give uniform
convergence as the support shrinks. Thus the final estimate passes to
the stated $C^8$ class, also for singular covariance.

The Ornstein--Uhlenbeck semigroup argument follows
\cite[Lemmas~1(3) and~2(1), and (10)]{Meckes}. We keep the calculation
in the common-point derivative norm to make its independence of the
list length explicit.

We now convert the summed identity once, using the aggregate norm at
the common Gaussian argument. Rescaling
Proposition~\ref{aggregate:prop:stein} gives the same identity for
\(\widetilde X=\sqrt2\widetilde Y\) with covariance \(\Sigma\).
The fixed factors from rescaling derivatives are absorbed in the constant;
they do not depend on the list length.
For \(Z\sim N(0,\Sigma)\), let
\[
 P_sF(x)=\E F(e^{-s}x+\sqrt{1-e^{-2s}}Z),\qquad
 v_T(x)=-\int_0^T(P_sF(x)-\E F(Z))\dd s.
\]
For every \(a\ge1\), differentiation of the Gaussian convolution,
followed by the triangle inequality with the common Gaussian argument,
gives
\begin{equation}
 \sup_x\sum_{j_1,\ldots,j_a}|\partial_{j_1}\cdots\partial_{j_a}v_T(x)|
 \le\frac1a\sup_x\sum_{j_1,\ldots,j_a}
             |\partial_{j_1}\cdots\partial_{j_a}F(x)|.
 \label{aggregate:eq:OUnorm}
\end{equation}
The generator identity is
\[
 (\Sigma:\nabla^2-x\cdot\nabla)v_T=F-P_TF.
\]
Apply the rescaled summed identity to \(f_i=\partial_i v_T\).
The vector derivative sums up to order seven become the scalar
derivative sums of orders one through eight. Bound
\eqref{aggregate:eq:OUnorm} is uniform in $T$ and uses no matrix norm
of $\Sigma$. Hence
\[
 |\E F(\widetilde X)-\E P_TF(\widetilde X)|
       \le C N^{-\eps/8}\cD_8(F).
\]
Send \(T\to\infty\) by bounded convergence. The linear term in the
Stein identity is integrable because \(\widetilde X\) is bounded
by \(CN\log N\). No inverse of \(\Sigma\) is used, and the argument
works even when coordinates repeat and \(\Sigma\) is singular.

Finally, the high-probability comparison between finite-integral and
actual fields is transferred by the deterministic inequality
\[
 |F(x)-F(y)|\le\norm{x-y}_\infty
       \sup_v\sum_i|\partial_iF(v)|.
\]
On the exceptional event use \(2\norm F_\infty\). The maximum coordinate
error multiplies the aggregate first-derivative bound only once. This proves
Theorem~\ref{aggregate:thm:au-comparison} for every fixed polynomial
list exponent $L$.

\subsection{Smooth maxima and the deterministic spatial mesh}
\label{aggregate:sec:maximum}

This subsection uses smooth maxima on a deterministic mesh to control the barrier throughout the disk during an initial fraction of scales.
We use the standard smooth-maximum approximation of
\cite[(1) and Lemmas~3--4]{CCK}. The categorical calculation below
records the derivative bounds through order eight required by
Theorem~\ref{aggregate:thm:au-comparison}.

We first transfer a maximum over any admissible finite list, then
choose a deterministic list for the spatial application. For deterministic
shifts \(b_i\), define
\[
 S_\beta(x)=\beta^{-1}\log\sum_{i=1}^d e^{\beta(x_i-b_i)},
 \qquad M(x)=\max_i(x_i-b_i).
\]
Then \(M\le S_\beta\le M+\beta^{-1}\log d\). If
\(p_i=e^{\beta(x_i-b_i)}/\sum_j e^{\beta(x_j-b_j)}\), the
\(a\)-th derivative tensor of \(S_\beta\) is
\(\beta^{a-1}\) times the order-\(a\) cumulant tensor of a categorical
vector with probabilities \((p_i)\). The partition formula gives
\begin{equation}
 \sup_x\sum_{i_1,\ldots,i_a}|\partial_{i_1}\cdots\partial_{i_a}S_\beta(x)|
 \le c_a\beta^{a-1},\qquad
 c_a=\sum_{k=1}^a S(a,k)(k-1)!,
 \label{aggregate:eq:softderivs}
\end{equation}
where \(S(a,k)\) is a Stirling number of the second kind.
For each partition, the sum of the absolute product of the corresponding
categorical moment tensors is one. Thus \(c_a\) depends on order only,
not on \(d\).

Let \(0\le\theta\le1\) be a fixed smooth function, zero on
\(( -\infty,0]\) and one on \([1,\infty)\). The chain rule and
\eqref{aggregate:eq:softderivs} yield
\[
 \cD_8\{\theta(S_\beta-B+1)\}\le C(1+\beta)^8.
\]
Take \(\beta=\max(1,\log d)\). The test equals one whenever
\(M>B\), and it can be nonzero only when \(M>B-2\).
Theorem~\ref{aggregate:thm:au-comparison} gives the finite-list estimate
\begin{equation}
 \Pp\{M(X)>B\}
 \le\Pp\{M(Z)>B-2\}
       +C N^{-\eps/8}(1+\log N)^8.
 \label{aggregate:eq:maxtransfer}
\end{equation}
The fixed height slack is allowed in barrier removal. There is no need
for Gaussian anti-concentration at a moving dimension.

\paragraph{Interpolation in the field scale and spatial point.}
For this application use the fixed entry law of Section~\ref{sp:common},
so $\kappa$ and the Gaussian family are fixed as $N\to\infty$.
Let \(K\) lie in the interior of \(K_+\), put
\(t_N=(\log N)/2\), \(J_N=\lfloor(1-\eps)t_N\rfloor\), and
\begin{align*}
 u_k&=\tanh k,\qquad
 \eta_k(z)=(1-u_k)\sqrt{u_k^{-1}-|z|^2},\\
 X_{N,k}(z)&=\sqrt2[\Phi_N(z,\eta_k(z))-E(z,\eta_k(z))].
\end{align*}
For fixed \(k_0\ge1\), uniformly on \(K_+\),
\[
 \eta_k(z)\ge c_{K_+}N^{-1+\eps},\qquad
 |\nabla_z\eta_k(z)|\le C_{K_+}\eta_k(z),
 \qquad k_0\le k\le J_N.
\]
The raw regularized log determinant has spatial and height derivatives
bounded deterministically by \(CN/\eta\). The MDE centering composed with the clock has spatial derivative
\(O(N)\) on a compact bulk set. Indeed, for $u=\tanh k$, its exact
potential is
\[
 \int\log(x^2+\eta_k(z)^2)\rho^z(x)\dd x
 =-u-\log u+|z|^2u^2,
\]
whose spatial derivative is bounded. The cumulant correction is smooth
and bounded. Therefore
\begin{equation}
 \sup_{z\in K_+}\max_{k_0\le k\le J_N}|\nabla_zX_{N,k}(z)|\le CN^2.
 \label{aggregate:eq:meshLip}
\end{equation}
Use a deterministic spatial mesh of spacing \(N^{-4}\) in a fixed
neighborhood of \(K\) contained in \(K_+\). It contains \(O(N^8)\)
points. Including the \(O(\log N)\) integer scales gives \(d_N\le N^9\)
for large \(N\). This is the $L=9$ application of the general comparison
theorem. The continuum maximum is at most the grid maximum plus
\(CN^{-2}\), deterministically.

Apply \eqref{aggregate:eq:maxtransfer} to the grid list with \(b_i=2k_i\).
By the clock identification in Section~\ref{aggregate:sec:theorem},
its Gaussian vector is precisely \((G_{k_i}(z_i))_i\) from the specified
Gaussian family. These probabilities are under the unweighted matrix law. For each fixed \(\eps\),
\begin{align}
 &\limsup_{N\to\infty}
 \Pp\left\{\sup_{z\in K}\max_{k_0\le k\le J_N}
               (X_{N,k}(z)-2k)>B\right\}\\
 &\hspace{18mm}\le
 \Pp\left\{\sup_{z\in K_+}\sup_{k\ge k_0}
               (G_k(z)-2k)>B-3\right\}.
 \label{aggregate:eq:gaussianreduction}
\end{align}
The extra unit absorbs the vanishing deterministic interpolation error.

\begin{corollary}[Uniform barrier bound on the initial scales]
\label{aggregate:cor:au-maximum}
Use the exact Gaussian maximum estimate proved in
Proposition {\hyperref[gaussian:gp2:max-prop]{\ref*{gaussian:gp2:max-prop}}} in Section~\ref{gaussian:gp2:maximum}:
\[
 \Pp\left\{\sup_{z\in K_+}\sup_{k\ge k_0}(G_k(z)-2k)>b\right\}
 \le C_{K_+}(1+b)e^{-2b},\qquad b\ge b_0.
\]
Then, for every fixed \(0<\eps<1\) and all sufficiently large \(B\),
\begin{equation}
 \limsup_{N\to\infty}
 \Pp\left\{\sup_{z\in K}\max_{k_0\le k\le J_N}
             (X_{N,k}(z)-2k)>B\right\}
 \le C_K(1+B)e^{-2B}=:p(B).
 \label{aggregate:eq:barrier}
\end{equation}
The function \(p(B)\) is independent of \(\eps\).
\end{corollary}

\begin{proof}
Combine \eqref{aggregate:eq:gaussianreduction} with the retained Gaussian bound
and absorb the fixed height shift into its constant. The matrix
comparison constant may depend on \(\eps\), but its error disappears
in the \(N\)-limit for each fixed \(\eps\). The surviving Gaussian
bound has no \(\eps\) dependence.
\end{proof}

This is the global barrier estimate through the stated initial fraction of
scales, under the unweighted matrix law, used in
Proposition~\ref{p34:global}. The Gaussian estimate is
Proposition~\ref{gaussian:gp2:max-prop}, which combines the sharp
fixed-scale maximum bound with the summable
\((1+B+\log k)k^{-3/2}e^{-2B}\) tail. The finite-$N$ comparison
retains its fixed-$\eps$ constants: only after taking $N\to\infty$
is the surviving bound independent of $\eps$. The final $\eps t_N$
scales remain outside this theorem, as does all-spatial comparison
at a merely logarithmic cutoff.

\endgroup

\begingroup
\newcommand{\E}{\mathbb E}
\newcommand{\PP}{\mathbb P}
\newcommand{\ind}{\mathbf1}
\newcommand{\dd}{\,\mathrm d}
\newcommand{\TV}{\mathrm{TV}}
\newcommand{\Var}{\operatorname{Var}}
\section{Loss from the final scales and barrier removal}\label{chap:barrier}

The purpose of this section is to calculate the mass lost on the final scales and use it to remove the history barrier in probability.
Indeed, we estimate the mass removed when the barrier is also
imposed on the final fraction of scales, as required by
Proposition~\ref{p34:loss}. We compare the Gaussian weighted laws while
retaining the barrier on a shorter initial sequence of integer scales.
We then calculate separately the expectations with this initial
restriction and with the full barrier. The resulting nonnegative loss supplies
Corollary~\ref{p34:joint-corollary}, whose proof gives barrier
removal in probability and joint convergence in distribution.

\subsection{Conventions and inputs}
\label{barrier:sec:setup}

Now we collect the fixed-barrier limits and the barrier bound through an initial fraction of scales needed to estimate the remaining mass loss.
We use the model, fixed disks, fields and exact Gaussian family of
Section~\ref{sp:common}. In particular, $t=\frac12\log N$,
$n=\lfloor t\rfloor$, the scale field is
\eqref{sp:fb:scale}, and the raw and corrected endpoints and marks
are those of \eqref{sp:fb:endpoints}--\eqref{sp:fb:roots}.
The deterministic correction cancels from the normalized density
$W_N$. At determinant zeros the weight vanishes and both marks
are set to zero, as in the common convention.

For an integer $j\ge k_0$, retain the barriers
\eqref{sp:barriers} and write
\[
 Z_N^{B,j}=\ind_{A_{N,j}^B}Z_N,\qquad Z_N^B=Z_N^{B,n}.
\]
The measures $\Xi_N$ and $\Xi_N^B$ use the raw mark $R_N$.
For the corrected mark, it is useful also to write
\begin{equation}
 \Xi_N^0(\dd z,\dd r)=Z_N(\dd z)\delta_{R_N^0(z)}(\dd r),
 \qquad
 \Xi_N^{B,0}(\dd z,\dd r)=Z_N^B(\dd z)\delta_{R_N^0(z)}(\dd r).
 \label{barrier:eq:corrected-marks}
\end{equation}
The second measure is exactly the one in \eqref{sp:fb:marked}.
All spatial measures are restricted to $K$, so their total mass
is tested by the continuous function one.

The covariance is $C_\kappa$ from \eqref{sp:fb:covariance}.
For a real root time $u>0$, the Gaussian tilted law and the
event that the barrier holds through a specified integer time are
\[
 \frac{\dd Q_{u,z}^G}{\dd\mathbb P_G}
   =e^{2G_u(z)-2\operatorname{Var}G_u(z)},\qquad
 A_j^B(z)=\bigcap_{k=k_0}^j\{G_k(z)\le2k+B\}.
\]
We keep the notation $\mu_s^{G,B}$ of
\eqref{sp:gaussian-measures}: its terminal weight is at real time
$s$, while its barrier ends at the integer $\lfloor s\rfloor$.

The following inputs have already been established independently
of the removal deduction.
\begin{enumerate}
\item Proposition~\ref{p1:proposition} gives
$\mu_s^{G,B}(f)\to Z_B(f)$ in $L^2$ for real $s\to\infty$,
where
\[
 Z_B=\ind_{A_\infty^B}M_\kappa|_K,\qquad
 A_\infty^B(z)=\bigcap_{k=k_0}^{\infty}\{G_k(z)\le2k+B\}.
\]
In particular, $\mathbb E_G Z_B(K)<\infty$ for each fixed $B$.
The maximum estimate of Proposition~\ref{gaussian:gp2:max-prop}
gives, on the fixed enlarged bulk disk $K_+$,
\begin{equation}
 \mathbb P_G\left\{\sup_{z\in K_+}\sup_{k\ge k_0}
                    (G_k(z)-2k)>B\right\}
 \le C_{K_+}(1+B)e^{-2B}
 \label{barrier:eq:br-gaussian-max}
\end{equation}
for sufficiently large $B$. Hence $Z_B\uparrow M_\kappa|_K$,
with eventual equality on each sample as $B\to\infty$.
\item The retained-history comparison in
Section~\ref{oneroot:retained-history-comparison-with-fixed-coarse-observations}
applies to smooth tests of the retained vector through a logarithmic
cutoff, with the stated derivative bounds and arbitrary prescribed
inverse-logarithmic accuracy. Its argument for removing the Gaussian
smoothing
also applies to a shorter, possibly $N$-dependent integer prefix.
The quantitative use needed here is given below; it is distinct
from a fixed early-prefix insertion in the full-history theorem.
\item Corollary~\ref{sp:fb:identification} gives
\begin{equation}
 (Z_N^B,\Xi_N^{B,0},C_N)
 \Longrightarrow (Z_B,Z_B\otimes\rho_{\mathrm R},C_G),
 \label{barrier:eq:br-fixed-limit}
\end{equation}
and the same fixed-barrier limit with the raw-marked $\Xi_N^B$.
Here the coarse vector is any admissible list in
\eqref{sp:coarse-vectors}, including height one, repeated observations
and prescribed sites in $K_*$ outside $K$. The Rayleigh law is
\[
 \rho_{\mathrm R}(\dd r)=r e^{-r^2/2}\ind_{\{r>0\}}\dd r.
\]
This convergence includes total mass and preserves the specified
Gaussian coupling. It does not assert an uncut matrix second moment.
\item Proposition~\ref{p34:global}, verified by
Corollary~\ref{aggregate:cor:au-maximum}, gives, for every fixed
$\epsilon\in(0,1)$ and $J=J_\epsilon=\lfloor(1-\epsilon)t\rfloor$,
\begin{equation}
 \limsup_N\mathbb P(E_N(B,\epsilon)^c)\le p(B),\qquad
 E_N(B,\epsilon)=\{A_{N,J}^B(z)\text{ for every }z\in K\},
 \label{barrier:eq:br-prefix-global}
\end{equation}
where $p(B)\le C_K(1+B)e^{-2B}$ is independent of $\epsilon$.
The finite-$N$ comparison rate may depend on fixed $\epsilon$.
This is an event under the unweighted matrix law.
\end{enumerate}
The spatial indicator defining $Z_B$ is not a global event
requiring all points to survive.

\subsection{Changing the Gaussian weighting scale while retaining the initial barriers}
\label{barrier:sec:root}

Next we shall change the Gaussian weighting scale without changing the initial barrier event, making its survival probability comparable at two terminal times.
We fix $B\in\mathbb R$, $\epsilon\in(0,1)$ and $A>2$, and put
\[
 J=\lfloor(1-\epsilon)t\rfloor,\qquad
 L=J-\lceil A\log(2+t)\rceil.
\]
All limits first send $N\to\infty$ with these parameters fixed.
Eventually $L\ge k_0$ and $L\asymp_\epsilon t$.

\begin{lemma}[Comparison of Gaussian weighting scales with the same initial barriers]
\label{barrier:lem:br-root}
Uniformly for $z\in K$,
\[
 \sqrt t\,\left|Q_{t,z}^G(A_J^B(z))
                  -Q_{J,z}^G(A_J^B(z))\right|\longrightarrow0.
\]
\end{lemma}
\begin{proof}
Under $Q_{u,z}^G$, the history mean is $2C_\kappa(k,z;u,z)$
and its centered covariance is unchanged. The exact same-site
kernel in \eqref{sp:fb:covariance} gives
\[
 |\partial_u C_\kappa(k,z;u,z)|\le Ce^{-2(u-k)},
 \qquad u\ge J,\quad k_0\le k\le L.
\]
Indeed, the denominator is comparable to $e^{-2k}+e^{-2u}$
on the fixed bulk disk, while differentiation in the deeper
root contributes $O(e^{-2u})$. The bounded cumulant term
satisfies the same bound. Integration gives
\begin{equation}
 \sup_{\substack{k_0\le k\le L\\z\in K}}
 |C_\kappa(k,z;t,z)-C_\kappa(k,z;J,z)|
 \le Ce^{-2(J-L)}.
 \label{barrier:eq:root-mean}
\end{equation}
Couple the two tilted arrays by the same centered Gaussian history.
Their scalar variances have a uniform positive lower bound,
including at $\kappa=-1$. More precisely,
$D(u_k,z;u_k,z)\le1-u_k^2$ and $p_z(u_k)^2\le u_k^2$ imply
\[
 \operatorname{Var}G_k(z)
 \ge\tfrac12[-\log(1-u_k^2)-u_k^2]\ge u_k^4/4,
 \qquad u_k=\tanh k,
\]
as in \eqref{gaussian:gpm:positive-variance}.
A union bound over scalar strips of width $Ce^{-2(J-L)}$
therefore gives
\begin{equation}
 \sup_{z\in K}|Q_{t,z}^G(A_L^B(z))-Q_{J,z}^G(A_L^B(z))|
 \le Ct e^{-2(J-L)}\le Ct^{1-2A}.
 \label{barrier:eq:early-root}
\end{equation}
This step uses neither a covariance inverse nor conditional
nondegeneracy of the history.

For either root $u\in\{J,t\}$ and $k\le J$, the mean of
$G_k-2k$ is uniformly bounded. For $L\le k\le J$, the centered
increment from $L$ has variance at most $C(k-L+1)$.
Fix $D>1/2$ and choose $R=C_{A,D}\log(2+t)$, with a sufficiently
large constant. Gaussian tails and a union bound over the
$O_A(\log t)$ coordinates yield
\begin{equation}
 \sup_{z\in K}Q_{u,z}^G\left\{
 \max_{L\le k\le J}[(G_k-2k)-(G_L-2L)]>R\right\}
 \le t^{-D}.
 \label{barrier:eq:increment-tail}
\end{equation}
No independence of increments is used. Outside this event,
survival through $L$ followed by a crossing by $J$ forces
$0\le B+2L-G_L\le R$.

For large $N$, $1\le R\le L$. Apply the one-root estimate
\eqref{gaussian:gpm:one-strip} in
Lemma~\ref{gaussian:gpm:strips} with observation horizon $L$
and external root $u=J$ or $u=t$. The external-clock convention
in Section~\ref{gaussian:gpm:endpoint} permits both choices,
uniformly in the gaps beyond $L$: the rooted means remain bounded
and the centered one-path covariance is unchanged. Thus
\begin{align}
 \sup_{z\in K}Q_{u,z}^G(A_L^B(z)\setminus A_J^B(z))
 &\le
 \frac{C_{B,D}R(1+\sqrt{\log(2+L)})
                    (1+R+\sqrt{\log(2+L)})}{L^{3/2}}
       +C_{B,D}L^{-D}+t^{-D}\nonumber\\
 &\le C_{B,\epsilon,A,D}
       \left((\log(2+t))^{5/2}t^{-3/2}+t^{-D}\right).
 \label{barrier:eq:block-loss}
\end{align}
The last line uses $L\asymp_\epsilon t$, including in the
$L^{-D}$ remainder. The barrier level remains $B$ throughout.

Since $A_J^B\subseteq A_L^B$, the triangle inequality and
\eqref{barrier:eq:early-root}--\eqref{barrier:eq:block-loss} give
\begin{align}
 &\sqrt t\sup_{z\in K}
    |Q_{t,z}^G(A_J^B(z))-Q_{J,z}^G(A_J^B(z))|\nonumber\\
 &\qquad\le C_{B,\epsilon,A,D}\left[
 t^{3/2-2A}+\frac{(\log(2+t))^{5/2}}{t}+t^{1/2-D}
 \right]\longrightarrow0.
 \label{barrier:eq:root-error}
\end{align}
Only the logarithmic final block has been discarded and restored.
The argument does not assume a Markov property for the exact
Gaussian scale family.
\end{proof}

\subsection{Expected mass lost on the final scales}
\label{barrier:sec:loss}

This subsection computes the limiting expected loss on the final scales and show that it vanishes as their fraction shrinks at fixed barrier height.
We  fix $B\in\mathbb R$ and $\epsilon\in(0,1)$, and keep
$J=J_\epsilon=\lfloor(1-\epsilon)t\rfloor$.
We first compare the growing prefix under the original determinant
root and the Gaussian root at real time $t$:
\begin{equation}
 \sup_{z\in K}\left|
 Q_{N,z}^{\mathrm{det}}(A_{N,J}^B(z))
       -Q_{t,z}^G(A_J^B(z))\right|=o(t^{-1/2}).
 \label{barrier:eq:rooted-prefix}
\end{equation}
This is an absolute error before critical normalization.
Because $J$ depends on $N$, we verify the quantitative invocation
of the retained-history comparison.

\paragraph{The retained vector and its smooth tests.}
In Section~\ref{oneroot:retained-history-comparison-with-fixed-coarse-observations},
denote the logarithmic history cutoff by
\[
 J_0=\left\lfloor t-(b/2+a)\log t-C_K\right\rfloor.
\]
All exponents will be fixed independently of $N$; hence $J\le J_0$
eventually at fixed $\epsilon$. Take the fixed insertion to be one.
For a smoothing width $h=t^{-s}$, upper and lower smooth barrier
products through $J$ are functions $F_N^\pm$ of the retained vector
through $J_0$ that ignore the later coordinates. With the
sum-of-suprema derivative convention
\[
 \mathcal D_r(F)=\|F\|_\infty+
  \sum_{q=1}^r\sum_{i_1,\ldots,i_q}
           \|\partial_{i_1}\cdots\partial_{i_q}F\|_\infty,
\]
they satisfy
\begin{equation}
 \|F_N^\pm\|_\infty\le1,\qquad
 \mathcal D_r(F_N^\pm)\le C_r t^{r(1+s)}.
 \label{barrier:eq:prefix-derivatives}
\end{equation}
Unused coordinates contribute zero derivatives.

Fix $P>1/2$, choose $P_*>P+3$, and then $s>P_*+2$.
Use the complete parameter order in
Section~\ref{oneroot:complete-parameter-order}, including the
local-factor, cap, moment and failure requirements. Its fixed
Stein derivative order $m$ precedes the choice of $b$; in particular,
the choices include
\begin{align*}
 b&>\max\{2[m(1+s)+P_*+4],\ m(1+s)+P_*+19\},\\
 a&>\max\{5,s+P_*+3\},
\end{align*}
with $a$ chosen after $b$. Thus the comparison error
\[
 C_b[t^{2-b/2}+t^{17-b}+N^{-c}]\mathcal D_m(F_N^\pm)
\]
is $o(t^{-P_*})$. The history-transport error times
$\mathcal D_1(F_N^\pm)$ is bounded by a constant times
\[
 [t^{1-a}+t^{1-2a}+t^{2-b-4a}+N^{-c}]t^{1+s}
     =o(t^{-P_*}).
\]
The remaining reconstruction, relative normalization, replacement
and errors from restoring the Gaussian parameters satisfy the required
bounds with the same parameter choices. These steps return to the actual
original determinant
density $W_N$, not merely to a normalized surrogate law.
Consequently the smooth products are compared at error $o(t^{-P})$,
uniformly over $z\in K$. The later choice of $a$ does not
change the already fixed reconstruction exponent $b$.

\paragraph{Removing the smoothing.}
Under $Q_{t,z}^G$, scalar variances are unchanged and comparable
to $k$ for integers $k\ge k_0$, uniformly on $K$, including
$\kappa=-1$. This follows from
\eqref{gaussian:gpm:positive-variance} and $v_k(z)=k+O(1)$.
The upper and lower products differ only on coordinate strips
of width $Ch$. Their total Gaussian probability is at most
\begin{equation}
 Ch\sum_{k=k_0}^{J}k^{-1/2}
       \le Ct^{1/2-s}=o(t^{-P}).
 \label{barrier:eq:prefix-desmoothing}
\end{equation}
The bound is uniform in the strip centers and rooted means.
The two smooth sandwiches therefore prove
\eqref{barrier:eq:rooted-prefix}; in fact they give $o(t^{-P})$
with the selected parameters. No matrix anti-concentration
estimate or fixed-prefix restriction is used.

\paragraph{Prefix and full-barrier expectations.}
By the exact rooted-density identity \eqref{sp:weighted-laws},
Tonelli's theorem, \eqref{barrier:eq:rooted-prefix} and
Lemma~\ref{barrier:lem:br-root},
\begin{align}
 \mathbb E Z_N^{B,J}(K)
 &=\sqrt t\int_K Q_{N,z}^{\mathrm{det}}(A_{N,J}^B(z))\dd A(z)
       \nonumber\\
 &=\sqrt t\int_K Q_{t,z}^G(A_J^B(z))\dd A(z)+o(1)
       \nonumber\\
 &=\sqrt t\int_K Q_{J,z}^G(A_J^B(z))\dd A(z)+o(1)
       \nonumber\\
 &=\sqrt{t/J}\,\mathbb E_G\mu_J^{G,B}(K)+o(1)
       \nonumber\\
 &\longrightarrow (1-\epsilon)^{-1/2}\mathbb E_G Z_B(K).
 \label{barrier:eq:prefix-first}
\end{align}
Both comparison errors are uniform on the fixed disk $K$.
The last identity is exact at the integer root time $J$.
Its factor reflects the normalizations $\sqrt t$ of the matrix
mass and $\sqrt J$ of the Gaussian prefix mass. The limit uses
Proposition~\ref{p1:proposition} with the constant test and
$t/J\to(1-\epsilon)^{-1}$.

Independently, Proposition~\ref{p2:matrix}\textup{(b),(e)},
with constant spatial and coarse tests, gives
\begin{equation}
 \mathbb E Z_N^B(K)
 =\mathbb E\widehat Z_N^B(K)
       +\mathbb E\|Z_N^B-\widehat Z_N^B\|_{\mathrm{TV}}
 \longrightarrow\mathbb E_G Z_B(K).
 \label{barrier:eq:full-first}
\end{equation}
For large $N$, $J\le n$ and
$A_{N,n}^B(z)\subseteq A_{N,J}^B(z)$, so the loss is
nonnegative pathwise. Subtracting \eqref{barrier:eq:full-first}
from \eqref{barrier:eq:prefix-first} proves
\begin{equation}
 \mathbb E[Z_N^{B,J}(K)-Z_N^B(K)]
 \longrightarrow
 \big[(1-\epsilon)^{-1/2}-1\big]\mathbb E_G Z_B(K).
 \label{barrier:eq:br-loss}
\end{equation}
This verifies Proposition~\ref{p34:loss} before the barrier is
removed in probability. Both expectations are barrier restricted.
The normalization without the barrier instead gives
$\mathbb E Z_N(K)=\sqrt t\,|K|$; no uniform integrability of
that mass or expectation convergence from weak convergence
has been assumed.

\subsection{Barrier removal in probability and the joint marked limit}
\label{barrier:sec:removal}

Now we combine the bound on early barrier failures with the expected loss on later scales to remove the barrier in the joint marked limit.
The calculation \eqref{barrier:eq:br-loss} completes the
verification of Proposition~\ref{p34:loss}. Together with the
global barrier estimate through the initial fraction of scales under the
unweighted matrix law in Proposition~\ref{p34:global}, it supplies the two
inputs to Corollary~\ref{p34:joint-corollary}. That corollary removes the
barrier in probability and gives the joint spatial, raw-marked and coarse
limit without the barrier. Its proof gives the complete deduction. In particular, \eqref{p34:ordered} uses only
$\mathbb E_G Z_B(K)<\infty$ at fixed $B$, with the order
\eqref{p34:order}: first $N\to\infty$, then
$\epsilon\downarrow0$ at fixed $B$, then $B\to\infty$.

The corrected-mark conclusion remains available as well. For
both mark conventions restriction discards exactly the spatial mass:
\begin{equation}
 \|\Xi_N^0-\Xi_N^{B,0}\|_{\mathrm{TV}}
 =\|\Xi_N-\Xi_N^B\|_{\mathrm{TV}}
 =\|Z_N-Z_N^B\|_{\mathrm{TV}}.
 \label{barrier:eq:both-marked-masses}
\end{equation}
Use \eqref{barrier:eq:br-fixed-limit} in the same approximation
argument of Corollary~\ref{p34:joint-corollary}. Its coarse
coordinate is unchanged, and its Gaussian approximation is still
$(Z_B,Z_B\otimes\rho_{\mathrm R},C_G)$ on the specified space.
Consequently, for every admissible fixed coarse list,
\begin{equation}
 (Z_N,\Xi_N^0,C_N)
 \Longrightarrow
 (M_\kappa|_K,M_\kappa|_K\otimes\rho_{\mathrm R},C_G)
 \label{barrier:eq:corrected-joint-limit}
\end{equation}
in the same finite-measure product topology as
\eqref{p34:joint-limit}. This also retains observations at
height one, repeated observations and permitted sites outside $K$.
The raw-mark passage at a fixed barrier is already included in
Corollary~\ref{sp:fb:identification}, so no second raw-mark
correction is needed here.

Appending original integral field tests is the separate transfer
argument of Proposition~\ref{p5:transfer} and Corollary~\ref{p5:field-corollary}.
The compact-window deduction additionally uses
Proposition~\ref{cw:inputs}. Neither conclusion is an input to
the expected-loss calculation or to barrier removal.

\endgroup

\begingroup
\newcommand{\E}{\mathbb E}
\newcommand{\PP}{\mathbb P}
\newcommand{\C}{\mathbb C}
\newcommand{\R}{\mathbb R}
\newcommand{\D}{\mathbb D}
\newcommand{\Ytest}{\mathcal Y}
\newcommand{\dd}{\,\mathrm d}
\newcommand{\br}[1]{\langle #1\rangle}
\newcommand{\one}{\mathbf1}
\newcommand{\Tr}{\operatorname{Tr}}
\section{Joint limits with averages of the original field}\label{chap:field}

In this section we recover averages of the original logarithmic field from regularized observations while preserving their joint dependence on the critical measure.

The marginal Gaussian limit of the original integrated logarithm follows
from \cite[Theorem 2.2 and Corollary 2.4]{CES}: spatial integration
expresses it as an eigenvalue linear statistic with a smooth
logarithmic-potential test, and an exterior cutoff is removed using
spectral confinement and the entry moment bounds, including in the
centering. This section verifies the joint transfer from fixed coarse
observations in Proposition~\ref{p5:transfer}, preserving the specified
limiting coupling with the additional random element left unchanged
in the approximation.

The first step approximates the spatially integrated logarithm by a field
at a positive cutoff, under the unweighted matrix law. The interval
reaching height zero needs separate treatment: for discrete entry laws,
a pointwise logarithm can have infinite negative expectation, although
its spatial integral is integrable. We therefore integrate first,
remove the microscopic interval, and only then use the resolvent
covariance estimates.

The final transfer leaves an arbitrary additional random element
unchanged. Its joint convergence with coarse observations is a
hypothesis, supplied separately in the application to critical measures.

\subsection{Definitions and resolvent estimates}
\label{field:sec:inputs}

We first define the centered spatial integrals and record the resolvent bounds that will control the removal of regularization.
We use the matrix $A_N$ and logarithms $L_{N,\eta}$ from
Section~\ref{sp:common}. Thus $A_N=N^{-1/2}(\xi_{ab})$ has complex
iid entries satisfying
\[
 \E\xi=\E\xi^2=0,\qquad \E|\xi|^2=1,\qquad
 \E|\xi|^p\le C_p\quad(p\ge1).
\]
For the final model the entry law is fixed. The estimates below are
uniform for families with one common collection of moment bounds.
Put $\kappa=\E|\xi|^4-2\ge-1$ and let $\dd A=\dd x\dd y$ denote
unnormalized planar area. For a prescribed finite family of real tests in $C_c^\infty(\D)$,
choose one closed bulk disk $K_*\Subset\D$ whose interior contains
all their supports and the measure disk $K$. The test supports need
not lie in $K$. Write $\varphi$ for one test in this family. For a fixed
positive continuous profile $a$ on $K_*$, or for $a\equiv0$, set
\begin{equation}
 \Ytest_{N,a}(\varphi)=\sqrt2\left\{
  \int_{K_*}\varphi(z)L_{N,a(z)}(z)\dd A(z)
  -\E\int_{K_*}\varphi(z)L_{N,a(z)}(z)\dd A(z)\right\}.
 \label{field:eq:of-definition}
\end{equation}
The spatial integral is taken before centering. In the common notation,
\[
 \Ytest_{N,0}(\varphi)=\Ytest_N(\varphi),
\]
where the right side is \eqref{sp:of:original}; the deterministic
potential in $Y_N$ cancels exactly. The symbol $\Ytest_{N,a}$ denotes
an integrated field test, whereas $W_N$ in \eqref{sp:fb:roots}
denotes the expectation-normalized critical density.

Write
\[
 H^z=\begin{pmatrix}0&A_N-zI\\(A_N-zI)^*&0\end{pmatrix},\quad
 G^z(i\eta)=(H^z-i\eta)^{-1},\quad
 \br A=(2N)^{-1}\Tr A,
\]
and $F_z(\eta)=\br{G^z(i\eta)-\E G^z(i\eta)}$.
The trace $F_z$ is purely imaginary. We use the following results of
Cipolloni--Erd\H{o}s--Schr\"oder \cite{CES},
under Assumption~2.1 of that paper.

First, its Lemma 4.4 gives an $L^1$ error tending to zero for the
centered integrated resolvent contribution over $[0,\eta_c]$ against
$\Delta f$, for every $f\in C_c^\infty(\C)$. Its proof splits this
interval at $\eta_0$, where
\[
 \eta_0=N^{-1-\delta_0},\qquad \eta_c=N^{-1+\delta_1}
\]
use the sufficiently small positive exponents provided by that proof;
in particular, choose $\delta_1<1/2$. We use the complete interval
$[0,\eta_c]$, not just its portion above $\eta_0$. The open set containing the support of $f$
may be any bounded open set containing the closed unit disk.
No density assumption on the entries is required: the density
condition in its Remarks 4.2 and 4.5 belongs only to an alternative
proof.

Second, Theorem 3.1 and the expectation estimate of Proposition 3.3
give, at any prescribed small fixed slack $\xi_{\rm loc}>0$,
\begin{equation}
 |F_z(\eta)|\le\frac{C N^{\xi_{\rm loc}}}{N\eta}
 \label{field:eq:of-local-law}
\end{equation}
with failure probability smaller than any prescribed inverse power of
the matrix size in the ranges used below.
Finally, its Proposition 3.3, equations (3.12)--(3.14), gives on a
fixed bulk compact
\begin{equation}
 \E F_z(\eta)F_w(\theta)
 =\frac{V(z,w;\eta,\theta)+\kappa U(z,\eta)U(w,\theta)}{2N^2}
     +O(\Psi),
 \label{field:eq:of-source-covariance}
\end{equation}
provided $\eta_*:=\min(\eta,\theta)\ge
N^{\xi_{\rm loc}-1}|z-w|^{-2}$, where
\begin{equation}
 \Psi\le
 \frac{C_{K_*}N^\epsilon}{(N\eta_*)^{1/2}}
 \frac{|z-w|^{-4}}{N^2\eta\theta}.
 \label{field:eq:of-source-error}
\end{equation}
Here $\epsilon>0$ is the fixed slack in the covariance error; it is
chosen together with $\xi_{\rm loc}$ in the proof of
Theorem~\ref{field:thm:of-approximation}.
The factors measuring distance from the unit circle in the cited
estimate are absorbed into $C_{K_*}$. In that notation,
\[
 V=\tfrac12\partial_\eta\partial_\theta\log D,
 \qquad U(z,\eta)=\frac{i}{\sqrt2}\partial_\eta m^z(i\eta)^2.
\]
Here $m$ is the scalar Hermitized MDE solution and $D$ is the mixed
deterministic stability denominator. Its explicit formula and the
bulk bounds in \cite[(4.21)]{CES} imply, for a fixed small $a_0>0$,
\begin{equation}
 |V(z,w;\eta,\theta)|
   \le\frac{C_{K_*}}{(|z-w|^2+\eta+\theta)^2},\qquad
 |U(z,\eta)|\le C_{K_*},\qquad 0\le\eta,\theta\le a_0.
 \label{field:eq:of-majorant}
\end{equation}
These are external inputs. The variable cutoffs and arbitrary spatial
weights below are deductions, not literal statements of the cited
resolvent theorem.

\subsection{Integrated logarithms and the microscopic interval}
\label{field:sec:microscopic}

This subsection shows that spatial integration makes the logarithmic singularities integrable, allowing us to remove the microscopic interval even for discrete entry laws.

\begin{lemma}[Spatial logarithmic integrability]
\label{field:lem:of-integrability}
For every fixed bounded measurable $Q\subset\C$, every matrix $A_N$,
and every measurable profile $a:Q\to[0,a_{\max}]$ with fixed
$a_{\max}<\infty$,
\begin{equation}
 \int_Q|L_{N,a(z)}(z)|\dd A(z)
 \le C_{Q,a_{\max}}N\{1+\log(2+\|A_N\|_{\rm F})\}.
 \label{field:eq:of-log-bound}
\end{equation}
In particular, a constant $C_Q$ suffices for all constant cutoffs
$0\le\eta\le1$.
Consequently \eqref{field:eq:of-definition} is well defined, also at zero,
for the stated matrix laws.
\end{lemma}
\begin{proof}
At zero, factor the characteristic polynomial into its $N$ eigenvalue
factors. The finitely many zeros have planar area zero. Uniformly in $\lambda\in\C$,
\[
 \int_Q\log_-(|z-\lambda|)\dd A(z)
 \le\int_{|u|\le1}-\log|u|\dd A(u)=\pi/2.
\]
The positive part is bounded by
$C_Q\log(2+|\lambda|)$, and $|\lambda|\le\|A_N\|_{\rm F}$.
Pointwise, $L_{N,a(z)}\ge L_{N,0}$, so the negative part
decreases. The positive part is bounded using
\[
 \|(A_N-zI)(A_N-zI)^*+a(z)^2I\|
 \le(\|A_N\|+|z|)^2+a_{\max}^2.
\]
This proves the display, including constant and zero cutoffs.
The entry moment assumptions give the required expectations and
all fixed polynomial moments of the right side, with at most
polynomial growth in $N$ under common moment bounds.
\end{proof}

\begin{remark}
A discrete entry law can give $\E L_{N,0}(z)=-\infty$ at a particular
deterministic $z$. Lemma~\ref{field:lem:of-integrability} does not assert
otherwise. It justifies centering the spatial integral. Moreover,
$L_{N,\eta_c}-L_{N,0}\ge0$ and, away from the finite eigenvalue set,
\[
 L_{N,\eta_c}(z)-L_{N,0}(z)
   =N\int_0^{\eta_c}\Im\br{G^z(i\eta)}\dd\eta.
\]
Its integral against $|\varphi|$ has finite expectation by the lemma.
Tonelli therefore applies to the nonnegative uncentered increment;
subtracting the expectation of its spatial integral gives the centered
identity. No pointwise expectation at height zero is needed.
\end{remark}

\begin{lemma}[Microscopic removal for arbitrary tests]
\label{field:lem:of-microscopic}
For the admissible cutoff $\eta_c$ above,
\begin{equation}
 \E|\Ytest_{N,\eta_c}(\varphi)-\Ytest_{N,0}(\varphi)|\longrightarrow0.
 \label{field:eq:of-micro}
\end{equation}
There is no condition that $\int\varphi\dd A=0$.
\end{lemma}
\begin{proof}
Define the smooth Newtonian potential
\[
 F=(2\pi)^{-1}\log|\cdot|*\varphi,\qquad \Delta F=\varphi.
\]
Choose $R>3$ larger by one than a fixed operator-norm bound that
holds with overwhelming probability. Let $\chi_R$ be a smooth cutoff equal to one on a disk
containing $K_*$, with derivative support in a fixed annulus outside
radius $R$, and put $f=\chi_RF$. Then $f\in C_c^\infty$ and
\[
 \Delta f=\varphi+h,\qquad
 h=2\nabla\chi_R\cdot\nabla F+F\Delta\chi_R,
\]
where $h$ is supported in the exterior annulus. This construction
allows $\int\varphi\dd A\ne0$: the extra term $h$ accounts for the
noncompact Newtonian potential. On
$\|A_N\|\le R-1$, all singular values of $A_N-z$ on this support
are at least one, whence
\[
 0\le L_{N,\eta_c}(z)-L_{N,0}(z)
 \le\tfrac12\sum_{i=1}^N\frac{\eta_c^2}{s_i(A_N-z)^2}
 \le N\eta_c^2/2=o(1).
\]
This fixed operator-norm bound holds with overwhelming probability,
for example, by
taking $z=0$ in \cite[Theorem 2.5]{CL} and using its Lemma 2.1(ii).
To check the expectation on the exceptional event, the integrated
logarithmic bound gives a dominating variable
\[
 H_N=C_{R,\varphi}N[1+\log(2+\|A_N\|_{\rm F})],
 \qquad \E H_N^2\le C_{R,\varphi}N^3.
\]
Here $\E\|A_N\|_{\rm F}^2=N$. For any fixed $D>3$, choose the
norm estimate with failure probability $O(N^{-D})$; Cauchy--Schwarz
then bounds the expected exceptional contribution by
$O(N^{3/2-D/2})=o(1)$. On the norm event the contribution is
$O(N\eta_c^2)=O(N^{-1+2\delta_1})=o(1)$. Centering its spatial
integral costs at most twice this $L^1$ bound. Apply \cite[Lemma 4.4]{CES} to $f$ and subtract this
centered exterior contribution. The resolvent coefficient in the cited
formula differs from the one here only by the fixed factor $-2\pi\sqrt2$.
This proves \eqref{field:eq:of-micro}. In particular, no compactly supported
solution of $\Delta f=\varphi$ has been assumed when the integral of
$\varphi$ is nonzero.
\end{proof}

\subsection{The small-cutoff estimate}
\label{field:sec:approximation}

Now we show that a small positive cutoff changes the original field average by little in probability.
The microscopic interval has been removed in $L^1$. Above that interval,
we use second moments: first discard a shrinking spatial diagonal,
then apply the separated-point resolvent covariance and restore the
deterministic diagonal. The cutoff may vary with the spatial point.

\begin{theorem}[Small positive cutoffs]
\label{field:thm:of-approximation}
For every fixed positive continuous $a$ on $K_*$, with
$a_*:=\sup_{K_*}a\le a_0$, and every $\delta>0$,
\begin{equation}
 \limsup_{N\to\infty}
 \PP\{|\Ytest_{N,0}(\varphi)-\Ytest_{N,a}(\varphi)|>\delta\}
 \le C_{K_*,\varphi}\delta^{-2}a_*.
 \label{field:eq:of-probability}
\end{equation}
The final constant does not depend on $\inf_{K_*}a$; the eventual size
of $N$ may depend on it.
\end{theorem}
\begin{proof}
Fix the profile $a$ before taking $N\to\infty$. Eventually
$\eta_c<\inf_{K_*}a$. Lemma~\ref{field:lem:of-microscopic} already gives
$\Ytest_{N,\eta_c}(\varphi)-\Ytest_{N,0}(\varphi)\to0$ in $L^1$.
We estimate the second moment only of the remaining increment
\[
 T_{N,a}:=\Ytest_{N,a}(\varphi)-\Ytest_{N,\eta_c}(\varphi)
 =\sqrt2N\int_{K_*}\varphi(z)\int_{\eta_c}^{a(z)}
                         \Im F_z(\eta)\dd\eta\dd A(z).
\]
\paragraph{The shrinking spatial diagonal.}
Fix $0<\nu<\delta_1/100$ and $\xi_{\rm loc},\epsilon<\nu/10$.
For each fixed $z$ and $\eta\in[\eta_c,a_0]$, apply
\eqref{field:eq:of-local-law} with polynomial failure exponent $D>2$.
On its exceptional event use $|F_z(\eta)|\le2/\eta$. Uniformly
in these parameters,
\[
 \E|NF_z(\eta)|^2
 \le C\eta^{-2}(N^{2\xi_{\rm loc}}+N^{2-D})
 \le C N^{2\xi_{\rm loc}}\eta^{-2}.
\]
Minkowski's integral inequality therefore gives
\[
 \E\left|N\int_{\eta_c}^{a(z)}F_z(\eta)\dd\eta\right|^2
 \le C N^{2\xi_{\rm loc}}\left(\int_{\eta_c}^{a(z)}\frac{\dd\eta}{\eta}\right)^2
 \le C N^{2\xi_{\rm loc}}\log^2N.
\]
Remove $|z-w|\le N^{-2\nu}$ from the spatial second-moment integral.
Pointwise Cauchy--Schwarz and the planar volume of the deleted
diagonal bound its contribution by
\begin{equation}
 CN^{2\xi_{\rm loc}-4\nu}\log^2N=o(1).
 \label{field:eq:of-random-diagonal}
\end{equation}
No common event over uncountably many spatial points is needed.

\paragraph{The separated region.}
On $|z-w|>N^{-2\nu}$, the condition of the cited covariance estimate follows from
$\delta_1>4\nu+\xi_{\rm loc}$. Multiplication of
\eqref{field:eq:of-source-error} by $2N^2$ and integration give at most
\begin{equation}
 CN^{\epsilon+8\nu-\delta_1/2}
       \left(\int_{\eta_c}^{a_*}\frac{\dd\eta}{\eta}\right)^2
 \le CN^{\epsilon+8\nu-\delta_1/2}\log^2N=o(1).
 \label{field:eq:of-integrated-error}
\end{equation}
Spatially dependent upper cutoffs only shorten these absolute
integrals. Since $F_z$ is purely imaginary,
$\E(\Im F_z)(\Im F_w)=-\E F_zF_w$. Thus
\begin{align}
 \lim_{N\to\infty}\E T_{N,a}^2
 ={}-&\int_{K_*^2}\varphi(z)\varphi(w)
       \int_0^{a(z)}\int_0^{a(w)}
       \{V(z,w;\eta,\theta)+\kappa U(z,\eta)U(w,\theta)\}
       \dd\eta\dd\theta\dd A(z)\dd A(w).
 \label{field:eq:of-limit-variance}
\end{align}
\paragraph{Restoring the deterministic integral.}
The deleted deterministic diagonal and lower endpoints can be
restored by dominated convergence using \eqref{field:eq:of-majorant}.
Indeed, for $s>0$,
\[
 \int_{K_*}\frac{\dd A(w)}{(|z-w|^2+s)^2}\le\frac\pi s,
 \qquad
 \int_0^{a_*}\int_0^{a_*}\frac{\dd\eta\dd\theta}{\eta+\theta}
       =2a_*\log2.
\]
Consequently
\begin{equation}
 \lim_{N\to\infty}\E T_{N,a}^2
 \le C_{K_*,\varphi}(a_*+a_*^2)\le C_{K_*,\varphi}a_*.
 \label{field:eq:of-variance-bound}
\end{equation}
Absolute values handle signed $\varphi$ and negative $\kappa$ in
this upper bound. The exact expression remains a variance.
Lemma~\ref{field:lem:of-microscopic} handles the microscopic term in
probability. Chebyshev's inequality handles $T_{N,a}$, with a
$\delta/2$ split, and proves \eqref{field:eq:of-probability}.
The bounds use $a_*$ but not $\inf_{K_*}a$; only the threshold ensuring
$\eta_c<\inf_{K_*}a$ depends on that minimum. The exact variance
limit above concerns $T_{N,a}$, without an $L^2$ assertion for the
zero-height microscopic term.
\end{proof}

\begin{remark}[Normalization and the limiting kernel]
The square of the prefactor $\sqrt2N$ cancels $1/(2N^2)$ in
\eqref{field:eq:of-source-covariance}. The minus sign in
\eqref{field:eq:of-limit-variance} comes from the purely imaginary
centered traces.

To identify the kernel, let $u_z(\eta)\in(0,1]$ be determined by
\[
 \eta=(1-u_z(\eta))\sqrt{u_z(\eta)^{-1}-|z|^2},
 \qquad u_z(0)=1,
\]
and express the common covariance \eqref{sp:fb:covariance} in height
variables as
\[
 \mathcal C_\kappa(z,w;\eta,\theta)
 =-\tfrac12\log D(u_z(\eta),z;u_w(\theta),w)
   +\tfrac\kappa2p_z(u_z(\eta))p_w(u_w(\theta)).
\]
Since $m^z(i\eta)^2=-p_z(u_z(\eta))$,
\[
 \partial_\eta\partial_\theta\mathcal C_\kappa
       =-(V+\kappa U U).
\]
For $z\ne w$, the zero-height value is exactly $K_\kappa(z,w)$
in \eqref{sp:kernel}: $D(1,z;1,w)=|z-w|^2$ and
$p_z(1)=1-|z|^2$. Consequently the inner integral, including its
minus sign, in \eqref{field:eq:of-limit-variance} is
\begin{align}
 &\mathcal C_\kappa(z,w;a(z),a(w))
   -\mathcal C_\kappa(z,w;a(z),0)\nonumber\\
 &\qquad-\mathcal C_\kappa(z,w;0,a(w))+K_\kappa(z,w).
 \label{field:eq:of-kernel-increment}
\end{align}
The preceding absolute majorant justifies its spatial integration.
No derivative of $a$ is used: the height integrations hold $z,w$
fixed. The area measure remains $\dd A=\dd x\dd y$ throughout;
there is no extra factor of $\pi$ in the limiting covariance.
\end{remark}

\paragraph{The smoothing-height profile used by Proposition~\ref{p5:transfer}.}
Take $a_S(z)=\eta_S(z)$ from \eqref{sp:fb:scale}. On the fixed
bulk disk, each fixed $S\ge1$ gives a continuous positive profile,
with $\sup_{K_*}a_S\le C_{K_*}e^{-2S}$. For every sufficiently
large fixed integer $S$, this supremum is below $a_0$, and
\eqref{field:eq:of-probability} gives
\begin{equation}
 \limsup_N\PP\{|\Ytest_N(\varphi)-I_{N,S}(\varphi)|>\delta\}
 \le C_{K_*,\varphi}\delta^{-2}e^{-2S}.
 \label{field:eq:of-clock-approximation}
\end{equation}
Here $I_{N,S}=\Ytest_{N,a_S}$ by \eqref{p5:centered-scale}.
The constant is independent of $S$, while the required size of $N$
may depend on the positive minimum of $a_S$. Thus $N$ tends to
infinity first at fixed $S$, and then the integer scale tends to
infinity. The theorem itself retains all continuous positive profiles.

\subsection{Exact cutoffs, moment convergence, and quadrature}
\label{field:sec:quadrature}

This subsection uses covariance convergence to approximate each regularized field integral by finitely many observations, with controlled mean-square error.
We use the fields, clocks and deterministic correction of
\eqref{sp:fb:Phi}--\eqref{sp:fb:scale}, and write $a_s=\eta_s$.
All corrections to $\sqrt2L_{N,a_S}$ are deterministic. Therefore
\begin{equation}
 \Ytest_{N,a_S}(\varphi)=I_{N,S}(\varphi)
 =\int_{K_*}\varphi(z)
       [X_{N,S}(z)-\E X_{N,S}(z)]\dd A(z).
 \label{field:eq:of-centered-coarse}
\end{equation}
This is actual expectation centering, as in \eqref{p5:centered-scale}.
For each fixed $S\ge1$, $a_S$ has a positive minimum on $K_*$,
its maximum is below one, and $\sup_{K_*}a_S\le C_{K_*}e^{-2S}$.
The reference is the common centered Gaussian family, with
\begin{equation}
 \E_G[G_s(z)G_q(w)]=C_\kappa(s,z;q,w),
 \label{field:eq:of-exact-covariance}
\end{equation}
where the complete formula is \eqref{sp:fb:covariance}.
Its covariance is positive semidefinite for every $\kappa\ge-1$
by Section~\ref{chap:gaussian}, including repeated observations.
The whole covariance is used when $\kappa<0$.

\begin{lemma}[Uniform moments at a fixed positive cutoff]
\label{field:lem:of-coarse-moments}
For fixed $S\ge1$, the unweighted finite-dimensional laws of $X_{N,S}$
converge to those of $G_S$. The first and second moments converge
uniformly over spatial arguments in $K_*$, including coincidences
and approaching points. In particular,
\begin{align}
 \sup_{z\in K_*}|\E X_{N,S}(z)|&\longrightarrow0,\nonumber\\
 \sup_{z,w\in K_*}
 |\operatorname{Cov}(X_{N,S}(z),X_{N,S}(w))
                 -C_\kappa(S,z;S,w)|&\longrightarrow0.
 \label{field:eq:of-moments}
\end{align}
For every fixed finite observation list, all prescribed fixed joint
polynomial moments have the same uniform convergence. The moment
bounds are uniform for entry laws with common moment bounds, as
in Section~\ref{field:sec:inputs}; the displayed limit uses the fixed
entry law and its fourth cumulant $\kappa$.
\end{lemma}
\begin{proof}
For fixed $S\ge1$, the continuous profile $a_S(z)=\eta_S(z)$ has a
positive minimum on $K_*$ and a maximum below one. Apply
Lemma~\ref{fixed:lem:positive-height} with this fixed lower bound.
Its first and second polynomial moments give the actual means and
covariances in \eqref{field:eq:of-moments}, uniformly in the spatial
arguments. Its general polynomial statement and higher absolute
moment bounds give all the remaining moment and uniform-integrability
assertions. No distributional limit alone is used to infer these
moments, and no critical determinant weight is present.
\end{proof}

\begin{proposition}[Fixed-cutoff quadrature]
\label{field:prop:of-quadrature}
Fix $S\ge1$ and put $\mu(\dd z)=\varphi(z)\dd A(z)$ on $K_*$.
Let $\nu_h$ be deterministic finite atomic real signed measures on
$K_*$, independent of $N$, such that
\[
 \sup_h\|\nu_h\|_{\mathrm{TV}}<\infty,
 \qquad \nu_h\Longrightarrow\mu\quad(h\downarrow0)
\]
against continuous functions on $K_*$. Define
\begin{align*}
 Q_{N,S,h}(\varphi)
   &=\int[X_{N,S}(z)-\E X_{N,S}(z)]\dd\nu_h(z),\\
 Q_{G,S,h}(\varphi)&=\int G_S(z)\dd\nu_h(z),\qquad
 I_{G,S}(\varphi)=\int_{K_*}\varphi(z)G_S(z)\dd A(z).
\end{align*}
Then
\begin{equation}
 \lim_{h\downarrow0}\limsup_{N\to\infty}
 \E|I_{N,S}(\varphi)-Q_{N,S,h}(\varphi)|^2=0,
 \label{field:eq:of-quadrature}
\end{equation}
and, on the specified Gaussian space,
\begin{equation}
 Q_{G,S,h}(\varphi)\longrightarrow I_{G,S}(\varphi)
       \quad\text{in }L^2\quad(h\downarrow0).
 \label{field:eq:of-gaussian-quadrature}
\end{equation}
\end{proposition}
\begin{proof}
Write $C_{N,S}$ for the actual centered covariance of $X_{N,S}$
and $\sigma_h=\mu-\nu_h$. Fubini's theorem gives the exact identity
\begin{equation}
 \E|I_{N,S}(\varphi)-Q_{N,S,h}(\varphi)|^2
 =\iint_{K_*^2}C_{N,S}(z,w)\dd\sigma_h(z)\dd\sigma_h(w).
 \label{field:eq:of-quadrature-variance}
\end{equation}
The moment bounds justify Fubini, and $\sup_h\|\sigma_h\|_{\mathrm{TV}}$
is finite. Replacing $C_{N,S}$ by $C_\kappa(S,\cdot;S,\cdot)$ costs
at most
\[
 \|\sigma_h\|_{\mathrm{TV}}^2
 \|C_{N,S}-C_\kappa(S,\cdot;S,\cdot)\|_\infty\longrightarrow0.
\]
At fixed $S$, the limiting kernel is continuous on $K_*^2$.
Approximate it uniformly by finite sums of products $f(z)g(w)$.
Signed weak convergence of $\sigma_h$ to zero handles each product,
and the bounded total variations control the uniform approximation
error. This proves \eqref{field:eq:of-quadrature}.
The same identity with the Gaussian covariance proves
\eqref{field:eq:of-gaussian-quadrature}; its integral is equivalently
the corresponding mean-square integral on the same Gaussian space.
\end{proof}

\paragraph{Common nodes for a finite test family.}
A convenient instance uses deterministic finite Borel partitions
$(E_{\ell,h})_{\ell=1}^{q(h)}$ of $K_*$, of mesh at most $h$,
and one node $z_{\ell,h}$ in each nonempty cell. Choose the partitions
and nodes independently of $N$ and of the test function, and set
\begin{equation}
 w_{\ell,h}(\varphi)=\int_{E_{\ell,h}}\varphi(z)\dd A(z),
 \qquad \nu_{h,\varphi}=\sum_\ell
            w_{\ell,h}(\varphi)\delta_{z_{\ell,h}}.
 \label{field:eq:of-common-nodes}
\end{equation}
Then $\|\nu_{h,\varphi}\|_{\mathrm{TV}}\le\|\varphi\|_{L^1}$.
Uniform continuity on $K_*$ shows that $\nu_{h,\varphi}$ converges
weakly to $\mu$. The proposition therefore applies simultaneously
to every test in the finite family.
At each fixed $S,h$, the observation list is finite and fixed before
$N\to\infty$; the mesh is refined afterwards. The proof uses
covariances, with no random spatial derivative estimate.

\subsection{Joint convergence with the original field}
\label{field:sec:transfer}

We now pass from finite observations to original field averages while leaving the other random coordinate unchanged throughout both approximations.
The preceding estimates allow two successive approximations: finite
quadrature at a fixed positive scale, followed by removal of that
scale. Keeping the other random coordinate unchanged in both steps
preserves the specified joint law.

\begin{proposition}[From finite coarse observations to the original field]
\label{field:prop:of-joint}
Fix a closed bulk disk $K_*$ and an integer $k_*\ge1$.
Let $\mathcal E$ be a Polish space, let $V_N$ be random elements of
$\mathcal E$, and let $V$ be defined on the same probability space
as the exact Gaussian family \eqref{field:eq:of-exact-covariance}.
Suppose that, for every fixed finite list of integer clocks
$s_i\ge k_*$ and sites $z_i\in K_*$,
\begin{equation}
 (V_N,(X_{N,s_i}(z_i))_i)
     \Longrightarrow (V,(G_{s_i}(z_i))_i)
 \label{field:eq:of-joint-premise}
\end{equation}
with this one specified limiting coupling. Then, for every fixed
finite family $\varphi_1,\ldots,\varphi_m\in
C_c^\infty(\operatorname{int}K_*;\R)$,
\begin{equation}
 (V_N,(\Ytest_N(\varphi_j))_{j=1}^m)
     \Longrightarrow (V,(G_\kappa(\varphi_j))_{j=1}^m)
 \label{field:eq:of-joint-conclusion}
\end{equation}
in $\mathcal E\times\R^m$.
For a prescribed test family, it suffices to require the premise
only at a fixed spatial family supporting deterministic quadratures
as in Proposition~\ref{field:prop:of-quadrature}. In particular, the
sites may be just the common nodes in \eqref{field:eq:of-common-nodes}.

Here $\Ytest_N=\Ytest_{N,0}$ and $G_\kappa$ is the original
Gaussian generalized field, with
\begin{equation}
 \E_G[G_\kappa(\varphi)G_\kappa(\psi)]
 =\iint_{K_*^2}\varphi(z)\psi(w)K_\kappa(z,w)\dd A(z)\dd A(w),
 \label{field:eq:of-generalized-covariance}
\end{equation}
where $K_\kappa$ is \eqref{sp:kernel}. If $V$ is the canonical
chaos, or another measurable object of this field, that coupling
is preserved. No moment or independence assumption on $V_N$ or
$V$ is required.
\end{proposition}
\begin{proof}
Fix the finite test family. Bold symbols below collect the corresponding
$m$ integrals or quadratures.

\paragraph{A metric with the first coordinate unchanged.}
Choose a bounded compatible metric $d_{\mathcal E}$ on $\mathcal E$
and use
\[
 d((v,x),(v',x'))=d_{\mathcal E}(v,v')+(1\wedge\|x-x'\|)
\]
on $\mathcal E\times\R^m$. When the first coordinate is unchanged,
for every $\delta>0$,
\begin{equation}
 \E d((V_N,x_N),(V_N,y_N))
 \le\delta+\PP\{\|x_N-y_N\|>\delta\}.
 \label{field:eq:of-unchanged-coordinate}
\end{equation}
The same inequality holds on the limiting Gaussian probability space
with $V$ unchanged. Thus the common approximation argument
\eqref{sp:approximation-proof} applies using only errors in the
field coordinates.

\paragraph{Joint convergence at a fixed positive scale.}
For each test choose deterministic finite atomic quadrature measures
$\nu_{h,j}$ permitted by Proposition~\ref{field:prop:of-quadrature}.
When all sites of $K_*$ are allowed, use the common construction
\eqref{field:eq:of-common-nodes}. Under the weaker spatial premise,
use the quadratures supported on that family. At each fixed $h$,
the union of their supports is one finite list of allowed sites.
Put $B_j=\sup_h\|\nu_{h,j}\|_{\mathrm{TV}}<\infty$.

Fix a sufficiently large integer $S\ge k_*$ and then a fixed $h$.
Apply \eqref{field:eq:of-joint-premise} at clock $S$ to this finite
list. The deterministic centering in the $j$th quadrature differs
from zero by at most
\[
 B_j\sup_{z\in K_*}|\E X_{N,S}(z)|\longrightarrow0
\]
by Lemma~\ref{field:lem:of-coarse-moments}. Taking the fixed linear
combinations and using these vanishing deterministic means gives
\begin{equation}
 (V_N,\boldsymbol Q_{N,S,h})
     \Longrightarrow(V,\boldsymbol Q_{G,S,h}).
 \label{field:eq:of-joint-quadrature}
\end{equation}
Sum the squared errors in \eqref{field:eq:of-quadrature} and
\eqref{field:eq:of-gaussian-quadrature} over the finite test family.
At this fixed $S$, these errors vanish as $h\downarrow0$, after
$N\to\infty$ on the matrix side. Apply
\eqref{sp:approximation-proof} with
\eqref{field:eq:of-unchanged-coordinate} to obtain
\begin{equation}
 (V_N,\boldsymbol I_{N,S})
     \Longrightarrow(V,\boldsymbol I_{G,S}).
 \label{field:eq:of-joint-scale}
\end{equation}
No observation list grows with $N$ in the finite-dimensional
comparison.

\paragraph{Identification with the same Gaussian field.}
Write $K_*=\{|z|\le r_*\}$, with $r_*<1$. The exact denominator
in \eqref{sp:fb:covariance} satisfies
\[
 (1-r_*^2)(1-uv)+uv|z-w|^2
 \le D(u,z;v,w)\le1-uv+uv|z-w|^2.
\]
For $S,T\ge1$, $u_S,u_T\ge\tanh1$. The cumulant term is bounded,
so, for $z\ne w$,
\begin{equation}
 |C_\kappa(S,z;T,w)|
 \le C_{K_*,\kappa}+\log_+|z-w|^{-1}.
 \label{field:eq:of-mixed-envelope}
\end{equation}
Off the diagonal the mixed kernel tends to $K_\kappa(z,w)$ as
$S,T\to\infty$. The logarithmic bound is integrable on $K_*^2$.
Apply dominated convergence to the two variances and the mixed
covariance in
\[
 \E_G|I_{G,S}(\varphi)-I_{G,T}(\varphi)|^2.
\]
Each integrated covariance tends to the same limit, so these
integrals are Cauchy in $L^2$. Section~\ref{gaussian:gp2:canonical}
defines the original field by this exact same limit:
\begin{equation}
 G_\kappa(\varphi)
      =L^2\!-\!\lim_{S\to\infty}I_{G,S}(\varphi).
 \label{field:eq:of-same-field}
\end{equation}
Its covariance is \eqref{field:eq:of-generalized-covariance}.
The limit is taken on the same probability space as $V$ and the
entire Gaussian family, retaining their specified coupling,
also at $\kappa=-1$.

\paragraph{Removal of the positive scale.}
Equation~\eqref{field:eq:of-clock-approximation} and a finite union
bound give
\begin{equation}
 \limsup_N\PP\{\|\boldsymbol{\Ytest}_N-
                           \boldsymbol I_{N,S}\|>\delta\}
 \le C_{K_*,\boldsymbol\varphi}\delta^{-2}e^{-2S}.
 \label{field:eq:of-vector-cutoff}
\end{equation}
For fixed $\delta$, the right side tends to zero as the fixed integer
scale $S$ tends to infinity. Equation~\eqref{field:eq:of-same-field}
gives the limiting approximation on the specified Gaussian space.
Apply \eqref{sp:approximation-proof} to
\eqref{field:eq:of-joint-scale}, with matrix target
$(V_N,\boldsymbol{\Ytest}_N)$ and approximation
$(V_N,\boldsymbol I_{N,S})$, and with limiting target
$(V,(G_\kappa(\varphi_j))_j)$. Inequality
\eqref{field:eq:of-unchanged-coordinate} controls both approximation
errors. Taking $N\to\infty$ first, then $S\to\infty$, and finally
$\delta\downarrow0$ proves \eqref{field:eq:of-joint-conclusion}.

We take the limits in the following order:
\[
 \text{fixed }S,h;\qquad N\to\infty;\qquad
 h\downarrow0\text{ at fixed }S;\qquad S\to\infty.
\]
Only sufficiently large integer clocks are used. Both approximations
leave $V_N$ and $V$ unchanged, which is why no extra assumptions on
those coordinates enter the proof. This verifies
Proposition~\ref{p5:transfer}.
\end{proof}

\paragraph{Application to the spatial and raw-marked measures.}
For the additional random element $V_N=(Z_N,\Xi_N)$, the complete application
is Corollary~\ref{p5:field-corollary} and its proof.
Corollary~\ref{p34:joint-corollary} supplies the required joint
measure--coarse premise on one enlarged bulk disk; no compact-window
premise is used. The conclusion retains total mass, the raw-marked
limit $M_\kappa|_K\otimes\rho_{\mathrm R}$, the normalization
$M_\kappa=\sqrt{2/\pi}\,D_\kappa$, and the specified Gaussian
coupling. The finite family of original field tests may have
supports outside the measure disk $K$.

\begin{remark}[Scope of the joint field comparison]
The proposition applies, for example, when $\mathcal E$ is a space
of finite spatial or marked measures, or a continuous measure-valued
path space, equipped with its stated Polish topology. It assumes the
corresponding joint coarse convergence; it does not establish it.
Original field tests may have support anywhere in the bulk disk,
not just in the spatial support of the measure, provided the premise
allows quadratures for those tests. No distributional field test is
inserted under a determinant tilt, and no coupling of all matrix sizes
to one Gaussian realization is asserted.
\end{remark}

\endgroup

\begingroup
\newcommand{\E}{\mathbb E}
\newcommand{\PP}{\mathbb P}
\newcommand{\dd}{\,\mathrm d}
\newcommand{\ind}{\mathbf1}
\newcommand{\D}{\mathcal D}
\newcommand{\csh}{\sqrt{2/\pi}}
\newcommand{\TV}{\mathrm{TV}}
\newcommand{\Var}{\operatorname{Var}}
\newcommand{\Cov}{\operatorname{Cov}}
\section{Uniform exponent comparison and the continuous near-critical window}\label{chap:window}

This section identifys the near-critical limit as a deterministic profile times the same critical measure on each fixed compact parameter interval.

The compact window requires the comparison of determinant-weighted laws
uniformly in the positive exponent, rather than only at the critical
value. We first establish this comparison on a fixed compact interval of
positive exponents, then prove the Gaussian survival estimate for both
signs of the near-critical parameter. The global mesoscopic event and the marked
convergence are separate inputs to the later measure deductions.

The comparison uses Theorems~\ref{scalar:thm:sr-sharp} and
\ref{scalar:thm:sr-regularization},
Proposition~\ref{scalar:prop:sr-moments},
Theorem~\ref{lf:thm:lfl-main} and
Proposition~\ref{lf:prop:lfl-uniform-charges}.
The retained-history passage uses
Lemma~\ref{calibration:lem:sc-history},
Theorem~\ref{replacement:thm:rgfc-main}, and
Theorem~\ref{signed:thm:ss-defect} with
Corollaries~\ref{signed:cor:ss-eight},
\ref{signed:cor:ss-normalizer} and~\ref{signed:cor:ss-ou}.
Each is used with its stated normalization and common moment bounds.

\subsection{Conventions and the precise assertions}
\label{window:sec:definitions}

We first fix the exponent conventions and state the uniform comparison needed to approach criticality from either side.
We use the model, exact density centering, clocks and corrected scale
field of \eqref{main:model}, \eqref{sp:fb:Phi} and
\eqref{sp:fb:scale}. Thus $t=\frac12\log N$,
$Y_N=\sqrt2\Phi_N(\cdot,0)$, and
$Y_N^0=Y_N+\frac{\sqrt2\kappa}{4}(1-|z|^2)^2$ as in
\eqref{sp:fb:endpoints}. The raw field defines the endpoint mark;
the corrected field $X_{N,k}$ defines the barrier
$A_{N,J}^B$ in \eqref{sp:barriers}. The deterministic correction
cancels from every normalized positive root.

The estimates below hold uniformly over complex iid laws with a
prescribed common profile of all fixed moments, including laws
with atoms. The enlarged profile includes the matched surrogates
and their normalized Gaussian convolutions, as in
\eqref{scalar:eq:sr-common}. Their actual fourth cumulants are used
at each stage. For the original law write $\kappa=\E|\xi|^4-2\ge-1$.
Fix the closed bulk disk $K$, the integer $k_0\ge1$, and a closed
bulk disk $K_+$ whose interior contains $K$. All input estimates are
applied on $K_+$. Since $\tau=t^b/N\to0$ for each fixed $b$,
the physical points $v=\sqrt{1+\tau}\,z$, $z\in K$, and the
auxiliary roots used below lie in $K_+$ for all sufficiently large
$N$. The disk is chosen before $b$; only the threshold in $N$
depends on $b$. Area remains $\dd A=\dd x\dd y$.

The same centered Gaussian scale family as in
\eqref{sp:fb:covariance} has covariance
\begin{align}
 C_\kappa(s,z;r,w)
 &=-\tfrac12\log D(u_s,z;u_r,w)
               +\tfrac\kappa2p_z(u_s)p_w(u_r),       \label{window:eq:cw-cov}\\
 D(u,z;v,w)
 &=1-uv-uv[(1-u)|z|^2+(1-v)|w|^2]+uv|z-w|^2.
 \nonumber
\end{align}
Use the complete positive-semidefinite covariance, including when
$\kappa=-1$; no independent negative-variance component is introduced.
Let $v_t(z)=C_\kappa(t,z;t,z)$ and
$A_J^B(z)=\{G_k(z)\le2k+B,\ k_0\le k\le J\}$.
For a positive exponent $\gamma$ of the scaled field, the two normalized roots are
\[
 \frac{\dd Q_{N,z}^{\gamma,\det}}{\dd\PP}
   =\frac{e^{\gamma Y_N(z)}}{\E e^{\gamma Y_N(z)}},
 \qquad
 \frac{\dd Q_{t,z}^\gamma}{\dd\PP_G}
   =e^{\gamma G_t(z)-\gamma^2v_t(z)/2}.
\]
The Gaussian root time $t$ is real. The determinant exponent $\beta$
multiplies $\Phi_N$, so the exponent of the scaled field is
$\gamma=\beta/\sqrt2$. Gaussian roots themselves are defined for
all real exponents, but matrix roots here use only positive exponents.

For each fixed compact interval $|\lambda|\le a$ and all
sufficiently large $N$, define, consistently with
\eqref{sp:window-measure},
\begin{align}
 \gamma_N(\lambda)&=2-\lambda/\sqrt t,&
 q_N(\lambda)&=\sqrt2\,\gamma_N(\lambda)
                  =2\sqrt2-\frac{2\lambda}{\sqrt{\log N}},\nonumber\\
 Z_N(\lambda)(\dd z)
   &=\sqrt t\,\frac{e^{\gamma_N(\lambda)Y_N(z)}}
                       {\E e^{\gamma_N(\lambda)Y_N(z)}}\dd A(z),&
 R_N(z)&=\frac{2t-Y_N(z)}{\sqrt t},                 \label{window:eq:cw-def}\\
 \Xi_N(\dd z,\dd u)&=Z_N(0)(\dd z)\delta_{R_N(z)}(\dd u).
 \nonumber
\end{align}
More precisely, take $t\ge4a^2$. Then
$\gamma_N(\lambda)\in[3/2,5/2]$ and
$\sqrt2\,\gamma_N(\lambda)\in[3\sqrt2/2,5\sqrt2/2]\Subset[2,4]$.
The scalar estimate below makes the normalizers finite and strictly
positive uniformly on this domain. At a determinant zero set $R_N=0$;
the positive-exponent density vanishes there. Thus $Z_N(0)=Z_N$
and $\Xi_N$ agrees with \eqref{sp:unbarred}. No negative
microscopic determinant moment is used.

\begin{theorem}[Uniform comparison for bounded positive field exponents]
\label{window:thm:cw-source}
Fix $\varepsilon\in(0,1)$, $B\in\mathbb R$ and
$J=\lfloor(1-\varepsilon)t\rfloor$. For every fixed $Q>0$, one
common choice of auxiliary parameters gives
\begin{equation}
 \sup_{\beta\in[2,4],\,z\in K}
 \left|Q_{N,z}^{\beta/\sqrt2,\det}(A_{N,J}^B(z))
       -Q_{t,z}^{\beta/\sqrt2}(A_J^B(z))\right|
       =o(t^{-Q}).
 \label{window:eq:cw-comparison}
\end{equation}
The determinant exponent $\beta$ multiplies $\Phi_N$. In particular,
for each fixed $a<\infty$ the comparison is $o(t^{-1/2})$
uniformly for $\gamma_N(\lambda)$, $|\lambda|\le a$.
\end{theorem}

The exponent interval and common moment bounds are part of the
uniform assertion. Its arbitrary fixed inverse power of $t$ comes
from the logarithmic local-factor estimates and the comparison
estimates for signed exponential coefficients;
no polynomial-in-$N$ rate is asserted for their combined conclusion.

\subsection{Uniformity of the one-root comparison}
\label{window:sec:comparison}

Next we make the comparison at one weighted point uniform over a fixed interval of positive exponents.

\subsubsection{Scalar, cap, and exceptional-event bounds}
\label{window:sec:scalar-inputs}

Now we shall choose common moment, cap, and exceptional-event bounds so that one set of estimates works for every exponent in the interval.
Theorem~\ref{scalar:thm:sr-sharp}, obtained from
\cite[Theorem 9.7]{CL}, gives,
uniformly for $\beta\in[2,4]$ and the common moment bounds,
\begin{equation}
 m_N(\beta,z):=\E e^{\beta\Phi_N(z,0)}
 =N^{\beta^2/8}c_\beta
   \exp\left\{\frac{\beta(\beta-2)}8
                    \kappa(1-|z|^2)^2\right\}
       [1+O(N^{-\alpha})],                           \label{window:eq:cw-scalar}
\end{equation}
where $c_\beta>0$ is continuous on the fixed interval.
In particular $cN^{\beta^2/8}\le m_N(\beta,z)\le
CN^{\beta^2/8}$, with uniform constants.
For $\beta\in[2,4]$, write $Q_{N,z}^{\beta,\Phi,\eta}$ for the law
with density $e^{\beta\Phi_N(z,\eta)}/\E e^{\beta\Phi_N(z,\eta)}$
with respect to $\PP$. At zero height this is
$Q_{N,z}^{\beta/\sqrt2,\det}$, since $Y_N=\sqrt2\Phi_N(\cdot,0)$.
At $\eta_{\rm m}=N^{-1}t^{-q}$,
Theorem~\ref{scalar:thm:sr-regularization}, in particular
\eqref{scalar:eq:sr-TV}, gives
\begin{equation}
 \sup_{\beta,z}
 \|Q_{N,z}^{\beta,\Phi,\eta_{\rm m}}
              -Q_{N,z}^{\beta,\Phi,0}\|_{\TV}
       \le C(t^{-q}+N^{-\alpha}).                    \label{window:eq:cw-proxy}
\end{equation}
Its exact Gaussian-minor inequality is valid at $\beta=2$ as well.
The matched and aged laws satisfy the same common moment bounds,
so these bounds include those triangular arrays.

For the compact interval $[2,4]$ of determinant exponents, the common positive
moment interval $[2,16]$ in
Proposition~\ref{lf:prop:lfl-uniform-charges} is fixed once;
it contains all moments used below, including the cap estimate at $r=12$.

One smooth odd nondecreasing $1$-Lipschitz truncation $c_A$
is fixed for all scaled coordinates. It is the identity on
$[-A\log N,A\log N]$ and bounded by $2A\log N$.
For positive scaled-field exponents
$\gamma=\beta/\sqrt2\in[\sqrt2,2\sqrt2]$, its lower-side
positive-weight error is at most $N^{-\sqrt2 A}$.
The upper-side error is bounded by a higher fixed positive moment,
uniformly in $\beta$. By Proposition~\ref{scalar:prop:sr-moments}, the scalar envelope
\[
 \E e^{r\Phi(z,\eta)}\le C_rN^{r(r+2)/8},\qquad r\ge2,
\]
at $r=12$ gives, writing $X=\sqrt2\Phi$ and
$\gamma=\beta/\sqrt2$,
\[
 \E[e^{\gamma X}\ind_{\{X>A\log N\}}]
 \le C N^{21-(12-\beta)A/\sqrt2}.
\]
Thus, by increasing $A$, this fixed moment interval gives any prescribed
polynomial bound on the upper-truncation error.
At mesoscopic coordinates the two-sided shell bound
\eqref{signed:eq:ss-tail}, with the choices in \eqref{signed:eq:ss-shell-budget},
is $N^{-c_0A^2}+N^{-D}$ and controls all errors from the caps and normalizers.
Fix the moment interval and the desired polynomial shell accuracy
first, choose $A$ next, and select the exceptional-event exponent $D$
under the unweighted law after fixing the cap and weight bounds. The same cap serves
every exponent in the interval. These are the cap-removal estimates
of Lemma~\ref{calibration:lem:sc-cap-removal}, with the bounds for positive
exponents made uniform as above.

For any unweighted exceptional event of probability $N^{-D}$,
Cauchy--Schwarz and the uniform scalar lower bound give
\[
 \sup_{\beta,z}Q_{N,z}^{\beta,\Phi,0}(\mathcal E^c)
       \le CN^{C_1-D/2}.
\]
The constant $C_1$ is fixed using the endpoints of the exponent and
moment intervals. No survival probability is a denominator.

\subsubsection{A local factor for every exponent in the interval}
\label{window:sec:local-factor}

Next we separate the local fluctuations from the initial field with an error controlled uniformly over the exponent interval.
We choose the matched normalized surrogate
$M=X_\tau/\sqrt{1+\tau}$ with initial Gaussian age comparable to
$\sigma=N^{-1/200}$ and $\tau=t^b/N$. Matching is imposed on $M$.
Use the physical height
$\eta_{\rm m}^{\rm phys}=\sqrt{1+\tau}\,\eta_{\rm m}$;
the local-factor theorem permits this bounded positive multiplier.
At the physical point $v=\sqrt{1+\tau}z$,
Theorem~\ref{lf:thm:lfl-main}, equation~\eqref{lf:eq:lfl-lf1},
gives
\[
 \Phi_\tau(v,\eta_{\rm m}^{\rm phys})
       =C_0+S+O(\rho_N),\qquad
 C_0=\Phi_0(v,B_v\tau),\quad B_v=\sqrt{1-|v|^2},
\]
outside an event whose probability is smaller than any prescribed
inverse power of the matrix size, with
$\rho_N=o(t^{-Q_*})$ for a chosen fixed $Q_*$.
To apply that theorem in its stated two-root form, append a
deterministic auxiliary root at distance $cN^{-1/4}$ from $v$
inside $K_+$, then discard its local factor. Its separation satisfies
the theorem's upper and lower restrictions for large $N$.
This is a construction for each prescribed root with uniform bounds,
not an independent family indexed by every point of the continuum.
The local family generating $S$ is independent of the entire
initial matrix. Fix a common cap event $\mathcal K$ depending only on
the local reference family,
and put $U_\beta=e^{\beta S}\ind_{\mathcal K}$.
Let $\Omega$ be the intersection of the local-factor reconstruction event,
$\mathcal K$, and the common scalar/local-law events used here. The
local-cap estimates in Section~\ref{lf:sec:lfl-weights} give $\PP(\Omega^c)\le C_DN^{-D}$
for every prescribed fixed $D$, with one choice uniform in the exponent.

The pointwise reconstruction does not depend on the exponent. On $\Omega$
the ratio of the two positive weights lies between
$e^{-4\rho_N}$ and $e^{4\rho_N}$ for every $\beta\in[2,4]$.
On its complement, uniform positive moments and the common local
cap give a fixed polynomial envelope. Indeed for a real $x$ and
$2\beta\in[4,8]$,
$e^{2\beta x}\le e^{4x}+e^{8x}$.
After dividing by \eqref{window:eq:cw-scalar}, this proves
\begin{equation}
 \sup_{\beta\in[2,4],\,z\in K}
 \frac{\E|e^{\beta\Phi_\tau(v,\eta_{\rm m}^{\rm phys})}
                      -U_\beta e^{\beta C_0}|}
       {\E e^{\beta\Phi_\tau(v,\eta_{\rm m}^{\rm phys})}}
       =o(t^{-Q_*}).                                 \label{window:eq:cw-L1}
\end{equation}
More explicitly, \eqref{lf:eq:lfl-cap-envelope} and the
positive moments on $[2,16]$ bound the weighted complement, after
normalization, by $CN^{C_{\rm cap}-D/2}$, where $C_{\rm cap}$
is fixed uniformly in $\beta,z$. Choose $D>2C_{\rm cap}+2$.
The good-event error is $O(\rho_N)$ directly. This is the
calculation of \eqref{lf:eq:lfl-one-bad}--\eqref{lf:eq:lfl-L1-weight}
with the uniform charges of Proposition~\ref{lf:prop:lfl-uniform-charges}.
The failure exponent under the unweighted law is chosen after the
uniform cap and moment bounds. This uses a uniform deterministic
bound, not a union over uncountably many exponent values.

Independence gives an estimate for
$\E U_\beta\,\E e^{\beta C_0}$ from \eqref{window:eq:cw-L1}.
For nonnegative weights $a,b$,
\[
 \E\left|\frac a{\E a}-\frac b{\E b}\right|
          \le\frac{2\E|a-b|}{\E a}.
\]
Therefore the tilted density on the joint probability space is uniformly
close to the product density
\[
 \frac{U_\beta e^{\beta C_0}}{\E U_\beta\,\E e^{\beta C_0}}.
\]
On every bounded initial-matrix-measurable test the local factor
cancels exactly. In corrected coordinates
$\sqrt2(\Phi-E)$, where $E$ is the deterministic mean correction,
the local factor becomes $U_\beta e^{\beta(E_0-E_{\rm m})}$.
This deterministic multiplier cancels from its normalized density
in the same way. All logarithms above use exact density centering;
variance scaling is the identity in Lemma~\ref{scalar:lem:sr-variance}.

\subsubsection{The retained smooth comparison}
\label{window:sec:retained}

Now the goal is to compare smooth tests of the exponentially weighted matrix and Gaussian histories, with one error bound across the exponent interval.
We fix the bulk disks, the common moment bounds, the interval $[2,4]$ of determinant exponents, the retained fraction $\varepsilon\in(0,1)$,
the barrier $B$, and a desired accuracy $Q>0$. Put
$J=\lfloor(1-\varepsilon)t\rfloor$ and use the retained vectors
\[
 \mathbf X_N(z)=(X_{N,k}(z))_{k_0\le k\le J},\qquad
 \mathbf G(z)=(G_k(z))_{k_0\le k\le J}.
\]
For a smooth bounded test on this $O(t)$-dimensional list, use
\[
 \D_m(F)=\sum_{r=0}^m\sum_{i_1,\ldots,i_r}
       \|\partial_{i_1}\cdots\partial_{i_r}F\|_\infty.
\]
Fix $m=8$ by Corollary~\ref{signed:cor:ss-eight}; this also
controls the order-five norm in
Theorem~\ref{replacement:thm:rgfc-main}. The derivative order is
fixed before any logarithmic time exponent. The root and all its
history coordinates remain at the same physical site until the
Gaussian comparison.

First we replace the height-zero root by its positive-height proxy
using \eqref{window:eq:cw-proxy}, then introduce the common caps
of Section~\ref{window:sec:scalar-inputs}. Apply
Theorem~\ref{replacement:thm:rgfc-main} to the original law and
the matched normalized surrogate of
Proposition~\ref{replacement:prop:rgfc-age}, with coefficient vector
$\gamma e_{r_z}$, $\gamma=\beta/\sqrt2\le2\sqrt2<6$.
Here $e_{r_z}$ denotes the coordinate of the microscopic root.
The centers are fixed and identical for the two replacement laws;
the return to their actual fourth cumulants is a deterministic
shift of size $O(\sigma+\tau)$ by
Lemma~\ref{replacement:lem:rgfc-centers} and the exact variance return.
Choose the replacement height slack $0<\delta<1/1200$ once.
For every fixed $q$, the microscopic height $N^{-1}t^{-q}$
and all retained heights lie in $[N^{-1-\delta},1]$ for large
$N$; the upper endpoint follows from $k_0\ge1$.
The theorem gives a polynomial-in-$N$ error times $\D_5(F)$
and a relative transfer of the scalar normalizer.
The positive cap lower bound is used before division.
Remove the microscopic positive cap after replacement, before
the history transport and local-factor reconstruction, using the
bounds at the full exponent in Section~\ref{window:sec:scalar-inputs}.
Normalization along each ray of coefficients takes place with capped
coordinates; no microscopic regularization bound is applied to determinant
exponents below $2$ along those rays.

Next we use exact variance scaling
(Lemma~\ref{scalar:lem:sr-variance}) to pass to the physical point
$v=\sqrt{1+\tau}\,z$ and physical heights. This identity cancels
the two $N\log(1+\tau)/2$ contributions. The least retained
height satisfies
\[
 \eta_J(z)\asymp N^{-(1-\varepsilon)},\qquad
 \tau/\eta_J(z)\le Ct^bN^{-\varepsilon}.
\]
It exceeds the lower height in
Lemma~\ref{calibration:lem:sc-history} for large $N$. That lemma,
including its exact clock, spatial and mean returns, gives
\begin{equation}
 \rho_{\rm hist}\le C_D\left[
 t^{1+b/2}N^{-\varepsilon/2}
 +t^{b+1}N^{-\varepsilon}
 +t^{b+2}N^{-2\varepsilon}+N^{-c}\right].
 \label{window:eq:cw-flow}
\end{equation}
The stochastic square-root term is retained. On the common event,
the test error is at most $\rho_{\rm hist}\D_1(F)$; on its
complement, the tilted exceptional-event bound in
Section~\ref{window:sec:scalar-inputs} applies. All transformed
physical points remain in the previously fixed $K_+$.

After this transport the test is measurable with respect to the
entire initial matrix. Equation~\eqref{window:eq:cw-L1} and
independence of the local factor therefore replace the surrogate
root by the normalized coarse root $e^{\beta C_0}/\E e^{\beta C_0}$.
Only the root density is reconstructed in this step, so its
$o(t^{-Q_*})$ error multiplies $\|F\|_\infty$, not a derivative
norm. For the coarse-root comparison use the mesoscopic caps;
\eqref{signed:eq:ss-tail}, above its fixed logarithmic height floor,
controls their removal after the estimate for signed exponential coefficients. The exponent
$b$ is chosen to satisfy that floor. The microscopic cap and proxy
are removed only at the full exponent $\beta\in[2,4]$.

Apply Theorem~\ref{signed:thm:ss-defect}, the resulting comparison of
normalizers along coefficient rays, and its Gaussian comparison,
Corollaries~\ref{signed:cor:ss-normalizer} and~\ref{signed:cor:ss-ou},
to this initial-matrix vector, using its actual cumulant $\kappa_0$.
The smooth comparison error is
\begin{equation}
 C_b[t^{2-b/2}+t^{17-b}+N^{-c}]\D_8(F).
 \label{window:eq:cw-Stein}
\end{equation}
All exponent coefficients are bounded on the fixed compact interval, so
one choice of caps and constants works throughout it.
The coarse root clock is
$\ell_v=t-(b/2)\log t+O(1)$, whereas the last retained clock
is $(1-\varepsilon)t+O(1)$. The exact same-site kernel gives
\[
 \sup_{k\le J}|C_{\kappa_0}(k,v;t,v)
                  -C_{\kappa_0}(k,v;\ell_v,v)|
 \le Ce^{-2(\ell_v-J)}\le Ct^bN^{-\varepsilon}.
\]
The rooted mean is linear in $\gamma$, so this estimate returns
the root clock to the real time $t$, uniformly in the exponent;
this is the prefix estimate \eqref{gaussian:gpm:prefix-decay}.
Since $\kappa_M=\kappa_0/(1+\tau)^2$ and
$\kappa_M-\kappa=O(\sigma)$, we have
$\kappa_0-\kappa=O(\tau+\sigma)$. Returning the common spatial
scaling, exact clocks and actual cumulant changes the observed means and covariance entries by
$O(\tau+\sigma)$. For the spatial part, at one site write
$D(u,z;v,z)=a_{u,v}-b_{u,v}|z|^2$, where $a_{u,v}=1-uv$ and
$b_{u,v}=uv(2-u-v)$ satisfy $0\le b_{u,v}\le a_{u,v}$. On $K_+$,
$D\ge(1-r_+^2)a_{u,v}$ with $r_+=\max_{z\in K_+}|z|<1$,
so $|\nabla_z\log D|\le C$ uniformly in the clocks. The
functions $p_z$ and their spatial derivatives are bounded too.
Together with the actual-cumulant change $O(\sigma+\tau)$, this
proves the stated uniform return bound. The common integer
observation labels are kept. Gaussian interpolation costs
$C(\tau+\sigma)(\D_1(F)+\D_2(F))$; it needs no inverse covariance.

For clarity, let
\[
 \Delta_N=N^{-c}+\rho_{\rm hist}+t^bN^{-\varepsilon}+\tau+\sigma.
\]
Absorb the fixed polynomial scalar, replacement, cap and exceptional
errors into $N^{-c}$, decreasing $c>0$ if necessary. The preceding
finite sequence of estimates yields
\begin{align}
 &\sup_{\beta\in[2,4],\,z\in K}
 \left|\E_{Q_{N,z}^{\beta/\sqrt2,\det}}F(\mathbf X_N(z))
       -\E_{Q_{t,z}^{\beta/\sqrt2}}F(\mathbf G(z))\right|\nonumber\\
 &\quad\le C_b[t^{2-b/2}+t^{17-b}+\Delta_N]\D_8(F)
       +[o(t^{-Q_*})+Ct^{-q}]\|F\|_\infty.
 \label{window:eq:cw-smooth-bound}
\end{align}
This displays which errors are multiplied by derivative bounds and
which only by the supremum norm of the test. All logarithmic exponents are fixed
before $N$ tends to infinity.

\subsubsection{Removing the smoothing and choosing the parameters}
\label{window:sec:desmoothing}

The next goal is to replace smooth tests by the retained barrier event, choosing parameters that preserve the required uniform comparison error.
The exact same-site covariance satisfies
\[
 C_\kappa(k,z;r,z)=\min(k,r)+O(1),\qquad
 ck\le v_k(z)\le Ck,\quad k\ge k_0.
\]
The lower bound includes $\kappa=-1$: for early times,
\[
 v_k(z)\ge\tfrac12[-\log(1-u_k^2)-u_k^2]>0.
\]
The Gaussian tilt changes the mean, but not the variance. Thus
for any deterministic centers $b_k$ and $0<h\le1$,
\begin{equation}
 Q_{t,z}^{\beta/\sqrt2}
   \{\exists k\in[k_0,J]:|G_k(z)-b_k|\le h\}
 \le C h\sum_{k=k_0}^J k^{-1/2}
 \le C h\sqrt t .
 \label{window:eq:cw-strips}
\end{equation}
The bound is uniform in $\beta,z$ and uses only one-coordinate
densities, with no nonsingularity assumption on the joint vector.

Take smooth nonincreasing one-dimensional cutoffs, equal to one
below a threshold and zero a distance $h$ above it. Their products
provide $F_-\le\ind_{A_J^B}\le F_+$, with transitions confined
to the strips in \eqref{window:eq:cw-strips} about $b_k=2k+B$.
Choose the common cap level $A>2$; for large $N$ it is the identity
at every transition $2k+B\pm h$, so the capped comparison does
not shift the event. For $h=t^{-s}$ the sum over coordinate
derivatives gives
\[
 \|F_\pm\|_\infty\le1,\qquad
 \D_m(F_\pm)\le C_m\sum_{r=0}^m t^r h^{-r}
                  \le C_m t^{m(1+s)}.
\]
Applying \eqref{window:eq:cw-smooth-bound} to both products and
using the strip bound brackets the difference of barrier
probabilities by
\begin{equation}
 C_b[t^{2-b/2}+t^{17-b}+\Delta_N]t^{m(1+s)}
      +o(t^{-Q_*})+Ct^{-q}+Ct^{1/2-s},\qquad m=8.
 \label{window:eq:cw-desmoothed}
\end{equation}
This is the first passage to the nonsmooth barrier event.

Here is one finite order of choices making this bound $o(t^{-Q})$.
Choose $Q_*>Q+2$ and $s>Q+2$. Fix the common exponent and moment
intervals, cap accuracies, cap levels and exceptional-event exponents,
as in Section~\ref{window:sec:scalar-inputs}, together with the
replacement slack and the fixed derivative order $m=8$.
Follow Section~\ref{lf:sec:lfl-parameters} with reconstruction
precision $P=Q_*$: first fix the weak-law exponents $(M,m_0,a_{\rm LF})$, where
$a_{\rm LF}$ is the exponent called $a$ in that schedule; then $q>Q_*+2$ and the lower cutoff exponent $c$ satisfying that
schedule. Choose $b$ above its stated floors, the mesoscopic
scalar floor, and
\[
 b>2[m(1+s)+Q+4],\qquad b>m(1+s)+Q+19.
\]
Choose the remaining local-factor exponents in their proved order,
$(j,k_1,R)$, then $\ell$, then $k$.
Increasing $b$ at its designated step does not alter an earlier
cutoff; later choices leave the fixed moment interval, cap, derivative
order and lower cutoff unchanged.

Indeed the two terms from the comparison for signed coefficients in
\eqref{window:eq:cw-desmoothed} are then $o(t^{-Q})$.
Every term in $\Delta_N$ is a fixed power of $t$ times a decaying
power of $N$, and hence still tends to zero after multiplication by
the displayed derivative bound.
The reconstruction, proxy and strip terms have exponents
$Q_*,q,s-1/2>Q$. This proves
Theorem~\ref{window:thm:cw-source} uniformly on the full exponent
interval. The conclusion is an arbitrary fixed inverse-logarithmic
accuracy, not a claimed polynomial-in-$N$ rate.

For $t\ge4a^2$, the near-critical determinant exponents belong to
$[3\sqrt2/2,5\sqrt2/2]\Subset[2,4]$. Thus the same construction
treats both signs of $|\lambda|\le a$ without a separate choice
of parameters depending on the exponent.

\subsection{Gaussian survival with both signs of the slope}
\label{window:sec:survival}

In this subsection we bound Gaussian barrier survival uniformly for both signs of the near-critical parameter, supplying the moment bound for restricted masses.

\begin{theorem}[Near-critical Gaussian survival]
\label{window:thm:cw-survival}
For fixed $a,B,k_0,K,\kappa$, uniformly for $|\lambda|\le a$,
$z\in K$, and integers $k_0+1\le J\le t$,
\[
 Q_{t,z}^{\,2-\lambda/\sqrt t}(A_J^B(z))
       \le C_{a,B}/\sqrt{1+J}.
\]
\end{theorem}

\begin{proof}
All constants may depend on the fixed $k_0,K,\kappa$ as well as
$a,B$. Under the root, the mean at time $k$ is
$(2-\lambda/\sqrt t)C_\kappa(k,z;t,z)$. Set
$v=\lambda/\sqrt t$ and subtract this mean. The centered vector,
denoted $H_k$, has covariance $C_\kappa(k,z;l,z)=\min(k,l)+O(1)$.
The event is contained in
\[
 \{H_k\le b_k:=B+C_a+vk,\ k_0\le k\le J\},
 \qquad |v|\sqrt J=|\lambda|\sqrt{J/t}\le a.
\]
Here the fixed compact slope range is used before any integral
estimate; the bound on the mean error is uniform even for bounded
$t$. We compare this vector to Brownian motion multiplicatively,
so that no additive error can mask a survival probability of order
$J^{-1/2}$.

\paragraph{Gaussian rectangle comparison.}
We choose $M_0\ge1$ bounding the covariance difference from
$\min(k,l)$, uniformly in the parameters. Add independent
$\xi_k\sim N(0,3M_0)$ to $H_k$, and compare with
\[
 B_k+Z+\eta_k,\qquad Z\sim N(0,2M_0),\qquad
 \eta_k\sim N(0,\operatorname{Var}H_k-k+M_0).
\]
All added variables and the Brownian motion $B$ are independent.
The last variances lie in $[0,2M_0]$. The two vectors have equal
diagonal variances, while the second has larger off-diagonal
covariances. Interpolating expectations of products of smooth
nonnegative decreasing cutoffs gives the rectangle order: the
mixed second derivatives are nonnegative and the diagonal
covariance derivatives vanish. Approximation of the cutoffs
then gives the same order for closed rectangles. This reasoning
also applies to singular covariance matrices, by adding a common
small independent Gaussian variance and taking it to zero.

We put $d_k=(k-k_0)\wedge(J-k)$ and
$a_k=A_0\sqrt{\log(2+d_k)}$, with $A_0$ fixed sufficiently large.
There are at most two indices at each distance from the endpoints.
Gaussian tails give
\[
 \sum_{k=k_0}^J\PP(\xi_k>a_k)\le C,
 \qquad \PP(\xi_k\le a_k\text{ for every }k)\ge c>0,
\]
where each individual factor is at least $1/2$, and the product
bound follows from independence and summable tails, uniformly in $J$. The variable
\[
 W=|Z|+\max_{k_0\le k\le J}(-\eta_k-a_k)_+
\]
is independent of $B$ and satisfies
$\PP(W>r)\le Ce^{-cr^2}$ for $r\ge0$, uniformly in $J,z$.
Indeed the same summation over distances bounds its maximum,
and the single variable $Z$ has a Gaussian tail.
Independence of $\xi$ from $H$, rectangle comparison, and the
inequality $-Z-\eta_k\le W+a_k$ now give
\begin{equation}
 \PP(H_k\le b_k\text{ for every }k)
 \le C\E_W\PP(B_k\le b_k+2a_k+W\text{ for every }k).
 \label{window:eq:cw-rectangle}
\end{equation}
The multiplicative constant is independent of the horizon.

\paragraph{From integer observations to a continuous curve.}
Conditional on the Brownian values at the integers, the centered
unit bridges are independent of that skeleton and of one another.
Their one-sided maximum has a Gaussian tail. Choose a bridge
bound of the same form $A_1\sqrt{\log(2+d)}$, where $d$ is the
distance of the unit interval from the nearer endpoint of
$[k_0,J]$. The probability that all bridges satisfy these bounds
is at least a constant $c_1>0$, by the same summable-tail product
argument. On that event the linear term $vs$ interpolates exactly.
For any fixed $0<\alpha<1/2$, the remaining logarithmic envelopes
are bounded by $C_2+C_1((s-k_0)\wedge(J-s))^\alpha$.
Consequently \eqref{window:eq:cw-rectangle} is bounded by
\begin{equation}
 C\E_W\PP\{B_s\le H+vs+
       C_1((s-k_0)\wedge(J-s))^\alpha,
                              \ k_0\le s\le J\},
 \label{window:eq:cw-Brownian-reduction}
\end{equation}
where $H=\max(1,B+C_a+C_2+W)$. This passage divides by $c_1$;
it adds no error. In particular $\E(1+H)^2\le C_{a,B}$.

\paragraph{The drifted endpoint integral.}
Fix $H$, put $\ell=J-k_0$, and first suppose $J\ge2k_0+2$.
Under the law of $U_s=B_s-vs$, the drift is $-v$. Relative to
centered Brownian motion on $[0,J]$, its likelihood is
$e^{-vU_J-v^2J/2}$. Condition in that centered law on
$U_{k_0}=y$ and $U_J=x$. The curved event forces $y,x\le H$.
Choose a fixed lift $h=H+C_3+1$.
The upper Brownian-bridge bound in
\cite[Proposition 2.1, (2.1)]{CHL}, applied after translation
and reflection, gives the following upper bound for the transition
density restricted by the curved event:
\begin{equation}
 \frac{C(h-y)(h-x)}{\ell}\varphi_\ell(x-y)\dd x,
 \qquad \varphi_s(r)=(2\pi s)^{-1/2}e^{-r^2/(2s)}.
 \label{window:eq:cw-killed}
\end{equation}
To match its curve parameters, choose $\delta\in(0,1/2)$ once
so that $1/2-\delta\ge\alpha$ and $\delta^{-1}\ge C_1$.
Take $C_3\ge C_1$ to cover distances below one. If
$d=s\wedge(\ell-s)$, the reflected curve is bounded below by
$-\delta^{-1}d^{1/2-\delta}$. Proposition 2.1 bounds the
probability of staying above that curve by $Cxy/\ell$ for
all endpoint deficits $x,y>\delta$ and all $\ell>0$. Here
the deficits are $h-y,h-x>1$. Its upper bound (2.1) imposes
no restriction on their product; the product restriction in
that proposition belongs to its separate lower bound (2.2).

Keeping the Gaussian density in \eqref{window:eq:cw-killed},
the probability inside \eqref{window:eq:cw-Brownian-reduction}
is at most
\begin{equation}
 \frac{Ce^{-v^2J/2}}{\ell}
 \int_{y\le H}\varphi_{k_0}(y)(h-y)
       \int_{x\le H}e^{-vx}(h-x)\varphi_\ell(x-y)\dd x\dd y.
 \label{window:eq:cw-drift-integral}
\end{equation}
The factor $e^{-vx}$ may grow on the negative half-line when
$v>0$, so it cannot be bounded by dropping the Gaussian density.
Instead complete the square:
\[
 e^{-vx}\varphi_\ell(x-y)
       =e^{-vy+v^2\ell/2}\varphi_\ell(x-y+v\ell).
\]
The inner integral is bounded by
$C e^{-vy+v^2\ell/2}(h+|y|+|v|\ell+\sqrt\ell)$.
For the remaining integral,
\[
 e^{-vy-v^2k_0/2}\varphi_{k_0}(y)
          =\varphi_{k_0}(y+vk_0).
\]
This Gaussian has uniformly bounded first two moments because
$k_0$ is fixed and $|v|\sqrt J\le a$. It follows that
\eqref{window:eq:cw-drift-integral} is at most
\begin{align*}
 \frac{C_a}{\ell}
    \big[(1+h)^2+(1+h)(|v|\ell+\sqrt\ell)\big]
 &\le \frac{C_a(1+H)^2}{\sqrt J}.
\end{align*}
Here $\ell\asymp J$, $|v|\ell\le a\sqrt J$, and $h-H$ is
fixed. The calculation uses $|v|$, so it covers both signs.
For the remaining bounded horizons $k_0+1\le J<2k_0+2$, the
probability bound $1$ is absorbed into $C_{a,B}/\sqrt{1+J}$.
Finally integrate over $W$, using the uniform bound on
$\E(1+H)^2$. This proves the stated estimate for every allowed
integer $J$ and both signs of $\lambda$. The Markov and bridge
properties used above belong only to the comparison Brownian
motion; none is assumed for the exact Gaussian scale family.
\end{proof}

\subsection{The global event and uniform marginal tightness}
\label{window:sec:masses}

This subsection shows that a common spatial barrier event gives marginal tightness uniformly over a fixed parameter interval.
The rooted comparison controls one prescribed spatial point. To
remove a barrier simultaneously over the measure disk, we use the
event under the unweighted matrix law from \eqref{p34:prefix}. For a fixed
$\varepsilon\in(0,1)$ and $J=\lfloor(1-\varepsilon)t\rfloor$,
write that event as
\[
 \mathcal E_N(B,\varepsilon)
 =\{X_{N,k}(z)\le2k+B\text{ for every }z\in K,
                                   \ k_0\le k\le J\}.
\]
The hypothesis needed for the following proposition is
\begin{equation}
 \limsup_{N\to\infty}\PP(\mathcal E_N(B,\varepsilon)^c)
       \le p(B),\qquad p(B)\longrightarrow0\quad(B\to\infty).
 \tag{GN}\label{window:eq:cw-global}
\end{equation}
Proposition~\ref{p34:global} supplies
this estimate for every fixed $\varepsilon$, with
$p(B)\le C_K(1+B)e^{-2B}$ for all sufficiently large $B$,
independently of
$\varepsilon$. This estimate is verified by
Corollary~\ref{aggregate:cor:au-maximum}. This is an unweighted
global estimate for the original matrix law; a rooted one-site
comparison alone does not imply it.

\begin{proposition}[Marginal tightness at nearby exponents]
\label{window:prop:cw-P}
Assume \eqref{window:eq:cw-global}. For every fixed $a<\infty$,
\begin{equation}
 \lim_{L\to\infty}\limsup_N\sup_{|\lambda|\le a}
       \PP\{Z_N(\lambda)(K)>L\}=0.
 \tag{P}\label{window:eq:cw-P}
\end{equation}
\end{proposition}

\begin{proof}
We fix the fraction $\varepsilon$ in \eqref{window:eq:cw-global}
and a barrier $B$. For large $N$, $J\ge k_0+1$ and the full
window has positive determinant exponents in $[2,4]$. Define the
measure restricted to the barrier by
\[
 Z_N^{B,J}(\lambda)(\dd z)
       =\ind_{A_{N,J}^B(z)}Z_N(\lambda)(\dd z).
\]
Theorem~\ref{window:thm:cw-source}, with any fixed $Q>1/2$,
and Theorem~\ref{window:thm:cw-survival} give, uniformly for
$|\lambda|\le a$ and $z\in K$,
\[
 Q_{N,z}^{\gamma_N(\lambda),\det}(A_{N,J}^B(z))
       \le \frac{C_{a,B}}{\sqrt{1+J}}+o(t^{-Q}).
\]
Tonelli's theorem and the exact normalized-root identity yield
\begin{align}
 \E Z_N^{B,J}(\lambda)(K)
 &=\sqrt t\int_K
   Q_{N,z}^{\gamma_N(\lambda),\det}(A_{N,J}^B(z))\dd A(z)
       \nonumber\\
 &\le |K|\sqrt t\left[\frac{C_{a,B}}{\sqrt{1+J}}
                       +o(t^{-Q})\right]
       \le C_{a,B,\varepsilon,K}+o(1).
 \label{window:eq:cw-stopped-moment}
\end{align}
Here $\sqrt t/\sqrt{1+J}=O_\varepsilon(1)$ and
$\sqrt t\,o(t^{-Q})=o(1)$. This bound concerns the restricted
mass and is uniform in the exponent.

On $\mathcal E_N(B,\varepsilon)$ the restricted and original
measures agree simultaneously for every $\lambda$: the event
depends only on the corrected scale field, not on the exponent.
Markov's inequality under the unweighted matrix law therefore gives
\begin{equation}
 \limsup_{N\to\infty}\sup_{|\lambda|\le a}
   \PP\{Z_N(\lambda)(K)>L\}
       \le p(B)+\frac{C_{a,B,\varepsilon,K}}L.
 \label{window:eq:cw-marginal-bound}
\end{equation}
For each fixed $B$, first let $L\to\infty$, then let
$B\to\infty$. The constant is allowed to depend on $B$.
Only one fixed retained fraction is used; no completion of the
last $\varepsilon t$ portion of the history is required.
\end{proof}

The raw expectations remain exactly
$\E Z_N(\lambda)(K)=\sqrt t\,|K|$.
Proposition~\ref{window:prop:cw-P} asserts uniform marginal
tightness. It does not bound these expectations or assert
curve tightness. The exponential-tail and path deductions use
the exact density identity separately.

\paragraph{The two masses in Proposition~\ref{cw:inputs}.}
Fix $b>0$ and take $\varepsilon=1/2$, so that
$J=\lfloor t/2\rfloor$. The common global event is available
by Proposition~\ref{p34:global}. Apply
\eqref{window:eq:cw-stopped-moment} to both $\lambda=b$ and
$\lambda=-b$ and add their restricted masses. On the same event
the sum equals $S_N(b)$ from \eqref{cw:endpoint-sum}. Thus
\begin{equation}
 \limsup_{N\to\infty}\PP\{S_N(b)>H\}
       \le p(B)+\frac{C_{b,B,K}}H.
 \label{window:eq:cw-two-masses}
\end{equation}
There is one exceptional-event term; no independence of the two
masses is needed. Let $H\to\infty$ with $B$ fixed, then
$B\to\infty$. This proves \eqref{cw:endpoint-tightness} while
retaining the stronger proposition above.

\subsection{Normalization and exponential tails of the marks}

This subsection  combines exact normalization with control of large marks to pass from the critical marked measure to the nearby exponent curve.

\subsubsection{The uniform normalizer correction}
\label{window:sec:normalizer}

Now we isolate the deterministic normalizer correction that turns a change of exponent into exponential reweighting of the raw height mark.
We use the raw field $Y_N$ and the definition of $\psi_N$ in
\eqref{sp:window-measure}. Set $\Gamma=[3/2,5/2]$.
The interval $\sqrt2\Gamma$ of determinant exponents lies inside $[2,4]$, so
\eqref{window:eq:cw-scalar} gives, uniformly on $\Gamma\times K$,
\begin{align}
 \psi_N(\gamma,z)
 &=\tfrac12\gamma^2t+H_\kappa(\gamma,z)+r_N(\gamma,z),
                                                        \label{window:eq:cw-psi}\\
 H_\kappa(\gamma,z)
 &=\log c_{\sqrt2\gamma}
    +\frac{\kappa\gamma(\gamma-\sqrt2)}4(1-|z|^2)^2,
 \qquad \epsilon_N:=\sup_{\Gamma\times K}|r_N|
       =O(N^{-\alpha}).\nonumber
\end{align}
Taking logarithms of the scalar relative error is legitimate
uniformly for large $N$, since its positive leading factor is
bounded away from zero after division by $N^{\gamma^2/4}$.
The correction $H_\kappa$ is bounded and jointly continuous on
this compact domain. For $Y_N^0$ the correction would instead be
$\log c_{\sqrt2\gamma}+\kappa\gamma^2(1-|z|^2)^2/4$.
The raw correction in \eqref{window:eq:cw-psi} is the one
consistent with the raw mark $R_N$.

\begin{lemma}
\label{window:lem:normalizer}
For every fixed $b>0$, the error in
\eqref{cw:normalizer-error} satisfies \eqref{cw:normalizer-control}.
More precisely, if $\omega_H$ is the uniform modulus of continuity
of $H_\kappa$ in its first variable on $\Gamma\times K$, then
for $t\ge4b^2$,
\begin{equation}
 \delta_N(b)\le\omega_H(b/\sqrt t)+2\epsilon_N
                   \longrightarrow0.
 \label{window:eq:cw-modulus}
\end{equation}
\end{lemma}

\begin{proof}
Write $\gamma=\gamma_N(\lambda)=2-\lambda/\sqrt t$.
The quadratic part of the normalizer difference is exactly
\[
 \tfrac12(4-\gamma^2)t
             =2\lambda\sqrt t-\lambda^2/2.
\]
We substitute \eqref{window:eq:cw-psi} into
\eqref{cw:normalizer-error}. These quadratic terms cancel, leaving
\[
 e_N(\lambda,z)=H_\kappa(2,z)-H_\kappa(\gamma,z)
                         +r_N(2,z)-r_N(\gamma,z).
\]
For $|\lambda|\le b$ and $t\ge4b^2$, both exponents lie in
$\Gamma$ and differ by at most $b/\sqrt t$. The modulus and
remainder bounds prove \eqref{window:eq:cw-modulus}. Only uniform
continuity is used; no derivative of an asymptotic error with respect
to the exponent is taken.
\end{proof}

The exact exponent-change identity was derived in
\eqref{cw:tilt}. With the same $e_N$ it reads
\begin{equation}
 Z_N(\lambda)(\dd z)
   =e^{\lambda R_N(z)-\lambda^2/2+e_N(\lambda,z)}
                                  Z_N(0)(\dd z).
 \label{window:eq:cw-tilt}
\end{equation}
It holds also at determinant zeros, where the positive densities
vanish. Together, Lemma~\ref{window:lem:normalizer} and
\eqref{window:eq:cw-two-masses} verify Proposition~\ref{cw:inputs}.

\subsubsection{Exponential mark tails and uniform mass bounds}
\label{window:sec:tails}

Now the goal is to show that two nearby masses, one on each side of criticality, control both exponential mark tails in probability.
The Rayleigh marked law and branching-window profile have the precedents
\cite[Theorem 6.7]{APS} and \cite[Theorem 1.1(ii)--(iii)]{Pain}
discussed in the Introduction. The Gaussian near-critical argument
already reweights critical endpoint marks exponentially and truncates
the resulting unbounded test; see
\cite[proof of Lemma~4.1, (4.3)--(4.4)]{APSCritical} and
\cite[Section~3.2, (3.17)--(3.18)]{PowellReview}.
For the present matrix measures, the tail bound in probability is
proved in \eqref{cw:tail-domination}--\eqref{cw:tails}, using the
uniform normalizer correction and tightness of both endpoint masses. The two endpoint window parameters
have opposite signs, and their absolute value must be larger than the
exponent in the mark test. In the present notation, for $a>0$ and
$b=a+1$ it gives
\begin{equation}
 \Xi_N(e^{a|u|}\ind_{\{|u|>L\}})
       \le C_b e^{-L}[Z_N(b)(K)+Z_N(-b)(K)]
 \label{window:eq:cw-tail-domination}
\end{equation}
for all large $N$. Consequently the uniform marginal estimate
\eqref{window:eq:cw-P} implies, for every $a,\delta>0$,
\begin{equation}
 \lim_{L\to\infty}\limsup_N
 \PP\{\Xi_N(e^{a|u|}\ind_{\{|u|>L\}})>\delta\}=0.
 \tag{T}\label{window:eq:cw-T}
\end{equation}
This is exactly the tail bound in probability in
Corollary~\ref{main:T}; its proof uses endpoint window parameters
of both signs. It asserts no expectation-level exponential uniform
integrability for the mass without the barrier.

Equations~\eqref{cw:weighted-mass} and~\eqref{cw:sup-mass} also
give tightness of $\Xi_N(e^{A|u|})$ and the pathwise estimate
\begin{equation}
 \sup_{|\lambda|\le A}Z_N(\lambda)(K)
       \le C_{A,b}[Z_N(b)(K)+Z_N(-b)(K)],\qquad b>A,
 \label{window:eq:cw-sup-mass}
\end{equation}
for all large $N$. The radius $A$ and the strict margin $b>A$
are fixed before $N$. These are consequences of the exact
reweighting identity; no supremum is moved through an expectation.
The full proof, including the decay factor $e^{-(b-A)L}$,
remains in Section~\ref{cw:subsection}.

\subsubsection{Convergence of the near-critical curve}
\label{window:sec:conditional}

Given the joint marked limit, we next recover the near-critical curve on each fixed compact parameter interval while preserving its dependence on the field.
We state the curve conclusion under the following hypotheses. Let $D_\kappa$ be
the derivative-normalized critical chaos of the specified Gaussian
field and $M_\kappa=\csh D_\kappa$. Suppose
\begin{equation}
 (\Xi_N,\mathcal Y_N)
   \ \Longrightarrow\
   (M_\kappa|_K(\dd z)\rho_R(\dd u),\mathcal Y_\kappa),
 \qquad
 \rho_R(\dd u)=u e^{-u^2/2}\ind_{\{u>0\}}\dd u.
 \tag{M}\label{window:eq:cw-M}
\end{equation}
Here $\mathcal Y_N$ is a prescribed finite vector of original
field tests, centered after spatial integration, and
$\mathcal Y_\kappa$ is the corresponding vector from the same
field generating $D_\kappa$. The marked measures converge weakly
as finite measures, including total mass. Since their spatial
marginals are $Z_N(0)$, this premise also gives joint convergence
with the critical spatial measure by continuous pushforward.

Corollary~\ref{p5:field-corollary} supplies \eqref{window:eq:cw-M}.
It combines the joint conclusion without the barrier in
Corollary~\ref{p34:joint-corollary} with the general
transfer in Proposition~\ref{field:prop:of-joint}. Fixed-barrier
identification alone does not supply this premise.
If a version of the premise initially uses the corrected mark
$R_N^0$, the uniform $O(t^{-1/2})$ mark displacement gives
bounded-Lipschitz distance at most $Ct^{-1/2}Z_N(0)(K)$.
The total mass is tight as part of that premise, so this distance
tends to zero in probability without changing a joint field
coordinate. Thus the conditional statement can use the raw mark.

\begin{theorem}[A continuous near-critical curve]
\label{window:thm:cw-W}
Assume \eqref{window:eq:cw-global} and \eqref{window:eq:cw-M}. For every fixed
$A<\infty$, jointly with the marked measure $\Xi_N$ and the same
field vector, whose joint limit is the one in \eqref{window:eq:cw-M},
\begin{equation}
 (Z_N(\lambda))_{|\lambda|\le A}
   \ \Longrightarrow\
   (F(\lambda)D_\kappa|_K)_{|\lambda|\le A}
 \quad\hbox{in }C([-A,A],\mathcal M_+(K)),             \label{window:eq:cw-W}
\end{equation}
with the supremum bounded-Lipschitz metric, where
\begin{equation}
 F(\lambda)=\csh e^{-\lambda^2/2}
       \int_0^\infty u e^{-u^2/2+\lambda u}\dd u
     =\csh e^{-\lambda^2/2}+2\lambda\Phi_{\rm nor}(\lambda).
 \label{window:eq:cw-F}
\end{equation}
In particular $F(0)=\csh$.
\end{theorem}

\begin{proof}
Fix $A$ and $b>A$. Under \eqref{window:eq:cw-global},
Proposition~\ref{window:prop:cw-P} gives tightness of each
endpoint mass at $b$ and $-b$, hence of their sum by the union
bound. This uses the fixed fraction supplied by that hypothesis;
it does not require the choice $\varepsilon=1/2$.
Lemma~\ref{window:lem:normalizer} supplies the other estimate,
\eqref{cw:normalizer-control}. Therefore the two estimates used
in the argument of Section~\ref{cw:subsection} hold under the present hypotheses.

Apply the proof in Section~\ref{cw:subsection}, beginning with
its truncated map \eqref{cw:truncated-map}, to the joint marked
premise \eqref{window:eq:cw-M}. Include the spatial marginal
$Z_N(0)$ by continuous pushforward as above. The map's continuity,
the uniform approximation \eqref{cw:curve-approximation}, and
the limiting approximation \eqref{cw:limit-approximation} give
the joint conclusion by \eqref{cw:joint-distances} and
\eqref{sp:approximation-proof}. That proof also establishes
finite-$N$ continuity and measurability of the curve. In both
approximations the spatial measure, marked measure and original
field vector are unchanged; only the appended curve is approximated.
The order is $N\to\infty$ at fixed mark cutoff, then removal
of that cutoff. The limiting argument uses only
$D_\kappa(K)<\infty$ almost surely and the Rayleigh law's fixed
exponential moments.

The profile calculation at the end of that same proof is
\eqref{sp:cw:profile}, reproduced as \eqref{window:eq:cw-F}.
Its integral representation is positive for every real $\lambda$,
and $F(0)=\csh$ retains the distinction between
$M_\kappa$ and $D_\kappa$. This proves the conditional conclusion
without requiring the full main theorem as an additional premise.
\end{proof}

\begin{remark}[Boundary of the conclusion]
The window radius and the retained fraction are fixed before $N$.
Both signs of the window are included. A growing $\lambda_N$ or
a fixed supercritical exponent is not covered. The global
no-crossing premise and the joint marked convergence in
\eqref{window:eq:cw-M}, with its specified Gaussian coupling, must be
proved in their separate sections. The deductions here preserve
that joint law and do not create it from marginal convergence.
\end{remark}

\endgroup

\clearpage
\phantomsection
\addcontentsline{toc}{section}{References}
\bibliographystyle{unsrturl}
\bibliography{IMW-references}

\end{document}